\documentclass{amsbook}

\usepackage{amsfonts,amssymb,verbatim,amsmath,amsthm,latexsym,textcomp,amscd}
\usepackage{latexsym,amsfonts,amssymb,epsfig,verbatim}
\usepackage{amsmath,amsthm,amssymb,latexsym,graphics,textcomp}
\usepackage{graphicx}
\usepackage{color}
\usepackage{url}
\usepackage{stmaryrd}
\usepackage{diagrams}
\usepackage{txfonts}

\usepackage[utf8]{inputenc}
\usepackage{glossaries}

\usepackage[T3,T1]{fontenc}
\DeclareSymbolFont{tipa}{T3}{cmr}{m}{n}
\DeclareMathAccent{\invbreve}{\mathalpha}{tipa}{16}

\makeindex

\usepackage{stackengine,scalerel}

\newtheorem{theorem}{Theorem}[chapter]
\newtheorem{thm}[theorem]{Theorem}
\newtheorem*{theorem*}{Theorem}

\newtheorem{prop}[theorem]{Proposition}
\newtheorem{lemma}[theorem]{Lemma}
\newtheorem{lem}[theorem]{Lemma}

\newtheorem{cor}[theorem]{Corollary}
\newtheorem{claim}[theorem]{Claim}

\theoremstyle{definition}

\newtheorem{defn}[theorem]{Definition}
\newtheorem{example}[theorem]{Example}

\theoremstyle{remark}

\newtheorem{rem}[theorem]{Remark}

\newtheorem{question}[theorem]{Question}
\newtheorem{notation}[theorem]{Notation}
\newtheorem{convention}[theorem]{Convention}
\newtheorem{conj}[theorem]{Conjecture}

\numberwithin{section}{chapter}
\numberwithin{equation}{chapter}

\usepackage{remreset}
\makeatletter
\@removefromreset{figure}{chapter}
\@removefromreset{table}{chapter}
\makeatother

\newcommand{\map}{\rightarrow}

\def\H{\mathbb H}
\def\N{\mathbb N}
\def\R{\mathbb R}
\def\D{\partial}

\newcommand{\C}{\mathcal C}
\newcommand{\CC}{\mathcal C}
\newcommand{\CCC}{\mathbb C}

\def\RA{\Rightarrow}

\def\al{\alpha}
\def\La{\Lambda}
\def\si{\sigma}

\def\Om{\Omega}
\def\A{{\mathfrak A}}
\def\Fl{{\mathfrak Fl}}
\def\HHH{{\mathfrak H}}
\def\L{{\mathfrak L}}
\def\QQQ{{\mathfrak Q}}
\def\T{{\mathfrak T}}
\def\X{{\mathfrak X}}
\def\Y{{\mathfrak Y}}
\newcommand\ZZZ{{\mathfrak Z}}

\def\AA{\mathcal A}

\newcommand\FF{{\mathcal F}}
\newcommand\GG{{\mathcal G}}

\newcommand\LL{{\mathcal L}}

\newcommand\NN{{\mathbb N}}

\newcommand\PP{{\mathcal P}}

\newcommand\RR{{\mathbb R}}

\newcommand\XX{{\mathcal X}}
\newcommand\YY{{\mathcal Y}}
\newcommand\ZZ{{\mathcal Z}}

\newcommand\flaring{{Corollary \ref{cor:super-weak flaring}}}
\def\Ga{\Gamma}
\def\Z{\mathbb Z}
\def\Si{\Sigma}

\def\Isom{\operatorname{Isom}}

\def\diam{\operatorname{diam}}
\def\dist{\operatorname{dist}}
\def\hull{\operatorname{Hull}}
\def\id{\operatorname{id}}
\def\Im{\operatorname{Im}}
\def\Hd{\operatorname{Hd}}

\renewcommand\int{\operatorname{int}}

\def\barycenter{\operatorname{center}}

\def\length{\operatorname{length}}

\def\mini{\scriptsize}

\def\acts{\curvearrowright}
\def\embed{\hookrightarrow}

\def\ga{\gamma}
\newcommand\la{\lambda}
\newcommand\eps{\epsilon}
\def\geo{\partial_{\infty}}

\newcommand{\cev}[1]{\reflectbox{\ensuremath{~\vec{\reflectbox{\ensuremath{#1}}}}}}

\makeglossaries

\begin{document}

\frontmatter

\title[Trees of hyperbolic spaces]{Trees of hyperbolic spaces}

\author{Michael Kapovich}
\address{University of California Davis, 1 Shields Avenue, Davis, CA 95616}
\email{kapovich@ucdavis.edu}

\author{Pranab Sardar}
\address{Indian Institute of Science Education and Research, Mohali, Sector 81, PB 140306, India}
\email{psardar@iisermohali.ac.in}

\date{\today}




\maketitle

\setcounter{page}{4}

\tableofcontents


\chapter*{Preface}

The goal of this book is to understand geometry of metric spaces $X$ which have structure of 
trees of hyperbolic spaces. The subject originates in the papers \cite{BF, BF-err} of Bestvina and Feighn, 
where they proved a {\em combination theorem}\footnote{Subsequently, alternative proofs of the group-theoretic version of this theorem were given by Kharlampovich  and Myasnikov in \cite{km-malnormal} and Gautero in \cite{Gautero-2003}, under certain extra assumptions.}, stating that under  certain conditions such $X$ itself is hyperbolic: 

\begin{theorem*}\label{thm:BF}
Suppose $\X= (\pi: X \map T)$ is a tree of hyperbolic metric spaces, where 
vertex and edge-spaces are uniformly hyperbolic, incidence maps of 
edge spaces into vertex spaces are uniformly quasiisometric and which satisfies 
 the hallway flaring condition. Then $X$ is a hyperbolic metric space. 
\end{theorem*}

In Chapter \ref{ch:trees} we give definitions clarifying the result. Informally, the hallway flaring condition means that 
two $K$-quasiisometric sections of $\X$ over a geodesic interval $I$ in $T$ ``diverge at a uniform exponential rate" 
as we move along $I$ in one of the two directions. The original proof of this theorem 
was by verifying that $X$ satisfies the linear isoperimetric inequality. In the book we 
give a new (and longer) proof under weaker flaring assumption than the one made by 
Bestvina and Feighn; we name the weakened condition {\em uniform} (or, in another version, {\em proper}) flaring. 
Informally speaking, instead of requiring the exponential divergence of sections, we only require {\em some} 
rate of divergence, given by a uniform proper function of the arc-length parameter of $I$. We refer the reader to Theorem \ref{thm:mainBF} for the precise statement.

The main benefit of our proof is that it is done by constructing a {\em slim combing} of $X$: We find a family of 
(uniformly quasigeodesic) paths $c(x,y)$ connecting pairs of points $x, y$ in $X$, satisfying {\em slim triangle} 
property: Given three points $x, y, z\in X$, any one of the three paths $c(x,y), c(y,z), c(z, x)$ is contained in 
a uniform neighborhood of the union of the two other paths. The description of the paths $c$ is a 6-step induction 
summarized in Chapter \ref{ch:description-of-geodesics}, starting with  paths in the trees of spaces of the simplest kind that we call 
{\em narrow carpets}: These are metric interval-bundles over geodesics in $T$ such that one of the interval-fibers 
has uniformly bounded length. We hope that our method of proof of Theorem \ref{thm:BF} has a potential to generalize this theorem  to complexes of hyperbolic spaces. 

The combing paths $c$ in $X$ are mostly concatenations of $K$-quasiisometric  sections of $\X$ 
over geodesics in $T$. Thus, we obtain (up to a uniformly bounded error) a description of geodesics in $X$ in terms 
of its structure as a tree of spaces, i.e. vertex-spaces and sections. 

As an application of this description of geodesics, we prove (Theorem \ref{thm:ECT}) 
the existence of Cannon--Thurston maps from 
Gromov-boundaries of subtrees of spaces $Y\subset X$ to $X$, extending an earlier result by  
Mitra \cite{mitra-trees}, who proved 
the existence of  Cannon--Thurston maps for the inclusion maps of vertex-spaces into $X$. Mitra's proof (as well as the subsequent work of Mj and Sardar, \cite{pranab-mahan}) was, in fact, 
a guideline for our description of geodesics in $X$. However, Mitra's description of geodesics stopped at geodesics connecting 
points in the same vertex-space (step 3 of our 6-step description), leaving much of the work to be done in general. Furthermore, we analyze in detail the {\em Cannon--Thurston laminations} of these maps. 

\medskip 
We also refer the reader to the related work of  Gautero describing uniform quasigeodesics in groups obtained via combination theorem in a special case in \cite{Gautero-2003} and in \cite{Gautero-2016} for trees of relatively hyperbolic spaces. However, we were unable to follow Gautero's proof.  
We note, furthermore, that various forms of the Bestvina--Feighn combination theorem for relatively hyperbolic groups and spaces were proven by Dahmani \cite{MR2026551}, Alibegovi\'{c} \cite{MR2131400}, Gautero and Lusztig 
\cite{Gautero-Lustig-2004, Gautero-Lustig-2007} and, in greatest generality, by 
Mj and Reeves \cite{mahan-reeves}.  In the book we did not attempt to describe (quasi)geodesics in the relatively hyperbolic
trees of spaces. However, we proved the existence of Cannon--Thurston maps for subtrees of spaces $Y\subset X$ in Chapter \ref{ch:CTR} using techniques of proof for the existence of Cannon--Thurston maps for subtrees of spaces in hyperbolic trees of spaces. 

\medskip 
{\bf Organization of the book.} In Chapter \ref{ch:basics} we review basic facts of coarse geometry and geometry of hyperbolic spaces. While most of the material of the chapter is standard and well-known, we included it for the ease of reference in 
the rest of the book. 

In Chapter \ref{ch:trees} we discuss  definitions of the theory of trees of metric spaces, state and compare different flaring conditions in trees of spaces, formulate our main theorem and prove it in some easier cases, e.g. for {\em quasiconvex  amalgamations} (Section \ref{sec:qcamalgam}).    
 
In Chapter \ref{ch:4 classes} we define a certain class of subspaces $Y$ in a tree of spaces $\X$, called {\em semicontinuous families}. These subspaces (each of which also has structure of a tree of hyperbolic spaces $\Y$) 
have the property that their intersections with vertex-spaces of $\X$ are uniformly 
quasiconvex and every point in $Y$ is connected to the intersection $Y_u= Y\cap X_u$ of $Y$ with 
a distinguished vertex-space $X_u$, by a $K$-quasiisometric section  of $\X$ over 
an interval in $T$. We prove that the subspaces $Y$ are coarse Lipschitz retracts of $X$, which is a generalization 
of the horocyclic projections to a geodesic in the hyperbolic plane; its existence was first proven by Mitra, \cite[Theorem 3.8]{mitra-trees} 
in the case of semicontinuous families called {\em flow-spaces} $Fl_K(X_u)$. Flow-spaces and three other types of semicontinuous families (ladders, carpets and bundles)  serve as key tools in our definition of combing paths $c$ in $X$. Ladders are certain (semicontinuous) families of intervals over subtrees in $X$, where semicontinuity (informally) 
means that the lengths of the intervals can shrink substantially as we move away one edge from a vertex $u$ (the center of the ladder). Bundles should be thought of as {\em continuous} families of quasiconvex subsets $Q_v$ of vertex-spaces of $\X$ with two (nonempty) vertex-spaces $Q_v, Q_w$ uniformly Hausdorff-close, whenever $v$ and $w$ span an edge of $T$.

Chapter \ref{sect:3} primarily deals with Steps 1--3 of our description of geodesics in $X$: 
We describe combing paths in carpets, carpeted ladders and general ladders and establish their hyperbolicity. 
Hyperbolicity of flow-spaces in proven in Chapter 
\ref{ch:flows}, which is technically the most difficult part of our work: We prove the slim triangle property for the combing paths $c$ 
by analyzing triples of ladders (with the common center $u$) in the flow-space $Fl_K(X_u)$ of a vertex-space $X_u$. 

Our last challenge is to connect by combing paths points in different vertex-flow-spaces $Fl_K(X_u)$, $Fl_K(X_v)$. This is done in 
Chapter \ref{sec:everything together}. The case of points in intersecting flow-spaces $Fl_K(X_u)$, $Fl_K(X_v)$ is handled in 
Section \ref{sec:union-of-two} where we primarily analyze  the case of {\em special intervals} $\llbracket u, v\rrbracket\subset T$, i.e. 
when $Fl_K(X_u)\cap X_v\ne \emptyset$. This covers Step 4 of our description of geodesics in $X$ and  is quite technical. 
The main trick is to introduce a certain generalization of flow-spaces of vertex-spaces and appeal to a special (and easy) case of Theorem \ref{thm:mainBF} proven earlier, the {\em quasiconvex amalgamation}, 
when the tree $T$ contains a single edge (Corollary \ref{cor:edge-spaces}). Once the case of special intervals 
is done, we complete easily Step 5 of our description of geodesics by considering points in flow-spaces 
$Fl_K(X_J)$ for subintervals $J\subset T$ represented as unions of three special subintervals: For the proof we use the quasiconvex amalgamation again. (A good example of such an interval $J$ is given by a {\em semispecial} interval  $\llbracket u, v\rrbracket$, where  the flow-spaces $Fl_K(X_u)$, $Fl_K(X_v)$ have nonempty intersection in $X$.)  Lastly, we conclude the 6-step description of geodesics in $X$ by appealing to the {\em horizontal subdivision} of geodesic intervals $J$ in $T$, 
so that the consecutive subdivision vertices $u_i, u_{i+1}$ define pairwise uniformly cobounded flow-spaces 
$Fl_K(X_{u_i})$, $Fl_K(X_{u_{i+1}})$ (their projections to the tree $T$ are disjoint), 
while each interval $\llbracket u_i, u_{i+1}\rrbracket$ between $u_i, u_{i+1}$ 
is a union of three special subintervals. This uniform coboundedness property allows us to 
reduce the problem of hyperbolicity of the flow-space $Fl_K(X_J)\subset X$ to the {\em pairwise cobounded quasiconvex chain-amalgamation} of hyperbolic spaces which is, again a special and easy case of Theorem \ref{thm:mainBF} proven earlier 
(Theorem \ref{thm:chain}). The combing paths $c(x,y)$ in $X$ are then defined as geodesics in flow-spaces 
$Fl_K(X_J)$, $x\in X_u, y\in X_v$, and $J=\llbracket u, v\rrbracket$.  
 Lastly, we verify the slim combing property for these paths $c(x,y)$ by considering flow-spaces $Fl_K(X_S)$ for 
 geodesic tripods $S\subset T$ and appealing to Theorem \ref{thm:chain} (or, more precisely, its consequence, Corollary 
 \ref{cor:finite-tree-hyp}) one last time. 

In Chapter \ref{ch:description-of-geodesics} we review the description of the combing paths $c(x,y)$ by putting together different steps of the descriptions scattered in the earlier parts of the book. We also prove an easy application of this description by giving  a 
 characterizations of geodesics $\al$ in vertex-spaces of $X$ 
which are quasigeodesics in $X$ itself, in terms of carpets bounded by subintervals in $\al$. Furthermore, assuming acylindricity, we give a simplified description of uniform quasigeodesics and quasigeodesic rays in $X$. We use this description to describe the ideal boundary of $X$ in terms of ideal boundaries of vertex spaces and of the tree $T$.  

In Chapter \ref{ch:CT} we apply our description of geodesics to prove Theorem \ref{thm:ECT}, establishing existence of 
Cannon--Thurston maps for  subtrees of spaces $\Y$ in a hyperbolic tree of spaces $\X$. The main technical result 
of the chapter is Theorem \ref{thm:cut-paste} which relates quasigeodesics $\phi$ in $X$ to quasigeodesics in $Y$ via a certain {\em cut-and-replace} procedure, replacing {\em detour subpaths} in $\phi$ by geodesics in vertex-spaces of $Y$. Along the way, we relate nearest-point projections to flow-spaces $Fl_K(X_u)$ taken in $X$ and in $Y$.  
We  give a necessary and sufficient condition for 
points to have equal images under these Cannon--Thurston maps (Theorem \ref{thm:CT-fibers}). In particular, such points have to belong to the Gromov boundary of the same vertex-space of $\Y$. In the following four sections of the chapter  
we discuss {\em Cannon--Thurston laminations} for the inclusion maps $X_v\to X$ in more detail. 
In the last section of the chapter, we discuss group-theoretic applications of our results. In particular, we construct examples of 
non-Anosov undistorted surface subgroups of $PSL(2, {\mathbb C})\times PSL(2, {\mathbb C})$ 
consisting entirely of semisimple elements.

In Chapter \ref{ch:CTR} we consider trees of relatively hyperbolic spaces and, generalizing the results of 
Chapter \ref{ch:CT},  prove existence of Cannon--Thurston maps in this context and establish some properties of the associated Cannon--Thurston laminations. 

There are many constants and functions used in the book. As a general rule, we label these using as the subscript the number of the theorem (or lemma, etc.) where these quantities 
are introduced.

\medskip
{\bf Acknowledgements.} During the work on this book the first author was partly supported by the NSF grant  DMS-16-04241 and by a Simons Fellowship, grant number 391602.

\mainmatter

\chapter{Preliminaries on metric geometry}\label{ch:basics}

\section{Graphs and trees}

Although we always work with unoriented metric graphs like Cayley graphs, we will also need  oriented graphs, 
to describe graphs of groups. The following definition 
is taken from \cite{serre-trees}:

\begin{defn}\index{graph} 
An {\em oriented graph} $\Gamma$ is a pair of sets $(V,E)$ together with two maps 
$$
E\map V\times V,\quad e\mapsto (o(e),t(e))$$ 
and
$$
E\map E, \quad e\mapsto \bar{e}$$
such that $o(\bar{e})=t(e), t(\bar{e})=o(e)$ and $\bar{\bar{e}}=e$ for all $e\in E$.
\end{defn}

We write $V(\Gamma)$  for $V$ and $E(\Gamma)$  for $E$. We refer to $V(\Gamma)$ as the {\em set of
vertices} of $\Gamma$ and $E(\Gamma)$ as the {\em set of edges} of $\Gamma$. We will almost always conflate a graph $\Gamma$ 
with its underlying space, i.e. its geometric realization as a 1-dimensional CW complex.

For an edge $e$ of a graph we refer to $o(e)$ as the {\em origin} and $t(e)$ as the {\em terminus} of $e$; the edge $\bar{e}$ is the same edge $e$ with opposite orientation. When $o(e)=v, t(e)=w$, we will use the notation $e=[v,w]$. While for general graphs this notation is ambiguous, for graphs which are  trees (and this is the case we are mostly interested in),  vertices $v, w$ uniquely determine the oriented edge $e$.

We shall denote by $|e|$ the edge $e$  of $\Gamma$ without any orientation, regarded as a subset of the underlying space of $\Gamma$. For each edge $e$ in a graph we define $\dot{e}$ to be $|e|$ with the end-points removed.

Given a subset $W\subset V$, we define the  {\em full subgraph} of $\Gamma$ 
spanned by $W$ as the maximal subgraph in $\Gamma$ with the vertex-set $W$. The {\em valence} or {\em degree} of a vertex $v\in V$ 
is the cardinality of the set $o^{-1}(v)\subset E$ (equivalently, $t^{-1}(v)\subset E$).

If $\Lambda\subset \Gamma$ is a subgraph, then a vertex $v\in V(\Lambda)$ is a {\em boundary vertex} 
 of $\Lambda$, if there is an edge $e=[v,w]\in V(\Gamma)$ such that $w\notin V(\Lambda)$.  The edge $e$ is then called a {\em boundary edge} of $\Lambda$ in $\Gamma$. We will use this notion only when $\Lambda$ is a subtree of a tree $\Gamma$. \index{boundary edge of a subtree} \index{boundary vertex of a subtree}

A {\em graph-morphism}, \index{graph-morphism} 
or a {\em morphism of graphs} $\phi: \Gamma\to \Gamma'$ is a pair of maps $\phi_V: V(\Gamma)\to V(\Gamma'), \phi_E: E(\Gamma)\to E(\Gamma')$, $v\mapsto v', e\mapsto e'$  such that the following diagrams commute for all oriented edges $e=[v,w]$:  

$$
 \begin{diagram}
 v              & \lTo & e & \rTo & w           \\
 \dTo_{\phi_V}   &     & \dTo_{\phi_E}    & &   \dTo_{\phi_V} \\
 v'              & \lTo & e' & \rTo & w'          
  \end{diagram}
$$
where the horizontal arrows are the origin/tail maps.

A {\em tree} is a simply-connected graph. \index{tree} 

\medskip 
We will use both the notation $uv$ and $\llbracket u,v\rrbracket$ for the (geodesic) segment, or {\em an interval}, in $\Gamma$ whose end-points are the vertices $u, v$. 
(Since we will be mostly working with graphs which are trees, this notation is unambiguous.) 
Given a segment $\llbracket u,v\rrbracket$ in a tree, we define 
$\rrbracket u, v\llbracket $ as the maximal subtree of $\llbracket u,v\rrbracket$ containing 
all the vertices of $\llbracket u, v\rrbracket $ except for $u$ and $v$. Similarly, we define 
subsegments $\rrbracket u, v\rrbracket $ and $\llbracket u, v\llbracket $. \index{intervals in trees}

\begin{convention}
\label{conv:order}
We will regard intervals $\llbracket u,v\rrbracket$ in simplicial trees as ordered sets with 
$u$ the smallest element and $v$ the largest. Accordingly, we will talk 
about supremums and infimums of subsets of $\llbracket u,v\rrbracket$ 
and $\sup(\emptyset)= u, \inf(\emptyset)=v$.  
\end{convention}

A {\em metric graph} is a connected graph $\Gamma$, every edge $e$ of which is assigned a positive real number $\ell(e)$ (its {\em length}). 
The vertex-set of $\Gamma$ then has a natural pseudometric $d_\ell$, where the distance between vertices  is defined to be the infimum of total lengths of edge-paths connecting these vertices. The metric $d_\ell$ extends to a pseudometric on the underlying space of $\Gamma$. Note, however, that the distance between the vertices of an edge $e$ of $\Gamma$ can be smaller than $\ell(e)$ even if the vertices are distinct. If $\ell: E(\Gamma)\to \RR_+$ takes only finitely many values, then $(\Gamma, d_\ell)$ is a complete geodesic metric space. 
If the function $\ell$ is bounded away from $0$, then $d_\ell$ is a metric, but, 
in general, metric graphs need not be complete nor geodesic and the pseudometric need not be a metric.

\begin{example}
Let $\Gamma$ be a graph with two vertices $v, w$ and edges $\{e_i\}_{i\in \NN}$ all of which connect $v$ to $w$. 

1. Take the function $\ell(e_i)= \frac{1}{i}$. Then $d_\ell(v,w)=0$. Hence,  $d_\ell$ is not a metric in this example. 

2. Take the function $\ell(e_i)= 1+ \frac{1}{i}$. Then $d_\ell$ is a metric but $\Gamma$ contains no geodesics between $v$ and $w$. 
\end{example}

\begin{example}
Consider the graph $\Ga$, which is the complete graph on the set $\{v_1, v_2, v_3\}$. Let $\ell([v_1,v_2])=\ell([v_2,v_3])=1, \ell([v_3,v_1])= 3$. Then 
$d_\ell(v_1, v_3)$ is $2$ rather than $3=\ell([v_3,v_1])$. 
\end{example}

Most of the time, unless stated otherwise, we will metrize connected graphs $\Ga$ by declaring that every edge has unit length: The distance between 
vertices equals the minimal number of edges in an edge-path connecting these vertices.  
 We refer to the resulting metric on $\Ga$ as the {\em graph-metric}.

\section{Coarse geometric concepts}

\subsection{Metric notions} \label{sec:Metric notions}

For a subset $A$ of a topological space $X$, $cl(A)$ will denote the closure of $A$ in $X$. By a {\em path} in a topological space $X$ we will always means a continuous map $I\to X$, where $I$ is an interval in $\RR$. Given a path $c: [a,b]\to X$, we denote by $\cev{c}$ the reverse path
$$
c(t)= c(a+b-t). 
$$

A path $c$ in a metric space $(X,d)$ is {\em geodesic} if it is an isometric embedding $I\to (X,d)$. 
We will frequently conflate paths and their images: 
Since we are primarily interested in geodesic and quasigeodesic paths, this conflation is mostly harmless. Accordingly, 
if $x, y$ are points in a path $c$ in $X$, then $c(x,y)$ will denote the subpath of $c$ between $x$ and $y$. This notation is, of course, slightly ambiguous since $c$ need not be injective, and a better notation would have been $c|_{[s,t]}$ where $c(s)=x, c(y)=t$. However, in practice, it will be always clear what the subpath $c(x,y)$ is. 
We will use the notation $c_1\star c_2$ to denote the concatenation of two paths. 

The {\em length} of a path $c: I\to (X,d)$, where $I=[a,b]$ is a finite closed interval in $\RR$, is
$$
\length(c)= \sup \sum_{i=1}^n d(c(t_i), c(t_{i+1})), 
$$ 
where the supremum is taken over all subdivisions of interval $I$:
$$
t_1=a\le t_2 \le ... \le t_n \le t_{n+1}=b. 
$$
A metric space $(X,d)$ is called {\em rectifiably connected} if every two points in it are connected by a path of finite length. 
A metric space $(X,d)$ is called a {\em path-metric space} and the metric $d$ a {\em path-metric} if for all points $x, y\in X$, \index{path-metric}
$$
d(x,y)= \inf_{c} \length(c),
$$
where the infimum is taken over all paths in $X$ connecting $x$ and $y$. Examples of path-metric spaces are given by metric graphs (assuming, of course, 
that $d_\ell$ is a metric and not merely a pseudometric). 

A metric space $(X,d)$ is called {\em geodesic} if for every two points $x, y\in X$ are connected by a geodesic path. 
We will use the notation  $xy$, or  $[xy]_X$, to denote a geodesic segment joining $x$ to $y$ in $X$. 
If $X$ is a simplicial tree and $u, v$ 
are vertices, then we will also use the notation $\llbracket u, v \rrbracket$ for this geodesic segment. 

For $x, y, z\in X$ we shall denote by $\Delta xyz$ a geodesic triangle with vertices $x, y, z$ which is the union 
of three geodesic segments $xy \cup yz \cup zx$. Similarly, a {\em geodesic quadrilateral} in $X$ with vertices $x, y, z, w$, 
denoted $\square xyzw$, is the union of four geodesics
$$
xy \cup yz \cup zw \cup wx. 
$$  
For any rectifiably--connected  subset $Y$ in a metric space $(X,d)$ we shall denote by $d_Y(\cdot , \cdot)$ the path-metric on $Y$ induced from $X$: The distance between 
two points in $Y$ is the infimum of lengths of paths in $Y$ between these points, where the length of a path is computed using the restriction of the metric $d$.

For $R\geq 0$ and a 
subset $A\subset X$,  
$$
N^X_R(A)=N_R(A):=\{x\in X: \, d(x,a)\leq R \, \mbox{ for some} \, a\in A\}$$ 
will denote the (closed) $R$-neighborhood of $A$ in $X$. A subset $A\subset X$ is 
said to be an {\em $R$-net} in $X$ if 
$$
N_R(A)=X. 
$$
\index{net}

\medskip

For subsets $Y, Z$ in a metric space $X$, $\Hd(Y, Z)\in [0,\infty]$ denotes the 
{\em Hausdorff distance} between $Y$ and $Z$:
$$
\Hd(Y, Z)= \inf \{R: Y\subset N_R(Z), Z\subset N_Z(Y)\}. 
$$
 We will use the notation 
$$
d(Y, Z)= \inf \{ d(y, z): y\in Y, z\in Z\}
$$
for the {\em minimal distance} between $Y$ and $Z$. (Note that, unlike the Hausdorff distance, the minimal distance, in general, fails to satisfy the triangle inequality.) We will sometimes add the subscript $X$ in this notation to emphasize that the distances and neighborhoods are taken in $X$.

For two maps $f, g: X\map Y$ between metric spaces, we define the distance between $f, g$ as 
$$
d(f, g)=\sup \{d(f(x),g(x)): x\in X\}.$$

\section{Group actions}\label{sec:Group actions}

Throughout the book, we will be only considering {\em left group actions}  of groups on sets (the notation for such an action is $G\times X\to X$ or $G\acts X$). For instance, if $X$ is a group and $G$ is a subgroup of $X$ then the action of $G$ on $X$ via left-multiplication
$$
L_g(x)= gx
$$
is a left action of $G$ on $X$.

Given a $G$-action on a set $X$ and a point $x\in X$, one defines the {\em orbit map} \index{orbit map} 
for the action to be the map $o_x: G\to X$, by $o_x(g)=gx$. We will be primarily interested in isometric  actions of discrete groups on metric spaces. Such an action is said to be {\em metrically\footnote{we will omit the adjective ``metrically'' in what follows} proper} if 
\index{metrically proper group action} 
for each bounded subset $B\subset X$ the subset 
$$
\{g\in G: gB\cap B\ne \emptyset\} 
$$
is finite. In other words, preimages of bounded subsets under orbit maps are finite. 
An isometric action is said to be {\em cobounded} \index{cobounded group action} 
if there exists a bounded subset $B\subset X$ such that $GB=X$.   An action is said to be {\em geometric}\index{geometric group action}  if it is both proper and cobounded. 

Suppose that we are given an isometric action $G\acts X$ and a subset $Y\subset X$. The {\em stabilizer} of $Y$ in $G$, denoted $G_Y$, is the subgroup of $G$ consisting of elements preserving $Y$ set-wise. 

\begin{defn}\label{defn:locally finite action} \index{locally finite orbit} 
One says that the {\em $G$-orbit} $GY$ of $Y$ is {\em locally finite} if for each $x\in X$ and $r\in \R_+$, there exist a finite subset $\{g_1,...,g_n\}\subset G$ such that 
$$
gY\cap B(x,r)\ne \emptyset \RA g\in g_i G_Y
$$
for some $i=1,...,n$.
\end{defn} 

 In order to see that this condition is natural, observe that for $h\in G_Y$,
$$
gY\cap B(x,r)\ne \emptyset \iff ghY \cap B(x,r)\ne \emptyset. 
$$

\begin{lemma}\label{lem:geo->lf} 
Suppose that $X$ is a finitely-generated group equipped with the word-metric and $Y< X$ is a subgroup. Then 
for each $G< X$, we have $G_Y=G\cap Y$ and the $G$-orbit $GY$ is locally finite.
\end{lemma}
\proof It suffices to prove the claim with $x=1$. Since the ball $B(1,r)$ is finite, there exist a finite set of pairs 
$(g_i, y_i)$, $i=1,...,n$, $g_i\in G, y_i\in Y$, such that whenever $y\in Y$, $g\in G$ satisfy 
$d_X(gy, 1) \le r$, we have $gy=g_iy_i$ for some $i$. Then
$$
h= g_i^{-1}g= y_i y^{-1}\in G\cap Y= G_Y
$$ 
and, hence, $g=g_ih$, i.e. $g\in g_iH$, as required. \qed

\begin{cor}
Suppose that $X$ is a geodesic metric space, $G'\acts X$ is a geometric action, $Y\subset X$ is a nonempty 
subspace whose $G'$-stabilizer $G'_Y$ also acts geometrically on $Y$. Then for each subgroup $G< G'$, the 
$G$-orbit $GY$ is locally finite. 
\end{cor}

\section{Length structures and spaces} \label{sec:length structures} 

Let $X$ be a topological space. A {\em length structure} on $X$ \index{length structure} 
 is a collection ${\mathcal P}$ of {\em admissible paths} 
(defined on closed intervals in $\R$) in $X$, together with a map
$$
\ell: {\mathcal P}\to \R_+
$$
(called a {\em length function}) satisfying the following axioms:

1. $\PP$ is closed under restrictions: The restriction of a path $c\in \PP$ to a subinterval again belongs to $\PP$. 

2. $\PP$ is closed under concatenations. 

3. $\PP$ is closed under linear reparameterizations.

4. $\ell(c_1\star c_2)= \ell(c_1) + \ell(c_2)$. 

5. $\ell$ is invariant under linear reparameterizations. 

6. For each $c: [A, B]\to X$, $c\in \PP$, the length $\ell(c|_{[a,b]})$ depends continuously on $a$. 

7. The length function $\ell$ is consistent with the topology of $X$ in the sense that for each $x\in X$ and each neighborhood $U$ of $X$
$$
\inf_{c\in \PP(x, X-U)} \ell(c) >0
$$
where $\PP(x, X- U)$ consists of all paths $c\in \PP, c: [a,b]\to X$, $c(a)=x, c(b)\in X-U$.   

8. For each pair of points $x, y\in X$ the subset $\PP_{x,y}$ consisting of paths $c\in \PP$ connecting $x$ to $y$ is nonempty. 

\medskip 
A {\em length space}\index{length space}
 is a topological space equipped with a length structure. Each length space $(X, \PP, \ell)$ has a canonical metric $d=d_\ell$ defined by 
$$
d(x,y)= \inf_{c\in \PP_{x,y}} \ell(c). 
$$
The topology defined by this metric is finer than the one of $X$; the metric $d_\ell$ is a path-metric (see Proposition 2.4.1 in \cite{BBI}).

\section{Coarse Lipschitz maps and quasiisometries}

Below, we let $X$, $Y,  Z$ denote metric spaces and let $L\geq 1, \,  \epsilon \geq 0$.

\begin{enumerate}
\item Suppose $Z$ is a set. A map $f: Z \rightarrow Y$ is said to be $D$-{\em coarsely surjective} \index{coarsely surjective map}
if $Y=N_D(f(Z))$, i.e. $f(Z)$ is a $D$-net in $X$.

\item Suppose $\{Z_{\alpha} \}$ and $\{ Y_{\alpha} \}$ are, respectively, a family of sets and a family of metric spaces.
A family of maps $f_{\alpha}: Z_{\alpha}\rightarrow Y_{\alpha}$ is said to be {\em uniformly coarsely surjective} 
if there is a constant $D\geq 0$, such that for all $\alpha$,  $Y_{\alpha}=N_D(f_{\alpha}(Z_{\alpha}))$. \index{coarsely surjective map}

\item A map $f: X\rightarrow Y$ is said to be  $(L,\epsilon)$-{\em coarsely Lipschitz} (or {\em coarse Lipschitz}) 
\index{coarsely Lipschitz}
if $\forall x_1, x_2\in X$ we have
$$d_Y(f(x_1), f(x_2))\leq L d_X(x_1,x_2) + \epsilon.$$ 
A map $f$ is coarsely Lipschitz if it is 
$(L,\epsilon)$-coarsely Lipschitz for some $L \geq 1, \epsilon\ge 0$. 
When $L=\epsilon$, we say that $f$ is $L$-{\em coarsely Lipschitz}. \index{coarsely Lipschitz map}

\item Let $\eta: S\subset \RR_+\to \RR_+$ be a  function. 
A map of metric spaces $f: X\to Y$ is called {\em $(\eta,L)$-proper} 
if $f$ is $L$-coarsely Lipschitz and 
$d(f(x_1), f(x_2))\le R$ implies that $d(x_1, x_2)\le \eta(R)$. The function $\eta$ is a {\em distortion function} of $f$.
\index{distortion function}  
We will frequently suppress the coarse Lipschitz constant $L$ (it will be often equal to $1$) and simply say that  
{\em $f$ is $\eta$-proper}. For instance, if $Y\subset X$ is a rectifiably connected subset of a path-metric space 
$(X, d_X)$, we say that $Y$ is $\eta$-properly embedded in $X$ if the inclusion map $(Y, d_Y)\to (X, d_X)$ is 
$\eta$-proper. \index{$\eta$-proper map} 

\item Similarly, suppose that $f_\al: (X_{\alpha}, d_{X_{\alpha}}) \to (Y_{\alpha}, d_{Y_{\alpha}})$, is a family of maps between 
metric spaces. If these maps are $(\eta,L)$-proper for some 
$\eta$ and $L$, then we will say that this family of maps is {\em uniformly proper}. \index{uniformly proper map}

\item  
A map $f: X\rightarrow Y$ is said to be
an $(L, \epsilon)$-{\em quasiisometric embedding} if   $\forall x_1, x_2\in X$ one has
$$ \frac{1}{L}d_X(x_1,x_2)  - \epsilon\leq d_Y(f(x_1), f(x_2))\leq L d_X(x_1,x_2) + \epsilon.$$
\index{quasiisometric embedding} 

\item 
A map $f: X\rightarrow Y$ will  be simply referred to as
a {\em quasiisometric embedding} 
if it is an $(L, \epsilon)$-quasi\-isometric embedding for some $L\geq 1$ and $\epsilon\geq 0$.

\item 
An $(L,L)$-quasiisometric embedding will  be referred to as an $L$-{\em quasiiso\-met\-ric embedding}. 

\item A map $f:X\rightarrow Y$ is said to be a $(L,\epsilon)$-{\em quasiisometry}\index{quasiisometry}
 if it is an $(L,\epsilon)$-quasiisomet\-ric embedding, which is also $\epsilon$-coarsely surjective. 
 If $L=\epsilon$ then we will refer to such $f$ as an  $L$-{\em quasiisometry}.

\item We will use the abbreviation {\em qi} for the word  {\em quasiisometric}. 

\item An $(L,\epsilon)$-{\em quasigeodesic} (resp. an $L$-{\em quasigeodesic}) in a metric space $X$ is 
a $(L,\epsilon)$-quasiisometric embedding (resp. a $L$-quasiisometric embedding) $\gamma: I\rightarrow X$, where
$I\subseteq \mathbb R$ is an interval.\index{quasigeodesic}

\item Given two maps $f: X\rightarrow Y$ and $g: Y\rightarrow X$, we say that $g$ is an $\epsilon$-{\em coarse left inverse} 
of $f$ if $d(f\circ g, \id_Y)\leq \epsilon$. Similarly one defines an $\epsilon$-{\em coarse right inverse}. 
If $g$ is both $\epsilon$-coarse left and right inverse then it is called an $\epsilon$-{\em coarse inverse} of $f$.

\item If $A\subset X$ and $i: A\map X$ is the inclusion map, then an $(L,\epsilon)$-{\em coarse retraction} of $X$ to $A$
\index{coarse retraction}
is a $(L,\epsilon)$-coarsely Lipschitz map $g:X\map A$ such that $g|_A=\id_A$. 
\end{enumerate}

\begin{rem}
More generally, one can define  an $(L,\epsilon)$-{\em coarse retraction} by requiring that 
$$
d(\id_A, g\circ i)\le \epsilon. 
$$
However,  in the book we only use the more restrictive definition. 
\end{rem}

\begin{example}
1. Let $G, H$ be finitely-generated groups equipped with word metrics and $\phi: G\to H$ is a homomorphism. Then $\phi$ is a coarse Lipschitz map. If $\phi$ has finite kernel, then it is also uniformly proper. 

2. Suppose that $G$ is a finitely-generated (discrete) group, $G\acts X$ is a proper isometric action on a metric space, then for each $x\in X$, the orbit map $o_x: G\to X$ is uniformly proper. 
\end{example}

\begin{defn}
One says that a finitely-generated subgroup $G$ of a finitely generated group $H$ has distortion at most $\eta$ if the inclusion map $G\to H$ is $\eta$-proper (when $G, H$ are equipped with word metrics as above). \end{defn} 

Thus, a subgroup $G$ of a group $H$ is at most {\em linearly distorted} if and only if it is qi embedded in $H$. In this case, one says that $G$ is {\em undistorted} in $H$.  One can make the notion of distortion independent of a generating set by working with a suitable equivalence relation on distortion functions.  For instance, one can talk about polynomial distortion, exponential distortion, etc. We refer the reader to \cite{Drutu-Kapovich} for further details.

\medskip 
We next discuss quasiisometries and qi embeddings of metric spaces. 
The following lemma is a direct calculation which we omit:

\begin{lemma}\label{lem:compositions} 
1. Suppose that $f_1: X_1\to X_2$ and $f_2: X_2\to X_3$ are, respectively, $(L_1,\eps_1)$ and 
$(L_2,\eps_2)$-coarse Lipschitz.  Then their composition 
is $(L_1L_2, L_2\eps_1+\eps_2)$-coarse Lipschitz. 

2. Suppose that $f_1: X_1\to X_2$ and $f_2: X_2\to X_3$ are, respectively, $(L_1,\eps_1)$ and $(L_2,\eps_2)$--qi embeddings.  Then their composition 
is an $(L_1L_2, L_2\eps_1+\eps_2)$--qi embedding. 
\end{lemma}

\begin{lemma}\label{lem:quasi-inverse}
Let $f: X\to Y$ be an $L$-quasiisometry. Then $f$ admits a coarse $3L^2$-inverse which is a $3L^2$-quasiisometry  $Y\to X$.  
\end{lemma}
\proof For $y\in Y$ define $g(y)=x$ such that $d(y, f(x))\le L$. Then
$$
L^{-1} d(y,y') - 3 \le d(g(y), g(y'))\le L d(y,y') +3L^2 
$$
and
$$
d(f\circ g(y), y)\le L, d(g\circ f(x), x)\le 2L^2. \qed 
$$

 \medskip
 
\begin{lemma}\label{lem:close->qi} 
Let $f_i: X_i\to Y$ be $k$-qi embeddings such that 
$$
\Hd(\Im(f_1), \Im(f_2))\leq r.
$$  
Define a map 
$$
g: X_1\to X_2
$$
sending $x_1\in X_1$ to a point $x_2\in X_2$ such that $d(f_1(x_1), f_2(x_2))\le r$. Then $g$ is a $K= K_{\ref{lem:close->qi}}(r,k)$-quasiisometry. 
\end{lemma}
\proof Let $x_1,y_1 \in X_1$. Then 
$$
-k+\frac{1}{k}d_{X_1}(x_1,y_1)\leq d_Y(f_1(x_1), f_1(y_1))\leq k+kd_{X_1}(x_1,y_1).
$$
Setting $x_2:=g(x_1)$ and $y_2:=g(y_1)\in X_2$, we obtain
$$
d_Y(f_1(x_1), f_2(x_2))\leq r \quad \hbox{ and } \quad d_Y(f_1(y_1), f_2(y_2))\leq r.$$
It follows that
$$
|d_Y(f_1(x_1), f_1(y_1))-d_Y(f_2(x_2), f_2(y_2))|\leq d_Y(f_1(x_1), f_2(x_2))+d_Y(f_1(y_1), f_2(y_2))\leq 2r.
$$
Hence, we get 
$$
-k-2r + \frac{1}{k}d_{X_1}(x_1, y_1)\leq d_Y(f_2(x_2), f_2(y_2))\leq 2r+k+kd_{X_1}(x_1, y_1).
$$
Since $f_2$ is a $k$-qi embedding, we have 
$$
-k+\frac{1}{k} d_{X_2}(x_2, y_2)\leq d_Y(f_2(x_2), f_2(y_2))\leq k+kd_{X_2}(x_2, y_2).
$$
Using these two sets of inequalities we obtain 
$$-\frac{2r+2k}{k} +\frac{1}{k^2} d_{X_1}(x_1, y_1)\leq d_{X_2}(x_2, y_2)\leq 2k^2+2rk+k^2d_{X_1}(x_1, y_1).
$$ 
Since $g(x_1)=x_2, g(y_1)=y_2$, it follows that $g$ is a $(2rk+2k^2)$-qi embedding.

Also, given any $x'_2\in X_2$, there is an $x'_1\in X_1$ such that $d_Y(f_1(x'_1), f_2(x'_2))\leq r$.
If $x''_2=g(x'_1)$. Then $d_Y(f_1(x'_1), f_2(x''_2))\leq r$. Hence, 
$$
d_Y(f_2(x''_2),f_2(x'_2))\leq d_Y(f_1(x'_1), f_2(x'_2))+ d_Y(f_1(x'_1), f_2(x''_2))\leq 2r.
$$ Since $f_2$ is
a $k$-qi embedding it follows that $d_{X_2}(x'_2, x''_2)\leq 2rk+k^2< 2rk+2k^2$.
Hence $g$ is a $K=(2rk+2k^2)$-quasiisometry. \qed

\begin{lemma}\label{lem: qi from lipschitz}
Suppose $f: X\map Y$ and $g: X\map Y$ are $(L,\epsilon)$-coarsely Lipschitz maps between metric spaces
such that $d(\id_X, g\circ f)\leq R$ and $d(\id_Y, f\circ g)\leq R$. Then $f$ as well as $g$
is an $(L,\epsilon+2R)$-quasiisometry.
\end{lemma}
\proof The proof of this lemma follows easily from  definitions. We refer to \cite[Lemma 1.1]{pranab-mahan}
for details. 
\qed

\medskip 
The next lemma follows immediately from definitions:

\begin{lemma}\label{qc plus prop emb}
Suppose that $(Y,d_Y)$ is a metric space, $X\subset Y$ is a subset equipped with a metric $d_X$ such that the inclusion map 
$(X,d_X)\to (Y,d_Y)$ is $L$-coarse  Lipschitz and admits an $L$-coarse Lipschitz retraction $(Y,d_Y)\to (X,d_X)$. 
Then the inclusion map $X\to Y$ is an $L$-qi embedding.  
\end{lemma}

\begin{rem}\label{rem:category}
1. Few categorical remarks are in order at this point. It is natural to consider the {\em coarse category} $\CC$, where objects are metric spaces and morphisms are equivalence classes of coarse Lipschitz maps (more generally, correspondences). Here two maps are declared to be equivalent if they are bounded distance apart. Isomorphisms in this setting are precisely quasiisometries of metric spaces. {\em Monomorphisms} or {\em monic morphisms} in this category are {\em uniformly left-cancellative} morphisms, meaning that $f: X\to Y$ is monic if for each pair of $(L,\eps)$-coarse Lipschitz maps $g_i: Z\to X, i=1,2$, if $d(f\circ g_1, f\circ g_2)\le D$ then 
$d(g_1, g_2)\le C(L,\eps, D)$. (Note the need for the uniform control on distances!) It is easy to verify that monic morphisms are precisely (equivalences classes of) uniformly proper maps; hence, monomorphisms are more general than qi embeddings.   
Epimorphisms are precisely the coarsely surjective maps. Most important examples of these, besides quasiisometries, are given by coarse retractions to subsets of metric spaces. The coarse retractions frequently used in the book 
are {\em Mitra's projections} (in the setting of subtrees of hyperbolic spaces) and nearest-point projections to 
quasiconvex subsets of hyperbolic spaces.  

2. An even more general formalism of coarse structures is developed by John Roe in \cite{Roe}. 

3. The above categorical notions, unfortunately, are not quite satisfactory for our purpose, since most of the time we have to keep track of various quantities such as distances between equivalent maps, coarse Lipschitz constants and distortion functions. For instance, when defining a graph (even a tree) of metric spaces, it is not quite enough to say that this is a functor from a graph (regarded as a category) to the category $\CC$, sending origin/terminus maps to monic morphisms of metric spaces: We will need uniform control of coarse Lipschitz constants and distortion functions. For a tree of hyperbolic spaces, we will need even more control, bounding hyperbolicity constants. For this reason, we will adopt a more pedestrian (and traditional) approach, and mostly refrain from using the categorical language. 
\end{rem}

\begin{lemma}\label{lem:qi-net} 
Let $Y$ be a path-metric space, $X\subset Y$ is an $\eta$-properly embedded 
subset equipped with the induced path-metric such that $X$ is an $R$-net in $Y$. Then the inclusion map $\iota: X\to Y$ is an $L$-qi embedding with $L= \eta(2R+1)$. 
\end{lemma}
\proof Take $x, x'\in X$, let $c=c_\eps: I=[0, T]\to Y$ be an arc-length parameterized path in $Y$ connecting $x$ to $x'$ whose length $T$ is $\le d(x,x')+\eps$. Subdivide the interval $I$ into $n+1$ subintervals $[t_i, t_{i+1}]$, $t_0=0$, such that $t_{i+1}-t_i=1$ except for $i=n$, $0\le r=t_{n+1}-t_n< 1$.  Let $P: Y\to X$ be a nearest-point projection. We apply $P$ to the sequence of points $y_i=c(t_i)$ and get a sequence $x_0=x, x_1=P(y_1),...,x_n=p(y_n), x_{n+1}=x'$ such that
$$
d_Y(x_i, x_{i+1})\le 2R+1, i=0,...,n. 
$$
Hence, $d_X(x,x')\le L(n+1)$, $L= \eta(2R+1)$. Therefore,
$$
d_X(x,x')\le L (d_Y(x,x')+1), 
$$
i.e. $\iota$ is an $L$-qi embedding. \qed 

\begin{rem}
Of course, in this situation, the map $\iota$ is also a quasiisometry $X\to Y$. 
\end{rem}

\begin{lemma}\label{close-spaces-qi}
Given a function $\eta: \RR_{+}\map \RR_{+}$ and $R\geq 0$, there is a constant $K=K_{\ref{close-spaces-qi}}(\eta, R)$
such that the following holds.

Suppose $Z$ is a  metric space and $Z_1, Z_2$ are two rectifiably-connected subsets in $Z$ such that, 
both $Z_1$ and $Z_2$ are $\eta$-properly embedded in $Z$. Assume $\Hd(Z_1,Z_2)\leq R$
and suppose $f: Z_1\map Z_2$ is any map that such that $d_Z(z,f(z))\leq R$ for all $z\in Z_1$.
Then $f$ is a $K$-quasiisometry.
\end{lemma} 

\proof  Since $\Hd(Z_1,Z_2)\leq R$ clearly there is a similar map $g: Z_2\map Z_1$.
We note that $d_Z(z, g\circ f(z))\leq 2R$ for all $z\in Z_1$. Hence, $d_{Z_1}(z, g\circ f(z))\leq \eta(2R)$
for all $z\in Z_1$, since $Z_1$ is $\eta$-properly embedded in $Z$, i.e. $d(\id_{Z_1}, g\circ f)\leq \eta(2R)$. 
Similarly,  $d(\id_{Z_2}, f\circ g)\leq \eta(2R)$.

We claim that $f$, $g$ are coarsely Lipschitz. Since the proofs are similar we shall show this only for $f$.
Since $(Z_1, d_{Z_1})$ is a path-metric space it is enough to show that if $z, w\in Z_1$
and $d_{Z_1}(z,w)\leq 1$, then $d_{Z_2}(f(z),f(w))$ is bounded by a constant independent of $z, w$. However, 
$$
d_Z(f(z),f(w))\leq d_Z(z,f(z))+d_Z(w,f(w))+d_Z(z,w)\leq 1+2R.$$
Hence $d_{Z_2}(f(z), f(w))\leq \eta(2R+1)$, since $Z_2$ is $\eta$-properly
embedded in $Z$. Now the claim follows from Lemma \ref{lem: qi from lipschitz}. \qed


\begin{lemma}\label{lem:quasigeodesic-paths} 
Given $D\geq 0$, $k\geq 1$ and $\eta:\RR_{\geq 0}\map \RR_{\geq 0}$, there is $K=K_{\ref{lem:quasigeodesic-paths}}(D,k,\eta)$  
with the following property:

Suppose $X$ is a metric space and $x,y\in X$ are arbitrary points. 
Suppose also that  $c$ is  a $1$-Lipschitz path  in a metric space $X$ joining points $x, y$ is $\eta$-properly embedded in $X$ 
and there is a continuous arc-length parametrized $k$-quasigeodesic $c_1$ joining $x, y$ in $X$ such that $\Hd(c,c_1)\leq D$.  
Then $c$ is a $K$-quasigeodesic. 
\end{lemma}

\proof Let $c: [0, l]\map X$, $c_1: [0, l_1]\map X$ be the given paths. Since $c$ is $1$-Lipschitz, for all $s,s'\in [0, l]$ 
we have $d_X(c(s), c(s'))\leq |s-s'|$. On the other hand, there are points $t,t'\in [0, l_1]$ such that $d_X(c_1(t), c(s))\leq  D$,
$d_X(c_1(t'), c(s'))\leq D$. Let 
$$
t_1=t\leq t_2\leq\cdots\leq t_n\leq t_{n+1}=t'$$
be points of $[0, l_1]$ such that 
$t_{i+1}-t_i=1$, $1\leq i\leq n-1$, and $t_{n+1}-t_n\leq 1$.
Then there are points 
$$
s_1=s, s_2,\cdots, s_n, s_{n+1}=s'$$
in $[0,l]$ such that $d_X(c_1(t_i), c(s_i))\leq D$. 
It follows that 
$$
d_X(c(s_i), c(s_{i+1}))\leq 2D+ d_X(c_1(t_i), c_1(t_{i+1}))\leq 2D+2k, \quad 1\leq i\leq n.
$$ 
Hence, $|s_i-s_{i+1}|\leq \eta(2D+2k)$ and, therefore, 
$$
|s-s'|\leq n\eta(2D+2k)\leq \eta(2D+2k)+\eta(2D+2k)|t-t'|.$$ 

However, $$|t-t'|\leq k^2+kd_X(c_1(t),c_1(t'))\leq k^2+2Dk+kd_X(c(s),c(s')).$$
It follows that
$$
-\frac{1+k^2+2Dk}{k}+\frac{1}{\max\{1, k\eta(2D+2k)\}}|s-t|\leq d_X(c(s),c(t))\leq |s-t|.
$$ 
Thus, we may take $K=\max\{1,\eta(2D+1), \frac{1+k^2+2Dk}{k}\}$. \qed

\medskip 

The next lemma follows immediately from the definition of a uniformly proper embedding: 

\begin{lemma}\label{lem:nested neighborhoods}
Suppose that $(X,d_X)$ is a path-metric space, $Y\subset X$ is rectifiably connected and 
$\eta$-properly embedded in $X$, i.e. the inclusion map $(Y,d_Y)\to (X,d_X)$ is $\eta$-proper. 
Then for every subset $Z\subset Y$ and $R\ge 0$ we have
$$
N_R^Y(Z)\subset N^X_R(Z)\subset N^Y_{\eta(R)}(Z). 
$$
Here the $R$-neighborhood $N^Y_R$ in $Y$ is taken with respect to the induced path-metric $d_Y$.  
\end{lemma}

In the following three lemmas, $I=[a, b]$, $I'=[a', b']$ denote nondegenerate intervals in $\RR$ (equipped with the standard metric).

\begin{lemma}
[Lipschitz approximation] \label{lem:approx}
Let $f: I\to I'$ be a coarse $(L,\epsilon)$-Lipschitz map. Then $f$ is within distance $\le 2(L+\epsilon)$ from a piecewise-linear 
$2(L+\epsilon)$-Lipschitz map $g: I\to I'$. 
\end{lemma}
\proof 1. First, assume that $b-a\ge 1$. We then subdivide the interval $I$ into subintervals 
$$
[a_0, a_1]= [a, a_1], [a_1, a_2], \ldots, [a_n, b]= [a_n, a_{n+1}]$$ 
of length at least $1/2$ and at most $1$. We replace the restriction of $f$ to each subinterval 
$I_i=[a_i, a_{i+1}]$ with a linear function $g_i$ such that $g_i(a_i)=f(a_i)$, $g_i(a_{i+1})=f(a_{i+1})$. Since 
\begin{equation}\label{eq:interp} 
|f(a_i) -f(a_{i+1})|=|g_i(a_i) - g_i(a_{i+1})| \le L |a_i - a_{i+1}| +\epsilon \le L+\epsilon, 
\end{equation}
it is easy to see that 
$$
d(f|_{I_i}, g_i)\le 2(L+\epsilon). 
$$
Combining the linear functions $g_i$ we obtain a piecewise-linear function $g: I\to I'$ such that 
$d(f, g)\le 2(L+\epsilon)$.  Since $a_{i+1}-a_i\ge 1/2$, the inequality \eqref{eq:interp} implies that the slope of each 
$g_i$ is at most $2(L+\epsilon)$. Hence, $g$ is   $2(L+\epsilon)$-Lipschitz. 

2. Suppose that $b-a<1$. Then we let $g$ be the constant function, equal $f(b)$. Since $|f(s)- f(t)|\le L+\epsilon$ for all 
$s, t\in I$, $d(f, g)\le L+\epsilon$. \qed

\begin{lemma}
[Coarse monotonicity of quasiisometries] \label{lem:c-mono} 
Set $D:= L(5\epsilon+4L)$. Suppose that $f: I\to I'$ is an $(L,\epsilon)$-qi embedding. 
Then $f$ is coarsely monotonic in the sense that if $r< s< t$ are in $I$ and $\min(s-r , t-s) > D$ 
then $f(s)$ is between $f(r)$ and $f(t)$. 
\end{lemma}
\proof Let $g$ be the Lipschitz approximation of $f$ as in Lemma \ref{lem:approx}. Since $f$ is an 
$(L,\epsilon)$-qi embedding and $d(f,g)\le 2(L+\epsilon)$, we conclude that $g$ satisfies
$$
L^{-1}|t- t'| - (5\epsilon+4L) \le |g(t)- g(t')|
$$
for all $t, t'\in I$. Suppose that, say, 
$$
f(s)> \max(f(r), f(t)). 
$$ 
Then
$$
\min( f(s)- f(r), f(s)- f(t)) \ge L^{-1}D -\epsilon.
$$
Once again since $d(f,g)\leq 2(L+\epsilon)$ we see that 
$$
\min( g(s)- g(r), g(s)- g(t)) \ge D'= L^{-1}D - (5\epsilon +4L). 
$$
For concreteness, assume that $g(r)\le g(t)$. By the Intermediate Value Theorem applied to the function $g$ (restricted to the interval 
$[r,s]$) there exists $t'\in [r,s]$ such that $g(t')=g(t)$. We have $t-t'\ge t-s>D$. Thus,
$$
L^{-1} D - (5\epsilon+4L)  < L^{-1}|t- t'| - (5\epsilon+4L) \le 0 \Rightarrow D< L(5\epsilon+4L) 
$$
which is a contradiction. If $f(s)< \min(f(r), f(t))$ then one can apply the same proof to the function $t\mapsto -f(t)$
to arrive at a contradiction. \qed

\begin{lemma}
[Approximating quasiisometries by homeomorphisms] \label{lem:approximation} 
Let $I=[a,b]$, $I'=[a', b']$ be nondegenerate intervals in $\RR$ (equipped with the standard metric). 
Suppose that $f: I\to I'$ is a $k$-quasiisometry sending $a$ to 
$a'$ and $b$ to $b'$. Then there exists  a (piecewise-linear) homeomorphism $\tilde{f}: I\to I'$ within 
distance $D_{\ref{lem:approximation}}=D_{\ref{lem:approximation}}(k)$ from $f$ which is
also a ${k}_{\ref{lem:approximation}}(k)$-quasiisometry.
\end{lemma}
\proof The proof is similar to that of Lemma \ref{lem:approx}. Set $L=\epsilon=k$. 

1. Suppose first that $b-a\ge 2D$, where $D= L(5\epsilon+4L)= 9k^2$ (as in Lemma \ref{lem:c-mono}). 
We subdivide the interval $I=[a, b]$ into subintervals 
$$
[a_0,a_1]=[a, a_1], [a_1,a_2], \ldots, [a_n, a_{n+1}]=[a_n, b], 
$$  
each of length greater than $D$ and at most $2D$. According to Lemma \ref{lem:c-mono}, $f$ restricted to the subset 
$$
J=\{a_0, a_1,..., a_n, a_{n+1}\}
$$
is strictly monotonic. We then let $\tilde{f}=g: I\to I'$ be the piecewise-linear function equal to $f$ on $J$ and linear on the complementary intervals. In particular,
$$
f(a)=a'=g(a), \quad f(b)=b'=g(b). 
$$
In view of monotonicity of $f|_J$, the function $g$ is also (strictly) monotonic, hence, a homeomorphism.  
For $s, t\in [a_i, a_{i+1}]$ we have 
$$
|f(s) - f(t)|\le 2D L + \epsilon= 18k^3 + k. 
$$
Therefore, for all $t\in [a_i, a_{i+1}]$ we have 
\begin{align*}
|f(t)- g(t)|\le |f(t)-f(a_i)|+|f(a_i)-g(t)|\leq\\ |f(t)-f(a_i)|+|f(a_i)-f(a_{i+1})| 
\leq 2(18k^3+k),
\end{align*}
since $f(a_i)=g(a_i), f(a_{i+1})=g(a_{i+1})$ and $g$ is monotonic. 

2. Suppose that $b-a\le 2D$. We then let $g=\tilde{f}$ be the linear function equal to $f$ on $\{a, b\}$; $g(a)=a'< g(b)=b'$. 
As in Case 1, $d(f, g)\le 36k^3+2k$. 
Hence we can choose $D_{\ref{lem:approximation}}=36k^3+2k$. Since $f$ is a $k$-quasiisometry
it follows that $\tilde{f}$ is a ${k}_{\ref{lem:approximation}}(k)=(k+D_{\ref{lem:approximation}})$-quasiisometry.
\qed

\begin{defn}\label{def:cob} \index{Lipschitz cobounded subsets}
Two subsets $Y, Z$ of a metric space $X$ are said to be $C$-{\em Lipschitz cobounded} if there exist 
$(L,\epsilon)$-coarse Lipschitz retractions $X\to Z, X\to Y$ whose restrictions 
$$
r_{Y,Z}: Y\to Z, \quad r_{Z,Y}: Z\to Y,
$$
satisfy: 

1. $r_{Y,Z}(Y)$ and $r_{Z,Y}(Z)$  have diameters $D_Y, D_Z$. 

2. $\max(L, \epsilon, D_Y, D_Z)\le C$. 
\end{defn}

In Section \ref{sec:Cobounded pairs of subsets} we will relate this definition to the more standard notion of cobounded subsets in a 
Gromov-hyperbolic space.

\begin{lemma}\label{lem:cobounded1}
If $Y, Z$ are $C$-Lipschitz cobounded, then for every $R$ there exists $D=D_{\ref{lem:cobounded1}}(R,C)$ such that if 
$$
a_i\in Y, b_i\in Z, i=1, 2
$$
are points satisfying $d(a_i, b_i)\le R$, $i=1,2$, then $d(a_1,a_2)\le D, d(b_1, b_2)\le D$. 
\end{lemma}
\proof One can take $D=2C(R +1) +C$. \qed

\section{Coproducts, cones and cylinders} \label{sec:cylinders}

In this section we discuss several purely topological notions used elsewhere in the book. Let $\{Z_\alpha: \alpha\in A\}$ be an indexed collection of topological spaces. Then the {\em coproduct} topology on the 
disjoint union 
$$
Z:= \coprod_{\alpha\in A} Z_\al
$$
is the finest topology on the disjoint union such that all the natural inclusion maps  $Z_\al\to Z$ are continuous. 
In particular, each $Z_\al\subset Z$ is a clopen subset homeomorphic to $Z_\al$.

\medskip 
In the following two constructions, the unit intervals $[0,1]$ are sometimes replaced by the half-intervals $[0,1/2]$ with $1/2$ playing the role of $1$. 

Let $X, Y$ be topological spaces, $f: X\to Y$ a continuous map. Then the {\em mapping cylinder} $Cyl(f: X\to Y)$, also denoted 
$X\cup_f Y$ and, sometimes, $Y\cup_f X$, 
is the quotient of $X\times [0,1] \sqcup Y$ (with the coproduct topology) by the equivalence relation 
$$
(x,1)\sim f(x), x\in X. 
$$

Lastly, the {\em cone} $C(a,X)$ over a topological space $X$ is the quotient of the product $X\times [0,1]$ by 
the subspace $X\times \{1\}$. 
The point $a$, the projection of $X\times \{1\}$ in $C(a,X)$, is called the {\em apex} of the cone. 
The projections of the intervals $\{x\}\times [0, 1]$ under the quotient map $q: X\times [0,1]\to C(a,X)$ will be called the 
{\em radial line segments} in $C(a,X)$. We will identify $X$ with the image of $X\times \{0\}$ in $C(a,X)$.

In the next section, we will metrize the cones $C(a,X)$.

\section{Cones over metric spaces}\label{sec:metric cones}  

Let $(X,d)$ be a path-metric space. We equip  $X\times [0, 1/2]$ with the product metric.  
We metrize the corresponding cone $C(a,X)$ as follows. Projecting rectifiable paths 
$$c: I=[0,1]\to X\times [0, 1/2]$$ we obtain a family of {\em admissible} paths in $C(a,X)$. For a rectifiable 
path $c$, we define the length  $L(q(c))$ to be equal to the length of $c$ minus the total length of the intersection $c(I)\cap X\times \{1/2\}$. In other words, the length of the path $q\circ c$ 
is the total length of $c|_J$, where $J= I \setminus (q\circ c)^{-1}(\{p\})$. 
For instance, the radial line segment connecting $x$ to $a$ has length $1/2$. 
Using this notion of length we path-metrize the cone $C(a,X)$. We leave it to the reader to check that this 
metric, denoted $\hat{d}$, metrizes the topology of $C(a,X)$ and that for each $x\in X$, 
$$
\hat{d}(x,a)=1/2. 
$$

\section{Approximation of metric spaces by metric graphs} \label{sec:approximation} 

In this section we discuss a generalization of path-metric spaces, called {\em quasi-path metric spaces}. 

\begin{defn}
A {\em finite $r$-path} in a metric space $(X ,d)$ is a map 
$$
c: [m, n]\subset \Z \map X$$
such that $d(c(i), c(i+1))\le r$ for all $m\leq i\leq n-1$. 
The {\em length} of such $c$ is defined as 
$$
\length(c)=\sum_{i=m}^{n-1} d(c(i), c({i+1})).$$
The finite path $c$ is said to {\em connect} the point $x=c(m)$ to the point $y=c(n)$. 

A metric space $(X,d)$ is called an {\em $r$-quasi-path metric space} for a constant $r>0$ if
for every pair of points $x, y\in X$ there exists a finite $r$-path $c$ connecting $x$ to $y$ such that $\length(c)\le d(x,y)+ r$. 
\end{defn}

For instance, every path-metric space is an $r$-quasi-path metric space for every $r>0$.

\begin{lemma}\label{quasi path sp qi to path sp}
Any $r$-quasi-path metric space is $(1,3r)$-quasiisometric to a path-metric space.
\end{lemma}
\proof Suppose $X$ is an $r$-quasi-path metric space. We construct a metric graph $Z$ with the vertex set
$V(Z)=X$ such that $x, y\in X$ are connected by an edge $e$ iff $x\neq y$ and $d_X(x,y)\leq r$, 
where the edge is assigned the length $\ell(e)=d_X(x,y)$.

Consider the inclusion map $\iota: X\map Z$. Suppose
$x,y\in X$ are arbitrary points. Then there is an $r$-path $x=x_0, x_1, \cdots, x_n=y$ in $X$ joining
$x$ to $y$ in $X$. By the definition of the graph $Z$, $x_i$'s  also form a sequence of vertices connected by edges in $Z$.
Hence, 
$$
d_Z(x,y)\leq \sum_{i}d_Z(x_i,x_{i+1})\leq \sum_{i}d_X(x_i,x_{i+1}) \leq d_X(x,y)+r.
$$ 
Thus, $\iota$ is $(1,r)$-coarse Lipschitz. Let $\rho: Z\map X$ be the following map. The restriction of $\rho$
on $V(Z)$ is simply the identity map and interiors of edges are mapped to one of the vertices.
Let $\alpha: I \map Z$ be any piecewise-linear path (see \cite{bridson-haefliger}, Chapter I.1, Section 1.9). Then
clearly, $\length(\alpha)$ and the length of the $r$-path $\rho\circ \alpha$ differ by at most $2r$.
Hence, $\rho$ is $(1,2r)$-coarsely Lipschitz. Moreover, it is clear that $d(\id_X, \rho\circ \iota)\leq r$
and $d(\id_Z, \iota\circ \rho)\leq r$. Hence, by Lemma \ref{lem: qi from lipschitz}, the maps $\iota, \rho$
are both $(1, 3r)$-quasiisometries. \qed

\begin{defn}[Rips graph] 
Let $(Y,d)$ be a metric space. For $R\ge 0$ the {\em $R$-Rips graph} of $(Y,d)$ is the graph $Z_R$ \index{Rips graph}
with the vertex-set $Y$ and edges $[y_1,y_2]$ for all pairs of distinct points $y_1, y_2\in Y$ 
such that $d(y_1,y_2)\le R$.  We will equip  $Z$ with its graph-metric (each edge has unit length). 
\end{defn}

Note that for a general metric space $Y$, the graph $Z_R$ is disconnected and  
the distance between points in different connected components is infinite. 
However, if $Y$ is a path-metric space, then each graph $Z_R$ is connected. 

\begin{defn}\index{coarsely connected metric space}
A metric space $(Y,d)$ is said to be {\em coarsely connected}\index{coarsely connected space} 
 if there exists $R<\infty$ such that the corresponding Rips graph $Z_R$ is connected. 
\end{defn}

The following fundamental result of geometric group theory is usually stated for {\em proper} metric spaces $Y$ and properly discontinuous actions, but, when the notion of metrically proper actions is used, properness of $Y$ is not needed 
and the proof is  the same as in the proper case, cf. \cite{Drutu-Kapovich}:

\begin{lemma}[Milnor--Schwarz Lemma] 
 \index{Milnor--Schwarz Lemma}\label{lem:Milnor--Schwarz Lemma}
Suppose that $(Y, d)$ is a (nonempty) metric space, $G$ is a discrete group and $G\acts Y$ is 
a geometric\footnote{i.e. isometric, metrically proper and cobounded} action. Then:

1. If $Y$ is coarsely connected, then $G$ is finitely generated. 

2. If $Y$ is a quasi-path metric space, then  for one (equivalently, every) $y\in Y$ the orbit map $G\to Gy\subset Y$ is a quasiisometry. 
\end{lemma}

\begin{lemma}\label{path sp qi to graph}
For a path-metric space $X$, let $Z=Z_1$ be the $1$-Rips graph of $X$. Then the inclusion map
$\iota: X\map Z$ is a $(1,1)$-quasiisometry with a $(1,3)$-qi inverse $\rho:Z\map X$. 
\end{lemma}
\proof The proof is very similar to that of the previous lemma: We let 
the map $\rho: Z\map X$ be identity on the $V(Z)=X$, etc.
Given $x, y\in X$, we join them in $X$ by an arc-length parametrized path 
$\gamma: [0,l]\map X$ such that $l\leq d_X(x,y)+\epsilon$, 
where $\epsilon>0$ is chosen in such a way that $d_X(x,y)<m+1$ where $m$
is the nonnegative integer determined by $m\leq d(x,y)<m+1$. 
Since $\length(\gamma)< m+1$,  it follows that
$d_Z(x,y)\leq m+1\leq d_X(x,y)+1$. Suppose $x,y\in X$ such that $d_Z(x,y)= n$. Let $x=x_0,x_1,\cdots, x_n=y$ the consecutive vertices on a geodesic in $Z$ joining $x,y$. Then we know that $d_X(x_i, x_{i+1})\leq 1$. 
Thus, 
$$
d_X(x,y)\leq \sum^n_{i=1}d_X(x_{i-1},x_i)\leq n
$$
and we get $d_X(x,y)\leq d_Y(f(x),f(y))\leq d_X(x,y)+1$. This proves the first statement of the lemma. 

Finally, $\iota(X)$ is a $1$-net in $Z$, and, hence, $\iota$ is coarsely $1$-surjective. 
The remaining parts of the proof follow 
from simple calculations and we leave details to the reader. \qed



\begin{cor}\label{cor: graph approx}
Suppose $X$ is an $r$-quasi-path metric space. Then there is a (connected) graph $X'$ equipped with the graph-metric 
 and a $(1, 3r+1)$-quasiisometry $\iota: X\map X'$ with a 
$(1,3r+3)$-qi inverse $\rho: X'\map X$ such that  $\rho(X')=X$. In particular, 
each quasi-path metric space is $(1, \eps)$-quasiisometric to a geodesic metric space. 
\end{cor}
\proof This is a straightforward consequence of Lemma \ref{quasi path sp qi to path sp}
and Lemma \ref{path sp qi to graph}.
\qed


\begin{lemma}
[Local-to-global principle for coarse Lipschitz maps from quasi\-path metric spaces]
\label{lem: proving coarse Lipschitz}
Suppose that $(X,d_X)$ is an $r$-quasi-path metric space, $(Y, d_Y)$ is any metric space
 and $f: X\to Y$ is a map such that for all $x_1, x_2\in X$, $d_X(x_1,x_2)\le r$ implies that 
$d_Y(f(x_1), f(x_2))\le r'$ for some $r'$ independent of $x_1, x_2$. 
Then $f$ is $(\frac{2r'}{r}, 3r')$-coarse Lipschitz. 
\end{lemma}
\proof Take $x, y\in X$. Suppose $n\in \mathbb N$ is the smallest integer such that
there is a finite sequence $x_0=x, x_1,\cdots, x_n=y$  in $X$
with $d_X(x_i, x_{i+1})\le r$ for all $i=0,...,n-1$ and $\sum_{i=0}^{n-1} d_X(x_i, x_{i+1})\le d(x,y)+ r$.
Then $d_X(x_i, x_{i+2})>r$ for $0\leq i\leq n-2$.
It follows that $r(n-1)/2<d_X(x,y)+r$. Hence, $n<3+\frac{2}{r} d_X(x,y)$.
On the other hand 
$$
d_Y(f(x), f(y))\leq \sum_{i=0}^{n-1} d_Y(f(x_i), f(x_{i+1}))\leq nr'< 3r'+ \frac{2r'}{r}d_X(x,y). \qed
$$


\medskip

Let $Y$ be a $D$-net in metric space $(X,d)$ with the inclusion map $\iota: Y\to  X$. 
Given $R>0$ we let $Z=Z_R$ be the full subgraph of the $R$-Rips graph 
of $(X,d_X)$ with the vertex set $Y$; we equip $Z$ with its graph-metric. We extend the map $\iota$ to the rest of $Z$, 
(the extension is still denoted by $\iota$), by taking an arbitrary orientation on $Z$ and 
sending all points of any open directed edge $[v,w]\setminus \{w\}= [v,w)$ in $Z$ to the point 
$v\in Y\subset X$. The next lemma (which is a form of the Milnor--Schwarz Lemma for metric spaces, 
cf. Theorem 8.52 in \cite{Drutu-Kapovich}) shows that, under some 
conditions, the graph $Z$ is connected and the map $\iota$ is a uniform quasiisometry $Z\to X$. This result generalizes 
Lemma \ref{path sp qi to graph}.

\begin{lemma}\label{lem:graph-approximation} 
Suppose $(X,d_X)$ is an $r$-quasi-path metric space and $Y\subset X$ is a $D$-net in $X$. 
If $R\ge r+2D$, then $Z$, as defined above, is a connected graph and the  map $\iota: (Z, d_Z)\to  X$ 
is a $(K_{\ref{lem:graph-approximation}}(r,R), \epsilon_{\ref{lem:graph-approximation}}(r,R))$-quasiisometry. 

Moreover, there is a $(1+r)$-coarse inverse $\rho:X\map Z$ to $\iota$ such that $\rho(x)=x$ for all $x\in Y$ which is
also a $(K_{\ref{lem:graph-approximation}}(r,R), \epsilon_{\ref{lem:graph-approximation}}(r,R))$-quasiisometry. 
\end{lemma}
\proof 1. Consider vertices $y, y'\in Y$. Since $(X,d)$ is a $r$-quasi-path metric space, there exists an $r$-path 
$$
y=x_0, x_1,..., x_n x_{n+1}=y'
$$
in $X$ from $y$ to $y'$. Since  $Y\subset X$ is a  $D$-net in $X$, there exist points $y_1,...,y_n\in Y$ satisfying 
$$
d_X(x_i, y_i)\le D, i=1,...,n. 
$$
By the triangle inequality, $d_X(y_i, y_{i+1})\le r+2D\le R$ for $i=0,...,n$ where $y_{n+1}=y'$, which implies that the vertices 
$x_0=y, y_1,...,y_n, y'$ are on an edge path in $Z$. This proves that $Z$ is connected.

Suppose $z\neq z'\in Z$ are any two  points. Let $\iota(z)=y, \iota(z')=y'$ and
let $z, y_0,\cdots, y_n, z'$ be a geodesic in $Z$ joining $z, z'$ where $y_i\in Y$.
We note that $d_Z(y,y_0)\leq 1, d_Z(y_n,y')\leq 1$ and $|d_Z(z, z') - n|\leq 2$. Now, 
\begin{align*}
d_X(\iota(z), \iota(z'))\leq d_X(y,y_0)+d_X(y_n, y')+\sum_{i=0}^n d_X(y_i, y_{i+1})\leq\\ 
2R+nR\leq 2R+R(d_Z(z, z')+2)=4R+Rd_Z(z, z').
\end{align*}
Hence the map $\iota:Z\map X$ is $(R,4R)$-coarse Lipschitz.

Next, we define a coarse inverse map $\rho$ to $\iota$ by defining it to be the identity map on $Y$ and 
sending  $x\in X\setminus Y$ to a point $\rho(x)=y\in Y$ such that $d_X(x, y)\le D$. 
Then $d_X(x, x')\le r$ implies that $d_X(\rho(x), \rho(x'))\le r+2D\le R$, i.e. $d_Z(\rho(x), \rho(x'))\le 1$. 
Since $(X,d)$ is a $r$-quasi-path metric space, the map $\rho$ is $(2/r, 3)$-coarse Lipschitz by Lemma \ref{lem: proving coarse Lipschitz}. 

By the construction, 
$$
d(\rho\circ \iota(t), t)\le 1, \quad d(\iota \circ \rho(x), x)\le r. 
$$
Hence, by Lemma \ref{lem: qi from lipschitz} we can take 
$$
K_{\ref{lem:graph-approximation}}(r,R)=\max\{R, 2/r\}, \quad  
\epsilon_{\ref{lem:graph-approximation}}(r,R)=\max\{ 4R,3\}+2\max\{1, r\}. \qedhere 
$$

Similarly to the approximation of coarse Lipschitz maps of line segments by piecewise-linear maps, one has a
uniform approximation of coarse Lipschitz maps of metric (simplicial) graphs by simplicial maps, i.e. 
maps sending every edge linearly to an edge or a vertex. 

\begin{lemma}\label{lem:simplicial-approximation}
Fix numbers $R, L, \epsilon$. 
Let $X, Y$ be connected metric graphs with the same edge-length $R$. Then, after subdividing edges
of $X$ in at most $n(R,L,\epsilon)$ equal length subsegments,  every $(L,\epsilon)$-coarse Lipschitz
map $f: X\to Y$ is within distance 
$D_{\ref{lem:simplicial-approximation}}(L,\epsilon,R)$ from a simplicial map $f_*: X\to Y$.  
\end{lemma}
\proof The proof is similar to that of Lemma \ref{lem:approx} and we omit it. \qed 

\medskip 
Since $Z$ constructed in this lemma is a complete geodesic metric space and since we are interested in
coarse geometric properties of metric spaces, we can always replace quasi-path metric spaces with
appropriate metric graphs. We will be using this in Section \ref{sec:trees-of-spaces} when constructing
total spaces of trees of spaces.


\section{Hyperbolic metric spaces}

Hyperbolic metric spaces are coarsifications of the classical hyperbolic $n$-space $\H^n$ and are characterized by a form of {\em thin triangle condition}. 
The most common notions of hyperbolicity for metric spaces are the one due to Rips (for geodesic metric spaces) and one due to Gromov (for general metric spaces). One drawback of Gromov's definition is that his notion of hyperbolicity is not qi invariant, although it is invariant under $(1,\eps)$-quasiisometries. One of the features (or bugs, depending on the perspective) of metric hyperbolicity is that it is stable under changes in metric below certain scale $\delta$ and that, accordingly, nothing can be said about general hyperbolic spaces below that scale. This also points to a limitation of Rips' notion of hyperbolicity since it applies only to geodesic metric spaces. This becomes somewhat important in the context of this book since metric spaces that we consider are frequently only {\em path-metric spaces}. One source where hyperbolicity along the lines of Rips definition is developed   for path-metric spaces $X$ is   V\"{a}is\"{a}l\"{a}'s paper \cite{MR2164775}: Instead of geodesics he considers {\em $h$-short paths}, which are rectifiable paths between points $x, y\in X$ whose length is at most $d(x,y)+h$. A drawback is that one is forced to carry an extra constant. Another possible approach is 
to extend Rips' definition to the class of {\em quasi-path metric spaces}. We will give basic definitions in Section \ref{sec:quasi-hyp} but will not pursue this direction much further beyond proving that such  metric spaces are $(1,\eps)$-quasiisometric to geodesic metric spaces and, hence, Gromov's notion of hyperbolicity in the context of quasi-path metric spaces is preserved by general quasiisometries, see Section \ref{sec:stability}. Yet, another possible approach is to work with path-metric spaces but instead of geodesics, work with sequences of paths whose lengths approximate distances between points. All the arguments appearing in the book will go through with constants unchanged comparing to the ones for geodesic metric spaces. A drawback is that this approach lengthens  the proofs (which are already long and technical  in chapters 
\ref{ch:4 classes}, \ref{sect:3}, \ref{ch:flows} and  \ref{sec:everything together}). Thus, for most of the book, we work with geodesic metric spaces. 
In this section we present various notions of hyperbolicity starting with the most familiar ones.

We assume that the reader is familiar with the basic definitions and facts about  hyperbolic metric spaces that can be found 
for instance in \cite{bridson-haefliger}, \cite{CDP}, \cite{Drutu-Kapovich}, \cite{gromov-hypgps}, \cite{GhH}, \cite{Shortetal}, \cite{MR2164775}. 
In this section we collect some of these to fix the notions and for later use.

\subsection{Hyperbolicity in the sense of Gromov} 

\begin{defn}[Gromov  product]\index{Gromov product} 
 Let $X$ be a metric space. Given points  $x, y, z\in X$, the 
 {\em Gromov  product} of $y, z$ with respect to $x$, denoted $(y.z)_x$, is defined as
$$
\frac{1}{2}\left( d(x,y)+d(x,z)-d(y,z) \right).$$
\end{defn}

\begin{defn}\index{hyperbolic space in the sense of Gromov} 
A metric space $X$ is said to be  {\em $\delta$-hyperbolic in the sense of Gromov} or simply $\delta$-{\em Gromov-hyperbolic} if 
 all $x, y, z, w\in X$ satisfy the inequality 
$$
(x.y)_w \geq \min\{ (x.z)_w, (y.z)_w\}- \delta.$$
A metric space $X$ is said to be {\em hyperbolic in the sense of Gromov} if it is $\delta$-{\em Gromov-hyperbolic} for some $\delta\in [0,\infty)$. 

A finitely-generated group $G$ is called {\em hyperbolic} if the metric space $(G, d_G)$ is hyperbolic in the sense of Gromov, where $d_G$ is the word metric on $G$ for some finite generating set. \index{hyperbolic group}  
\end{defn}

\begin{example}
Consider the graph $G\subset \RR^2$ of the function $y=|x|$. We equip $G$ with the restriction of the standard Euclidean metric on $\RR^2$. 
For $n\in \N$ consider points 
$$
o=(0,0), p= (-n,n), q=(n,n), z=(2n,2n). 
$$ 
Then $(p, q)_o= (\sqrt{2}-1)n$, $(p,z)_o= 2\sqrt{2}n$, $(q,z)_o=\sqrt{2}n$. The difference 
$$
(p.q)_o -  \min\{ (p.z)_o, (q.z)_o\}= (\sqrt{2}-1)n - \sqrt{2}n = -n 
$$
diverges to $-\infty$ as $n\to\infty$. Thus, $G$ is not  hyperbolic in the sense of Gromov.  On the other hand, $G$ is qi to the real line via the map 
$x\mapsto (x, |x|)$ and the real line is $0$-hyperbolic. 
\end{example}

Thus, Gromov-hyperbolicity is not preserved by quasiisometries even for {\em quasigeo\-desic metric spaces}, i.e. metric spaces where all points are connected by uniform quasigeodesics. On the other hand, the following lemma is straightforward from the definition of Gromov-hyperbolicity: 

\begin{lemma}\label{lem: gromov hyp qi inv}
Suppose $X, Y$ are metric spaces and $f: X\map Y$ is a $(1,\epsilon)$-quasi\-iso\-met\-ry for some
$\epsilon\geq 0$. Then $X$ is Gromov-hyperbolic iff $Y$ is. More precisely, if $X$ $\delta$-hyperbolic in the sense of Gromov, 
then $Y$ is $\delta+3\eps$-hyperbolic in Gromov's sense. 
\end{lemma}

\subsection{Hyperbolicity in the sense of Rips} 

Suppose now that $X$ is a geodesic metric space. 

\begin{defn}
Consider a geodesic triangle $\Delta= \Delta x_1x_2x_3\subset X$  with the vertices $x_1$, $x_2$, $x_3$,
and let  $\delta\geq 0$. 

\begin{enumerate}
\item The triangle $\Delta$ is said to be  $\delta$-{\em slim} \index{slim triangle}
if each side of  $\Delta$ is contained in the
$\delta$-neighborhood of the union of the other two sides.  

\item For all $i\neq j\neq k\neq i$, let $c_k \in x_ix_j$ be such that
$d(x_i,c_j)=d(x_i,c_k)$. The points $c_i$ are called the {\em internal points} of $\Delta$. \index{internal points} 
Note that, for all $i\neq j\neq k\neq i$, 
$$
d(x_i,c_j)=\frac{1}{2}\left(d(x_i,x_j)+d(x_i,x_k)-d(x_j,x_k)\right)= (x_j . x_k)_{x_i}. 
$$

\item If $X$ is a tree, then $p= c_1=c_2=c_3$ and in this case we shall refer to the point $p$ as the {\em center} of the 
$\Delta$.

\item The diameter of the set $\{c_1, c_2, c_3\}$ will be referred to as the {\em insize} of the triangle $\Delta$.\index{insize} 

\item The triangle $\Delta$ is said to be $\delta$-{\em thin} if for all $i\neq j\neq k\neq i$ and \index{thin triangle} 
$p\in x_ic_j\subset x_ix_k$, $q\in x_ic_k\subset x_ix_j$ with $d(p,x_i)=d(q,x_i)$, one has 
$$
d(p,q)\leq \delta.$$
\end{enumerate}
\end{defn}

The next lemma is clear from the definitions:

\begin{lemma}
If $\Delta$ is $\delta$-thin, then it is also $\delta$-slim and its insize is $\le \delta$. 
\end{lemma}

\begin{defn}[Rips hyperbolicity] \index{hyperbolic space in the sense of Rips} 
A geodesic metric space $X$ is said to be $\delta$-hyperbolic in the sense of Rips if each geodesic triangle in $X$ is $\delta$-slim. A geodesic metric space is said to be {\em Rips-hyperbolic} if it is $\delta$-hyperbolic in the sense of Rips for some $\delta<\infty$.  
\end{defn}

\begin{lemma}[Proposition 2.1 in  \cite{Shortetal}]\label{hyp-defn} 
Suppose $X$ is a $\delta$-hyperbolic metric space in the sense of Rips. Then the following hold:
\begin{enumerate}
\item All the geodesic triangles in $X$ have insize at most $4\delta$.
\item All the geodesic triangles in $X$ are $6\delta$-thin.
\end{enumerate}
\end{lemma}

It follows that a geodesic metric space is hyperbolic in the sense of Rips if and only if all geodesic triangles in $X$ are uniformly thin.

The following  lemmata  are also very standard and follow easily from definitions, see for instance  \cite{Drutu-Kapovich}, \cite{Shortetal}, or \cite{MR2164775}: 

\begin{lemma}\label{lem:gromov iff rips}
A geodesic metric space is Rips-hyperbolic if and only it is Gromov-hyperbolic.  More precisely:
\begin{enumerate}
\item If a metric space $X$ is $\delta$-hyperbolic in the sense of Rips, then $X$ is $3\delta$-hyperbolic in Gromov's sense. 

\item If $X$ is geodesic and  $\delta$-hyperbolic in Gromov's sense, then $X$ is $2\delta$-hyperbolic in the sense of Rips. 
\end{enumerate}
\end{lemma}

In view of this lemma, when talking about hyperbolicity for geodesic metric spaces, we will always mean hyperbolicity in the sense of Rips. 

\medskip 
For geodesic hyperbolic spaces, the Gromov-product $(y.z)_x$ ``almost equals'' the distance from $x$ to the geodesic $yz$: 

\begin{lemma}\label{lem:gromov product}
Suppose that $X$ is a $\delta$-hyperbolic space in the sense of Rips. Then for each triple $x, y, z\in X$,
$$
 (y.z)_x - 2\delta \le  (y.z)_x \le d(x, yz).  
$$ 
\end{lemma}

\begin{lemma}\label{lem:slim-quad}
Every geodesic quadrilateral $\square=xyzw$ in a  $\delta$-hyperbolic metric space $X$ is $2\delta$-slim, i.e. the side $xy$ of $\square$ is contained in the $2\delta$-neighborhood of the union of the other three sides of $\square$. Similarly, each geodesic $n$-gon in $X$ is 
$(n-2)\delta$-slim. 
\end{lemma}

\subsection{Hyperbolicity for path-metric spaces}\label{sec:quasi-hyp}

The form of Rips-hyperbolicity discussed in this section is a mild generalization of Rips-hyperbolicity for geodesic metric spaces. 
We refer the reader to \cite{MR2164775} for a more general discussion. 

\begin{defn}
A rectifiable  path  $c$ connecting $x$ to $y$ in a metric space $X$ is called {\em $\eps$-short} \index{$\eps$-short path}
if 
$$
\length(c) \leq d(x, y)+\eps.$$
\end{defn}

\begin{defn}[Triangles formed by paths]
Suppose $X$ is any metric space.
 Given any three points $x,y,z\in X$ and three (continuous or finite) paths   
$c(x,y), c(x,z), c(y,z)$ joining these points, the {\em triangle} formed by these paths is the set
$\{ c(x,y), c(x,z), c(y,z)\}$ and the members of this set will be called the {\em sides} of the triangle.  
\end{defn}

\begin{defn}[Slimness constant for a path-family] \label{defn:slimness constant} 
(1) Suppose $X$ is a metric space and $x, y, z\in X$. We shall say that a triangle 
formed by three paths $c(x,y), c(x,z), c(y,z)$ is {\em $\delta$-slim} for some $\delta\geq 0$, if each side of the triangle is contained in the union of  $\delta$-neighborhoods of the remaining two sides. 

(2) Given a family $\CC$ of paths $c(x,y)$ 
connecting points $x, y$ in a metric space $X$,  
for all $x\neq y\in X$ we define the {\em slimness constant} of $\CC$ as
$$
\delta_s(\CC):= \sup_{ c(x,y), c(y,z), c(z,x)\in \CC} \inf \{r:\mbox{ the triangle}\,\, 
\{c(x,y), c(y,z),c(z,x)\}\, \mbox{is} \,\,r\mbox{-slim}\}.
$$
\end{defn}

We are now ready to define a form of Rips-hyperbolicity for path-metric spaces:  

\begin{defn}[Rips hyperbolicity of path-metric spaces]  \label{defn:Rips-hyperbolic} 
If $(X,d)$ is a path-metric space, and $\epsilon\geq 0$, let $\FF_\eps$ be the family of all $\eps$-short paths $c(x,y)$ in $X$. 
We say that $X$ is {\em $\delta$-hyperbolic in the sense of Rips} if
$$
\limsup_{\eps\to 0+} \delta_s(\FF_\eps) \le \delta. 
$$
\end{defn}

\begin{rem}
\begin{enumerate}
\item $(X,d)$ is a geodesic metric space if and only if the set $\FF_0$ is nonempty.
Elements of $\FF_\eps$ are  $(1,\eps)$-quasigeodesics in $X$.  
\item For a  geodesic metric space, Definition \ref{defn:Rips-hyperbolic} is equivalent to 
the standard notion of $\delta$-hyperbolicity in the sense of Rips (with a slightly different hyperbolicity constant). 
There,  a space is $\delta_0$-hyperbolic (in the sense of Rips) if $\delta(\FF_0)=\delta_0$.
\item If for a path-metric space $X$ we have $\delta(\FF_\eps)< \infty$ for some $\epsilon \geq 0$, then $X$ is 
$\delta(\FF_\eps)$-hyperbolic in the sense of Rips. 
\end{enumerate}
\end{rem}

\subsection{Stability of quasigeodesics and qi invariance of hyperbolicity}\label{sec:stability} 

One of the most fundamental facts about hyperbolic spaces is that quasigeodesics are uniformly close to geodesics. This fact is also known as the (hyperbolic) Morse lemma, as it first appeared in a work of Morse on geodesics in the hyperbolic plane, \cite{Morse}. Morse did not have the notion of quasigeodesics and he was interested in how geodesics on a surface change with a change of its hyperbolic metric. Since changing a Riemannian metric on a compact manifold results in a quasiiisometric change of the metric on its universal cover, Morse's result can be interpreted as stability of quasigeodesics. Morse's proof was quite general and most modern proofs of stability of quasigeodesics follow the same line of reasoning.

The next lemma is a converse to Lemma \ref{lem:quasigeodesic-paths} in the setting of 
 hyperbolic spaces. 

\begin{lemma}[Morse Lemma or  stability of quasigeodesics] 
\label{stab-qg} \label{Morse}\index{Morse Lemma} 
There is a function \newline $D_{\ref{stab-qg}}=D_{\ref{stab-qg}}(\delta, k)$ defined for 
 $\delta\geq 0$ and $k\geq 1$,  such that
the following holds:

Suppose $X$ is a $\delta$-hyperbolic geodesic metric space. Then for every $k$-quasigeodesic $\phi: [a,b]\to X$, the Hausdorff distance between the image of $\phi$ 
and that of the geodesic $\phi^*$ connecting the end-points of $\phi$,  is $\le D_{\ref{stab-qg}}$. 

More precisely, according to  \cite{GS}, for a $(k,\epsilon)$-quasigeodesic $\phi$ in $X$,
$$
\Hd(\phi, \phi^*)\le 92k^2(\eps + 3\delta). 
$$ 
Thus, for a $k$-quasigeodesic one can take 
$$
D_{\ref{stab-qg}}(\delta, k)= 92k^2(k + 3\delta). 
$$
\end{lemma}

With minor modifications, the proofs go through for path-metric spaces, when geode\-sics are replaced with $\eta$-short paths $\phi^*_\eta$ 
(cf. \cite{MR2164775}).  One obtains an estimate $D_{\ref{stab-qg}}(\delta, k, \eta)$ and, hence,  
$$
D_{\ref{stab-qg}}(\delta, k)=\lim_{\eta\to 0+} D_{\ref{stab-qg}}(\delta, k, \eta). 
$$

As a consequence:

\begin{lemma}\label{lem:sub-close} 
There exists a function $D=D_{\ref{lem:sub-close}}(\delta,k,r)\ge r$ such that the following holds. 
If $X$ be a $\delta$-hyperbolic geodesic space, and $\phi_i: I_i=[a_i,b_i]\to X$ are $k$-quasigeodesics 
satisfying 
$$
d(x_1, x_2)\le r, d(y_1, y_2)\le r, x_i=\phi_i(a_i), y_i= \phi_i(b_i), i=1, 2,
$$
then the images $\phi_1(I_1), \phi_2(I_2)$ are  $D$-Hausdorff close.  
\end{lemma}
\proof Let $\phi_i^*$ be geodesics  connecting the end-points $x_i, y_i$ of $\phi_i$, $i=1,2$. Then, since quadrilaterals in $X$ as $2\delta$-slim, 
$$
\Hd(\phi_1^*, \phi_2^*)\le 2\delta+r. 
$$ 
Applying Lemma \ref{stab-qg}, we conclude that
$$
\Hd(\phi_1, \phi_2)\le D_{\ref{lem:sub-close} }(\delta,k,r)= 
\max(2D_{\ref{stab-qg}}(\delta, k) +\delta, D_{\ref{stab-qg}}(\delta, k) +\delta +r). 
$$
More explicitly, since $D_{\ref{stab-qg}}(\delta, k)= 92 k^2(k+3\delta)$, we get:
$$
D_{\ref{lem:sub-close} }(\delta,k,r)= 
\max(184 k^2(k+3\delta)  +\delta, 92 k^2(k+3\delta) +\delta +r). \qed 
$$

\begin{lemma}\label{lem:qi-preserves-paths}\label{lem:qi-preserves}
Suppose that $Y, X$ are path-metric spaces, $X$ is  $\delta$-Rips-hyperbolic and $f: Y\to X$ is a $(K,\eps)$-qi embedding. 
Then $Y$ is also $\delta_{\ref{lem:qi-preserves-paths}}(\delta, K, \eps)$-hyperbolic. 

In particular, (Rips) hyperbolicity is qi invariant among path-metric spaces: \newline If $f:~Y\to X$ is a $(K, \eps)$-quasiisometry and $X$ is $\delta$-hyperbolic 
in the sense of Rips, then $Y$ is $\delta'_{\ref{lem:qi-preserves-paths}}(\delta, K, \eps)$-hyperbolic in the sense of Rips.
\end{lemma}
\proof Consider a triple of points $x, y, z\in Y$ and for $\eta\ge 0$ take the triangle $\Delta_\eta$ in $Y$ formed by 
$\eta$-short arc-length parameterized paths $c(x,y), c(y,z), c(z,x)$ connecting these points. These paths are $(1,\eta)$-quasigeodesics in $Y$.  

Then $f(c(x,y)), f(c(x,z)), f(c(y,z))$ are $(K, \eps')$-quasigeodesics in $X$, for 
$$
\eps'= K(1+\eta)+ \eps,$$ 
see Lemma \ref{lem:compositions}. 
Hence, by Lemma \ref{stab-qg}, the quasigeodesic triangle formed by the quasigeodesic paths $f(xy), f(xz), f(yz)$ is  
$D=(2\cdot 92K^2(\eps' + 3\delta) +\delta)$-slim. 
This implies that the triangle $\Delta_\eta$ is $K(D+\eps)$-slim. 
Sending $\eta$ to $0$, we conclude that $Y$ is $\delta_{\ref{lem:qi-preserves-paths}}(\delta, K, \eps)$-hyperbolic 
for
$$
\delta_{\ref{lem:qi-preserves-paths}}(\delta, K, \eps)= K (184K^2(K + \eps + 3\delta) +\delta+\eps).  
$$ 
The second statement of the lemma follows from the first, combined with Lemma \ref{lem:quasi-inverse}. \qed

\begin{rem}
For $\eps=K$ we will use then notation $\delta_{\ref{lem:qi-preserves-paths}}(\delta, K, \eps)= \delta_{\ref{lem:qi-preserves-paths}}(\delta, K)$ and 
 $\delta'_{\ref{lem:qi-preserves-paths}}(\delta, K, \eps)= \delta'_{\ref{lem:qi-preserves-paths}}(\delta, K)$. Thus, for $\delta_{\ref{lem:qi-preserves-paths}}(\delta, K)$ we can take 
the number 
$$
K (184K^2(2K  + 3\delta) +\delta+ K).$$ 
\end{rem}

\medskip 
Lemma \ref{lem:qi-preserves-paths} combined with Lemma \ref{lem:gromov iff rips} immediately imply: 

\begin{cor}\label{cor:Rips-qi}
Hyperbolicity is qi invariant among  path-metric spaces. 
\end{cor}

\begin{cor}
Gromov-hyperbolicity is  qi invariant among quasi-path metric spaces.
\end{cor}
\proof Suppose that $X, Y$ are quasiisometric quasi-path metric spaces and $X$ is Gromov-hyperbolic. 
By Corollary \ref{cor: graph approx}, there are metric graphs $X', Y'$ which are geodesic metric spaces, and
are $(1,\epsilon)$-quasiisometric to $X$ and $Y$ respectively for a suitable
$\epsilon \geq 0$. Thus:

(1) by Lemma \ref{lem: gromov hyp qi inv} $X'$ is Gromov hyperbolic. 

(2) $X', Y'$ are quasiisometric geodesic metric spaces. 

However, then by Corollary \ref{cor:Rips-qi}, $X'$ is Rips-hyperbolic since it is a geodesic metric space. 
Then, by Lemma \ref{lem:qi-preserves}, $Y'$ is also Rips hyperbolic. 
Since $Y'$ is a geodesic metric space, again by Lemma \ref{lem:gromov iff rips},  
it is also Gromov-hyperbolic. Finally, $Y$ is also Gromov-hyperbolic by 
Lemma \ref{lem: gromov hyp qi inv}. \qed

\medskip As another corollary we get:  

 \begin{cor}\label{cor:hyp-th}
 Suppose that $Y$ is a path-metric space, $X\subset Y$ is a rectifiably connected $R$-net in $Y$ ($N_R(X)=Y$),
  equipped with a path-metric such that the inclusion map $X\to Y$ is  $\eta$-uniformly proper,  and 
 $X$ is $\delta$-hyperbolic. Then $Y$ is $\delta_{\ref{cor:hyp-th}}(\delta,\eta(2R+1),R)= 
 \delta_{\ref{lem:qi-preserves}}(\delta, \eta(2R+1))$-hyperbolic.  
 \end{cor}
 \proof According to Lemma \ref{lem:qi-net}, the inclusion map $X\to Y$ is an $L$-qi embedding, $L=\eta(2R+1)$, hence, an $L$-quasiisometry. Now the corollary follows immediately from Lemma \ref{lem:qi-preserves}. 
   \qed


 \section{Combings and a characterization of hyperbolic spaces} 

One of the key tools in our work is a characterization of hyperbolicity in terms of {\em slim combings} due to Bowditch. The idea is that if $X$ is a $\delta$-hyperbolic geodesic 
metric space, then for each pair of points we have a (typically non-unique) geodesic path $xy$ between these points and these paths satisfy the $\delta$-slim triangle property. 
Bowditch's characterization reverses this definition. 

\medskip
Let ${\mathcal P}(X)$ denote the space of paths in a topological space $X$.

\begin{defn}\label{defn:fellows} 
1. Two paths in a metric space are said to {\em Hausdorff $D$-fellow-travel} \index{Hausdorff fellow-traveling paths} 
if their images are $D$-Hausdorff-close. 

2. A {\em combing} ${\mathcal C}$ of a metric space $X$ is a map \index{combing} 
$$
c: X_0\times X_0\to {\mathcal P}(X)
$$
sending each pair $(x,y)\in X_0^2$ to a path $c_{x,y}$ in $X$ 
(also frequently denoted $c(x,y)$) connecting $x$ to $y$, where $X_0\subset X$ is a $D$-net for some $D$. 

3. For a function $C=C(r)$, a combing ${\mathcal C}$ is said to satisfy the 
{\em $C(r)$-Hausdorff fellow-traveling property}\index{fellow-traveling property} 
 if for every triple of points $x, y, z\in X_0$ with $d(y,z)\le r$,  
$$
\Hd(c_{x,y}, c_{x,z})\le C(r) 
$$
and $\diam(c_{x,x})\le C(0)$. 
\end{defn} 

While we define combings as maps,  we will think of each combing as a subset ${\mathcal C}= c(X^2)\subset {\mathcal P}(X)$.

\begin{lemma}\label{lem:fellow->proper}
Suppose that we have a combing ${\mathcal C}$  in a metric space $X$ such that:

1. Every path $c\in {\mathcal C}$ is a concatenation of at most $n$ subpaths each of which is in $\CC$ and is $\eta$-proper. 

2. The family $\CC$ satisfies the $C$-Hausdorff fellow-traveling property. 

3. The family $\CC$ is {\em consistent} in the sense that for each triple of points $x, y, z\in X$ such that $y\in c= c_{x,z}$, 
the subpath $c(x, y)$ of $c$ between $x$ and $y$ is $R$-Hausdorff close to the path $c_{x,y}$.

Then the family $\CC$  is uniformly proper. More precisely,   there is a function 
$$
\zeta_{\ref{lem:fellow->proper}}(r, C, R, \eta, n)$$ 
such that 
for each path $c: I\to X$ in $\CC$, for all $s, t\in I$, 
$$
d(c(s), c(t))\le r\Rightarrow d(s, t)\le \zeta(r, C, R, \eta, n). 
$$
\end{lemma}
\proof Consider a path $c=c_1\star ... \star c_n$ in $\CC$, $p=c(s), q=c_m(t), s<t$ and $d(p,q)\le r$. 
We assume that the domain $[0,T_{n+1}]$ of $c$ is the union of subintervals $[T_i, T_{i+1}]$, $i=1...,n$. 
By the consistency condition 3, the path $c_{p,q}$ is within the Hausdorff-distance $D$ from the subpath $c|_{[s,t]}$. By the definition of a combing, the path 
$c_{p,p}$ has diameter $\le C(0)$, while by the fellow-traveling property, 
$$
\Hd(c_{p,q}, c_{p,q})\le C(r). 
$$
Hence, by the $\eta$-properness assumption on the subpaths $c_i, i=2,...,n-1$, their domains have lengths $\le \eta(C(r)+D+R)$. For the same reason, 
$$
T_1- s\le R, t-T_n\le \eta(C(r)+C(0)+R). 
$$
 We, thus, obtain the required estimate  $t-s \le n\eta(C(r)+ C(0)+ R)=: \zeta(r, C, R, \eta, n)$.  \qed

\medskip

The following characterization of hyperbolicity is due to Bowditch, \cite[Proposition 3.1]{bowditch-unif-hyp}.

\begin{thm}
 \label{prop:bowditch-unif-hyp}
Given $h \geq 0$, there is  $k=k_{\ref{prop:bowditch-unif-hyp}}(h)<\infty$ with the following property.
Suppose that $X$ is a connected graph, and that for all $x, y \in V (X)$, we have
associated a connected subgraph, $Y_{xy} \subset X$ containing both $x$ and $y$, 
satisfying  the
following properties.

(1) For all $x, y \in V (X)$ with $d(x, y) \leq 1$, the diameter of $Y_{xy}$ in $X$ is at most $h$.

(2) For all $x, y, z \in V (X)$, $Y_{xy} \subset N_h ( Y_{xz} \cup Y_{zy})$.

Then $X$ is $k$-hyperbolic. In fact, we can take any $k \geq (3m- 10h)/2$, where $m$ is
any positive real number satisfying $2h(6 + \log_2 (m + 2)) \leq m$. Moreover, for all
$x, y \in V (X)$, the Hausdorff distance between $Y_{xy}$ and any geodesic from $x$ to $y$
is bounded above by $m-4h$.
\end{thm}

In the book we will be using the following corollary of Bowditch's characterization. We will refer to the family ${\mathcal C}$ 
of paths $c(x,y)$ appearing in the corollary as a {\em slim combing} of $X$. \index{slim combing} 
Note that Property (2) below amounts to the condition that the family $\CC$ is 
$D_2$-slim, i.e. $\delta_s(\CC)\le D_2$, see Definition \ref{defn:slimness constant}.

\begin{cor}\label{cor:bowditch-1}
Given $D_0\ge 0$, $D_1>0$, $D_2>0$ and a coarse Lipschitz function $\eta:\RR_{\geq 0}\map \RR_{\geq 0}$ 
there are $\delta=\delta_{\ref{cor:bowditch-1}}(D_0, D_1, D_2)\geq 0$, and $R=R_{\ref{cor:bowditch-1}}(D_0, D_1, D_2)\geq 0$ such that:

Suppose $X$ is a path-metric space and $X_0\subset X$ is a $D_0$-net  such that for each pair of points $x,y\in X_0$ we are 
given a rectifiable path $c(x,y)$ in $X$ joining $x,y$ with the following properties. 

Property (1). For all $x, y \in X_0$ with $d(x, y) \leq 1+2D_0$, the length of $c(x,y)$ in $X$ is at most $D_1$.

Property (2). For all $x, y, z \in X_0$, $c(x,y) \subset N_{D_2} ( c(x,z) \cup c(z,y))$.

\noindent Then $X$ is $\delta$-hyperbolic in the sense of Rips. Moreover, for all $x,y\in X_0$ and all $1$-quasigeodesics 
$\gamma_{xy}$ in $X$ joining $x,y$, $\Hd(\gamma_{xy},c(x,y))\leq R$.
\end{cor}
\proof 
Let $R=1+2D_0$. Since path metric spaces are $1$-quasi-path metric spaces, by Lemma \ref{lem:graph-approximation} 
there is a $(K_{\ref{lem:graph-approximation}}(1,R),\epsilon_{\ref{lem:graph-approximation}}(1,R))$-quasiisometry
$\iota: Z\map X_0\subset X$ which is the identity on the vertex set of $Z$, where $Z$ is the full subgraph of the
$R$-Rips graph of $X$ with the vertex set $X_0$. Also there is coarse inverse $\rho:X\map Z$ to $\iota$ which is the identity 
map when restricted to $X_0$ and and which is also a 
$(K_{\ref{lem:graph-approximation}}(1,R),\epsilon_{\ref{lem:graph-approximation}}(1,R))$-quasiisometry.
Let $k_0=K_{\ref{lem:graph-approximation}}(1,R)+\epsilon_{\ref{lem:graph-approximation}}(1,R)$. Then $\iota$ and
$\rho$ are $k_0$-quasiisometries.

Now we shall verify that $Z$ satisfies the conditions of Proposition \ref{prop:bowditch-unif-hyp} 
by  finding a suitable set of edge-paths joining each pair of vertices. Since hyperbolicity of path-metric spaces 
is invariant under quasiisometry it will then follow that $X$ itself is uniformly hyperbolic. Suppose $x,y\in V(Z)$ 
and $\alpha:[0,l]\map X$ is the arc-length parametrization of the path $c(x,y)$.
Let $n=\lfloor l\rfloor$. Consider the points $y_i=\alpha(i)$, $0\leq i\leq n$ and $y_{n+1}=\alpha(l)=y$ if $l>nD_0$. 
For each $i$, there is a point $x_i\in X_0=V(Z)$ such that $d_X(y_i,x_i)\leq D_0$. 
Then 
$$
d_X(x_i, x_{i+1})\leq  1+2D_0. 
$$
In particular, there is an edge in $Z$ joining $x_i, x_{i+1}$ for $0\leq i\leq n$. 
This defines an edge-path $\beta=\beta(x,y)$ in $Z$ joining $x,y\in V(Z)$ and the length of $\beta(x,y)$ is at most  $l+1$.

We take the subgraph $Y_{xy}\subset Z$ to be the image of the path $\beta(x,y)$. This is clearly a connected subgraph 
in $Z$. We now verify the conditions of Proposition \ref{prop:bowditch-unif-hyp} for this family of subgraphs.  

(1) Suppose $x,y\in V(Z)$ with $d_Z(x,y)\leq 1$. This implies $d_X(x,y)\leq R$. Hence the 
length of $c(x,y)$ is at most $D_1$ by property $(1)$. Thus the diameter of $Y_{xy}$ is at most $D_1+1$.

(2) If $V(\beta)$ denotes the vertex set $\{x_i\}$ of $\beta$ then we note that 
$$
\Hd(c(x,y), V(\beta))\leq 1+D_0.
$$
Thus, clearly, for all $x, y, z\in V(Z)$, we have 
$$
\beta_{xy}\subset N_{2+D_0+D_2}(\beta_{xz}\cup \beta_{zy}).
$$

Thus if we take $h=\max\{D_1+1, 2+D_0+D_2\}$, then $Z$ is $k_{\ref{prop:bowditch-unif-hyp}}(h)$-hyperbolic in the sense of Rips.
Since $\rho: X\map Z$ is a $k_0$-quasiisometry, $X$ is 
$\delta=\delta_{\ref{lem:qi-preserves}}(k_{\ref{prop:bowditch-unif-hyp}}(h),k_0)$-hyperbolic 
in the sense of Rips by Lemma \ref{lem:qi-preserves}.

For the second part of the corollary we set $m$ to be the least positive number satisfying $2h(6 + \log_2 (m + 2)) \leq m$.
Let $x\neq y\in X_0$ be arbitrary points.
Then by Theorem \ref{prop:bowditch-unif-hyp} any geodesic $\gamma_{xy}$ in $Z$ joining $x,y$ we have
$\Hd(\gamma_{xy}, \beta_{xy})\leq m-4h$. Suppose $\alpha_{xy}$ is a $1$-quasigeodesic joining $x,y$ in $X$.
Then clearly $\rho(\alpha_{xy})$ is a $2k_0$-quasigeodesic in $Z$. Thus, 
$\Hd(\gamma_{xy}, \rho(\alpha_{xy}))\leq D_{\ref{stab-qg}}(k_{\ref{prop:bowditch-unif-hyp}}(h),2k_0)$.

Hence, $\Hd(\beta_{xy}, \rho(\alpha_{xy}))\leq m-4h+D_{\ref{stab-qg}}(k_{\ref{prop:bowditch-unif-hyp}}(h),2k_0)$ and, therefore, 
$$
\Hd(V(\beta_{xy}), \rho(\alpha_{xy}))\leq 1+m-4h+D_{\ref{stab-qg}}(k_{\ref{prop:bowditch-unif-hyp}}(h),2k_0)=R_0.$$ 
But $V(\beta_{xy})=\rho(V(\beta_{xy}))$ since $\rho$ is the identity map when restricted to $X_0$. Hence, 
$$\Hd(\rho(V(\beta_{xy})),\rho(\alpha_{xy}))\leq R_0.$$
Since $\rho$ is a $k_0$-qi embedding we have $\Hd(V(\beta_{xy}), \alpha_{xy})\leq k_0(R_0+k_0)$.
The inequality $\Hd(c(x,y),V(\beta_{xy}))\leq D_0+1$ implies $\Hd(c(x,y), \alpha_{xy})\leq 1+D_0+k_0(R_0+k_0)$.
Hence we may set $R=1+D_0+k_0(R_0+k_0)$.
\qed

\begin{cor}\label{cor:bowditch}
Given $D_0\ge 0$, $D>0$, and a coarse Lipschitz function $\eta:\RR_{\geq 0}\map \RR_{\geq 0}$,  
there are $\delta=\delta_{\ref{cor:bowditch}}(D_0, D, \eta)\geq 0$, and $K=K_{\ref{cor:bowditch}}(D_0, D, \eta)<\infty$ such that:

Suppose $X$ is a path-metric space and $X_0\subset X$ is a $D_0$-net  such that for each pair of points $x,y\in X_0$ we are 
given a rectifiable arc-length parametrized path $c(x,y)$ in $X$ joining $x,y$ with the following properties:  

Property (1). For all $x, y \in X_0$ the path $c(x,y)$ is  $\eta$-proper.

Property (2). For all $x, y, z \in X_0$, $c(x,y) \subset N_{D} ( c(x,z) \cup c(z,y))$.

\noindent Then $X$ is $\delta$-hyperbolic in the sense of Rips and moreover, the paths $c(x,y)$ are $K$-quasigeo\-desics in $X$.
\end{cor}
\proof This follows from Corollary \ref{cor:bowditch-1}. Property 1 of  Corollary \ref{cor:bowditch-1} is verified
with $D_1=\eta(2D_0+1)$ and property 2 is verified with $D_2=D$. Thus we can take
$\delta= \delta_{\ref{cor:bowditch}}(D_0,D,\eta)=\delta_{\ref{cor:bowditch-1}}(D_0, \eta(2D_0+1), D)$.

For the second part of the corollary suppose $\gamma_{xy}$ is an arc-length parametrized path in $X$ joining $x,y\in X_0$
such that $\length(\gamma_{xy})\leq 1+d(x,y)$. Then clearly it is a $1$-quasigeodesic. By the second part of Corollary
\ref{cor:bowditch-1} we have 
$$
\Hd(\gamma_{xy}, c(x,y))\leq R_{\ref{cor:bowditch-1}}(D_0, \eta(2D_0+1),D)=R.
$$
Then, by Lemma \ref{lem:quasigeodesic-paths}, the paths $c(x,y)$ are 
$K_{\ref{lem:quasigeodesic-paths}}(R,1,\eta)$-quasigeodesics. \qed

\begin{rem}
1. Suppose that the assumption in Part (b) of this lemma holds. Then each path $c(x,x)$ has length $\le D_1:=\eta(0)$, 
i.e. the condition (1) in Part (a) necessarily holds as well.  

2. In  the proofs of hyperbolicity of various spaces given in this books, based on Corollary \ref{cor:bowditch}, we first verify the uniform properness of the 
paths $c(x,y)$ and, frequently, also verify that they satisfy the {\em Hausdorff fellow-traveling condition}, before 
proving that (a2) holds. 

3. In our proofs, instead of  using arc-length parameterizations of the paths $c(x,y)$ we will be using some uniformly quasiisometric reparameterizations of these paths. 
Clearly, uniform properness of one implies uniform properness of the other. 

4. We will refer to a family of paths $c$ as satisfying the assumptions of the corollary as a {\em slim combing} of $X$. 
\end{rem}

Lastly, we generalize this corollary to the case of discrete paths $c: [m, n]\cap \Z\to X$. We will be using the extension $\tilde{c}$ of maps $c$ to the real interval $[m, n]$ defined by sending the open interval $(i, i+1)$ to $c(i)$. 

\begin{cor}
Given $r>0$, $D_0\ge r$, $D_1>0$, and a function $\eta:\RR_{\geq 0}\map \RR_{\geq 0}$ 
there are $\delta=\delta(D_0, D_1,\eta)\geq 0$, and $K=K(D_0, D_1, \eta)<\infty$ such that:

(a) Suppose $X$ is a $r$-quasi-path metric space and $X_0\subset X$ is a $D_0$-net 
such that for each pair of points $x,y\in X_0$ we are 
given a discrete $r$-path $c(x,y)$ in $X$ joining $x,y$ with the following properties:  

Property (a1). For all $D\geq 0$ and $x, y \in X_0$ with $d(x, y) \leq D$, the length of $c(x,y)$ in $X$ is at most $\eta(D)$.

Property (a2). For all $x, y, z \in X_0$, $c(x,y) \subset N_{D_1} ( c(x,z) \cup c(z,y))$.

\noindent Then $X$ is $\delta$-hyperbolic in the sense of Gromov.

(b) Moreover, if paths $c(x,y)$ are $\eta$-proper, then each map $\tilde{c}$ is 
a $K$-quasigeo\-desic in $X$.  
\end{cor}

\proof By Corollary \ref{cor: graph approx}, there is a $(1,\epsilon)$-quasiisometry $\rho:X\map X'$, 
where $X'$ is a metric graph all whose each edges are of length $r$ and 
$\epsilon=\epsilon_{\ref{cor: graph approx}}(r)$. Therefore,
$X'$ is a geodesic metric space. Now, define $X'_0= \rho(X_0)$. We take any map
$g: X'_0\map X_0$ such that $x=\rho(g(x))$ for all $x\in X'_0$. Then for all $x,y\in X'_0$
we define a path $c'(x,y)$ as follows. Let $c_{xy}: I_{xy}\cap \mathbb Z\map X$ be the
parametrization of $c(g(x),g(y))$. For any two consecutive points $s,t\in I_{xy}\cap \mathbb Z$
we join $\rho\circ c_{xy}(s), \rho\circ c_{xy}(t)$ by a geodesic in $X'$. 
Concatenation of these forms a path $c'(x,y)$.  We leave it to the reader to verify that the two properties of Corollary
\ref{cor:bowditch} hold for this family of paths in $X'$. This implies that $X'$ is uniformly Rips-hyperbolic
and, hence, uniformly Gromov-hyperbolic by Lemma \ref{lem:gromov iff rips}. Since $X$ is $(1,\epsilon)$-qi to
$X'$, it is also uniformly Gromov-hyperbolic by Lemma \ref{lem: gromov hyp qi inv}. The last part of the corollary follows from 
Lemma \ref{lem:quasigeodesic-paths}. 
\qed

\section{Hyperbolic cones} 
\label{sec:hyperbolic cone}

Suppose that $(Z, d_Z)$ is a path-metric space. In this section we define the {\em hyperbolic cone}
\index{hyperbolic cone}
 $Z^h$ over $Z$. This definition will be used in Chapter 9 when discussing relatively hyperbolic spaces

As a topological space, $Z^h$ is the product $Z\times [1,\infty)$, where we identify $Z$ with $Z\times \{1\}$. We equip $Z^h$ with the length structure (imitating the description of the hyperbolic metric on horoballs in the real-hyperbolic space): Paths in this length structure are concatenation of vertical and horizontal path, with respect to the product decomposition of $Z^h$. Given two points $y_1, y_2\in [1,\infty)$, and $z\in Z$, we let the length of the interval between $(z,y_1), (z,y_2)$ in  $\{z\}\times [1,\infty)$ equal 
$|\log(y_2/y_1)|$. We let the length of each horizontal path, contained in the ``horosphere'' $Z\times \{y\}$
equal $y^{-1}$ times the length of the corresponding path in $Z$. This length structure defines a path-metric on $Z^h$.

\begin{rem}
We refer the reader to the book of Roe \cite[2.5]{Roe} and the paper by Bowditch \cite{bowditch-relhyp} for alternative definitions of hyperbolic cones. For instance, Roe's construction 
works for general metric spaces $(Z, d_Z)$. 
\end{rem}

  \begin{figure}[tbh]
\centering
\includegraphics[width=70mm]{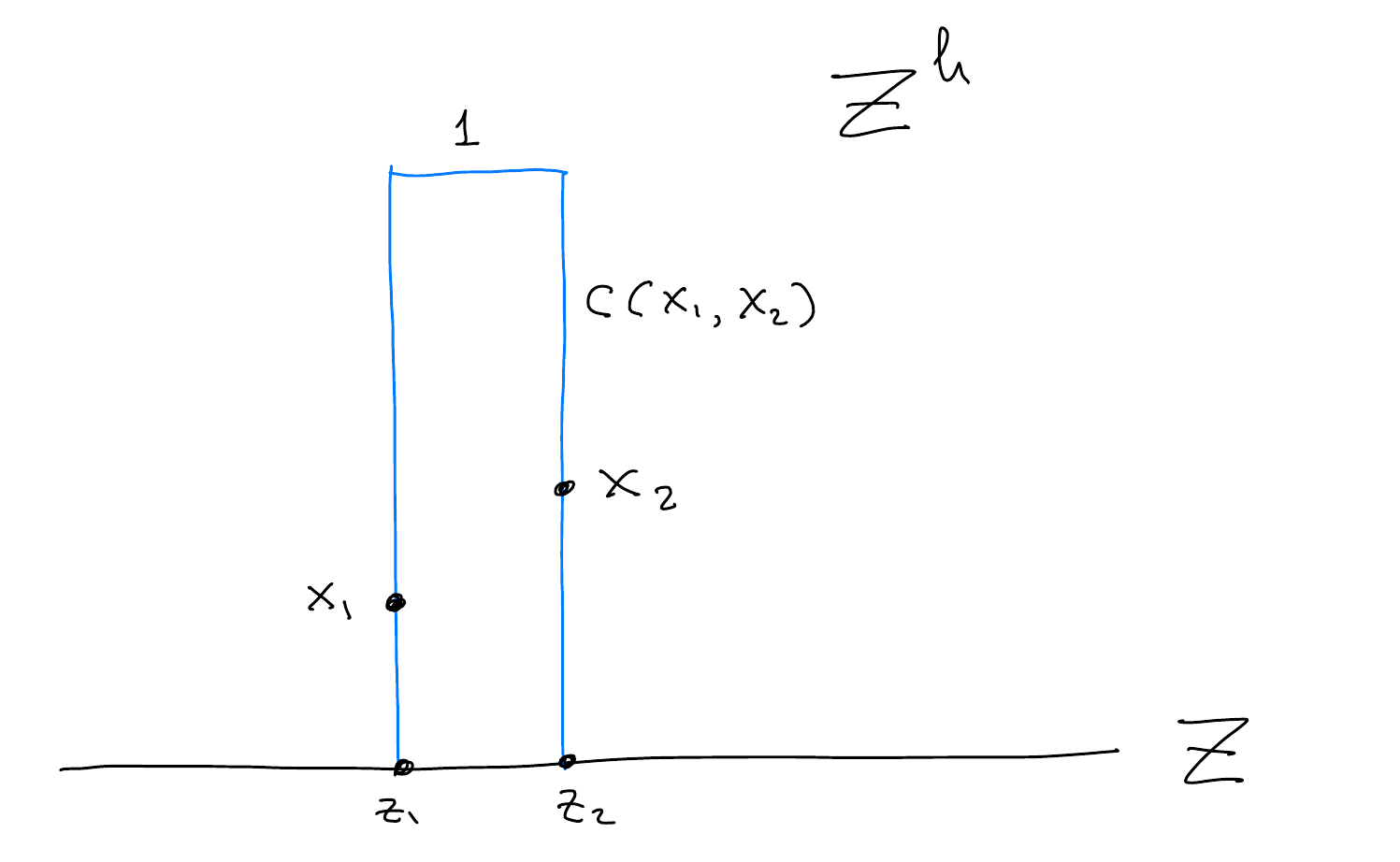}
\caption{Combing of hyperbolic cone}
\label{fig:hyp cone}
\end{figure}

\begin{prop}\label{prop:hyper-cone}
The metric space $Z^h$ is $\delta$-hyperbolic for some uniform constant $\delta$. 
\end{prop}
\proof We will describe a slim combing on $Z^h$. 
Each combing path will be a concatenation of at most two vertical paths and at most one horizontal path. Consider two points 
$x_1=(z_1, y_1), x_2=(z_2, y_2)$ in $Z^h$. If $z_1=z_2$ then the map $c({x_1,x_2})$ connecting $x_1$ to $x_2$ will be the unique vertical interval connecting these points. 
Suppose that $y_1\le y_2$. Find the smallest $y\ge y_2$ such that 
$$
y^{-1} d_Z(z_1, z_2)\le 1. 
$$
Set $x'_i:= (z_i, y)$, $i=1, 2$. Then  $c({x_1,x_2})$ is the concatenation 
$$
[x_1 x'_1]\star [x'_1 x'_2] \star [x'_2 x_2],
$$
where the first and the last segments are vertical intervals between $x_i, x'_i$, $i=1, 2$, and the middle segment $[x'_1 x'_2]$ is any path in  $Z\times \{t\}$ connecting $x_1'$ to $x_2'$ and having length $\le 1$.

We now verify the slim combing properties of the paths $c$, as required by Corollary \ref{cor:bowditch}. 

{\bf 1. Uniform properness.} We define two projections, $\pi_1: Z^h\to Z, \pi_2: Z^h\to [1,\infty)$,
$$
\pi_1((z,y))=z, \quad \pi_2((z,y))=y. 
$$
Then the composition $\log\circ \pi_2$ is $1$-Lipschitz, while the first projection has the property that if $\beta: [0,1]\to Z^h$ is a path such that $\pi_2(\beta([0,1]))\le t_0$, then 
$\length(\pi_1\circ \beta))\le \exp(t_0) \length(\beta)$.   We will use these two observations to verify uniform properness of the paths $c$. 

First, we note that it suffices to estimate from above the length of the path $c=c({x_1,x_2})$ in terms of the distance between $x_1=(z_1, y_1), x_2=(z_2, y_2)$. 

a. Suppose, first that  $y_1\ne y_2$. Then, since $\log\circ \pi_2$ is  $1$-Lipschitz, 
$$d(x_1, x_2)\ge d(\pi_2(x_1), \pi_2(x_2))\ge \length(c).$$ 

b. Thus, we only have to consider the case $y_1=y_2$. Furthermore, if 
$$y_1^{-1}d_Z(z_1,z_2)\le 1,$$ then 
$\length(c)\le 1$ as well. Hence, without loss of generality,  $y_1^{-1}d_Z(z_1,z_2)\ge 1$. 
Setting $D:= d_Z(z_1, z_2)$, we get $y=D$ and $c$ is the concatenation 
$$
[x_1 x'_1]\star [x'_1 x'_2] \star [x'_2 x_2],
$$
where $x'_i:= (z_i, y)$, $i=1, 2$. The overall length of $c$ is 
\begin{equation}\label{eq:length-of-c}
1+ 2\log(y/y_1)=1+ 2\log(D/y_1).
\end{equation}
 Let $\beta$ be a path connecting $x_1, x_2$. Define 
$$
y_0:=\max \pi_2\circ \beta.$$
Since $\log\circ \pi_2$ is  $1$-Lipschitz, 
$$
\log(y_0/y_1)\le \length(\beta), \quad \log(y_0) \le \length(\beta) +\log(y_1). 
$$
Using the projection $\pi_1$, we see that 
$$
D\le \length( \pi_1(\beta))\le y_0 \length(\beta)\le y_0 \le y_1 \exp(\length(\beta)).  
$$  
Hence, 
$$
\length(c)\le 1+ 2\log(D/y_1) \le 1+ 2 \length(\beta).  
$$
It follows that $\length(c)\le 1+ 2 d(x_1, x_2)$, as required.

\medskip
{\bf 2. Slim triangle condition.} Consider three points $x_i=(z_i,y_i)\in Z^h$, $i=1, 2, 3$. We let $y_{ij}$ denote the maximum of $\pi_2\circ c(x_i, x_i)$, $i\ne j$. After relabeling the points, we can assume that 
$$
y_{12}\le y_{23}\le y_{31}. 
$$
Define $x_1':=(z_1, y_{12})$ and $x'_2:= (z_2,y_{12})$. 

Replacing the points $x_1, x_2$, respectively, with $x_1'=(z_1, y_{12})$ and $x'_2= (z_2,y_{12})$ we see that the 
$c(x_1, x_2), c(x_2, x_3), c(x_3,x_1)$ is $\delta$-slim if and only if 
the triangle formed by the paths $c(x'_1, x'_2), c(x'_2, x_3), c(x_3,x'_1)$ is $\delta$-slim. Thus, without loss of generality, $x_1'=x_1, x_2'=x_2$ and the path $c(x_1,x_2)$ 
is horizontal,  of length $\le 1$. We claim that in this situation, the paths  $c(x_1, x_3), c(x_2, x_3)$ are uniformly Hausdorff-close, which, of course, will imply the uniform slimness. Denoting 
$$
x''_i:= (z_i,y_{23}), i=1, 2, 3, 
$$
we see that the parts $\Hd(c(x_1, x''_1), c(x_2, x_2''))\le 1$, while the parts between $x_3''$ and $x_3$ of 
 $c(x_1, x_3), c(x_2, x_3)$ are equal.    Thus, we need to bound the Hausdorff distance between 
$c(x''_1, x''_3), c(x''_2, x''_3)$. 

We have:
$$
1\ge y^{-1}_{12} d(z_1,z_2)\ge y^{-1}_{23} d(z_1,z_2). 
$$
Applying the triangle inequality in $(Z,d_Z)$, we obtain
$$
2=1+1 \ge y^{-1}_{23}( d(z_1,z_2) + d(z_2,z_3))\ge y^{-1}_{23} d(z_3,z_1).  
$$
It follows that $y_{23}\le y_{31}\le 2 y_{23}$. Hence, 
$$
d(x_1'', (z_1, y_{31}))\le \log(2). 
$$ 
It follows that the length of $c(x''_1, x''_3)$ is $\le 1+ 2\log(2)$, while the length of $c(x''_2, x''_3)$ is $\le 1$. Clearly, these path are at the Hausdorff-distance $\le 2+ 2\log(2)$. 
\qed

\begin{rem}\label{rem:qg-in-horoball}
By Corollary \ref{cor:bowditch}, the slim combing paths $c(x_1, x_2)$ in $Z^h$ defined in this proof are  $k$-quasigeodesics in $Z^h$ with the quasigeodesic constant $k$ independent of $Z$. 
\end{rem}

\begin{lem}\label{lem:strictly-convex} 
Suppose that $z\in Z$ is within distance $\le C$ from a segment $[z_1 z_2]_{Z^h}$ with the end-points $z_1, z_2$ in $Z$. 
Then  $d_Z(z, \{z_1, z_2\})\le C'$, where $C'$ depends only on $C$.  
\end{lem}
\proof The combing path $c(z_1, z_2)$ is uniformly Hausdorff-close to the geodesic  $[z_1 z_2]_{Z^h}$. Hence, it suffices to prove the lemma with $[z_1 z_2]_{Z^h}$ replaced by the combing path $c(z_1, z_2)$. Then the inequality $d_{Z^h}(z, c(z_1, z_2))\le C$ implies that 
$$
d_Z(z, \{z_1, z_2\})\le 2C+1. 
$$
\qed

\begin{prop}\label{prop:exp-distortion}
The inclusion map $Z\to Z^h$ is $\eta$-uniformly proper, where $\eta(t)=a\exp(at)$ and $a$ is a universal constant. 
\end{prop}
By the above remark, there is a universal constant $L$ such that for any two points $z_1, z_2\in Z$ the combing path 
$c=c(z_1, z_2)$ is an $L$-quasigeodesic in $Z^h$. By the equation \eqref{eq:length-of-c}, unless $D=d_Z(z_1,z_2)\le 1$, 
the length of this path is $1+ 2\log (D)$ since in our situation $y_1=1$. In the exceptional case when $D\le 1$, 
the path $c$ is horizontal and has length $D$. Thus,
$$
d_{Z^h}(z_1, z_2)\ge L^{-1} \length(c) \ge \min(2 L^{-1}\log (D), L^{-1} D), 
$$
and, hence,
$$
D\le \exp(\frac{1}{2}L d_{Z^h}(z_1, z_2)) + Ld_{Z^h}(z_1, z_2)\le a \exp (a d_{Z^h}(z_1, z_2)) 
$$
for $a=\max(1, \frac{1}{2}L)$.  \qed 

\section{Geometry of hyperbolic triangles} 

Informally speaking, triangles in hyperbolic spaces resemble triangles in trees, i.e. tripods. The {\em comparison map} to trees makes this precise and allows one to reduce proofs of various statements about hyperbolic triangles to that of tripods in trees. Below, by a {\em tree} we mean a regular simplicial tree $T$ of valence $\ge 3$, 
equipped with the standard graph-metric. However, any real tree (a $0$-hyperbolic geodesic metric space) not isometric to an interval, will work just as well. 
For any three numbers $a_1, a_2, a_3$ satisfying the triangle inequalities $a_i\le a_j + a_k, \{i, j, k\}=\{1, 2, 3\}$, 
there exists a triangle $\Delta\subset T$ with the side-lengths $a_1, a_2, a_3$. Accordingly, for each triangle 
$\Delta= \Delta x_1x_2x_3$ in a metric space $X$, we define its {\em comparison triangle} \index{comparison triangle}
$$
\bar{\Delta}= \Delta \bar x_1 \bar x_2 \bar x_3
$$
in  $T$, as a triangle in $T$ such that
$$
d(x_i, x_j)= d_T(\bar{x}_i, \bar{x}_j), \quad 1\le i, j\le 3. 
$$

For each point $p$ in the side $x_ix_j$ of $\Delta$, define its {\em comparison point} \index{comparison point} 
$\bar{p}\in \bar{x}_i \bar{x}_j\subset \bar\Delta$ by the condition
$$
d(p, x_i)= d_T(\bar{p}, \bar{x}_i). 
$$
Thus, we get the {\em comparison map} \index{comparison map}
$\theta: \Delta\to \bar\Delta$, $\theta(p)=\bar{p}$, which   restricts to an isometry on each side of $\Delta$. 
The internal points of $\Delta$ are the points of the $\theta$-preimage of the center of $\bar\Delta$. A triangle $\Delta$ is $\delta$-thin if and only if the diameters of all fibers of $\theta$ are $\le \delta$. For each pair of points  $p, q\in \Delta$ in a 
$\delta$-thin triangle $\Delta\subset X$, triangle inequalities imply: 
$$
d(p, q)- 2\delta \le d_T(\bar{p}, \bar{q})\le d(p, q)+\delta.
$$
Thus, the map $\theta$ is a $(1, 2\delta)$-quasiisometry $\Delta\to \bar\Delta$ for $\delta$-thin triangles $\Delta\subset X$.

\begin{defn}
Let $(X,d)$ be a metric space. A point $p\in X$ is said to be a $C$-{\em center} of a geodesic triangle $\Delta$ if  
$p$ lies within distance $C$ from all three sides of the triangle. \index{$C$-center of a triangle} 
A $0$-center  is simply called a center of $\Delta$. We will use this definition  when $(X,d)$ is a tree, in which case the $0$-center is unique. 
\end{defn}

Note that if a geodesic metric space $X$ is $\delta$-hyperbolic, every geodesic triangle $\Delta$ 
in $X$ has a $\delta$-center, e.g. a point on one side of $\Delta$, within distance $\delta$ from the two other sides. The internal points of a $\delta$-thin triangle are mapped via the comparison map $\theta$ to the center of the comparison triangle in the tree and each internal points is a $\delta$-center of $\Delta$. It follows directly from the definition that every side of every $\delta$-slim triangle $\Delta$ contains a $\delta$-center of $\Delta$.

\begin{defn}\label{defn:tripod} \index{$C$-tripod} 
{\em A $C$-tripod} in $X$  is the union $T_p(xyz)$ of three geodesic segments $px\cup py \cup pz$, where $p$ is a $C$-center of a geodesic triangle $\Delta xyz$. The points $x, y, z$ are called the {\em extremities} of the tripod and the segments 
$px, py, pz$ are the {\em legs} of the tripod. If $X$ is a tree and $C=0$ then 
by a {\em tripod} in $X$ we mean a $0$-tripod and the center of a tripod means its $0$-center. 
\end{defn}

As noted above, if a triangle $\Delta=\Delta x y z$ is $\delta$-slim, then there exists a point  $p\in x y$ such that the union $px\cup py \cup pz$ is a $\delta$-tripod.

The next lemma follows immediately from the definitions:

\begin{lemma}\label{lem:C-tripods}
If $X$ is $\delta$-hyperbolic in the sense of Rips, then each $C$-tripod  $T_p(xyz)\subset X$ is $C+\delta$-Hausdorff close to the triangle $\Delta xyz$. 
\end{lemma}

\begin{lemma}\label{lem:K-centers}
If $X$ is $\delta$-hyperbolic in the sense of Rips, then for each $C$-tripod  $T_p(xyz)$, 
the point $p$ is a $C+2\delta$-center of every triangle $\Delta x' y' z'$, where $x'\in px, y'\in py, z'\in pz$. 
\end{lemma}
\proof Consider a point $u\in xy$ within distance $C$ from $p$. By the $2\delta$-slimness of the quadrilateral $xx'y'y$, either 
there exists a point $u'\in x'y'$ at distance $\le 2\delta$ from $u$, or there is a point $v\in xx'\cup yy'$ at distance $\le \delta$ from $u$. In the latter case, for some $u'\in \{x', y'\}$, we get  $d(p,u')\le C+\delta$. Thus, in each case, there is $u'\in x'y'$ within distance $C+2\delta$ from $p$. \qed

\begin{lemma}\label{lem:centers} 
Suppose that $X$ is a $\delta$-hyperbolic geodesic metric space. 
If $p, q$ are $C$-centers of the same geodesic triangle, 
then $d(p,q)\le D_{\ref{lem:centers} }(\delta,C)$. 
\end{lemma}
\proof We will use the comparison map $\theta: \Delta\to \bar\Delta$. Since $X$ is $\delta$-hyperbolic, the triangle $\Delta$ will be $6\delta$-thin (Lemma \ref{hyp-defn}) and, hence, $\theta: x\mapsto \bar{x}$ satisfies the inequalities
$$ 
d(a, b)- 12\delta \le d_T(\bar{a}, \bar{b})\le d(a, b)+ 6\delta,
$$
 $a, b\in \Delta$. We now prove the lemma. Let $p$ be a $C$-center of $\Delta$ and $p_1, p_2, p_3$ be the points on the sides of $\Delta$ within distance $C$ from $p$. 
 Then $d(p_i, p_j)\le 2C$, $1\le i< j\le 3$. We will estimate the distances from the points $p_i$ to the internal points $c_i$ of $\Delta$, where we label the points so that $p_i, c_i$ lie on the same side of $\Delta$. Let $\bar{c}\in \bar\Delta$ denote 
 the center of $\bar\Delta$, $\bar{c}=\theta(c_i)$, $i=1, 2, 3$.     Then 
 $$
d_T(\theta(p_i), \theta(p_j))\le 2C + 6\delta.
 $$
 It is easy to see that all three points $\theta(p_i)$ cannot lie on the same leg of the tripod $\bar\Delta$, unless one of them equals to the center of 
 $\bar\Delta$. Thus, there exists $i$ such that $d_T(\theta(p_i), \bar{c})\le \frac{1}{2}(2C + 6\delta)= C + 3\delta$.  
 Since $\theta$ is an isometry on each side of $\Delta$, we then obtain 
 $$
 d(p_i, c_i)\le C + 3\delta 
 $$
and, hence, $d(p, c_i)\le 2C + 3\delta$ for one of the internal points of $\Delta$. Since the internal points are distance $\le 6\delta$ apart, we conclude that for any two $C$-centers $p, q$ of $\Delta$ 
$$
d(p,q)\le D_{\ref{lem:centers} }(\delta,C):= 2C + 9\delta. \qed 
$$

\begin{lemma}\label{lem:almost-geodesic} 
Suppose that $X$ is $\delta$-hyperbolic. For a geodesic triangle $\Delta xyz$ suppose that 
$x'\in xz, y'\in yz$ are equidistant from $z$ and satisfy the inequality $d(x',y') > 2\delta$. 
Then the path $xx' \star x'y' \star y'y$ is 
$C_{\ref{lem:almost-geodesic}}(\delta)$-quasigeodesic in $X$. 
\end{lemma}
\proof First we show that $x'\in N_{\delta}(xy)$. 
By $\delta$-hyperbolicity of the $\Delta xyz$, $x'\in N_{\delta}(xy\cup yz)$. If possible suppose
$d(w,x')\leq \delta$ for some $w\in yz$. Then $|d(z,w)-d(x',w)|\leq d(w,x')\leq \delta$.
This implies $d(y',w)\leq \delta$, since $d(z,x')=d(z,y')$. It follows that $d(x',y')\leq d(w,x')+d(w,y')\leq 2\delta$, 
 a contradiction. Hence, $x'\in N_{\delta}(xy)$. Let $x''\in xy$ be any point with $d(x',x'')\leq \delta$.

Next we claim that $y'\in N_{5\delta}(x''y)$.
Since the quadrilateral $\square x'x''yz$ is $2\delta$-slim, 
$$
y'\in N_{2\delta}(x'x''\cup x''y\cup x'z).
$$
If $y'\in N_{2\delta}(x'x'')$ then clearly $d(x'',y')\leq 3\delta$. As in the previous paragraph, 
if $y'\in N_{2\delta}(x'z)$ then $d(x',y')\leq 4\delta$ and thus $d(y', x'')\leq 5\delta$.
If neither of these happen then $y'\in N_{2\delta}(x''y)$. This proves our claim. Let
$y''\in x''y$ be any point such that $d(y',y'')\leq 5\delta$.

Let $\alpha:[0,l]\map X$ be the arc-length parametrization
of $xx'* x'y'* y'y$. Suppose $s, t\in [0,l]$ and $s<t$. Let $x_1=\alpha(s), y_1=\alpha(t)$. 
We already have $d(x_1,y_1)\leq t-s$. 

If both $x_1, y_1$ are on the same segment then there is nothing to prove. So assume otherwise.
Suppose $x_1\in xx'$ and $y_1\in x'y'$.
Since the quadrilateral $\square x'y'y''x''$ is $2\delta$-slim,  there is a point $y'_1\in x''y''$ such that 
$d(y_1, y'_1) \leq 7\delta$. Also by the $\delta$-slimness of $\Delta xx'x''$, there is a point
$x'_1\in xx''$ such that $d(x_,x'_1)\leq 2\delta$.
It follows that 
\begin{align*}
d(x_1,y_1)\geq d(x'_1,y'_1)-d(x_1,x'_1)-d(y_1,y'_1)\geq\\ 
d(x'_1,x'')+d(x'',y'_1)-9\delta
\geq \\
d(x_1,x')-d(x_,x'_1)-d(x',x'')+d(x',y_1)-d(x',x'')-d(y_1,y'_1)-9\delta
\geq\\
d(x_1,x')+d(x',y_1)-20\delta=t-s-20\delta. 
\end{align*}
 In the same way, if
$x_1\in x'y', y_1\in y'y$ we have $d(x_1,y_1)\geq t-s-20\delta$. Next suppose $x_1\in xx'$ and $y_1\in y'y$.
Then there are points $x'_1\in xx''$ and $y'_1\in y''y$ such that $d(x_1,x'_1)\leq 2\delta$ and
$d(y_1,y'_1)\leq 6\delta$ by the $\delta$-slimness of the triangles $\Delta xx''x'$ and $\Delta y'y''y$ respectively.
It follows that 
\begin{align*}
d(x_1,y_1)\geq d(x'_1,y'_1)-d(x_1,x'_1)-d(y_1,y'_1)
\geq \\
d(x'_1,x'')+d(x'',y'')+d(y'',y'_1)-8\delta
\geq \\
d(x_1,x')- d(x_1,x'_1)- d(x',x'')+ d(x',y')- d(x',x'')- \\ 
d(y',y'')+ d(y',y_1)- d(y',y'')- d(y_1,y'_1)- 8\delta
\geq \\
d(x_1,x')+ d(x',y')+ d(y',y_1)-28\delta= t- s- 28\delta. 
\end{align*}

Hence, the arc-length parametrization of $xx' \star x'y' \star y'y$ is a
$(1,28\delta)$-quasigeodesic. In particular, we can take
$C_{\ref{lem:almost-geodesic}}(\delta)=1+28\delta$.  \qed

\medskip 
As an application of the lemma we will prove that two geodesics $xz, zy$ either are uniformly close to each other 
on ``long subintervals'' $x'z, y'z$, or their concatenation $xz\star zy$ is a uniform quasigeodesic in $X$:

\begin{lemma}\label{lem:g-concat}
Suppose that $R\le \min(d(x,z), d(y,z))$ is such that points $x'\in xz, y'\in yz$ at the distance $R$ from $z$ satisfy $d(x',y')>2\delta$. Then the concatenation    $xz\star zy$ is an $L_{\ref{lem:g-concat}}(R,\delta)$-quasigeodesic in $X$. 
\end{lemma}
\proof According to the previous lemma, the concatenation $xx' \star x'y' \star y'y$ is a 
$C_{\ref{lem:almost-geodesic}}(\delta)$-quasigeodesic in $X$. We regard $xz\star zy$, 
$xx' \star x'y' \star y'y$ as paths $c, c'$ on intervals $[0, T]$, $[0,T']$, parameterized by arc-length and connecting $x$ to $y$, 
where $T=d(x,z)+d(z,y)]\le T'\le T+ 2R$.  Then for each $t\in [0,T]$
$$
d(c(t), c'(t))\le 2R.
$$
Since $c'$ is a $C_{\ref{lem:almost-geodesic}}(\delta)$-quasigeodesic, Lemma follows.   \qed 

\medskip

The next lemma is a consequence of Lemma \ref{lem:approximation}, but we will give a direct proof: 

\begin{lemma}
Suppose that  $\Delta xyz$ is a geodesic triangle in a $\delta$-hyperbolic space $X$ and $d(y,z)\le C$. 
Then there is a monotonic map $f: xz\to xy$ such that $d(f, \id)\le 2(\delta+C)$. 
\end{lemma}
\proof We define $f$ as follows. 

1. Points $p\in xz$ such that $d(x,p)\le \min(d(x,z), d(x,y))$, will be sent to $\bar{p}\in xy$ such that  
$d(x,\bar{p})= d(x,p)$. 

2.  Points $p$ such that $d(x,p)\ge d(x,y)$ will be all sent to $y$. 

\smallskip 
\noindent Thus, the map $f: xz\to xy$ is monotonic but possibly discontinuous at $z$. We next estimate the distances $d(p, f(p))$. Since $X$ is $\delta$-hyperbolic, there exists a point $p'\in xz \cup yz$ within distance $\delta$ from $p$. 

Case 1: $p'\in xy$ and  $d(x,p)< \min(d(x,y), d(x,z))$. Then by the triangle inequalities, 
$$
d(x,p)-\delta \le d(x, p')\le d(x,p)+\delta. 
$$
In particular, $d(p',\bar{p})\le \delta$ and, thus, $d(p, f(p))\le 2\delta$. 

 Case 2: $d(x,p)\ge \min(d(x,z), d(x,y))$, which implies that $d(p,z)\le C$, $d(p,y)\le 2C$. Thus, $d(p, f(p))\le 2C$. 
 
 Case 3. $p'\in yz$ and  $d(x,p)< \min(d(x,y), d(x,z))$. Since $p'\in yz$ and $d(y,z)\le C$, we get $d(p,y)\le \delta+C$. By the same argument as in Case 1, 
 $$
d(x, \bar{p})-C-\delta \le d(x, y)\le d(x,\bar{p})+\delta+C, 
$$
 with the point $y$ playing the role of $p'$ in Case 1. Thus, $d(y,\bar{p})\le \delta+C$, which implies
 $$
 d(p, f(p))\le 2(\delta+C). \qed 
 $$
 
 \medskip 

We will apply this construction in the following situation. Let $\Delta= \Delta x_1 x_2 x_3$ be a geodesic triangle in $X$, $z$ be its $\delta$-center, $y\in x_1x_2$ a point within distance $\delta$ from $z$. 
We then use the maps $f_1: x_1z\to x_1y, f_2: x_2z\to x_2y$ as in the lemma. These maps combine to define a map 
$f: x_1z \cup zx_2\to x_1x_2$. We parameterize concatenation $\al=x_1z \star zx_2$ by its arc-length. 
Then the map $f$  is monotonic (with respect to the parameterization), fixes the endpoints of $\al$ and satisfies
$$
d(p, f(p))\le 4\delta, p\in  \al. 
$$

We, thus, obtain:

\begin{cor}\label{lem:proj}
Suppose that $z$ is a $\delta$-center of a geodesic triangle $\Delta= \Delta x_1 x_2 x_3$ in a $\delta$-hyperbolic space $X$. Then there exists a monotonic map 
$$
f: x_1z \star zx_2 \to x_1x_2
$$
such that $d(f, \id)\le 4\delta$. 
\end{cor}

\section{Ideal boundaries}\label{sec:ideal boundaries}

In the book we will be mostly working with Gromov's notion of ideal boundary (Gromov-boundary)
 $\geo Z$  of a space $Z$, which is hyperbolic in the sense of Gromov. The elements of $\geo Z$ are equivalence classes $[z_n]$ of {\em Gromov-sequences} $(z_n)$ in $Z$ (see e.g. \cite{Drutu-Kapovich, MR2164775}):

A sequence $(z_n)$ in $Z$ is called a {\em Gromov-sequence} if \index{Gromov-sequence}
$$
\lim_{m, n\to\infty} (z_m. z_n)_z=\infty 
$$ 
for some (equivalently, every) $z\in Z$. 

Two Gromov-sequences $(w_m), (z_n)$ are equivalent if 
$$
\lim_{m, n\to\infty} (w_m. z_n)_z=\infty. 
$$ 
One extends the definition of the Gromov-product to the elements $\xi, \zeta$ of the Gromov-boundary 
(equivalence classes of Gromov-sequences)  
$\xi, \zeta$ by  
$$
(\xi. \zeta)_z:= \inf  \liminf_{m, n\to\infty} (w_m. z_n)_z, 
$$
where the infimum is taken over all sequences $(w_n), (z_n)$ representing $\xi, \eta$. (Here we follow the definition in \cite{CDP, Shortetal, MR2164775}, which differs (by $\le 2\delta$) from the one in \cite{bridson-haefliger} and \cite{GhH}, where the supremum is taken instead of the infimum.) 
 
Similarly, one defines 
$$
(w. \zeta)_z:= \inf  \liminf_{n\to\infty} (w. z_n)_z, 
$$
for $w\in Z, \zeta\in \geo Z$. One topologizes the space $\bar{Z}:= Z\cup \geo Z$ so that a neighborhood basis of $z\in \geo Z$ in $\bar{Z}$ 
is given by the subsets (with fixed $p\in Z$) 
$$
U_{z, \eps}= \{w \in Z\cup \geo Z: (z.w)_p<\eps\}. 
$$
In particular, a sequence $(z_n)$ in $Z$ converges to $\zeta\in \geo Z$ if and only if $(z_n)$ is a Gromov-sequence representing $\zeta$. 

\begin{defn}\index{Gromov-boundary} \label{not:rel-bdry} 
Suppose that $Z$ is a Gromov-hyperbolic metric space. We will use the notation $\geo Z$ for the Gromov-boundary of $Z$, equipped with the above topology. For a subset $Y\subset Z$, we will use the notation
$$
\geo(Y,Z) 
$$
for the accumulation set of $Y$ in $\geo Z$, the {\em relative ideal boundary} of $Y$ in $Z$. \index{relative ideal boundary}
 
We also define $\geo^{(2)}Z\subset (\geo Z \times \geo Z)/\tau$, the set of unordered pairs of {\em distinct} elements of $\geo Z$: The involution $\tau$ swaps the two factors of $\geo Z \times \geo Z$. 
\end{defn}

\medskip 
Another common definition of the {\em visual boundary} of a hyperbolic geodesic metric space $Z$ uses equivalence classes of geodesic rays in $Z$: Two rays are equivalent if they are at finite Hausdorff distance from each other. The two definitions agree if $Z$ is a proper metric space, see e.g. \cite[III.H.3]{bridson-haefliger}. For non-proper spaces one can also use quasigeodesic rays, see e.g. \cite[Section 5]{Bonk-Schramm}, \cite{MR2164775}. 

\medskip 
For our purpose, it suffices to observe that if $\gamma: \R_+\to Z$ is a quasigeodesic ray in a geodesic hyperbolic space $Z$, for each sequence $(t_n)$ in $\R_+$ diverging to $\infty$, the sequence $(\gamma(t_n))$ is a Gromov-sequence in $Z$, and any two such sequences are equivalent. Thus, $\gamma$ defines a point $\xi\in \geo Z$. We will say that 
$\gamma$ is {\em asymptotic} of $\xi$ and that $\gamma$ {\em joins} $p=\gamma(0)$ and $\xi$. We will use the notation $\gamma=p\xi$ if $\gamma$ is a geodesic ray joining $p$ and $\xi$. If $\gamma$ is a biinfinite quasigeodesic, then it defines two 
quasigeodesic rays $\gamma_\pm$ (the restrictions of $\gamma$ to $\R_+$ and to $\R_-$) and these are asymptotic to points  $\xi_\pm$, also denoted $\gamma(\pm\infty)$. A hyperbolic space $Z$ is said to be a {\em visibility space} if any two distinct ideal boundary points are connected by a biinfinite geodesic. For instance, each proper geodesic hyperbolic space is a visibility space. Even if $Z$ is a non-proper geodesic metric space, each 
point in $X$ can be joined to each point in $\geo Z$ by a quasigeodesic ray and any two distinct points in 
$\geo Z$ are connected by a biinfinite quasigeodesic, see e.g. \cite[Section 5]{Bonk-Schramm}. 

A {\em generalized geodesic triangle} in a $\delta$-hyperbolic geodesic metric space $Z$ is defined by taking a triple of points $z_1, z_2, z_3\in \bar{Z}$ (such that no two ideal boundary points in this triple are equal) and connecting them by geodesics in $Z$, the {\em sides} of the triangle $\Delta z_1 z_2 z_3$. An ideal triangle is a generalized triangle with all three vertices in $\geo Z$. 

The next lemma is an 
application of the slim triangle property in $Z$, cf. \cite[section 11.11]{Drutu-Kapovich}: 

\begin{lemma}
[Slim generalized triangle property] Suppose that $Z$ is a $\delta$-hyperbolic geodesic metric space. Then:

1. Every generalized geodesic triangle $\Delta$ in $Z$ with two non-ideal vertices is 
$2\delta$-slim: Each side of $\Delta$ is contained in the $2\delta$-neighborhood of the union of the two other sides. 

2. Every generalized geodesic triangle $\Delta$ in $Z$ with one non-ideal vertex is $3\delta$-slim. 

3. Every ideal triangle in $Z$ is $4\delta$-slim.     
\end{lemma}
\proof 1. Let $\Delta=\Delta xy\zeta$, where $x, y\in Z, \zeta\in \geo Z$. Take diverging sequences $p_n\in x\zeta$, $q_n\in y\zeta$, where $x\zeta, y\zeta$ are the infinite sides of  $\Delta$, such that $d(p_n, q_n)\le C$, where $C$ is a constant. Since the quadrilateral $xp_nq_ny$ is $2\delta$-slim (Lemma \ref{lem:slim-quad}), it follows that for each point $z\in y\zeta$, if $n$ is sufficiently large, then $z\in yq_n$ and, furthermore, $d(z, xp_n\cup xy)\le 2\delta$. (The point $z$ can be $2\delta$-close to the side $p_nq_n$ only for finitely many values of $n$.) In particular, $z$ lies in the $2\delta$-neighborhood of $xy\cup x\zeta$. The same argument proves that each point $z\in xy$ also lies in the $2\delta$-neighborhood of $xy\cup x\zeta$. 

2. Suppose that $\Delta = \Delta x \eta \zeta$, where $\eta, \zeta$ are in $\geo Z$ and $x\in Z$. As before,  consider diverging sequences $p_n\in x\eta$, $q_n\in x\zeta$. Find points $p'_n\in \eta\zeta$, $q'_n\in \eta\zeta$ within distance $C$ from $p_n, q_n$ respectively. Since the pentagon 
with the vertices $x, p_n, p'_n, q'_n q_n$ is $3\delta$-slim, it follows that each point $z\in x\zeta$ lies in the $3\delta$-neighborhood of the union $\eta\zeta\cup x\eta$.  
 
3. The proof for ideal triangles is similar (we use $4\delta$-slimness of hexagons in $X$) 
and is left to the reader. \qed

\begin{lemma}
Suppose that $Z$ is a $\delta$-hyperbolic geodesic metric space, $\gamma_i: \R_+\to Z, i=1,2$, are geodesic rays within finite Hausdorff distance from each other. Then there exist $T_1, T_2$ (depending on $\gamma_1, \gamma_2$) such that 
$$
\Hd(\gamma_1([T_1,\infty)), \gamma_2([T_2,\infty)))\le 4\delta 
$$
and, moreover, $d(\ga_1(T_1), \ga_2(T_2))\le 2\delta$. 
\end{lemma}
\proof Since the rays are at finite Hausdorff distance, $\ga_1(\infty)=\ga_2(\infty)=\xi$ for some $\xi\in \geo Z$. 
Set $p:= \gamma(0)$, $D:=d(p_1,p_2)$.  
Consider the generalized   geodesic triangle $\Delta$ in $Z$ with the vertices $p_1= \gamma_1(0), p_2=\gamma_2(0)$ and the third vertex at infinity, $\xi$. Then by the slim triangle property, every point $x_i\in \gamma_i(\R_+)$ is within distance 
$2\delta$ from $\gamma_{3-i}(\R_+)\cup p_1p_2, i=1, 2$.

Take $T_1:= D+2\delta$. The triangle inequality implies that $x_1:= \ga_1(T_1)$ is within distance $2\delta$ from some 
point $x_2=\ga_2(T_2)$.  Since the generalized triangle with the vertices $x_1, x_2, \xi$ is $2\delta$-slim, it follows that 
every point of $\gamma_1([T_1,\infty))$ is within distance $4\delta$ from a point of $\gamma_2([T_2,\infty))$ and vice versa. 
 \qed 

\begin{cor}\label{cor:asymptotic-quasirays} 
Suppose that $c_1, c_2$ are $k$-quasigeodesic rays in $Z$ such that $c_1(\infty)=c_2(\infty)$. 
Then there exist $T_1, T_2$  such that 
$$
\Hd(\gamma_1([T_1,\infty)), \gamma_2([T_2,\infty)))\le 4\delta+ 2D_{\ref{Morse}}(\delta,k). 
$$
Moreover,  $d(\ga_1(T_1), \ga_2(T_2))\le 2\delta+ 2D_{\ref{Morse}}(\delta,k)$. 
\end{cor}

A similar result holds for biinfinite quasigeodesics: 

\begin{lem}\label{lem:gen-bigons}
Suppose that $X$ is a $\delta$-hyperbolic space, $\al, \beta$ are biinfinite $L$-quasi\-geo\-de\-sics in $X$ such that
$$
\xi_\pm= \al(\pm\infty)=\beta(\pm\infty)\in \geo X
$$
Then $\Hd(\al,\beta)\le D_{\ref{lem:gen-bigons}}(L,\delta)$. 
\end{lem}
\proof Using the above corollary, we find subrays $\al_\pm$ in $\al$ and $\beta_\pm$ in $\beta$ which are asymptotic to, respectively, $\xi_\pm$, such that the respective initial points $x_\pm$ (of $\al_\pm$) and $y_\pm$ (of $\beta_\pm$) satisfy
$$
d(x_\pm, y_\pm)\le r:=2\delta+2D_{\ref{Morse}}(\delta,k). 
$$
Removing  the above rays from the quasigeodesics $\al, \beta$, we are left with two finite quasigeodesic subsegments 
$\al_0\subset \al$ (between the points $x_\pm$) and  $\beta_0\subset \beta$ (between the points $y_\pm$). 
Now, applying Lemma \ref{lem:sub-close} to $\al_0, \beta_0$, we get:
$$
\Hd(\al_0, \beta_0)\le D_{\ref{lem:sub-close}}(\delta, L, r). 
$$ 
Hence, taking $D_{\ref{lem:gen-bigons}}(L,\delta): \max(r, D_{\ref{lem:sub-close}}(\delta, L, r))$, concludes the proof of the lemma. \qed

Each qi embedding $f: X\to Y$ of geodesic hyperbolic spaces has an extension to a map of Gromov-boundaries, a topological embedding $\geo f: \geo X\to \geo Y$. The combined map
$$
f\cup \geo f: X\cup \geo X\to Y\cup \geo Y
$$
is continuous at each point of $\geo X$, see e.g. \cite[Exercise 11.109]{Drutu-Kapovich}. In Chapter \ref{ch:CT} we will discuss the existence of such an extension in the case of more general coarse Lipschitz maps of hyperbolic spaces.

\section{Quasiconvex subsets} 

\begin{defn}
[Quasiconvex subset] Let $X$ be a geodesic metric space, $Y\subseteq X$ and let $\lambda\geq 0$.
We say that $Y$ is $\la$-quasiconvex in $X$ if every geodesic with end-points
in $Y$ is contained in the $\la$-neighborhood of $Y$. A subset $Y\subset X$
is said to be quasiconvex if it is $\la$-quasiconvex for some $\la\geq 0$. \index{quasiconvex subset} 
\end{defn}

This definition generalizes in the setting of path-metric spaces where $Y\subset X$ is $\la$-quasiconvex if there is a function $\eta: \R_+\to \R_+$ 
converging to $0$ as $\eps\to 0+$, such that  each $\eps$-geodesic $\ga\subset X$ with end-points in $Y$, satisfies
$$
\ga\subset N_{\al+ \eta(\eps)}(Y). 
$$

\medskip 
\noindent Convex subsets of, say, Hadamard spaces, have several characteristic properties such as:

\begin{itemize}
\item Intersections of convex subsets are again convex. 

\item The nearest-point projection to a nonempty closed convex subset is well-defined and 1-Lipschitz. 
\end{itemize}

Below we discuss analogues of these properties for quasiconvex subsets of geodesic 
hyperbolic spaces.  

\medskip 
The following two lemmata are straightforward and we omit the proofs:

\begin{lemma}
Every geodesic triangle in a $\delta$-hyperbolic geodesic metric space is $\delta$-quasiconvex. 
\end{lemma}

\medskip 

\begin{lemma}\label{lem:nbds-of-qc}
1. Suppose that $A$ is a $\la$-quasiconvex subset of a $\delta$-hyperbolic space $X$. 
Then the $R$-neighborhood $N_R(A)$ is $(R+2\delta+\la)$-quasiconvex in $X$.

2. Suppose that $A$ is a $\la$-quasiconvex subset of a $\delta$-hyperbolic space $X$ and $B\subset X$ is such that $\Hd(A,B)\le R$. Then $B$ is $2R+2\delta+\la$-quasiconvex in $X$. 
\end{lemma}

\medskip 
We next discuss the relation between the quasiconvexity and qi embeddings.

\begin{lemma}\label{lem:qc-qp} 
If $A$ is a $\la$-quasiconvex subset of a  geodesic hyperbolic metric space $(X,d)$, then:

1. The metric space 
$(N_\la(A), d)$ is a quasi-path metric space.

2.  When we equip $N_\la(A)$ with the path-metric $d_p$ induced from $X$, the inclusion map  $N_\la(A)\to X$ is a qi embedding.   
\end{lemma}
\proof If $\la=0$, then $(A,d)$ is actually a geodesic metric space isometrically embedded in $(X,d)$ 
and there is nothing to prove. Thus, we assume that $\la>0$. 

1. Consider points $x, y\in N_\la(A)$ and $a\in A, b\in A$ be points within distance $\la$ from $x, y$ respectively. Let $\ga=ab$ be a geodesic  in $X$ contained in $N_\la(A)$. Then divide $\ga$ into subsegments $x_i x_{i+1}, 1\le i\le n-1$ of equal length $\le \la$. Then the finite sequence
$$
x=x_0, a=x_1, x_2,..., x_n=b, x_{n+1}=y
$$ 
is a $\la$-path connecting $x$ to $y$ and 
$$
\sum_{i=0}^n d(x_i, x_{i+1})\le d(x,y) + 4\la. 
$$

2. By the same argument, for any two points $x, y\in N_\la(A)$,
$$
d_p(x,y)\le d(a,b) + 2\la \le d(x,y) + 4\la. 
$$
Therefore, the inclusion map $(N_\la(A), d_p)\to (X,d)$ is a $(1, 4\la)$-qi embedding. \qed

\begin{lemma}\label{lem:qi-preserves2}
Suppose that $Y, X$ are geodesic metric spaces, $X$ is  $\delta$-hyperbolic and $f: Y\to X$ is a $K$-qi embedding. 
Then $f(Y)$ is $\la_{\ref{lem:qi-preserves2}}(\delta, K)=D_{\ref{stab-qg}}(\delta,K)$-quasiconvex in $X$. 
\end{lemma}
\proof For a geodesic $\ga=y_1y_2\subset Y$,  $f(\ga)$ is a $K$-quasigeodesic in $X$ connecting $x_i=f(y_i), i=1,2$. Therefore, by Lemma \ref{stab-qg}, 
$$
\Hd(f(\ga), x_1 x_2)\le \la=D_{\ref{stab-qg}}(\delta,K)=92K^2(K + 3\delta),
$$
which implies that $x_1x_2\subset N_\la(f(Y))$, i.e. $f(Y)$ is $\la$-quasiconvex.  \qed 

\medskip 
The following is a converse to Lemma \ref{lem:qi-preserves2}:

\begin{lemma}\label{lem:up+qc->qi} 
Suppose that $(X,d)$ is hyperbolic, $Y$ is a geodesic metric space and $f: Y\to X$ is a uniformly proper map with quasiconvex image. Then $f$ is a qi embedding.  
\end{lemma}
\proof By Lemma \ref{lem:qc-qp}, if $A:= f(Y)$ is $\la$-quasiconvex, then $(N_\la(A), d_p)$ is qi embedded in $X$. 
Since $f$ is uniformly proper, Lemma \ref{lem:qi-net} then implies that $f: Y\to (N_\la(A), d_p)$ is a quasiisometry. Therefore, $f: Y\to X$ is a qi embedding as a composition of two qi embeddings. \qed 


 \section{Quasiconvex hulls}

A common source of quasiconvex subsets in hyperbolic spaces is given by {\em quasiconvex hulls} which we discuss in this section. 

\begin{defn}\label{defn:hull}
For a subset $U$ of a geodesic metric space $X$ we define the {\em quasiconvex hull} \index{quasiconvex hull} 
$\hull(U)$ as the union of all geodesics in $X$ connecting all pairs of points in $U$. 
The $\eps$-{\em quasiconvex hull} $\hull_{\eps}(U)$ of $U$ in $X$ is defined as the closed 
$\eps$-neighborhood of the quasiconvex hull $\hull(U)$. Thus, $cl(\hull(U))= \hull_0(U)$. 
\end{defn}

By the construction, if $U$ is $\la$-quasiconvex, then 
\begin{equation}\label{eq:qc-nbd} 
\hull_{\eps}(U)\subset N_{\eps+\la}(U). 
\end{equation}

\begin{lemma}\label{lem:delta-hull}
If $X$ is a $\delta$-hyperbolic geodesic metric space, then the inclusion map $Y=\hull_{\delta}(U)\to X$ is a $(1,\epsilon)$-quasiisometric embedding with 
$\epsilon=6\delta$ and $Y$ equipped with the path-metric induced from $X$. 
\end{lemma}
\proof Suppose first that $y_1, y_2\in \hull(U)$ belong to geodesics $a_1b_1, a_2b_2$ with $a_i, b_i\in U, i=1,2$. Since $X$ is $\delta$-hyperbolic, by considering geodesic triangles 
$\Delta a_1 b_1 b_2$, $\Delta a_1 a_2b_2$, we see that $y_1$ is within distance $\delta$ from 
a point $z_1 \in a_1b_1 \cup b_1b_2$. Similarly, $y_2$ is within distance $\delta$ from a point $z_2 \in a_1b_1 \cup a_1a_2$. If both $z_1, z_2$ belong to the same geodesic $a_1b_2$ then we connect $y_1, y_2$ by a path of length $d(z_1, z_2) + 2\delta$ in $Y$. Hence, $d_Y(y_1, y_2)\le d_X(y_1,y_2)+4\delta$. The cases when one of the points $z_1, z_2$ is on the geodesic $a_1b_2$ are similar and left to the reader. 

Suppose that $z_1\in b_1b_2, z_2\in a_1a_2$. We then switch to the triangles $\Delta a_1 a_2 b_1$, $\Delta a_1 b_1b_2$. Again, it suffices to consider the case when $y_1$ is within distance $\delta$ from $w_1\in a_1a_2$ or $y_2$ is within distance $\delta$ from $w_2\in b_1b_2$. But then we can connect $y_1, y_2$ by a path in $Y$ of length 
$\le 2\delta+ d(w_1, z_2)$ or $\le 2\delta+ d(z_1, w_2)$. Again,  $d_Y(y_1, y_2)\le d_X(y_1,y_2)+4\delta$. 

Lastly, any two points $x_1, x_2\in Y$ are within distance $\delta$ from $\hull(U)$ and, therefore, for $\epsilon=6\delta$,
$$
d_Y(x_1, x_2)\le d_X(x_1,x_2)+\epsilon.   \qed 
$$

\begin{cor}\label{cor:delta-hull} Suppose $X$ is a $\delta$-hyperbolic geodesic metric space, 
 $U\subset X$ is a $\la$-quasiconvex subset and $R \ge \la+\delta$. Then the inclusion map $N_R(U)\to X$ is a   
$(1,6\delta+R)$-quasiisometric embedding, where $N_R(U)$ equipped with the path-metric induced from $X$. 
\end{cor}

\begin{lemma}\label{lem:QCunion}
Let $X$ be $\delta$-hyperbolic geodesic metric space; assume that $U, V\subset X$ are $\la$-quasiconvex subsets with nonempty intersection. Then:

1. For every $R\ge \la+\delta$ the subset 
$$
Y= N_R(U\cup V)\subset X
$$
is $(R+\delta)$-quasiconvex. 

2. Furthermore, for each $R\ge 2\la+4\delta$,  the inclusion map $Y\to X$ is a $(1,6\delta+R)$-quasiisometric embedding where 
 we equip $Y$ with the path-metric induced from  $X$. 
\end{lemma}
\proof 1. For $a\in U, b\in V, c\in U\cap V$ we have
$$
ab\subset N_{\delta}(ac\cup bc) \subset N_{\la+\delta}(U\cup V)\subset Y=N_R(U\cup V).
$$
For $x\in N_R(U), y\in N_R(V)$ there exist $a\in U, b\in V$ satisfying $xy\subset N_{R+\delta}(ab)$. Since $ab\subset Y$, we obtain 
  $xy\subset N_{R+\delta}(Y)$, i.e. $Y$ is $(R+\delta)$-quasiconvex. 

2. Observe that
$$
N_R(U\cup V)= N_{R-(\la+\delta)}(N_{\la+\delta}(U\cup V)). 
$$
Thus, the second part of the lemma follows from the first in combination with Corollary \ref{cor:delta-hull}. \qed

\begin{lemma}\label{lem:qc-hull}
Suppose that $X$ is a $\delta$-hyperbolic geodesic metric space. Then for every $Y\subset X$, $\hull_\eps(Y)$ is $3\delta+\eps$-quasiconvex in $X$. \end{lemma}
\proof Take $x_1, x_2\in \hull_\eps(Y)$; these points are within distance $\eps$ from $y_1, y_2\in \hull(Y)$. 
As in the proof of Lemma \ref{lem:delta-hull}, there exist $z_1, z_2$ lying on a common geodesic $a_1a_2$ connecting points in $Y$, such that 
$$
d(y_i, z_i)\le 2\delta, i=1, 2. 
$$
Hence, $d(x_i, z_i)\le \eps+ 2\delta$, which implies that for $\la=3\delta+\eps$, 
$$
x_1 x_2\subset N_\la(z_1 z_2) \subset N_\la(\hull(Y)). \qedhere  
$$

\section{Projections} 

We now discuss properties of quasiconvex subsets in hyperbolic metric spaces. 
For a complete metric space $(X,d)$, a closed subset $Y\subset X$ and a point $x\in X$, we will be using {\em nearest point projections} $\bar{x}\in Y$ of the point $x$ to the subset $Y$. We note that a nearest point projection always exists (but need not be unique). 

\begin{rem}
The assumptions that $X$ is complete and $Y$ is closed are essentially irrelevant: If they fail, then instead of nearest-point projections $\bar{x}\in Y$ of $x\in X$, 
we can use $\eps$-nearest point projections, i.e. points $\bar{x}\in Y$ such that  for all $y\in Y$, 
$$
d(x, y)\ge d(x, \bar{x}) -\eps. 
$$
Since we are interested in coarse-geometric aspects of metric spaces, such ``almost'' nearest point projections will work just as well. In what follows, we, therefore, will use the same fudge as for path-metric spaces: We will talk about nearest point projections without imposing the assumption that a subset is closed and the ambient space is 
complete. 
\end{rem}

{\bf Notation for nearest point projections:} Suppose that $Y$ is a closed subset in a complete metric space $X$. 
We shall denote by $P_{X,Y}$ or $P_Y$ (when the choice of $X$ is clear) or sometimes simply by $P$ 
(when the choices of $X$ and $Y$ are clear) a nearest-point projection 
 $X\to Y$.   We will see below (Lemma \ref{lip-proj}) that in hyperbolic spaces, 
 nearest-point projections to quasiconvex subsets 
 are ``coarsely well-defined,'' thus, justifying the notation.

 \begin{lemma}\label{lem:citerion-of-projection} 
Suppose that $X$ is a geodesic metric space, $Y\subset X$, $x\in X$ and $\hat{x}\in Y$ is such that each geodesic $xy, y\in Y$, passes within distance 
$r$  from $\hat{x}$. Then $d(\bar{x}, \hat{x})\le 2r$, where $\bar{x}=P_Y(x)$. 
  \end{lemma}
 \proof We take $y=\bar{x}$ and let $z\in xy$ be a point within distance $r$ form $\hat{x}$. 
 Since $d(x, \hat{x})\ge d(x, y)$, it follows that $r\ge d(z, \hat{x})\ge d(z,y)$. Therefore, $d(\hat{x},y)\le 2r$. 
  \qed 
 
 \medskip 
We now turn to quasiconvex subsets of hyperbolic spaces.

\begin{lemma}
[See e.g. Lemma 11.53 in \cite{Drutu-Kapovich}] \label{lip-proj}
For each $\la\geq 0$ there is a constant $L_{\ref{lip-proj}}=L_{\ref{lip-proj}}(\delta,\la)$
such that the following holds: 

If $X$ is a $\delta$-hyperbolic metric space and $Y$ is a $\la$-quasiconvex subset, then the projection map 
$P=P_{X,Y}: X\rightarrow Y$ is coarsely $L_{\ref{lip-proj}}(\delta,\la)$-Lipschitz, i.e. for all $x,y\in X$ we have 
$d_X(P(x),P(y))\leq L_{\ref{lip-proj}}(\delta,\la) (d_X(x,y)+ 1)$. 
In particular, $P$ is coarsely well-defined: For different choices $\bar{x}_1, \bar{x}_2$ of points in $Y$ nearest to $x$ we have $d_X(P(x),P(y))\leq L_{\ref{lip-proj}}(\delta,\la)$. \end{lemma}

\begin{rem}\label{rem:lip-proj}
We will frequently use this lemma when $Y$ is a geodesic in $X$, in which case $\la=\delta$ 
and $d(P(x), P(y))\le d(x,y)+ 12\delta$, i.e. 
$L_{\ref{lip-proj}}(\delta,\delta)=12\delta$, compare Lemma \ref{proj-geod}. In general:
$$
L_{\ref{lip-proj}}(\delta,\la)= \max(2, 2\la+ 9\delta), 
$$
see Lemma 11.53 in \cite{Drutu-Kapovich}. 
\end{rem}

\medskip 
The next lemma is a converse to Lemma \ref{lip-proj}:  

\begin{lemma}\label{lem:qc-criterion} 
Suppose that a metric space $(X,d_X)$ is $\delta$-hyperbolic, $Y\subset X$ is a rectifiably connected subset (equipped with the induced path-metric $d_Y$) 
such that there exists a  $k$-coarse-Lipschitz retraction 
$$
P: X\to Y. 
$$
Then $Y$ is $\la=\la_{\ref{lem:qc-criterion}}(k,\delta)$-quasiconvex in $X$. 
\end{lemma}
\proof The inclusion map $i: (Y, d_Y)\to (X,d_X)$ is $1$-Lipschitz. Combining this with the existence of the  retraction $P$, we conclude that 
$i$ is a $k$-quasiisometric embedding. For each $\eps>0$ any two points $y_1, y_2\in Y$ can be connected by a path $\al$ whose length is 
$\le d_Y(y_1, y_2)+\eps$, in particular, $\al$ is $(1,\epsilon)$-quasigeodesic. We will take $\eps=1$. 
The composition $i\circ \al$ in $X$ is $L=k(1+\epsilon)$-quasigeodesic.  By the Morse Lemma (Lemma \ref{stab-qg}), the image of  $i\circ \al$ is within distance $\la=D_{\ref{stab-qg}}(\delta, L)$ from a geodesic $\ga=y_1 y_2$ in $X$ connecting $y_1$ to $y_2$. Hence,  $\ga$ is contained in the $\la$-neighborhood of $Y$. \qed

\medskip 
Lemma \ref{lem:projection-1} below is a converse to Lemma \ref{lem:citerion-of-projection} in the context of hyperbolic spaces; it will be used repeatedly in the book: 

\begin{lemma}\label{lem:projection-1}
Suppose that $X$ is $\delta$-hyperbolic geodesic metric space, $Y\subset X$ is a $\la$-quasiconvex subset of $X$. Let $\bar{x}\in Y$ be a nearest-point projection of $x\in X$ to $Y$ and let $y\in Y$ be an arbitrary point. Then:

(i) $\bar{x}$ lies within distance $\la+2\delta$ 
from (any) geodesic $xy\subset X$. 

(ii) $\Hd(xy, x\bar{x}\cup \bar{x} y)\leq \la+3\delta$.

(iii)  The concatenation $x\bar{x}\star \bar{x} y$ is a  
$2(\la+2\delta)$-quasigeodesic. 
\end{lemma}
\proof (i) We consider a geodesic triangle $\Delta x  \bar{x} y$. Since this triangle is $\delta$-slim, there exist points $z\in xy, z'\in x\bar{x}, z''\in y \bar{x}$ such that $d(z, z')\le \delta, d(z', z'')\le \delta$. Then $d(z'', Y)\le \la$. 
Since $\bar{x}$ is a nearest point to $x$ in $Y$, it follows that $d(z', \bar{x})\le \la+\delta$, hence, $d(\bar{x}, z)\le \la+ 2\delta$.  

(ii) As in the proof above, let $z\in xy$ be such that $d(x,z)\leq \la+2\delta$. Then by $\delta$-hyperbolicity of $X$, 
it is clear that $\Hd(\bar{x}y, zy)\leq \la+3\delta$ and $\Hd(x\bar{x}, xz)\leq \la+3\delta$.
From these it follows that $\Hd(xy, x\bar{x}\cup \bar{x} y)\leq \la+3\delta$. 

(iii) Taking $z\in xy$ as above, we see that
$$
d(x,y)\le d(x, \bar{x})+ d(\bar{x}, y) + 2(\la+2\delta). 
$$
Hence, $x\bar{x}\star \bar{x} y$ is a $(1,2(\la+2\delta))$-quasigeodesic. \qed 

\medskip 
Below are several corollaries of the lemma: 

\begin{cor}\label{cor:center-projection}
If $Y$ is a geodesic in $X$, then  $\bar{x}$ lies within distance $3\delta$ 
from (any) geodesic $xy\subset X$. In particular, for each geodesic triangle $\Delta=\Delta xyz$ in $X$, the projection $P_{yz}(x)$ is 
a $3\delta$-center of $\Delta$. 
\end{cor}

\begin{cor}
[Almost nearest-point projection] \label{cor:almost-npp}
Suppose that $X, Y$ are as above and $y\in Y$ satisfies
$$
d(x, y)\le d(x, \bar{x}) + C. 
$$
Then $d(\bar{x}, y)\le C + 2\la + 4\delta$. 
\end{cor}
\proof Let $z\in xy$ be a point within distance $\la+2\delta$ from $\bar{x}$ and let $D$ denote $d(y, z)$. Then 
$$
d(x, \bar{x}) + C  -D \ge  d(x,y) - D = d(x,z)\ge d(x, \bar{x}) - (\la+2\delta),
$$
which implies that $D\le C+ \la+2\delta$. Hence, $d(\bar{x}, y)\le d(\bar{x},z)+d(y,z)\le C + 2\la + 4\delta$.\qed 
 
 \medskip 

\begin{cor}\label{cor:proj-to-close-subsets} 
Let $U, V\subset X$ be a pair of closed $\la$-quasiconvex subsets such that $\Hd(U,V)\le D$. 
Then $d(P_U, P_V)\le D_{\ref{cor:proj-to-close-subsets}}(\delta,\la,D)$. 
\end{cor} 
 \proof Let $u\in U, v\in V$ denote $P_{X,U}(x)$ and $P_{X,V}(x)$ respectively. Since  $\Hd(U,V)\le D$, there exist points $u'\in V, v'\in U$ such that 
 $$
 d(u,u')\le D, \quad d(v,v')\le D. 
 $$
 Hence,
 $$
d(x, u') -D \le d(x, u) \le  d(x, v')\le d(x, v)+ D
 $$
By Corollary \ref{cor:almost-npp}, $d(v, u') \le 2D + 2\la +4\delta$. Hence, $d(u, v)\le D_{\ref{cor:proj-to-close-subsets}}(\delta,\la,D):=3 D + 2\la +4\delta$. 
  \qed

\begin{cor}\label{cor:projection-1}
Suppose that $Y\subset X$ is a $\la$-quasiconvex subset of a $\delta$-hyperbolic space $X$ and $x, \hat{x}$ are points in $X, Y$ respectively connected by a geodesic $c=x \hat{x}$ 
such that:

There exists a function $R\mapsto \hat{R}$ satisfying $y\in Y$, $d(y, c)\le R\Rightarrow d(y, \hat{x})\le \hat{R}$.  

\noindent Then $d(\hat{x}, P_Y(x))\le \widehat{\la +2\delta}$. 
\end{cor} 
 \proof Let $\bar{x}= P_Y(x)$. Then, by Lemma \ref{lem:projection-1}(1), $d(\bar{x}, c)\le \la+2\delta=:R$. It then follows that 
 $d(\bar{x}, \hat{x})\le \hat{R}=\widehat{\la +2\delta}$. \qed 
 
 \begin{cor}\label{cor:projection-2}
Suppose that $V\subset U\subset X$ are two $\la$-quasiconvex subsets of $X$. Then 
$$
d(P_{X,V}, P_{U,V}\circ P_{X,U})\le C_{\ref{cor:projection-2}}(\delta,\la). 
$$
 \end{cor} 
\proof The proof is by repeated use of Lemma \ref{lem:projection-1}. For $x\in X$ set $x_1:= P_U(x), x_2:= P_V(x), x_3:= P_V(x_1)$.  
Consider the triangle $\Delta x x_1x_2$. By Lemma \ref{lem:projection-1}(i), there is a point $x'_1\in xx_2$ such that $d(x_1, x'_1)\leq \la +2\delta$
because $x_2\in U$ and $x_1$ is a nearest point projection of $x$ to $U$. Now we note that
$x_2$ is a nearest point projection of $x'_1$ on $V$. Hence there is a point $x'_2\in x'_1x_3$ such that
$d(x_2, x'_2)\leq \la+2\delta$. Since $d(x_1, x'_1)\leq \la+2\delta$ and the $\Delta x_1x'_1x_3$
is $\delta$-slim, there is a point $x'_3\in x_1x_3$ such that $d(x'_2, x'_3)\leq \la+3\delta$. Hence, 
$$
d(x_2, x'_3)\leq d(x_2,x'_2)+d(x'_2, x'_3)\leq 2\la+5\delta.$$ 
Finally, we note that $x_3$ is a nearest point
projection of $x'_2$ on $V$ too. Hence, $d(x'_3,x_3)\leq d(x'_3,x_2)\leq 2\la+5\delta$.
Thus 
$$
d(x_2,x_3)\leq C_{\ref{cor:projection-2}}(\delta,\la):=d(x_2,x'_3)+d(x'_3,x_3)\leq 4\la+10\delta. \qed 
$$

\medskip The following  variation on this corollary will be used in the proof of Corollary \ref{cor:proj-to-qc-subgroup} and 
Theorem \ref{thm:mitras-projection}: 

\begin{lemma}\label{lem:two-projections} 
Suppose that  $X$ is $\delta$-hyperbolic, $U, V\subset X$ are, respectively, $\la$ and $\mu$-quasiconvex subsets. Set $\eps:= \mu+2\delta$. Then for every $r\ge r_{\ref{lem:two-projections}}(\la,\mu,\delta)=\eps+\delta+\la$, if $W=N_r(U)\cap V$ is nonempty, then the distance between the restrictions to $V$ of the projections 
$P_{X,W}, P_{X,U}$,  
is at most $\eps+\delta+r$.  
\end{lemma} 
 \proof Take $x\in V$ and consider points $\bar{x}=P_{X, W}(x), \hat{x}:= P_{X,U}(x)$. Since $\bar{x}\in N_r(U)$, there exists a point $x'\in U$ such that 
 $d(x', \bar{x})\le r$. By Lemma \ref{lem:projection-1}, 
 there exists a point $y\in x x' \cap B(\hat{x}, \eps)$. Since the triangle  $\triangle x x' \bar{x}$ is $\delta$-slim, the 
 point $y$ is within distance $\delta$ either from a point $z\in \bar{x} x'$ or from a point $z\in \bar{x} x$. In the former case, 
 $$
 d(\bar{x}, \hat{x})\le \eps+ \delta+r,
 $$
 as claimed. Suppose, therefore, that $z\in \bar{x} x$.  Since $V$ is $\la$-quasiconvex, there exists $w\in V$ within distance $\la$ from $z$. Thus,
 $$
 d(\hat{x}, w)\le \eps+\delta+\la. 
 $$
 The assumption that $r\ge \eps+\delta+\la$ implies that $w\in N_r(U)\cap V= W$. In particular,
 $$
 d(x, w)\ge d(x, \bar{x}). 
 $$ 
 Setting $t=d(x, z), s=d(z, \bar{x})$, we obtain $d(x, \bar{x})=t+s$ and, hence,
 $$
t+s=d(x, \bar{x})\le d(x, w)\le t+\la.
 $$ 
 It follows that $s\le \la$ and, therefore, $d(\hat{x}, \bar{x})\le \eps+\delta+s\le \eps+\delta+\la \le r$. Thus, in both cases, 
 $d(\bar{x}, \hat{x})\le \eps+ \delta+r$.  \qed 
  
 \medskip 
Below is another variation on Corollary \ref{cor:projection-2} and Lemma \ref{lem:projection-1}(iii): 

\begin{lem}\label{lem:two-projections-2}
Let $X$ be a $\delta$-hyperbolic space, and let $Y_1\supset Y_2 \supset ... \supset Y_n$ be a chain of $\la$-quasiconvex subsets 
of $X$. For a point $x\in X$ define inductively points $y_1= P_{X,Y_1}(x),..., y_n=P_{Y_{n-1},Y_n}(y_{n-1})$. Then:

1. The concatenation 
$$
xy_1 \star .... \star y_{n-1} y_n
$$
is an $L_{\ref{lem:two-projections-2}}(\delta,\la,n)$-quasigeodesic in $X$. 

2. $d_X(P_{X,Y_n}(x), x_n)\le D_{\ref{lem:two-projections-2}}(\delta,\la,n)$.  
\end{lem} 
 \proof The proof is by induction on $n$. If $n=1$, there is nothing to prove. Suppose that the claim holds for $n-1$, we will prove it for $n$.  
 
 By Lemma \ref{lem:projection-1}(iii), the concatenation $xy_1\star y_1 y_n$ is a  $2(\la+2\delta)$-quasigeodesic in $X$. By the induction hypothesis, the concatenation $y_1 y_2\star ... \star y_{n-1} y_n$ is a $L(\delta,\la,n-1)$-quasigeodesic. Now, Part 1 follows from the stability of quasigeodesics in hyperbolic spaces.  
 
 Part 2 of the lemma follows from the inductive application of Corollary \ref{cor:projection-2}.  \qed 

\medskip 
We conclude the section with two technical lemmata that 
will be used in Section \ref{sec:def-of-flow-spaces}. 
The lemmata generalize the obvious fact that if $Z\subset Y\subset X$ are inclusions of subsets in a 
metric space $X$ such that $Z$ is convex in $X$, then $Z$ is also convex in $Y$, while if 
$Z$ is convex in $Y$ and $Y$ is convex in $X$ then $Z$ is convex in $X$. 

\begin{lemma}\label{lem:triple1}
Suppose that $Z\subset Y\subset X$ are inclusions of metric spaces such that $Y$ is a $\delta$-hyperbolic geodesic metric space, $X$ is a geodesic metric space, $Z\subset X$ is $\la$-quasiconvex, and 
the inclusion map $Y\to X$ is an $L$-quasiisometric embedding, 
where $\la\ge \frac{3}{2}\delta$. Then $Z\subset Y$ is 
$\la'$-quasiconvex with $\la'=1500 (L\la)^3$.
\end{lemma}
\proof Pick any pair of points $p, q\in Z$, let $pq$ be a geodesic in $X$ contained in $N_\la(Z)$, 
let $\ga: I=[a,b]\to X$ be its arc-length parameterization.
Pick a maximal finite sequence of 1-separated points $x_i=\ga(t_i)$ and project these points to a sequence 
$z_i= P_Z(x_i)\in Z$. Then the maps $[t_i, t_{i+1})\to \{z_i\}$ define a $k=2L\la$-quasigeodesic 
$\beta: I\to Z\subset Y$.  Let $\beta^*$ denote a geodesic in $Y$ connecting $p$ to $q$. By Theorem \ref{stab-qg}, 
$$
\Hd(\beta, \beta^*)\le 92k^2(k + 3\delta)\le 92\cdot 16 (L\la)^3 \le 1500  (L\la)^3. \qedhere 
$$

\begin{rem}
One can give a faster proof of the lemma using the restriction to $Y$ 
of the coarse Lipschitz retraction $P_{X,Z}$ (and quoting Lemma \ref{lem:qc-criterion}), but 
we prefer to get an explicit estimate.  
\end{rem}

\begin{lemma}\label{lem:triple2}
Suppose that $Z\subset Y\subset X$ are inclusions of metric spaces such that $X, Y$ are geodesic metric spaces (with the path-metric on $Y$ is induced by that of $X$), $X, Y$ are $\delta_1$ and 
$\delta_2$-hyperbolic respectively, $Y$ is $\la_1$-quasiconvex in $X$, $Z\subset Y$ is $\la_2$-quasiconvex, then $Z$ is $\la_{\ref{lem:triple2}}(\la_1,\la_2,\delta_1,\delta_2)$-quasiconvex in $X$. 
\end{lemma}
\proof The projections $P_1=P_{X,Y}: X\to Y, P_2=P_{Y,Z}: Y\to Z$ are $L_{\ref{lip-proj}}(\delta_i,\la_i)$-coarse Lipschitz retractions, $i=1,2$. Hence, their composition $P_2\circ P_1: X\to Z$ is an $k$-coarse Lipschitz retraction with  
$$
k= L_{\ref{lip-proj}}(\delta_1,\la_1)\cdot L_{\ref{lip-proj}}(\delta_2,\la_2). 
$$
Therefore, by Lemma \ref{lem:qc-criterion}, $Z\subset X$ is  $\la$-quasiconvex for 
$$
\la= \la_{\ref{lem:qc-criterion}}(k,\delta_1). \qedhere 
$$

\section {Images and preimages of quasiconvex subsets under projections}

 
 In this section we discuss the extent to which images and preimages 
 of projections of quasiconvex subsets to quasiconvex subsets are again quasiconvex. 
 

\begin{lemma}\label{lem:small-projection}
Let $\al=xy\subset X$ be a geodesic in a $\delta$-hyperbolic geodesic metric space, $Y\subset X$ a $\la$-quasiconvex subset and $\bar{x}=P_Y(x), \bar{y}=P_Y(y)$.
Then $\bar\al=P_Y(\al)$ is $D_{\ref{lem:small-projection}}(\delta,\la)$-Hausdorff close to the geodesic $\bar x \bar y$. 
\end{lemma}
\proof i. Take $z\in \al$. Then $z$ lies within distance $2\delta$ from a point $w\in x \bar x \cup \bar x \bar y \cup \bar y y$. If $w\in \bar x \bar y$ then $d(z, Y)=d(z, \bar{z})$ ($\bar z=P_Y(z)$) satisfies $d(z, Y)\le \la+2\delta$. In particular, $d(\bar z, w)\le 4\delta +\la$. Suppose that $w\in x \bar x$. Then, $d(w, Y)= d(w, \bar{x})$ and, without loss of generality, $P_Y(w)= \bar{x}$. Lemma  \ref{lip-proj} implies that $d(P(z), \bar{x})= d(P(z), P(w))\le   (2\delta+ 1) L_{\ref{lip-proj}}(\delta,\la)$. The case $w\in y \bar y$ is handled by relabelling $x$ and $y$. To conclude:
$$
\bar\al\subset N_{D}(\bar x \bar y), \quad  D= \max( (2\delta+ 1) L_{\ref{lip-proj}}(\delta,\la), 4\delta +\la).  
$$ 

ii. Consider a point $z\in \bar x \bar y$. 
The point $z$ is within distance $2\delta$ from some 
$w\in xy \cup x \bar{x}\cup y \bar{y}$. If $w\in xy$ then $d(w, P_Y(w))\le 2\delta +\la$ and, hence, 
$d(z, P(w))\le 4\delta +\la$. Suppose that $w\in x \bar{x}$. Then $d(z, \bar{x})\le 2\delta$.
Similarly in case $w\in y \bar{y}$ we have $d(z, \bar{y})\le 2\delta$.
Hence, in either case, we have $d(z, \bar\al)\le  2\delta +\la\le D$.

To conclude, for 
$$
D=D_{\ref{lem:small-projection}}(\delta,\la):=  \max( (2\delta+ 1) L_{\ref{lip-proj}}(\delta,\la), 4\delta +\la),
$$
 the geodesic $\bar x \bar y$ and the set $\bar\al$ are $D$-Hausdorff close.   \qed 

\begin{cor}\label{cor:small-projection}
If $\diam\{\bar x, \bar y\}\le D$ then $\diam(\bar \al)\le 2D_{\ref{lem:small-projection}}(\delta,\la)+D$. 
\end{cor}

\begin{rem}
If $\delta\ge 1$, then 
$$
D_{\ref{lem:small-projection}}(\delta,\la)\le 3\delta(2\la+9\delta).$$ 
\end{rem}

\begin{cor}\label{cor:npp-of-qgeodesic}
Let $X$ be a $\delta$-hyperbolic space, $Y\subset X$ a $\la$-quasiconvex subset and 
$\ga\subset X$ be a $K$-quasigeodesic connecting points $x$ and $y$ and $P=P_{X,Y}$. 
Then $P(\ga)$ is $C_{\ref{cor:npp-of-qgeodesic}}(\delta,\la,K)$-Hausdorff close to 
the geodesic segment $P(x)P(y)$. 
\end{cor}
\proof This corollary follows from the Morse Lemma, Lemmata \ref{lip-proj} and   \ref{lem:qc-projection}.

 \medskip 
As another application of the lemma, we obtain:

\begin{lemma}\label{lem:qc-projection}
Suppose that $Y, Z$ are $\la_Y, \la_Z$-quasiconvex subsets respectively in a $\delta$-hyperbolic space $X$. Then the projection 
$P_Y(Z)$ is $D_{\ref{lem:qc-projection}}(\delta,\la_Y, \la_Z)$-quasiconvex in $X$. 
\end{lemma}
\proof Take points $y_i=P_Y(z_i), z_i\in Z, y_i\in Y, i=1,2$. By Lemma \ref{lem:small-projection},
$$
\Hd(y_1 y_2, P_Y(z_1z_2))\le D_{\ref{lem:small-projection}}(\delta,\la_Y). 
$$
Thus, for every $y\in y_1y_2$ there exists $x\in z_1 z_2$ such that 
$d(y, P_Y(x))\le D_{\ref{lem:small-projection}}(\delta,\la_Y)$. Since $Z$ is $\la_Z$-quasiconvex in $X$, there exists $z\in Z$ such that $d(x,z)\le \la_Z$. By Lemma \ref{lip-proj}, 
$$
 d(P_Y(x), P_Y(z))\le (\la_Z+1) L_{\ref{lip-proj}}(\delta,\la_Y).  
 $$
 Putting together these inequalities, we obtain:
 \begin{align*}
 d(y, P_Y(Z))\le d(y, P_Y(z))\le
 D_{\ref{lem:qc-projection}}(\delta, \la_Y, \la_Z):=  \\
 (\la_Z+1) L_{\ref{lip-proj}}(\delta,\la_Y) + D_{\ref{lem:small-projection}}(\delta,\la_Y). \qedhere
 \end{align*}
 
 \begin{rem}\label{rem:qc-projection}
 If $\delta\ge 1, \la\ge 1$ then
 $$
 D_{\ref{lem:qc-projection}}(\delta,\la,\la)\le    (2\la +3\delta)(2\la + 9\delta). 
 $$
 \end{rem} 

 Thus, projections (to uniformly quasiconvex subsets) in hyperbolic spaces send uniformly quasiconvex subsets to uniformly quasiconvex subsets. The next lemma, which we add only for the completeness of the picture and which is not used elsewhere otherwise, establishes a similar statement for preimages. 
 We need to warn the reader that preimages of quasiconvex subsets under projections need not be quasiconvex, the true statement is more subtle. 
 
 
\begin{lemma}\label{lem:preimages} 
Suppose that $X$ is geodesic, $\delta$-hyperbolic,  $Y, Z$ are $\la_Y, \la_Z$-quasicon\-vex subsets in $X$ respectively, such that $Z\subset Y$. Then 
$$
\hull(P_Y^{-1}(Z)) \subset P_Y^{-1}(N_D(Z))
$$
for some $D= D_{\ref{lem:preimages}}(\delta, \la_Y, \la_Z)$. 
\end{lemma} 
\proof Take two points $x, y\in P_Y^{-1}(Z)$  and a geodesic $\al=xy\subset X$ connecting these points. Then, for $\bar{x}=P_Y(x), \bar{y}=P_Y(y)$, 
by Lemma \ref{lem:small-projection}, $P_Y(\al)$ is   $D_{\ref{lem:small-projection}}(\delta,\la_Y)$-Hausdorff close to the segment $\beta=\bar x \bar y$. Since 
$\bar{x}, \bar{y}$ are in $Z$ and $Z$ is $\la_Z$-quasiconvex in $X$, $\beta\subset N_{\la_Z}(Z)$. Thus,
$$
\al\subset P_Y^{-1}(N_D(Z)), \quad D= D_{\ref{lem:preimages}}(\delta, \la_Y, \la_Z)= D_{\ref{lem:small-projection}}(\delta,\la_Y) + \la_Z. \qed 
$$

 \medskip
 For a pair of $\la$-quasiconvex subsets $U, V$ in a hyperbolic space $X$ there  is a basic dichotomy: Either $P_U(V)$ has uniformly bounded diameter 
 (in terms of $\la$ and $\delta$) or it is uniformly close to a quasiconvex subset of $V$. We will discuss this and related issues in more detail in Sections 
 \ref{sec:Projections and coarse intersections} and \ref{sec:Cobounded pairs of subsets}; for now, we prove this statement in the context of projections of geodesics to 
 geodesics:

\begin{lemma}\label{proj-geod}
For any $\delta\geq 0$ and $\la\geq 0$ there is are constants $D=D_{\ref{proj-geod}}(\delta,\la)$
$R=R_{\ref{proj-geod}}(\delta,\la)$ such that the following holds:

1. Suppose $X$ is a $\delta$-hyperbolic geodesic metric space and $Y\subset X$ is $\la$-quasiconvex.
Let $x,y\in X$ and let $\bar{x},\bar{y}\in Y$ be respectively their nearest-point projections to $Y$.
If $d(\bar{x},\bar{y})\geq D$ then $\bar{x}\bar{y}\subset N_R(xy)$. 
One can take $D=2\la+7\delta$ and $R=\la+5\delta$. 

2. When $Y$ is a geodesic, $\la=\delta$ and we can take: $D_{\ref{proj-geod}}(\delta,\delta)=8\delta$, 
$R_{\ref{proj-geod}}(\delta,\delta)=6\delta$.   
\end{lemma}
\proof  We prove Part (1) and leave computations in Part (2) (as a special case) to the reader. 
By Lemma \ref{lem:projection-1}(i) there is a point $z\in x\bar{y}$ such that $d(\bar{x},z)\leq \la+2\delta$.
Now we consider the geodesic triangle $\Delta xy\bar{y}$. Then $z\in N_{\delta}(xy\cup y\bar{y})$.
Suppose $z\in N_{\delta}(y\bar{y})$ and let $w\in y\bar{y}$ such that $d(z,w)\leq \delta$. Then
$d(w,\bar{x})\leq \la+3\delta$. Since $\bar{y}$ is a nearest point of $Y$ from $w$ it follows that
$d(w, \bar{y})\leq \la+3\delta$. Thus it follows that $d(\bar{x}, \bar{y})\leq 2\la+6\delta$. Hence, 
if $d(\bar{x}, \bar{y})\geq \delta+2\la+6\delta$ then $\bar{x}, \bar{y}\in N_{\la+3\delta}(xy)$.
Since geodesic quadrilaterals in $X$ are $2\delta$-slim, in that case it follows 
that $\bar{x}\bar{y}\subset N_{\la+5\delta}(xy)$. Hence, we may take
$D=2\la+7\delta$ and $R=\la+5\delta$. \qed 

 
 \section{Modified projection}

 In Section \ref{sec:ubiquity} we will need a minor modification of the projection $P_{X,Y}$ to quasiconvex subsets $Y\subset X$ 
 in the setting when $Y=T=T_p(xyz)$ is a $C$-tripod and also in the setting of quasigeodesic tripods.

 \begin{defn}\label{def:modified projection} \index{modified projection $\bar{P}$}
 Suppose that $T$ is a tripod as above in a $\delta$-hyperbolic geodesic space $X$. For a subset $U\subset X$, we define its 
 {\em modified projection} $\bar{P}_{X,T}(U)= 
 \bar{P}_T(U)\subset T$ as 
\begin{equation}\label{eq:barP} 
\bar{P}_{T}(U)= 
\hull_0(P_T(U)), 
\end{equation}
where the (closed) hull is taken  with respect to the intrinsic path-metric of $T$.
\end{defn}

If the nearest-point projection were continuous and $U$ were compact, then, of course, $\bar{P}_{T}(U)= P_T(U)$. 

\begin{lemma}\label{lem:proj-to-tripod} 
Suppose that $U$ is a $\la$-quasiconvex subset of a $\delta$-hyperbolic space $X$, $T= T_p(xyz)\subset X$ is a $C$-tripod. Then 
$$
\Hd_X(P_T(U), \bar{P}_T(U))\le C_{\ref{lem:proj-to-tripod}}(\delta,\la,C)=D_{\ref{lem:qc-projection}}(\delta,\delta,\la) + C+3\delta. 
$$
\end{lemma}
\proof First of all, each tripod $T\subset X$ is $\delta$-quasiconvex. Therefore, according to Lemma \ref{lem:qc-projection}, 
the projection $P_T(U)$ is a $\la'=D_{\ref{lem:qc-projection}}(\delta,\delta,\la)$-quasiconvex subset of $X$. Hence, for points $x'\in xp, y'\in yp$ which belong to $P_T(U)$, the segment $x'y'\subset X$ is contained in the $\la'$-neighborhood of $P_T(U)$. Since $p$ is a $C+2\delta$-center of the tripod $T_p(x'y'z)$ (Lemma \ref{lem:K-centers}), the geodesic segment $\ga\subset T$ connecting $x'$ to $y'$ is $C+3\delta$-Hausdorff close to $x'y'$. Thus, each 
 point $u\in \ga$ is within distance $C+3\delta$ from some $u'\in x'y'$ and, by the $\la''$-quasiconvexity of $P_T(U)$, there exists 
 $v\in P_T(U)$ such that $d(v,u')\le   \la'$. By the triangle inequality,
 $d(u,v)\le \la'+ C+3\delta$. Thus, $\bar{P}_T(U)\subset N_{\la'+ C+3\delta}(P_T(U))$. 
 Since $P_T(U)\subset \bar{P}_T(U)$, lemma follows. \qed 

\begin{rem}
Assuming that $\delta\ge 1$, we can take: 
\begin{align*}
C_{\ref{lem:proj-to-tripod}}(\delta,\la,C)= C+3\delta + 
(L_{\ref{lip-proj}}(\delta,\delta)+1)\la + D_{\ref{lem:small-projection}}(\delta,\delta)\le \\
C+3\delta + (12\delta+1)\la + 3\delta(2\la+9\delta)\le C+ 30\delta^2 +19\delta\la. 
\end{align*}
\end{rem}

\begin{rem}\label{rem:barP}
 As a special case, we will use the modified projection when the tripod $T$ is a single geodesic segment and $U$ is also a geodesic segment. 
 The estimate on the Hausdorff-distance in this situation is better; we leave it to the reader to verify 
 (analogously to the proof of Lemma  \ref{proj-geod})  that 
 $$
 C_{\ref{lem:proj-to-tripod}}(\delta,\delta,0)= 4\delta. 
 $$
 This estimate will be used in Section \ref{sec:proj-ladders}. 
\end{rem}

We will need a generalization of $\bar{P}$ and the lemma in the following setting: We let $Y$ be a union of three rectifiable arcs $\al\cup\beta\cup \ga$ in $X$, connecting points $x, y, z$ in $X$ to a certain point $p\in X$ and parameterized by their arc-length. Thus, $Y$, equipped with its intrinsic path-metric $d_Y$, is an abstract tripod. 
We assume that the inclusion map  $(Y, d_Y)\to X$ is a $K$-qi embedding. Note that if $\al,\beta,\ga$ are all geodesics in $X$ and $Y$ is a $C$-tripod $T$, 
then the inclusion map  $(T,d_T)\to X$ is a $4\delta$-qi embedding. As before, for $U\subset X$ we define $\bar{P}_Y(U)$ as the closed convex hull of  $P_Y(U)$ with respect to the metric $d_Y$. The next lemma follows from the fact that $Y$ is uniformly Hausdorff-close to the geodesic tripod $T= T_p(xyz)$; we leave a proof of the reader:

\begin{lemma}\label{lem:proj-to-tripod1} 
Suppose that $U$ is a $\la$-quasiconvex subset of a $\delta$-hyperbolic space $X$. Then 
$$
\Hd_{X}(P_T(U), \bar{P}_T(U))\le C_{\ref{lem:proj-to-tripod1}}(\delta,\la,K). 
$$
\end{lemma}

\section{Projections and coarse intersections}\label{sec:Projections and coarse intersections}

A basic fact of convex geometry is that intersections of convex subsets are again convex. In the context of 
quasiconvex subsets $U, V$ of a hyperbolic space, one needs to modify the notion of intersection. 
The most esthetically pleasing way to do so is to intersect $R$-neighborhoods of $U$ and $V$. 
However, most useful for us will be  {\em asymmetric} coarse intersections  $N_R(U)\cap V$. 
In this section  we discuss these in conjunction with the projections $P_U(V)$ and $P_V(U)$.

\begin{lemma}
[Coarse intersections with quasiconvex subsets are quasiconvex] 
\label{lem:coarse-intersections-are-qc-2} 
Suppose that $Y_i\subset X$ are $\lambda_i$-quasiconvex subsets in a $\delta$-hyperbolic space $X$, $i=1,2$. Then for every $\eps\ge \la_1+ \la_2+2\delta$, 
the intersection
$$
Y= N_\eps(Y_1)\cap Y_2
$$
is $\la_{\ref{lem:coarse-intersections-are-qc-2}}(\eps,\delta)$-quasiconvex, with $\la_{\ref{lem:coarse-intersections-are-qc-2}}= \eps+2\delta$. 
\end{lemma}
\proof Take two points $x, y\in Y$ and let $x_1, y_2\in Y_1$, be points within distance $\eps$ from $x, y$ respectively.
In view of the $2\delta$-slimness of quadrilaterals in $X$, for each $z\in xy$ either $d(z, x_1y_1)\le 2\delta$ or 
$d(z, x)\le 2\delta+\eps$  or $d(z, y)\le 2\delta+\eps$. In the last two cases, $z\in N_{2\delta+\eps}(Y)$. Suppose, therefore, 
that $d(z, x_1y_1)\le 2\delta$. By the $\la_i$-quasiconvexity of $Y_i$, there exist points $z_i\in Y_i$, such that $d(z,z_1)\le \la_1+2\delta$ and 
$d(z,z_2)\le \la_2$. In particular, $d(z_1,z_2)\le \la_1+\la_2+2\delta\le \eps$, i.e. $z_2\in Y$. Since 
$$
d(z,z_2)\le \la_2\le \eps+2\delta,
$$
we conclude that $Y$ is $(\eps+2\delta)$-quasiconvex. \qed

 \medskip 
The next lemma will be used in the proof of Theorem \ref{thm:mitras-projection}: 

\begin{lemma}\label{lemma0-flow-space}
Define the function
$$
R':=R_{\ref{lemma0-flow-space}}(R,\la,\delta)=  2\la+3\delta +R. 
$$
Let $U_1, U_2$ be $\la$-quasiconvex subsets of a $\delta$-hyperbolic space $X$ such that $d(U_1, U_2)\le R$. 
Then 
$$
P_{U_2}(U_1)\subset N_{R'}(U_1)\cap U_2
$$
and 
$$
\Hd( P_{U_1}(U_2), P_{U_2}(U_1))\le R'.$$
\end{lemma}
\proof The key is to  show that for every $a\in U_1$ its nearest-point projection $b=P_{U_2}(a)$ lies in the $R'$-neighborhood of $U_1$. 

Suppose $a_1\in U_1, b_1\in U_2$ are such that $d(a_1,b_1)\leq R$. By Lemma \ref{lem:projection-1}, 
there exists a point $c\in ab_1$ within distance $\la+2\delta$ from $b$. 
Since $d(a_1,b_1)\le R$, the $\delta$-slimness of the triangle $\Delta aa_1b_1$ implies existence of a point $c_1\in aa_1$ within distance $R+\delta$ from $c$. 
In view of the $\la$-quasiconvexity of $U_1$, $c_1$ belongs to the $\la$-neighborhood of $U_1$. Thus,
$$
b\in N_{2\la+3\delta +R}(U_1)= N_{R'}(U_1). 
$$
Therefore, the distance from $b$ to $P_{U_2}(b)$ is at most $R'$, verifying the inclusion
$$
P_{U_2}(U_1)\subset N_{R'}(P_{U_1}(U_2)). 
$$
The reverse inclusion is proven by switching the roles of $U_1$ and $U_2$. \qed 

\medskip
Continuing with the notation of the lemma: 

\begin{cor}\label{cor:0-flow-space}
If $d(U_1,U_2)\le R$, then 
$$\Hd(P_{U_1}(U_2), N_{R'}(U_2)\cap U_1)\le 2R' \hbox{~~and~~} 
\Hd(N_{R'}(U_1)\cap U_2, N_{R'}(U_2)\cap U_1)\le R'.$$ 
\end{cor}
\proof 1. According to the lemma, $P_{U_1}(U_2)\subset N_{R'}(U_2)\cap U_1$. Conversely, given 
$x\in N_{R'}(U_2)\cap U_1$, there exists $y\in U_2$ with $d(x, y)\le R'$. Hence, 
$d(y, P_{U_1}(U_2)(y))\le R$, implying 
$$
d(x, P_{U_1}(U_2)(y))\le 2R'.  
$$
2. The second claim is clear and holds for arbitrary $R'\ge 0$ and arbitrary subsets of arbitrary metric spaces. \qed

\medskip
Thus, we proved that if  $d(U_1,U_2)\le R$, then 
all four subsets 
$$P_{U_1}(U_2), P_{U_2}(U_1), N_{R'}(U_1)\cap U_2, 
N_{R'}(U_2)\cap U_1$$ are within Hausdorff distance $2R'$ from each other.

 \section{Quasiconvex subgroups and actions}

In this section we discuss quasiconvexity in the context of subgroups of hyperbolic groups and, more generally, group actions.

\begin{defn}
\label{defn:qc-subgroup} 
A subgroup $H$ of a hyperbolic group $G$ is said to be {\em quasiconvex} \index{quasiconvex subgroup}  \index{quasiconvex action} 
if it is a quasiconvex subset of a Cayley graph of $G$ for a finite generating set. More generally, 
a (metrically) proper isometric action of a discrete group $H$ on a geodesic hyperbolic metric space $X$ is {\em quasiconvex} if one (equivalently, every) $H$-orbit in $X$ is a quasiconvex subset in $X$. 
\end{defn}

\begin{lemma}\label{lem:qc action}
If the action of $H$ on $X$ is {\em quasiconvex} then $H$  is  finitely generated, the orbit map 
$H\to H\cdot x\subset X$ is a qi embedding and $H$ is a hyperbolic group. 
\end{lemma}
\proof 1. Quasiconvexity of $Hx\subset X$ implies that $Hx$ is coarsely connected. Hence, by the Milnor--Schwarz Lemma, $H$ is finitely generated. We, thus, equip $H$ with a word metric corresponding to a finite generating set. 

2. Metric properness  of the action implies that the orbit map $o_x: H\to Hx\subset X$ is 
uniformly proper. Since the image of this map is a quasiconvex subset of $X$, the orbit map is a qi embedding (see Lemma \ref{lem:up+qc->qi}). 

3. Since $X$ is assumed to be hyperbolic, in view of Lemma 
\ref{lem:qi-preserves-paths}, the existence of a qi embedding $o_x$ implies hyperbolicity of $H$.  \qed

\medskip
We now discuss the notion of coarse intersection in relation to quasiconvex subgroups and actions.

For general quasiconvex subsets $U, V$ of hyperbolic spaces $X$, coarse intersections  $N_R(U)\cap V$ might  
not be Hausdorff-close to the actual intersections $U\cap V$: For instance, $U\cap V$ might be empty while for 
some $R>0$ the intersection $N_R(U)\cap V$ might be unbounded. As a specific example, consider $X=\R$ 
(which is $0$-hyperbolic) and $1$-quasiconvex subsets $U, V$ consisting of odd/even integers respectively. Then 
$N_1(U)\cap V=V$, while $U\cap V=\emptyset$. Nevertheless, in the group-theoretic setting we have

\begin{lemma}\label{lem:qc-subgroups} 
Let $G$ be a hyperbolic group, $U, V$ be quasiconvex subgroups in $G$ with $W:= U\cap V$. 
Then for every $r>0$ there exists $R=R_{\ref{lem:qc-subgroups}}(G, r)$ such that 
$$
W\subset W_r:=V\cap N_r(U)\subset N_R(W). 
$$
\end{lemma}
\proof The proof is quite standard, cf. \cite[pp. 164-165]{gromov-ai} or Lemma 2.6 in \cite{MR1389776}. 
Suppose that $u\in U, v\in V$ satisfy $d_G(u,v)\le r$, i.e. $u^{-1}v\in B_G(1, r)$. Since the ball $B_G(1, r)$ is finite, 
there exists a finite set of pairs $(u_i, v_i), i=1,...,n$ such that for any pair $u\in U, v\in V$ within distance $r$ from 
each other,  the product  $u^{-1}v$ equals  $u_i^{-1}v_i$ for some $i\in \{1,...,n\}$. We have 
$$
u^{-1}v=u_i^{-1}v_i \Rightarrow u u_i^{-1}= v v_i^{-1}=w\in W=U\cap V.
$$
Hence, $u= w u_i, v=wv_i$ and, therefore, for 
$$
R:= \max \{|v_i|: i=1,...,n\}
$$
we have $d_G(v, W)\le R$. \qed

\begin{cor}\label{cor:proj-to-qc-subgroup}
In the setting of Lemma \ref{lem:qc-subgroups}, the distance between the restrictions to $V$ 
of the projections $P_{X,W}, P_{X,U}$,  is at most $C_{\ref{cor:proj-to-qc-subgroup}}(G, \delta,\la)$, where $\delta$ is the hyperbolicity constant of $G$ and $\la$ is the  maximum of the quasiconvexity constants of $U, V$ in $G$. 
\end{cor}
\proof We take $r:= r_{\ref{lem:two-projections}}(\la,\la,\delta)$. According to 
Lemma \ref{lem:two-projections}, the 
restrictions to $V$ of the projections 
$$
P_{G, W_r},  P_{G,U}
$$ 
are within distance $\mu+3\delta+r$. By Lemma \ref{lem:qc-subgroups}, the subsets 
$W, W_r\subset X$ are $R=R_{\ref{lem:qc-subgroups}}(G, r)$-Hausdorff-close.  Therefore, by Corollary 
\ref{cor:proj-to-close-subsets}, the distance between the projections $P_{G, W_r}, P_{G, W}$ is 
$\le D_{\ref{cor:proj-to-close-subsets}}(\delta, \la, R)$. Thus, we can take 
$$
C_{\ref{cor:proj-to-qc-subgroup}}(G, \delta,\la):= D_{\ref{cor:proj-to-close-subsets}}(\delta, \la, R) + \mu+3\delta+r. \qed 
$$

As an immediate consequence we obtain the standard result on  intersections of quasiconvex subgroups of hyperbolic groups (see e.g. \cite{short} or \cite[Proposition 4.13]{bridson-haefliger}):

\begin{cor}\label{cor:qc-in} 
If $G$ is a hyperbolic group and $U, V$ are quasiconvex subgroups of $G$, then $U\cap V$ is also a quasiconvex subgroup of $G$, $H$ and $V$. 
\end{cor}
\proof This is a combination of Lemmata  \ref{lem:qc-subgroups}, \ref{lem:coarse-intersections-are-qc-2} and \ref{lem:nbds-of-qc}(2). 
 \qed

\medskip 
Essentially the same proofs as above work in the more general setting, when we have a quasiconvex  
 metrically proper action of a hyperbolic group $G$ on a $\delta$-hyperbolic geodesic  metric space $X$, and 
 $Y\subset X$ is a quasiconvex subset with locally finite $G$-orbit (see Definition \ref{defn:locally finite action}).

\begin{prop}\label{prop:proj-to-qc-action}
Let $H< G$ denote the stabilizer of $Y$ in $G$, $x\in Y$. Then: 

a. There exists a function $R=R_{\ref{prop:proj-to-qc-action}}(x, r)$ such that 
$$
Hx\subset Gx\cap N_r(Y)\subset N_R(Hx)
$$

b. $Hx$ is a quasiconvex $\mu$-subset of $X$, $\mu=\mu_{\ref{prop:proj-to-qc-action}}(x, \delta, \la)$, where $\la$ is the maximum of quasiconvexity constants (in $X$) of $Gx$ and $Y$.  

c. $H$ is a quasiconvex subgroup of $G$. 

d. The restrictions of $P_{X,Y}$ and $P_{X, Hx}$ to the orbit $Gx\subset X$ are within distance 
$C=C_{\ref{prop:proj-to-qc-action}}(x,\delta,\la)$.  
\end{prop} 
\proof a. Suppose that $d_X(gx,Y)\le r$ for some $g\in G$; equivalently, $g^{-1}Y\cap B(x,r)\ne \emptyset$. By the definition of a local finiteness of the $G$-orbit $GY$, there exist  
$h^{-1}\in H$ and $g_i^{-1}\in S$, where $S\subset G$  is a finite subset depending only on $x$ and $r$, such that $g^{-1}= g_i^{-1}h^{-1}$. We let  $R=R(x,r)$ be the maximum of distances $d(x, g_i^{-1}(x))$ taken over $g_i^{-1}\in S$. Then $d(hx, gx)\le R$. This proves (a). 

b. We take $r:= 2\la + 2\delta$. By Lemma \ref{lem:coarse-intersections-are-qc-2}, the coarse intersection $Gx\cap N_r(Y)$ is 
$\la_{\ref{lem:coarse-intersections-are-qc-2}}(r,\delta)$-quasiconvex in $X$. On the other hand, 
by Part (a),
$$
\Hd(Hx, Gx\cap N_r(Y))\le R.
$$
Therefore, by Lemma \ref{lem:nbds-of-qc}(2), the subset $Hx$ is $\mu$-quasiconvex in $X$ with 
$$
\mu= 2R+2\delta+\la_{\ref{lem:coarse-intersections-are-qc-2}}(r,\delta).$$

c. Since the actions of $G$ and $H$ on $X$ are quasiconvex, the orbit maps $o_x: G\to X$, $o_x: H\to X$ are qi embeddings (see Lemma \ref{lem:qc action}). From this, we conclude that $H$ is qi embedded and, hence, is quasiconvex in $G$. 

d. The proof of this part is identical to that of Corollary \ref{cor:proj-to-qc-subgroup} and we omit it. 
\qed 

\medskip
We assume now that $X$ is hyperbolic and that for each point $x\in X$ and each ideal boundary point $\xi\in \geo X$, there exists a geodesic ray $x\xi$ connecting $x$ to $\xi$ (e.g. $X$ is a proper geodesic metric space). 

\begin{defn}\label{defn:conical-limit} \index{conical limit point} \index{limit point} 
Suppose that $G$ acts isometrically and properly on $X$. A point $\xi\in \geo X$ is called a {\em limit point} of this action (or, simply, a limit point of $G$) if there exists a sequence $g_i\in G$ such that for some (equivalently, every) $x\in X$, the sequence $(g_i(x))$ converges to $\xi$. A limit point $\xi$ is called {\em conical} if the sequence $(g_i)$ can be chosen so that 
for some (equivalently, all) $x\in X, y\in X$, there exists a constant $R$ such that $d(g_i y, x\xi)\le R$ for all $i$.   
\end{defn}

The proof of the following result (a {\em Beardon--Maskit criterion} for quasiconvexity) can be found for instance in Swenson's paper \cite{MR1804703} (cf. also \cite{MR1317633, bowditch-cgnce,MR1637829}); we will only need the easier direction (every limit point of a quasiconvex action is conical):  

\begin{thm}
Suppose that $X$ is a proper geodesic hyperbolic metric space. Then a proper isometric action of a discrete group $G\acts X$ is quasiconvex if and only if every limit point of $G$ is conical. 
\end{thm}

\section{Cobounded pairs of subsets}\label{sec:Cobounded pairs of subsets}

Recall that in Definition \ref{def:cob} we defined Lipschitz-cobounded pairs of subsets of general metric spaces. Below, we establish a characterization of cobounded pairs in hyperbolic spaces.

\begin{prop}[Characterizations of cobounded pairs] 
\label{prop:cobounded2} 
The following are equivalent for $\la$-quasiconvex subsets $Y, Z\subset X$ in a $\delta$-hyperbolic geodesic metric space $X$:

\begin{enumerate}
\item $Y, Z$ are $C_1$-Lipschitz cobounded. 

\item For every $R$ there exists $D=D(R)$ such that if 
$$
a_i\in Y, b_i\in Z, i=1, 2, 
$$
satisfy $d(a_i, b_i)\le R$, $i=1,2$, then $d(a_1,a_2)\le D, d(b_1, b_2)\le D$. 

\item The diameters of nearest-point projections $P_{X,Y}(Z)$, $P_{X,Z}(Y)$ are $\le C_2$. 
\end{enumerate} 

Moreover, once a constant $C_i$ (or a function $D(R)$) in one of the items is chosen, this, together with $\delta$ and $\la$, 
determines the constant/function in the other two items. 
\end{prop}
\proof The implication (1)$\Rightarrow$(2) is proven in Lemma \ref{lem:cobounded1} for arbitrary subsets of arbitrary metric spaces, with 
$$
D=2C_1(R +1) +C_1
$$

For the implication (2)$\Rightarrow$(3), consider points $a_i\in Y, b_i\in Z$ such that $a_i\in P_{X,A}(b_i), i=1,2$. 
By Lemma \ref{proj-geod}, if $d(a_1,a_2)\geq D_{\ref{proj-geod}}(\delta, \la)$ then there exists 
$R=R_{\ref{proj-geod}}(\delta,\la)$ such that
$$
a_1a_2\subset N_R(b_1b_2).
$$
In that case, there are points $b_i'\in Z$ within distance $\le R+\la$ from $a_i$, $i=1,2$.  Then, by (2), 
$$
d(a_1, a_2)\le D(R+\la). 
$$ 
Hence, we can take $C_2= \max\{D_{\ref{proj-geod}}(\delta, \la), D(R+\la)\}$. 

For the implication (3)$\Rightarrow$(1), we can take the retractions 
$$
r_A:= P_{X,A},  \quad r_B:= P_{X,B}.   \qed
$$

\medskip 
In view of this proposition, for quasiconvex subsets of hyperbolic spaces we will adopt the following terminology:

\begin{defn}\label{defn:hyp-cobounded} \index{$C$-cobounded subsets in a hyperbolic space}
A pair of subsets $Y, Z\subset X$ in a hyperbolic space $X$ is $C$-cobounded if the diameters of the  projections $P_{X,Y}(Z), P_{X,Z}(Y)$ are $\le C$.  
\end{defn}

\begin{lemma}\label{cobdd-cor}
Given $\delta\geq 0$ and $\la\geq 0$ there are constants $R=R_{\ref{cobdd-cor}}(\delta, \la)$ and 
$D=D_{\ref{cobdd-cor}}(\delta, \la)$ such 
that the following  holds:

Suppose $X$ is a $\delta$-hyperbolic metric space and $Y, Z\subset X$ are two $\la$-quasiconvex and $R$-separated subsets.
Then $Y, Z$ are $D$-cobounded. In fact, one can take $D=2\la+7\delta$ and $R=  2\la+5\delta$. 
\end{lemma}
\proof We will show that the choice of $D=D_{\ref{proj-geod}}(\delta,\la)=2\la+7\delta$ and 
$$
R=\la+R_{\ref{proj-geod}}(\delta, \la)=  2\la+5\delta$$
 works.
Let $R_1=R_{\ref{proj-geod}}(\delta, \la)$, so that $R=\la+R_1$.
Suppose the diameter of $P_{X,Z}(Y)$ is greater than or equal to $D$. Let $x,y\in Y$ be such that
$d(P_{X,Z}(x), P_{X,Z}(y))\geq D$. Then by Lemma \ref{proj-geod} 
$P_{X,Z}(x)\in N_{R_1}(xy)$. But $Y$ is $\la$-quasiconvex and $x,y\in Z$. It follows that 
$P_{X,Z}(x)\in N_R(Y)$. Thus if $Y,Z$ are $R$-separated then the diameter of $P_{X,Y}(Z)$ and $P_{X,Z}(Y)$ are both
less than $D$. 
\qed

\medskip 
A consequence of this lemma allows one to simplify the verification that two subsets are cobounded; namely, it suffices to check only that one projection is bounded:

\begin{cor}\label{cor:cob-char} 
Suppose that $U, V\subset X$ are $\la$-quasiconvex subsets in a $\delta$-hyperbolic space. 

a. If $\diam(P_U(V))\le D$, then $\diam(P_V(U))\le C=C_{\ref{cor:cob-char}}(\la,\delta,D)$, where $D\le C$. In particular, 
the pair $(U, V)$ is  $C$-cobounded. 

b. If the pair $U, V\subset X$ is not $D_{\ref{cobdd-cor}}(\delta, \la)$-cobounded then 
$$
\Hd(P_{U}(V), P_{V}(U))\le D_{\ref{cor:cob-char}}(\delta,\la)= 
R_{\ref{cor:cob-char}}(\delta, \la)=2\la+3\delta + R_{\ref{cobdd-cor}}(\delta, \la). 
$$
\end{cor}
\proof a. There are two cases to consider:

1. If $d(U,V)\ge R=R_{\ref{cobdd-cor}}(\delta, \la)$, then the pair $(U,V)$ is $D_{\ref{cobdd-cor}}(\delta, \la)$-cobounded by Lemma \ref{cobdd-cor}. 

2. Suppose that $d(U,V)\le R=R_{\ref{cobdd-cor}}(\delta, \la)$. Then by Lemma \ref{lemma0-flow-space}, 
$$
\Hd(P_{U}(V), P_{V}(U))\le R'=2\la+3\delta +R.$$
Since $\diam(P_U(V))\le D$, it follows that 
$$
\diam(P_V(U))\le D+R'. 
$$
Taking $C:= \max(D+R', D_{\ref{cobdd-cor}}(\delta, \la))$, concludes the proof of a. 

\begin{rem}
Note that $C= \max(D+4\la+8\delta, D_{\ref{proj-geod}}(\delta, \la))$. 
\end{rem}

b. By the argument in Part a1, since the pair $U, V\subset X$ is not $D_{\ref{cobdd-cor}}(\delta, \la)$-cobounded, 
$d(U,V)< R=R_{\ref{cobdd-cor}}(\delta, \la)$. Thus, as in Part a2, 
$$
\Hd(P_{U}(V), P_{V}(U))\le R'=2\la+3\delta +R= 2\la+3\delta + R_{\ref{cobdd-cor}}(\delta, \la). \qedhere 
$$

\begin{rem}\label{rem:cbb-geodesics} 
If $U_1, U_2$ are geodesics in $X$, $\la=\delta$ and, 
by Lemma \ref{proj-geod}, one can take $D_{\ref{cobdd-cor}}(\delta, \delta)= 8\delta$ and 
$R_{\ref{cor:cob-char}}(\delta, \la)=12\delta$. 
\end{rem}

\medskip 
Another application of Lemma \ref{cobdd-cor} is:

\begin{cor}\label{cor:noncbd} 
Suppose that $\la$-quasiconvex subsets $U_1, U_2\subset X$ are not $D=D_{\ref{cobdd-cor}}(\delta, \la)$-cobounded. Then 
$$
P_{U_2}(U_1)\subset N_{4\la+8\delta}(U_1)\cap U_2. 
$$
\end{cor}
\proof By Lemma \ref{cobdd-cor}, since $U_1, U_2$ are not $D$-cobounded, then
$$
d(U_1,U_2)\le R=R_{\ref{cobdd-cor}}(\delta, \la)=2\la+5\delta. 
$$
According to Lemma \ref{lemma0-flow-space}, 
$$
P_{U_2}(U_1)\subset N_{R'}(U_1)\cap U_2,
$$
where $R'=2\la+3\delta +R=4\la+8\delta$.  \qed

\medskip 

\begin{lemma}\label{cobdd-lem1}
Given $\delta\geq 0$, $\la\geq 0$ and $C\ge 0$, there exists a constant $D= D_{\ref{cobdd-lem1}}(\delta, \la, C)$ 
such that the following  holds: 

Suppose $X$ is a $\delta$-hyperbolic metric space and $U,V\subset X$ are two nonempty $\la$-quasicon\-vex and $C$-cobounded subsets. 
Then there are points $x_0\in U_0=P_U(V)\subset U$, $y_0\in V_0=P_V(U)\subset V$,    
such that $x_0y_0\subset N_D(xy)$, for all $x\in U$ and $y\in V$. 
\end{lemma} 
\begin{proof} 
 Since the pair $(U, V)$ is $C$-cobounded, 
$$
\diam(V_0)\le C, \quad \diam(U_0)\le C. 
$$ 
Choose any pair of points $x_0\in U_0$, $y_0\in V_0$. Take $x\in U, y\in V$ and consider the points $\bar{x}= P_V(x)\in V_0, \bar{y} = P_U(y)\in U_0$. By Lemma 
\ref{lem:projection-1}, the points $\bar{x}, \bar{y}$ are within distance $\la+2\delta$ from $xy$. Therefore,
$$
\max( d(x_0, xy), d(y_0, xy)) \le \la+2\delta+C 
$$
and, hence, we can take $D=\la+4\delta+C$. \end{proof}

\begin{cor}\label{cobdd-cor1}
Given $\delta\geq 0$ and $\la\geq 0$, there are constants $R=R_{\ref{cobdd-cor1}}(\delta, \la)$ and $D=D_{\ref{cobdd-cor1}}(\delta, \la)$ such 
that the following  holds:

Suppose $X$ is a $\delta$-hyperbolic metric space and $U,V\subset X$ are two $\la$-quasiconvex and $R$-separated subsets.
Then there are points $x_0\in U$, $y_0\in V$ such that $x_0y_0\subset N_D(xy)$, for all $x\in U$ and $y\in V$.
\end{cor}
\begin{proof} 
By Lemma \ref{cobdd-cor}, there exists $R =R_{\ref{cobdd-cor}}$ such that the pair $(U,V)$ is 
$C=D_{\ref{cobdd-cor}}$-cobounded 
whenever $U,V$ are  $R$-separated. Now, the claim follows from Lemma \ref{cobdd-lem1}.  
\end{proof}

\chapter{Graphs of groups and trees of metric spaces}\label{ch:trees}

\section{Generalities} \label{sec:generalities}

We presume that the reader is familiar with the Bass--Serre theory. However, we  briefly recall some of the concepts
that we shall need. For details we refer the reader  to Section 5.3 of Serre's book \cite{serre-trees}. 

\begin{defn}[Graph of groups] \index{graph of groups} 
A {\bf graph of groups} $(\mathcal G,\Ga)$ consists of the following data: 

(1) A connected graph $\Ga$. 

(2) An assignment to each vertex $v\in V(\Ga)$ (and edge $e\in E(\Ga)$) of a group $G_v$ (respectively $G_e$) together with injective homomorphisms $\phi_{e,o(e)}: G_e\map G_{o(e)}$ and $\phi_{e,t(e)}: G_e\map G_{t(e)}$ for all $e\in E(\Ga)$, 
such that the following conditions hold: 

(i) $G_e=G_{\bar{e}}$,

(ii) $\phi_{e,o(e)}=\phi_{\bar{e},t(\bar{e})}$ and $\phi_{e,t(e)}=\phi_{\bar{e},o(\bar{e})}$. 
\end{defn}

We shall refer to the maps $\phi_{e,v}$ as the {\em canonical maps} of the graph of groups.
We shall refer to the groups $G_v$ and $G_e$, $v\in V(\Ga)$ and $e\in E(\Ga)$ as {\em vertex groups} and {\em edge groups} respectively.
For topological motivations of graph of groups and the following definition of the fundamental group of a graph of groups
one is referred to \cite{scott-wall} or \cite{altop-hatcher}. In the terminology of \cite{bridson-haefliger}, a graph of groups is a covariant functor from the graph $\Ga$ (regarded as a small category with set of objects $E\sqcup V$ and the set of 
morphisms consisting of the maps $o$ and $t$) to the category of groups, sending morphisms $o, t$ to group-monomorphisms.  Functorially, in the case when $\Ga$ is a tree, one can define the group $G$, the {\em fundamental group} of $(\mathcal G, \Ga)$, or the {\em pushout} of the diagram ${\mathcal G}$, by a universal property. Namely, there exist monomorphisms $G_e\to G, G_v\to G$ such that the diagrams 
$$
\begin{diagram}
G_e & \rTo & G_v\\
& \rdTo & \dTo \\
       & \rdTo&  G
\end{diagram}
$$
commute, and, whenever we have a group $H$ and a compatible collection of homomorphisms $G_e\to H, G_v\to H$ forming commutative diagrams
$$
\begin{diagram}
G_e & \rTo & G_v\\
& \rdTo & \dTo \\
       & \rdTo&  H
\end{diagram}
$$
there is a unique homomorphism $G\to H$ forming commutative diagrams
$$
\begin{diagram}
G_e & \rTo & G\\
& \rdTo & \dTo \\
       & \rdTo&  H
\end{diagram} \quad \hbox{and}\quad  
\begin{diagram}
G_v & \rTo & G\\
& \rdTo & \dTo \\
       & \rdTo&  H
\end{diagram} 
$$

The general definition is more complicated: 

\begin{defn}[Fundamental group of a graph of groups]\index{fundamental group of a graph of groups} \label{defn:fundamental group of graph of groups}
Suppose $(\mathcal G, \Ga)$ is a graph of groups and 
let $S\subset \Ga$ be a maximal (spanning) subtree.
Then the fundamental group $G=\pi_1(\mathcal G,\Ga, S)$ of $(\mathcal G,\Ga)$ is defined in terms of generators
and relators as follows:

The generators of $G$ are the elements of the disjoint union of the generating sets of the vertex groups $G_v$, $v\in V(\Ga)$
and the set $E(\Ga)$ of {\em oriented} edges of $\Ga$. 

The relators are of four types: 

(1) Those coming from the vertex groups; 

(2) $\bar{e}=e^{-1}$ for all edge $e$;  

(3) $e=1$ whenever $|e|$ is a unoriented edge of $S$; 

(4) $e\phi_{e,t(e)}(a)e^{-1}=\phi_{e,o(e)}(a)$ for all oriented edges $e$ and $a\in G_e$.
\end{defn}

The group $G$ does not depend on the choice of $S$ and it will be denoted  $G=\pi_1(\mathcal G)$ in what follows. We will also frequently suppress the letter $\Ga$ in the notation of a graph of groups.

\begin{defn}[Bass--Serre tree of a graph of groups]\index{Bass--Serre tree}
Suppose $(\mathcal G, \Ga)$ is a graph of groups 
and let $S$ be a maximal tree in $\Ga$ as in the above definition.
Let $G=\pi_1(\mathcal G, \Ga, S)$ be the fundamental group of this graph of groups. 
The {\em Bass--Serre tree}, denoted $T$, is the tree with the 
vertex set 
$$\coprod_{v\in V(\Ga)} G/G_v$$ 
and the edge set 
$$
\coprod_{e\in E(\Ga)} G/G^e_e$$ where $G^e_e=\phi_{e,t(e)}(G_e)<G_{t(e)}$. The origin/terminus maps  are given by
$$
t(gG^e_e)=g G_{t(e)},\,o(gG^e_e)=gG_{o(e)}.$$
Note that whenever $|e|$ is a unoriented edge of $S$, then we have $e=1$ in $G$. The group $G$ acts on $T$ via left multiplication. 
\end{defn}

Conversely, given an action {\em without inversions}\footnote{which means that if $g\in G$ preserves an edge $[v, w]$ of $T$, then it also fixes both $v$ and $w$} of a group $G$ on a tree $T$, there exists a  graph of groups ${\mathcal G}$ with $\pi_1({\mathcal G})\cong G$ such that $T$ is equivariantly isomorphic to the Bass--Serre tree of $\mathcal G$, see \cite{serre-trees}.

Since our main motivation comes from geometric group theory and, hence, finitely generated groups, we observe that for $G=\pi_1(\mathcal G, \Ga, S)$ to be finitely generated, it suffices (not not necessary!) to assume that each vertex group $G_v$ is finitely generated and the graph $\Ga$ is finite. On the other hand, the edge groups $G_e$ need not be finitely generated. Natural examples of the latter situation are given by amalgams 
$$
G=G_v \star_{G_e} G_w, 
$$
where $G_e$ is an infinite rank free subgroup in two finitely-presented groups $G_v, G_w$: Such groups $G$ are finitely generated but not finitely presentable. In the context of combination theorems for hyperbolic groups, one assumes that the graph $\Ga$ is finite, each vertex/edge group is hyperbolic and the monomorphisms $\phi$ are qi embeddings, i.e. have quasiconvex images. 

Returning to the general setting with finitely generated vertex groups and finite graph $\Ga$, we note that while it is meaningless to assume that the canonical maps $\phi$ are uniformly proper (as edge-groups do not have canonical qi classes of metrics), 
nevertheless, if we equip $G_e$ with the pull-back a word metric from $G_{o(e)}$, while $G_{t(e)}$ has a word metric coming from a finite generating set, then the monomorphism $G_e\to G_{t(e)}$ is uniformly proper. Since the graph $\Ga$ is finite, we conclude that 
each edge group has a left-invariant proper metric, such that the homomorphisms $\phi_{e,o(e)}$ and $\phi_{e,t(e)}$ are 
$(\eta,L)$-proper for some uniform function $\eta$ and a constant $L$. 

\medskip
A {\em morphism} of graphs of groups, $\Psi: {\mathcal G}\to {\mathcal G}'$, consists of a morphism of the underlying graphs $\psi: \Ga\to \Ga'$ together with a collection of group homomorphisms 
$$
\Psi_v: G_v\to G_{\psi(v)}, v\in V(\Ga), \quad \Psi_e: G_e\to G_{\psi(e)}, e\in E(\Ga)
$$
such that the following diagrams are commutative for $v=o(e)$ and $w=t(e)$ and their respective images $v'=\psi(v), w'=\psi(w), e'=\psi(e)$:
$$
\begin{diagram}
G_e & \rTo & G_{e'}\\
\dTo_{\phi_{e,v}} &  & \dTo_{\phi_{e',v'}} \\
G_v & \rTo & G_{v'}\\
\end{diagram}, 
 \begin{diagram}
G_e & \rTo & G_{e'}\\
\dTo_{\phi_{e,w}} &  & \dTo_{\phi_{e',w'}} \\
G_w & \rTo & G_{w'}\\
\end{diagram}
$$

Given a graph of groups $({\mathcal G}',\Ga')$ and a graph-morphism $\psi: \Ga\to \Ga'$ from a connected graph $\Ga$, there is a canonical {\em pull-back} graph of groups $({\mathcal G}, \Ga)$ and a 
morphism of graphs of groups $\Psi: {\mathcal G}\to {\mathcal G}'$, such that the underlying 
morphism of graphs $\Ga\to \Ga'$ equals $\psi$. In the special case when $\Ga$ is a connected subgraph 
of $\Ga'$, the graph of groups $({\mathcal G}, \Ga)$ is called the {\em restriction} of 
${\mathcal G}'$ to $\Ga$ (see \cite[2.15]{Bass-93}). In this case, the Bass--Serre tree $T$ of $({\mathcal G}, \Ga)$ admits a 
$G=\pi_1({\mathcal G})$-equivariant embedding in the Bass--Serre tree $T'$ of $({\mathcal G}', \Ga')$ and $G$ equals the stabilizer of $T\subset T'$ in $G'$.  We refer the reader to \cite{Bass-93} for further discussion of morphisms of graphs of groups.

In the book, on several occasions we will use the following definition from the theory of group actions on trees:

\begin{defn}\label{defn:acylindrical-action} \index{$k$-acylindrical group action}
An action of a group $G$ on a tree $T$ is said to be {\em $k$-acylindrical} if whenever a 
nontrivial\footnote{In the literature, acylindricity is sometimes defined by requiring only that $G$-stabilizers of intervals of length $\ge k$ are finite, rather than trivial, subgroups.} 
element $g\in G$ fixes element-wise an interval $J\subset T$, then $J$ has length $\le k$. 
\end{defn}

 This terminology 
originates in Sela's paper \cite{Sela97}. The definition of acylindrical actions on trees was later coarsified and generalized by Bowditch in \cite{Bowditch-08}; 
we will not use his generalization.

\section{Trees of spaces} \label{sec:trees-of-spaces}

Each graph of groups  yields a ``tree of metric spaces'' over its Bass--Serre tree; this was first formalized and used by Bestvina and Feighn in \cite{BF}. Below is our version  
of their definition. 

We start with the simpler concept of a {\em tree of topological spaces}. 
One can regard a (simplicial) tree $T$ (or a general graph) as a small category with object sets equal to $V(T)\sqcup E(T)$ and morphisms given by origin/terminus arrows. 
Then a tree of topological spaces over a tree $T$ is a functor $\X$ from $T$ to the category of topological spaces. More explicitly:

\begin{defn}\label{defn:top-tree}\index{tree of topological spaces} 
A tree of topological spaces over a tree $T$  is a   collection $\X$ of nonempty topological  
spaces (vertex and edge-spaces) $X_v, v\in V(T), X_e, e\in E(T)$, together a collection of 
continuous  {\em incidence maps} $f_{ev}: X_e\to X_v$ defined for each oriented edge $e=[v,w]$. The {\em total space} $X$ of $\X$ is the mapping cylinder of the collection of the maps $f_{ev}$, i.e. the quotient of the disjoint union 
$$
\coprod_{v\in V(T)}  X_v \sqcup \coprod_{e\in E(T)} X_e\times [0,1]
$$
by the equivalence relation
$$
(x,0)\sim f_{ev}(x), (x,1)\sim f_{ew}(x), e= [v,w]\in E(T). 
$$
 \end{defn}

We will use trees of topological spaces in Section \ref{sec:CT-lamination}. For most of the book, we will work with trees of metric spaces defined below.

\begin{comment}

\begin{defn}\label{defn:retractive tree} 
A tree of spaces $\X= (\pi: X\to T)$ is called 
 $L'$-{\em retractive} if with the last (uniform properness) condition can be replaced by the following stronger property:  
For every oriented edge $e=[v,w]$, the restriction of the attaching map $f_e|_{ X_e \times\{v\}}$, admits an  
$L'$-coarse Lipschitz retraction  
$$
f_{ve}: X_{v}\to X_e. 
$$
The number $L'$ is the {\em retractivity constant} of $\X$. 
\end{defn}

In particular, for such a tree of spaces the attaching maps are  $K$-quasiisometric embeddings to $(X_{v}, d_{X_{v}})$, where
$$
K= \max(L, 2L').  
$$

We will see below (Proposition \ref{unif-emb-subtree}) 
that uniform properness of the inclusion maps $X_v\to X$ is a consequence of the 
existence of coarse $L'$-retractions $X_{v}\to X_e$.

Later, we will impose further {\em hyperbolicity} restrictions on trees of spaces.  But first, we relate trees of metric spaces to a slightly different concept, which we call {\em an abstract tree of spaces}.  In the second definition, instead of starting with a metric space equipped with a certain collection of maps, we start with a collection of metric spaces and {\em incidence maps} and from that produce a tree of spaces.

\end{comment}

Again, regarding a tree $T$ as a small category, to {\em some degree}, a tree of metric spaces $\X$ over a tree $T$ is a functor from $T$ to the coarse category $\CC$, see Remark \ref{rem:category}. The actual definition is somewhat more restrictive:

\begin{defn}[Abstract tree of spaces] \label{defn:abstract-tree-of-spaces} 
An {\em abstract tree of (metric) spaces} $\X$ over a simplicial tree $T$, is a   collection of nonempty  
metric  spaces (vertex and edge-spaces) 
$X_v, v\in V(T), X_e, e\in E(T)$, together a collection of $\psi$-uniformly proper coarse $L$-Lipschitz 
 {\em incidence maps} $f_{ev}: X_e\to X_v$ defined for each oriented edge $e=[v,w]$. The constant $L$ and the function $\psi$ are the {\em parameters} of the abstract tree of spaces $\X$. The tree $T$ is the {\em base} of $\X$. 
 \end{defn}
 
 Throughout the book, we will be assuming that all vertex-spaces $X_v$ are path-metric spaces.

 \medskip
 In view of the approximation lemmata (Lemma \ref{lem:graph-approximation} and Lemma \ref{lem:simplicial-approximation}), one can replace general path-metric spaces $X_v$ and  incidence maps $f_{ev}$ by (connected) metric graphs (equipped with graph-metrics) and simplicial  incidence maps. Below  
 we define the {\em total space} $X$ of a tree of spaces and a projection $\pi: X\to T$. Thus, we will frequently refer to trees of spaces as $\X= (\pi: X\to T)$, since the map $\pi$ records the most important information about $\X$.

In important class of trees of spaces consists of {\em metric bundles}. We refer to \cite{pranab-mahan} for the general definition; for the purpose of this book the following will suffice:

\begin{defn}\label{defn:bundle} 
An abstract tree of spaces $\X= (\pi: X\to T)$ is a {\em metric bundle} if the incidence maps $f_{ev}$ are uniform quasiisometries, i.e. there exists $\eps\ge 0$ such that for each edge $e=[v,w]\in E(T)$, the image $f_{ev}(X_e)$ (and, hence, $f_{ew}(X_w)$, by reversing the orientation on $e$) is $\eps$-dense in $X_v$.  
\end{defn}
 
 While the main motivation for trees of spaces comes from graphs of groups, the main group-theoretic 
 examples of metric bundles over trees are short exact sequences 
 $$
 1\to K\to G \to H\to 1,
 $$ 
 where $K$ is a finitely generated group and $H$ is a free group of finite rank.

 \begin{defn}\label{defn:total space} 
 The {\em total space}, or the {\em push-out}, of a tree of spaces $\X$ is a metric space $X$ admitting a collection of $L'$-coarse Lipschitz maps $X_e\to X$, $e\in E(T)$, $X_v\to X, v\in V(T)$, and satisfying the following universal property: For every metric space $Y$ and a compatible collection of $L_1$-coarse Lipchitz maps $X_e\to Y, X_v\to Y$, there exists a unique, up to a uniformly bounded error, $L_2$-coarse Lipschitz map $X\to Y$ forming diagrams which commute up to a uniform error $C$:  
 $$
\begin{diagram}
X_e & \rTo & X\\
& \rdTo & \dTo \\
       & \rdTo&  Y
\end{diagram}\quad \hbox{and}\quad 
\begin{diagram}
X_v & \rTo & X\\
& \rdTo & \dTo \\
       & \rdTo&  Y
\end{diagram} 
$$
Here $L_2$ and $C$ depend on $L_1$.  
 \end{defn}
 
This definition implies uniqueness (up to a quasiisometry) of the total space $X$. We will prove the existence of $X$ below (Theorem \ref{thm:existence-of-trees}). 
 
 \begin{defn}\label{defn:retractive tree} 
 An abstract tree of spaces is said to be {\em retractible} (or {\em retractive}), if there exists a collection of (uniformly) $L$-coarse Lipschitz maps (retractions) $f_{ve}: X_v\to X_e$ defined 
  for oriented edges $e=[v,w]$, which are uniformly coarse left-inverses to the incidence maps $f_{ev}$, i.e. 
  $$
  \dist(f_{ve}\circ f_{ev}, \id_{X_e})\le \eps,
  $$
  for some uniform constants $L\ge 1, \eps\in [0,\infty)$.  
 \end{defn}

Under the retractibility assumption, the incidence maps are not only uniformly proper but are also uniformly quasiisometric embeddings. While the definition is general, in this book, vertex and edge-spaces {\em mostly} will be uniformly hyperbolic, images of edge-spaces in vertex spaces 
will be uniformly quasiconvex and the retractions $f_{ve}$  will be given by nearest-point projections $P_{X_v, X_e}: X_v\to X_{ev}$.  

\medskip
{\bf Morphisms.} Let $\X, \X'$ be abstract trees of spaces over trees $T, T'$ respectively with the respective vertex/edge spaces $X_v, X'_{v'}, X_e, X'_{e'}$. 
A {\em morphism} of abstract trees of spaces $\X\to \X'$ is a graph-morphism 
$T\to T'$, $v\mapsto v', e\mapsto e'$,  together with a collection of 
\index{morphism of trees of spaces} 
uniformly coarse Lipschitz maps, 
between respective vertex and edge-spaces 
$$
h_v: X_v\to X'_{v'}, \quad h_e: X_e\to X'_{e'}
$$
such that the diagrams (where the horizontal arrows are the incidence maps) 
$$
 \begin{diagram}
 X_e              & \rTo & X_v           \\
 \dTo_{h_e}   &     & \dTo_{h_v} \\
  X'_{e'}              & \rTo & X'_{v'} 
  \end{diagram}
$$
commute up to uniformly bounded errors. An {\em isomorphism} of abstract trees of spaces is an invertible morphism, equivalently, it is an isomorphism of trees $T\to T'$ and a collection of uniform quasiisometries $X_v\to X'_{v'}, X_e\to X'_{e'}$. 

\begin{rem}
In this book we will be only considering {\em monic} morphisms of trees of spaces, i.e. ones for which the graph-morphism $T\to T'$ is injective 
and the maps $X_v\to X_{v'}$, $X_e\to X'_{e'}$ are $\zeta$-proper for some uniform function $\zeta$. 
\end{rem}

\begin{example}
The most common examples of morphisms of trees of spaces used in this book  are {\em subtrees of spaces}. Namely, let $S\subset T$ is a subtree, $\X$ is a trees of spaces over $T$. Then the {\em pull-back} of $\X$ over $S$ is a tree of spaces $\Y$ such that $Y_v=X_v, Y_e=X_e$, $v\in V(S), e\in E(S)$. The collection of identity maps 
$Y_v\to X_v, Y_e\to X_e$ defines a morphism of trees of spaces $\Y\to \X$. We will use the notation $X_S$ for the total space of the tree of spaces $\Y$. In the case when $S$ is an interval (resp. tripod) in $T$, we will refer to $X_S$ as an {\em interval-space} (resp. {\em tripod-space}). 
\end{example}

While the above definition is the main definition used in this book, we  now 
connect the  notion of an abstract trees of spaces to the notion of a tree of spaces as defined by Mitra in \cite{mitra-trees}. According to Mitra's definition, a 
tree of metric spaces is a path-metric space equipped with a certain auxiliary data, such as a map to a simplicial tree and a collection of maps to $X$ from certain spaces. 
Below, we use the $\ell_1$-metric on $X_e\times [v, w]$. Recall that for edges  $e$ of $T$,  
$\dot{e}$ denotes the edge minus its end-points; below we will use the notation $m(e)$ for  the midpoint of $e$.

\begin{defn}[Tree of spaces] \label{defn:tree-of-spaces} 
A {\em tree  of  metric spaces}, denoted ${\mathfrak X}$, is a path-metric space $(X, d)$ equipped with a $1$-Lipschitz surjective map  
$\pi : X \rightarrow T$ onto a simplicial tree $T$, satisfying the following:

\begin{enumerate}
\item For each $v\in V(T)$,  the corresponding {\em vertex-space} $X_v := \pi^{-1} (v)\subset X$ is rectifiably connected.

\item For every edge $e\in E(T)$, the {\em edge-space} $X_e:= \pi^{-1}(m(e))$  is rectifiably connected.  Every oriented edge $e=[v,w]$ 
comes equipped with an $L$-Lipschitz\footnote{The  Lipschitz condition is absent in \cite{mitra-trees}, but it holds in all natural examples. 
On the other hand, Mitra assumes that each restriction $f_{e} |_{ X_e \times\dot{e}}$ is  an isometry onto $\pi^{-1}(\dot{e})$, equipped with its path-metric induced from $X$. We find this assumption unnecessarily restrictive.}  
$\eta$-proper map 
$$
f_{e}: X_e\times [v,w]\to X, 
$$
such that $f_{ev}(X_e\times \{v\})\subset X_{v}$.
\end{enumerate}
\end{defn}

By abusing the notation, we will denote a tree of spaces by $\pi: X\to T$. We will use the notation $f_{ev}$ for the composition
$$
X_e\to X_e \times \{v\} \stackrel{f_{e}}{\longrightarrow} X_v
$$
and $X_{ev}:= f_{ev}(X_e)$. 

\begin{rem}
Mitra also assumes that  inclusion maps $X_v\to X$ to be $\zeta$-proper for some function $\zeta$. We will see below that this is a consequence of uniform properness 
of the maps $X_e\to X_v$. 
\end{rem}

\begin{comment}

\begin{defn}
We refer to the constants and functions appearing  in this definition as {\em parameters} of the tree of spaces $\pi: X\to T$. 
The space $X$ is the {\em total space} of the tree of spaces $\pi: X\to T$. 
\end{defn}

\begin{defn}
A {\em morphism} of trees of spaces $\pi: X\to T, \pi': X'\to T'$ consists of a pair of maps $\iota: T\to T'$, $h: X\to X'$, where $\iota$ is a simplicial embedding and $h$ 
is a coarse Lipschitz map $h: X\to X'$ such that $\pi'\circ h= \iota\circ \pi$. Since $\iota$ is a simplicial embedding, we will frequently identity $T$ with its image $\iota(T)$. 
Two morphisms $(h,\iota), (h',\iota')$ are coarsely inverse to each other if $\iota'=\iota^{-1}$ and $h'$ is coarsely inverse to $h$. 
Similarly, one defines coarse right/left inverse morphisms. Two trees of spaces $\pi: X\to T, \pi': X'\to T$ are (coarsely) isomorphic 
if there exist their morphisms which are  coarse inverse to 
each other.  In particular, such morphisms are quasiisometries between $X$ and $X'$. 
\end{defn}

{\mini 
Note also that for every subtree $S\subset T$ and a tree of spaces ${\mathfrak X}= (\pi: X\to T)$, the preimage $X_S:= \pi^{-1}(S)$ carries a natural structure of a tree of spaces 
${\mathcal S}= (X_S\to S)$ for which the inclusion maps $S\to T, X_S\to X$ define a morphism of trees of spaces. }

\end{comment}

\medskip 
Observe that each tree of spaces ${\mathfrak X}$ yields naturally an abstract tree of spaces ${\mathfrak X}^{ab}$, 
the {\em abstraction} of ${\mathfrak X}$, with the incidence maps $f_{ev}$. 
The next theorem is a converse to this {\em abstraction} procedure.  

\begin{theorem}
[An existence theorem for  trees of spaces] \label{thm:existence-of-trees}
For each abstract tree of spaces $\X$ over a tree $T$, there exists
a (unique up to an isomorphism) tree of spaces $(\pi: X\to T)$, called a {\em concretization} of ${\Y}$, such that 
${\X}$ is isomorphic to the abstraction ${\mathfrak X}^{ab}$ of  $(\pi: X\to T)$. The total space $X$ of $\X$ satisfies the universal property in the 
Definition \ref{defn:total space}. 
\end{theorem}
\proof Our proof mimics the definition of the underlying topological space (equipped with the weak topology) of a cell complex, where the latter is defined via an inductively defined collection of attaching maps. 
We let $X$ denote the topological space 
obtained by attaching the products 
$X_e\times [0,1]$ to the disjoint union 
$$
{\mathcal X}=\coprod_{v\in V(T)} X_v  
$$
via the attaching maps $f_{ev}: X_e\times \{v\}\to X_{v}, f_{ew}: X_e\times \{w\}\to X_w$, $e=[v,w]\in E(T)$. In other words, $X$ is the mapping cylinder 
$$  
Cyl(f: {\mathcal X}_E \to \XX_V)$$
of the map 
$$
f: {\mathcal X}_E:= \coprod_{e\in E(T)} X_e \to {\mathcal X}_V:= \coprod_{v\in V(T)} X_v,  
$$
given by the collection of  incidence maps $f_{ev}$. 
For each edge $e$ of $T$ we will identify $X_e\times \dot{e}$ with its image in $X$. 

We define {\em admissible} paths in $X$ (see Section \ref{sec:length structures}) to be  the  continuous maps $c: [a, b]\to X$ which are concatenations of {\em vertical} paths, which are rectifiable (with respect to the metrics on vertex-spaces) paths contained in the vertex-spaces of $X$ and {\em horizontal} paths, which are rectifiable paths contained in the intervals of the form $x\times [0,1]$, $x\in X_e, e\in E(T)$. For every admissible path $c$, we let $length(c)$ 
be the sum of measures of lengths of its vertical and horizontal components.  We leave it to the reader to verify that this defines a 
{\em length-structure} on $X$ and, hence, a path-metric $d$. {\em  We retopologize $X$ using this path-metric}. 
By the construction, each inclusion map $X_e\times \dot{e}\to X$ is an isometry to its image, each 
vertex space is rectifiably--connected in $X$, each inclusion map $X_v\to X$ is $1$-Lipschitz and the projection map 
$$
\pi: X\to T
$$
is  $1$-Lipschitz as well. The verification that the space $X$ satisfies the universal property is rather straightforward. Given a collection of compatible coarse $L$-Lipschitz maps $h_v:X_v\to Y, h_e: X_e\to Y$ to a metric space $Y$, we define a map $h: X\to Y$ by sending each open interval $\{x\}\times (0,1)\subset X_e\times (0,1)$ to the point 
$h_e(x)$.  The uniqueness of $h$ (up to a bounded error) follows from the fact that the union 
$$
\coprod_{v\in V(T)} X_v \sqcup \coprod_{e\in E(T)} X_e
$$ 
forms a $1/2$-net in $X$. We will leave it to the reader to check that $\X$ is isomorphic to the abstraction of $(\pi: X\to T)$. \qed 

\medskip

\begin{rem}
1. A definition similar  to our abstract tree of spaces and a construction analogous to the one in the proof of Theorem 
\ref{thm:existence-of-trees} appear in  the work of Cashen and Martin \cite[2.4]{Cashen-Martin}. However, they work in the category of proper metric spaces and the metric spaces they produce do not satisfy all the properties in Definition \ref{defn:tree-of-spaces} and, hence, we cannot directly use their work. 

2. Throughout the book, we will work with  geometric realizations of trees of spaces constructed in the proof of Theorem \ref{thm:existence-of-trees}. In particular, the inclusion maps
$$
X_v\to X
$$
are $1$-Lipschitz. For every edge $e=[u,v]\in E(T)$ we will be frequently using 
the path-metric spaces 
$$
X_{uv}=X_{\llbracket u, v\rrbracket} =\pi^{-1}(\llbracket u, v\rrbracket). 
$$
 The inclusion maps $X_u\to X_{uv}\leftarrow X_v$ are also $1$-Lipschitz.   
 
3. One drawback of our construction is that even if vertex and edge-spaces are complete and geodesic, the tree of spaces we construct in the proof of Theorem  \ref{thm:existence-of-trees} is a only a path-metric space and, a priori, is not a geodesic metric space and need not be complete. There are two ways to rectify this issue which we describe below.
\end{rem}

a. In the book, when we say ``a geodesic" we really mean a path which is $\eps$-short for 
a suitably chosen sufficiently small $\eps>0$. Similarly, when dealing with nearest-point projections, we frequently project to non-closed subsets. 
Then a nearest-point projection of $x\in X$ to $Y\subset X$ is a point $\bar{x}\in Y$ such that for a suitable chosen, sufficiently small $\eps>0$,
$$
d(x, \bar{x}) \le d(x, Y) +\eps. 
$$ 

b. For a reader uncomfortable with such a fudge, we describe an alternative approach to rectifying the issue with geodesics and nearest-point projections. 

First of all, as we noted earlier (Lemmata \ref{lem:graph-approximation}, \ref{lem:simplicial-approximation}) without loss of generality, we may assume that all vertex spaces and edge-spaces $X_v, X_e$ are connected graphs equipped with standard graph-metrics. We will replace each $X_e$ with its vertex-space. Then the space $X$ defined in 
the proof of Theorem \ref{thm:existence-of-trees}  is a connected graph and the path-metric on this graph defined in the proof is the standard graph-metric.  
 The drawback of this approach is the need to keep track of combinatorial issues which, are, ultimately, irrelevant.

\medskip 
From now on, we will work with abstract trees of spaces $\X$ and their concretizations $\pi: X\to T$. The metric space $X$ is the {\em total space} of $\X$. There is nothing particularly canonical about our choice of $X$ in this construction, it is just something we find convenient to work with. The reader could alternatively work for instance with, say, the $\ell_1$-metric coming from the products $X_e\times [v,w]$ in the mapping cylinder $X$. In fact, most of our arguments deal with vertex-spaces and pull-backs $X_{vw}$: We will be using the fact that the natural inclusion maps $X_v\to X_{vw}\leftarrow X_w$ are $1$-Lipschitz and either uniformly proper or, for trees of hyperbolic spaces, uniform qi embeddings. 

\begin{example}\label{ex:BS-tree} 
One motivation for our construction of $X$ comes from Cayley graphs of fundamental groups $G$ of  graphs  of groups. We assume that $({\mathcal G}, Y)$ is a finite 
graph of finitely generated groups, $S\subset Y$ is a spanning tree, and $G=\pi_1({\mathcal G}, Y, S)$ is the fundamental group. 
We will identify $S$ with a subtree in the Bass--Serre tree $T$ of  $({\mathcal G}, Y)$. 
 Then form a graph $\Gamma$ using the generators of $G$ as described in Definition \ref{defn:fundamental group of graph of groups}, except:
 
 (a) We fix an orientation of  the edges of $Y$ and use only one generator per each edge (not two). 
 
 (b) We use the given generating sets of the vertex-groups $G_v$ instead of the entire $G_v$. 
 
 Thus, in the graph $\Gamma$  there are vertical edges (corresponding to translates $\Gamma_v, v\in V(T)$, 
 of Cayley graphs of vertex groups) and horizontal edges (corresponding to the generators coming from the edges of $Y$). The vertex-spaces $X_v$ are, then 
 the graphs $\Gamma_v$. The edge-spaces are the translates of the edge-groups, $gG_e$, $g\in G$, $e\in E(Y)$. The incidence maps $f_{ev}, f_{ew}$ for the oriented 
 edges $e=[v,w]$ in $S$ come from the monomorphisms $\phi_{e,o(e)}$ and $\phi_{e,t(e)}$. For general edges $e$ of $T$ (which are translates of the edges $e'\in S$), 
 the incidence maps are obtained by composing with the action of $G$ by left multiplication. Thus, we obtain a tree of spaces $\X$ over $T$ with vertex spaces isometric to Cayley graphs of the vertex-groups $G_v$, $v\in V(Y)$, and edge-spaces isometric to edge-groups $G_e, e\in E(Y)$, with metrics obtained via pull-backs of word-metrics on 
 the incident  vertex-groups $G_v$, $v=t(e)$. Note that in the Cayley graph of $G$ there are no edges corresponding to generators of the edge-groups. This is consistent to our use of only horizontal paths over the edges of $T$ in the construction of the total space $X$ in the proof of Theorem \ref{thm:existence-of-trees}. We leave it to the reader to check that the Cayley graph $\Gamma$ as above is $G$-equivariantly isometric to the total space $X$ of the tree of spaces $\X$ defined in 
the proof of Theorem \ref{thm:existence-of-trees}.  
\end{example}

\begin{prop}\label{unif-emb-subtree}
There exists a continuous function $\eta_{\ref{unif-emb-subtree}}$ depending on the parameters of an abstract tree of spaces $\X$, such that for every subtree 
$S\subset T$, the inclusion map 
$$
X_S\to X
$$
is an $\eta_{\ref{unif-emb-subtree}}$-uniformly proper embedding. 
\end{prop}
\proof The key case to understand is when $T$ has a single edge $e=[u,v]$ and $S=\{u\}$. We let $Y$ denote the total space of the corresponding tree of spaces. 
It suffices to estimate (from below, in terms of $d_{X_u}(x,x')$) lengths of paths $c$ in $Y$ connecting $x=x_1, x'=x_n\in X_u$, such that $c$ is a concatenation of the form
$$
c(x_1, y_1) \star c(y_1, z_1) \star c(z_1,z_2) \star c(z_2, y_2) \star c(y_2, x_2) \star c(x_2, x_3) \star  ... \star c(y_n, z_n),
$$
where $x_i=f_{eu}(y_i), z_i= f_{ev}(y_i)$ and paths $c(x_i,y_i), c(y_i,z_i)$ are horizontal, while the paths $c(x_j, x_{j+1}), c(z_k, z_{k+1})$ are vertical geodesics in the 
vertex-spaces $X_u, X_v$. The lengths of this path is
$$
\length(c)= \sum_{i=\hbox{even}}  d_{X_u}(x_i, x_{i+1}) + n + \sum_{j=\hbox{odd}}  d_{X_v}(z_j, z_{j+1}). 
$$
Assume that $\length(c)\le D$. Then $n\le D$ and $d_{X_v}(z_j, z_{j+1})\le D$ for each odd index $j$. 
We have (for $j$ odd):
$$
L^{-1} d_{X_u}(x_j, x_{j+1}) \le d_{X_e}(y_j, y_{j+1})\le \psi(d_{X_v}(z_j, z_{j+1}))
$$
and, hence,
$$
d_{X_u}(x_j, x_{j+1})  \le L \psi( d_{X_u}(z_j, z_{j+1})) \le L \psi(D).  
$$
Thus, the concatenation $c_u$ of vertical geodesics $[x_i x_{i+1}]_{X_u}$ connecting $x$ to $x'$ has total length $\length(c_u)$ satisfying 
\begin{align*}
d_{X_u}(x,x')\le \length(c_u)= \sum_{j=\hbox{odd}}  d_{X_u}(x_j, x_{j+1}) + \sum_{i=\hbox{even}}  d_{X_u}(x_i, x_{i+1})  \le \\
L \sum_{j=\hbox{odd}}  \psi(D) +  \sum_{i=\hbox{even}}  d_{X_u}(x_i, x_{i+1}) \le LD \psi(D) + D. 
\end{align*} 
It follows that $d_{X_u}(x,x')\le LD \psi(D) + D$ and, hence, the inclusion map $X_u\to X_{uv}$ is $\eta$-proper, for $\eta(D):= D(L\psi(D) + 1)$. 

\begin{rem}\label{rem:linear} 
Assuming that the map $f_{ev}: X_e\to X_v$ is an $L$-qi embedding (which will be eventually our assumption for trees of hyperbolic spaces), 
we obtain a better estimate:  
\begin{align*}
d_{X_u}(x,x')\le \length(c_u)= \sum_{j=\hbox{odd}}  d_{X_u}(x_j, x_{j+1}) + \sum_{i=\hbox{even}}  d_{X_u}(x_i, x_{i+1})  \le \\
\sum_{i=\hbox{even}}  d_{X_u}(x_i, x_{i+1})  + 
L^2 \sum_{j=\hbox{odd}}  d_{X_v}(z_j, z_{j+1}) + L^3 n \le L^3 \length(c). 
\end{align*} 
Thus, we conclude that each inclusion map $X_u\to X_{uv}$ in this case is an $(L^3, 0)$-qi embedding. 
\end{rem}

 We now deal with the general case. Consider an admissible  path $\beta: [0,1]\to X$ connecting $x, y\in X_S$. 
 The projection $\pi\circ \beta$ is a path $p$ in $T$ connecting $\pi(x)$ to $\pi(y)$ whose length is $\le \length(\beta)$. Without loss of generality, we may assume that $\pi(x), \pi(y)$ are vertices in $S$ and $p$ is a simplicial path in $T$. We now construct inductively a sequence of paths 
$$
\beta_0= \beta, \beta_1,...,\beta_n
$$
in $X$ with simplicial projections to $T$, all connecting $x$ to $y$, such that:

(1) $\beta_n$ is a path in $X_S$. 

(2) The length of $\pi\circ \beta_{i+1}$ is at most 
$\length(\pi\circ \beta_{i}) -2$.

(3)   
$$
\length(\beta_{i+1}) \le \eta(\length(\beta_i)) 
$$
where $\eta(D):= D(L\psi(D) + 1)$ as above.

Assume that $\beta_{i}$ is defined. If this path is contained in $X_S$, then $n=i$ we are done. Otherwise, there exists an edge $e=[v,w]$ in the tree $\pi \beta_i([0,1])$ such that $\beta$ contains a subpath $\beta'$ connecting points $x', y'\in X_v$ and contained in the subspace $X_{vw}$. We then replace $\beta'$ with a geodesic in $X_v$ connecting 
the end-points $x', y'$ of $\beta'$. By the above estimate in $X_{vw}$, 
$$
d_{X_v}(x',y')\le \eta(\length(\beta')) 
$$
and, hence, the new path $\beta_{i+1}$ satisfies the required conditions. 

Clearly, $n\le \length(\beta)$, hence, 
$$
\length(\beta_n)\le \eta^{(n)}(\length(\beta)), 
$$
where $\eta^{(n)}$ is the $n$-fold iteration of the function $\eta$. Therefore, for $\eta_{\ref{unif-emb-subtree}}= \eta^{(n)}$, 
$$
n=\lceil  d_X(x,y) \rceil$$ 
we obtain 
$$
d_{X_S}(x,y)\le \eta_{\ref{unif-emb-subtree}} ( d_X(x,y)).  \qed 
$$ 

Applying the arguments of the proof of the proposition with the linear estimate in Remark \ref{rem:linear} we obtain: 

\begin{cor}\label{cor:exp-dist} 
If each incidence map $f_{ev}$ is an $L$-qi embedding, then each $X_S$ is at most exponentially distorted in $X$, i.e. is 
$\eta$-uniformly properly embedded in $X$ with $\eta(t)= \exp(at)$ for some $a\ge 1$ depending only on $L$. 
\end{cor}

We omit the proof of this corollary since it is straightforward and the result is not used elsewhere. 

\begin{defn}\label{defn:lift} \index{$K$-qi section} 
Let $\X=(\pi: X\to T)$ be a tree of spaces. 

\begin{enumerate}

\item By a $K$-{\em qi section} (or a $K$-{\em qi lift} of $S$) over a subtree $S\subset T$ we mean a map 
$\si: S\to X$ such that for each vertex $v\in S$, $\si(v)\in X_v$, for any pair of adjacent vertices $u, v\in S$, 
we have $d_{X_{uv}}(\sigma(u), \sigma(v))\leq K$ and the restriction of $\si$ to the interval $uv$ is a parameterization of a geodesic 
 $\si(u) \si(v)$ in $X_{uv}$. 

\item $K$-qi lifts of geodesic segments of $T$ will be referred to as {\em $K$-qi leaves} in $X$ 
and denoted by $\ga$ or $\ga_x$ or $\ga_{xy}$, provided they start at $x$ and end at $y$. We will refer to such $\gamma$'s as 
{\em horizontal} paths in $X$.

\item A {\em vertical path} in $X$ is a path contained in one of the vertex-spaces.

\item If $Y$ is a subset of $X$ then the {\em fiberwise neighborhood} of $Y$ in $X$ (denoted $N_r^{fib}(Y)$) is the union
$$
\bigcup_{v\in V(T)} N_r(Y\cap X_v),
$$
where the latter neighborhood is taken with respect to the (intrinsic) metric of $X_v$. 
\end{enumerate}
\end{defn}

\medskip

Let ${\mathfrak X}$ be an abstract  tree of spaces. A {\em subtree of spaces} in ${\mathfrak X}=(\pi: X\to T)$ is a tree of spaces ${\mathfrak X}'=(\pi': X'\to T')$ whose base tree is a subtree $T'\subset T$, and  vertex/edge spaces $X'_v, X'_e$ are rectifiably connected uniformly properly embedded subsets of  $X_v, X_e$ respectively, so that the incidence maps of $\X'$ are uniformly close to 
restrictions of incidence maps of $\X$.

\section{Coarse retractions}

In this section we  prove a general existence theorem of coarse Lipschitz left-inverses ({\em retractions}) for morphisms of trees of spaces.

Let $T'$ be a subtree of $T$ and let ${\mathfrak X}=( \pi: X\to T)$, ${\mathfrak X}'=(\pi': X'\to T')$ be trees of spaces. We say that a morphism $h: X\to X'$ of these trees of spaces is a {\em relative $K$-qi embedding} if for each $v\in V(T'), e\in E(T')$, the map $h_v: X'_v \to X_v, h_e: X'_e\to X_e$ is a $K$-qi embedding. 
Similarly, one can define a {\em relatively retractible} morphism of trees of spaces (a morphism which admits a relative $L$-coarse Lipschitz retraction) as a morphism $h$ such that for each $v\in V(T'), e\in E(T')$ 
the maps $h_v: X'_v \to X_v, h_e: X'_e\to X_e$ admit $L$-coarse Lipschitz left-inverses $h'_v: X_v\to X'_v$,
$h_e: X_e\to X'_e$. If ${\mathfrak X}, {\mathfrak X}'$ are trees of $\delta$-hyperbolic spaces then the two notions are equivalent and, moreover, the subspaces $h_v(X'_v)\subset X_v, h_e(X'_e)\subset X_e$ are  $\la$-quasiconvex for $\la=\la(L,\delta)$. 
Our goal is to prove that, under some conditions, a relatively retractive morphism is {\em absolutely retractive}, i.e. admits a coarse left-inverse $h': X\to X'$. 
(Recall that the morphism $h'$ is a collection of maps $h'_v: X_v\to X'_v, h'_e: X_e\to X'_e$ satisfying certain compatibility properties.) 
 This result is motivated by Mitra's construction of a coarse retraction in \cite[Theorem 3.8]{mitra-trees}. For relatively retractible morphisms of trees, by abusing the notation, we will identify the vertex/edge spaces $X'_v, X'_e$ of ${\mathfrak X}'$ with their images $h_v(X'_v)\subset X_v$ and $h_e(X'_e)\subset X_e$ respectively.

 \medskip 
 The following theorem is inspired by Mitra's coarse retraction theorem in 
 \cite[Theorem 3.8]{mitra-trees} and its proof closely follows Mitra's argument.

\begin{theorem}[Existence of a retraction]\label{thm:left-inverse}
Suppose that for some constants $C, D$, a relatively retractive morphism of trees of spaces $h: {\mathfrak X}'\to {\mathfrak X}$ satisfies the following  conditions:

(i) For every boundary edge $e$ of $T'$, $e=[v,w], v\in V(T'), w\in V(T)-V(T')$,
$$
\diam_{X'_v} (h'_v \circ f_{ev}(X_e))\le D. 
$$

(ii) For every  edge $[v,w]=e\in E(T')$ 
$$
\dist_{X'_v}( h'_v \circ f_{ev}, f'_{ev}\circ h'_e)\le C. 
$$
Then the map $h: X'\to X$ admits a coarse $L_{\ref{thm:left-inverse}}$-Lipschitz retraction $h': X\to X'$ whose restriction to $X_v$ 
equals $h'_v$ for each $v\in V(T')$.  Here $L_{\ref{thm:left-inverse}}$ depends only on $C, D$, coarse Lipschitz constants of the maps $h'_v, h'_e$, and 
the parameters of trees of spaces 
$\X, \X'$. 
\end{theorem} 
\proof  We let $K$ denote the maximum of Lipschitz constants of the projections $\pi: X\to T, \pi': X'\to T'$. 
For each  $v\in V(T')$ then we let $h'(x):= h'_v(x)$. Let $p: T\to T'$ denote the nearest-point projection.  

Suppose $x\in X_w$, $w\in V(T)\setminus V(T')$; then $v=p(w)\in T'$ is the vertex nearest to $w$.  
Let $e\in E(T)$ be the edge incident to $v$ and contained in the geodesic $wv$. Thus, $e$ is a boundary edge of the subtree $T'\subset T$. 
By the assumption (i), 
the projection $h'_v(X_{ev})\subset X'_v$ has the diameter $\le D$. We let $h'(x)$ be any point $x'$ 
of this projection (we will use the same point $x'$ for all vertices $w$ in each component of $T-T'$). 

In order to verify that $h'$ is (uniformly) coarse Lipschitz it suffices to find a uniform upper bound on distances $d(h'(x), h'(y))$ for points 
$$
x, y\in {\mathcal X}= \coprod_{v\in V(T)} X_v
$$
which are within  distance $K$ from each other. If $x, y$ belong to the same vertex space $X_v$, then $d(h'(x), h'(y))\le L$, the upper bound for coarse Lipschitz constants of the maps $h'_v:  X_v\to X'_v$. Suppose that $x, y$ belong to $X_v, X_w$, $v, w\in V(T')$ are vertices spanning an edge $e\in E(T')$. Then, necessarily, 
$$
x\in X_{ev}, y\in X_{ew}, x= f_{ev}(z), y=f_{ew}(z)
$$ 
for some $z\in X_e$. The condition (ii) then implies the estimates
$$
d(h'_v(x), f'_{ev}\circ h'_e(z))\le C, \quad d(h'_w(y), f'_{ew}\circ h'_e(z))\le C,
$$
hence $d(h'(x), h'(y))\le 2C$. 

If $v, w\in V(T)- V(T')$ then the inequality $d_T(v,w)\le 1$ implies that $p(v)=p(w)=u\in V(T')$ and there is a common boundary edge $e$ of $T'$ contained in the geodesics $uv, uw\subset T$. In particular, 
both $h'(x), h'(y)$ belong to the subset 
$$
h'_u \circ f_{eu}(X_e)\subset X'_u
$$  
and, hence, $d(h'(x), h'(y))\le D$ by the condition (i). 

Lastly, consider the case when $x\in X_v, y\in X_w$, where $v\in V(T'), w\in V(T) - V(T')$ and $v, w$ span a boundary edge $e$ of 
$T'$. Since $d(x,y)\le 1$, it follows that $x\in X_{ev}, y\in X_{ew}$. Then $p(w)=v$ and, by the definition of $h'$, 
$h'(x), h'(y)\in h_v'(X'_{ev})$ and, therefore, 
$$
d_{X'_v}(h'(x), h'(y))\le D. \qed 
$$

An easy corollary of Theorem \ref{thm:left-inverse} is:  

\begin{cor}\label{cor:r'}
Suppose that ${\mathfrak X}=(\pi: X\to T)$ is a retractive tree of spaces. 
For every edge $e=[u,v]\in E(T)$ there exists an $r=r_{\ref{cor:r'}}$-coarse retraction 
$X_{uv}\to X_u$, where $r$ depends only on the parameters of ${\mathfrak X}$ and its {retractivity constant}. 
\end{cor}
\proof We have a retractive tree of spaces ${\mathfrak Y}= (\pi: X_{uv}\to \llbracket u,v\rrbracket )$. In ${\mathfrak Y}$ we have a subtree of spaces $\pi': {\mathfrak Y}'= (Y'\to \llbracket u,v\rrbracket )$, whose vertex spaces are 
$Y'_u=Y_u=X_u, Y'_e= Y_e=X_e, Y'_v= X_{e}, f'_{ev}=\id, f'_{eu}=f_{eu}$, and the morphism $h: {\mathfrak Y}'\to {\mathfrak Y}$ is defined by using $\id: Y'_u\to Y_u$ and  
$f_{ev}: Y'_v\to Y_v$.  Since ${\mathfrak X}$ is retractive, the morphism $h$ is relatively retractive.  Hence, by Theorem \ref{thm:left-inverse}, 
the identity map $X_u\to X_u$ and the retraction $X_v\to X_e$ define a coarse Lipschitz retraction $Y_{uv}\to Y'$. Since $Y'$ is Hausdorff-close to $X_u$, we obtain a coarse   
Lipschitz retraction $Y_{uv}\to X_u$. The reader will verify that the coarse Lipschitz bound for this retraction depends only on the parameters of ${\mathfrak X}$ and its retractivity constant. \qed 

\medskip 
Another useful application of Theorem \ref{thm:left-inverse} is in the setting of trees of hyperbolic spaces (which we will discuss in more detail in the next section): 

\begin{cor}\label{cor:projection}
Suppose that the trees of spaces ${\mathfrak X}, {\mathfrak X}'$ and a morphism ${\mathfrak X}'\to {\mathfrak X}$ which is a fiberwise $L$-qi embedding, have the following properties:

1. For some $\delta$, all vertex and edge-spaces $X_v, X_e$ are $\delta$-hyperbolic. 
(Accordingly, the images vertex and edge-spaces $h_v(X'_v)\subset X_v$,  $h'_e(X'_e) \subset X_e$ 
are $\la$-quasiconvex subsets in $X_v, X_e$ respectively, where $\la_{\ref{lem:qi-preserves}}(\delta, L)$.) 

2. The retractions $h'_v, h'_e$ are ``nearest-point projections'' in the sense that 
$$
h'_v= P_{X_v,X'_v} \circ h_v, h'_e= P_{X_e,X'_e} \circ h_e. 
$$ 

3. There is constant $K$, for every edge $e=[v,w]\in T'$ 
the Hausdorff distances  $\Hd_{X'}(X'_v, X'_e)$ and $\Hd_{X'}(X'_w, X'_e)$ are $\le K$. 

4. $T'=T$. 

Then the fiberwise nearest point projections $h'_v, h'_e$ extend to an
$L_{\ref{cor:projection}}(\delta,L,K)$-coarse retraction $h': X\to X'$, where (without loss of generality) 
$$
L_{\ref{cor:projection}}(\delta,L,K) \ge \max(L,K). 
$$ 
\end{cor}
\proof For  vertices $v$ in $T$ incident to an edge $e$, the images $h_v(X'_v)$, $f_{ev}\circ h_e(X'_e)$ are 
uniformly Hausdorff-close to each other. Therefore, the nearest-point projections (in $X_v$) to these uniformly quasiconvex subsets 
are also uniformly close to each other (see Corollary \ref{cor:proj-to-close-subsets}). 
Now, the claim follows from Theorem \ref{thm:left-inverse}. 
 \qed

\section{Trees of hyperbolic spaces} \label{sec:hyperbolic trees}

We now introduce {\em hyperbolicity conditions} for trees of spaces. 

\begin{defn}\index{Axiom H} \label{defn:axiom H}
A  tree of spaces ${\mathfrak X}$ satisfies Axiom {\bf H} if there are constants $\delta_0$ and $L_0$ such that: 

(1) Each vertex/edge space $X_v, X_e$ of 
${\mathfrak X}$ is a $\delta_0$-hyperbolic geodesic metric space. 

(2) Each incidence map $f_{ev}: X_e\to X_v$ is an $L_0$-qi embedding.

\noindent We will refer to such $\X$ as a {\em  tree of hyperbolic spaces}. 

A finite graph of finitely generated groups ${\mathcal G}$ satisfies Axiom {\bf H} if the the corresponding tree of spaces does. In other words, all vertex and edge-groups have to be hyperbolic and edge-groups are quasiconvex in the incident vertex-groups. 
\end{defn}

A word of caution: Our terminology does not mean that a  tree of hyperbolic 
spaces  ${\mathfrak X}=(\pi: X\to T)$ 
has $\delta$-hyperbolic total space $X$. 
Simple examples are given by Euclidean plane and Cayley 
complexes of Baumslag--Solitar groups. 
One needs to add a suitable {\em flaring condition} on $\X$ to ensure hyperbolicity of $X$, as discussed in Section \ref{sec:flare}.  
Note also that our terminology requires not only uniform hyperbolicity of vertex and edge-spaces but also uniform qi embedding condition 
for the incidence maps.

\begin{defn}\index{parameters of a  tree of hyperbolic spaces} 
We will refer to $\delta_0$ and $L_0$ as the {\em primary parameters} of a  
tree of hyperbolic spaces  ${\mathfrak X}$. 
\end{defn}

In general, throughout the book, we will suppress the dependence of various constants and functions on the parameters of  
${\mathfrak X}$.

\begin{defn}
Suppose that $\X$ is a tree of hyperbolic spaces. Let $A, B\subset X$ with $\pi(A)\subset \pi(B)$. If $X_v\cap A$ and $X_v\cap B$ are uniformly quasiconvex in $X_v$ for all $v\in \pi(A)$, we  define the nearest projections in $X_v$ of 
$A\cap X_v$ to $B\cap X_v$. This gives us a map $A\map B$. We refer to this map as the {\em fiberwise projection} of $A$ to $B$. 
\end{defn}

\medskip 
It is immediate that for every tree of hyperbolic spaces, for every edge $e=[v,w]\in E(T)$, the subset $X_{ev}\subset X_v$ 
is $\la_0= \la_{\ref{lem:qi-preserves2}}(\delta_0,L_0)$-quasiconvex. In particular,  
every  tree of hyperbolic spaces  is retractive (see Definition \ref{defn:retractive tree}) with retractions 
$$
f_{ve}: X_v\to X_e, e=[v,w], 
$$
given by the nearest-point projections $P=P_{X_v,X_{ev}}$ 
to the quasiconvex subsets $X_{ev}=f_{ev}(X_e)\subset X_v$; more precisely:
$f_{ve}(x)$ is defined to be an arbitrary point in $f_{ev}^{-1}(P(x))$.

\medskip 
 As an application of Remark \ref{rem:linear} or, alternatively, of Corollary \ref{cor:r'}, we obtain: 

\begin{lemma}\label{lem:L'0}
Suppose that ${\mathfrak X}$ is a  tree of hyperbolic spaces with the primary parameters $\delta$ and $L$. Then for every edge $e=[u,v]\in E(T)$, 
the inclusion maps $X_u\to X_{uv}, X_v\to X_{uv}$ are $L'_0= L_{\ref{lem:L'0}}(\delta,L)$-qi embeddings where $L'_0$ is the maximum of $2$ and of the 
coarse Lipschitz constant for a retraction $X_{uv}\to X_v$ (see Corollary \ref{cor:r'}).    
\end{lemma}

\begin{rem}
In this lemma we ensured that $L'_0\ge 2$. This, somewhat artificial, convention will be used in the proof of  Lemma \ref{lem:growth-of-flare} below.  
\end{rem}

\medskip

Suppose that $\X'=(\pi: X'\to T)$ is a tree of hyperbolic spaces, $G< \Isom(X')$ is a subgroup acting by automorphisms of $\X$, such that the quotient graph $T/G$ is finite and for every vertex $v\in V(T)$ (resp. edge $e\in E(T)$) the action of the corresponding stabilizer $G_v< G$ (resp. $G_e< G$)  on $X'_v$ (resp. $X_e$) is quasiconvex (see Definition \ref{defn:qc-subgroup}).   Thus, the group $G$ also has structure of a graph of finitely generated groups ${\mathcal G}$ (with the underlying graph $T/G$); in particular, $G$ is finitely generated. (Note that we are not assuming hyperbolicity of the space $X$.) 

Since $G$ acts via automorphisms of $\X'$, for each edge $e=[v,w]\in E(T)$ the subspace 
$X'_{ev}\subset X'_v$ is $G_e$-invariant. We will also assume that for each $v\in V(T)$ the  $G_v$-orbit of $X'_{ev}$ is {\em locally finite} in $X'_v$ (see Definition \ref{defn:locally finite action}). Note that the local finiteness assumption is automatic for instance if there exists a larger discrete group $G'_v$ (containing $G_v$) acting  on $X'_v$  geometrically (Lemma \ref{lem:geo->lf}). 

\begin{prop}\label{prop:retract-so-subgroup} 
Under the above assumptions, there exists a coarse Lipschitz retraction $X'\to G x$ for each $G$-orbit in $X'$. 
In particular, each orbit map $G\to Gx\subset X'$ is a  qi embedding. 
\end{prop}
\proof The proof is similar to that of Corollary \ref{cor:projection}. Let $\X= (\pi: X\to T)$ denote the tree of hyperbolic spaces corresponding to the graph of groups ${\mathcal G}$. The isometric action of $G$ via automorphisms of $\X'$ defines a morphism of trees of spaces $\X\to \X'$. This morphism is relatively retractive in view of the quasiconvexity assumption for the actions  $G_v\acts X'_v, G_e\acts X'_e$. In view of Proposition \ref{prop:proj-to-qc-action}(4), the local finiteness assumption 
implies that for $y\in X'_v$, the restriction to $X'_e$ of the nearest-point 
projection $P_{X'_v,G_vy}$ is within uniformly bounded distance from the projection $P_{X'_{ev},G_e y}$. Thus, 
Theorem \ref{thm:left-inverse} applies and the coarse Lipschitz retractions $X'_v\to X_v, X'_e\to X_e$ 
 together give rise to a coarse Lipschitz retraction $X'\to X$. Since $G$ acts cocompactly on $X$, we, thus, obtain a coarse Lipschitz retraction $X'\to G x$. \qed

\begin{cor}
Suppose that ${\mathcal G}'$ is a finite, connected graph of hyperbolic groups satisfying Axiom {\bf H}, with 
$\pi_1({\mathcal G}')=G'$ and let $T$ be the Bass--Serre tree of ${\mathcal G}'$.  Let  $G< G'$ be a subgroup such that:

1. For every  vertex $v$ (resp. edge $e$) of $T$, the $G$-stabilizer $G_v< G$ of $v$ (resp.  
the $G$-stabilizer $G_e< G$ of $e$) is a quasiconvex subgroup of the $G'$-stabilizer $G'_v< G'$ of $v$ 
(resp. of the $G'$-stabilizer $G'_e< G'$ of $e$). 

2. The quotient-graph $T/G$ is   finite. 

There exists a coarse Lipschitz retraction $G'\to G$. In particular, the subgroup $G$ is qi embedded in $G'$. 
\end{cor}

\begin{example}
Let $H=\pi_1(S_1)\star \pi_1(S_2)$, where $S_1, S_2$ are closed connected hyperbolic surfaces, and let $\phi_i:  \pi_1(S_i)\to \pi_1(S_i), i=1, 2$, be automorphisms. Then $\phi_1, \phi_2$ define an automorphism $\phi: H\to H$ and 
we obtain subgroups  $G_i=\phi_1(S_i) \rtimes_{\phi_i} \Z$ in $G'=H\rtimes_\phi \Z$. The subgroups $G_i< G'$ clearly satisfy the assumptions of the corollary  
(where $T$ is the line) which implies that they are  coarse Lipschitz retracts of $G'$. Note, furthermore, that if $\phi_1, \phi_2$ are induced by 
pseudo-Anosov homeomorphisms of the surfaces $S_1, S_2$, then the group $G'$ is isomorphic to the amalgam of hyperbolic groups 
$G_1\star_{\Z} G_2$, where $\Z$ is a malnormal subgroup of both $G_1, G_2$. Hence, the group $G'$ is hyperbolic (see Corollary \ref{cor:malnormal-amalgam} below) and
 the subgroups $G_1, G_2$ are quasiconvex in $G'$.  
\end{example}

Below, $H$ is a quasiconvex subgroup of $G'_v$ for some vertex $v$ in $T$. 

\begin{lemma}\label{lem:qc-stabs} 
For each edge $e$ and vertex $w$ in $T'$, the $H$-stabilizer of $e$ (resp. $w$) is a quasiconvex subgroup of $H$ and $G'_e$ (resp. $G'_w$). 
\end{lemma}
\proof Consider first an edge $e=[v,w]$. Then $H_e=H\cap G'_e$ is the intersection of two quasiconvex subgroups of $G'_v$, hence, is quasiconvex in $G'_v$, $G'_e$ and $H$  
(see Corollary \ref{cor:qc-in}). The general case follows from induction on the edge-path connecting $e$ (resp. $w$) to $v$. \qed

\begin{lem}\label{lem:finite stabs} 
Suppose, in addition, that the $H$-stabilizers of edges incident to $v$ are all finite. 
Then for each $R\ge 0, x\in X'_v$ and all 
the edges $e$ incident to $v$, the coarse intersections intersections $Hx \cap N_R(X'_{ev})$ are uniformly bounded, with bound independent of $e$.  In other words, the 
pairs $Hx, X'_e$ are uniformly cobounded. 
\end{lem}
\proof By properness of the action, there exists an edge $e$  incident to $v$ such that the diameter of the intersection $Hx \cap N_R(X'_{ev})$ is maximal. The subset 
$X'_{ev}$ is Hausdorff-close to the orbit $G'_e x$. According to Proposition \ref{prop:proj-to-qc-action}, the coarse intersection $Hx \cap N_R(X'_{ev})$ is 
Hausdorff-close to the orbit $H_e x$, where $H_e=H\cap G'_e$. 
Since, by the hypothesis of the lemma, the subgroup $H_e$ is finite, the coarse intersection $Hx \cap N_R(X'_{ev})$  is bounded. \qed 

\begin{cor}\label{cor:finite stabs}
If the pair  $Hx, X'_e$ is not cobounded then the intersection $H\cap G'_e$ is infinite.  
\end{cor}

We will prove in Corollary \ref{cor:finite-tree-hyp} that there exists a function $\delta(n)$ (depending also on the constants $\delta_0$ and $L_0$)  such that 
for each interval $J$ of  length $n$, the pull-back space $X'_J$ (with its intrinsic path-metric) is $\delta(n)$-hyperbolic. Thus, we can talk about 
cobounded pairs of subspaces in vertex-spaces $X'_v, X'_w$, $v, w\in V(J)$.

\medskip 
Continuing with the notation of Lemma \ref{lem:finite stabs}, and applying Corollary \ref{cor:finite stabs} inductively (with Lemma \ref{lem:qc-stabs}),  
we obtain:

\begin{lemma}\label{lem:trivial stabs} 
Suppose that for some vertex $w\in T, x\in X'_v$, the subsets $Hx, X'_w$ are not cobounded in $X'_J$, $J= \llbracket v, w\rrbracket$. 
Then the $H$-stabilizer of the segment $J$ is an infinite subgroup of $H$. 
\end{lemma}

\medskip 
Below are few more easy consequences of Axiom {\bf H} for trees of spaces. 

\begin{lemma}\label{lem:edge-spaces}
Assume that ${\mathfrak X}$ is a  tree of hyperbolic spaces. 
Then 
for every edge $e=[v_1,v_2]$ of $T$, if $\al_{i}= [x_i y_i]_{X_{v_i}} \subset X_{v_i}$ are vertical geodesics such that 
$$
d_{X_{v_1 v_2}}(x_1, x_2)\le C, \quad d_{X_{v_1 v_2}}(y_1, y_2)\le C,
$$
then the Hausdorff distance between these vertical geodesics in $X_{v_1 v_2}$ is at most $C_1=C_{\ref{lem:edge-spaces}}(C)$. 
\end{lemma} 
\proof Geodesics $x_1x_2, y_1y_2$ have to cross $X_e$ (separating $X_{vw}$) at some points  
$x, y\in X_e$. Since both $X_{ev_i}\subset X_{v_i}$ are $\la_0$-quasiconvex, it follows that geodesics $\al_i$ lie in 
$N_{\la_0}(X_{ev_i})$, $i=1,2$. Lemma \ref{lem:sub-close} applied to the geodesic $\al_i$ and the $L_0$-quasigeodesic 
$\al'_i= f_{ev_i}([xy]_{X_e})$  implies that 
$$
\Hd_{X_{v_i}}(\al_i, f_{ev_i}(\al))\le  D=D_{\ref{lem:sub-close}}(\delta_0,L_0,C). 
$$ 
Since 
$$
\Hd_{X_{v_1v_2}}(\al, f_{ev_i}(\al))\le 1,
$$
we conclude:
$$
\Hd_{X_{v_1v_2}}(\al_1, \al_2)\le 2(1+ D). \qed 
$$

\begin{lemma}\label{lem:growth-of-flare} 
Let $I=\llbracket v,w\rrbracket \subset T$ be a subinterval, we denote its consecutive vertices $v_0=v, v_1, ..., v_n=w$. 
Let $\gamma_0, \gamma_1$ be $K$-qi sections over  $I$. Then the function
$$
\ell(i):= d_{X_i}(\ga_{0,k}(v_i), \ga_{1,k}(v_i)), i\in [0, n]\cap \Z,
$$
satisfies
$$
\ell(n)\le a^n  \ell(0) + \frac{a^n-1}{a-1}b< a^n  (\ell(0) + b), 
$$
where $a= L'_0$, $b=2L_0' K$. 
\end{lemma}
\proof Consider an edge $e=[v_i, v_{i+1}]\subset I$. The points 
$$
\ga_{0}(v_{i+1}), \ga_{1}(v_{i+1})\in X_{v_i} 
$$
are connected by a path of length $\le 2K+ \ell(i)$ in $X_{\llbracket v_i, v_{i+1} \rrbracket}$, 
obtained by concatenating a vertical geodesic 
$[\ga_{0}(v_i)\ga_{1}(v_i)]_{X_{v_i}}$ with two geodesics of length $\le K$. Since $X_{v_{i+1}}$ is $L'_0$-qi embedded in 
$X_{\llbracket v_i, v_{i+1} \rrbracket}$, we have
$$
\ell(i+1)\le L'_0( 2K+ \ell(i)). 
$$
Then
$$
\ell(n)\le a^n \ell(0) + (a^{n-1}+\ldots +1) b = a^n  \ell(0) + \frac{a^n-1}{a-1}b< a^n ( \ell(0) + b). \qed 
$$

\begin{cor}\label{cor:contraction} 
If $\ell(0)\ge a (M+b)$, then for all 
$$
n\in [0, N], N= \left\lfloor  \log_a  \left(   \frac{\ell(0) }{M +b}  \right) \right\rfloor, 
$$
we have
$$
\ell(n)>M. 
$$
\end{cor}
\proof We first reverse the role of $\ell(0)$ and $\ell(n)$ and obtain from the lemma that 
$$
\ell(n)<  a^{-n} \ell(0) - b, n\in \NN. 
$$
The inequality $\ell(n)> M$ then follows from 
$$
n\le N \le \log_a  \left(   \frac{\ell(0) }{M +b}  \right).
$$
The assumption that $\ell(0)\ge a (M+b)$ ensures that $N\ge 1$. \qed 

\medskip
Another corollary (or, rather, a special case of the lemma) is

\begin{cor}\label{bdd-flaring}
For every edge $e=[u,v]$ in $T$, any pair points $x, y\in X_u$, and a pair of $K$-qi sections $\ga_0, \ga_1$ over the interval $uv$, we have 
$$
d_{X_v}(\ga_0(v), \ga_1(v))\le D_{\ref{bdd-flaring}}(K, d_{X_u}(x,y))= L'_0( 2K+ d_{X_u}(x,y)). 
$$
\end{cor}

\section{Flaring}\label{sec:flare}

Geodesics (and, hence, uniform quasigeodesics) in hyperbolic spaces diverge (exponentially fast). Since $k$-qi leaves in hyperbolic trees of spaces $\pi: X\to T$ are uniform quasigeodesics,  they should also diverge if $X$ is hyperbolic. In this section we discuss several divergence 
conditions, called {\em flaring conditions}\footnote{Flaring conditions do not require Axiom {\bf H}.}, one can impose on qi-leaves  
in trees of spaces. These conditions involve pairs $\Pi=(\ga_0, \ga_1)$ of $k$-sections $\gamma_0, \gamma_1$ over a common  geodesic segment $J=\llbracket t_{-n}, t_n\rrbracket \subset T$ of length $2n$ 
and prescribe the nature of growth of the vertical distances 
$$
d_{X_{v_i}}(\gamma_0(v_i), \gamma_1(v_i))
$$
for $i>0$ or $i<0$. The {\em girth} $\Pi_0$ of the pair $(\gamma_0,\gamma_1)$ is the vertical distance 
$$
d_{X_0}(\gamma_0(0), \gamma_1(0)). 
$$

\begin{rem}
$\Pi_0$ need not be equal to 
$$
\min_{v\in V(J)} d_{X_v}(\ga_0(v), \ga_1(v)). 
$$
\end{rem}

We will  frequently use the notation $\Pi_{\max}$ for the {\em maximal separation of the ends}  of the pair $\Pi=(\ga_0, \ga_1)$, 
$$
\Pi_{\max}:= \max \left(d_{X_{t_{-n}}}(\ga_0(t_{-n}), \ga_1(t_{-n})), d_{X_{t_{n}}}(\ga_0(t_{n}), \ga_1(t_{n})) \right), 
$$
describing the rate of growth of the above vertical distances (in one of the directions).

\begin{figure}[tbh]
\centering
\includegraphics[width=60mm]{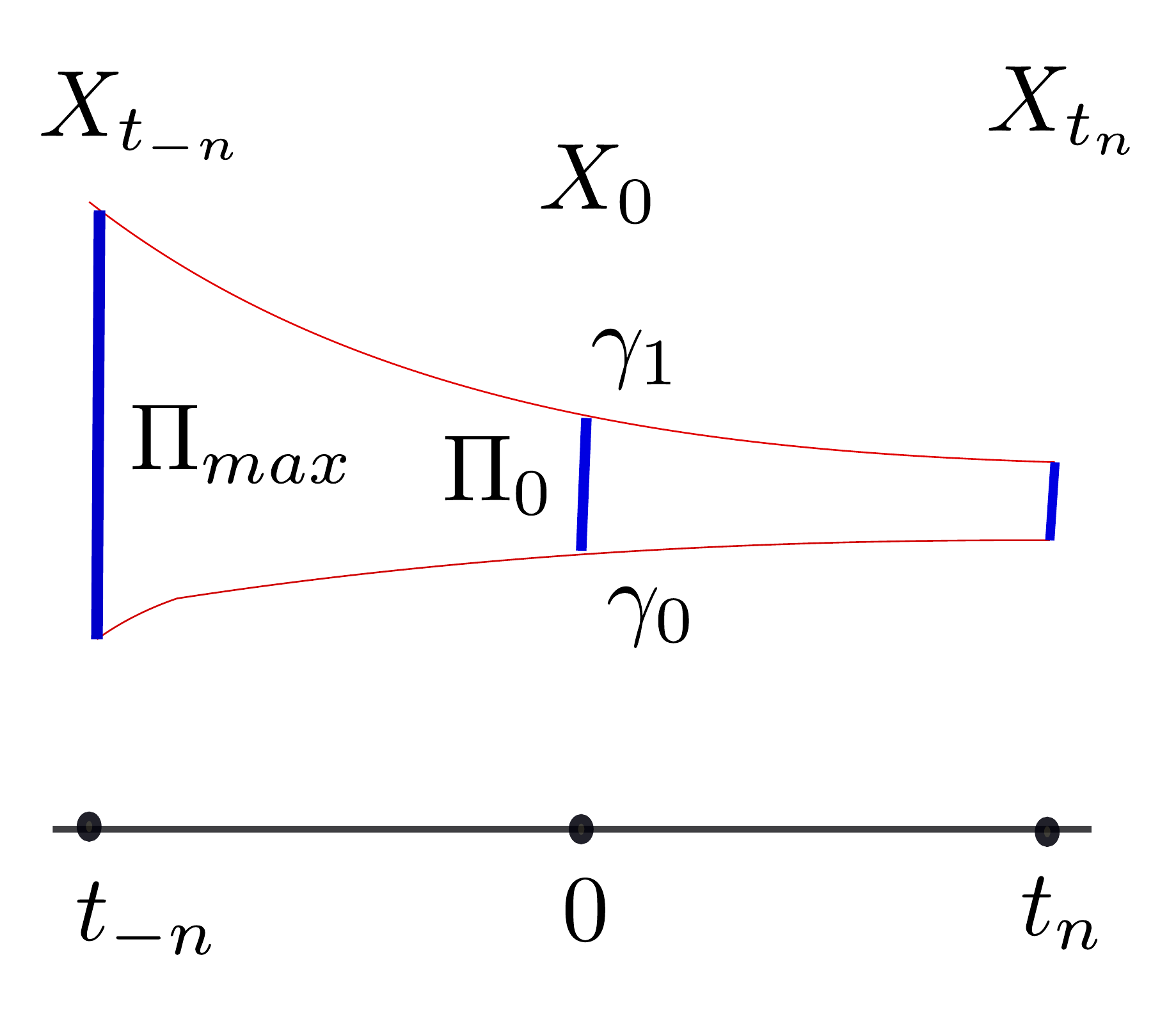} 
\caption{Flaring.}
\label{flaring.fig}
\end{figure}

\subsection{Proper  and uniform flaring conditions}  

The {\em proper flaring} condition requires $k$-qi sections 
over the same geodesic in $T$ to diverge at 
some uniform rate in at least one direction.  More precisely:

\begin{defn}
[Proper flaring] \index{proper flaring}
A tree of spaces ${\mathfrak X}=(\pi: X\to T)$ is said to satisfy the {\em proper $\kappa$-flaring condition} if there exists $m_\kappa\ge 0$ and a positive proper function $\phi_{\kappa}: {\mathbb N}\to {\mathbb R}_+$ such that for every  pair $\Pi$ of $\kappa$-qi sections $\gamma_0,\gamma_1$ of girth $> m_\kappa$, over an interval of length $2n$ in $T$,  
we have
$$
\Pi_{\max}\ge \phi_\kappa(n). 
$$
\end{defn}

In other words, $\kappa$-sections have to diverge uniformly fast but the rate of divergence is allowed to be, say, sublinear (unlike the the Bestvina--Feighn flaring condition where one has  an exponential rate of flaring).

\medskip 
It is clear from the definition that if ${\mathfrak X}$ satisfies the proper $K$-flaring condition, then it also satisfies the $\kappa$-flaring condition for all $\kappa\in [1, K]$: We simply take $\phi_\kappa:=\phi_K$ and $m_\kappa:=m_K$.

Note also that it would be too much to ask for 
$$
\Pi_{\min}=\min\left( d_{X_{-m}}(\ga_0(-m), \ga_1(-m)),  d_{X_{n}}(\ga_0(n), \ga_1(n)) \right) \ge \phi_\kappa(n),  
$$
for some (uniform) proper function $\phi_\kappa$.

\begin{defn}\label{defn:signed flare} 
We will say that a pair $\Pi=(\ga_0, \ga_1)$ of sections over an interval $\llbracket -N, N\rrbracket$ 
 in $T$ is {\em flaring in the positive/negative direction} if, respectively,
$$
d_{X_{n}}(\ga_0(n), \ga_1(n)) \ge \phi_\kappa(n),
$$
or 
$$
d_{X_{-n}}(\ga_0(-n), \ga_1(-n)) \ge \phi_\kappa(n) 
$$
for all $n\in \NN\cap [1, N]$. 
\end{defn}

We will see in Lemma \ref{lem:exp} and Corollary \ref{cor:signed flare}, that proper flaring (for all $\kappa\ge 1$) implies proper flaring in positive or negative direction 
(after a possible change of the function $\phi_\kappa$).

\medskip 
In the book we will be mostly using an alternative form of the proper flaring condition established in the next proposition. 
For the ease of the notation, in this section we will identify a geodesic, say $\llbracket v,w\rrbracket$, of length $\ell$ in $T$, 
where $v,w\in V(T)$ with an interval $[a, b]\subset \RR$ of length $\ell$, where $a, b\in \mathbb Z$ through an implicit isometry 
$[a,b]\map \llbracket v,w\rrbracket$; this means, in particular, that integers correspond to the vertices in $\llbracket v, w\rrbracket$.

\begin{prop}\label{prop:weak flaring} 
The following are equivalent:

1. A tree of spaces $\X=(\pi: X\to T)$ satisfies the proper $\kappa$-flaring condition. 

2. There exist $M_\kappa$ such that for all $D\ge 0$, there is {$\tau=\tau_{\ref{prop:weak flaring}}(\kappa,D)$} satisfying the following:

For every pair $\Pi$ of $\kappa$-qi sections $\gamma_0, \gamma_1$ over a geodesic interval $[-m, n]\subset T$ ($m, n\in \NN$),  
if 
$$
d_{X_i}(\gamma_1(i),\gamma_2(i)) >M_\kappa, \forall i\in [-m+1, n-1],$$
and
$$
\Pi_{max}=\max\left( d_{X_{-m}}(\ga_0(-m), \ga_1(-m)),  d_{X_{n}}(\ga_0(n), \ga_1(n)) \right) \le D,
$$
then 
$$
n+m\leq \tau.$$
\end{prop}
\proof First of all, we leave it to the reader to check that if (2) holds for all $m=n$ then 
it holds for all $m, n\in \NN$. Therefore, in what follows, in (2) we will be always assuming that $n=m$.

i. Assume that the proper flaring condition holds. Take $M_\kappa:=m_\kappa$ and consider a pair $\Pi=(\ga_0, \ga_1)$ 
of $\kappa$-qi sections over an interval $[-n, n]\subset T$ as in part (2). In particular, girth of $(\gamma_0,\gamma_1)$ is $>M_\kappa$. 
By the proper flaring condition, we have 
$$
D\ge \Pi_{\max} \ge \phi_\kappa(n). 
$$
Since $\phi_\kappa(t)$ is proper, the preimage $\phi_\kappa^{-1}([0,D])$ is contained in an interval $[0, t_{\kappa,D}]$. Then we take 
$$
\tau(\kappa,D):= 2t_{\kappa,D}. 
$$

ii. Conversely, suppose that (2) holds but proper flaring fails. 
Then  there exist a constant $D>0$ and a sequence $\Pi^m$ of pairs of $\kappa$-qi sections $\gamma^m_{0}, \gamma^m_{1}$
over some intervals $\llbracket s_m, t_m\rrbracket\subset T$ of length $2n_m$ with the midpoint vertex $r_m$ such that
$\Pi^m_0\to\infty, n_m\to\infty$, but
$$
\Pi^m_{\max}= \max\left( d_{X_{t_m}}(\ga^m_{0}(t_m), \ga^m_{1}(t_m)),  d_{X_{s_m}}(\ga^m_{0}(s_m), \ga^m_{1}(s_m)) \right) \le D. 
$$ 
We will isometrically parameterize the geodesic $[s_m, t_m]$ by the interval $[-n_m, n_m]\subset \Z$ so that 
$r_m$ corresponds to $0$. 
Set $\tau:= \tau(\kappa,M)$ where
$$
M=\max(D, M_\kappa). 
$$
Define the function
$$
\ell_m(i):= d_{X_i}(\ga^m_{0}(i), \ga^m_{1}(i)), i\in [-n_m, n_m]; \ell_m(0)= \Pi^m_0.  
$$
Then for sufficiently large $m$ we have 
$$
\Pi^m_0=\ell_m(0)> a(M+ b); a= L'_0, b=2L_0' \kappa, 
$$
and, hence, according to Corollary \ref{cor:contraction}, for all $n\in [-N_m+1, N_m-1]$, 
we have
$$
\ell_m(n)>M. 
$$
Here 
$$
N_m= \left\lfloor \log_a \left(   \frac{\Pi^m_0}{M +b}  \right) \right\rfloor. 
$$
Observe that the right hand side diverges to infinity as $m\to \infty$. Therefore, for sufficiently large $m$, $N_m>\tau/2$. 
Thus, we obtain a  contradiction with (2). \qed 

\medskip
While the proper flaring condition is quite natural, it is the condition (2) in the proposition 
that we will use throughout the book. 

\begin{defn}\index{uniform flaring}\label{uniform flaring}
We will say that a tree of spaces $\X$ satisfies the {\em uniform $\kappa$-flaring} condition with the parameter $M_\kappa$ if the condition (2) in the proposition holds. 
\end{defn}

\begin{convention}
In what follows, unless indicated otherwise,  $\kappa$-flaring  always means uniform $\kappa$-flaring. 
\end{convention}

\begin{lemma}\label{lem:hyp->uniform flaring} 
Suppose that $\X=(\pi: X\to T)$ is a tree of hyperbolic spaces with 
$\delta$-hyperbolic total space $X$.  Then $\X$ satisfies the uniform $\kappa$-flaring condition for all $\kappa\ge 1$. 
\end{lemma}
\proof As noted earlier, it suffices to consider the case $n=m$. Since $\ga_0, \ga_1$ are $\kappa$-quasigeodesics in $X$, they are within Hausdorff 
distance $D_{\ref{Morse}}(\delta,\kappa)$ from  geodesics $\ga_i^*$ in $X$ connecting the endpoints of $ga_0, \ga_1$ respectively. 
Take  $x_0=\ga_0(0)$ and $x_0^*\in \ga_0^*$ a point within distance $D_{\ref{Morse}}(\delta,\kappa)$ from $x_0$. The projections to $T$ of the geodesics $[\ga_0(-n) \ga_1(-n)]_X$, $[\ga_0(n) \ga_1(n)]_X$ each have length $\le D$. Thus, 
$$
d(x_0^*, [\ga_0(\pm n) \ga_1(\pm n)]_X)\ge D. 
$$

Suppose for a moment that $n-D> 2\delta$. By the slim quadrilateral property, there is a point $x_1^*\in \ga_1^*$ within distance $2\delta$ from $x_0^*$. (A priori, this could have been a point on one of two other sides of the geodesic quadrilateral with the vertices 
$\ga_i(\pm n), i=0,1$, but this possibility is ruled out by our assumption that  $n-D> 2\delta$.)  Thus, we find a point $x_1\in \ga_1\cap X_v$, within distance  
$$
D_0=2(\delta + D_{\ref{Morse}}(\delta,\kappa))+\kappa
$$
from $x_0$. While $v$ need not be equal to the vertex $0\in \llbracket -n, n\rrbracket\subset T$, we still have
$$
d_T(0, v)\le D_0. 
$$
In particular,
$$
d_{X_0}(\gamma_0(0), \ga_1(0))=\Pi_0\ge d_X(\ga_0(0), \ga_1(0))\le D_1=D_0(\kappa+1). 
$$ 
We, therefore, set
$$
M_\kappa= D_1 
$$
and $\tau(\kappa,D)= \delta+ \frac{1}{2}D$. Since in the uniform $\kappa$-flaring property, it is assumed, in particular, that 
$$
d_{X_0}(\gamma_0(0), \ga_1(0))> M_\kappa, 
$$
we obtain a contradiction with the above estimates, unless the inequality $n-D\ge 2\delta$ is violated, i.e. unless $n\le D+2\delta$, equivalently, the length of the interval $\llbracket -n, n\rrbracket$ is at most $\tau$, as required. \qed

\medskip 
The uniform flaring condition  has an immediate consequence that we will use on few occasions:

\begin{lemma}[Three flows lemma] \label{lem:3-flows} 
Suppose that $\pi: X\to T$ satisfies the uniform $K$-flaring condition. 
Suppose that $\ga_1, \ga_2, \ga_3$ are $K$-qi sections of $\pi: X\to T$ over an interval $\llbracket s,t\rrbracket $ such that for all $r\in \rrbracket s,t\llbracket$, 
$$d_{X_r}(\ga_1(r), \ga_3(r))> M_K$$
 while 
$$
\max_{i,j} d_{X_s}(\ga_i(s), \ga_j(s))\le C, \quad \max_{i,j} d_{X_t}(\ga_i(t), \ga_j(t))\le C. 
$$
Then the length of the interval $\llbracket s,t\rrbracket $ is uniformly bounded, i.e is 
$\le \tau_{\ref{lem:3-flows}}(K,C)$.  
\end{lemma}

The property appearing below will be also used quite often in our book: 

\begin{defn}
We say that a tree of spaces $\X$ {\em satisfies the $R(K,C)$-thin $K$-bigon property} if  there is a function $R(K,C)$ such that 
for every pair $\Pi=(\ga_1, \ga_2)$ of  $K$-qi sections of $\pi: X\to T$ over any interval $I=\llbracket v,w\rrbracket$, 
$$
\Pi_{max}\le C \Rightarrow   \quad   \forall t\in V(I), ~~
d_{X_t}(\ga_1(t), \ga_2(t))\le R(K,C).$$
\end{defn} 

Here, as before,  
$$
\Pi_{max}= \max\left( d_{X_v}(\ga_1(v), \ga_2(v)),   d_{X_w}(\ga_1(w), \ga_2(w)) \right). 
$$

\begin{cor}\label{cor:super-weak flaring} 
Aa tree of spaces $\X=(X\to T)$  satisfies the uniform $K$-flaring condition if and only if it satisfies the  
$R(K,C)$-thin  $K$-bigon property for some $R(K,C)=R_{\ref{cor:super-weak flaring}}(K,C)$. 
\end{cor} 
\proof 1. Assume that $\X$ satisfies the uniform $K$-flaring condition. Consider a pair of $K$-qi sections over 
an interval $I\subset T$. If for every vertex $t\in I$, $d_{X_t}(\ga_1(t), \ga_2(t))\le M_K$, then we are done. Otherwise, let 
$I'=\llbracket v', w'\rrbracket \subset \llbracket v,w\rrbracket $ be a maximal subinterval such that for all vertices $t\in I'$ we have 
$$d_{X_t}(\ga_1(t), \ga_2(t))> M_K.$$
Then there are edges $[v'',v'], [w',w'']$ in $I$ (not contained in $I'$) such that 
$$
d_{X_s}(\ga_1(s), \ga_2(s))\le C':=\max(M_K, C, 3\delta_0), s\in \{v'', w''\}. 
$$ 
By Lemma \ref{lem:3-flows} applied to $K$-qi sections $\ga_1, \ga_2=\ga_3$, restricted to $I'':= \llbracket v'', w''\rrbracket $, we obtain:
$$
d_T(v'', w'')\le \tau:=\tau_{\ref{lem:3-flows}}(K, C'). 
$$
By Lemma \ref{lem:growth-of-flare}, we get that for all $t\in V(I')$, 
$$
d_{X_t}(\ga_1(t), \ga_2(t))\le R_{\ref{cor:super-weak flaring}}(K,C):=a^{\tau} \left(C' + \frac{b}{a-1}\right),
$$ 
with $a=L'_0, b=2L'_0K$. (Recall that $L'_0\ge 2$.) 

\medskip 
2. We argue as in the proof of Proposition \ref{prop:weak flaring}. Suppose that proper $\kappa$-flaring fails. 
Then  there exist a constant $D>0$ and a sequence $\Pi^m$ of pairs of $\kappa$-qi sections $\gamma^m_{0}, \gamma^m_{1}$ 
over some intervals $J_m=\llbracket s_m, t_m\rrbracket\subset T$ of length $2n_m$ with the midpoint vertex $r_m$ such that
$\Pi^m_0\to\infty, n_m\to\infty$, but
$$
\Pi^m_{\max}= \max\left( d_{X_{t_m}}(\ga_{0,m}(t_m), \ga_{1,m}(t_m)),  d_{X_{s_m}}(\ga_{0,m}(s_m), \ga_{1,m}(s_m)) \right) \le D. 
$$ 
Setting $C:=D$, the hypothesis in Part 2 of the corollary means that  
$$d_{X_t}(\ga_1(t), \ga_2(t))\le R(K,C)$$
 for all vertices $t\in J_m$. This contradicts  
$\Pi^m_0\to\infty$. \qed 

\subsection{Acylindrical trees of spaces} \label{sec:acyl trees}

An easy, and frequently occurring, {\em sufficient} condition for uniform $\kappa$-flaring is {\em acylindricity}:

\begin{defn}\label{defn:acylindrical} \index{acylindrical tree of spaces} 
Fix constants $\kappa\ge 1$ and $\tau\ge 1$. A tree of spaces $(\pi: X\to T)$ is {\em $(M,\kappa,\tau)$-acylindrical} if for every pair of $\kappa$-sections $\ga_0, \ga_1$ over an interval $J\subset T$ of length $\ge \tau$, we have
$$
d_{X_v}(\ga_0(t), \ga_1(t))\le M, \forall t\in V(J). 
$$  


\end{defn}



\medskip 

We give few geometric examples of acylindrical trees of spaces in Section \ref{sec:qcamalgam}. In order to see that acylindrical trees of spaces satisfy uniform flaring, we take $M_\kappa:=M$ and $\tau(\kappa,D):= \tau+2$. Then, regardless 
of $D$, if $\Pi=(\gamma_0,\gamma_1)$ is a pair of $\kappa$-qi sections over an interval $J=\llbracket u, v\rrbracket\subset T$ and 
$$
d_{X_i}(\gamma_1(i),\gamma_2(i)) >M, i\in V(\rrbracket u, v\llbracket),
$$
then the length of $\rrbracket u, v\llbracket$ is $<\tau$ and, hence, $J$ has length $< \tau(\kappa,D)= \tau+2$.

The terminology {\em acylindrical} has its origin in 3-dimensional topology: A compact oriented 3-dimensional manifold with incompressible boundary $M$ is called (homotopically) {\em acylindrical} if every  map of an annulus $(A, \partial A)\to (M, \partial M)$ is homotopic (rel. $\partial A$) to a map $A\to \partial M$. Algebraically speaking, this condition means that if two elements of $\pi_1(\partial M,m)$ are conjugate in $\pi_1(M,m)$, then they are conjugate in $\pi_1(\partial M,m)$. If one glues two connected 
acylindrical 3-manifolds $M_1, M_2$ along their boundary surfaces to form a 3-manifold $M$, then every subgroup of 
$\pi_1(M)$ isomorphic to $\Z^2$ is contained (up to conjugation) in $\pi_1(M_1)$ or in $\pi_1(M_2)$. 
Algebraically speaking, topological acylindricity corresponds to acylindricity in the sense of group actions on trees (Definition \ref{defn:acylindrical-action}) as follows. 
The decomposition $M=M_1\cup M_2$ yields 
graph-of-groups decomposition of the fundamental $G= \pi_1(M)$. Let $G\times T\to T$ denote the action of $G$ on 
 the Bass--Serre tree $T$ corresponding to  this decomposition of $G$. Then the action of $G$ on $T$ is 1-acylindrical if and only if both manifolds $M_1, M_2$ are acylindrical. Suppose again that $G$ is the fundamental group of a finite graph of finitely generated groups  $({\mathcal G}, Y)$; let  $G\times T\to T$ be the corresponding $G$-action on the Bass--Serre tree and 
  $\X=(X\to T)$ the tree of spaces with $X$ equal to the Cayley graph of  $G$ as discussed in Example \ref{ex:BS-tree}. 
 We will see in Proposition \ref{example: acylindrical} that in this setting the tree of spaces $\X$ is $(\kappa,\tau)$-acylindrical provided that the action of 
 $G$ on $T$ is $k$-acylindrical for suitable values of $\kappa, \tau$ and $k$.

\subsection{Group-theoretic examples}

The following proposition was proved by Ilya Kapo\-vich \cite{ikap-acyl}; below, we give a different proof.

\begin{prop}\label{example: acylindrical}
Suppose $(\mathcal G, Y)$ is a finite graph of hyperbolic groups satisfying Axiom {\bf H} 
 and $G:=\pi_1(\mathcal G)$. 
If the $G$-action on the Bass--Serre tree $T$ of $\mathcal G$ is $R$-acylindrical in the
sense of Sela \cite{Sela97}, then for all $\kappa\geq 1$ there is a constant 
$M_{\kappa}$ such that the induced 
tree of metric spaces $\X=(\pi: X\map T)$ is $(M_{\kappa}, \kappa, R)$-acylindrical. 
In particular, in view Theorem \ref{thm:mainBF}, $G$ is hyperbolic.
\end{prop}
\proof The first part of the proof follows in the arguments in \cite[Section 3]{ps-limset}. 
We will  need  some properties of the tree of spaces $\X$ listed below. 

(1) The vertex-spaces of $\X$ are metric graphs which are isometric copies of various cosets of $G_y$'s in $G$, where $y\in V(Y)$. The map $\pi: X\map T$ is $G$-equivariant. 
The $G$-action on $X$ is proper and cocompact, and the stabilizer of each $v\in V(T)$ acts
on $V(X_v)$ transitively.

(2) Suppose that $\Gamma$ is a Cayley graph of $G$ with respect to a finite generating set. 
Let $f: G \map X$ be an orbit map. 
We know that for each $y\in V(Y)$ and $g\in G$, $gG_y$ is a vertex of $T$. 
 We have $\Hd(X_{gV_y}, f(gG_y))\leq D$, where
$D$ is a constant independent of $g\in G, y\in V(Y)$.
  
Suppose that the claim of the proposition  fails for some $\kappa$. 
Then there is a sequence of pairs of $\kappa$-qi sections
$\ga_{0,n}, \ga_{1,n}$ over geodesic intervals 
$$
\beta_n: [0, R+1]\map T$$
of length $R+1$, such that for some integer $t\in [0,R+1]$ we have 
$$
d_{X_{\beta_n(t)}}(\ga_{0,n}(t), \ga_{1,n}(t))\geq n, \quad \forall n\in \mathbb N.$$ 
Note that for all integers $s\in [0,R+1]$  
$$
d_X(\ga_{0,n}(s), \ga_{1,n}(s))\geq d_X(\ga_{0,n}(t), \ga_{1,n}(t))- d_X(\ga_{0,n}(s), \ga_{0,n}(t))-d_X(\ga_{1,n}(s), \ga_{1,n}(t)).
$$
Since $\ga_{0,n}, \ga_{1,n}$ are $\kappa$-qi sections, we have
\begin{align*}
d_X(\ga_{0,n}(s), \ga_{1,n}(s))\geq d_X(\ga_{0,n}(t), \ga_{1,n}(t))- 2\kappa|s-t|-2\kappa\geq \\ d_X(\ga_{0,n}(t), \ga_{1,n}(t))-2(R+1)\kappa-2\kappa.
\end{align*}
Since vertex-spaces of $\X$ are uniformly properly embedded in the ambient space $X$ 
we see that $d_X(\ga_{0,n}(t), \ga_{1,n}(t))\map \infty$ as
$n\map \infty$. Thus, $d_X(\ga_{0,n}(t), \ga_{1,n}(t))\map \infty$, which in turn implies that
$d_{X_{\beta_n(s)}}(\ga_{0,n}(s), \ga_{1,n}(s))\map \infty$ for all $s\in [0,R+1]$. 
Thus, passing to subsequence, if necessary, we may assume that
$d_{X_{\beta_n(s)}}(\ga_{0,n}(s), \ga_{1,n}(s)) \geq n$ for all $n\in \mathbb N$ and $s\in [0,R+1]$. Also, since the group $G$ acts on $T$ cocompactly, 
we can assume, by passing to subsequence if necessary, that $\beta_n(0)$ is a fixed vertex 
$v$ and $\ga_{0,n}(0)$ is a fixed point $x\in X_v$.
Since $X$ is quasiisometric to $G$, by passing to a further subsequence, if necessary, we may assume that $\beta_n$ is a fixed geodesic $vw$ in
$T$, where $d_T(v,w)=R+1$. We note that since 
$d_X(\ga_{0,n}(v), \ga_{1,n}(w))\leq \kappa +\kappa R$, by Lemma \ref{lem:edge-spaces} 
we have 
$$
\Hd([\ga_{0,n}(v) \ga_{1,n}(v)]_{X_v}, [\ga_{0,n}(w) \ga_{1,n}(w)]_{X_v})\leq C_{\ref{lem:edge-spaces}}(\kappa+ (R+1)\kappa).
$$
Now, by (2) above we have a constant $D_1$ and  $y, y'\in V(Y)$, $g, g'\in G$ such that 

(i) $v=gG_y$, $w=g'G_{y'}$ and 

(ii) the diameter of $N_{D_1}(gG_{y})\cap g'G_{y'}$ is infinite in $\Gamma$. 

\medskip 
The rest of the argument is borrowed from \cite[Theorem 4.6]{mitra-ht}. 
Let $\{h_n\}\subset gG_{y}$ and $\{h'_n\}\subset g'G_{y'}$ be sequences of distinct elements
such that $d_{\Gamma}(h_n, h'_n)\leq D_1$ for all $n\in \mathbb N$. 
Hence, $d_{\Gamma}(1, h^{-1}_nh'_n)\leq D_1$.
But there are only finitely many elements of $G$ inside $B(1; D_1)$. Hence, passing to a subsequence, we may assume that the sequence
$\{h^{-1}_nh'_n\}$ is constant. Let $x=h^{-1}_nh'_n$. Consider the equations 
$x=h^{-1}_mh'_m=h^{-1}_nh'_n$; whence $h_mx=h'_m, h_nx=h'_n$. Thus, we have
$$x^{-1}h^{-1}_mh_nx={h'}^{-1}_mh'_n \Rightarrow h^{-1}_mh_n=x{h'}^{-1}_mh'_nx^{-1}$$
$$ \Rightarrow h_m(h^{-1}_mh_n) h^{-1}_m=(h_mx){h'}^{-1}_mh'_n(h_mx)^{-1}=h'_m({h'}^{-1}_mh'_n){h'}^{-1}_m.$$
Clearly, $h^{-1}_mh_n\in G_y$ and, hence, $h_m(h^{-1}_mh_n) h^{-1}_m\in h_mG_yh^{-1}_m=gG_yg^{-1}$, since $h_m\in gG_y$.
Similarly, $h'_m({h'}^{-1}_mh'_n){h'}^{-1}_m\in g'G_{y'}g'^{-1}$. This implies that 
$$
h_m(h^{-1}_mh_n) h^{-1}_m=h'_m({h'}^{-1}_mh'_n){h'}^{-1}_m\in gG_yg^{-1}\cap g'G_{y'}g'^{-1}.
$$ 
However, $gG_yg^{-1}$ is the stabilizer of the vertex $v=gG_y$ and $g'G_{y'}g'^{-1}$ is the stabilizer of $w=g'G_{y'}$.
Since $\{h_n\}$ and $\{h'_n\}$ are sequences of distinct elements in $gG_y$ and $g'G_{y'}$ respectively, the intersection $G_v\cap G_w$ is
infinite. Since $d_T(v,w)=R+1$ this contradicts the $R$-acylindricity of the $G$-action.
\qed

\begin{rem}
1. The proof of Proposition \ref{example: acylindrical} also works even if we assume that $G_v\cap G_w$ finite whenever $d_T(v,w)\geq k+1$.

2. In fact, to conclude hyperbolicity of $G$ in the proposition, one does not need the full power of the Combination Theorem, Theorem \ref{thm:mainBF}, one can derive the result from the cobounded quasiconvex chain-amalgamation, Theorem \ref{thm:chain}. 
\end{rem}


For the next corollary, we recall that a subgroup $H$ in a group $G$ is {\em weakly malnormal} if for every $g\in G\setminus H$ the intersection
$$
gHg^{-1} \cap H
$$
is finite. 

\begin{cor} \label{cor:malnormal-amalgam} 
{\em (\cite[Theorem 2]{km-malnormal})}
If $G_1, G_2$ are hyperbolic groups and $H$ is a common quasiconvex subgroup which is weakly malnormal
in either in $G_1$ or in $G_2$, then $G=G_1*_H G_2$ is hyperbolic.
\end{cor}
\proof We claim that the action of $G$ on the Bass--Serre tree for the given amalgam decomposition is $3$-acylindrical. Without loss of generality let us assume that $H<G_1$ is weakly malnormal.
If $T$ is the Bass--Serre tree and $v,w\in V(T)$ with $d_T(v,w)\geq 4$, then there are is a sequence of consecutive vertices on $vw$ 
of the form $v_1=xG_2, v_2=yG_1, v_3=zG_2$, where $x,y,z\in G$. Then $G_{v_1}\cap G_{v_3}$ is equal to the intersection of  
the stabilizers of the two edges:

(i) one edge connecting $v_1,v_2$, and 

(ii) the one edge connecting $v_2,v_3$. 

However, these are two distinct conjugates of $H$ in the stabilizer of $v_2=yG_1$, i.e. they are 
 of the form $ygHg^{-1}y^{1}, yg'H{g'}^{-1}y^{-1}$ in $yG_1y^{-1}$ where $g,g'\in G_1$.
Since  
$$
ygHg^{-1}y^{1}\cap yg'H{g'}^{-1}y^{-1}=y(gHg^{-1}\cap g'H{g'}^{-1})y^{-1}$$ 
(weak) malnormality of $H$ in $G_1$
proves our claim. Then the hyperbolicity follows from Proposition \ref{example: acylindrical}. \qed

An example analogous to the situation of the corollary above in the context of a tree of spaces is discussed in Section \ref{sec:qcamalgam}. 

\subsection{Exponential flaring (Bestvina--Feighn  flaring condition)}

\begin{defn}[Exponential flaring condition]\label{exp-flare}\index{exponential flaring}
We say that a tree ${\mathfrak X}$ of metric spaces $\pi:X\rightarrow T$ 
it satisfies the {\em Bestvina--Feighn} $\kappa$-{\em flaring  condition} 
or the {\em exponential} $\kappa$-{\em flaring  condition}, if there exist  $\lambda_\kappa>1, M_\kappa>0$ 
and  $n_\kappa\in \mathbb N$ such that 
the following holds:

For every pair $\Pi=(\ga_0,\ga_1)$ of $\kappa$-qi sections of ${\mathfrak X}$ over a length $2n_\kappa$ geodesic interval 
$\llbracket -n_\kappa,n_\kappa\rrbracket \subset T$, if  the girth $\Pi_0$ of the pair $(\ga_0, \ga_1)$ is $\ge M_\kappa$,  then 
$$
\lambda_\kappa\cdot \Pi_0 \le \Pi_{\max}.
$$
\end{defn}

A form of this flaring condition first appeared in the paper  \cite{BF} of Bestvina and Feighn. 
Actually, the original Bestvina--Feighn flaring condition was a bit different from the exponential flaring condition above 
as they required not just two qi sections but a 1-parameter family of $\kappa$-qi sections interpolating these two, 
i.e. a $\kappa$-hallway, see Definition \ref{defn:hallway}. 
The existence of such a family (with a different but uniform qi constant $\kappa'$) 
follows from \cite{pranab-mahan}. It will be also proven in Lemma \ref{lem:E-ladder-structure}(b).

We will see below that the exponential flaring implies proper flaring with an exponential function $\phi_\kappa$ 
and that if $X$ is hyperbolic, then $\X$ satisfies the exponential  
$\kappa$-flaring condition for all $\kappa\ge 1$. 
Note that while in their first paper  \cite{BF}  Bestvina and Feighn imposed the exponential flaring condition 
for all $\kappa\ge 1$, in the addendum \cite{BF-err} to their paper, 
the flaring condition was required only for some value of $\kappa$, cf. the statement of our main result, 
Theorem \ref{thm:mainBF}.

\begin{lemma}\label{lem:exp} 
Bestvina--Feighn $\kappa$-flaring implies exponential proper $\kappa$-flaring. Moreover, the proper flaring condition holds either in the negative or in the positive direction (see Definition \ref{defn:signed flare}). 
\end{lemma}
\proof We fix $\kappa$ and set $n:= n_\kappa, \la:=\la_\kappa$. Suppose that $\Pi=(\ga_0,\ga_1)$ is a pair of $\kappa$-qi sections over a geodesic interval $I$ of length $N= s n$ and of girth $\Pi_0\ge m_\kappa$. For concreteness, we assume that 
$$
d_{X_{n}}(\gamma_0(n), \gamma_1(n))\ge \la d_{X_{0}}(\gamma_0(0), \gamma_1(0))=\Pi_0. 
$$
Then, applying the flaring inequality to the subinterval in $I$ of length $2n$ centered at $n$, we obtain
$$
\max( d_{X_{2n}}(\gamma_0(2n), \gamma_1(2n)), \Pi_0)\ge \la d_{X_{n}}(\gamma_0(n), \gamma_1(n)). 
$$
Since $\la>1$, the maximum in this inequality is attained by $d_{X_{2n}}(\gamma_0(2n), \gamma_1(2n))$ and, thus,
$$
d_{X_{2n}}(\gamma_0(2n), \gamma_1(2n)) \ge \la d_{X_{n}}(\gamma_0(n), \gamma_1(n)).
$$
Applying this argument inductively, we obtain:
$$
\lambda^s_\kappa\cdot \Pi_0\leq d_{X_{sn}}(\gamma_0(sn), \gamma_1(sn))\le \Pi_{\max}(sn).$$
By reducing $\la$ to $\mu>1$ if necessary and using Lemma \ref{lem:growth-of-flare}, we also get
$$
\Pi_{\max}(m)\ge  d_{X_{m}}(\gamma_0(m), \gamma_1(m))\ge \mu^m \Pi_0, \forall m\ge n. 
$$
Since the function $m\mapsto \mu^m, m\in \NN$, is proper, 
the exponential proper $\kappa$-flaring condition for ${\mathfrak X}$ follows.  \qed

\begin{prop}\label{hyp to lin flaring}
If $\X$ satisfies the proper $\kappa$-flaring condition for all $\kappa\ge 1$, then  $\X$ also satisfies satisfies an exponential $\kappa$-flaring condition for all $\kappa\geq 1$. In particular, if $X$ is hyperbolic, then $\X$ satisfies satisfies an exponential $\kappa$-flaring condition for all $\kappa\geq 1$.
\end{prop}
\proof Since $X$ is hyperbolic, the tree of spaces $\X=(X\map T)$ satisfies both the proper $\kappa$-flaring and 
the property obtained in Corollary \ref{cor:super-weak flaring} for all $\kappa\ge 1$. 
We will use both of these in the proof. 
The proof is inspired by, but is conceptually simpler than \cite[Proposition 5.8]{pranab-mahan}.
For each $K\geq 1$, we inductively define $K_0:=K$ and $K_i:=\max\{K_{i-1}, C_{\ref{lem:edge-spaces}}(K_{i-1})\}$, 
$i\geq 1$. Given $\kappa \geq 1$ we set  
$$
L:=\eta_{\ref{unif-emb-subtree}}(2\kappa_3), \epsilon=3\eta_{\ref{unif-emb-subtree}}(2\kappa_3), \quad 
R:=\max\{1, m_{\kappa_3}, L(5\epsilon+4L)\}$$
and $D:=\max\{R, R_{\ref{cor:super-weak flaring}}(\kappa_3,R)\}$.
Let $n=n_{\kappa}$ be any integer such that $\phi_{\kappa_3}(n)\geq 12D$; set $\lambda_{\kappa}:=2$
and $M_{\kappa}:=D+1$. 

If $\Pi=(\ga_0,\ga_1)$ is a pair of $\kappa$-qi sections over an interval $J=[-n,n]\subset T$, $\Pi_0\geq M_{\kappa}$,  
then we form a metric bundle $\Y= (Y\to J)$:

The vertex-spaces $Y_i$ of $Y$ are geodesic segments in $X_i$ 
 joining  $\ga_0(i),\ga_1(i)$. The edge-spaces $Y_e$, $e=[i, i+1]$, of $\Y$ are geodesic segments in $X_e$ with end-points 
 within distance $\kappa$ from the respective end-points of $Y_i$. The incidence maps of $\Y$ are obtained by composing the incidence maps of $f_{ev}$, $v\in V(J)$, composed with the nearest-point projections to $Y_v$ (taken in $X_v$).

After that, the idea is to first decompose this interval-bundle into a finite number of subbundles by constructing
qi sections in $Y$ (cf. \cite[Proposition 3.14]{pranab-mahan}, also Proposition \ref{vertical subdivision}), where the
subbundles intersect along the qi sections. We then use proper flaring to prove that the qi sections bounding
each of the subbundles flare in at least one direction. Finally, as in the last step of the proof of 
\cite[Proposition 5.8]{pranab-mahan}, we verify that at least half of these will flare in the same direction
to finish the proof. 

{\bf Step 1: Construction of qi sections in $Y$.}
We note that through any point of the metric bundle formed by two $K_{i-1}$-qi sections, there is a 
$K_i$-qi section, $i\geq 1$. Let $\alpha_i=Y\cap X_i$, $i\in V(J)$. For two consecutive integers $i, j$ we have
a map $h_{ij}: \alpha_i\map \alpha_j$ such that for all $x\in \alpha_i$, $d_{X_{ij}}(h_{ij}(x),x)\leq \kappa_3$. 
This map is clearly $\eta_{\ref{unif-emb-subtree}}(2\kappa_3)$-coarsely Lipschitz, with a similarly defined
$\eta_{\ref{unif-emb-subtree}}(2\kappa_3)$-coarse inverse $h_{ji}:\alpha_j\map \alpha_i$, which is also an 
$\eta_{\ref{unif-emb-subtree}}(2\kappa_3)$-coarsely Lipschitz map. Hence, by Lemma \ref{lem: qi from lipschitz}, 
the maps $h_{ij}, h_{ji}$ are both $(\eta_{\ref{unif-emb-subtree}}(2\kappa_3), 3\eta_{\ref{unif-emb-subtree}}(2\kappa_3))$-quasiisometries.
By Lemma \ref{lem:c-mono}, if $x, y, z\in \alpha_i$ and $y$ is between $x,z$ with $d_{X_i}(x,y)\geq L(5\epsilon+4L)$,
and $d_{X_i}(y,z)\geq L(5\epsilon+4L)$, then $h_{ij}(y)$ is between $h_{ij}(x)$ and $h_{ij}(z)$.
In particular, this is true if $d_{X_i}(x,y)\geq R$ and $d_{X_i}(y,z)\geq R$ by the choice of $R$.

Suppose $l(\alpha_0)=l$ and let $\alpha_0$ also denote the parametrization of this geodesic in $X_0$ so that
$\alpha_0(0)=\ga_0(0), \alpha_0(l)=\ga_1(0)$. Next, we inductively construct a sequence of numbers 
$s_0=0, \cdots, s_n=l$ and a sequence of $\kappa_1$-qi sections $\ga_0=\beta_0,\beta_1,\cdots, \beta_n=\ga_1$ 
in $Y$ such that each $\beta_{i+1}$ is contained in the  metric bundle formed by $\beta_i$ and $\beta_n$, 
$0\leq i\leq n-2$ as follows. Suppose $s_0,\cdots, s_i$ are chosen and so are $\beta_0,\cdots, \beta_i$ and
$s_i<l$. To construct $s_{i+1}$ and $\beta_{i+1}$ consider the subset $S\subset (s_i,l]$ consisting of $s$ such that 
there is a $\kappa_2$-qi section $\beta$ through $s$ satisfying
$$
\min_jd_{X_j}(\beta(j), \beta_i(j))\leq R.
$$

If $S=\emptyset$ then define $s_{i+1}=s_n=\ga_1$. Assume now that $S$ is nonempty. 
Suppose there is $s\in S$ and a $\kappa_1$-qi section $\beta$ in $Y$ through $\alpha_0(s)$ 
such that $\min_j d_{X_j}(\beta(j), \beta_i(j))= R$. 
In this case, if $s\neq l$, then define 
$s_{i+1}=s$, $\beta_{i+1}=\beta$. Otherwise, if $s=l$, then we define $s_{i+1}=s_n=l$ and $\beta_n=\ga_1$. 
Suppose there is no such $s\in S$. Then let $s_{i+1}=\min\{l, 1+\sup S\}$. If $s_{i+1}\neq l$, then let $\beta_{i+1}$ be 
any $\kappa_1$-qi section in $Y$ passing through $s_{i+1}$. Otherwise, define $s_n=s_{i+1}$ and $\beta_{i+1}=\ga_1$.
We note that $s_{i+1}-s_i\geq R$ unless $s_{i+1}=s_n$.

\medskip 
{\bf Step 2: Verification of the properties of the qi sections.}
Let $\Pi^i=(\beta_i, \beta_{i+1})$ and let $Y^i$ denote the interval-bundle over $J$ formed by these qi sections. 
We claim that
$Y^i\cap Y^j=\emptyset$, unless $|i-j|\leq 1$, and $Y^i\cap Y^{i+1}=\beta_{i+1}$ for all permissible $i$. 
Both  claims follow from Lemma \ref{lem:c-mono}, cf. Lemma 3.12 of \cite{pranab-mahan}. 

\medskip 
{\bf Step 3: Flaring of $\Pi^i=(\beta_i, \beta_{i+1})$.}
 We know that there is a $\kappa_2$-qi section
 $\bar{\beta}_i$ through either $s_{i+1}$ or $s_{i+1}-1$ inside the subbundle $Y^i$, such that
$$
d_{X_j}(\beta_i(u_i), \bar{\beta}_i(u_i))\leq R$$ 
for some $u_i\in V(J)$. Without loss of generality, we may assume that 
$u_i<0$. Now we construct a new set of $\kappa_3$-qi sections inside the bundle formed by $\beta_i$ and $\bar{\beta}_i$ 
as follows. Let $r=\lfloor (s_{i+1}-s_i-1)/D\rfloor$. Let $\beta'_j$, $0\leq j\leq r$ be arbitrary $\kappa_3$-qi sections 
in the bundle formed by $\beta_i, \bar{\beta}_i$ such that $\beta'_0=\beta_i$, and for $j\neq 0$
$\beta'_j$ passes through $\alpha_0(s_i+jD)$. It follows from Lemma \ref{cor:super-weak flaring} that for all $j, k\geq 0$
$d_{X_j}(\beta'_k(j), \beta'_{k+1}(j))\geq R$ and, thus, as in step 2, by Lemma   \ref{lem:c-mono} we see that
$$
\beta_i(j)=\beta'_0(j),\cdots, \beta'_m(0), \beta_{i+1}(j)$$
is a monotonic sequence in the geodesic interval $[\beta_i(j) \beta_{i+1}(j)]_{X_j}$
for all $j\geq 0$. Thus,  
\begin{align*}
d_{X_n}(\beta_i(n), \beta_{i+1}(n))\geq \sum^{r-1}_{j=0} d_{X_n}(\beta'_j(n), \beta'_{j+1}(n)) \geq \\
\sum^{r-1}_{j=0}12D
=\sum^{r-1}_{j=0} 12d_{X_0}(\beta'_j(0), \beta'_{j+1}(0)) = 12d_{X_0}(\beta'_0(0), \beta'_r(0))\end{align*}
by the proper flaring and by the choice of $n$. However,  
$$
d_{X_0}(\beta_i(0), \beta_{i+1}(0))=d_{X_0}(\beta'_0(0), \beta'_r(0))+d_{X_0}(\beta'_r(0), \beta_{i+1}(0)),$$ 
where $d_{X_0}(\beta'_r(0), \beta_{i+1}(0)\leq D+1$. It follows that
$$
d_{X_n}(\beta_i(n), \beta_{i+1}(n))\geq \frac{12D}{2D+1}d_{X_0}(\beta_i(0), \beta_{i+1}(0))\geq 
4d_{X_0}(\beta_i(0), \beta_{i+1}(0)),$$ 
since $D\geq 1$.

\medskip 
{\bf Step 4: Exponential flaring of $\ga_0, \ga_1$.}
We know that each pair $\Pi^i=(\beta_i, \beta_{i+1})$ exponentially flares in at least one direction (say, in the positive direction) 
by Step 3. Then there is a subset of induces 
${\mathcal I}\subset \{1, 2,...\}$, such that 
$$
\sum_{i\in {\mathcal I}} \Pi^i_0 \ge \frac{1}{2} \Pi_0.  
$$
It follows that $\Pi$ flares exponentially in the positive direction with $\lambda_{\kappa}=2$. \qed

\begin{cor}\label{cor:signed flare} 
If $\X$ satisfies the proper $\kappa$-flaring condition for all $\kappa\ge 1$, then (again, for all $\kappa\ge 1$) it satisfies proper flaring either in positive or in the negative direction.  
\end{cor}

\section{Hyperbolicity of  trees of hyperbolic spaces}

\subsection{The combination theorem} 

We are now ready to state our version of the  combination theorem of Bestvina and Feighn \cite{BF}: 

\begin{theorem}\label{thm:mainBF}
There exist $K_*=K_{\ref{thm:mainBF}}(\delta_0,L_0)$ and $\delta_*=\delta_{\ref{thm:mainBF}}(\delta_0,L_0)$, 
depending only on $\delta_0$ and $L_0$,  
such that the following holds. Suppose $\X=(\pi: X\map T)$ is a tree of hyperbolic spaces 
(with primary parameters $\delta_0, L_0$) satisfying 
the uniform   $K_*$-flaring condition. 
Then $X$ is a $\delta_*$-hyperbolic metric space. 
\end{theorem}

The constants $K_*$ and $\delta_*$ are computable. In Remark \ref{rem:K*}  we will give a formula for $K_*$, which is inductive in nature, as it relies upon earlier computations of various constants and functions scattered throughout the book. (We will not attempt to write a formula for $\delta_*$.)   
The reader unwilling to keep track of such computations, can simply assume that ${\mathfrak X}$ 
satisfies the uniform $\kappa$-flaring condition {\em for all} $\kappa\ge 1$.

\subsection{Cobounded quasiconvex chain-amalgamation}\label{sec:qcamalgam}

In the book  we will be frequently using the following 
 very special case of  Theorem \ref{thm:mainBF} which is much easier to prove, see e.g. 
 \cite[Proposition 1.51]{pranab-mahan}. This special case was motivated by a result of Hamenstadt, \cite[Lemma 3.5]{hamenst-word}. Although Hamenstadt
used much stronger assumptions, it is clear that the proof of Hamenstadt goes through with the weaker hypothesis as well. 
We include a proof along the lines of Hamenstadt's arguments for the sake of completeness  
and also since we want a description of geodesics.

We assume that $X$ is a path-metric space that can be represented as a union of a finite chain of rectifiably-connected subsets equipped with induced path-metrics 
$$
Y= Q_0\cup Q_1\cup  ... \cup Q_{n},
$$
such that for some constants $C$ and $\delta$ the following hold: 

\begin{enumerate}
\item Each $Q_i$ is $\delta$-hyperbolic.  

\item For each $i<n$ the intersection $Q_{i,i+1}= Q_{i}\cap Q_{i+1}$ is rectifiably connected  and $L$-qi embedded 
in $Q_i, Q_{i+1}$. 

\item  Each $Q_{i,i+1}$ separates (in $Y$) $Q_{i}$ from $Q_{i+1}$ in the sense that every path $c$ in $X$ connecting $Q_i$ to $Q_{i+1}$ has to cross $Q_{i,i+1}$. 

\item Each pair of intersections $Q_{i-1,i}, Q_{i,i+1}$ is $C$-cobounded in $Q_i$.

\item $d_{Q_{i}}(Q_{i-1,i}, Q_{i,i+1})\ge 1$. 
\end{enumerate}

\medskip 
We will say that such $X$ is a {\em cobounded quasiconvex chain-amalgam} of $Q_i$'s.  If $n=1$, we 
will refer to $X=Q_0\cup Q_1$ simply as a {\em quasiconvex amalgam}. 
\index{cobounded quasiconvex chain-amalgam}

Clearly, the collection  $Q_i$'s in a cobounded quasiconvex chain-amalgam 
gives $X$ structure of a tree of hyperbolic spaces with vertex-spaces $Q_i$ 
and edge-spaces $Q_{i,i+1}$, such that the  tree $T$ is isometric to the interval $J$ of length $n+1$ in $\RR$ with integer vertices. Conversely, consider a tree of hyperbolic spaces $\X$ over an interval $T$ such that for each vertex $v$ with the incident edges $e_\pm$, the corresponding subsets $X_{e_\pm v}$ are $C'$-cobounded in $X_v$. Then $\X$ 
yields a cobounded quasiconvex chain-amalgam with 
subsets $Q_i=Q_v$, $v=v_i$, equal to the unions 
$$
X_{e_{\scriptstyle -}}\times \left[\frac{1}{2},1 \right] \cup_{f_{e_-v}} X_v \cup_{f_{e_+v}} X_{e_{\scriptstyle +}}\times \left[0,\frac{1}{2}\right], 
$$ 
$$
Q_{i-1,i}= X_{e_{\scriptstyle -}} \times \frac{1}{2}, 
$$
see Section \ref{sec:cylinders} for the definition of mapping cylinders.

For each $i$ pick  points 
$$
x_i^-\in N_{C'} (P_{Q_i,Q_{i-1,i}}(Q_{i,i+1})) \cap Q_{i-1,i}, 
x_{i}^+\in N_{C'} (P_{Q_i,Q_{i,i+1}}(Q_{i-1,i}))\cap Q_{i,i+1},$$
where the $C'$-neighborhoods are taken with respect to the metric of $Q_{i}, Q_{i}$. 
Since both projections of $Q_{i-1,i}$ to  $Q_{i,i+1}$ and of $Q_{i,i+1}$ to $Q_{i-1,i}$ have diameters $\le C$, we obtain
\begin{equation}\label{eq:min-cobounded} 
d_{Q_i}(Q_{i-1,i}, Q_{i,i+1})\le d_{Q_i}(x_i^-, x_i^+)\le d_{Q_i}(Q_{i-1,i}, Q_{i,i+1})+2(C+C'),
\end{equation}
i.e. the pair of points $x_i^-, x_i^+$ ``almost'' realizes the minimal distance in $Q_i$ between the subsets $Q_{i-1,i}$,  
$Q_{i,i+1}$.  

We will simultaneously prove 
hyperbolicity of $X$ and describe uniform quasigeodesics connecting points in $X$. For this description, given points 
$x\in Q_{i-1}, x'\in Q_{k+1}$, it will be convenient to name their nearest-point projections 
(in $Q_{i-1}, Q_{k+1}$) to $Q_{i-1,i}, Q_{k,k+1}$ as $\bar{x}, \bar{x}'$, respectively. Suppose, furthermore, that 
$$
c(x_{i}^+, x_{i+1}^-), $$
are $L'$-quasigeodesic paths in $Q_{i,i+1}$ connecting $x_{i}^+$ to $x_{i+1}^-$ 
and 
$$
c(x_i^-, x_{i}^+), c(x,\bar{x}),  c(\bar{x}', x')
$$
are $L'$-quasigeodesic paths in $Q_{i}$ connecting $x_i^-$ to $x_{i}^+$, etc. We let $c^*(\cdot, \cdot)$ denote the corresponding 
geodesics paths in $Q_{i,i+1}$, $Q_{i}$, connecting the respective points.

\begin{thm}\label{thm:chain}\label{thm:hyp-tree}
Under the above assumptions, $X$ is $\delta_{\ref{thm:chain}}(\delta,L,L',C,C')$-hyperbo\-lic. Moreover, the following paths $c(x,x')$ 
are $K=K_{\ref{thm:chain}}(\delta,L,C,D,L')$-quasigeodesics in $X$ connecting $x\in Q_{i-1}, x'\in Q_{k+1}$, $i\le k$: 

1. If  $x, x'$ belong to $Q'_i=Q_i\setminus (Q_{i-1,i}\cup Q_{i,i+1})$ for some $i$, then we assume that $c(x,x')$ is an $L'$-quasigeodesic in $Q_i$ connecting $x$ to $x'$.  

2. Otherwise,  
$$
c(x,x')= c(x,\bar{x})\star c(\bar{x}, x_{i}^-) \star c(x_i^-, x_{i}^+) \star 
c(x_{i}^+, x_{i+1}^-) \star ... \star c(x^+_k,\bar{x}') \star c(\bar{x}', x'). 
$$
\end{thm}
\proof This theorem is proven by verifying the assumptions of Corollary \ref{cor:bowditch}, i.e. axioms of a slim combing. 

(a1)  
We will need to estimate the length of $c(x,x')$ in terms of $d(x,x')$. First of all, 
$$
\length(c(x,x'))\le L' (\length (c^*(x,x')) +1), 
$$
hence, it suffices to get an estimate for $c^*$. 

Let $\gamma$ be any geodesic in $X$ connecting $x$ to $x'$. In view of the separation property (3) in the theorem, for each 
$i\le j\le k$, $\ga$ will contain subpaths $\ga_j\subset Q_j$ (necessarily a geodesic in $Q_j$) 
connecting a point $p_j^-\in  Q_{j-1,j}$ to some $p_j^+\in  Q_{j,j-1}$. 

Let $P_-, P_+$ denote the projections $Q_j\to Q_{j-1,j}$, $Q_j\to Q_{j,j+1}$ respectively. 
According to Lemma \ref{lem:projection-1}, $\ga_j$ contains points $y_j^\pm$ satisfying 
$$
d_{Q_j}(y_j^\pm, P_\pm(p_j^\mp))\le 2\delta+\la, 
$$
hence,
$$
d_{Q_j}(y_j^\pm, x_j^\pm)\le D:=C+C'+2\delta+\la,  
$$
where $\la=\la_{\ref{lem:qi-preserves2}}(\delta,L)$ 
is the quasiconvexity constant of $Q_{j,j\pm 1}$ in $Q_j$. Thus, 
$$
\length( [p_j^- x_j^-]_{Q_j} \star [x_j^- x_j^+]_{Q_j} \star     [x_j^+ p_j^+]_{Q_j} ) \le \length(\ga_j) + 4D.  
$$
Since $Q_{j-1,j}, Q_{j,j+1}$ are $L$-qi embedded in $Q_j$ we also obtain 
\begin{equation}\label{eq:gammaj}
\length(c(p_j^-, p_j^+))\le L \cdot \length(\ga_j) + 4D (L+1). 
\end{equation}
Since 
$$
\length(\ga)\ge  d(x, p^-_{i}) + \sum_{j=i}^{k} \length(\ga_j) + d(p_{k}^+, x'), 
$$
by combining the inequalities \eqref{eq:gammaj}, we get: 
$$
 \length(c^* (x^-_{i}, x_{k}^+)) \le  \length(c (p^-_{i}, p_{k}^+)) \le L  \sum_{j=i}^{k} \length(\ga_j)  + 4D (L+1)(k-i+1).
$$
To estimate the term $4D (L+1)m$, $m=k-i+1$, note that $d(x,y)\ge m$ (in view of the assumption 5 in the theorem). Thus,
$$
 \length(c^* (x^-_{i}, x_{k}^+)) \le L d(x,x') + 4D (L+1) d(x,x') = (L+ 4D (L+1)  )d(x,x'). 
$$

Lastly, we deal with $d(x, x^-_{i})$ and $d(x', x_{k}^+)$. Recall that the metric space $(Q_i, d_{Q_i})$ 
is $\eta$-properly embedded in $X$ (Proposition \ref{unif-emb-subtree}). We obtain:
\begin{align*}
\length c^*(x,  x^-_{i}) = d_{Q_i}(x, \bar{x}) + d_{Q_{i-1, i}}(\bar{x},  x^-_{i})\le \\
d_{Q_{i-1}}(x, \bar{x}) + L + L d_{Q_{i-1}}(\bar{x},  x^-_{i+1}) \le L + 2D + L d_{Q_{i-1}}(x,  x^-_{i}) \le \\
   L + 2D + L  \eta( d(x,  x^-_{i}) ) \le \\
    L + 2D +  L  \eta( d(x,  y^-_{i}) +D) \le  L + 2D +  L  \eta( d(x,  x') +D). 
\end{align*}
 Similarly, 
$$
\length c^*(x',  x^+_{k}) \le L + 2D +  L  \eta( d(x,  x') +D). 
$$
Combining the inequalities, we obtain: 
\begin{align*}
\length(c^*(x,x'))= \length c^*(x,  x^-_{i}) +  \length(c (x^-_{i}, x_{k}^+) + \length c^*(x',  x^+_{k}) \le \\
\eta_{\ref{thm:chain}}(d(x,x')) := 2(L + 2D +  L  \eta( d(x,  x') +D)) + (L+ 4D (L+1)  )d(x,x').  
\end{align*}

\medskip

(a2) Consider a triple of points $x\in Q_i, y\in Q_j, z\in Q_k$, $i\le j\le k$, and the ``triangle''
$$
\Delta_c=c(x,y)\cup c(y,z)\cup c(z,x). 
$$
By the definition of the paths $c$ in $X$, the paths $c(x,y), c(y,z)$ coincide away from $Q_j$, the same applies to the pair of paths 
$c(y,z), c(z,x)$. Therefore, it suffices to consider the case when $i=j=k$.  

We will use the notation $pq$ for geodesics in $Q_i$. Our goal is to verify that each of the paths $c(p,q)$ connecting points $p, q\in Q_i$ are uniformly Hausdorff-close to a geodesic $pq$: Once we are done with this, then uniform slimness of 
$\Delta_c$ will follow from $\delta$-hyperbolicity of $Q_i$. 

If both $p, q$ are in $Q'_i$ or in $Q_{i,i+1}$ or in $Q_{i-1,i}$, 
there is nothing to prove. Hence, up to permutation of the points $p, q$ and reversing the order in the interval $[0, n]$, there are two cases to consider, depending on the position of the points $p, q$ with respect to the subsets $Q_{i-1,i},   Q_{i,i+1}$:

Case 1. Suppose that $p\notin Q_{i-1,i}$ and $q\in Q_{i,i+1}$. We will be using the notation 
$\bar{p}=P_{Q_i,Q_{i,i+1}}(p)$.  Then 
$$
c(p,q)= c(p, \bar{p})  \cup c(\bar{p}, p) 
$$ 
Since $Q_{i,i+1}$ is $L$-qi embedded in $Q_i$, this path  is  $D_{\ref{stab-qg}}(\delta, LL')$-Hausdorff close to the union $p \bar{p} \cup \bar{p} q$. According to Lemma \ref{rem:lip-proj}, 
$$
\Hd_{Q_i}(p \bar{p} \cup \bar{p} q, pq)\le \la+2\delta,$$
concluding the proof in this case. 

Case 2. Suppose that $p\in Q_{i-1,i}$ and $q\in Q_{i,i+1}$. In view of the assumption that  $Q_{i-1,i}, Q_{i,i+1}$ are 
$L$-qi embedded in $Q_i$, we will work with $Q_i$-geodesics connecting pairs of points points in $Q_{i-1,i}$ and 
pairs of points in $Q_{i,i+1}$ instead of the $c$-paths in $Q_{i-1,i}$ and   $Q_{i,i+1}$. 
Continuing with the notation of Case 1, 
$$
\Hd_{Q_i}(pq, p \bar{p} \cup \bar{p} q)\le \la+3\delta 
$$
(see Lemma \ref{lem:projection-1}) and 
$$
d(\bar{p}, x_i^+)\le C.
$$
Thus,
$$ 
\Hd_{Q_i}(pq, p x_i^+ \cup x_i^+ q)\le C+\la+4\delta. 
$$
Similarly,
$$
\Hd_{Q_i}(p x_i^+, px_i^- \cup x_i^- x_i^+)\le C+\la+4\delta. 
$$
Combining the inequalities, we obtain 
$$
\Hd_{Q_i}(pq, p x_i^+  \cup x_i^- x_i^+ \cup x_i^+ q)\le 2(C+\la+4\delta). \qed 
$$

\begin{rem}
The flexibility of working with concatenations of quasigeodesics points $x_i^\pm$ uniformly close to nearest-point projections will be critical in several places in the book, e.g. in Chapter \ref{ch:CT}. 
\end{rem}

\begin{cor}\label{cor:chain}
Assuming that $X\to T$ is a tree of hyperbolic spaces satisfying the assumptions of Theorem \ref{thm:chain}, for every subinterval $S\subset T$, the inclusion map
$$
X_S\to X
$$  
is an $L_{\ref{cor:chain}}(\delta,L,C)$-qi embedding. 
\end{cor} 
\proof This is an immediate consequence of the of the fact that for any pair of points $x, y'\in S_X$, the path $c_{X_S}(x,x')$ equals 
the path $c_{X}(x,x')$ where the subscript denotes the space in which we define the combing. \qed


\begin{cor}\label{cor:edge-spaces}
Suppose that ${\mathfrak X}=(\pi: X\to T)$ is a tree of spaces satisfying Axiom {\bf H}. Then, 
for every edge $e=[v,w]$ of $T$ the space 
$X_{vw}$ 
equipped with its natural path-metric, is $\delta'_0$-hyperbolic with the 
hyperbolicity constant $\delta'_0$ depending only on the primary parameters of the tree of hyperbolic spaces 
${\mathfrak X}$. 
\end{cor}

We now return to the general quasiconvex chain-amalgamation and relate this class of trees of spaces to acylindrical trees of spaces. 
Suppose that $\gamma$ is a $\kappa$-qi section of the tree of spaces $\X=(\pi: X\to J)$, defined on an interval $I\subset J$. Thus, for each integer $i$ we have a point $x_i\in Q_i$ and $d_Y(x_i, x_{i+1})\le K$ for some $K$ depending on $\kappa$.  
If the length of $I$ is $\ge 3$, it follows that for each triple of indices $i-1, i, i+1$ the point $x_i$ is within uniform distance $D=D(K)$ from both $Q_{i-1,i}$ and $Q_{i,i+1}$: There exist $y^-_i\in Q_{i-1,i}$, $y^+_i\in Q_{i,i+1}$ such that
$$
d(x_i, y^\pm_i)\le K. 
$$
Such a point $x_i$ might not even exist, which would mean that each $\kappa$-qi section $\gamma$ of $\X$ is defined only on an interval of length $\le 2$, and that would definitely ensure acylindricity of $\X$. In general, one can say that 
$$
d(P_{Q_i, Q_{i-1,i}}(y_i^+), x^-_i)\le 2K, \quad d(P_{Q_i, Q_{i,i+1}}(y_i^-), x^+_i)\le 2K,
$$
and, hence, 
$$
d(x_i, x_i^\pm)\le 3K. 
$$
It follows that any two $\kappa$-qi sections $\gamma_0, \gamma_1$ defined on $I$ satisfy
$$
d(\gamma_0(i), \gamma_1(i))\le 6K, i\in V(I),
$$
thereby ensuring $(6K,\kappa,3)$-acylindricity of $\X$.

\subsection{Hyperbolicity of finite trees of hyperbolic spaces}

We will also need the following version of Theorem \ref{thm:chain} 
in the situation when the coboundedness condition is dropped:

\begin{cor}\label{cor:finite-tree-hyp} 
Suppose that $T$ is a finite tree, $\X=(\pi: X\to T)$ is a tree of hyperbolic spaces (satisfying Axiom {\bf H}). Then $X$ is $\delta$-hyperbolic with 
$$
\delta=\delta_{\ref{cor:finite-tree-hyp}}(\delta_0, L_0, |V(T)|),$$  
i.e. $\delta$ depends only on the parameters of $\X$ and the cardinality of $|V(T)|$. 
\end{cor}
\proof The proof is induction on $|V(T)|$. For $|V(T)|=1$, there is nothing to prove. For $n=2$, the corollary is nothing but Theorem \ref{thm:chain} for the quasiconvex amalgam of pairs. We, thus, assume that the corollary holds for all 
trees $S$ with $|V(S)|= n\ge 2$. Let  $\X=(\pi: X\to T)$ be a tree of hyperbolic spaces with $|V(T)|= n+1$. Pick a valence $1$ 
vertex $w$ of $T$ and let $e=[v,w]$ be the incidence edge. Set $Y_v:=X_{vw}$. We then form a new tree of spaces $\Y=(Y\to S)$, where $S$ is obtained from $T$ by removing $w$ and $e$, hence, $S$ has one less vertex than $T$. For vertices  of $S$ which are distinct
from $v$ and edges which are not incident to $v$,  we use the same incidence maps for $\Y$ as we had for $\X$. 
For edges $e_i=[v_i,v]$ incident to $v$ we use incidence maps 
$Y_{e_i}=X_{e_i}\to Y_v=X_{vw}$ equal to the corestrictions of the incidence maps $X_{e_i}\to X_v$.  The new tree of spaces still satisfies 
the assumptions of the corollary since $Y_v=X_{vw}$ is $\delta_1=\delta(\delta_0, L_0, 2)$-hyperbolic and  incidence maps 
$$
f_{e_iv}: X_{e_i}\to Y_v= X_{vw}
$$
are $L_1=L_0\cdot L'_0$-qi embeddings, where $L'_0=L_{\ref{lem:L'0}}(\delta,L_0)$. 
Now, $\delta=\delta(\delta_0, L_0, n)$-hyperbolicity of $X$ follows from the induction hypothesis. 
\qed 

\subsection{Secondary parameters of trees of hyperbolic spaces} 
\label{not:edge-space-constants}\label{not:K0}\label{not:secondary}

In addition to the primary parameters of trees of hyperbolic spaces $\X= (X\to T)$, we will be frequently using {\em secondary parameters}, which are functions of the primary parameters. Since these secondary parameters are used so often, we will give them 
special names. First of all, we recall some constants defined earlier, namely, $\la_0$, the quasiconvexity constant of the images $X_{ev}$ of incidence maps $X_e\to X_v$ and  $L'_0\ge 2$, the quasiisometry constant for the inclusion maps $X_v\to X_{vw}$, where $e=[v,w]$ runs over all edges of $T$ (Lemma \ref{lem:L'0}). Also, $\delta'_0$ is the supremum of hyperbolicity constants of the spaces 
$X_{uv}= X_{\llbracket u, v \rrbracket}$ (Corollary \ref{cor:edge-spaces}). Let $\la'_0$ denote an upper bound on  
the quasiconvexity constants for the images in $X_{vw}$ of $4\delta_0$-quasiconvex subsets in 
$X_v, X_e$ (in particular, each $X_v, X_e$ is $\la'_0$-quasiconvex in $X_{vw}$).  Explicitly, one can take 
$$
\la'_0= 92(L'_0)^2(L'_0 + 3\delta'_0).$$
We will also use the notation 
$L'_1$ for an upper bound of coarse Lipschitz constants of projections $P=P_{X_{uv},X_v}: X_{uv}\to X_v$, 
$L'_1=(L'_0+1)\cdot D_{\ref{lip-proj}}(\delta'_0,\la'_0)$ (Lemma \ref{lip-proj}). 
Last, but not least, we define the constant \index{$K_0$}
\begin{equation}\label{K0}
K_0:= (15(2\la'_0 + 5\delta'_0)L'_0)^3. 
\end{equation}
The importance of this constant will become clear in Chapter \ref{ch:4 classes} during the discussion of flows of quasiconvex 
subsets of vertex-spaces of $\X$.  This constant will be critical in computing the flaring constant $K_*$ in Theorem 
\ref{thm:mainBF}.  

\section{Flaring for semidirect products of groups}\label{sec:semi-direct-flaring} 

The purpose of this section is to illustrate the concept of  flaring in the case of semidirect products of groups, 
$G=H\rtimes \Z$.

Suppose $H$ is a nonelementary finitely generated group (which we will eventually assume to be hyperbolic) with a finite generating set $S$ and the 
corresponding word-metric $d_H$. Recall that the word-length of an 
element $h\in H$, denoted $|h|_H$ or $|h|_{S}$, when the generating set is to be
stressed, is related to $d_H$ by $|h|_H=d_H(1, h)$. 

Let  $f: H\map H$ be an automorphism and 
 $G=H\rtimes_f \langle t\rangle$  the corresponding semidirect product. Let $S_G= S\cup \{t\}$ be a generating set
of $G$, where $t$ is the stable letter corresponding to the infinite cyclic factor in the semidirect product. 
Let $X$ be the Cayley graph $\Gamma(G,S)$; define the linear tree $T=\Gamma(\Z, 1)$. 
Then we have a tree of metric spaces $\pi: X\map T$, where the vertex spaces are various left  
cosets $X_i:=t^iH, i\in \Z$, of $H$ in $G$. (Strictly speaking, $X$ is only quasiisometric to the 2-dimensional 
complex which is the total space of the abstract tree of spaces whose vertex spaces are isometric copies of the Cayley graph of $H$.) 
We shall denote the standard metric on $X$ by $d_G$ and the metrics on
the left cosets $t^iH\subset G$ by $d_{t^iH}$; the latter are isometric to the word-metric on $H$ corresponding to the generating set $S$.

Given $m\in \Z, n\in \N$, a $\kappa$-qi section over the interval $\llbracket m-n, m+n\rrbracket$ in $T$ is a sequence\footnote{For further computations, we find it notationally  convenient to write elements of $t^iH$ as $h_it^i, h_i\in H$, which is possible since $H$ is normal in $G$.} $\{h_i t^i\}$, $m\leq i\leq n$, such that for each $i\in [m-n, m+n]\cap \Z$, $d_{X_{i,i+1}}(h_it^i, h_{i+1}t^{i+1})\leq \kappa$, where we identify 
integers $i\in [m-n, m+n]$ with the corresponding vertices of $T$.  This inequality is satisfied, in particular, when $d_{X_i}(h_it^i, h_{i+1}t^{i})\leq \kappa-1$. Since the vertex spaces $X_i, X_{i+1}$ 
  are qi embedded in $X_{i,i+1}$, after changing $\kappa$ if necessary, we can (and will) identify $\kappa$-qi sections with 
 sequences $\{h_i t^i\}$ satisfying the inequality 
$$
d_{X_i}(h_it^i, h_{i+1}t^{i})=d_H(1, t^{-i}h^{-1}_i h_{i+1}t^{i})=d_H(1,f^{-i}(h^{-1}_i h_{i+1}))= 
|f^{-i}(h^{-1}_i h_{i+1})|_H\le \kappa,$$
equivalently,
\begin{equation}\label{eq:fkappa}
d_H(f^{-i}(h_i), f^{-i} (h_{i+1}))\le \kappa . 
\end{equation}

Here is an explicit example:

\begin{example}
Fix $h\in H$. Then $i\mapsto ht^i, i\in \Z$, is a $1$-qi section over $T$. 
\end{example}
 
Now, let us see what, respectively, exponential and proper flaring conditions in this context mean in group-theoretic terms.
Suppose $\gamma,\gamma'$ are two $\kappa$-qi sections over $\llbracket m-n,m+n\rrbracket$, 
where $m,n\in \Z$, given by maps $i\mapsto a_it^i$ and $i\mapsto b_it^i$. Then for
each integer $i\in [m-n, m+n]$, the fiber-distance equals 
$$
d_{t^iH}(\gamma(i),\gamma'(i))=d_{t^iH}(a_it^i, b_it^i)=d_H(t^{-i}a_it^i, t^{-i}b_it^i)=
d_H(1, t^{-i}a^{-1}_ib_it^i)= |f^{-i}(a^{-1}_ib_i)|_H.$$ 
If we denote the pair of sections $(\gamma, \gamma')$ by $\Pi$, then 
$$\Pi_{max}=\max \{|f^{-m+n}(a^{-1}_{m-n}b_{m-n})|_H, |f^{-m-n}(a^{-1}_{m+n}b_{m+n})|_H\}.$$
In the special case, when $m=0, n=1$,  
\begin{equation}\label{eq:Pimax-auto} 
\Pi_{max}=\max \{|f(a^{-1}_{-1}b_{-1})|_H, |f^{-1}(a^{-1}_{1}b_{1})|_H\}. 
\end{equation}

\begin{example}\label{simplest qi section}
If $\gamma, \gamma'$ are given by the maps $i\mapsto t^i$ and $i\mapsto ht^i$ respectively, where $h\in H$, then
$$\Pi_{max}=\max\{|f^{-m+n}(h)|_H, |f^{-m-n}(h)|_H\}.$$
\end{example}

Since $G$ acts on itself isometrically via the  left multiplication, in order to
formulate flaring conditions,  without loss of generality, we may assume that the qi sections $\gamma, \gamma'$ are
defined over intervals of the form $\llbracket -n,n\rrbracket$ (i.e. $m=0$) and $\gamma(0)=1$. 

\medskip
One can also reformulate the above conditions and quantities using the notion of  {\em pseudo-orbits} coming from the theory of dynamical systems. 

\begin{defn}\index{pseudo-orbit}\label{defn:pseudo-orbit}
Let $(Y,d)$ be a metric space and $\phi: (Y,d)\to (Y,d)$ be a homeomorphism. For a number $K$, a {\em $K$-pseudo-orbit} of $\phi$ in $Y$ is a biinfinite sequence $(y_i)_{i\in \Z}$ in $Y$ such that for each $i$
$$
d(y_{i+1}, \phi(y_i))\le K. 
$$
For instance, if $K=0$ then $0$-pseudo-orbits are just orbits of $\phi$ (or, more precisely, the cyclic group generated by $\phi$) in $Y$. The element $y_i$ is called the $i$-th member of the pseudo-orbit $(y_i)_{i\in \Z}$. 

A {\em partial} $K$-pseudo-orbit is the restriction of a $K$-pseudo-orbit sequence to an interval in $\Z$. 
\end{defn}

Given an automorphism $f$ of the group $H$, set $\phi:= f^{-1}$.  We will use $(H, d_H)$ as our metric space $(Y,d)$. 
For a section $\gamma$,  $\gamma(i)= h_it^i$, we define $g_i:= \phi^i(h_{i+1})$. In particular, $g_0=h_1$. 
Then the inequality \eqref{eq:fkappa} is equivalent to 
$$
d_H(\phi(g_i), g_{i+1})\le \kappa,  
$$
the $\kappa$-pseudo-orbit condition. In other words, instead of working with $\kappa$-sections, we can work with (partial) $\kappa$-pseudo-orbits. 
The case of a 1-qi section corresponds to the case when $(g_i)$ is the (partial) $\phi$-orbit of $g_0$. 
Given two sections $\ga, \ga'$, we note that the corresponding partial pseudo-orbit sequences $(g_i), (g'_i)$, satisfy  
$$
d_{t^iH}(\gamma(i),\gamma'(i))=d_H(\phi(g_i), \phi(g'_i)). 
$$
In particular, for fixed $\phi$, a uniform bound on $d_{t^iH}(\gamma(i),\gamma'(i))$ is equivalent to a (possibly different) uniform bound on $d_H(\phi(g_i), \phi(g'_i))$. 

\medskip

We can now restate various flaring conditions:

\medskip
(1) The linear (Bestvina--Feighn) $\kappa$-flaring condition is equivalent to: 

There exist constants $M_\kappa\ge 0$, $\la_\kappa>1$ and $n_\kappa\in \N$ such that 
for every pair of  maps $i\mapsto a_i\in H$ and $i\mapsto b_i\in H$, $i\in [-n_\kappa, n_\kappa]\cap \Z$, satisfying:

(a) $a_0=1$, $|b_0|_H\ge M_\kappa$,

(b) $d_H(f^{-i}(a_i), f^{-i}(a_{i+1}))\le \kappa$, $d_H(f^{-i}(b_i), f^{-i}(b_{i+1}))\le \kappa$, $i\in [-n,n]$, 

we have 
$$
\max \{d_H(f^{n}(a_{-n}), f^{n} (b_{-n})),  d_H(f^{-n}(a_{n}), f^{-n} (b_{n}))\} \ge \la |b_0|_H. 
$$

\medskip 
(2) The proper $\kappa$-flaring condition is equivalent to:

There exists a constant $M_\kappa\ge 0$ and a proper function $\phi_\kappa: \N\to \R_+$ such that 
for every pair of  maps $i\mapsto a_i\in H$ and $i\mapsto b_i\in H$, $i\in \Z$,  satisfying:

(a) $a_0=1$, $|b_0|_H\ge M_\kappa$,

(b) $d_H(f^{-i}(a_i), f^{-i}(a_{i+1}))\le \kappa$, $d_H(f^{-i}(b_i), f^{-i}(b_{i+1}))\le \kappa$, $i\in [-n,n]$, 

we have 
$$
\max \{d_H(f^{n}(a_{-n}), f^{n} (b_{-n})),  d_H(f^{-n}(a_{n}), f^{-n} (b_{n}))\} \ge \phi_\kappa(n). 
$$

It is also useful to spell out the negation of the proper $\kappa$-flaring condition, which is most apparent as the negation of the bigon property from Corollary \ref{cor:super-weak flaring}: 

There exists elements $g, h\in H$ and pairs sequences  of partial 
$\kappa$-pseudo-orbits $(g_{i,n})_{n\in \N}$, $(g'_{i,n})_{n\in \N}$ of $f$ in $H$ defined for $i\in [0, N_n]$,  
such that:

(a)  For all $n$, $g_{0,n}=1, g'_{N_n,n}=g$, $g'_{N_n}=h g_{N_n}$.

(b) $\lim_{n\to \infty} \max_{i\in [0, N_n]} d_H(g_{i,n}, g'_{i,n})= \infty$.   

Note that, by possibly increasing $\kappa$ to $K:=\kappa+C$, where $C=\max\{ |g|, |h|\}$, and working with partial $K$-pseudo-orbits, we can even ensure that $g=1, h=1$ and, hence,  
$g'_{0,n}=1$, $g_{N_n,n}=g'_{N_n,n}$. 
  
  \medskip

\begin{comment}
\medskip 
(3) The uniform flaring ....

\medskip
(4) Maybe also the bigon property from Corollary \ref{cor:super-weak flaring}. It also might be convenient to state the negation of the thin bigon property:
\end{comment}

We next relate flaring to {\em hyperbolicity properties} of the automorphism $f$.

\begin{defn}(Bestvina, Feighn, \cite{BF})\label{defn hyp autom}
Suppose $H$ is a finitely-generated group and $S$ is finite generating set for $H$.
Suppose $f: H\map H$ is an automorphism. We say that $f$ is {\em weakly hyperbolic}
if there is $m\in \NN$, $\lambda>1$ and a finite subset $E\subset H$, such that
for all $h\in H\setminus E$ 
we have $$\lambda |h|\leq \max \{ |f^m(h)|, |f^{-m}(h)|\}.$$
We say that $E$ is the {\em exceptional} subset. 
An automorphism is called {\em hyperbolic} if the above inequality holds with $E=\emptyset$. 
\end{defn}

Some remarks are in order regarding this definition. 

\begin{rem}
\begin{enumerate}
\item  The notion of {\em hyperbolicity} of an automorphism   was introduced by 
Bestvina and Feighn in \cite{BF} (the exceptional subset $E$ was absent). They also proved hyperbolicity of semidirect products $H\rtimes_f \langle t\rangle$ of hyperbolic automorphisms 
of hyperbolic groups, see Corollary in \cite[section 5]{BF}.  
However, the original Bestvina--Feighn  definition  is 
too restrictive if the purpose is to conclude hyperbolicity of semidirect products. Lemma \ref{lem:BF-auto} below 
(which is already present in  Gersten's paper \cite[Corollary 6.9]{MR1650363})  
shows that hyperbolicity 
of the semidirect product is equivalent to the weak hyperbolicity of the automorphism. 

\item If $f: H\to H$ is a weakly 
hyperbolic automorphism, then for any nontrivial finite group $H_1$,
the automorphism $f'=f \times \id$ of $H'=H\times H_1$ is also weakly hyperbolic but fixes the subgroup $H_1$ element-wise and, hence, is not hyperbolic. 

\item In Corollary \ref{cor:wh-auto} we will prove that for automorphisms of torsion-free hyperbolic groups weak hyperbolicity is equivalent to hyperbolicity. 

\item  The only hyperbolic groups which admit weakly hyperbolic automorphisms are the ones commensurable to free products of surface groups and free groups, as 
follows for instance from \cite{rips-sela}. 

\item The hyperbolicity inequality trivially holds for the trivial element $h=1\in H$. 
Suppose that the exceptional subset $E$ is a ball $B(1, r)\subset H$ and that (with $\la>1$) for $h\notin E$, 
$|f^m(h)|\ge \la |h|$. Then $f^m(h)\notin E$ and, thus, we can apply the same hyperbolicity inequality to $f^m(h)$. Clearly,
$$
 |f^{2m}(h)|= \max \{ |f^{2m}(h)|, |f^{-m+m}(h)|\}\ge \la^2|h|.  
 $$
Repeating this argument, we see that for each $i\ge 1$, 
$$
|f^{im}(h)|\ge \la^i |h|. 
$$
\end{enumerate}
\end{rem}

\begin{lemma} \label{lem:hyp-independent} 
Hyperbolicity and weak hyperbolicity of an automorphism $f: H\map H$ are 
 independent of the finite generating set of $H$. 
\end{lemma}
\proof We will verify the claim for the weak hyperbolicity property since the proof for hyperbolicity is identical (with $E=\emptyset$). 

Suppose $f$ is weakly hyperbolic with respect to a finite generating set $S$, i.e. there exist $m\in \NN$, $\lambda>1$  and a finite subset $E\subset H$ such that 
$$
\lambda |h|_S\leq \max \{ |f^m(h)|_S, |f^{-m}(h)|_S\}$$ 
for all $h\in H\setminus E$. 

 Suppose $S'$ is another finite generating 
set of $H$. For any $h\in H$ let $|h|_S$ and $|h|_{S'}$ denote the word-lengths of $h$ with respect to 
$S$ and $S'$ respectively. Then there is a constant $C>0$ such that 
$\frac{1}{C}|h|_S\leq |h|_{S'}\leq C|h|_S$ for all $h\in H$. Also, we note that for all
$r\in \NN$ we have $\lambda^r |h|_S\leq \max\{|f^{mr}(h)|_S, |f^{-mr}(h)|_S\}$ for all 
$h\in G\setminus E_r$ where 
$$E_r=\bigcup_{-(r-1)\leq i\leq r-1} f^{im}(E).$$
Hence, 
$$\lambda^r |h|_{S'}\leq C\lambda^r |h|_S \leq C\max\{|f^{mr}(h)|_S, |f^{-mr}(h)|_S\}
\leq C^2 \max\{|f^{mr}(h)|_{S'}, |f^{-mr}(h)|_{S'}\}$$ for $h\in H\setminus E_r$.
It follows that $C^{-2}\lambda^r |h|_{S'}\leq \max\{|f^{mr}(h)|_{S'}, |f^{-mr}(h)|_{S'}\}$
for $h\in H\setminus E_r$.
Thus, if we choose $r=r_1$ large enough, we get $\lambda_1=C^{-2}\lambda^{r_1}>1$ and setting  
$m_1=r_1m$, we obtain  
$$\lambda_1 |h|_{S'}\leq \max \{ |f^{m_1}(h)|_{S'}, |f^{-m_1}(h)|_{S'}\}$$
for all $h\in H\setminus E_{r_1}$.  Whence $f$ is weakly hyperbolic with respect to the generating set $S'$. \qed

\begin{comment}
{\mini 
Thus, one is led to the following:

\begin{question} 
Is the hyperbolicity of $f$ (in the sense of Bestvina--Feighn) equivalent to the weak hyperbolicity 
for torsion-free hyperbolic groups? 
\end{question}

The following results provide  some partial evidence to the positive answer.
}
\end{comment}

\begin{example} \label{lemma: hyp autom}
(1) Suppose  $H=\pi_1(\Sigma)$ where $\Sigma$ is a closed connected hyperbolic surface. 
Then an automorphism $f$ of $H$ is hyperbolic if and only if it is induced by a pseudo-Anosov automorphism of $\Sigma$, if and only if it has no nontrivial periodic conjugacy classes, if and only if the semidirect product $H\rtimes_f \Z$ is hyperbolic.  (The equivalence of the last three properties is due to William Thurston, see e.g. \cite{Casson, Otal-book}. The equivalence with hyperbolicity of $f$ can be seeing either as a consequence of the pseudo-Anosov property or 
of the combination of Lemma \ref{lem:BF-auto} and Corollary \ref{cor:wh-auto}.)

(2) If $H={\mathbb F}_n$, $n\geq 2$, then any automorphism $f\in Aut(H)$ with no periodic conjugacy classes is hyperbolic
(in the sense of Bestvina--Feighn). See  
Theorem 5.1 in \cite{BFH-lam}.  
\end{example}

We refer the reader to \cite{MR1800064, MR3720349, MR4237417} for other results in this direction.

\begin{lemma}[\cite{BF, MR1650363}] 
 \label{lem:BF-auto} If $f$ is weakly hyperbolic, then the tree of metric spaces 
$X=\Gamma(H\rtimes_f \Z, S_G)\map T=\Gamma(Z,1)$--- as constructed in the
beginning of this subsection--- satisfies the Bestvina--Feighn  flaring condition. The converse is also true. 
\end{lemma}
\proof  1. Suppose $f$ is weakly hyperbolic. 
Let $R=\max\{d_H(1,h): h\in E\}$ where $E\subset H$ is a finite exceptional subset as in Definition \ref{defn hyp autom}. Then for all $x\in H$ with $|x|_H> R$ we have
$$\lambda |x|\leq \max\{ |f^m(x)|, |f^{-m}(x)|\}.$$ 
First of all, since the Bestvina--Feighn  flaring condition is equivalent to hyperbolicity of $G$ and the latter is equivalent to the hyperbolicity of the semidirect product 
$H\rtimes_{f^m} \Z$ (commensurable to the original group $G$), it suffices to consider the case when $m=1$. We will also assume that in the definition of a weakly hyperbolic automorphism the maximum is attained by $\phi=f^{-1}$ rather than $f$ (otherwise, we replace $f$ with $f^{-1}$). Then the weak hyperbolicity inequality reads 
\begin{equation}\label{eq:weak}
\lambda |x|\leq  |\phi(x)|
\end{equation} 
for all $x\in H\setminus E$. 

Take $\kappa \geq 1$. 
As we noted earlier, 
it suffices to verify the Bestvina--Feighn flaring condition for pairs of $\kappa$-qi sections $\Pi=(\gamma, \gamma')$ 
 defined over intervals of the form $\llbracket -m,m\rrbracket$, satisfying  $\gamma(0)=1, h:=\gamma'(0)$, where, as above, $m=1$, such that $\Pi_0=|h|\ge M_\kappa$ for a suitable uniform constant $M_\kappa$. Pick any $\la'$ in the open interval $(1,\la)$; for concreteness, we take $\la'= \frac{1}{2}(\la+1)$. We claim that 
 $$
\la' |h|= d_H(\gamma(0), \gamma'(0))\leq d_{tH}(\gamma(m), \gamma'(m)).  
 $$
Set $\gamma(1)=h_{1}t, \gamma'(1)=h'_{1}t$. Then (as we noted earlier, after changing $\kappa$) the $\kappa$-qi section condition for $\gamma, \gamma'$ 
 over the interval $[-1,0]$ is equivalent to the inequalities 
 \begin{equation}\label{eq:kappa1}
 d_H(\phi(h_{1}),1)\le \kappa, \quad d_H(\phi(h), \phi(h'_{1}))\le \kappa.
 \end{equation} 
 We will estimate from below the distance (in $H$) between $\gamma(1), \gamma'(1)$; according to the computation in \eqref{eq:Pimax-auto}, we need to estimate from below the distance 
 $$
 d_H(\phi(h_{1}), \phi(h'_{1})). 
 $$
 By the triangle inequality, combined with the inequalities \eqref{eq:weak} and \eqref{eq:kappa1}, we have
 $$
  d_H(\phi(h_{1}), \phi(h'_{1})) \ge d_H(1, \phi(h)) -2\kappa \ge \la |h| -2\kappa. 
 $$
 Then the desired inequality $\la' |h| \leq d_{tH}(\gamma(1), \gamma'(1))$ is equivalent to 
 $$
 (\la-\la')|h|\ge 2\kappa 
 $$
 By out choice of $\la'$, this amounts to 
 $$
|h| \ge \frac{2\kappa}{\la - 1}. 
$$
Therefore, taking $M_\kappa$ equal to the maximum of $\frac{2\kappa}{\la - 1}$ and  $1+\diam_H(E)$, we ensure the Bestvina--Feighn flaring condition.  

\medskip 
2. For the converse one applies the flaring condition to $\kappa=1$-qi sections. More precisely, we work with pairs of sections
as in Example \ref{simplest qi section}. We take the finite set in the definition of weak hyperbolicity 
to be $E=\{h\in H: |h|_H\leq M_1\}$. The rest is straightforward and hence left as an exercise 
for the reader. \qed

\begin{lemma}
Suppose that $H$ is a hyperbolic group, $f: H\to H$ is weakly hyperbolic. Then the exceptional subset of $E$ can be chosen to contain only finite order elements. 
\end{lemma}
\proof Let $E\subset H$ be an exceptional subset of $f$, i.e. there exist $m\in \NN$, $\lambda>1$  such that 
$$
\lambda |h|\leq \max \{ |f^m(h)|, |f^{-m}(h)|\}$$ 
for all $h\in H\setminus E$. After replacing $f$ with $f^m$, we can assume that $\lambda |h|\leq \max \{ |f(h)|, |f^{-1}(h)|\}$ unless $h\in E$. 
After enlarging $E$ is necessary, we can assume that it equals to the ball of certain radius $r$ in $H$ (centered at $1\in H$). 
We define $E'\subset E$, the subset consisting of infinite order elements. 

Suppose that $h\in E$ is such that for infinitely many values of $m\ge 1$,
$$
|f^m(h)|\le r.   
$$
We claim that $h$ has finite order. Indeed, then there exist two numbers $n> m\ge 1$ such that 
$$
f^m(h)=f^n(h), f^{n-m}(h)=h.
$$
It follows that in the group $G=H\rtimes_f \Z$ we have 
$$
t^{n-m} h= h t^{n-m}. 
$$
Since $f$ is weakly hyperbolic, the semidirect product $G$ is a hyperbolic group, see Lemma \ref{lem:BF-auto}. Hence, the abelian subgroup $A< G$ generated by $t^{n-m}$ and $h$ is virtually cyclic, i.e. $h$ has finite order. 

Thus, as noted above, for each $h\in E'$ there exists a smallest natural number 
$n=n_h$ such that $|f^n(h)|>r$, which, in particular, implies that $f^n(h)\notin E'$. Thus, 
$$
\max \{|f\circ f^n(h)|, |f^{-1}\circ f^n(h)|\}\ge \la |f^n(h)|\ge \la |h|.
$$
Since $n_h$ was chosen to be smallest, it follows that the above inequality holds for $f\circ f^n(h)=f^{n+1}(h)$. By the same argument as in the proof of Lemma \ref{lem:exp}, we see 
that  for each $m\ge \max\{n_h: h\in E'\}$, and $h\in E'$, 
$$
\lambda |h|\leq  |f^m(h)|.$$ 
Since the hyperbolicity inequality holds for all $h\in H\setminus E$, we conclude that $f$ satisfies the hyperbolicity condition except for the subset $E''\subset E$ consisting of torsion elements. \qed

\begin{cor}\label{cor:wh-auto}
If $H$ is torsion-free then each weakly hyperbolic automorphism of $H$ is hyperbolic. 
\end{cor}

\chapter{Flow-spaces, ladders and their retractions} \label{ch:4 classes} 

In this chapter we introduce and analyze four classes of subtrees of spaces in hyperbolic trees of spaces:

\begin{itemize}
\item Ladders

\item Metric bundles 

\item Carpets

\item Flow-spaces
\end{itemize}

These spaces play key role in proving hyperbolicity and describing geodesics in trees of spaces $(\pi: X\to T)$: Uniform quasigeodesics in $X$ will be inductively described as concatenations of uniform quasigeodesics in carpets, ladders and flow-spaces. The main result of this and the next chapter is that all ladders, carpets and flow-spaces are hyperbolic and admit coarse Lipchitz retractions from $X$. We note that our definitions of ladders and flow-spaces are 
inspired by the {\em ladder} construction of Mitra, \cite{mitra-trees}, while the notion of metric bundles is adapted from \cite{pranab-mahan}.

\section{Semicontinuous families of spaces}

All four classes of spaces discussed in this (and the next) chapter are special cases of {\em semicontinuous families} 
of spaces, which are certain subtrees of spaces 
${\mathfrak Y}\subset {\mathfrak X}$. 
In what follows, given a subtree of spaces 
$\Y= (\pi: Y\to S)\subset \X=(\pi: X\to T)$, it will be notationally convenient to extend $\Y$ to a tree of spaces (still denoted $\Y$) 
over the entire tree $T$ by declaring $Y_v=\emptyset, Y_e=\emptyset$ for each $v\in V(T)-V(S)$ and $e\in E(T)-E(S)$. 

\begin{defn}\label{defn:scfamily} \index{semicontinuous family} 
Suppose that  ${\mathfrak X}=(\pi: X\to T)$ is a tree of hyperbolic spaces. Fix constants 
$\la\in [0,\infty)$, $E, K\in [1, \infty), D\in [0,\infty]$. 
We will say that a subtree of spaces  ${\mathfrak Y}= (\pi: Y\to S), S\subset T$, in ${\mathfrak X}$ forms a {\em $(K,D,E,\la)$-semicontinuous family}, 
relative to a vertex $u\in V(S)$, called the {\em center} of $\Y$, if the following conditions hold:

1. Each vertex/edge space $Y_v\subset X_v, Y_e\subset X_e$, $v\in V(S), e\in E(S)$, is $\la$-quasiconvex. 

2. Each $y\in \YY$ 
is connected to $Y_u$ by a $K$-leaf $\ga_y$ in $Y$.  

3. For each edge $e=[v,w]\in E(T)$ 
we define the (possibly empty!) projection
\begin{equation}\label{Yvw}
Y^v_w:=P_{X_{vw}, X_w}(Y_v)\subset X_w.
\end{equation}
We require that whenever $e=[v,w]\in E(S)$ is oriented away from $u$, 
\begin{equation}\label{eq:E-in} 
\Hd_{X_{vw}}(Y^v_w, Y_w)\le E
\end{equation}
 and 
 \begin{equation}\label{eq:K-in} 
 \Hd_{X_{vw}}(Y_w, Y_e)\le K
 \end{equation}  

4. For every edge $e=[v,w]$ such that $v\in V(S), w\notin V(S)$, we require the pair 
of quasiconvex subsets $(Y_v, X_w)$ in $X_{vw}$ to be $D$-cobounded. 
\end{defn}

\begin{figure}[tbh]
\centering
\includegraphics[width=80mm]{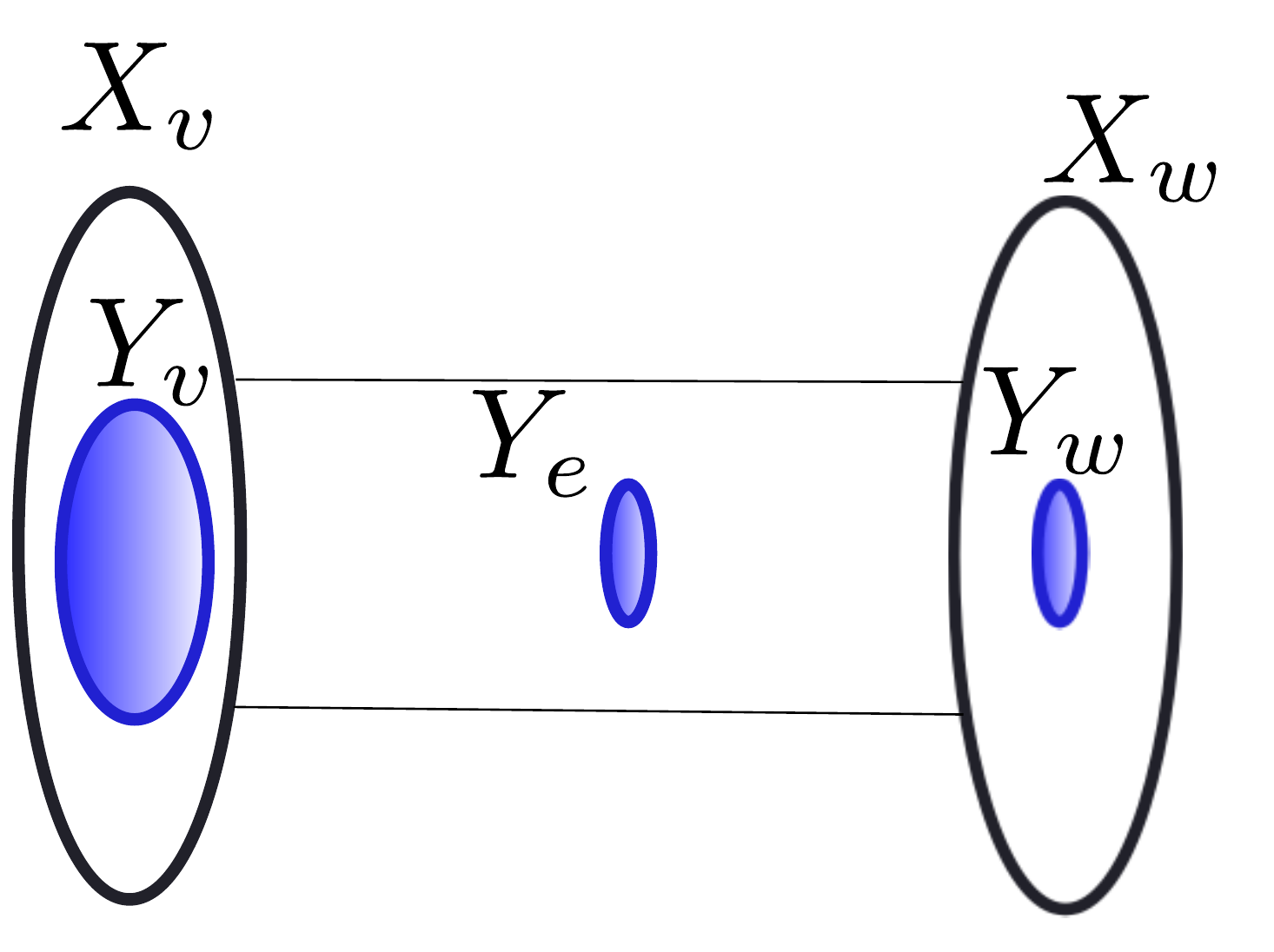}  
\caption{A semicontinuous family}
\label{scf.fig}
\end{figure}

\begin{rem}\label{rem:semico}
i. Condition 2 implies that $Y_w\subset N^e_{K}(Y_v)$, hence, $Y_w\subset N_{\eta_0(2K)}(Y^v_w)$ with respect to the metric of $X_w$. (Recall that $\eta_0$ is the distortion function of $X_v$ in $X$, hence, in $X_{uv}$.)

ii. Condition 2 in this definition ensures that $Y_w$ cannot be ``much larger'' than $Y_v$, while Condition 3 ensures a certain lower bound on $Y_w$.  Thus, as we move away from $X_u$, vertex spaces of $\Y$ can shrink substantially (even disappear) but they cannot substantially increase. 

iii. In most examples in our book, $\la=4\delta_0$, hence, we will be suppressing the dependence on this parameter and record only the triple of numbers $(K,D,E)$.

iv. To ensure uniform coboundedness of the pairs $(Y_v, X_w)$ in Axiom 4, it suffices to 
get a uniform upper bound $C$ on the diameters of $Y^v_w$: It will then follow that the pair $(Y_v, X_w)$ 
is $D$-cobounded for some $D=D(\delta'_0,\la'_0, C)$, see Corollary \ref{cor:cob-char}. 

v. We do not insist on the converse implication in Axiom 4: There will be important situations when we have to consider $\Y$'s with uniformly bounded fibers over non-boundary vertices of $S$. 

vi. Axioms 3 and 4 will be needed in order to have a uniform coarse Lipschitz retraction $X\to Y$, see 
Theorem \ref{thm:mitras-projection}. 

vii. The edge-spaces $Y_e$ of  are largely irrelevant for our discussion. 

viii. The projections $P_{X_{vw},X_w}$ restricted to $X_v$ are at uniformly bounded distance (as measured  in $X_{vw}$) from the projections $P_{X_v,X_{ev}}$. The same, of course, applies to restrictions of the projections $P_{X_{vw},X_v}$. However, we decided to work with the projections $P_{X_{vw},X_w}$ as computations tend to be simpler in this setting. 
\end{rem}

\begin{theorem}\label{thm:mitras-projection}
Suppose that $\Y$ is a $(K,D,E,\la)$-semicontinuous family of spaces with $D<\infty$. Then there exists an 
$L_{\ref{thm:mitras-projection}}(K,D,E,\la)$-coarse Lipschitz retraction $\rho_{\Y}: X\to Y$. 
\end{theorem}
\proof We will verify that the subtree of spaces $\Y$  is (uniformly) retractible;  we use Theorem \ref{thm:left-inverse} as follows. 

(i) We let $h'_v: X_v\to Y_v$ denote the restriction of the nearest-point projections $P_{X_{vw},Y_v}$. According to  Theorem \ref{thm:left-inverse}(i), we need to bound the diameter of the image of $X_{ev}$ in $Y_v$ under $h'_v$. However,  $X_{ev}$ is contained in the unit neighborhood of $X_w$ (taken in $X_{vw}$);
thus, we need to bound the diameter of the projection  (in $X_{vw}$)  of $X_w$ to $Y_v$, which is 
the content of Axiom 4 of the definition of a semicontinuous family of spaces: This diameter is $\le D$. 

\medskip 
(ii) Assuming that $e=[v,w]$ is an edge of $S$ oriented away from $u$, we need to get a uniform bound 
$$
\dist_{X'_v}( h'_v \circ f_{ev}, f'_{ev}\circ h'_e)\le Const. 
$$

The subsets $Y_{ev}, Y_e \subset X_{vw}$ are within unit Hausdorff distance, while 
 $$
 \Hd_{X_{vw}}(Y_e,Y_{w})\le E \hbox{~~and~~} 
 \Hd_{X_{vw}}(Y^v_w, Y_w)\le E.$$
  Since $d_{X_{vw}}(Y_w, X_v)\le K$, by applying 
  Lemma \ref{lemma0-flow-space} to the subsets $U_1=X_v, U_2=Y_w$ in $X_{vw}$, we conclude that 
  $$
  \Hd_{X_{vw}}(P_{X_{vw},Y_v}(X_w), P_{X_{vw},X_w}(Y_v))\le 2\la_0'+3\delta'_0 + K. 
  $$
Taking into account that the projection $P_{X_{vw},Y_v}$ is uniformly coarse Lipschitz, 
we conclude that $Y_{ev}$ 
 is uniformly close to the image of $X_{e}$ under the nearest-point projection 
$P_{X_{vw},Y_v}$ and, accordingly, the nearest-point projection $h'_e: X_e\to Y_e$ is uniformly close to the restriction of the nearest-point projection $X_{vw}\to P_{X_v,Y_v}(X_e)$ (see Lemma \ref{lem:two-projections}). 
Taking also into account that the map $f'_{ev}$ moves points by distance $\le K$ in $X_{vw}$, we 
 can replace $f'_{ev}\circ h'_e$ with the restriction of the projection $X_{vw}\to P_{X_v,Y_v}(X_e)$ 
 to $X_e$. Similarly, the map $f_{ev}$ moves points distance $\le 1$  in $X_{vw}$ and, hence (in view of the uniform coarse Lipschitz property of $h'_v$), we can replace the composition $h'_v \circ f_{ev}$ with the restriction 
 of the nearest-point projection $P_{X_{vw},Y_v}$ to $X_e$. But now, the projections 
$X_{vw}\to P_{X_v,Y_v}(X_e)$ and $P_{X_{vw},Y_v}$ are uniformly close to each other according to 
Corollary \ref{cor:projection-2} applied to the $\la_0'$-quasiconvex subsets $Y_v$ and $X_e$
in the ambient hyperbolic space $Z=X_{vw}$. \qed

\begin{cor}\label{cor:mitras-projection}
If $\Y=(\pi: Y\to S)$ is a $(K,D,E,\la)$-semicontinuous family of spaces with $D<\infty$, then the inclusion map 
$Y\to X$ is an $L_{\ref{thm:mitras-projection}}(K,D,E,\la)$-qi embedding. 
\end{cor}

\medskip 
We next describe a class of semicontinuous subtrees of spaces, called {\em metric bundles}.  
The theory of metric bundles was developed in \cite{pranab-mahan} in a more general setting when the base is allowed to be an arbitrary geodesic metric space 
but we will not need that in our book. The following definition of {\em metric bundles} is adapted from \cite{pranab-mahan} in a form suitable for our purposes. It is easy to verify
that the two definitions (ours and that of \cite{pranab-mahan}) 
are equivalent when the base is a tree. The reader should also compare this definition with 
the notion of an abstract metric bundle given in Definition \ref{defn:bundle}: 
Each metric bundle defined below is also an abstract metric bundle.

\begin{defn}[Metric bundles] \label{defn:flow0} \index{metric bundle} 
 A subtree of spaces $\Y =(\pi: Y\to S) \subset \X=(\pi: X\to T)$ is called a {\em $K$-metric bundle} if:
 
 1. Every vertex/edge space of $\Y$ is $\la$-quasiconvex in the respective vertex/edge space of $\X$. 
 
 2. For every  vertex $u\in V(S)$ and edge $e=[v,w]\in E(S)$ (directed away from $u$),  and $x\in Y_w, x\in Y_e$, 
 there exist $K$-qi sections $\ga_x$ in $\Y$ on $\llbracket w, u\rrbracket$, such that $\ga_x(w)=x$. 
 \end{defn} 

It follows immediately 
that each $K$-metric bundle forms  a $(K,\infty,E,\la)$-semicontinuous family of spaces in $\X$ (relative to any vertex $u\in S$), with $E=\eta_0(2K)$. 
The reader uncomfortable with using $D=\infty$ here can simply restrict $\X$ to $S$, then one can take $D=0$.

\medskip
While Theorem \ref{thm:mitras-projection} does not directly apply to metric bundles $\Y\subset X$ (unless $S=T$), we will see in Corollary \ref{cor:X-to-bundle} that under certain extra assumption weakening condition 4 in Definition 
\ref{defn:scfamily}, these too admit uniform coarse Lipschitz retraction from $\X$.

\section{Ladders} \label{sec:Ladders} 

\index{ladder}

Ladders are certain semicontinuous subtrees of spaces  $\L=(\pi: L\to S)\subset {\mathfrak X}$ whose fibers are geodesic segments. 
However, in addition to the properties of a semicontinuous family of spaces, we will impose 
a certain extra structure on $\L$.

\medskip 
Each ladder $\L=(\pi: L\to S)$ comes equipped with certain parameters and  two pieces of extra data: An orientation of the fibers (hence, ladders generalize oriented line bundles) and a canonical choice of a  maximal $K$-qi section $\Sigma_x\subset L$  through each point $x\in L$. The choice of $\Sigma_\bullet$ can be regarded as a ``connection'' on $\L$. Thus, ladders can be regarded as ``oriented line semi-bundles equipped with connections.''  

We will be primarily interested in ladders such that $\LL$ is contained in the $4\delta_0$-fiberwise neighborhood of a $k$-flow space ${\mathcal Fl}_k(Q_u)\subset X$ (these will be defined in Section 
\ref{sec:flow-spaces}). For the ease of notation, we will be ignoring the flow-spaces for now; formally speaking, one can regard the flow-space ${\mathfrak Fl}_k(Q_u)$ as a tree of hyperbolic spaces satisfying a uniform flaring condition.

\medskip 
We now begin with an axiomatic definition of ladders. Let ${\mathfrak X}=(\pi: X\to T)$ be a tree of hyperbolic spaces satisfying Axiom {\bf H}. Fix  positive numbers $K, D, E$ and a vertex $u\in T$; 
these are the {\em parameters} of a ladder $\L$.   
A {\em ladder}  with these parameters (a {\em $(K,D,E)$-ladder centered at $u$}) is a subtree of oriented geodesic intervals in $\X$, 
$\L= (\pi: L\to S)$, $S=\pi(L)$ which satisfies further axioms described below. Each fiber $L_v:= L\cap X_v, v\in V(T), L_e= L\cap X_e, e\in E(T)$, of $\L$ 
is an  {\em oriented} geodesic segment denoted   $[x_v y_v]_{X_v}$ or $[x_e y_e]_{X_e}$.

Furthermore, we fix once and for all a family $\Sigma_\bullet$ of {\em maximal partial $K$-qi sections} $\Sigma_x$ of $\pi: L\to S$, 
whose domains $T_x$ are subtrees in $S$ containing the vertex $u$ (and $\pi(x)$, of course). Maximality here is understood in the sense that if $\Sigma'_x$ is another partial $K$-qi section containing $\Sigma_x$, then $\Sigma_x=\Sigma'_x$.  The subscript $x$ in $\Sigma_x$ indicates that $x\in \Sigma_x$. We will assume that $x$ belongs to a vertex-space 
of $\X$.

\medskip 
{\bf Axiom L0}.  
We will require the family of 
sections $\Sigma_\bullet$ to be {\em consistent} in the sense that whenever $v=\pi(y)$ is between $u$ and $w=\pi(x)$, the sections $\Sigma_y$ and $\Sigma_x$ agree on the interval $\llbracket u, v\rrbracket\subset T$. 

\begin{defn}
Let $\L$ be a ladder centered at the vertex $u$, $L_u= [x_u y_u]_{X_u}$. We 
 will refer to the subsets  $\Sigma_{x_u}=bot(\L), \Sigma_{y_u}=top(\L)$ as, respectively, the {\em bottom} and the {\em top} of the ladder $\L$. 
\end{defn}

Thus, $\Sigma_\bullet$ defines a family of maps
$$
\Pi_{w,v}: L_w\to L_v,
$$
for every vertex $v\in V(S)$ between $u$ and $w\in V(S)$: 
$$
\Pi_{w,v}(x)= \Sigma_x\cap L_v, x\in \LL_w. 
$$
(Note that this intersection is nonempty since $\pi(\Sigma_x)$ is a subtree containing both $w$ and $u$, hence, also containing $v$.) 
These maps can be regarded analogues of {\em parallel transport maps} in the conventional theory of connections on bundles. 
Consistency of sections implies the following {\em semigroup property}: 
$$
\Pi_{w_1,w_3}= \Pi_{w_2,w_3}\circ \Pi_{w_1,w_2}  
$$
whenever $w_2, w_3$ belong to the interval $\llbracket u,w_1\rrbracket$ and appear there in the following order:
$$
u\le w_3\le w_2\le w_1. 
$$
As we will see below, axioms of a ladder require each map $\Pi_{w,v}$ to be either constant or an orientation-preserving topological embedding. The maps $\Pi_{w,v}$ need not be surjective; 
for every oriented edge $e=[v,w]$ in $S$ (oriented away from $u$) we have an oriented subsegment
$$
L'_v= [x'_v y'_v]_{X_v}:= \Pi_{w,v}(L_w)
$$
in $L_v$. Here $x'_v=\Pi_{w,v}(x_w), y'_v:= \Pi_{w,v}(y_w)$. The orientation of the segment  
$L'_v$ is then consistent with that of $L_v$ (since $\Pi_{w,v}$ is orientation-preserving).
Observe also that the mere existence of $K$-qi sections $\Sigma_x$ 
implies {\em some  semicontinuity} of the ladder $\L$: For every edge $e=[v,w]\subset S$ (oriented away from $u$) 
\begin{equation}\label{eq:semicontinuity} 
L_w \subset N^e_K(L'_v)\subset N^e_K(L_v), 
\end{equation}
 where the $K$-neighborhood is taken in the subspace $X_{vw}$ (which is what the superscript $e$ indicates). 
 However, $L_w$ can be much smaller than $L_v$.

\medskip 
For $\L$ (equipped with $\Sigma_\bullet$) to be a ladder, it has  to satisfy three further axioms listed below. Note, however, that the assumption that the fibers $L_v, L_e$ of $\L$ are geodesic segments ensures Property  1 in Definition \ref{defn:scfamily} with $\la=\delta_0$, while Axiom L1 implies Property 2  in that definition, thus making Axiom L3 somewhat redundant. 

We now fix $K\in [1,\infty), E\ge 1$ and  $D\in [0,\infty]$. While all other two parameters 
in the triple $(K,D,E)$ are real numbers, as with general semicontinuous subtrees of spaces, 
it is convenient to allow for infinite values of the parameter $D$.

\begin{figure}[tbh]
\centering
\includegraphics[width=60mm]{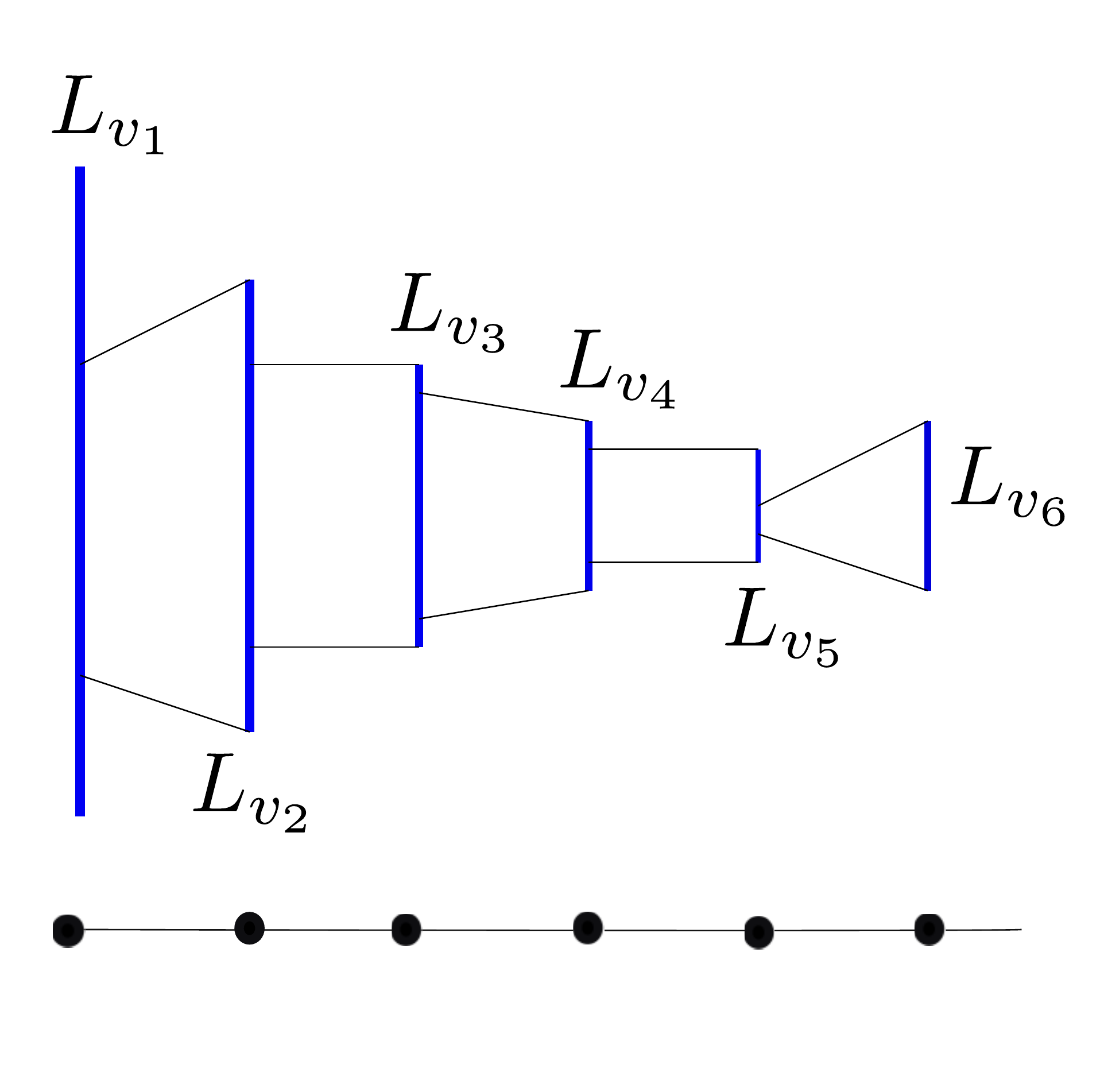} 
\caption{Ladder}
\label{ladder.fig}
\end{figure}

\medskip 
{\bf Ladder Axioms}: 

\begin{itemize}
\item{{\bf L1}} Each $x\in L$ belongs to some $\Sigma_x$.  

\item{{\bf L2}} Each map $\Pi_{w,v}$ is either constant or is an orientation-preserving  topological embedding. 

\item{{\bf L3}} $\L$ is  
a $(K,D,E,\delta_0)$-semicontinuous family of spaces. 
\end{itemize}

\begin{rem}
In general, $\pi(top(\L))$ and $\pi(bot(\L))$ are smaller than the base  $S$. If $$\pi(top(\L))=\pi(bot(\L))=S$$ then 
$\L$ will be a $K$-metric bundle. This will happen in the important case of {\em carpets}. \index{carpet}
\end{rem}

\begin{defn}
A $(K,D,E)$-{\em ladder} in $\X$ is a subtree of spaces $$\L=(\pi: L\to S)\subset {\mathfrak X}=(\pi: X\to T)$$ whose vertex and 
edge-spaces are oriented geodesic segments, equipped with a family of $K$-qi sections $\Sigma_\bullet$ and satisfying Axioms L0---L3.  
\end{defn}

\begin{example}
Let $x$ be a point in $X_u$ and let $\Sigma_x$ be a $K$-qi section of $\pi: X\to T$ defined over a subtree $S\subset T$ such that $x\in \Sigma_x$. 
Then $L=\Sigma_x$ is the total space of a $(K,0, \eta_0(2K))$-ladder. 
\end{example}

\begin{defn}\label{defn:subladder}
A {\em subladder} in $\L$ is a ladder $\L'=\L(\al')\subset \L=\L(\al)$ with the same center $u$ as $\L$, such that the sections $\Sigma'_\bullet$ of $\L'$ are restrictions of the sections $\Sigma_\bullet$ of $\L$. In particular, top and the bottom of $\L'$ are contained in sections of $\L$ through the end-points of $\al'$. 
\end{defn}

 In what follows, given a ladder $\L=\L_K(\al)$, 
 $\al\subset X_u$,  for each point $x\in \LL$ let $\ga_x\subset \Sigma_x\subset L$ denote the section over the interval $\llbracket u, \pi(x)\rrbracket$, connecting $x$ to a point in $\al$. Similarly, 
 given two points $x, y\in \LL$, if $\Sigma_x\cap \Sigma_y\ne \emptyset$ and the restriction of $\pi$ to $\Sigma_x\cup \Sigma_y$ is 1-1, then  
there exists a unique $K$-leaf $\ga_{x,y}$ in 
$\Sigma_x\cup \Sigma_y\subset \L$  connecting $x$ to $y$. 
 
 \medskip 
We omit a proof of the next lemma as it is straightforward: 

\begin{lemma}[Bisecting a ladder]\label{bisecting a ladder}
Suppose $u\in V(T)$, $\al= [xy]_{X_u}\subset X_u$ and we are given a ladder $\L=\L_{K,D,E}(\al)$. 
Then for every point $z\in [xy]_{X_u}$ the $K$-qi section $\Sigma_z\subset \LL$ decomposes $\L$ into two $(K,D,E)$-subladders $\L^+, \L^-$ such that 
\begin{enumerate}
\item $L^+_u=[zy]_{X_u}\subset \al$, $L^-_u=[xz]_{X_u}\subset \al$, 
\item 
$$
top(\L^-)= \Sigma_z= bot(\L^+). 
$$
\end{enumerate}
\end{lemma}

Applying this lemma twice, we obtain: 

\begin{cor}[Trisecting a ladder]\label{cor:trisecting a ladder}
Suppose $u\in V(T)$, $\al= [xy]_{X_u}\subset X_u$ and we are given a ladder $\L=\L_{K,D,E}(\al)$. 
Then for every subsegment $\al'=[x'y']_{X_u}\subset \alpha$ there exists a subladder $\L'=\L_{K,D,E}(\al')\subset \L$ bounded by the $K$-qi sections $\Sigma_{x'}, \Sigma_{y'}\subset \L$ (its bottom and top respectively).  
\end{cor}

Since a $(K,D,E)$-ladder $\L=(\pi: L\to S)$  
is a $(K,D,E,\delta_0)$-semicontinuous subtree of spaces in $\X$, 
as an application of the retraction Theorem \ref{thm:mitras-projection} we obtain:

\begin{cor}
[Retraction to ladders] \label{cor:ladder-retraction}
For every $(K,D,E)$-ladder $\L=(\pi: L\to S)$ there exists a 
$L_{\ref{thm:mitras-projection}}(K,D,E,\delta_0)$-coarse Lipschitz retraction $\rho_{\L}: X\to L$. 
\end{cor}

  \medskip 
 We next define {\em carpets} which are both ladders and metric bundles. While in Axiom L3 
 of a ladder we assume that fibers over {\em all} boundary vertices of $S$ have uniformly bounded diameter when projected to adjacent vertex-spaces $X_w, w\notin V(S)$,  in the definition of carpets  (where the base $S$ is an oriented interval $\llbracket u, w\rrbracket$)  
 we will allow one of the boundary vertices of $S$ (namely the vertex $u$) to violate this property (which is why $D=\infty$).  
 However, instead, we will add a stronger requirement on the other boundary vertex $w$ and a requirement on the top and the bottom.   
  
\begin{defn}
A $(K,\infty, \eta_0(2K))$-ladder $\A=(\pi: A\to S)\subset \X$ is called a $(K,C)$-{\em narrow carpet} or just a {\em $(K,C)$-carpet} if:

1.  $S$ is an interval $\llbracket u,w\rrbracket$ and, furthermore, the top and the bottom of $\A$ connect the two end-points of $A_u$ to that of  $A_w$, i.e.
$$
\pi(top(\A))=\pi(bot(\A))=S. 
$$
 In this case, we will say that $\A$ is bounded by the $K$-qi sections $\ga_1=bot(\L), \ga_2=top(\L)$  of the carpet. We will  refer to  $\beta=A_w$ 
 as the {\em (narrow) end} of the carpet and will say that $\A$ is {\em from $\al=A_u$ to $\beta=A_w$}. 

2. $A_w$ has length $\le C$. 

\noindent We will use the notation $\A= \A_K(\al)$ for such carpets to indicate the two key parameters. 
\end{defn}  
  
\begin{defn}\index{hallway} \label{defn:hallway} 
A $(K,\infty)$-carpet is called a $K$-hallway. 
\end{defn}

  \begin{rem}
 1.  Every $K$-hallway is a $K$-metric bundle. 
 
 2. The pair of sections $\ga_1, \ga_2$ determines a hallway $\A$ ``coarsely uniquely'': 
 The ambiguity in the definition comes from the choice of the vertical geodesics $A_t, t\in V(S)$,  
 and, hence, is uniformly controlled. Therefore, in what follows, we will ignore this ambiguity. 
\end{rem}

\medskip 
The next lemma establishes existence of ladder and hallway structures on subsets in $X$ 
which are unions of vertical geodesic segments.

\begin{lemma}\label{lem:E-ladder-structure} 
Suppose that ${\mathfrak X}$ is a tree of hyperbolic spaces. 
There exists a function $K'=K'_{\ref{lem:E-ladder-structure}}(K)$ such that the following holds:  

a0. Suppose that $\LL\subset \XX$ is a subset whose projection to $T$ is the vertex-set of a subtree $S\subset T$ containing  a vertex $u$ satisfying: 

a1. Every fiber $L_v=\LL\cap X_v, v\in V(S)$,  is an oriented geodesic segment $[x_v y_v]_{X_v}$ in $X_v$. 

a2. $\LL$ satisfies Property 4 of a semicontinuous family of spaces with the parameter $D$. 
Furthermore, in line with Property 3, for every oriented away from $u$ edge $e=[v,w]\in E(S)$, 
$\Hd_{X_{vw}}(L^v_w,L_w)\le E$, where, as before, 
$$
L^v_w= P_{X_{vw},X_w}(L_v)\subset X_w. 
$$ 

a3. Points $x_w, y_w$ are within distance $K$ (in $X_{vw}$) from points $x'_v, y'_v\in L_v$ respectively, so that 
$$
x_v\le x'_v\le y'_v\le y_v
$$
on the oriented segment $L_v$. 

Then $\LL\subset X$ is the union of vertex-spaces  of  a $(K',D,E)$-ladder $\L\subset \X$ centered at $u$. 

b0. Suppose that $\AA\subset \XX$ is a subset whose projection to $T$ is the vertex-set of an interval $S=\llbracket u,w\rrbracket\subset T$ such that:  

b1. Every fiber $A_v, v\in V(S),$ of $\AA$,  is an oriented geodesic segment $[x_v y_v]_{X_v}$ in $X_v$.

b2. For every edge $[v_1,v_2]$ in $S$, $d_{X_{v_1v_2}}(x_{v_1}, x_{v_2})\le K$,   $d_{X_{v_1v_2}}(y_{v_1}, y_{v_2})\le K$.  

Then $\AA$ is the union of vertex-spaces of a $K'$-hallway $\A\subset X$.  
\end{lemma}
\proof Our first goal is to define the function $K'$. 

We let $r':=D_{\ref{lem:sub-close}}(\delta_0', L'_0, K)$ be given by Lemma \ref{lem:sub-close}. For 
$k=K_{\ref{lem:close->qi}}(r', L'_0)$ given by Lemma \ref{lem:close->qi}, we let $\la'={k}_{\ref{lem:approximation}}(k)$ be given by  Lemma \ref{lem:approximation}. Lastly, set  
$$
K':= r'+ D_{\ref{lem:approximation}}(k).
$$

\noindent We now prove the lemma. 

a. We  define inductively 
the projections $\Pi_{v_1,v_2}$ (where $e=[v_1,v_2]$ is an edge in $S$ oriented away from $u$), as well as the edge-spaces $L_e$.

Suppose that for the subtree $B_n\subset S$, which is the closed $n$-ball centered at $u$, we 
defined (partial) $K$-qi sections $\Sigma$ and maps $\Pi$ satisfying all the requirements of a ladder 
with respect to the parameter $K'$. 

We extend the definitions of these sections and maps to the vertices in the ball $B_{n+1}\subset S$ 
as follows. Let $e=[v_1,v_2]$ be an edge of $S$ with $v_1\in  B_n, v_2\notin B_n$. 
Let $L'_{v_1}$ denote the oriented subsegment of $L_{v_1}$ bounded by $x'_{v_1}, y'_{v_1}$ respectively. 
Similarly, we define the edge-space $L_e$ as the oriented geodesic segment 
in $X_e$ spanned by the nearest-point projections of the end-points $x_{v_2}, y_{v_2}$ of $L_{v_2}$. 

According to Lemma \ref{lem:sub-close}, 
we have
$$
\Hd_{X_{v_1v_2}}(L'_{v_1}, L_{v_2})\le r'=D_{\ref{lem:sub-close}}(\delta_0', L'_0, K)\le K',
$$ 
$$
\Hd_{X_{v_1v_2}}(L_{e}, L_{v_2})\le r'. 
$$
These conditions ensure Property 3 of Definition \ref{defn:scfamily}, i.e. Axiom L3 of a ladder. 
   
 By  Lemma \ref{lem:close->qi}, we  extend the map $x_{v_2}\mapsto x'_{v_1}, y_{v_2}\mapsto y_{v_1}'$ to a 
 $k=K_{\ref{lem:close->qi}}(r',L'_0)$-quasiisometry of geodesic segments $q: L_{v_2}\to L'_{v_1}$,  
 which moves each point distance $\le r'$ (with respect to the metric of $X_{v_1v_2}$).  
 Applying Lemma \ref{lem:approximation}, we can replace the quasiisometry $q$ by an  
 increasing homeomorphism $\tilde{q}$ (or a constant function) within distance $D_{\ref{lem:approximation}}(k)$  from $q$, such that  
 $\tilde{q}$ is a ${k}_{\ref{lem:approximation}}(k)$-quasiisometry.

Since $q$ was moving every point of $L_{v_2}$ at most distance $r'$, it follows that $\tilde{q}$ moves every point within distance $K'=r'+D_{\ref{lem:approximation}}(k)$ in $X_{v_1v_2}$. 
We set 
$$
\Pi_{v_2,v_1}:= \tilde{q}. 
$$
Thus, we obtain maps $\Pi_{w,v}: L_{v_2}\to L_{v_1}$ 
for oriented edges $[v_1,v_2]$ of the tree $S=\pi(\LL)$, such that $d_{X_{v_1v_2}}(x, \Pi_{v_2,v_1}(x))\le K', x\in L_{v_2}$.  
For vertices $v_1, v_2$ of $S$ such that $v_1$ is between $u$ and $v_2$ we define 
the  map $\Pi_{v_2,v_1}: L_{v_2}\to L_{v_1}$ as the composition of  maps defined for the sequence of edges  connecting $v_2$ to $v_1$. 
If $\Pi_{v_2,v_1}$ is injective, then for $z\in L'_{v_1}$ we define the section $\Sigma_z\cap L_{v_2}$ as
$$
\Pi_{v_2,v_1}^{-1}(z). 
$$
If the map is not injective, i.e. constant, we choose an arbitrary point in  $L_{v_2}$ as  $\Sigma_z\cap L_{v_2}$. 

b. The proof of this part is exactly the same as of Part a, except that we use $x'_{v_1}=x_{v_1}, y'_{v_1}=y_{v_1}$. \qed

\section{Flow-spaces}\label{sec:flow-spaces}

\subsection{$K$-flow spaces and Mitra's retraction} \label{sec:def-of-flow-spaces}

Suppose that ${\mathfrak X}=(\pi: X\to T)$ is a tree of hyperbolic spaces. 
We fix a vertex $u\in T$, the {\em center of the flow} and orient all edges $e=[v, w]$ of $T$ so that $v$ is closer to $u$ than $w$. For each 
$4\delta_0$-quasiconvex subset  $Q_u\subset X_u$ we will define the $K$-{\em flow-space} 
$$
{\mathfrak Fl}_K(Q_u)= (\pi: {Fl}_K(Q_u) \to S) \subset {\mathfrak X},
$$
which, unlike ladders and carpets, depends only on $Q_u$ and on $K$, and 
which will be a $(K,D,E,4\delta_0)$-semicontinuous family of spaces (relative to the vertex $u$), with the parameter 
$E$ depending only on $K$ and $D=D_0$, where 
\begin{equation}\label{eq:D0}
D_0=D_{\ref{cobdd-cor}}(\delta'_0, \la'_0)
\end{equation}
is independent of $K$.

  However, for the construction to work, the parameter $K$ has to be sufficiently large, specifically, $K\ge K_0$, where $K_0$ (which depends only on the parameters of the tree of spaces $\X$) 
is given by the equation \eqref{K0}.    As before, we will use   ${\mathcal Fl}_K(Q_u)$ to denote the union of 
vertex-spaces of ${\mathfrak Fl}_K(Q_u)$. We first compute the auxiliary parameter $E$ and a certain parameter $R$ (depending on $K$) which will be used to define the $K$-flow. 

Suppose that $\la\ge \frac{3}{2}\delta_0$. 
Recall (Lemma \ref{lem:triple1}) that if the image of a  subset $Q$ of $X_v$ is 
$\la$-quasiconvex in $X_{vw}$ then $Q$ is $\hat{\la}$-quasiconvex in $X_v$ with 
\begin{equation}\label{hatla}
\hat{\la}= 1500 (L'_0\la)^3. 
\end{equation}

Take   
\begin{equation}\label{eq:R-ineq}
R\ge  R_0:= \max(2(\la'_0+ \delta'_0), R_{\ref{cobdd-cor1}}(\delta'_0, \la'_0))= 
2\la'_0+ 5\delta'_0. 
\end{equation}
Set (cf. \eqref{hatla}) 
\begin{equation}\label{eq:lambda'} 
\la':= 1500 (L'_0(R+2\delta'_0))^3, 
\end{equation}
\begin{equation}\label{eq:E-R-flow} 
E:=2(2\la'_0+3\delta'_0 +R)+ (\la'+\delta_0). 
\end{equation}
While our proofs will work whenever 
\begin{equation}\label{eq:K-R-ineq}
K\ge R+\la'+\delta_0, 
\end{equation} 
concretely, we will use 
\begin{equation}\label{eq:vee} 
K = R^\wedge:=(15 L'_0 R)^3, \hbox{~~~i.e.~~~} R= K^\vee:= \frac{1}{15L'_0} K^{1/3}  
\end{equation}
(The reader will verify that this $K$ satisfies the inequality \eqref{eq:K-R-ineq}.) 
Thus, the inequality \eqref{eq:R-ineq} translates into the inequality 
\begin{equation}\label{eq:K-ineq}
K\ge K_0= 15^3 (2\la'_0+ 5\delta'_0)^3 (L'_0)^3. 
\end{equation}
Note also that \eqref{eq:E-R-flow} makes $E$ a function of $K$ 
\begin{equation}\label{eq:E-K-flow} 
E=E_{\ref{eq:E-K-flow}}(K),
\end{equation}
while and $R$ also becomes a function of $K$.

We inductively define $4\delta_0$-quasiconvex subsets $Q_v\subset X_v, Q_e\subset X_e$, $v\in V(T), e\in E(T)$,  and, at the same time, verify conditions of Definition \ref{defn:scfamily} for the collection of subsets $Q_v, Q_e$, aiming eventually to Prove Theorem \ref{thm:semicontinuity-of-flows}. 
Assuming that for $v\in V(T)$ a $4\delta_0$-subset  
$Q_v\subset X_v$ is defined, for the oriented edge $e=[v,w]$ of $T$ 
(oriented away from $u$)  
we set 
$$
Q^v_w:= P_{X_{vw},X_w}(Q_v), \quad Q'_w:= N^e_R(Q_v)\cap X_w.
$$
According to Corollary \ref{cor:0-flow-space}, 
$$
\Hd_{X_{vw}}(Q^v_w, Q'_w)\le 2(2\la'_0+3\delta'_0 +R). 
$$
Note that both $X_w, Q_v$ are $\la'_0$-quasiconvex in $X_{vw}$.

Furthermore, by Lemma \ref{lem:coarse-intersections-are-qc-2}, 
since  $R\ge R_0\ge 2\la'_0+ 2\delta'_0$, 
the intersection $Q'_w:= N^e_R(Q_v)\cap X_w$ is 
$\la_{\ref{lem:coarse-intersections-are-qc-2}}= R+2\delta'_0$-quasiconvex in $X_{vw}$.  
Hence, $Q'_w$ is $\la'=\widehat{R+2\delta'_0}$-quasiconvex in $X_w$, where 
$$
\la'= 1500 (L_0' (R+ 2\delta'_0))^3, 
$$
see Lemma \ref{lem:triple1}.

Therefore, by \eqref{eq:qc-nbd}, the $\delta_0$-hull, taken in $X_w$,  
$$
Q_w:=\hull_{\delta_0}(Q'_w)
$$
is $(\la'+\delta_0)$-Hausdorff close to $Q'_w$, thus,
$$
\Hd_{X_{vw}}(Q^v_w, Q_w)\le E=2(2\la'_0+3\delta'_0 +R)+ (\la'+\delta_0), 
$$
verifying the condition \eqref{eq:E-in} in Part 3 of a semicontinuous family of spaces (in the case when $Q_w'\ne \emptyset$, equivalently, $Q_w\ne \emptyset$). 

We define the edge-space $Q_e$ as the $\delta_0$-null (in $X_e$) of the projection
$$
P_{X_{vw},X_e}(Q_w). 
$$
Thus,
$$
\Hd_{X_{vw}}(Q^v_w, Q_e)\le \delta_0+1. 
$$
At the same time, since each point of $Q'_w$ is within distance $R$ from $Q_v$, each point of 
$Q_w$ is within distance  
$$
R+\la'+\delta_0$$ 
 from $Q_v$, where both distances are computed in $X_{vw}$. Since
 $$
 K=  (15 L'_0 R)^3\ge  R+\la'+\delta_0,  
 $$
we conclude that each point of $Q_w$ is within distance $K$ from $Q_v$. From this, since $Q_e$ was defined as the projection of $Q_w$ to $X_e$, it also follows that $\Hd_{X_{vw}}(Y_w, Y_e)\le K$. Thus, we verified Part 3 of Definition 
\ref{defn:scfamily} (for the edge $e$). Since the subsets $Q_w, Q_e$ were defined as $\delta_0$-hulls in $\delta_0$-hyperbolic spaces, we conclude that $Q_e\subset X_e, Q_w\subset X_w$ are $4\delta_0$-quasiconvex, verifying Part 1 of Definition \ref{defn:scfamily}. 

Lastly, we turn to Part 4 of Definition \ref{defn:scfamily}. As we noted earlier, $Q_w=\emptyset$ if and only if 
$Q'_w=N^e_R(Q_v)\cap X_w=\emptyset$. In other words, the $\la'_0$-quasiconvex subsets $Q_v, X_w\subset  X_{vw}$ are $R$-separated. Since $R$ was chosen to satisfy 
$$
R\ge R_0=R_{\ref{cobdd-cor}}(\delta'_0, \la'_0)=2\la'_0+5\delta_0,$$
Corollary \ref{cobdd-cor} implies that subsets $Q_v, X_w\subset X_{vw}$ are 
$D=D_{\ref{cobdd-cor}}(\delta'_0, \la'_0)$-coboun\-ded.  This verifies Part 4 of Definition \ref{defn:scfamily}.

\begin{defn}\label{defn:flow-space} \index{flow-space $Fl_K(Q)$}
We define the {\em $K$-flow space ${\mathfrak Fl}_K(Q_u)$} of $Q_u$ as the subtree of spaces in 
${\mathfrak X}$ as follows. The nonempty subsets $Q_v, Q_e$ defined by the inductive procedure above will be the vertex/edge spaces of ${\mathfrak Fl}_K(Q_u)$. 
The incidence maps $g_{ev}$ of  ${\mathfrak Fl}_K(Q_u)$ are the compositions of the incidence maps $f_{ev}$ with fiberwise nearest-point 
projections in $X_v$ to $Q_{v}$. The vertex and edge-spaces of  ${\mathfrak Fl}_K(Q_u)$ are equipped with path-metrics induced from the ambient 
path-metrics on vertex and edge-spaces of 
${\mathfrak X}$.  We let $Fl_K(Q_u)\subset X$ denote the total space of  ${\mathfrak Fl}_K(Q_u)$, set $S:= \pi(Fl_K(Q_u))$; 
we will use the notation ${\mathcal Fl}_K(Q_u)$ for the disjoint union 
$$
\coprod_{v\in V(S)} Q_v,
$$
which is the union of vertex-spaces of ${\mathfrak Fl}_K(Q_u)$. We will equip $Fl_K(Q_u)$ with the standard 
path-metric of a tree of spaces.  

Sometimes it will be convenient to restrict flow-spaces to subtrees $T'\subset T$. We will denote such ``subflows'' by
$$
{\mathfrak Fl}^{T'}_K(Q_u). 
$$
\end{defn}

\begin{rem}\label{rem:flows}
1. The $\delta_0$-neighborhoods in the definition of  flow-spaces  are taken in order to ensure that the each inclusion map $Q_w\to X_w$ is a $(1,C_{\ref{lem:delta-hull}}(\delta_0))$-quasiisometric 
embedding, where $Q_w$ is equipped with the path-metric induced from $X_w$, see Lemma \ref{lem:delta-hull}. 

2.  In general, it is not true that for Hausdorff-close subsets $A, B\subset X_u$, the $K$-flow spaces are Hausdorff-close to each other. However, if 
$$
B\subset N^{fib}_{r}(Q_u)\subset X_{u}$$ 
then (by the very definition of a flow-space) 
$$
Fl_K(B)\subset Fl_{K+ r}(Q_u). 
$$
Similarly,
$$
N^{fib}_r(Fl_K(Q_u)) \setminus N^{fib}_r(Q_u) \subset Fl_{K+ r}(Q_u). 
$$
\end{rem}

The discussion preceding the definition of flow-spaces proves: 

\begin{theorem}\label{thm:semicontinuity-of-flows}
For every $K\ge K_0$, the flow-space ${\mathfrak Fl}_K(Q_u)$ is a $(K,D,E)$-semicontinuous family of spaces 
in $\X$, where $D=D_0=D_{\ref{cobdd-cor}}(\delta'_0, \la'_0)$ and $E=E_{\ref{thm:semicontinuity-of-flows}}(K)$ is given by the equation \eqref{eq:E-R-flow}.  
In particular, every $x\in Fl_K(Q_u)$ belongs to a $K$-leaf $\ga_x$ in $Fl_K(Q_u)$ connecting $x$ to $Q_u$. 
\end{theorem}

Combining with with the existence of uniform coarse Lipschitz retractions to semicontinuous subtrees of spaces 
(Theorem \ref{thm:mitras-projection}), we conclude:

\begin{theorem}[Mitra's Retraction] \label{mjproj}\index{Mitra's retraction $\rho$} 
Suppose that ${\mathfrak X}$ is a tree of hyperbolic spaces.
Then for each $K\ge K_0$,  there exists an $L_{\ref{mjproj}}(K)$-coarse Lipschitz retraction,  
called {\em Mitra's retraction}, $\rho=\rho_{Fl_K(Q_u)}: X\to Fl_K(Q_u)$,  where  
$$
L_{\ref{mjproj}}(K)= L_{\ref{thm:mitras-projection}}(K,D_0,E_{\ref{thm:semicontinuity-of-flows}}(K),4\delta'_0). 
$$
\end{theorem}

Below we collect several consequences of  Theorem \ref{mjproj}. 

\begin{cor}\label{cor:mjproj}
The inclusion map $Fl_K(Q_u)\to X$ is an $L_{\ref{mjproj}}(K)$-qi embedding.
\end{cor}

\begin{cor}[M.~Mitra, \cite{mitra-trees}] 
 If $X$ were a hyperbolic metric space then for all $u\in V(T)$ and $4\delta_0$-quasiconvex subsets $Q_u\subset X_u$,
the flow-spaces $Fl_{K}(Q_u)$ would be uniformly quasiconvex subsets in $X$. 
\end{cor}

\begin{cor}\label{cor:disjoint->cobounded}
If $\pi(Fl_K(X_{u_1}))\cap \pi(Fl_K(X_{u_2}))=\emptyset$, then the flow-spaces $Fl_K(X_{u_1})$, $Fl_K(X_{u_2})$ are 
$L_{\ref{mjproj}}(K)$-Lipschitz-cobounded in $X$ (cf. Definition \ref{def:cob}). 
\end{cor}
\proof We will be using Mitra's projections
$$
\rho_i= \rho_{Fl_K(X_{u_i})}, i=1,2. 
$$
Since these projections are $L_{\ref{mjproj}}(K)$-coarsely Lipschitz, it suffices to show that diameters of $\rho_{i}(Fl_K(X_{u_{3-i}}))$, $i=1, 2$, are uniformly bounded. 
Let $v_i\in \pi(Fl_K(X_{u_i}))$, $i=1, 2$, denote the vertices realizing the minimal distance between these subtrees of $T$. 
Let $e_1=[v_1,w_1]$, $e_2=[v_2,w_2]$ be the edges incident to $v_1, v_2$ and contained in the interval 
$\llbracket v_1, v_2\rrbracket $ (it is possible that $w_1=v_2$, $w_2=v_1$). Then, by  the definitions of Mitra's 
projection (see the proof of Theorem \ref{thm:left-inverse}) and the $K$-flow, for $i=1, 2$, 
$$
\rho_{3-i}(Fl_K(X_{u_i}))\subset \rho_{3-i} (X_{v_i})= \{x_{3-i}\}\subset Fl_K(X_{u_{3-i}})\cap X_{v_{3-i}},
$$
i.e., $\rho_{3-i}(Fl_K(X_{u_i}))$ is the singleton $\{x_{3-i}\}$. Thus, 
the flow-spaces $Fl_K(X_{u_1})$, $Fl_K(X_{u_2})$ are $L_{\ref{mjproj}}(K)$-Lipschitz-cobounded in $X$, see Definition \ref{def:cob}. \qed 

\begin{example}
One can realize the hyperbolic plane as the total space (up to a quasiisometry) of a metric line bundle over a line. Namely, 
let $T=\RR$, where the vertices are the integer points. We will identify  $T$ with the $y$-axis in the upper half-plane model of the hyperbolic plane (of course, we parameterize the $y$-axis in $\H^2$ by the hyperbolic arc-length). The projection $\pi: \H^2\to T$ is given by the $y$-coordinate of the points in $\H^2$. Examples of 1-qi sections of the bundle $\pi: \H^2\to T$ are given by  hyperbolic geodesics in $\H^2$ which are vertical rays in the half-plane model. For $K=1$, the $K$-flow space of a singleton $Q_u=\{q\}$ is then 
such a hyperbolic geodesic through this point.  Mitra's retraction to $L=Fl_K(Q)$ in this example is within finite distance from the horocyclic projection $\H^2\to L$. While it is $1$-Lipschitz, it is not close to the nearest-point projection $\H^2\to L$. For a general $K\ge 1$, the $K$-flow of a singleton $Q_u=\{q\}$ is Hausdorff-close to a geodesic in $\H^2$, equivalently, to a $K$-qi section through $q$. 
\end{example}

\subsection{Basic properties of flow-spaces} \label{sec:basics-of-flows} 

Most of the time, besides Mitra's retractions and the fact that each flow-space forms a semicontinuous subtree of spaces in $\X$, instead of the definition of flow-spaces we will use their properties summarized in the next proposition (we recall the equations \eqref{eq:R-ineq} and \eqref{eq:vee} defining the constants $K_0, R_0$ and the function $R\mapsto R^\wedge=K$):

\begin{prop}\label{prop:flow-prop} 
Suppose that $Q_u\subset X_u$ is a $4\delta_0$-quasiconvex subset. 
\begin{enumerate}
\item 
Suppose that a vertex $w$ lies between vertices $u$ and $v$. Then for every $r\ge 0$, and all $K\ge K_0$, 
\begin{equation}\label{i}
N^{fib}_r(Fl_K(X_u))\cap N^{fib}_r(Fl_K(X_v)) \subset N^{fib}_r(Fl_K(X_w)), 
\end{equation}
and 
\begin{equation}\label{ii}
N_r(Fl_K(X_u))\cap N_r(Fl_K(X_v)) \subset N_r(Fl_K(X_w)). 
\end{equation}

\item 
Suppose that $R\ge R_0$ and let $\ga$ is a $R$-qi leaf in $X$ emanating from some  
$\ga(u)\in Q_u$. Then  $\ga\cap \XX$   
is contained in ${\mathcal Fl}_{K}(Q_u)$, where $K=R^{\wedge}$.

\item 
For all $K\ge R_0$ and $Q_v:= Fl_{K}(Q_u)\cap X_v$, we have
$$
Q_u\subset Fl_{K^{\wedge}}(Q_v). 
$$

\item 
For every boundary edge $e=[v,w]$ of $S= \pi(Fl_{K}(Q_u))$, the subsets $Q_v, X_w$ are $D_0$-cobounded in $X_{vw}$, where $D_0$ is given by \eqref{eq:D0}. 

\end{enumerate}
\end{prop}
\proof (1)  The first containment follows from the inductive nature of the definition of flows. 
 We now prove the second inclusion. Let 
 $x\in Fl_K(X_u), y\in Fl_K(X_v)$ and $z\in \bar{B}(x,r)\cap \bar{B}(y,r) \subset X$. 
Up to relabeling, there are two cases: 

a. $w$ lies in the geodesic segment $\pi(y) v\subset T$. Then $y\in Fl_K(X_w)$ (according to \eqref{i} with $r=0$), which implies that $z\in N_r(Fl_K(X_w))$. 

b. The vertex $w$ is  not in  $u\pi(x) \cup v\pi(y)$, hence, $w$ separates $\pi(x), \pi(y)$ in $T$. Therefore, after relabeling, $w$ separates $\pi(x)$ from $\pi(z)$. 
In particular the geodesic $xz\subset X$ crosses $X_w$ and, hence, $z\in N_r(X_w)\subset N_r(Fl_K(X_w))$.

(2) The proof is by induction on the distance from vertices of $\pi(\ga)$ to $u$. We will use the notation 
Let  $Q_v, Q'_w$, etc.   as in the definition of the flow-space  
$Fl_{K}(Q_u)$.

The base of induction (for the vertex $u$) is clear. Suppose that $e=[v,w]$ is an edge in $\pi(\ga)$ oriented away from $u$ such that   $x= \ga(v)\in Q_v, y=\ga(w)\in X_w$. 
  Arguing inductively, we 
assume that the claim holds for the point $x$, 
i.e. 
$$
x\in  Q_v={\mathcal Fl}_{K}(Q_u) \cap X_v.$$ 
 Since $\ga$ is an $R$-qi leaf in $X$, 
 $$
 d_{X_{vw}}(x, y)\le R,  
 $$
i.e. by the definition of the function $R\mapsto R^{\wedge}=K$, 
$$
y\in N^e_R(Q_v)\cap X_w= Q'_w\subset Q_w.
$$
  
Part (3) is an immediate consequence of (2).   \qed

\begin{notation} In what follows, we will refer to $K$-qi leaves in $Fl_K(Q_u)$ as {\em leaves} in  $Fl_K(Q_u)$ and denote them $\ga$ with some subscripts. We will use the notation $\ga_{x}$, 
$\ga_{x,y}$ for such leaves provided that $\ga$ has $x$ as one of its end-points (in the former case) or $\ga$ connects $x$ and $y$. 
\end{notation}

 \begin{lemma}\label{max-qi-section}
Given any point $x\in Fl_K(X_u)$, there is a maximal\footnote{Here maximal means we can not find a $K$-qi 
section $\Sigma'_x$ over a subtree $T'_x\subset T$ containing $T_x$ such that 
$\Sigma_x\subsetneq \Sigma'_x$.  } 
$K$-qi section $\Sigma_x$ in $Fl_K(X_u)$ over a subtree $T_x$ of $T$ such that $x\in \Sigma_x$ and $u\in T_x$. 
\end{lemma}
\proof Let $w=\pi(x)$. By Theorem \ref{thm:semicontinuity-of-flows}, there exists a $K$-qi section $\Sigma_{\llbracket u, w\rrbracket,x}$ 
over the geodesic 
$\llbracket u, w\rrbracket \subset T$ connecting $u$ and $w$ such that $x\in \Sigma_{\llbracket u, w\rrbracket,x}$. 
We define a poset consisting of  
$K$-qi sections $\Sigma_{S,x}$ over subtrees $S\subset T$ and containing $\Sigma_{\llbracket u, w\rrbracket,x}$. Define the partial order 
$\Sigma_{S,x}\le \Sigma_{S',x}$ if the qi section $\Sigma_{S',x}$ extends $\Sigma_{S,x}$. This poset is clearly nonempty. 
The existence of  a maximal element in this poset follows from Zorn's lemma. \qed

\medskip
In chapter \ref{ch:CT} we will also use the following property of flow-spaces: 

\begin{lemma}\label{lem:intersection-with-flow-nbd}
Suppose that $Q_w=Fl_K(X_u)\cap X_w\ne \emptyset$. Then $N_R(Fl_K(X_u))\cap X_w$ is contained in 
$$
N^{fib}_{D}(X_u),
$$
$D=D_{\ref{lem:intersection-with-flow-nbd}}(R,K)$. 
\end{lemma}
\proof Let $S$ denote the subtree $S=\pi(Fl_K(X_u))\subset T$. Since $Y= X_S$ is $\eta_0$-properly embedded in $X$, it suffices to measure distances and consider geodesics in $Y$. 
Let $v$ be a vertex in $S$ and consider a point $x\in X_v$ such that
$$
d_Y(x, Q_v)\le r
$$
and suppose that $y\in X_w$ is within distance $R$ from $x$. 
Our goal is to bound $d_{X_w}(y, Q_w)$ in terms of $R$ and $r$. Lemma will then follow by considering the case $r=0$, as
$$
d_{X_w}(y, Q_w)\le \eta_0(  d_(y, Q_w) ) 
$$
Since the length of the projection of $[xy]_Y$ to $S$ is bounded from above by $R$, the proof can be done by induction 
on $d_S(v,w)$. The base-case $n=0$ is clear as 
$$
d_Y(y, Q_w)\le R+r. 
$$
Suppose that there is a function $\theta(n,r,R)$ such that 
if $d_S(v,w)= n$, then 
$$
d_Y(y, Q_w)\le \theta(n,r,R). 
$$
Consider now points $x, y$ such that $d_S(v,w)= n+1$.

Let $[zy]_Y$ be a subsegment in $[xy]_Y$ connecting $z\in X_t$ to $y$, such that 
$d_S(t,w)=1$.   We will consider the more difficult case when 
$t$ lies between $u$ and $w$ and leave the case when $w$ is between $u$ and $t$ to the reader.  
By the induction hypothesis (applied to the points $x, z$), there exists 
$p\in Q_t$ such that $d(p,z)\le  \theta(n,r,R)$.  
Let $\bar{p}$ denote the projection (taken in $X_{tw}$) of $p$ to 
$X_w$. Then 
$$
d_Y(\bar{p}, Q_w)\le d_{X_{tw}}(\bar{p}, Q_w)\le E= E_{\ref{thm:semicontinuity-of-flows}}(K). 
$$
Since $y\in X_w$ and 
$$
d_{X_{tw}}(p, y)\le L'_0(d(p, y)+1) \le L'_0(R+\theta(n,r,R)+1), 
$$
it follows that $d(p, \bar{p})\le L'_0(R+\theta(n,r,R)+1)$ as well. By combining the inequalities, we obtain:
$$
d_Y(y, Q_w)\le \theta(n+1,r,R):= E +  2L'_0(R+\theta(n,r,R)+1). \qed 
$$


\subsection{Generalized flow-spaces}\label{sec:gfs}

In this section we discuss several generalizations of flow-spaces and retractions to these, generalizing Mitra's projection.  We recall (see Definition \ref{defn:flow-space}) that superscript in the form of a subtree $S\subset T$ in the notation of a flow-space,  
means that the flow-space is taken in the subtree of spaces ${\mathfrak X}_S\subset {\mathfrak X}$ 
over  the subtree $S\subset T$.

\begin{defn}\label{def:subtree-flows} 
Assume that $K\ge K_0$. 

1. 
We define flow-spaces of 
 quasiconvex subsets $Q_e\subset X_e$ of edge-spaces, $e=[v,w]$. Define subtrees $T_v$ (resp, $T_w$) in $T$ as the maximal subtree in $T$ containing $v$ (resp. $w$) and disjoint from $w$ (resp. $v$).   We define the flow-space ${\mathfrak Fl}_K(Q_e)$ of such $Q_e$ 
 so that its union of vertex-spaces is 
 $$
{\mathcal Fl}_K(Q_e):= {\mathcal Fl}^{T_v}_K(\hull_{\delta_0}f_{ev}(Q_e)) \cup {\mathcal Fl}^{T_w}_K(\hull_{\delta_0} f_{ew}(Q_e)). 
$$

2. Similarly, we define flow-spaces in $X$ of $X_S\subset X$, where  $S\subset T$ is a subtree. For each $w\in V(S)$ we 
define the maximal subtree $T_w\subset T$ whose intersection with $S$ equals $\{w\}$. 
We then define the $K$-flow-space  ${\mathfrak Fl}_K(X_S)$ in $X$ so that $X_S\subset {\mathfrak Fl}_K(X_S)$ and 
the union of vertex-spaces of   ${\mathfrak Fl}_K(X_S)$ equals   
\begin{equation} \label{eq:flow-of-subtree}
{\mathcal Fl}_K(X_{S})= \bigcup_{v\in V(S)} {\mathcal Fl}^{T_v}_K(X_v).  
\end{equation}
\end{defn}

\begin{cor}\label{cor:qi-embedding-of-flow-of-edge-space}
For  oriented edges $e=[u,v]$ of $T$ the inclusion maps  
$$
Fl_K(X_e)\to Fl_K(X_{u}) 
$$
are uniformly quasiisometric embeddings. \end{cor}
\proof Apply Theorem \ref{mjproj} to the $\delta_0$-quasiconvex hull of $Q_u:=f_{eu}(X_{e})$ and observe that $Fl_K(X_e)\subset Fl_K(Q_u)\subset Fl_K(X_{u})$. \qed 

\begin{prop}\label{prop:uniform embeddings for subtree flows} 
There exists $L=L_{\ref{prop:uniform embeddings for subtree flows}}(K)$ such that 
for every subtree $S\subset T$ and $K\ge K_0$, the flow-space $Fl_K(X_{S})$ is $L$-qi embedded in $X$. 
\end{prop}
\proof We continue with the notation introduced in Definition \ref{def:subtree-flows}.   For each $w\in V(S)$ we define Mitra's retraction $\rho_w= \rho_{X_w}:  X_{T_w}\to Fl^{T_w}_K(X_w)$.   
Hence, the collection of maps $\rho_w$ is uniformly coarsely Lipschitz. We then obtain a (uniformly) coarse  Lipschitz retraction 
$$
X\to Fl_K(X_{S})
$$
whose restriction to $\pi^{-1}(S)$ is the identity and whose restriction to each $X_{T_w}$ equals 
$\rho_w$. \qed

\medskip 
We now define flow spaces of metric bundles. Let $S\subset T$ be a subtree and ${\mathfrak Q}=(\pi: Q\to S)$ be a metric $K_1$-bundle in $X$ whose vertex spaces $Q_v\subset X_v$ are $4\delta_0$-quasiconvex subsets. Let  ${\mathcal Q}=\bigcup_v Q_v$ denote the union of vertex-spaces of ${\mathfrak Q}$. 
For each vertex $v\in V(S)$, as above, we have the associated subtree $T_v\subset T$, which is the maximal subtree in $T$ such that 
$T_v\cap S=\{v\}$. Accordingly, we have subtrees of spaces $X_{T_{v}}\subset X$. 

For each vertex $v\in S$ we take the $K_2$-flow-space  
$$
{\mathcal F}_v:={\mathcal Fl}_{K_2}(Q_v)\cap X_{T_v}  
$$
inside $X_{T_v}$. Lastly, set 
$$
{\mathcal Fl}_{K_2}(Q_u):= \bigcup_{v\in V(S)} {\mathcal F}_v. 
$$
Then, as in the case when ${\mathcal Q}$ was the union of vertex spaces $X_v, v\in V(S)$ (see   
in Definition \ref{def:subtree-flows}(2)), 
${\mathcal Fl}_{K_2}(Q_u)$ is the union of vertex-spaces of a tree of spaces 
${\mathfrak Fl}_{K_2}({\mathfrak Q})$. 

\begin{defn}\label{defn:generalized flow-space}
We will refer to ${\mathfrak Fl}_{k}({\mathfrak Q})$ (and its total space $Fl_k(\QQQ)$) as a {\em generalized $k$-flow-space}, or the {\em $k$-flow-space of a metric bundle}. \index{generalized flow-space}  
\end{defn}

\medskip 

The following is an extension of Mitra's theorem to such generalized flow-spaces: 

\begin{theorem}\label{thm:generalized retraction} 
For every $K_2\ge K_0$, there exists an $L=L_{\ref{thm:generalized retraction}}(K_1,K_2)$-coarse Lipschitz retraction $X\to  Fl_{K_2}({\mathfrak Q})$. 
\end{theorem}
\proof The proof is similar to that of Proposition \ref{prop:uniform embeddings for subtree flows}. Over subtrees $T_w$, $w\in V(S)$, we use Mitra's retractions 
$X_{T_w}\to Fl_{K_2}({\QQQ}) \cap X_{T_w}$. The fact that these maps define a uniformly coarse  Lipschitz retraction $X_S\to Q$ (over the three $S$) follows from 
the  assumption that  ${\mathfrak Q}= (\pi: Q\to S)$ is a $K_1$-bundle, cf. Corollary \ref{cor:projection}.  
 \qed 

This theorem will be used in  Section \ref{sec:generalized hallways}.

\subsection{Boundary flows}\label{sec:Boundary flows}

In this section we will define {\em ideal boundary flows}, which are ideal counterparts of flow-spaces discussed above. Ideal boundary flows will be used only in chapters \ref{ch:description-of-geodesics} (specifically, Section \ref{sec:brdy-acyl-case}) and \ref{ch:CT}. The definition of  boundary flows given below   
 was latent in \cite[Definition 4.3]{ps-limset}.

\medskip \index{boundary flow} 
Fix a vertex $u\in V(T)$ and consider a subset  $Z_u\subset \geo X_u$. We will define the {\em flow-space} ${\mathfrak Fl}(Z_u)$ of $Z_u$ in $\geo X$, analogously to the definition of flow-spaces ${\mathfrak Fl}_K(Q_u)$ of quasiconvex subsets $Q_u\subset X_u$ given in Section \ref{sec:flow-spaces}. Each  ${\mathfrak Fl}(Z_u)=\ZZZ$ will be a tree of topological spaces over a subtree $S\subset T$, with vertex-spaces  $Z_v\subset \geo X_v$ and edge-spaces $Z_e\subset \geo X_e$. 

Our definition of subsets  $Z_v\subset \geo X_v$, $Z_e\subset \geo X_e$, is inductive on $d_T(u,v)$, analogously to the definition of ${\mathfrak Fl}_K(Q_u)\subset X$.  
We assume that subsets $Z_v$,  $Z_e$ are defined for all vertices and edges of $T$ contained in the ball (in $T$) of radius $n$ centered at $u$.  Consider an edge $e=[v,w]$ in $T$ oriented away from $u$, such that $d_T(u,v)=n$. We have qi embeddings $f_{ev}: X_e\map X_v$ and $f_{ew}: X_e\map X_w$. They induce topological embeddings ({\em boundary maps}) 
$$
\geo f_{ev}: \geo X_e\map \geo X_v, \quad \geo f_{ew}: \geo X_e\map \geo X_w, 
$$
see Section \ref{sec:ideal boundaries}. Define 
$$
Z_e:=(\geo f_{ev})^{-1}(Z_v), \quad \hbox{and}\quad Z_w:= \geo f_{ew}(Z_e).
$$ 
The boundary maps $\geo f_{ev}, \geo f_{ew}$ provide incidence maps $Z_e\to Z_v, Z_e\to Z_w$. We continue inductively. Define $S$ as the subtree in $T$ spanned by the vertices $v\in T$ such that $Z_v\ne \emptyset$. Thus, we obtain a tree $\ZZZ=(Z\to S), Z=Fl(Z_u)$, of topological spaces with the vertex-spaces $Z_v, v\in V(S)$, the edge-spaces 
$Z_e, e\in E(S)$ and incidence maps  $\geo f_{ev}$ as above, see Definition \ref{defn:top-tree}. 
We will use the notation 
${\mathcal Fl}(Z_u)$ for $\ZZ$, the union of vertex-spaces of $\ZZZ$. Note that since the maps 
$\geo f_{ev}$ are 1-1, the set ${\mathcal Fl}(Z_u)$ breaks as a disjoint union of flow-spaces 
${\mathcal Fl}(\{z\})$ of singletons $\{z\}\subset Z_u$. For each $v\in V(T)$, the intersection  
$$
Fl_v(\{z\}):={\mathcal Fl}(\{z\})\cap \geo X_v$$
is either empty or is a singleton. Similarly, for each edge $e=[v,w]\in E(T)$, the intersection $Fl_e(\{z\}):= Fl(\{z\})\cap Z_e$ is either empty or is a singleton.  In view of the inductive definition of 
$Fl(Z_u)$, if $Fl_w(\{z\})\ne \emptyset$ for a vertex $w$, then for every vertex $v\in uw$,  
$Fl_v(\{z\})=\{z'\}$ is also nonempty and we have
$$
Fl_v(\{z'\})= Fl_w(\{z\}). 
$$
Furthermore, for each edge $e$ in the interval $uv$, $Fl_e(\{z\})$ is nonempty as well. For 
 $z\in \geo X_u$ we define a tree $T_z\subset S$  such that $V(T_z)$ consists of vertices $v$ for which $Fl_v(\{z\})$ is nonempty.

\medskip 
The next lemma is partially proven in \cite[Lemma 4.4]{ps-limset}; we include a proof for the sake of completeness: 

\begin{lemma}\label{lem: flow condition}
Suppose that $\X=(\pi: X\to T)$ is a tree of hyperbolic spaces. Consider a point 
 $z\in \geo X_u$ and a geodesic ray $\alpha_u\subset X_u$ asymptotic to $z$. 
Then for a vertex $v\in V(T)$ the following are equivalent:

1. $Fl_v(\{z\})$ is nonempty.

2.  There is a geodesic ray $\al_v\subset X_v$ such that $\Hd(\alpha_u, \al_v)<\infty$. In this case $Fl_v(\{z\})=\{\al_v(\infty)\}$.

3. $\al_u$ is Hausdorff-close to a subset of $X_v$. 
\end{lemma}
\proof It suffices to prove the equivalence in the case when $u, v$ span an edge $e=[v,w]$ in $T$ (the general case is proven  by induction on $d_T(u,v)$). As we noted earlier, $Fl_v(\{z\})=\{z_v\}\ne \emptyset$ if and only if $Fl_e(\{z\})=\{z_e\}\ne \emptyset$. The latter holds if and only if  there is a geodesic ray $\alpha_e$ in $X_e$ such that $\alpha_e(\infty)=z_e$ and
$$
\geo f_{eu}(z_e)=z, \quad \geo f_{ev}(z_e)=z_v. 
$$
The paths  $f_{eu}(\alpha_e)$, $f_{ev}(\alpha_e)$  are quasigeodesic rays in $X_u$ and $X_v$ respectively, which are asymptotic, respectively, to the points $\{z\}, \{z_v\}$.  By the stability of quasigeodesics we have $\Hd(\alpha_u,f_{eu}(\alpha_e))<\infty$ and there is a geodesic ray $\al_v$ in $X_v$ Hausdorff-close to $f_{ev}(\alpha_e)$. This proves the equivalence of (1) and (2). 

The implication (2) $\Rightarrow$ (3) is clear. If (3) holds, then the image of $\al_u$ under 
the projection $P=P_{X_{uv},X_v}$ is again a quasigeodesic in $X_v$ (which is Hausdorff-close to $\al_u$). By the stability of quasigeodesics, $P(\al_u)$ is Hausdorff-close to a geodesic $\al_v$ in $X_v$. \qed

\begin{cor}
The following is an equivalence relation on $\coprod_{v\in V(T)} \geo X_v$: 

Two points $z_i\in \geo X_{u_i}$, $i=1,2$ 
are related iff $Fl_{u_2}(\{z_1\})= \{z_2\}$. 
\end{cor}

We next relate boundary flow-spaces to the flow-spaces in $\X$.  

\begin{lemma}\label{lem:flow-flow}
There is a constant $K$ depending on the parameters of $\X$ such that 
the following holds:

(1)  If $\alpha_u\subset X_u$ is a geodesic ray asymptotic to $z\in \geo X_u$, then 
$Fl_v(\{z\})= \geo(Fl_K(\alpha_u)\cap X_v, X_v)$ for all $v\in V(T_z)$. 

(2) $Fl_v(\geo X_u)=\geo (Fl_K(X_u)\cap X_v, X_v)$ for all $v\in V(T)$.
\end{lemma}
\proof (1) The inclusion $\geo(Fl_K(\alpha_u)\cap X_v, X_v)\subset Fl_v(\{z\})$ is clear, we will  prove the opposite inclusion. Suppose first that $v_1, v_2$ are  vertices in the tree of $T_z$ satisfying 
$$
d_T(u,v_2)=d_T(u,v_1)+1,
$$
and let $\alpha_i\subset X_{v_i}$, $i=1,2$, be geodesic rays such that 
$Fl_{v_i}(\{z\})=\{\alpha_i(\infty)\}$, $i=1,2$. According to Lemma \ref{lem: flow condition}, 
$$
\Hd_{X_{v_1v_2}}(\alpha_1,\alpha_2)<\infty.$$ 

Since $X_{v_1v_2}$ is $\delta'_0$-hyperbolic and $\al_1, \al_2$ are $L'_0$-quasigeodesics in 
$X_{v_1v_2}$, it follows that there are subrays in $\al'_1\subset \al_1, \al'_2\subset \al_2$ 
which are $R$-Hausdorff close in $X_{v_1v_2}$, where $R$ depends only on the parameters of $\X$ (see Corollary \ref{cor:asymptotic-quasirays}). Set $K:=R^\wedge$. Then,  
if $\al_1$ is contained in $Fl_K(\al_u)$, then so does $\al'_2$ (see Proposition \ref{prop:flow-prop}(2)). Now, Part (1) of lemma follows by induction on $d_T(u,v)$. 

Part (2) follows immediately from (1). \qed

\subsection{Geometry of the flow-incidence graph}\label{sec:flow-incidence graph}

In this section we will assume that $K\ge K_0$. We will analyze (to some degree) the intersection  pattern of projections 
to $T$ of the flow-spaces $Fl_K(X_u)$, $u\in V(T)$. Some of these results will be important (specifically, Subdivision Lemma and its corollary, Corollary \ref{cor:subdivision}) in Chapter \ref{sec:everything together}. 

\begin{defn}
An interval $J=\llbracket u , v\rrbracket\subset T$ is {\em special} (more precisely, $K$-{\em special}) if one of its end-points (say, $u$) has the property that $J\subset \pi(Fl_K(X_u))$. In this case, the vertex $u$ is said to be {\em special} in $J$. 
\end{defn}

For instance, for every edge $e=[u,v]\in E(T)$, the interval   $J=\llbracket u , v\rrbracket$ is special. 

\begin{rem}
The notion of a special interval can be refined and one can say that an {\em oriented} interval 
$J=\llbracket u , v\rrbracket\subset T$ is special if $J\subset \pi(Fl_K(X_u))$. We will not need this refinement.  
\end{rem}

The importance of special intervals comes from the following simple fact:

\begin{lemma}\label{lem:special-sub}
Suppose that $u, v\in V(T)$ are such that $\pi(Fl_{K}(X_u))\cap \pi(Fl_{K}(X_v))$ contains a vertex $t$. Then the center $w$ of the triangle $\Delta tuv$ subdivides the interval $I= \llbracket u , v\rrbracket$ in two special subintervals
$$
J= \llbracket u , w\rrbracket, J'= \llbracket w , v\rrbracket.
$$
\end{lemma}
\proof Since $w$ separates $u$ and $v$ from $t$, Proposition \ref{prop:flow-prop}(1) implies that 
$$
J\subset \pi(Fl_K(X_u)), J'\subset \pi(Fl_K(X_v)). \qed 
$$

Of prime importance for us will be whether intersections as in the lemma are empty or not, since empty intersection will imply that 
the pair of subsets $Fl_K(X_u)$, $Fl_K(X_v)$  of $X$ is $L_{\ref{mjproj}}(K)$-Lipschitz cobounded, see Corollary \ref{cor:disjoint->cobounded}. 
This observation motivates the following definition:

\begin{defn}\label{defn:flow-graph}
For each $K$ we define the {\em flow-incidence graph} $\Ga=\Ga_K$. Its vertex-set is $V(T)$. Two vertices $u, v$ of $\Ga$ are connected by an edge $e\in E(\Ga)$ 
if and only if  $\pi(Fl_{K}(X_u))\cap \pi(Fl_{K}(X_v))\ne \emptyset$. 
\end{defn}

Lemma \ref{lem:special-sub} implies: 

\begin{cor}
$d_\Ga(u, v)\le 1$ if and only if the interval $\llbracket u , v\rrbracket$ is the union of two special subintervals. 
\end{cor}

While the graph $\Ga$ is not necessarily a tree,  we will see that it is a {\em quasi-tree}. 
Recall that a geodesic metric space $Y$ is called a {\em quasi-tree} (see Manning's paper \cite{Manning}) if it is quasiisometric to a simplicial tree. According to 
\cite{Manning}, $Y$ is a quasi-tree if and only if there exists a constant $r$ such that for every geodesic segment $xy\subset Y$, the (closed) ball $B(m,r)$ centered at the 
midpoint $m$ of $xy$ separates $x$ and $y$. (Here and below, separation is understood in the sense that every path connecting $x$ and $y$ has to pass through $B(z,r)$.) 
Alternatively, one characterizes quasi-trees by the existence of $r$ such 
that for every geodesic segment $xy\subset Y$ and every $z\in xy$, the ball  $B(z,r)$ separates $x$ and $y$.

\begin{lemma}\label{lem:flow-separation}
Suppose that $p$ is a point\footnote{not necessarily a vertex} in an interval $J=\llbracket u , v\rrbracket\subset T$. Then:

1. The closed unit ball $B_\Ga(p,1)\subset \Ga$ separates $u$ and $v$. 

2. $d_\Ga(u,v)\ge d_\Ga(u,w)$.   
\end{lemma}
\proof Both parts of the lemma are proven by the same argument. Suppose that $w_0=u, w_1,...,w_n, w_{n+1}=v$ is a vertex-path in $\Ga$ connecting $u$ to $v$. Since $p$ separates $u$ from $v$, there exists $i\le n$ such that $w_i, w_{i+1}$ are separated 
by $p$ in $T$ and, of course, $d_\Ga(w_i, w_{i+1})=1$. 
Lemma \ref{lem:special-sub} implies that 
$$
d_\Ga(w_i, p)\le 1, d_\Ga(w_{i+1}, p)\le 1. 
$$
This proves Part (1) of the lemma. To prove Part (2) we take the above vertex-path to be geodesic in $\Ga$ and observe that
$$
d_\Ga(p, u)\le d_\Ga(p, w_i) + d_\Ga(w_i, u)\le 1+ d_\Ga(u, w_i) \le d_\Ga(u,v). \qed 
$$

Lemma \ref{lem:flow-separation} thus implies that for every $K$ and a tree of hyperbolic spaces 
spaces $\X$, the graph $\Ga=\Ga_K$ is a quasi-tree with the constant $r=1$. 
 
 \medskip 
We  are now ready to prove the {\em horizontal subdivision lemma}, which will play an important role 
 in Chapter \ref{sec:everything together}, when we establish uniform hyperbolicity of $K$-flows of interval-spaces $X_J$ in $\X$.

\begin{figure}[tbh]
\centering
\includegraphics[width=60mm]{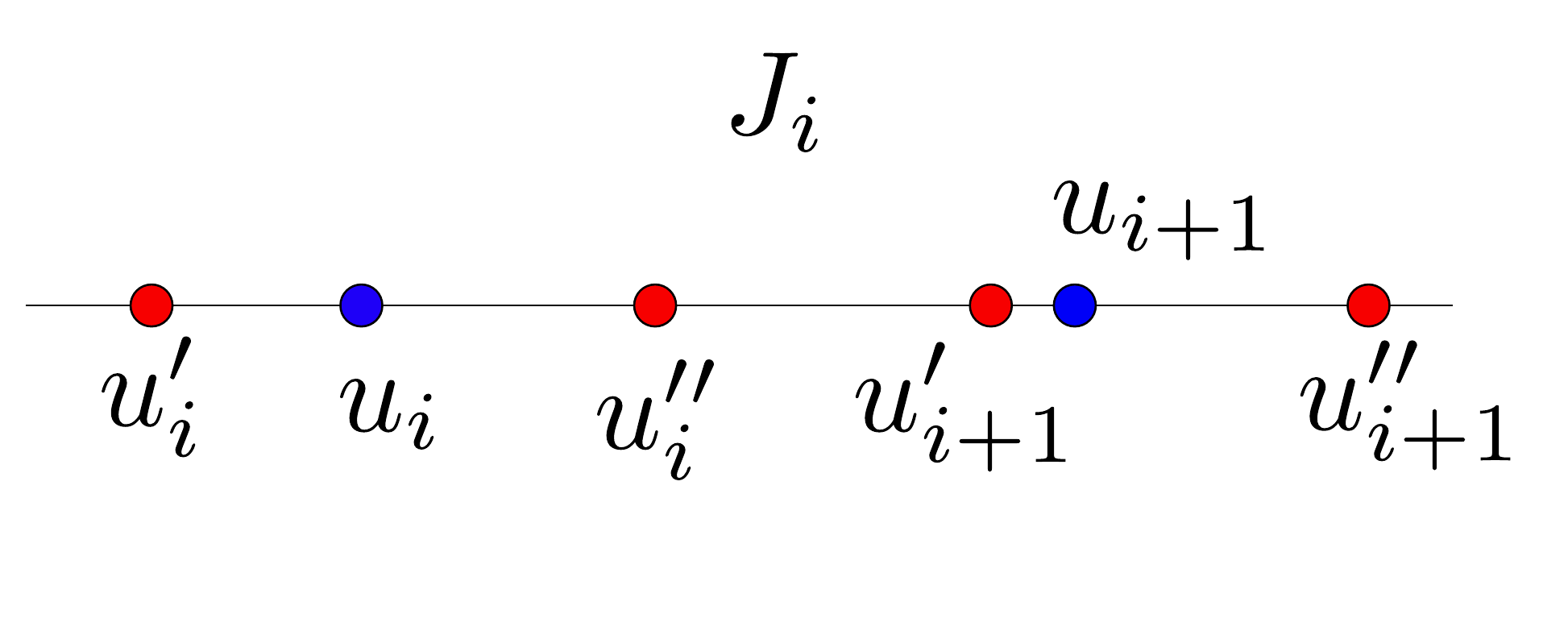} 
\caption{Horizontal subdivision}
\label{h-sub.fig}
\end{figure}

\begin{lemma}
[Horizontal subdivision lemma] 
\label{lem:subdivision}
For any pair of distinct vertices $u, v\in T$, the interval $J=\llbracket u,v\rrbracket$, can be subdivided into 
nondegenerate subintervals 
$$
J = \llbracket u_0,u_1\rrbracket \cup ....\cup \llbracket u_{n-1},u_{n}\rrbracket  \cup \llbracket u_{n},u_{n+1}\rrbracket , \quad u=u_0, v=u_{n+1},
$$
and each $J_i=\llbracket u_i, u_{i+1}\rrbracket $ can be further subdivided into subintervals (some of which could be degenerate),  
$$
\llbracket u_i, u_i''\rrbracket \cup \llbracket u_i'', u'_{i+1}\rrbracket  \cup \llbracket u'_{i+1}, u_{i+1}\rrbracket, 
$$
so that the following hold for all $i\le n$:

\begin{enumerate}
\item\label{I2} 
$$
\pi(Fl_K(X_{u_i}))\cap J_{i}= \llbracket u_i,u_i''\rrbracket,$$
(i.e. the interval $\llbracket u_i,u_i''\rrbracket$ is special) and 
$$
u_i''\notin \pi(Fl_K(X_{u_{i+1}})),$$
unless $i=n$ in which case we {\em could} have $u_i''\in \pi(Fl_K(X_{u_{i+1}}))$. 

\item\label{I3} The interval $\llbracket u''_{i},u'_{i+1}\rrbracket$ is special, it is contained in $\pi(Fl_K(X_{u'_{i+1}}))$.  

\item The interval $\llbracket u_i'', u'_{i+1}\rrbracket $ is nondegenerate unless $i=n$. 

\item $d_T(u'_{i+1}, u_{i+1})\le 1$, thus, each interval $\llbracket u'_{i+1},u_{i+1}\rrbracket$ is special. 
\end{enumerate}
\end{lemma}
\proof We find the subdivision vertices inductively. Set $u_0:=u$. 
Inductively, we assume that $u_i$ is defined. If $u_i=v$, we set $n+1=i$ and terminate the induction. 
Suppose, therefore, that this is not the case. We then define $u_i''$, $u_{i+1}'$ and $u_{i+1}$:

We choose a vertex $u''_i\in \rrbracket u_i,v\rrbracket $ to be the farthest  from $u_i$ such that 
$$
Fl_K( X_{u_i})\cap X_{u''_i} \neq \emptyset.$$ 
Note that such a vertex always exists since for the edge $[u_i,v_i] \in E( \llbracket u_i, v\rrbracket )$ we have
$$
Fl_K( X_{u_i})\cap X_{v_i} \neq \emptyset.
$$
If it so happens that $u''_i=v$, we set $n=i$, and $u'_{i+1}:=u_{i+1}=v$; this will conclude the induction.
Suppose that this is not the case.

Then consider the vertices $s\in \rrbracket u''_i,v\rrbracket $ such that 
$$
Fl_K(X_s)\cap X_{u''_i}=\emptyset.$$ 
If such a vertex does not exist, then we set $n=i$, $u'_{i+1}=u''_i$ and $u_{i+1}=v$, and again conclude the subdivision process. 
Assume that this is not the case. Then we define $u_{i+1}\in \rrbracket u''_i,v\rrbracket $ to be the closest vertex to $u''_i$ such that 
$Fl_K(X_{u_{i+1}})\cap X_{u''_i}=\emptyset$.  We define $u'_{i+1}$ in this case to be the vertex in 
$J_i= \llbracket u_i, u_{i+1}\rrbracket$ adjacent to $u_{i+1}$, i.e. $d_T(u'_{i+1}, u_{i+1})=1$. Then, by the definition of $u_{i+1}$,
$$
Fl_K(u'_{i+1})\cap X_{u_i''}\ne\emptyset. 
$$
Hence, the vertices $u_i'', u_{i+1}', u_{i+1}$ satisfies  requirements of the lemma and we continue inductively. \qed

\begin{cor}\label{cor:subdivision}
1. Each interval $J_i$ as above is the union of (at most) three special intervals and the sequence 
$$
... u_i, u'_{i+1}, u_{i+1},...
$$ 
is a vertex-path in $\Ga$. 

2. For any two consecutive vertices $u_i, u_{i+1}$, $i\le n-1$, 
$$
d_\Ga(u_i, u_{i+1})=2, 
$$
while $1\le d_\Ga(u_n, u_{n+1})\le 2$.

3. For each pair of indices $i, j$, if $0\le i+1<j\ne n$, then 
$$
d_\Ga(u_i, u_j)\ge 2. 
$$
In particular, if $|i- j|\ge 2$, 
the flow-spaces $Fl_K(X_{J_i}), Fl_K(X_{J_j})$ are $L_{\ref{mjproj}}(K)$-Lipschitz cobounded in $X$. 
\end{cor}

\section{Retractions to bundles}\label{sec:retractions}

The main goal of this section is to prove Theorem \ref{thm:flow-to-bundle}, which is an analogue of Theorem \ref{mjproj}, constructing coarse Lipschitz retractions from flow-spaces to certain 
$K'$-metric bundles ${\mathfrak Y}= (\pi: Y\to S)\subset  {\mathfrak X}$ with $\la$-quasiconvex fibers $Y_v\subset X_v$. In Theorem \ref{thm:flow-to-bundle} we will impose a stronger assumption on 
${\mathfrak X}$, namely the 
$\kappa$-uniform flaring condition for a certain constant $\kappa\ge K$ (see \eqref{eq:key-kappa} for the definition of this constant, which depends on $K, K'$, on a quasiconvexity constant $\la$ of $Y_v\subset X_v$,  and on a constant $D$ which is an upper bound on the diameter of $Y_w$ for some $w\in V(S)$), that was not needed in Theorem \ref{mjproj}. While the $\kappa$-flaring condition implies $k$-flaring for all $k\in [1,\kappa]$, it will be notationally convenient to also have the constants $M_k, k\le \kappa$, at our disposal, hence, we will be {\em assuming} the uniform $k$-flaring condition for all $k\in [1,\kappa]$.  

The retractions 
$\rho_{\Y}: Fl_K(Q_u)\to Y$ will be defined on flow-spaces $Fl_K(Q_u)$ whose 
$4\delta_0$-fiberwise neighborhoods contain $Y$, but composing $\rho_{\Y}$ with Mitra's retraction $\rho: X\to Fl_K(Q_u)$, we then obtain retractions defined on the entire $X$. In the special case, when $Fl_K(Q_u)$ is 
$\delta$-hyperbolic, the retraction $\rho_{\Y}$ will be uniformly close to the nearest-point projection $Fl_K(Q_u)\supset Y$ (see Proposition \ref{prop:rho-npp}). 

\begin{rem}
1. The condition that $\rho_{\Y}$ is a retraction should be understood coarsely since $Y$ is not quite contained in 
$Fl_K(Q_u)$: We can only guarantee that $\rho_{\Y}$  fixes all points in $\YY\cap {\mathcal Fl}_K(Q_u)$; the rest of the points of $Fl_K(Q_u)$ lie in the $\max(K,4\delta_0)$-neighborhood of $Y$ and $\rho_{\Y}$ 
can move them only by a uniformly bounded amount. 

2. In the case when ${\mathfrak Y}$ is a $K'$-carpet, which is of the main interest,  
$\la=\delta_0$. 
\end{rem}

\begin{lemma}\label{bundle-proj}
Fix $\la$ and $K'$. Suppose that a subtree of spaces ${\mathfrak Y}\subset  {\mathfrak X}$ is a $K'$-metric bundle over a subtree $S= \pi(Y)\subset T$. Assume also that vertex-spaces $Y_v= Y\cap X_v, v\in V(S)$, are  $\lambda$-quasiconvex in $X_v$. Then the fiberwise nearest-point projection $X_{S}\map Y$ is a 
$D_{\ref{bundle-proj}}(\lambda, K')$-coarse Lipschitz retraction. In particular, if $\gamma$ is a $C$-qi section over some interval $J\subset S$,  
then the fiberwise projection $\bar\ga$ of $\ga$ to $Y$ is a 
$K_{\ref{bundle-proj}}(\lambda, K', C)= C D_{\ref{bundle-proj}}(\lambda,K')$-qi section over $J$, where
$$
K_{\ref{bundle-proj}}(\lambda, K', C) \ge \max(C, K'). 
$$
\end{lemma}

\proof The lemma is an immediate corollary of  Corollary \ref{cor:projection}. \qed

\begin{notation}\label{barK and barbarK} 
1. For the rest of this subsection we will use the notation $\bar{K}$ for 
$K_{\ref{bundle-proj}}(\lambda, K', K)$. We also set  
$\bar{\bar{K}}:=K_{\ref{bundle-proj}}(\delta_0, K'_{\ref{lem:E-ladder-structure}}(K), \bar{K})$. Note that 
$$
\bar{\bar{K}}\ge \max(\bar{K}, K'_{\ref{lem:E-ladder-structure}}(K)). 
$$
In these notation we suppress the dependence on $\la$ and $K'$: In the most useful for us case, when ${\mathfrak Y}$ is a $K$-carpet, we will have $\la=\delta_0$. 

2. Define 
\begin{equation}\label{eq:key-r}
r:=r_{\ref{eq:key-r}}=
3\delta_0 + \la + R_{\ref{cor:super-weak flaring}}(\bar{\bar{K}}, 1) + R_{\ref{cor:super-weak flaring}}(\bar K, M_{\bar K})
\end{equation}
\begin{equation}\label{eq:key-k}
k= K'+ R_{\ref{cor:super-weak flaring}}(\bar K, r)
\end{equation}
\begin{equation}\label{eq:key-kappa}
\kappa:= \kappa_{\ref{eq:key-kappa}}(\la,K, K'):= \max(k,  \bar{\bar K}) 
\end{equation}
\end{notation}
Observe that $\kappa\ge \bar{\bar{K}}\ge \bar K\ge K$. The proof of the following theorem will need uniform 
$\bar{\bar{K}}$-flaring (in numerous places) as well as the uniform $(K'+ R_{\ref{cor:super-weak flaring}}(\bar K, r_i))$-flaring for some numbers $r_1, r_2, r_3$ (subcases (i), (ii) and (iii) respectively in the proof of Proposition \ref{lem:proj-to-bundle}); the  constant $r$ above is chosen to be the maximum of the numbers $r_1, r_2, r_3$. 

\begin{theorem}\label{thm:flow-to-bundle}
Fix constants $K, K'$ and $\la$ and assume that, $\Y=(\pi: Y\to S)\subset \X$ is a $K'$-metric bundle with $\la$-quasiconvex fibers. We assume, furthermore, that 

1. There exists a vertex $w\in S$ such that $\diam_{X_w}(Y_w)\le D'$. 

2. The tree of spaces ${\mathfrak X}$ satisfies the uniform flaring condition for all parameters in the interval $[1,\kappa]$, in particular, for $\bar{K}, K'$ and $\bar{\bar{K}}$; see 
\eqref{eq:key-kappa} for the definition of $\kappa$ which depends on $K, K'$ and $\la$. 

3. We assume that $\Y$ is either a $(K',D')$-carpet $\A(\al')$ contained in a $K$-ladder\footnote{The other two parameters, $D, E$, of $\L$ play no role in this theorem.} ${\mathfrak Z}=\L_K(\al)$, 
$$
\al'\subset \al\subset X_u, \length(\al')\ge \length(\al)- M_{\bar{K}},
$$
or $\Y$ is a general $K'$-metric bundle contained in the fiberwise $4\delta'_0$-neighborhood of 
a $K$-flow-space ${\mathfrak Z}= {\mathfrak Fl}_K(Q_u)$, such that $Y_u=Q_u$. 

In both cases, we let $Z$ denote the total space of ${\mathfrak Z}$ and $\ZZ:= Z\cap \XX$. 

Then there exists a coarse $L_{\ref{thm:flow-to-bundle}}(\lambda,K,K',D')$-Lipschitz retraction 
$$
\rho=\rho_{Y}: \ZZ\to \YY.$$ 
\end{theorem}
\proof The proof of this theorem is quite long and technical; it occupies most of the rest of this section. 

{\bf Step 1:} {\bf Construction of the map $\rho: \ZZ\map \YY$.} 

For each $x\in \ZZ$ with $v=\pi(x)\in V(T)$, we 
fix a $K$-qi section $\gamma_x$ of $\pi: Z\to \pi(Z)\subset T$ over 
$\llbracket u, v\rrbracket $, connecting $x$ to some point in $Z_u$ once and for all. 

Let $b=b_x$ be the nearest point projection of $v$ to $S$ in $T$.  We define the following important points:

\begin{itemize}
\item We let $t=t_x\in \llbracket u, b\rrbracket $ be the vertex farthest from $u$ such that there exists a point $\tilde{x}\in \ga_x(t)$ for which 
$$
d_{X_t}(\tilde{x}, Y_t)\le M_{\bar{K}}. 
$$
(Note that  it is possible that $t=u$ and $\tilde{x}\in Z_u$.)  

\item Let $\bar{x}\in Y_t$ be a nearest-point projection to $\tilde{x}$ to $Y_t$
in the vertex space $X_t$. 
\end{itemize}

Thus, 
\begin{equation}\label{eq:tildex}
d_{X_{t_x}}(\tilde{x}, \bar{x})\le M_{\bar{K}}
\end{equation}
and if $x\in Y_v$, then $t_x=v$ and $\bar{x}=\tilde{x}=x$.

\begin{figure}[tbh]
\centering
\includegraphics[width=60mm]{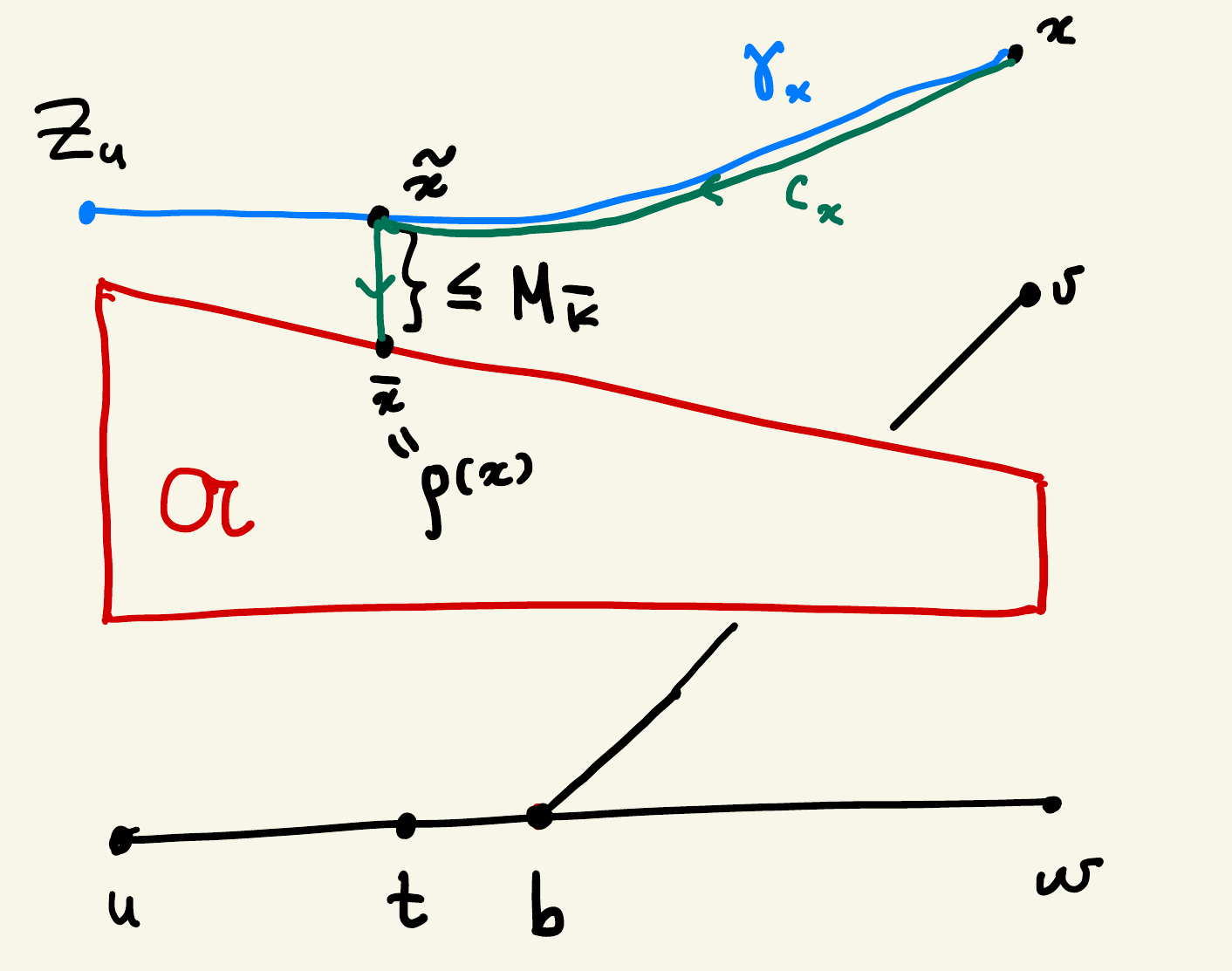}  
\caption{Projection to a bundle which is a  carpet}
\label{proj-to-carpet.fig}
\end{figure}

\begin{defn}\label{defn:retraction-to-bundle} 
1. We define the retraction $\rho=\rho_{\Y}: \ZZ \to \YY$ by $\rho(x):= \bar{x}$. 
We extend this map to $Z$ using the fact that 
${\mathcal Z}$ is a $K$-net in 
$Z$ (to define the extension we compose a nearest-point projection with $\rho$). See Figure \ref{proj-to-carpet.fig}. 

2. We define  a path $c_x= c_{x,Y}$ connecting $x$ to $\rho(x)\in Y$  
as the concatenation 
$$
\ga_{x,\tilde{x}}\star [\tilde{x} \bar{x}]_{X_{t}},$$
where $t=t_x$, 
and 
$$
\ga_{x,\tilde{x}}=  \ga_x|_{\llbracket v, t\rrbracket }, 
$$
is  the subpath of $\ga_x$ connecting $x$ to $\tilde{x}$. We will refer to the latter as the 
{\em horizontal part} of  $c_x$. The {\em vertical} part of $c_x$ is the geodesic    $[\tilde{x} \bar{x}]_{X_{t}}$;  
it is a path of uniformly bounded length (see \eqref{eq:tildex}) connecting $\tilde{x}$ to $\bar{x}$ and contained in the $4\delta_0$-neighborhood of $Z_t= Z\cap X_t$. 
\end{defn}

\medskip

{\bf Step 2: Verification of the properties of $\rho$.} It suffices to verify the coarse Lipschitz property for the restriction of $\rho$ to $\ZZ$. We note further that it is enough to get a uniform 
upper bound on $d(\rho(x), \rho(y))$ for two types of 
pairs $(x, y)$:

a. $x, y\in Z_v$, $d_{X_v}(x,y)\le 1$. 

b. The vertices $v_1=\pi(x), v_2=\pi(y)$ span an edge in $T$ and $d_{X_{v_1v_2}}(x,y)\le K$.

\noindent These two cases are treated in Proposition \ref{lem:proj-to-bundle} and Lemma \ref{lem:proj-to-bundle1} respectively. The former  is the longest and hardest part of the proof. 

For the following proposition we observe that, according to the 3rd assumption of Theorem \ref{thm:flow-to-bundle}, 
$\X$ satisfies the uniform $k$-flaring condition 
for 
$$
k= K'+ R_{\ref{cor:super-weak flaring}}(\bar K, r),
$$
where
$$
r= \max( M_{\bar{K}}, 
6\delta_0 +\la+ R_{\ref{cor:super-weak flaring}}(\bar{\bar{K}}, 1),  
3\delta_0+ \la + R_{\ref{cor:super-weak flaring}}(\bar{\bar{K}}, 1) + 
 R_{\ref{cor:super-weak flaring}}(\bar{K}, M_{\bar{K}}) ), 
$$
where 
$$
r= 
3\delta_0+ \la + R_{\ref{cor:super-weak flaring}}(\bar{\bar{K}}, 1) + 
 R_{\ref{cor:super-weak flaring}}(\bar{K}, M_{\bar{K}}), 
$$
since $R_{\ref{cor:super-weak flaring}}(\bar{K}, M_{\bar{K}})\ge \max(3\delta_0, M_{\bar{K}})$.

\begin{prop}\label{lem:proj-to-bundle}
Suppose that $x, y\in \ZZ\cap X_v$. If $x,y$ are within  
distance $1$ from each other in $X_v$, then $d_X(\rho(x), \rho(y))\le C_{\ref{lem:proj-to-bundle}}(\la,K,K',D')$. 
The bound is independent of the choice of the paths $\ga_x$, $\ga_y$ as above.  
\end{prop}
 The proof of this proposition is long and is done through analyzing several cases and subcases. 
 We use the notation preceding Definition \ref{defn:retraction-to-bundle} and note that in the setting of the proposition, 
$b_x=b_y$; we will denote this vertex simply by $b$.

Now we define certain auxiliary  objects and make  general remarks to be used in the proof, especially in Cases 2 and 3 below. We let $z$ denote the nearest point projection (in $T$) of $b$ to $\llbracket u,w\rrbracket $, i.e. $z$ is the center of the triangle $\Delta uwb$. 

\begin{rem}\label{rem:bundle-proj}
\begin{enumerate}
\item For each vertex $s\in T$, every geodesic $\alpha\subset X_s$  is $\delta_0$-quasiconvex in $X_s$.
Hence, the nearest-point projection (in $X_s$) to $\alpha$ is coarsely $L_{\ref{lip-proj}}(\delta_0, \delta_0)$-Lipschitz. 

\item 
 We let $\bar{\gamma}_x$ and $\bar{\gamma}_y$ denote the fiberwise projections (to $Y$) of the restrictions to $\llbracket b,u\rrbracket $ of 
$\gamma_x$ and $\gamma_y$ respectively. These are $\bar{K}$-qi sections over $\llbracket b,u\rrbracket $. 
(See Lemma \ref{bundle-proj} and Notation \ref{barK and barbarK}.) Then, by the definition of $t_x, t_y$,  
$$
d_{X_s}(\ga_x(s), \bar\ga_x(s))> M_{\bar{K}}, \forall s\in V(\rrbracket t_x, b\rrbracket ), 
d_{X_{t_x}}(\ga_x(t_x), \bar\ga_x(t_x))\le M_{\bar{K}}, 
$$
$$
d_{X_s}(\ga_y(s), \bar\ga_y(s))> M_{\bar{K}}, \forall s\in V(\rrbracket t_y, b\rrbracket )
d_{X_{t_y}}(\ga_y(t_y), \bar\ga_y(t_y))\le M_{\bar{K}}.
$$

\item 
Since ${\bar{K}}\ge {K}$ and $\ga_x(u)=\bar\ga_x(u)$, by \flaring~we have  
$$
d_{X_s}(\gamma_x(s), \bar{\gamma}_x(s))\leq R_{\ref{cor:super-weak flaring}}(\bar{K}, M_{\bar{K}}),$$
 for all $s\in V(\llbracket t_x,u\rrbracket )$ and, similarly, 
$$
d_{X_s}(\gamma_y(s), \bar{\gamma}_y(s))\leq R_{\ref{cor:super-weak flaring}}(\bar{K}, M_{\bar{K}}), \quad \forall s\in V(\llbracket t_y,u\rrbracket ).
$$
\item 
The carpet ${\mathfrak A}$ bounded by $\gamma_x$ and $\gamma_y$ (with the narrow end $[\ga_x(v) \ga_y(v)]_{X_v}$) 
is a $K'_{\ref{lem:E-ladder-structure}}({K})$-carpet over $\llbracket u,v\rrbracket$; in particular, it 
is a $K'_{\ref{lem:E-ladder-structure}}({K})$-metric bundle 
whose fibers are $\delta_0$-quasiconvex in the corresponding vertex spaces. 
We let $\bar{\bar{\gamma}}_x$ and $\bar{\bar{\gamma}}_y$ denote, respectively, the 
 fiberwise projections of $\bar{\gamma}_x$ and $\bar{\gamma}_y$ to $\A$.  
 Thus, by 
Lemma \ref{bundle-proj}, both $\bar{\bar{\gamma}}_x$ and $\bar{\bar{\gamma}}_y$ are 
$\bar{\bar{K}}=K_{\ref{bundle-proj}}(\delta_0, K'_{\ref{lem:E-ladder-structure}}(K), \bar{K})$-qi sections over $\llbracket u,b\rrbracket $.

\item 
Since $\A$ is a $K'_{\ref{lem:E-ladder-structure}}(K)$-carpet, 
we can join the points $\bar{\bar{\gamma}}_x(b)$ and $\bar{\bar{\gamma}}_y(b)$ to some points of $[xy]_{X_v}$ in  $\A$  via $K'_{\ref{lem:E-ladder-structure}}(K)$-qi sections over $\llbracket v,b\rrbracket $. 
Since $\bar{\bar{K}}\ge K'_{\ref{lem:E-ladder-structure}}(K)$, the concatenation of these qi sections with sections 
$\bar{\bar{\gamma}}_x$ and $\bar{\bar{\gamma}}_y$ (over $\llbracket u,b\rrbracket $) 
are both $\bar{\bar{K}}$-qi sections over 
$\llbracket u, v\rrbracket $, joining $\gamma_x(u)$ and $\gamma_y(u)$ to some points of $[xy]_{X_v}$. 
We retain the notation $\bar{\bar{\gamma}}_x$ and $\bar{\bar{\gamma}}_y$ for these concatenations.

\item Notice that $\bar{\bar{\gamma}}_x(u)= \bar\ga_x(u)= \ga_x(u)$ and $\bar{\bar{\gamma}}_y(u)= \bar\ga_y(u)= \ga_y(u)$. At the same time, the other pair of 
end-points (namely, points in $[xy]_{X_v}$) of $\gamma_x$ and $\bar{\bar{\gamma}}_x$, respectively, of 
$\gamma_y$ and $\bar{\bar{\gamma}}_y$, are within  distance $\le 1$ in $X_v$. Therefore, by \flaring, 
we have 
$$
d_{X_s}(\gamma_x(s), \bar{\bar{\gamma}}_x(s))\leq R_{\ref{cor:super-weak flaring}}(\bar{\bar{K}}, 1), \forall s\in V(\llbracket v,u\rrbracket )$$
and 
$$
d_{X_s}(\gamma_y(s), \bar{\bar{\gamma}}_y(s))\leq R_{\ref{cor:super-weak flaring}}(\bar{\bar{K}}, 1), \forall s\in V(\llbracket v,u\rrbracket ).$$ 
\end{enumerate}
\end{rem}

\begin{figure}[tbh]
\centering
\includegraphics[width=60mm]{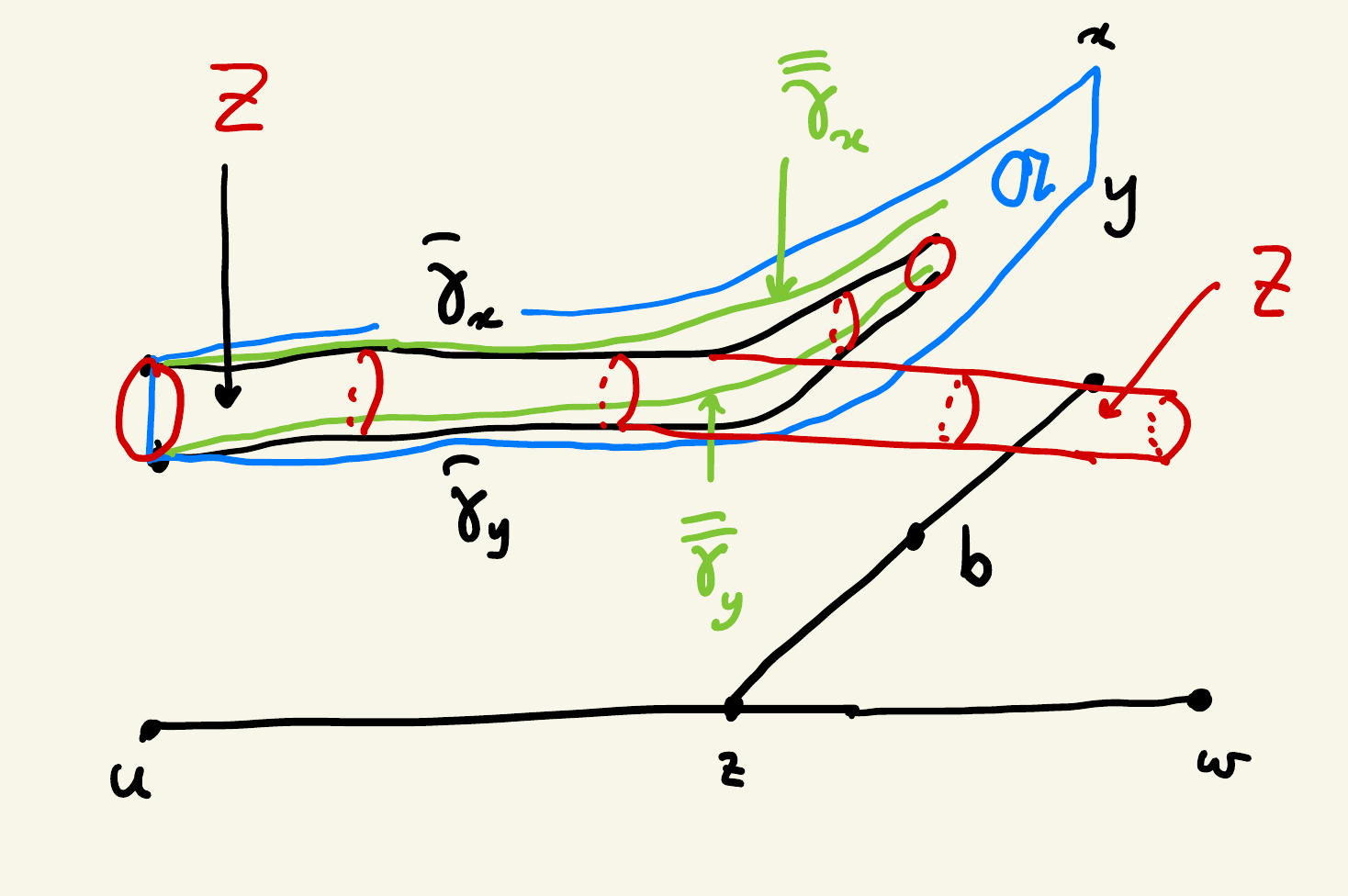}  
\caption{Projections $\bar\ga_x, \bar\ga_y, \bar{\bar{\ga}}_x,  \bar{\bar{\ga}}_y$.}
\label{four-projections.fig}
\end{figure}

Before proving the proposition we will need a technical lemma: 

\begin{lemma}\label{case3-bundle-proj}
Suppose that $r$ is such that ${\mathfrak X}$ satisfies the uniform $k$-flaring condition 
with  
$$
k= k_{\ref{case3-bundle-proj}}= K'+ R_{\ref{cor:super-weak flaring}}(\bar K, r).$$
Then the following holds. 

Suppose that there are vertices $v_1\in \llbracket t_x,b\rrbracket \cap \llbracket z,b\rrbracket$ and 
$v_2\in \llbracket t_y, b\rrbracket \cap \llbracket z,b\rrbracket$  such that
\begin{equation}\label{eq:case3-bundle-proj}
d_{X_{v_1}}(\gamma_x(v_1), \bar{\gamma}_x(v_1))\leq r ~~\hbox{and} ~~d_{X_{v_2}}(\gamma_y(v_2), \bar{\gamma}_y(v_2))\leq r . 
\end{equation}
Then:

(i) $d_T(v_1,v_2)\leq \tau_{\ref{case3-bundle-proj}}=\tau_{\ref{case3-bundle-proj}}(K,r,D')$. 

(ii) $d(\gamma_x(v_1), \gamma_y(v_2))\leq R_{\ref{case3-bundle-proj}}=R_{\ref{case3-bundle-proj}}(K,r,D')$. 

(iii)  
$$
 d_{X_s}(\gamma_x(s), \gamma_y(s))\leq 
 R_{\ref{cor:super-weak flaring}}(k, \max(1, D')), 
 \forall s\in \llbracket v,z\rrbracket. 
$$
\end{lemma}
\proof (i) Taking into account the fact that $K\le \bar{K}\le \kappa$, 
$\ga_x(u)=\bar\ga_x(u)$, $\ga_y(u)=\bar\ga_y(u)$, as well as the inequalities \eqref{eq:case3-bundle-proj}, we see that   \flaring  \ applied to $\gamma_x|_{\llbracket u,v_1\rrbracket }$ and $\bar{\gamma}_y|_{\llbracket u,v_2\rrbracket }$ and the vertex 
$z\in \llbracket u, v_1\rrbracket  \cap \llbracket u, v_2\rrbracket $,  
implies that 
$$
d_{X_z}(\gamma_x(z), \bar{\gamma}_x(z))\leq R_{\ref{cor:super-weak flaring}}(\bar{K}, r)$$
and 
$$d_{X_z}(\gamma_y(z), \bar{\gamma}_y(z))\leq R_{\ref{cor:super-weak flaring}}(\bar{K}, r).$$

We  join $\bar{\gamma}_x(z)$ and $\bar{\gamma}_y(z)$
to $Y_w$ by two $K'$-qi sections over $\llbracket z,w\rrbracket $ contained in $Y$, $\gamma_{x,1}$ and $\gamma_{y,1}$ respectively.

We let $z_1$ denote the vertex in $\llbracket z,w\rrbracket $ adjacent to $z$, assuming that $z\ne w$. We connect $\gamma_x(z)$ and $\gamma_{x,1}(z_1)$ by a geodesic path $\gamma_{1,z,z_1}$ 
 in $X_{zz_1}$, similarly, connect   $\gamma_y(z)$ and $\gamma_{y,1}(z_1)$ by a geodesic path $\gamma_{2,z,z_1}$ in $X_{zz_1}$. 
 Both paths have length $\le k=K'+  R_{\ref{cor:super-weak flaring}}(\bar{K}, r)$.

\begin{figure}[tbh]
\centering
\includegraphics[width=60mm]{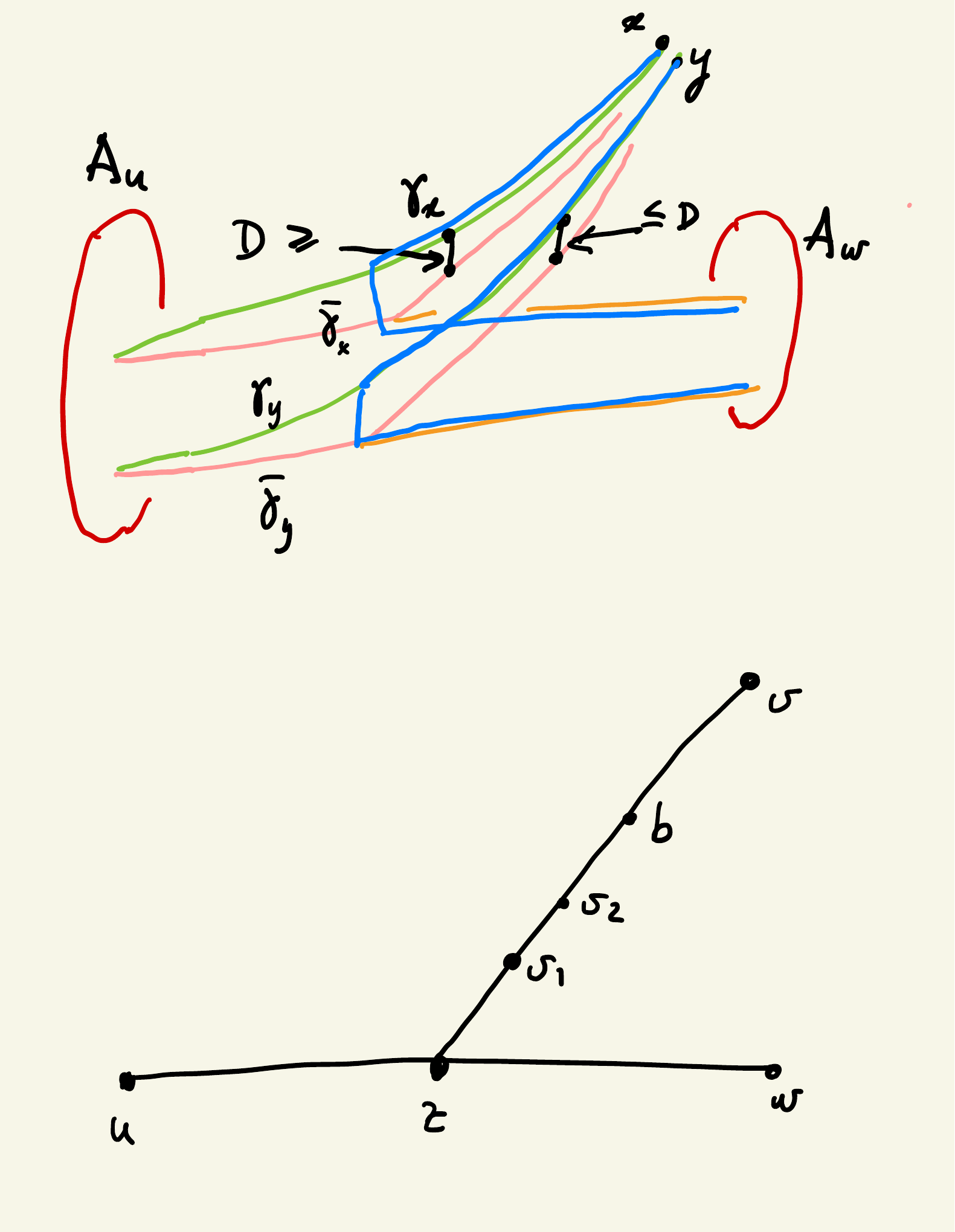}
\caption{Illustration of the proof of Lemma \ref{case3-bundle-proj}}
\label{lemma2.51.fig}
\end{figure}

 Then we get that the concatenations 
 $$
 \ga'_x:= \gamma_x|_{\llbracket v,z\rrbracket }\star \gamma_{1,z,z_1} \star \gamma_{x,1}|_{\llbracket z_1,w\rrbracket }$$
  and 
  $$
\ga'_y:=  \gamma_y|_{\llbracket v,z\rrbracket }\star \gamma_{2,z,z_1}\star \gamma_{y,1}|_{\llbracket z_1,w\rrbracket }.$$ 
These are 
$k$-qi sections over $\llbracket v,w\rrbracket $. (See the bold paths in Figure \ref{lemma2.51.fig}.) 
Their end-points are at a distance at most $\max(1, D')$ of the respective vertex-spaces
 (since $d_{X_v}(x,y)\le 1$ and $Y_w$ has diameter $\le D'$).

 Hence (since ${\mathfrak X}$ is assumed to satisfy the uniform $k$-flaring condition), by \flaring, 
 $$
 d_{X_s}(\gamma'_x(s), \gamma'_y(s))\leq R_{\ref{cor:super-weak flaring}}(k, \max(1, D')), \forall s\in \llbracket v,w\rrbracket 
 $$
 and, restricting to the subinterval $\llbracket v,z\rrbracket$, we obtain 
 \begin{equation}\label{eq:case3-bundle-proj2}
 d_{X_s}(\gamma_x(s), \gamma_y(s))\leq R_{\ref{cor:super-weak flaring}}(k, \max(1, D')), 
 \forall s\in \llbracket v,z\rrbracket. 
 \end{equation}
In the case $z=w$, we will use the paths 
$$
\ga'_x:= \gamma_x|_{\llbracket v,z\rrbracket }, \ga'_y:=  \gamma_y|_{\llbracket v,z\rrbracket }$$
and obtain the same inequality \eqref{eq:case3-bundle-proj2}. 

Thus, we established Part (iii) of the lemma. 

Since the fiberwise projection to $Y_s, s\in \llbracket b,z\rrbracket$ is 
$L_{\ref{lip-proj}}(\delta_0, \la)$-Lipschitz, we get
 \begin{equation}\label{eq:case3-bundle-proj3}
d_{X_s}(\bar\gamma_x(s), \bar\gamma_y(s))\leq L_{\ref{lip-proj}}(\delta_0, \la)(R_{\ref{cor:super-weak flaring}}(k, \max(1, D'))+1),  \forall s\in \llbracket b,z\rrbracket. 
 \end{equation}
Without loss of generality (by switching the roles of $x$ and $y$ if necessary), 
we may assume that $v_2$ is a vertex in $\llbracket v_1,b\rrbracket$.  

Combining the second inequality in  \eqref{eq:case3-bundle-proj} with the inequalities 
\eqref{eq:case3-bundle-proj2}, \eqref{eq:case3-bundle-proj3} applied to $s=v_2$, by the triangle inequality in $X_{v_2}$ we obtain:
 \begin{align}\label{eq:case3-bundle-proj4}
 \begin{split}
d_{X_{v_2}}(\gamma_x(v_2),\bar{\gamma}_x(v_2))\leq\\ R_1:=r+ 
R_{\ref{cor:super-weak flaring}}(k, \max(1, D'))+ 
L_{\ref{lip-proj}}(\delta_0, \la)(R_{\ref{cor:super-weak flaring}}(k, \max(1, D'))+1).
\end{split}
 \end{align}
Hence, by taking into account the fact that 
$$
d_{X_s}(\gamma_x(s),\bar{\gamma}_x(s))> M_{\bar{K}}
$$
for all vertices $s$ of $\rrbracket t_x, b\rrbracket $ (and the reverse inequality at $t_x$ and the 
inequality \eqref{eq:case3-bundle-proj4}) and using the 
$\bar{K}$-uniform flaring property of the sections $\ga_x, \bar\ga_x$ over the interval 
$\llbracket t_x, v_2\rrbracket$, we obtain 
$$
d_T(v_1, v_2)\leq d_T(t_x, v_1)\le \tau_{\ref{case3-bundle-proj}}(K, r, D'):= 
\tau_{\ref{prop:weak flaring}}(\bar{K}, R_1).$$ 
This concludes the proof of Part (i).

\medskip 
(i)$\Rightarrow$(ii): $d(\gamma_x(v_1),\gamma_y(v_2))$ is bounded by the length of the concatenation of the
paths $\gamma_x|_{\llbracket v_1,v_2\rrbracket }$ (whose length is estimated by (i)) and 
$[\gamma_x(v_2)\gamma_y(v_2)]_{X_{v_2}}$ 
(whose length is estimated by \eqref{eq:case3-bundle-proj2} since $v_2\in \llbracket b,z\rrbracket $) 
is which therefore at most
$$
R_{\ref{case3-bundle-proj}}(K,r,D'):= 
K\cdot \tau_{\ref{case3-bundle-proj}}(K,r,D')+R_{\ref{cor:super-weak flaring}}(k, \max(1, D')). \qedhere$$

We will be using this lemma in the proof of Proposition \ref{lem:proj-to-bundle} with $r=r_1, r=r_2, r=r_3$ (defined below), where $r_2$ is the largest of the three parameters. The uniform flaring condition made in the statement of the Proposition ensures that the uniform $k$-flaring condition in the lemma is satisfied for $r=r_2$ and, hence, for 
the two other (smaller) values of $r$. 
 
\medskip
{\em Proof of Proposition \ref{lem:proj-to-bundle}.} There are several cases to consider depending on mutual position of various vertices in $T$. 

\medskip
{\bf Case 1:} Suppose $v\in V(\llbracket u,w\rrbracket )$, in which case $b=v$. 
Without loss of generality we may assume that   $t_y\in V(\llbracket t_x,v\rrbracket )$. 

\medskip
Consider the subset 
$$
W=\{s\in V(\llbracket t_y,v\rrbracket ): d_{X_{s}}(\bar{\gamma}_x(s), \bar{\gamma}_y(s))\geq D'_{\ref{proj-geod}}(\delta_0,\delta_0)\}.$$
If $W\neq \emptyset$, let $v_1$ be the farthest vertex from $t_y$ in $W$. If $W=\emptyset$ then define $v_1$ to be $t_y$. In other words, $v_1=\sup W$  in the (oriented) interval 
$\llbracket t_y, v\rrbracket $.

\begin{figure}[tbh]
\centering
\includegraphics[width=60mm]{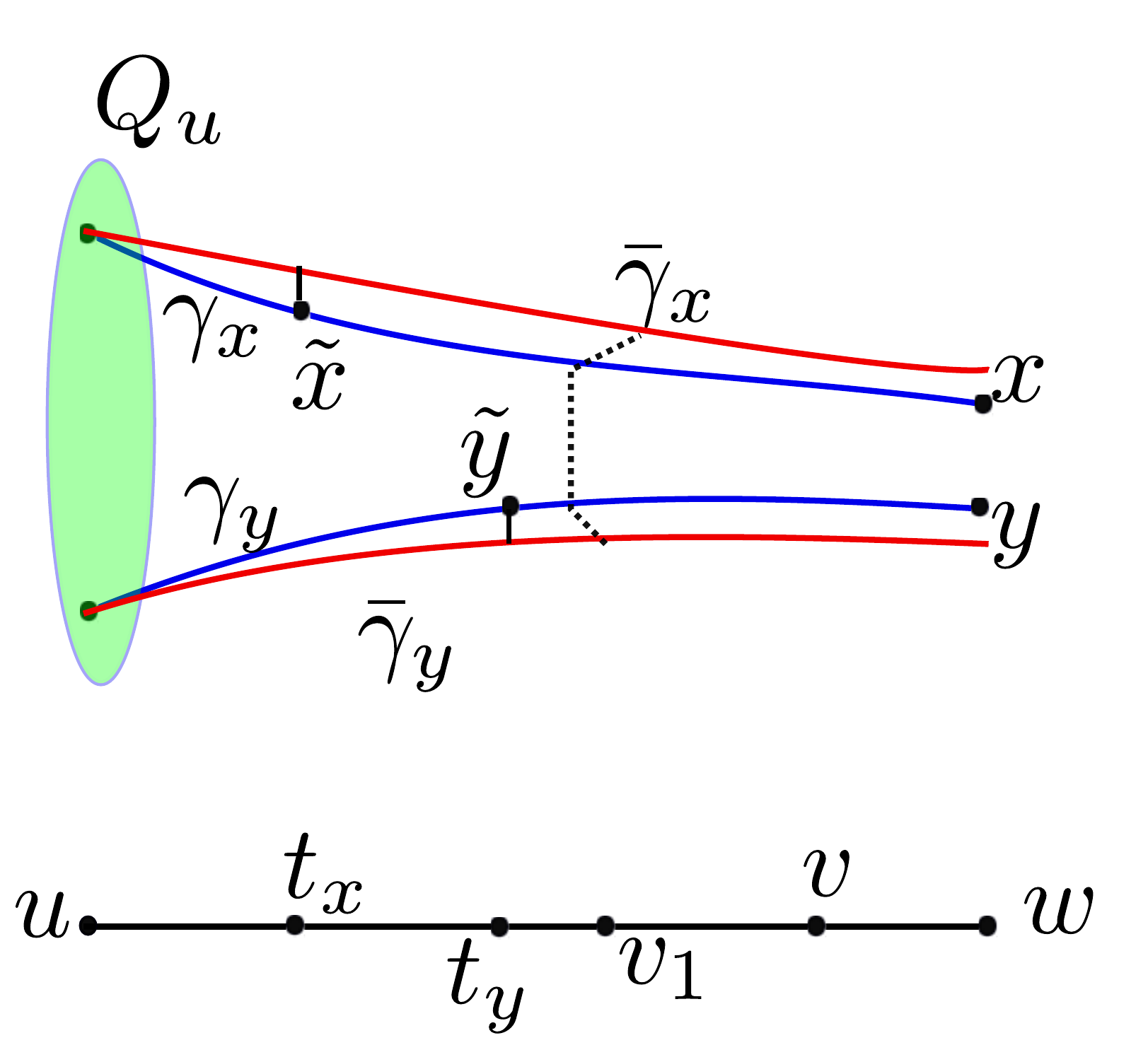}     
\caption{Case 1.}
\label{case1-prop2.50.fig}
\end{figure}

\begin{claim}\label{claim:Case1} 
All three distances 
$$
d_{X_{v_1}}(\gamma_x(v_1), \gamma_y(v_1)), 
d_{X_{v_1}}(\gamma_x(v_1), \bar{\gamma}_x(v_1)), 
d_{X_{v_1}}(\gamma_y(v_1), \bar{\gamma}_y(v_1))$$ 
are bounded above by 
\begin{align*}
C_{\ref{claim:Case1}}= 
2R_{\ref{cor:super-weak flaring}}(\bar{\bar{K}}, 1)+ \\
\max( M_{\bar{K}}+  \delta_0(21 + 72 \delta_0),  9 \delta_0(1+ D_{\ref{bdd-flaring}}(\bar{K},9\delta_0)) ).
\end{align*}

See bold curves in Figure \ref{case1-prop2.50.fig}. 
\end{claim} 
\proof The proof is divided into two subcases. 

{\bf Subcase (i):} Suppose that $W=\emptyset$, thus, $v_1=t_y\in \llbracket u, v\rrbracket$. In this case 
$$
d_{X_{v_1}}(\bar{\gamma}_x(v_1), \bar{\gamma}_y(v_1))< D_{\ref{proj-geod}}(\delta_0,\delta_0).
$$
By Remark \ref{rem:bundle-proj}(1),  
$$
d_{X_{v_1}}(\bar{\bar{\gamma}}_x(v_1)), \bar{\bar{\gamma}}_y(v_1))\leq L_{\ref{lip-proj}}(\delta_0, \delta_0)(1 + D_{\ref{proj-geod}}(\delta_0,\delta_0)).
$$ 
By combining this inequality with the two inequalities 
in Remark \ref{rem:bundle-proj}(6), we obtain  (by the triangle inequality in $X_{v_1}$) 
\begin{align}\label{eq:xy}
\begin{split}
d_{X_{v_1}}(\gamma_x(v_1), \gamma_y(v_1))\leq 2R_{\ref{cor:super-weak flaring}}(\bar{\bar{K}}, 1)+ 
L_{\ref{lip-proj}}(\delta_0, \delta_0)(1 + D_{\ref{proj-geod}}(\delta_0,\delta_0))= \\
2R_{\ref{cor:super-weak flaring}}(\bar{\bar{K}}, 1)+ 
12\delta_0(1 + 9\delta_0),
\end{split}
\end{align}
 c.f. Remark \ref{rem:lip-proj} and Lemma \ref{proj-geod}(2).

Since $v_1=t_y$, 
\begin{equation}\label{eq:y}
d_{X_{v_1}}(\gamma_y(v_1), \bar{\gamma}_y(v_1))\leq M_{\bar{K}}
\end{equation}
and  it  follows from 
the triangle inequality applied to the quadrilateral in $X_{v_1}$ with
the vertices 
$$\gamma_x(v_1),\gamma_y(v_1),\bar{\gamma}_y(v_1), \bar{\gamma}_x(v_1),$$
that 
\begin{align*}\label{eq:x} 
\begin{split}
d_{X_{v_1}}(\gamma_x(v_1), \bar{\gamma}_x(v_1))\leq 
2R_{\ref{cor:super-weak flaring}}(\bar{\bar{K}}, 1)+ 
L_{\ref{lip-proj}}(\delta_0, \delta_0)(1 + D_{\ref{proj-geod}}(\delta_0,\delta_0)) + 
M_{\bar{K}}\\  +D_{\ref{proj-geod}}(\delta_0,\delta_0)= \\
2R_{\ref{cor:super-weak flaring}}(\bar{\bar{K}}, 1)+ M_{\bar{K}}+  \delta_0(21 + 72 \delta_0)
\end{split}
\end{align*} 
(see Remark \ref{rem:lip-proj}). By combining this inequality with  \eqref{eq:y} and \eqref{eq:xy}, we obtain the upper bound
\begin{align*}
\max \{d_{X_{v_1}}(\gamma_x(v_1), \gamma_y(v_1)),  d_{X_{v_1}}(\gamma_x(v_1), \bar{\gamma}_x(v_1)), d_{X_{v_1}}(\gamma_y(v_1), \bar{\gamma}_y(v_1)) \} \le \\
2R_{\ref{cor:super-weak flaring}}(\bar{\bar{K}}, 1)+ M_{\bar{K}}+  \delta_0(21 + 72 \delta_0)\le C_{\ref{claim:Case1}}. 
\end{align*} 

This proves the  inequality in the claim in the subcase (i).

\medskip 
{\bf Subcase (ii):} Suppose $W\neq \emptyset$; hence, 
$$d_{X_{v_1}} (\bar{\gamma}_x(v_1), \bar{\gamma}_y(v_1))\ge D_{\ref{proj-geod}}(\delta_0,\delta_0)=9\delta_0.$$ 
By Lemma \ref{proj-geod}, 
$$
\max\left( d_{X_{v_1}}(\bar{\gamma}_x(v_1), \bar{\bar{\gamma}}_x(v_1)), d_{X_{v_1}}(\bar{\gamma}_y(v_1), \bar{\bar{\gamma}}_y(v_1)) \right) 
\leq R_{\ref{proj-geod}}(\delta_0,\delta_0)=6\delta_0.   
$$
Combining this inequality with  Remark \ref{rem:bundle-proj}(6) and the triangle inequality, we obtain 
\begin{align*}
d_{X_{v_1}}(\bar{\gamma}_x(v_1), \gamma_x(v_1))\leq d_{X_{v_1}}(\bar{\gamma}_x(v_1), \bar{\bar{\gamma}}_x(v_1))+
d_{X_{v_1}}(\bar{\bar{\gamma}}_x(v_1),\gamma_x(v_1))\\
\leq 6\delta_0 + R_{\ref{cor:super-weak flaring}}(\bar{\bar{K}}, 1)
\end{align*} 
and, similarly, 
$$
d_{X_{v_1}}(\bar{\gamma}_y(v_1), \gamma_y(v_1))\leq 6\delta_0 + R_{\ref{cor:super-weak flaring}}(\bar{\bar{K}}, 1)\le C_{\ref{claim:Case1}}.$$
This establishes two out of three bounds in the claim. 

Lastly, we get a bound on $d_{X_{v_1}}(\gamma_x(v_1), \gamma_y(v_1))$. If $v_1=v$ then the bound definitely holds since 
$d_{X_{v}}(\gamma_x(v), \gamma_y(v)) =d_{X_v}(x,y)\le 1$ by the assumption of the proposition. 

Otherwise, by the definition of $v_1$, if $v_2$ is the vertex in $\llbracket v_1, v\rrbracket$ adjacent to $v_1$, then 
$$
d_{X_{v_2}}(\bar{\gamma}_x(v_2), \bar{\gamma}_y(v_2))<  9\delta_0=D_{\ref{proj-geod}}(\delta_0,\delta_0).$$ 
Hence, by  Corollary \ref{bdd-flaring}, 
we have 
$d_{X_{v_1}}(\bar{\gamma}_x(v_1), \bar{\gamma}_y(v_1))\leq D_{\ref{bdd-flaring}}(\bar{K},9\delta_0)$.
By Remark \ref{rem:bundle-proj}(1)
$$ 
d_{X_{v_1}}(\bar{\bar{\gamma}}_x(v_1)), \bar{\bar{\gamma}}_y(v_1))\leq 9\delta_0(1+ D_{\ref{bdd-flaring}}(\bar{K},9\delta_0)).
$$ 
Hence, using Remark \ref{rem:bundle-proj}(6), we obtain 
$$
d_{X_{v_1}}(\gamma_x(v_1), \gamma_y(v_1))\leq \\
2R_{\ref{cor:super-weak flaring}}(\bar{\bar{K}}, 1)+ 9 \delta_0(1+ D_{\ref{bdd-flaring}}(\bar{K},9\delta_0))\le C_{\ref{claim:Case1}}.
$$
This completes the proof of the claim. \qed

We now prove  the proposition in Case 1.
By \flaring  ~applied to the $\bar{K}$-qi sections ${\gamma}_x$ and $\bar{\gamma}_x$ over 
the interval $\llbracket t_x, v_1\rrbracket$, in view of the bound in Claim \ref{claim:Case1}, since
$$
M_{\bar{K}}\le C_{\ref{claim:Case1}},
$$
we obtain: 
$$
d_T(t_x,v_1)\leq R_{\ref{cor:super-weak flaring}}(\bar{K},C_{\ref{claim:Case1}})$$
and 
$$
d_T(t_y,v_1)\leq R_{\ref{cor:super-weak flaring}}(\bar{K},C_{\ref{claim:Case1}}).$$

 By the triangle inequality, $d_X(\rho(x),\rho(y))$ is bounded by 
\begin{align*}
d_{X_{t_x}}(\rho(x), \gamma_x(t_x))+ d_X(\gamma_x(t_x), \gamma_x(v_1))+
d_{X_{v_1}}(\gamma_x(v_1), \gamma_y(v_1))+\\
d_X(\gamma_y(v_1), \gamma_y(t_y))+ 
 d_{X_{t_y}}(\gamma_y(t_y), \rho(y)).
\end{align*}

This, in turn, is
bounded by 
\begin{align*}
M_{\bar{K}}+K d_T(t_x,v_1)+Kd_T(t_y,v_1)+d_{X_{v_1}}(\gamma_x(v_1), \gamma_y(v_1)) + K d_T(t_x,v_1)
+  M_{\bar{K}}\le \\
2(M_{\bar{K}} +  R_{\ref{cor:super-weak flaring}}(\bar{K},C_{\ref{claim:Case1}})) + C_{\ref{claim:Case1}}.
\end{align*}
This concludes the proof of the proposition in Case 1.

{\bf Case 2:} Assume that $v\not\in \llbracket u,w\rrbracket$. Let $z$ denote 
the nearest point projection of $b$ to $\llbracket u,w\rrbracket $ in $T$. 
There are three subcases to be dealt with.

\medskip 
{\bf Subcase (i):} Suppose $t_x, t_y\in \llbracket b,z\rrbracket $. 
Letting $v_1=t_x, v_2=t_y$ and applying
Lemma \ref{case3-bundle-proj}(ii) with $r=r_1=M_{\bar{K}}$ 
we obtain 
$$
d(\gamma_x(t_x),\gamma_y(t_y))\leq R_{\ref{case3-bundle-proj}}(K,M_{\bar{K}},D'),$$
 and, hence,
$$
d(\rho(x),\rho(y))\leq 2M_{\bar{K}}+R_{\ref{case3-bundle-proj}}(K,M_{\bar{K}},D'). 
$$

{Observe that we need  uniform $k$-flaring for 
$$
k= k_{\ref{case3-bundle-proj}}= K'+ R_{\ref{cor:super-weak flaring}}(\bar K, r_1), r_1=M_{\bar{K}}. 
$$
}

\medskip 
{\bf Subcase (ii):} Suppose $t_x\in \llbracket u,z\rrbracket $ and $t_y\in \llbracket b,z\rrbracket $, or vice versa.

Without loss of generality, we may assume that $t_x\in \llbracket u,z\rrbracket $ and 
$t_y\in \llbracket b,z\rrbracket $. We first get a uniform upper bound on 
$d_{X_z}(\bar{\gamma}_x(z), \bar{\bar{\gamma}}_x(z))$.

 By Lemma 
 \ref{lem:projection-1},  since $Y_z$ is $\la$-quasiconvex in $X_z$, $\bar{\gamma}_x(z)$ lies  in the $(\la+2\delta_0)$-neighborhood of 
$[\gamma_x(z) \bar{\gamma}_y(z)]_{X_z}$ in $X_z$.

\begin{figure}[tbh]
\centering
\includegraphics[width=60mm]{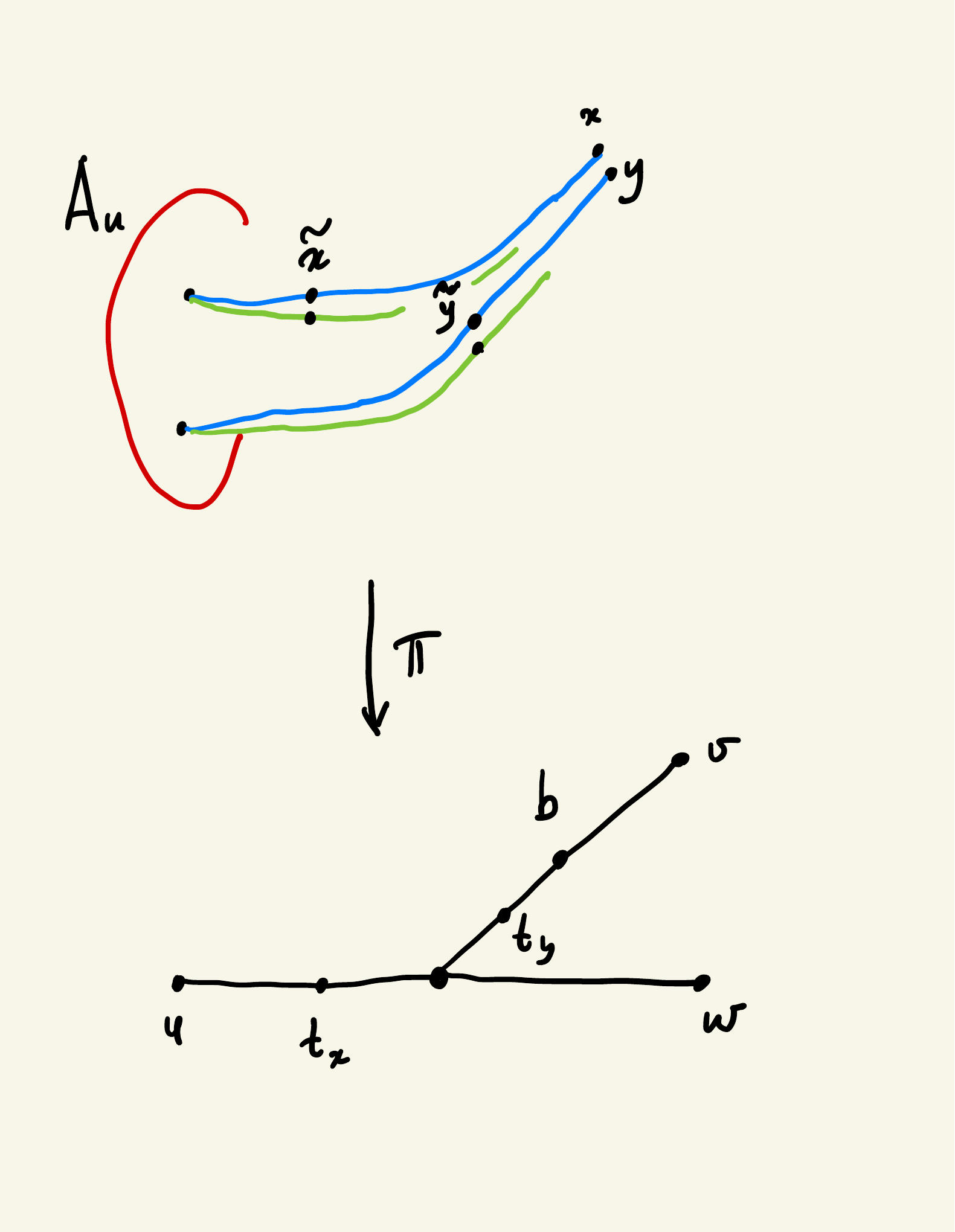}
\caption{Case 3(ii)}
\end{figure}

Since 
$$d_{X_z}(\gamma_y(z), \bar{\gamma}_y(z))\leq R_{\ref{cor:super-weak flaring}}(\bar{K}, M_{\bar{K}}),$$
 by Remark \ref{rem:bundle-proj}$(3)$, it follows from the $\delta_0$-hyperbolicity of $X_z$ 
that the Hausdorff distance between the geodesics $[\gamma_x(z) \gamma_y(z)]_{X_z}$ and $[\gamma_x(z) \bar{\gamma}_y(z)]_{X_z}$ in $X_z$ 
is at most $\delta_0+R_{\ref{cor:super-weak flaring}}(\bar{K}, M_{\bar{K}})$. 

Combining this with the earlier observation that 
$$
\bar{\gamma}_x(z)\in N^{fib}_{\la+2\delta_0}(
[\gamma_x(z) \bar{\gamma}_y(z)]_{X_z}),$$ 
we conclude that  $\bar{\gamma}_x(z)$ belongs to the 
$(3\delta_0+ \la + R_{\ref{cor:super-weak flaring}}(\bar{K}, M_{\bar{K}}))$-neighborhood of
$[\gamma_x(z) \gamma_y(z)]_{X_z}$ in $X_z$. 
Therefore,
$$
d_{X_z}(\bar{\gamma}_x(z), \bar{\bar{\gamma}}_x(z))\leq 
3\delta_0+ \la + R_{\ref{cor:super-weak flaring}}(\bar{K}, M_{\bar{K}}).
$$

\medskip

Since, by Remark \ref{rem:bundle-proj}(6), 
$$
d_{X_z}(\gamma_x(z), \bar{\bar{\gamma}}_x(z))\leq R_{\ref{cor:super-weak flaring}}(\bar{\bar{K}}, 1),$$ 
by triangle inequality we obtain 
$$
d_{X_z}(\gamma_x(z), \bar{\gamma}_x(z))\leq r_2:= 
R_{\ref{cor:super-weak flaring}}(\bar{\bar{K}}, 1) + 
3\delta_0+ \la + R_{\ref{cor:super-weak flaring}}(\bar{K}, M_{\bar{K}}).$$

We will apply Lemma \ref{case3-bundle-proj}(ii) with $r=r_2$ and $v_1=z$, $v_2=t_y$.

\begin{rem}
Note that in order to apply Lemma \ref{case3-bundle-proj}(ii) we need the uniform $k$-flaring condition for 
$$
k= k_{\ref{case3-bundle-proj}}= K'+ R_{\ref{cor:super-weak flaring}}(\bar K, r_2), 
r_2=3\delta_0+ \la + R_{\ref{cor:super-weak flaring}}(\bar{\bar{K}}, 1)  
 + R_{\ref{cor:super-weak flaring}}(\bar{K}, M_{\bar{K}}), 
$$
which is ensured by the choice of the constant $r$ in \eqref{eq:key-r}, $k$ in  \eqref{eq:key-k} and the uniform 
$\kappa$-flaring assumption in Theorem \ref{thm:flow-to-bundle}. 
\end{rem}

In view of the $\bar{K}$-uniform flaring condition, applied to $\ga_x, \bar\ga_x$, we have $d_T(t_x, z)\leq \tau_{\ref{prop:weak flaring}}(\bar{K}, r_2)$. Hence, 
\begin{equation}\label{eq:gx-fla}
d(\gamma_x(t_x),\gamma_x(z))\leq K\tau_{\ref{prop:weak flaring}}(\bar{K}, r_2).
\end{equation}
 
Since
$$
d_{X_{t_y}}(\ga_y(t_y), \bar\ga_y(t_y))\le M_{\bar{K}}\le r_2
$$
and $$
d_{X_z}(\gamma_x(z), \bar{\gamma}_x(z))\leq r_2,$$
we can apply Lemma \ref{case3-bundle-proj}(ii) with $v_1=z, v_2=t_y$ and obtain:  
$$
d(\gamma_x(z), \gamma_y(t_y))\leq R_{\ref{case3-bundle-proj}}(K,r_2,D').
$$ 
Combining this estimate with the inequality \eqref{eq:gx-fla} we get 
$$
d(\gamma_x(t_x),\gamma_y(t_y)) \leq 
K\tau_{\ref{prop:weak flaring}}(\bar{K}, r_2)+R_{\ref{case3-bundle-proj}}(K,r_2,D')$$ 
This, in turn, implies that 
$$d(\rho(x), \rho(y))\leq K\tau_{\ref{prop:weak flaring}}(\bar{K}, r)+R_{\ref{case3-bundle-proj}}(K,r)+2M_{\bar{K}}.$$

\begin{figure}[tbh]
\centering
\includegraphics[width=60mm]{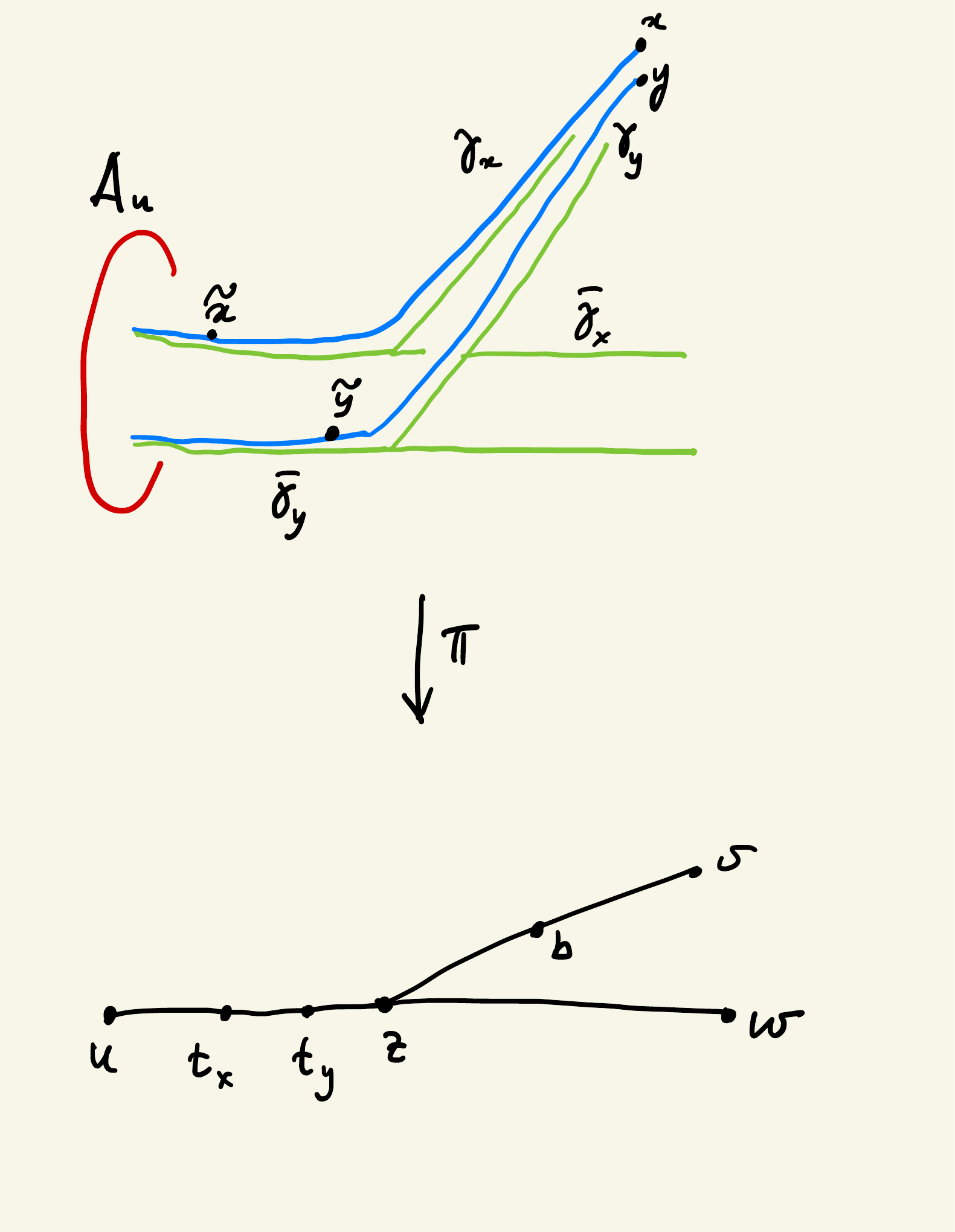}
\caption{Case 3(iii)}
\label{Case 3(iii)}
\end{figure}

\medskip 
{\bf Subcase (iii):} Suppose $t_x, t_y\in \llbracket u,z\rrbracket $.  Without loss of generality, after swapping $x$ and $y$ if necessary, 
we may assume that $t_y\in \llbracket t_x, z\rrbracket $, see Figure \ref{Case 3(iii)}. 
We will  show that $d_{X_z}(\gamma_x(z), \gamma_y(z))$ is uniformly bounded in this subcase.

Suppose first that 
$$
d_{X_z}(\bar{\bar{\gamma}}_x(z), \bar{\bar{\gamma}}_y(z)) \geq 9\delta_0=D_{\ref{proj-geod}}(\delta_0, \delta_0).
$$ 
Then (by Lemma \ref{proj-geod}) both points 
$\bar{\bar{\gamma}}_x(z), \bar{\bar{\gamma}}_y(z)$ belong to the 
$6\delta_0=R_{\ref{proj-geod}}(\delta_0, \delta_0)$-neighbor\-hood of 
$[\bar{\gamma}_x(z) \bar{\gamma}_y(z)]_{X_z}$ in $X_z$.

However, by Remark \ref{rem:bundle-proj}(6), 
$$d_{X_z}(\bar{\bar{\gamma}}_x(z),\gamma_x(z)) \leq R_{\ref{cor:super-weak flaring}}(\bar{\bar{K}}, 1) 
\hbox{~~and~~} d_{X_z}(\bar{\bar{\gamma}}_y(z),\gamma_y(z)) \leq R_{\ref{cor:super-weak flaring}}(\bar{\bar{K}}, 1).$$
Thus, since $Y_z$ is $\la$-quasiconvex in $X_z$ and both endpoints 
of $[\bar{\gamma}_x(z) \bar{\gamma}_y(z)]_{X_z}$ are in $Y_z$, the above inequalities imply: 
$$
d_{X_z}(\ga_x(z), Y_z)\le  r_3:=\la+ R_{\ref{cor:super-weak flaring}}(\bar{\bar{K}}, 1)+ 6\delta_0, 
$$
and 
$$
d_{X_z}(\ga_y(z), Y_z)\le  r_3.  
$$

Since $\bar\ga_x(z)$ is a nearest-point projection (in $X_z$) of $\ga_x(z)$ to $Y_z$, it follows that 
$$
d_{X_z}(\gamma_x(z), \bar{\gamma}_x(z))\leq r_3,$$
and, similarly, 
$$
d_{X_z}(\gamma_y(z), \bar{\gamma}_y(z))\leq r_3.$$ 
We are now again in position to apply Lemma \ref{case3-bundle-proj}(iii) with $r=r_3, v_1=v_2=z$ and conclude that
$$
 d_{X_z}(\gamma_x(z), \gamma_y(z))\leq 
 R_{\ref{cor:super-weak flaring}}(k, \max(1, D')), 
$$
where $k= K'+ R_{\ref{cor:super-weak flaring}}(\bar K, r_3)$. 

\begin{rem}
Observe that in order to apply Lemma \ref{case3-bundle-proj}(iii), we need the uniform $k$-flaring condition for this $k$ and 
$$r_3:=\la+ R_{\ref{cor:super-weak flaring}}(\bar{\bar{K}}, 1)+ 6\delta_0,$$ which is ensured by 
the choice of the parameter $r$ in \eqref{eq:key-r} and the uniform $\kappa$-flaring assumption in  Theorem 
\ref{thm:flow-to-bundle}. 
\end{rem}

The uniform $\bar{K}$-flaring condition in $\X$ applied to the pairs of $\bar{K}$-qi sections 
$(\ga_x, \bar\ga_x)$, $(\ga_y, \bar\ga_y)$, then implies the inequality   
$$
\max( d_T(t_x, z), d_T(t_y, z)) \leq \tau_{\ref{prop:weak flaring}}(\bar{K}, r_3). 
$$ 
Hence,
$$
\max( d(\rho(x), \ga_x(z)),   d(\rho(y), \ga_y(z)) ) 
\le \bar{K} \tau_{\ref{prop:weak flaring}}(\bar{K}, r_3) + M_{\bar{K}}. 
$$
By the triangle inequality we get
$$
d(\rho(x), \rho(y))\le 2(\bar{K} \tau_{\ref{prop:weak flaring}}(\bar{K}, r_3) + M_{\bar{K}}) + 
 R_{\ref{cor:super-weak flaring}}(k, \max(1, D')). 
$$

\noindent This concludes the argument in subcase (iii) and, hence,  the proof of Proposition \ref{lem:proj-to-bundle}. \qedhere

\medskip 
The following is an immediate corollary of the proposition:

\begin{cor}\label{cor3.32} 
For each vertex $v\in \llbracket u,w\rrbracket$  the restriction of the retraction $\rho$ to the subspace $Z_v$ is 
$C_{\ref{lem:proj-to-bundle}}(\la, K,K',D')$-coarse Lipschitz. This 
bound is independent of the choice of the paths $\ga_x$, $\ga_y$ as above.  
\end{cor}

\begin{lemma}\label{lem:proj-to-bundle1}
If $\pi(x)=v_1, \pi(y)=v_2$, $d_T(v_1, v_2)=1$ and $d_{X_{v_1v_2}}(x,y)\le K$,  
 then $d(\rho(x), \rho(y))\leq C_{\ref{lem:proj-to-bundle1}}(\la,K, K',D')$.
\end{lemma}

\proof Without loss of generality, we may assume that $d(u,v_2)=d(u,v_1)+1$. 
Let $y_1=\ga_y(v_1)$. Since $d(x, y_1)\leq 2K$,  we also have 
$$
d_{X_{v_1}}(x,y_1)\le \eta_0(2K). 
$$
Applying Corollary \ref{cor3.32} to the points $x, y_1\in X_{v_1}$, we obtain:
$$
d(\rho(x), \rho(y_1))\le C_{\ref{lem:proj-to-bundle}}(\la,K,K',D') (\eta_0(2K) +1). 
$$

We next note that, without loss of generality we may assume that $\ga_{y_1}$ is chosen to be the restriction of 
$\ga_y$ to the subinterval $\llbracket u, v_1\rrbracket$ since the Lipschitz bound in 
Corollary \ref{cor3.32} holds regardless of the choice of the sections $\ga_y$. 
Hence, 
$$
d(\rho(x), \rho(y))\le C_{\ref{lem:proj-to-bundle}}(\la,K,K',D') (\eta_0(2K) +1).  \qedhere
$$

\medskip 
This completes the proof of Theorem \ref{thm:flow-to-bundle}.  \qed 

\medskip 
The next corollary is immediate from Theorem \ref{thm:flow-to-bundle}: 

\begin{cor}
The map $\rho$ in Theorem \ref{thm:flow-to-bundle} is ``coarsely independent'' of the choice of paths $\ga_x$ used in its construction. More precisely, if $\rho, \rho'$ are two projections defined using different choices of paths $\ga_x$, then 
$$
d(\rho(x), \rho'(x))\le L_{\ref{thm:flow-to-bundle}}(\la,K,K',D').
$$
\end{cor}

Recall that for $x\in \ZZ$ we defined a path $c_x$ in $Z$ connecting $x$ to $\rho(x)$, see Definition 
\ref{defn:retraction-to-bundle}. 

\begin{cor}\label{cor:paths-to-bundle}
Under the assumptions of Theorem \ref{thm:flow-to-bundle}, 
for any two points $x, y\in Z$ within distance $C$ from each other, the Hausdorff 
distance between the paths $c_x, c_y$ is $\le D_{\ref{cor:paths-to-bundle}}(\la,K, K', D', C)$. 
\end{cor}
\proof As in the proof of Theorem \ref{thm:flow-to-bundle}, it suffices to verify the claim in two cases: 

{\bf Case 1:} Suppose that $x, y\in Z_v$ and $d(x,y)\le C$.  
Without loss of generality we may assume that $t_y\in \llbracket t_x,v\rrbracket $. By Theorem \ref{thm:flow-to-bundle} we have 
$$
d_T(\pi(t_x), \pi(t_y))\le d(\rho(x), \rho(y))\leq L_{\ref{thm:flow-to-bundle}}(\la,K,K',D') \cdot (C+1).$$ 
Hence, the length of the portion of $c_x$ between $\gamma_x(t_y)$ and $\rho(x)$  is at most 
$$
M_{\bar{K}}+ K\cdot L_{\ref{thm:flow-to-bundle}}(\la,K,K',D') \cdot (C+1).$$ 
It follows that 
$$
d_X(\gamma_x(t_y),\gamma_y(t_y))\leq R_1:=L_{\ref{thm:flow-to-bundle}}(\lambda,K,D')+2M_{\bar{K}}+
K\cdot L_{\ref{thm:flow-to-bundle}}(\la,K,K',D') \cdot (C+1).$$
Since $X_{t_y}$ is $\eta_0$-uniformly properly embedded in $X$, we also obtain 
$$
d_{X_{t_y}}(\gamma_x(t_y),\gamma_y(t_y))\leq R_2:=\eta_0(R_1).$$
By Corollary \ref{cor:super-weak flaring}, we obtain that for all vertices $s\in V(\llbracket t_y,v\rrbracket)$
$$
d_{X_s}(\gamma_x(s), \gamma_y(s))\leq R_3:=R_{\ref{cor:super-weak flaring}}(K, \max(R_2, C)). 
$$
By combining this with the earlier estimate on the  length of the portion of $c_x$ between $\gamma_x(t_y)$ and $\rho(x)$, we obtain that the Hausdorff distance between $c_x,c_y$ is at most 
$$
R_4:=M_{\bar{K}}+ K\cdot L_{\ref{thm:flow-to-bundle}}(\la,K,K',D') \cdot (C+1) + M_{\bar{K}} + R_3.$$
This concludes the proof in Case 1.

\begin{rem}\label{rem:diff-paths} 
This argument also proves that if in the definition of $\rho(x)$ and $c_x$ we use different $K$-qi sections $\ga'_x, \ga''_x$, then the resulting paths $c'_y, c''_y$  are within Hausdorff distance $R_4$ from each other. 
\end{rem}

\medskip 
{\bf Case 2:} Suppose that $\pi(x)=v_1, v_2=\pi(y)$, $d_T(v_1, v_2)=1$ and $d(x,y)\le K$. 
Without loss of generality we may assume that $d(u,v_2)=d(u,v_1)+1$. Setting $y_1:= \ga_y(v_1)$, according to Case 1 and the above remark, we obtain the bound $R_4$ on the Hausdorff distance between $c_x$ and $c_{y_1}$. It follows that the Hausdorff distance between $c_x$ and $c_y$ is $\le K+R_4$. \qed

\medskip

\begin{prop}\label{prop:rho-npp}
Under the assumptions of Theorem \ref{thm:flow-to-bundle}, 
assuming, in addition, that $Z$ is $\delta$-hyperbolic,  
there exists $C=C_{\ref{prop:rho-npp}}(\delta,\la,K,K',D')$ such that 
$d(\rho_Y, P)\le C$, where $P=P_{Z,Y}$ is the nearest-point projection in $Z$. 
\end{prop}
\proof By Theorem \ref{thm:flow-to-bundle}, there exists a coarse  $L=L_{\ref{thm:flow-to-bundle}}(\la,K,K',D')$-Lipschitz  retraction $\rho: Z\to Y$. By 
Lemma \ref{lem:qc-criterion}, the subset $Y$ is $\la'=\la_{\ref{lem:qc-criterion}}(L,\delta)$-quasiconvex in $Z$.

By the construction of $\rho$, the path $c_x$ connecting $x\in Q_v$ to $\bar{x}=\rho_Y(x)$ is the concatenation of the $K$-qi section $\ga_{x,\tilde{x}}$ and the vertical geodesic $[\tilde{x} \bar{x}]_{X_t}$ of length $\le M_{\bar{K}}$, $t=t_x$. Hence, considering the geodesic triangle $\Delta x \tilde{x} \bar{x}$ in $Z$ and using the $\delta$-hyperbolicity of $Z$, we conclude that 
$$
\Hd_Z(\ga_{x,\tilde{x}}, [x \bar{x}]_{Z})\le \delta+ M_{\bar{K}} + D_{\ref{stab-qg}}(\delta, K).
$$ 
Therefore, since $Y$ is $\la'$-quasiconvex in $Z$, 
in order to prove the proposition, according to Corollary \ref{cor:projection-1}, it suffices to verify that 
there exists a function $R\mapsto R'$ (depending on $\delta,\la'$) such that 
for all $v\in V(\pi(Z)), x\in Q_v, y\in \YY$ 
$$
d_{Z}(y, \ga_{x,\tilde{x}} \cap \XX)\le R\Rightarrow d_Y(y, \bar{x})\le R'.$$
Our goal then is to define such a function.

Take $y\in \YY$ such that  $d_{Z}(y, \ga_{x,\tilde{x}})\le R$. Pick a point 
$$
x'\in \ga_{x,\tilde{x}}\cap X_{v'}, v'\in \llbracket u, v\rrbracket,$$
such that $d_{Z}(y, x')\le R$. In particular, $d_T(v', \pi(y))\le R$.  

\begin{figure}[tbh]
\centering
\includegraphics[width=60mm]{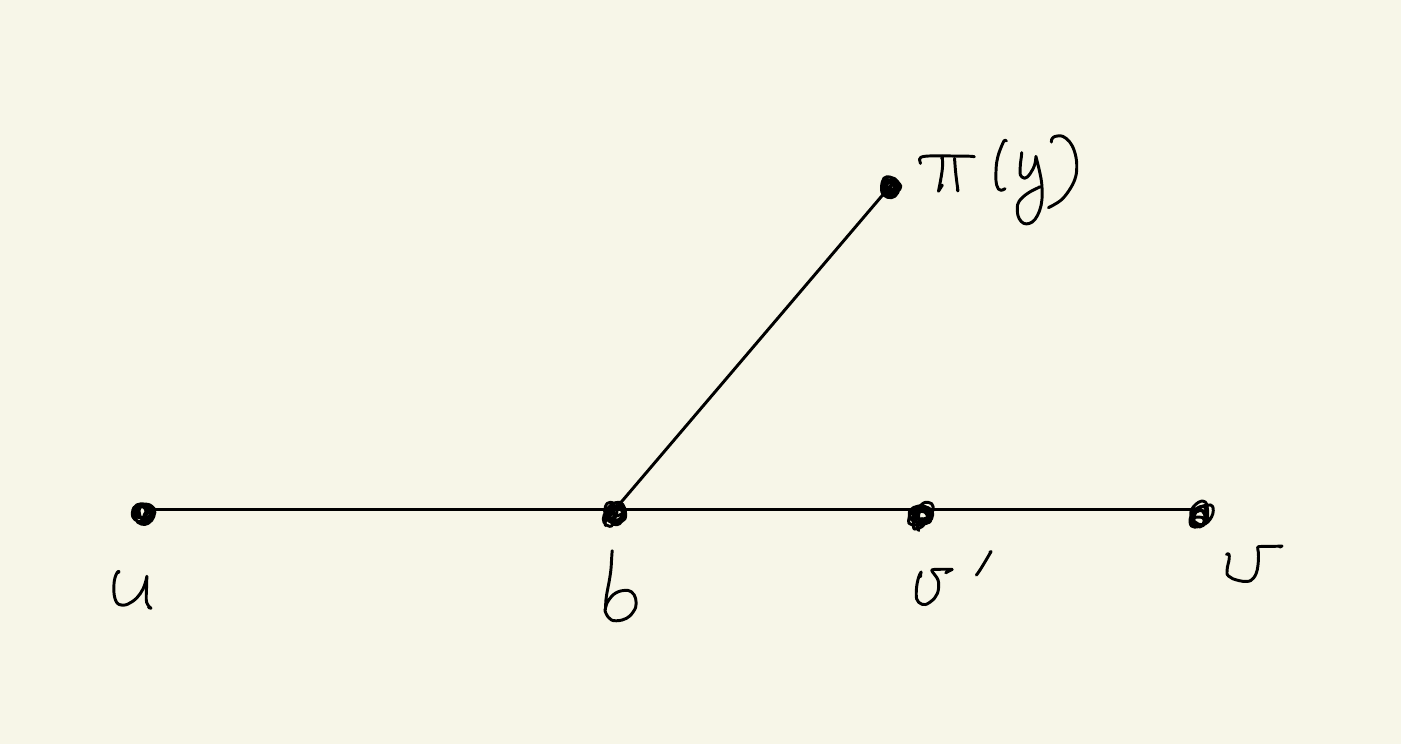}
\label{Trpd.fig}
\end{figure}

Let $b$ denote the center of the tripod $\Delta  v' \pi(y) u \subset T$, 
$$
b\in \llbracket u, v'\rrbracket \subset \llbracket u, v\rrbracket,$$
see Figure \ref{Trpd.fig}. Then 
\begin{equation}\label{eq:proj-to-Y} 
d_T(v',b)+ d_T(b,\pi(y))= d_T(v', \pi(y))\le R.  
\end{equation}
Set $x'':= \ga_x(b)$ (this point is defined since $b\in \llbracket u, v\rrbracket$ and $\ga_x$ is a section over 
that interval) and let $\ga''$ denote $\ga_{x'',\tilde{x}}$, the restriction of $\ga_x$ to 
$\llbracket b, t\rrbracket$. (Note that the order in which the vertices $t, b$ appear in the interval 
$\llbracket u, v'\rrbracket$ is unclear.)

Let $\bar{\ga}''$ denote the fiberwise projection of $\ga''$ to 
${\YY}_{t b}$; without loss of generality, $\bar{\ga}''(b)=\bar{x}$.  
Then $\bar\ga$ is a $\bar{K}$-qi leaf in $Y$ 
connecting $\bar{x}$ to $\bar{x}''= P_{X_b, Y_b}(x'')$.

Furthermore, let $\ga_y$ be a $K'$-qi section in $Y$ over $\llbracket \pi(y), u\rrbracket$, 
connecting $y$ to $Y_u$ and let $\ga''_{y}$ be its restriction to the interval $\llbracket \pi(y), b\rrbracket$
(recall that $b\in \llbracket \pi(y), u\rrbracket$). Set $y'':= \ga_y(b)=\ga''_y(b)$.

Now, consider the quadruple of points $x'', x', y, y''$: We have a $K$-qi section $\ga_{x'',x'}$ connecting $x''$ to $x'$, a geodesic $[x' y]_Z$ in $Z$  connecting $x'$ to $y$ and the $K'$-qi section $\ga''_y$ connecting $y$ to $y''$. 
Therefore, 
$$
d_Z(y,y'')\le K' d_T(\pi(y),b)
$$
and 
$$
d_Z(x'', y'')\le K d_T(b,v') + R + K' d_T(\pi(y),b). 
$$
Taking into account the inequality \eqref{eq:proj-to-Y}, we obtain 
\begin{align}\label{eq:x"y'}
\begin{split}
d_Z(x'', y'')\le  R +  \max(K, K') d_T(\pi(y), v')  \le \\
R + R \max(K, K') = R( 1+ \max(K, K')),   
\end{split}
\end{align}
\begin{equation}\label{eq:yy''}
d_Z(y,y'')\le K' d_T(\pi(y), v')  \le K'R. 
\end{equation}

Since $x'', y''\in X_b$, we also get 
\begin{equation}\label{eq:x"y''}
d_{X_b}(x'', y'')\le \eta_0R( 1+ \max(K, K'))). 
\end{equation}
Since $y''\in Y_b$ and $\bar{x}''$ is the projection of $x''$ to $Y_b$, it follows that
\begin{equation}\label{eq:x''barx''}
d_{X_b}(x'', \bar{x}'')\le \eta_0(R( 1+ \max(K, K'))). 
\end{equation}

 The uniform $\bar{K}$-flaring condition applied to the restrictions of $\ga''$ and $\bar\ga$ to
 $\llbracket t, b \rrbracket$,  and the inequalities  
 \begin{align}
 \begin{split}
 d_{X_b}(\ga''(b), \bar\ga''(b)) \le \eta_0(R( 1+ \max(K, K'))),\\
  d_{X_t}(\ga''(t), \bar\ga''(t))= d_{X_t}(\tilde x, \bar{x})\le M_{K'},
 \end{split}
 \end{align}
 imply  that
$$
d_T(t, b)\le \tau_{\ref{prop:weak flaring}}(\bar{K},  \max(\eta_0R( 1+ \max(K, K'))), M_{K'}). 
$$
In particular,
\begin{equation}\label{eq:barxbarx''}
d(\bar{x}, \bar{x}'') \le \bar{K} \tau_{\ref{prop:weak flaring}}(\bar{K},  \max(\eta_0R( 1+ \max(K, K'))), M_{K'}). 
\end{equation}

Putting together the inequalities \eqref{eq:barxbarx''}, \eqref{eq:x''barx''}, \eqref{eq:x"y''}, \eqref{eq:yy''}, by the triangle inequality, we get:
\begin{align*}
d_Z(\bar{x}, y)\le d_Z(\bar x, \bar{x}'') + d_Z(\bar{x}'', {x}'') + d_Z({x}'', y'') + d_Z({y}'', y) \le \\
R':= \\
\bar{K} \tau_{\ref{prop:weak flaring}}(\bar{K},  \max(\eta_0R( 1+ \max(K, K'))), M_{K'}) +\\
\eta_0R( 1+ \max(K, K'))+\\
R( 1+ \max(K, K'))+ K'R.  
\end{align*}
This proves the proposition. \qed

\medskip 
Combining Theorem \ref{thm:flow-to-bundle} with the existence of uniformly coarse Lipschitz retractions  
$X\to Z$, where $Z=Fl_K(Q_u)$ or $Z=L_K(\al)$, (Theorem \ref{mjproj}, Corollary \ref{cor:ladder-retraction}), 
we obtain:

\begin{cor}\label{cor:X-to-bundle} 
Suppose that $Y\subset N^{fib}_{4\delta_0}(Fl_K(Q_u))\subset X$ or 
$A_{K',D'}(\al')\subset L=L_{K,D,E}(\al)$ 
satisfy the assumptions of  Theorem \ref{thm:flow-to-bundle} with parameters $\lambda, K, D$, $E, K'$ and $D'$. 
Then there exists  a coarse  $L_{\ref{cor:X-to-bundle}}(\lambda,K,K',D,D',E)$-Lipschitz retraction 
$\rho_Y: X\to Y$.  
\end{cor}

\begin{rem}\label{carpet-rem}
1. It follows that if $X$ were hyperbolic, then the total space $A$ of each $(K',D')$-carpet $\A_{K',D'}(\al')$  
would be uniformly quasiconvex  in $X$. 

2. Unlike the existence theorem for coarse Lipschitz retractions to semicontinuous subtrees of spaces, in order to get a uniform coarse Lipschitz retractions to carpets (and  to bundles) we do not need lower bounds on $K$ besides the obvious inequality $K\ge 1$.  
\end{rem}

\chapter{Hyperbolicity of  ladders}\label{sect:3}

In this chapter we prove uniform hyperbolicity of $(K,D,E)$-ladders in $X$. The proof is divided in three main steps. We first prove (section \ref{sec:Hyperbolicity of carpets}) hyperbolicity of carpets by exhibiting slim combings of carpets (combings satisfying the conditions of Corollary \ref{cor:bowditch}). 
We use these paths, in conjunction with the retractions to carpets (see Corollary \ref{cor:X-to-bundle}) 
to construct combings of {\em carpeted ladders}, which are ladders $\L(\al)$ containing carpets $\A(\al')$ with $\al'\subset \al$ whose length almost equals that of $\al$. 
This is done in Section \ref{sec:narrow ladders}. Lastly, in Section \ref{sec:hyperbolicity-of-ladders} we prove hyperbolicity of general ladders subdividing these 
(a ``vertical subdivision'') into carpeted ladders and then using quasiconvex amalgamation to prove hyperbolicity. 
Hyperbolicity of ladders is a key step for proving hyperbolicity of flow-spaces, which is done in the next chapter.

\section{Hyperbolicity of carpets} \label{sec:Hyperbolicity of carpets}

The proof of the following simple proposition will serve as a model for more complex proofs of 
hyperbolicity of certain subspaces of $X$. 

\begin{prop}\label{prop:easy-one}
For every $K\ge 1$, every $(K,C)$-narrow carpet ${\mathfrak A}=(\pi: A\to \llbracket u,w\rrbracket)$ 
in ${\mathfrak X}$, equipped with its intrinsic metric, is 
$\delta_{\ref{prop:easy-one}}(K,C)$-hyperbolic, provided that ${\mathfrak X}$ 
satisfies the uniform {$K$-flaring condition}. 
\end{prop}
\proof   
Let $\beta=\AA\cap X_w$ denote the $C$-narrow end of ${\mathfrak A}$. 
For each $x\in \AA$ we have the $K$-qi section $\ga_x\subset \Sigma_x$ of 
$\pi: {A}\to  \llbracket u, w\rrbracket$ over $\llbracket \pi(x),w\rrbracket$, connecting $x$ to $\beta$. 
For any two such sections $\ga_x, \ga_y$ we let $t_{xy}\in V( \llbracket w, u\rrbracket)$ denote the supremum of
$$
\{ t\in V(\llbracket w,u\rrbracket): d_{X_t}(\ga_x(t), \ga_y(t))\le M_{K}\}. 
$$
In particular, if this subset is empty, then $t_{x,y}=w$. Set $t:= t_{xy}$. We then define a path $c(x,y)$ as the concatenation of the section $\ga_x$ restricted to $\llbracket \pi(x), t\rrbracket$ with the vertical segment 
$[\ga_{x}(t) \ga_y(t)]_{X_t}$, followed by the concatenation with the restriction of the section $\ga_y$ to the subinterval $\llbracket t, \pi(y)\rrbracket$. The assumption that $\AA$ is $(K,C)$-narrow implies that the length of 
$[\ga_{x}(t) \ga_y(t)]_{X_t}$ is at most 
$$
\max\left(C, M_{K}\right). 
$$
We claim that this family of paths in $A$ 
satisfies the conditions of Corollary \ref{cor:bowditch} with constants depending only on $K$ and $C$.  The assumption that the paths $c(x,y)$ are uniformly coarse Lipschitz 
 follows from the fact that each path $c(x,y)$ is a concatenation of $K$-qi sections and of vertical geodesics. 
 
 \begin{lemma} \label{lem:easy-up} 
 The family of paths $c(x,y)$ is uniformly proper in $A$, with distortion function depending only on $K$ and $C$.  
 \end{lemma}
 \proof Let $x, y\in \AA$ be such that $d(x,y)\le r$.  Set $v_1:=\pi(x), v_2:=\pi(y)$. Then 
 $d_T(v_1, v_2)\le r$ as well. Without loss of generality, on the oriented interval $\llbracket w,u\rrbracket $ we have
 $$
 w\le t_{xy}\le v_1\le v_2 \le u. 
 $$
 Let $y_1\in A_{v_1}$ be the intersection point with $\ga_y$. Then the subpath $\ga_{yy_1}$ in $\ga_y$ between $y, y_1$ has length $\le K r$ and $d(y_1, x)\le r_1=r(K+1)$. 
 Furthermore, $d_{X_{v_1}}(y_1, x)\le \eta_0(r_1)$. 
 By the uniform $K$-flaring condition, 
 $$
 d_T(v_1, t_{xy})\le \tau_{\ref{prop:weak flaring}}(K, \max(M_{K},C, \eta_0(r_1))). 
 $$
 Therefore, the overall length of $c(x,y)$ is
 $$
 \le \max(M_{K}, C) + 2K\tau_{\ref{prop:weak flaring}}(K, \max(M_{K},C, \eta_0(r_1))) + K r. \qed 
 $$

 \begin{rem}\label{rem:dependence on C}
 This lemma is the only place where the constant $C$ plays any role in the proof of the proposition. 
 \end{rem}

 We next verify the  condition (a2) of   Corollary \ref{cor:bowditch} for the family of paths $c(x,y)$. 
 First of all, in the special case of (a2) when $d(y, z)$ is uniformly bounded,   the paths $c(x,y), c(x,z)$ are uniformly close (fellow-travel), as follows from \flaring. 
 (The bounds  depend only on $K$ and on $d(y,z)$.)

Consider now three points $x, y, z\in \AA$ and the triangle $\Delta$ formed by the paths $c(x,y)$, $c(y,z)$, $c(z,x)$ connecting them. 
After relabelling the points $x, y, z$ we can assume that the vertices $t_{xy}$, $t_{yz}$, $t_{zx}$ appear in the interval $\llbracket u,w\rrbracket$ in the following order:
$$
u\le t_{xy}\le t_{yz}\le t_{zx}\le w. 
$$

{\bf Case 1.} Suppose first that $t:= t_{xy}< w$. Then
$$
d_{X_{t}}(\ga_x(t), \ga_y(t))\le M_{K}. 
$$
Therefore, we replace $x$, $y$ with $x':= \ga_x(t), y':= \ga_y(t)$ respectively; $d_{X_t}(x',y')\le M_{K}$, i.e. the length of the path $c(x',y')$ is $\le M_{K}$. 
Thus, $\delta_{\ref{prop:easy-one}}(K)$-slimness of the triangle $\Delta'= \Delta x' y' z$ formed by the paths $c(x', y'), c(y',z), c(z, x')$ follows from the uniform fellow-traveling property of  the paths $c$ in $\A$. Since, without loss of generality we may assume that 
$c(y',z), c(z,x')$ are subpaths in $c(y,z), c(z,x)$, we conclude the $\delta_{\ref{prop:easy-one}}(K)$-slimness of the original triangle $\Delta$. 

\medskip 
{\bf Case 2.} It remains to consider the case $t_{xy}= t_{yz}= t_{zx}= w$. Then, as in Case 1, we replace the points $x, y, z$ with the points $x':= \ga_x(w), y':= \ga_y(w), z'=\ga_z(w)$. The triangle $\Delta'$ formed by the paths $c(x',y'), c(y',z'), c(z',x')$ is contained in the geodesic segment $\beta$; hence $\Delta'$ is $\delta_0$-slim. 
We conclude that $\Delta$ is $\delta_{\ref{prop:easy-one}}(K)$-slim as well. \qed

\begin{rem}
Lemma \ref{lem:easy-up} establishes that the paths $c(x,y)$ are uniformly proper (with distortion depending only on $K$ and $C$). Hence, by Corollary \ref{cor:bowditch}(b), the paths $c(x,y)$ are uniformly 
(in terms of $K$ and $C$) close to geodesics in $A$. 
\end{rem}

\section{Hyperbolicity of carpeted ladders}\label{sec:narrow ladders}

\begin{defn}\label{defn:carpeted ladder} 
Let $\bar{K}$ be defined as in Notation \ref{barK and barbarK}, with $K'=K,\la=\delta_0$, 
\begin{equation}\label{eq:barK.4}
\bar{K}:=K_{\ref{bundle-proj}}(\delta_0, K, K). 
\end{equation}
Set 
\begin{equation}\label{eq:kappa.4}
\kappa:= \kappa_{\ref{defn:carpeted ladder}}(K)=\kappa_{\ref{eq:key-kappa}}(\delta_0,K, K).
\end{equation}
A $(K,D,E)$-ladder $\L(\al)$ containing a $(K,C)$-carpet $\A=\A_{K,C}(\al')$, 
$\al'\subset \al$, as a subladder\footnote{see Definition \ref{defn:subladder}},   
satisfying
$$
\length(\al')\ge \length(\al)- M_{\bar{K}},
$$
 will be called {\em carpeted}.
 \end{defn}

In this section, we will prove that carpeted ladders are uniformly hyperbolic. Only the parameters $K$ and $C$ will play a role in the proof, the parameters $D$ and $E$ will be irrelevant, just as in the proof of the existence of a coarse Lipschitz retraction $L_K(\al)\to A_{K,C}(\al')$.

\begin{thm}
[Hyperbolicity of carpeted ladders] \label{small-ladder} 
Carpeted ladders in $\X$ are  hyperbolic. 
More precisely: Fix $K\ge 1$ and 
suppose that that $\X$ satisfies the uniform $\kappa$-flaring condition, 
where $\kappa=\kappa(K)$ is as in Definition \ref{defn:carpeted ladder}, \eqref{eq:kappa.4}. 
Let  $\L=\L_K(\al)$ be a $K$-ladder containing a $(K,C)$-carpet $\A(\al')$ as in Definition 
\ref{defn:carpeted ladder}. Then $L_K(\al)$ (with its intrinsic metric) is 
$\delta_{\ref{small-ladder}}(K,C)$-hyperbolic. 
\end{thm}
\proof  The proof of this theorem is divided in two steps: We first define a family of paths $c(x,y)$ in $L_K(\al)$ 
connecting points $x, y\in \LL$, {\em using the family of uniform quasigeodesics for the carpet $\A$ as a blackbox},  and then  check that these paths satisfy the conditions of Corollary \ref{cor:bowditch}. We let $\llbracket u, w\rrbracket\subset S\subset T$ denote the base-interval of the carpet $\A$, where $S=\pi(L_K(\al))$ is the base of $\L= (\pi: L_K(\al)\to S)$. Then $u$ is the center of both $\A$ and $\L$, 
$\al'=A_u$ and $\beta= A_w$ is the narrow end of $\A$. 

\medskip
For each $x\in \LL$ we let $b_x\in V(T)$ 
denote the center of the triangle $\Delta u w \pi(x)\subset T$.

\medskip 
According to Theorem \ref{thm:flow-to-bundle},  there exists  
$\rho=\rho_{\L,\A}: \LL\to \AA$, a coarsely $k$-Lipschitz retraction 
with 
$$k=L_{\ref{thm:flow-to-bundle}}(\delta_0,K,K,C).$$

 We first review some facts about  $\rho(x)$.  The point $\rho(x)=\bar{x}$ belongs 
 to the interval $A_s\subset \AA$, for a certain vertex 
 $$
 s\in \llbracket u, \pi(x)\rrbracket \cap \llbracket u, w\rrbracket, 
 $$ 
 such that 
 $\bar{x}$  is within vertical distance $M_{{\bar K}}$ (in the vertex space $X_s$)   
from a point $\tilde{x}= \ga_x(s)\in X_s$. Here  
$\ga_x\subset \Sigma_x$ is the $K$-section in $\LL$ over $\llbracket u, \pi(x)\rrbracket$ connecting 
$x$ to $\al$. (Such point $\bar{x}$ always exists in view of the assumption on the length of $\al'$.)

The paths $\ga'_x=\ga_{x\tilde{x}}\star [\tilde{x} \bar{x}]_{X_s}$ defined in the proof 
of  Theorem \ref{thm:flow-to-bundle},  connecting $x$ to $\bar{x}$,  
satisfy the Hausdorff  fellow-traveling property with respect to variations of $x$ (see Corollary 
\ref{cor:paths-to-bundle}). Here  $\ga_{x\tilde{x}}$ is a subpath of $\ga_x$ connecting $x$ to $\tilde{x}$.

\begin{figure}[tbh]
\centering
\includegraphics[width=60mm]{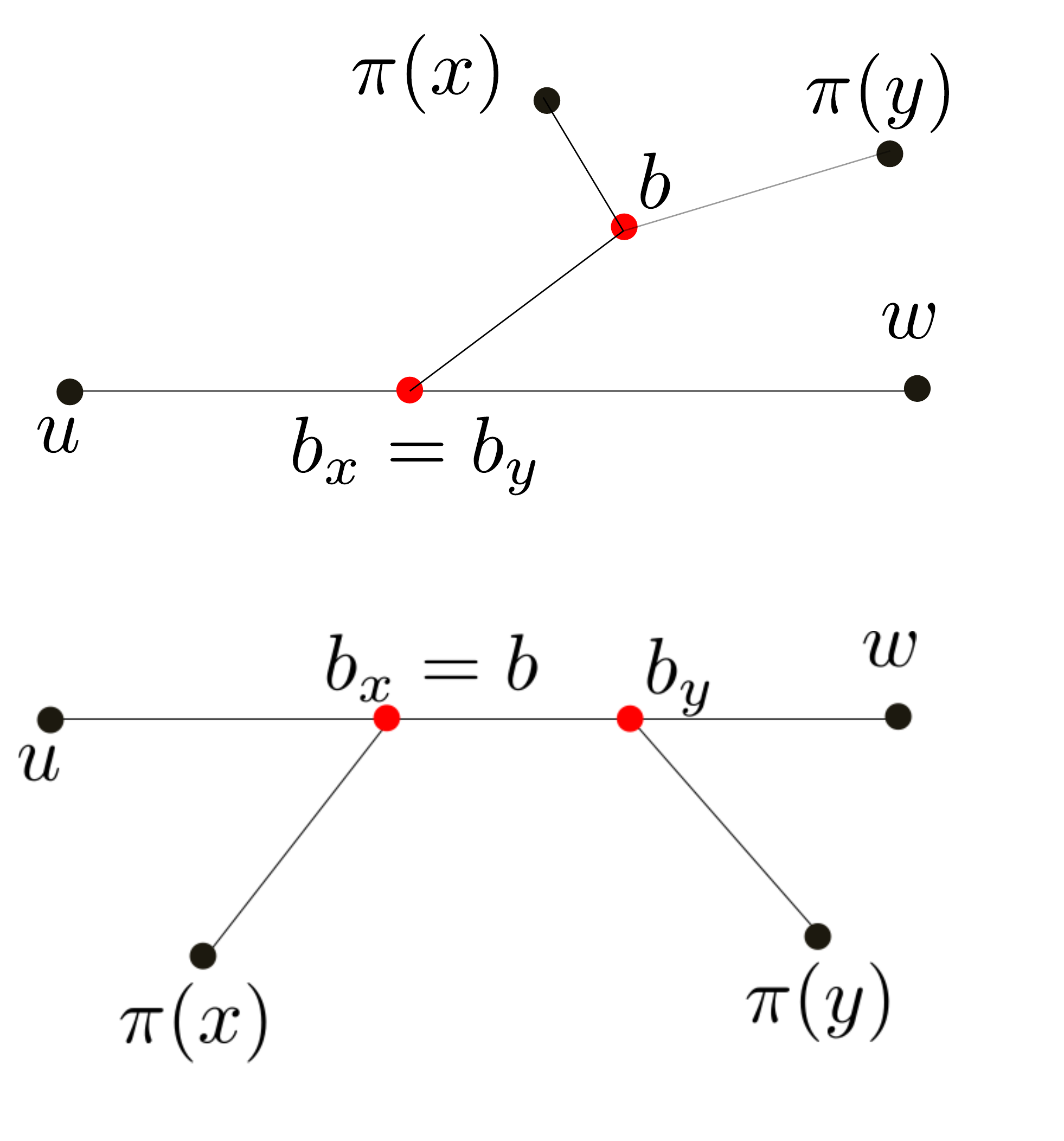}
\caption{Up to relabeling $x, y$ there are two possible configurations of the points $u, w, \pi(x), \pi(y)$.}
\end{figure}

\medskip
{\bf Step 1:} {\bf  Definition of the paths $c_{x,y}=c(x,y)$.}

\smallskip

For $x,y \in \LL$ we let $b=b_{xy}$ be the center of the triangle $\triangle u\pi(x)\pi(y)\subset T$.
There are two cases to consider, depending on which we get two types of paths $c(x,y)$.

\medskip

{\bf Paths of type 1:} There exists $t\in V(\llbracket \pi(x), \pi(\bar{x})\rrbracket \cap \llbracket \pi(y), \pi(\bar{y})\rrbracket)\subset V(\llbracket u,b\rrbracket)$ such that 
$$
d_{X_t}(\gamma_x(t), \gamma_y(t))\leq M_{\bar{K}},$$
i.e. the paths $\ga_x, \ga_y$ ``come sufficiently close'' in some common vertex-space.

\begin{figure}[tbh]
\centering
\includegraphics[width=60mm]{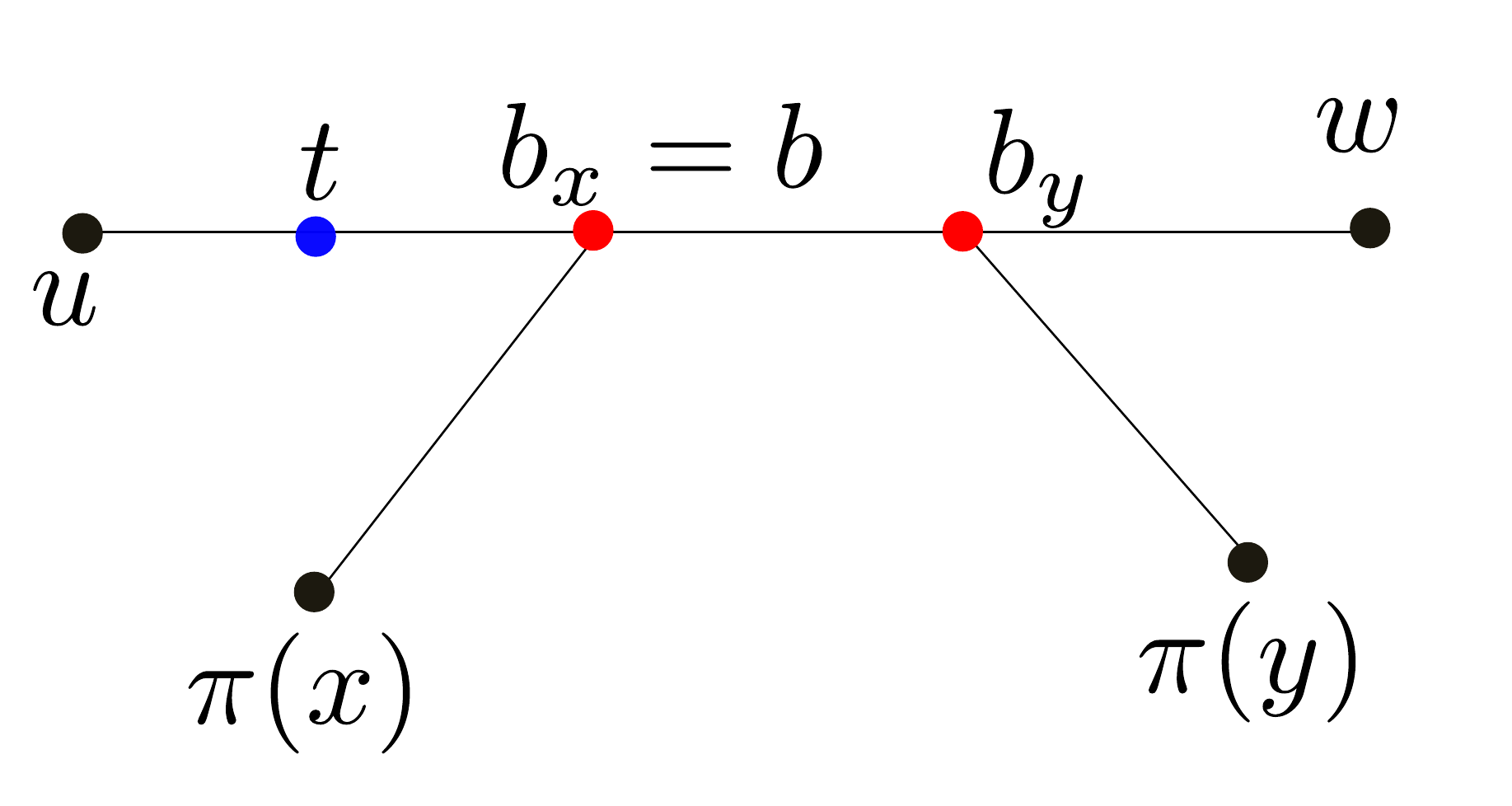} 
\caption{Paths of type 1: $t=t_{x,y}$.}
\end{figure}

Let $t_{xy}$ be the maximal vertex in $\llbracket u,b\rrbracket$ with this property. Then define $c(x,y)$ to be the concatenation of the portions of $\gamma_x$ and (the reverse of) $\gamma_y$ over $\llbracket t_{xy}, \pi(x)\rrbracket$ and $\llbracket t_{x,y}, \pi(y)\rrbracket$ respectively with the subsegment of $L_{t_{xy}}$ joining their end points.

\medskip 
{\bf Paths of type 2:} Suppose $t$ as in type 1 does not exist, i.e. the paths $\ga_x, \ga_y$ ``stay far apart'' in every vertex-space they both enter.  

Then define $c(x,y)$ to be the concatenation of $\ga'_x$ and the reverse of $\ga'_y$ with 
a geodesic $[\bar{x} \bar{y}]_A$ in $A=A_{(K,C)}(\al')$ connecting $\bar{x}$ to $\bar{y}$. 

\begin{rem}
More precisely, instead of geodesics $[\bar{x} \bar{y}]_A$ one should use uniform quasigeodesic paths in $A$ defined in the previous section. Since the two families are uniformly close to each other, we will work with geodesics for the ease of the notation. We will do the same repeatedly later in the book. 
\end{rem}

\medskip 
We observe that each path of type 1 is a concatenation of (at most) three uniformly 
quasigeodesics paths (the middle one of which has) uniformly bounded length $\le M_{\bar{K}}$,  while each path 
of type 2 is a concatenation of (at most) five uniformly quasigeodesics paths (two of which have uniformly bounded 
lengths $\le M_{\bar{K}}$). 

\medskip 
Our next task is to establish a uniform Hausdorff fellow-traveling property of the family of paths $c(x,y)$ (Lemma \ref{wd00}). Even though is property is not required by Corollary \ref{cor:bowditch}, it will play key role in 
verifying the other conditions of the corollary.

\begin{figure}[tbh]
\centering
\includegraphics[width=60mm]{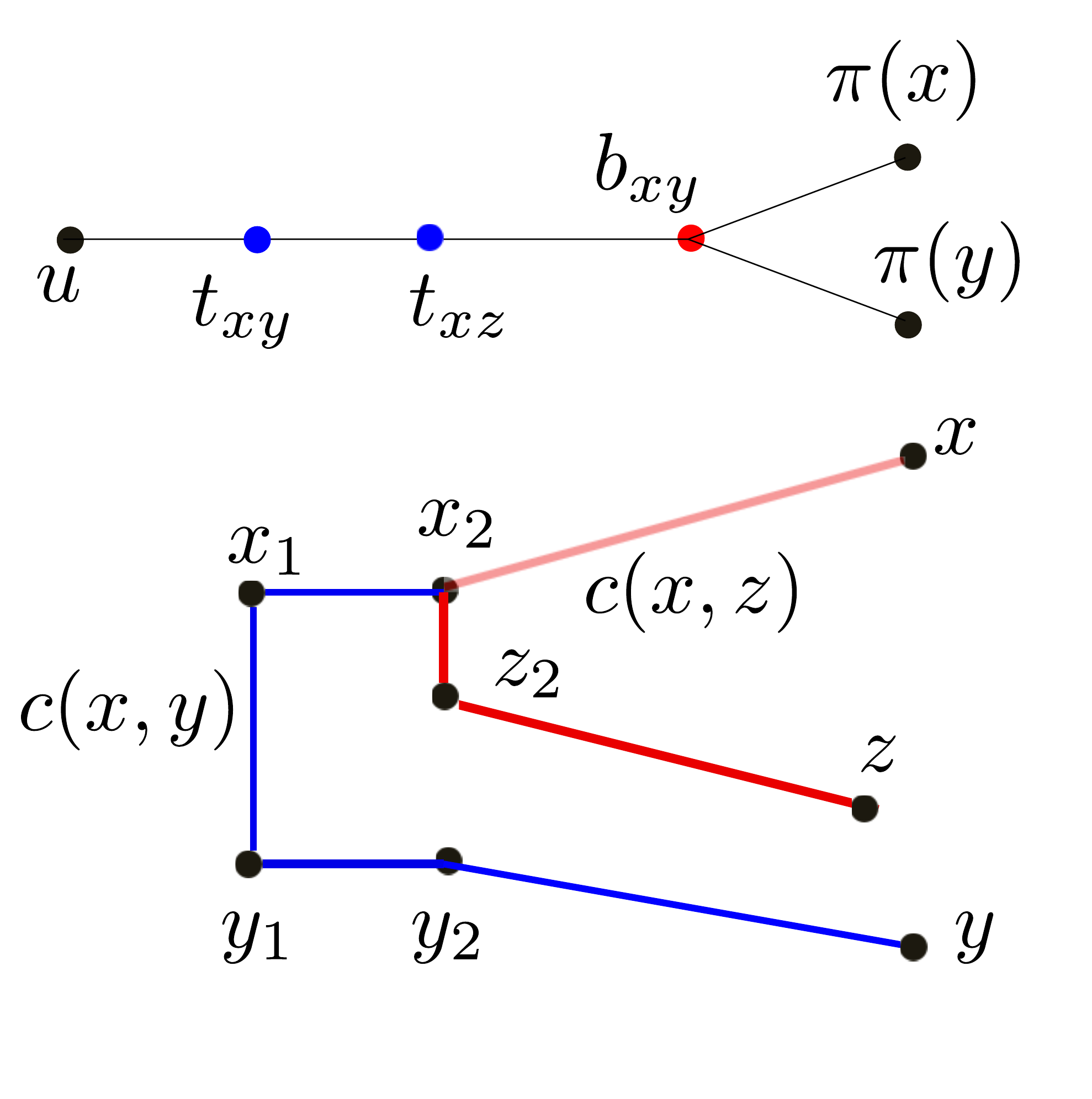} 
\caption{Case 1-1.}
\label{fig:utxy}
\end{figure}

\begin{lemma}\label{wd00}
The paths $c(x,y)$ satisfy the Hausdorff fellow--traveling property,  
i.e. if $y, z$ are uniformly close to each other (in the total space $L$ of the ladder $\L$), so are the images of the 
paths $c(x,y), c(x,z)$. More precisely, 
there is a function $D_{\ref{wd00}}(C,K,r)$ such that 
if $d_L(y, z)\le r$, then 
$$
\Hd(c_{x,y}, c_{x,z}) \le D_{\ref{wd00}} (C, K, r), 
$$ 

\end{lemma}
\proof 
As in the proof of Theorem \ref{thm:flow-to-bundle}, there are two cases to consider: $\pi(y)=\pi(z)$ 
(see part of the proof covered in Proposition \ref{lem:proj-to-bundle}) 
and $d(\pi(y), \pi(z))$ (see part of the proof covered in Lemma \ref{lem:proj-to-bundle1}). 
The second case  follows from the first one just as in the proof of Lemma \ref{lem:proj-to-bundle1}, so we assume that $\pi(y)=\pi(z)$. 
There are three cases to check depending on the types of the paths $c(x,y), c(x,z)$.

\medskip 
{\bf Case 1-1}: Both paths $c(x,y), c(x,z)$ have type 1. The paths $c(x,y), c(x,z)$ agree over the interval 
$\llbracket t, x\rrbracket$ where $t\in \{t_{xy}, t_{xz}\}$ is the vertex closer to $b_{xy}$; after swapping the roles of $y$ and $z$ we may assume that $t=t_{xz}$, see Figure \ref{fig:utxy}.

Define the points
$$
x_1:= \ga_x(t_{xy}), y_1:= \ga_y(t_{xy}),  x_2:= \ga_x(t_{xz}), y_2= \ga_y(t_{xz}), 
z_2:= \ga_z(t_{xz}).
$$
They satisfy the inequalities
$$
d_{X_{t_{xy}}}(x_1, y_1)\le M_{\bar{K}}, \quad d_{X_{t_{xz}}}(x_2, z_2)\le M_{\bar{K}}
$$
Except for the point $y_2$, these are the points where the paths $c(x,y), c(x,z)$ switch from vertical to horizontal. 
The points  $x_1, y_1, y_2$ lie in the image of $c(x,y)$, while the points $x_2, z_2$ lie in the image 
of $c(x,z)$. 

We will show that the length of the interval $\llbracket t_{xy}, t_{xz}\rrbracket$ is uniformly bounded (in terms of $r$ and the other parameters in the theorem). Since
$$
d_{X_{t_{xy}}}(x_1,y_1)\le M_{\bar{K}},
$$
it will follow, according to Lemma \ref{lem:growth-of-flare} that  the part of $c(x,y)$ 
lying between $x_2$ and $y_2$ is uniformly Hausdorff-close to the vertical segment $[x_2 y_2]_{X_{t_{xz}}}$. 
In particular, the length of that segment will be uniformly bounded. 
However, the points $x_2$ and $z_2$ are within vertical distance $M_{\bar{K}}$ from each other. Hence, 
the vertical distance between $z_2$ and $y_2$ is also uniformly bounded. Since $d_L(y,z)\le r$, by \flaring, it will follow that the vertical distance between $\ga_y, \ga_z$ over the interval $\llbracket t_{xz}, \pi(y)\rrbracket$ is also uniformly bounded, thereby, concluding the proof. Thus, it remains to bound the length of 
the interval $\llbracket t_{xy}, t_{xz}\rrbracket$.

\medskip

By the definition of the projection $\rho: \LL\to \AA$, 
$$
\rho(x_1)= \rho(x_2)=\bar{x}, \rho(y_1)=\bar{y}, \rho(z_2)=\bar{z}.
$$
 Since $\rho$ is $k$-coarse Lipschitz, we have 
$$
d_{A}(\bar{x}, \bar{y})\le k (M_{\bar{K}}+1), \quad d_{A}(\bar{x}, \bar{z})\le k (M_{\bar{K}}+1), 
\quad d_{A}(\bar{y}, \bar{z})\le k (r+1). 
$$
Lemma \ref{lem:3-flows} now implies that 
$$
d_T(t_{xy}, t_{xz})\le \tau_{\ref{lem:3-flows}}(K, \max(k (M_{\bar{K}}+1), k (r+1) )).
$$  
This concludes the proof in Case 1-1.

\medskip
{\bf Case 1-2}: The path $c(x,y)$ is of type 1 while $c(x,z)$ is of type 2. 

Since $c(x,y)$ has type 1, we define the vertex 
$$
t=t_{xy}\in \llbracket u, \pi(x)\rrbracket \cap \llbracket u, v\rrbracket, v= \pi(y)=\pi(z). 
$$
At this vertex 
\begin{equation}\label{gaxgayt}
d_{X_t}(\ga_x(t), \ga_y(t))\le M_{\bar{K}}. 
\end{equation}
All three paths $\ga_x, \ga_y, \ga_z$ enter the same vertex space $X_t$, $t=t_{xy}$ at points $x_1, y_1, z_1$ respectively.  

Since $\rho$ is $k$-coarse Lipschitz, and $\bar{x}=\rho(x)=\rho(x_1), \bar{y}=\rho(y)=\rho(y_1)$, we have 
$$
d_{A}(\bar{x}, \bar{y})\le k (M_{\bar{K}}+1).$$
Similarly, since $d_L(y,z)\le r$, we obtain
$$d_{A}(\bar{y}, \bar{z})\le k (r+1).
$$ 
Define $v_y:= \pi(\rho(y)), v_z:= \pi(\rho(z))$. Thus, 
$$
d_T(v_y,v_z)\le k(r+1)
$$
$$d(\ga_y(v_y), \ga_z(v_z))\le k(r+1) + M_{\bar{K}}. 
$$
Let $\bar{v}\in \{v_y, v_z\}$ denote the vertex further away from $u$. Thus,
$$
d(\ga_y(\bar{v}), \ga_z(\bar{v}))\le k(r+1) + 2M_{\bar{K}}  + 2kK(r+1). 
$$
By \flaring, the fiberwise distance between $\ga_y, \ga_z$ over  the interval $\llbracket \bar{v}, v\rrbracket$ is 
at most
$$
\tau_{\ref{cor:super-weak flaring}}(K, \max(r, k(r+1) + 2M_{\bar{K}}  + 2kK(r+1)))= 
\tau_{\ref{cor:super-weak flaring}}(K, k(r+1) + 2M_{\bar{K}}  + 2kK(r+1)). 
$$

In particular, this inequality holds at the vertex $t$ since it belongs to  
the interval $\llbracket \bar{v}, v\rrbracket$:
$$
d_{X_t}(\ga_y(t), \ga_z(t))\le \tau_{\ref{cor:super-weak flaring}}(K, k(r+1) + 2M_{\bar{K}}  + 2kK(r+1)). 
$$
Hence, by \eqref{gaxgayt}, 
$$
d_{X_t}(\ga_x(t), \ga_z(t))\le M_{\bar{K}} + \tau_{\ref{cor:super-weak flaring}}(K, k(r+1) + 2M_{\bar{K}}  + 2kK(r+1)). 
$$
Define a vertex $v'$ by 
$$
\llbracket \pi(\bar{x}), t\rrbracket\cap \llbracket \pi(\bar{z}),t\rrbracket = \llbracket v', t\rrbracket.  
$$
Since the path $c(x,z)$ has type 2, for all 
$s\in V(\llbracket v', t\llbracket)$ 
we have the inequality
$$
d_{X_s}( \ga_x(s), \ga_z(s))> M_{\bar{K}}. 
$$
Therefore, as in Case 1-1,  Lemma \ref{lem:3-flows} 
 implies a uniform upper bound in the lengths of the intervals 
 $\llbracket \pi(\bar{x}), t\rrbracket, \llbracket \pi(\bar{z}),t\rrbracket$. 
Thus, just as in Case 1-1, we obtain a uniform upper bound on the distances 
$d(\bar{x},  x_1), d(\bar{z}, z_1)$ and, hence, the paths $c(x,y), c(x,z)$ uniformly Hausdorff fellow-travel.

\medskip 
{\bf Case 2-2}: Both paths $c(x,y), c(x,z)$ have type 2. The points $\bar{y}, \bar{z}\in \AA$ are within distance 
$\le k(r+1)$ from each other and, hence, by \flaring, the paths 
$\ga'_y, \ga'_z$ uniformly fellow--travel. 
The same holds for geodesics $[\bar{y} \bar{x}]_A, [\bar{z} \bar{x}]_A$ since 
$A$ is $\delta_A=\delta_{\ref{prop:easy-one}}(K,C)$-hyperbolic. 
Hence, the paths $c(x,y), c(x,z)$ uniformly fellow-travel as well. \qed

Combining this Lemma with Lemma \ref{lem:fellow->proper} and the fact that each path $c(x,y)$ is a concatenation of at most five uniformly quasigeodesic paths, we obtain: 

\begin{cor}\label{wd0}
The paths $c(x,y)$ are uniformly proper. More precisely, there is a function $\zeta_{\ref{wd0}}(r, K, C)$ such that 
for each path $c(x,y)$ defined above, for any two points $x', y'$ on $c(x,y)$, 
if $d(x',y')\le r$ then the length of the portion of $c(x,y)$ between $x', y'$ is $\le \zeta_{\ref{wd0}}(r, K, C)$. 
\end{cor}

\smallskip
{\bf Step 2:}  We shall now check that the paths $c$ defined on Step 1 
satisfy the conditions of Corollary \ref{cor:bowditch} to conclude the proof of the theorem.

\medskip 

{\bf Condition (a1):} This is an immediate consequence of Lemma \ref{wd0}.

\medskip 

{\bf Condition (a2):} The verification of the condition of uniform slimness of the triangle $\Delta$ in $X$ formed by the paths $c(x,y), c(y,z), c(z,x)$ 
is broken into several cases depending on the relative positions of the three points $x,y,z$
and the types of the paths $c(x,y),c(y,z), c(z,x)$ (type 1 or type 2). 
The trick is to reduce the proof to a simpler case by replacing $x,y,z$ with some other
suitable points, analogously to the proof of Proposition \ref{prop:easy-one} and then appeal to Proposition \ref{small-ladder}. For instance, suppose that there is a constant $r=r(K,C)$ such that a triangle $\Delta$ as above is $r$-thin, i.e. there exists a point $x\in X$ within distance $r$ from all three sides of $\Delta$. Then the Hausdorff-fellow-traveling condition (Lemma \ref{wd00}) will imply that $\Delta$ is  $D_{\ref{wd00}}(C, K, 2r)$-thin.

\medskip 
{\bf Case 1:} Suppose all three paths are of type 2. Then we  replace $x,y,z$ by their $\rho$-projections to $A$: The points $\bar{x}, \bar{y},\bar{z}$ respectively. The subtriangle in $\Delta=c(x,y)\cup c(y,z)\cup c(z,x)$ which is the union $\bar{\Delta}= c(\bar{x},\bar{y})\cup c(\bar{y},\bar{z})\cup c(\bar{z},\bar{x})\subset A$ is $\delta_{\A}$-hyperbolic, where $\delta_{\A}$ (depending only on $K$ and $C$) is a uniform bound on the hyperbolicity constant  of the carpet $\A$ (Proposition \ref{prop:easy-one}). Thus, $\Delta'$ is $r=\delta_{\A}$-thin, therefore, 
as we noted above, $\Delta$ itself is $D_{\ref{wd00}} (C, K, 2r)$-slim.

\medskip 
{\bf Case 2:} Suppose we have a triangle $\Delta$ formed by three paths exactly two of which 
are of type $2$; say, $c(x,y),c(y,z)$ are of type $2$ and $c(x,z)$ is of type $1$. 
Since $c(x,z)$ has type 1,  the vertex $t=t_{xz}$ satisfies the inequalities 
$$
d_T(\pi(x), \pi(\bar{x}))\ge d_T(\pi(x), t), \quad  d_T(\pi(z), \pi(\bar{z}))\ge d_T(\pi(z), t).
$$  
 Thus, we can replace the points $x, z$ by $x':= \ga_x(t)$ and $z':= \ga_z(t)$ respectively (as they belong to $c(x,y)$ and $c(z,y)$ respectively).  
 Now, by  the definition of $t=t_{xz}$, $d(x',z')\le d_{X_t}(x', z')\le M_{\bar{K}}$. 
Hence, the subtriangle $\Delta'\subset \Delta$ formed by the paths $c(x',y), c'(y,z'), c(z',x')$ is 
$r=M_{\bar{K}}$-thin, which, in turn, implies that $\Delta$ is $D_{\ref{wd00}} (C, K, 2r)$-slim.

\begin{figure}[tbh]
\centering
\includegraphics[width=40mm]{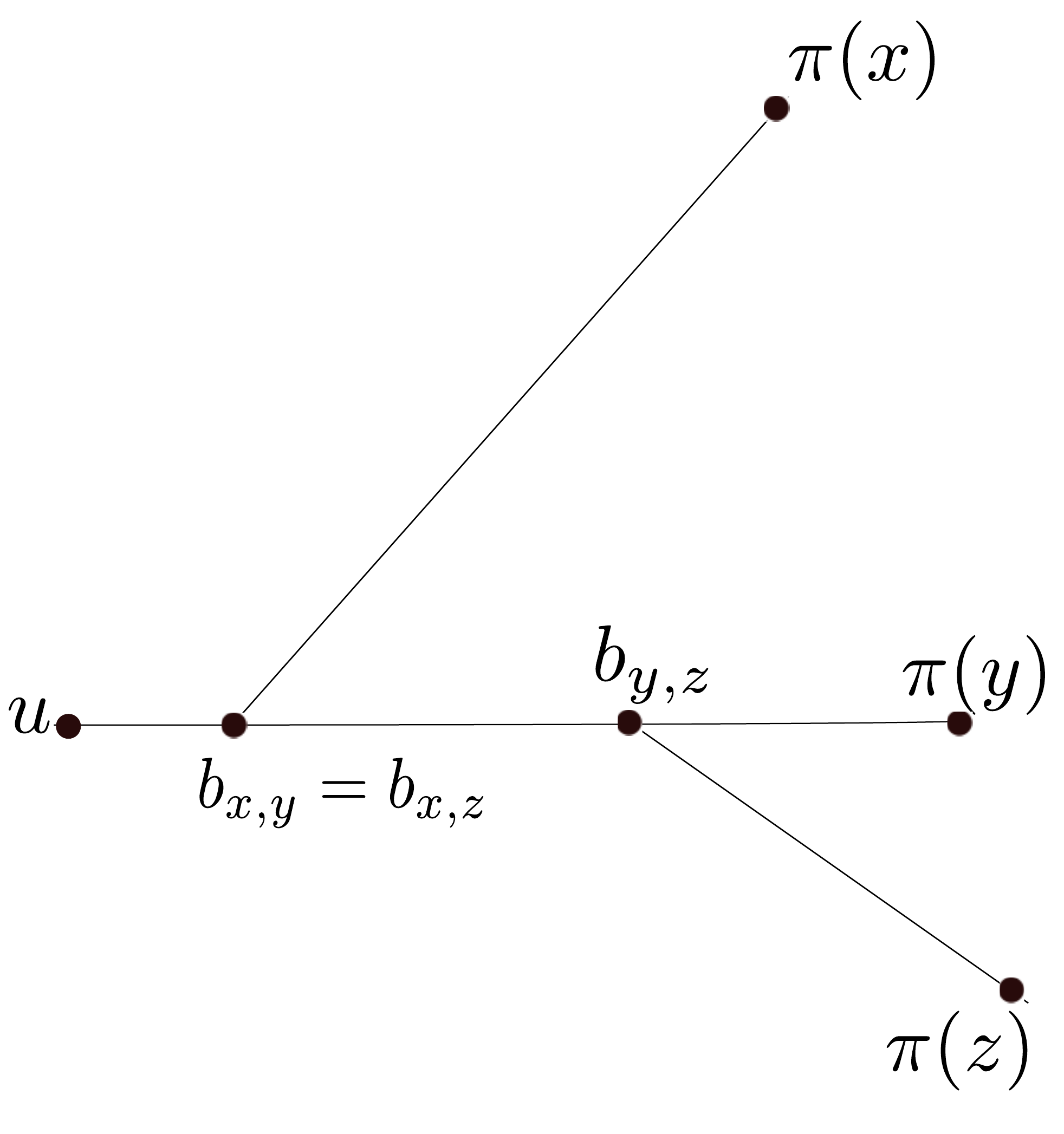}
\caption{Three centers}
\label{three_centers.fig}
\end{figure}

\medskip 
{\bf Case 3:} Suppose we have a triangle $\Delta$ formed by three paths exactly one which, 
say, $c(x,y)$ is of type 2. After swapping $x$ and $y$ we can assume that on the interval $\llbracket u, \pi(z)\rrbracket$ the following vertices appear in this order:
$$
u\le \pi(\bar{z}) \le t_{yz}\le t_{xz}\le \pi(z). 
$$
(Since $c(y,z), c(z,x)$ are of type 1, the vertex $\pi(\bar{z})$ is closer to $u$ than both $t_{yz}, t_{xz}$.)  
As in Case 2, for $t=t_{x,z}$ we replace $x, z$ by the points 
$$
x':= \ga_x(t), \quad z':= \ga_z(t)
$$
respectively. This defines a subtriangle $\Delta'\subset \Delta$ formed by the paths $c(x',y)$, $c(y, z')$, 
$c(z',x')$. By the definition of type 1 paths, $d_{X_t}(x',z')\le M_{\bar{K}}$. 
Thus, again, the subtriangle  $\Delta'\subset \Delta$  is  $r=M_{\bar{K}}$-thin and, 
therefore, $\Delta$ is $D_{\ref{wd00}} (C, K, 2r)$-slim.

\medskip 
{\bf Case 4:} Suppose all three paths $c(x,y), c(y,z), c(z,x)$ forming a triangle $\Delta$ 
are of type 1. This is the most interesting of the four cases.

Without loss of generality we may assume that  in the tree $T$, 
$$
(\pi(x).\pi(y))_{u}\leq (\pi(y).\pi(z))_{u}. 
$$ 
In particular,
$$
b_{x,y}=b_{x,z},
$$
see Figure \ref{three_centers.fig}.

\begin{figure}[tbh]
\centering
\includegraphics[width=80mm]{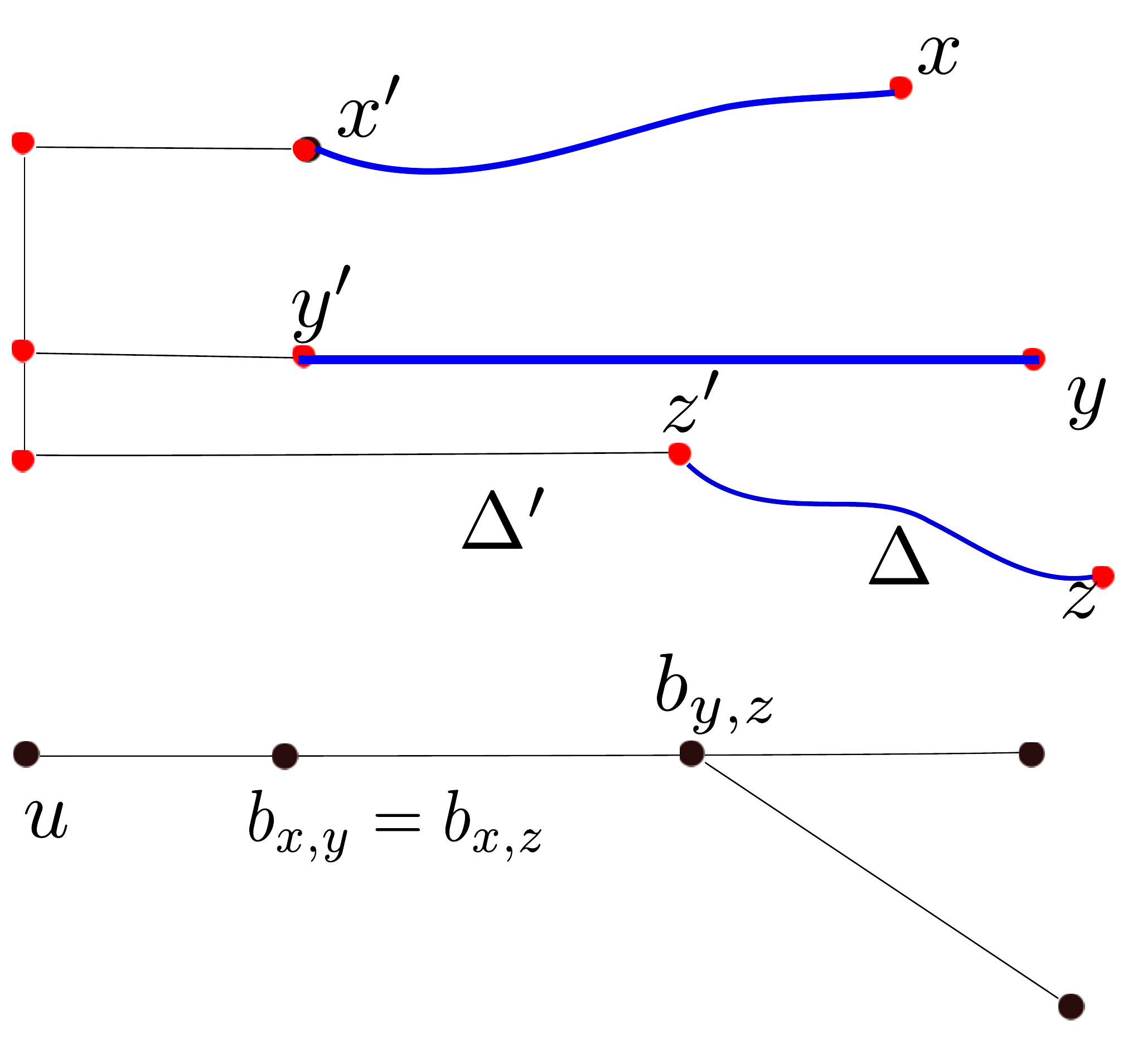}
\caption{Triangles $\Delta, \Delta'$}
\label{fig14.fig}
\end{figure}

Consider the points
$$
x'\in c(x,z)\cap X_{b_{x,z}}, y'\in c(x,y)\cap X_{b_{x,y}}, z'\in c(x,z)\cap X_{b_{y,z}}
$$
which are closest to, respectively, $x, y, z$ along the above paths, see Figure \ref{fig14.fig}. 
The point $z'$ then equals the point of the intersection
$$
c(y,z)\cap  X_{b_{y,z}}
$$
which is closest to $z$ along the path $c(y,z)$. 

Moreover, again by the definition of the paths $c$, the triangle 
$\Delta:= c(x,y)\cup c(y,z)\cup c(z,x)$ is obtained from $\Delta':= c(x',y')\cup c(y',z')\cup c(z',x')$ 
by attaching (to its vertices) the following segments of the $K$-flow-lines:
$$\ga_{x,x'}\subset \ga_x, \ga_{y,y'}\subset \ga_y, \ga_{z,z'}\subset \ga_z.$$
Hence, the $r$-slimness of the triangle $\Delta'$ will imply  $r$-slimness of the triangle $\Delta$.

\begin{figure}[tbh]
\centering
\includegraphics[width=80mm]{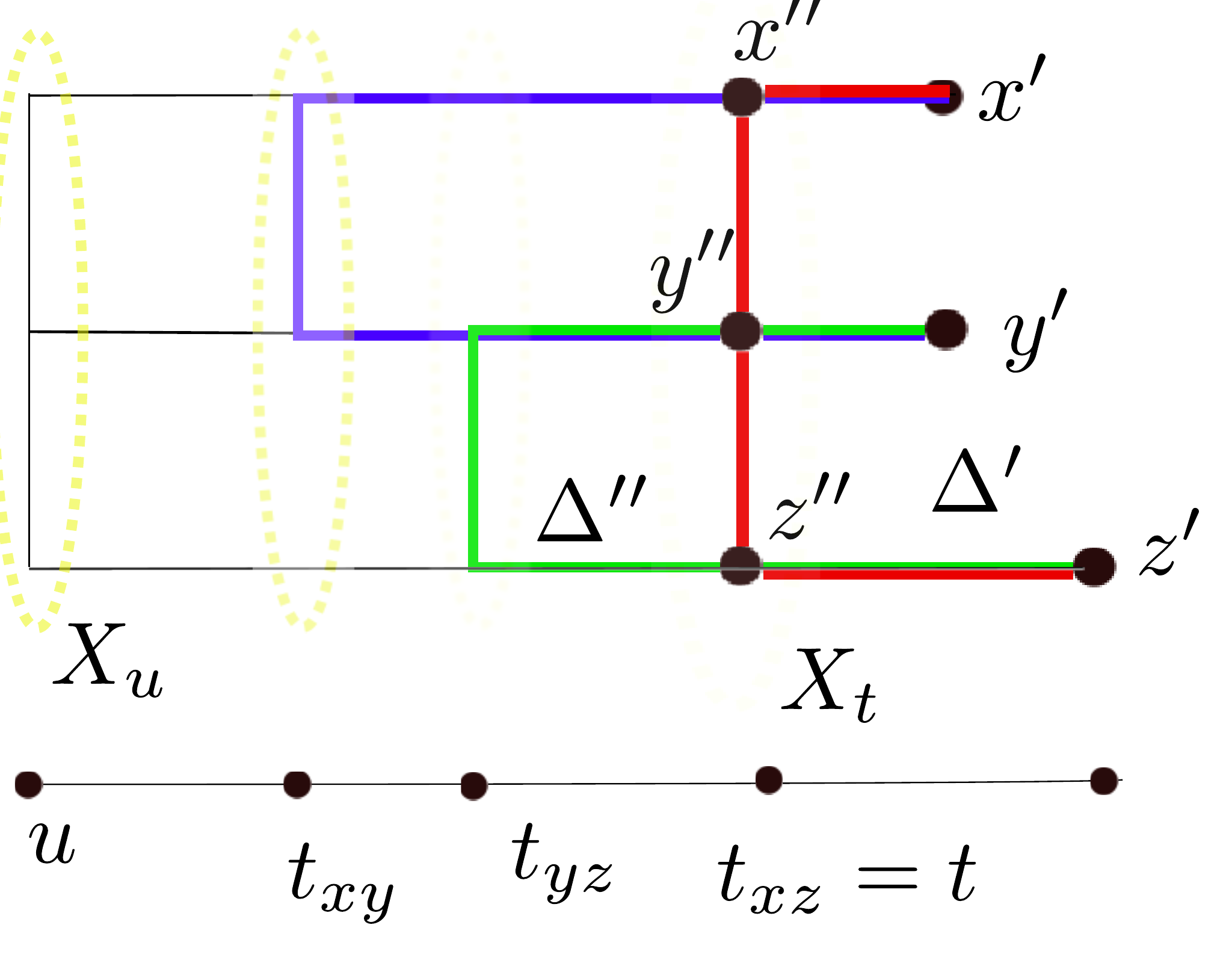}
\caption{Triangles $\Delta', \Delta''$.}
\label{fig15.fig}
\end{figure}

Thus, after replacing $x\to x', y\to y', z\to z'$, it suffices to consider the case when 
$v'=\pi(x)=\pi(y)$ and $v'\in \llbracket u, \pi(z)\rrbracket$.  
Below, {\em we will not be using the property that $\pi(x)=\pi(y)$}, only that all three projections $\pi(x), \pi(y), \pi(z)$ belong to a common  oriented interval $J=\llbracket u,u'\rrbracket\subset T$. Therefore, 
$$
\{t_{zy}, t_{yx}, t_{xz} \}\subset J. 
$$
After a permuting the points $x, y, z$, we can assume that 
$$
t_{zy}\le t_{yx} \le t_{xz} 
$$
on the oriented interval $J$. Therefore, all three paths $\ga_x$, $\ga_y$, $\ga_z$ contains the subpaths 
$\ga_{xx''}$, $\ga_{yy''}$, $\ga_{zz''}$ with $x'', y'', z''\in X_t$, $t=t_{xz}$ and these subpaths are in the respective sides $c(x,y), c(y,z), c(z,x)$ of $\Delta$, see Figure \ref{fig15.fig}.

\medskip 

Thus, we perform one more reduction, replacing the triangle $\Delta$ with the subtriangle $\Delta''$ formed by the paths $c(x'',y''), c(y'',z''), c(z'',x'')$, such that $x'', y'', z''\in X_t, t=t_{xz}$. By the definition of the vertex $t_{xz}$,
$$
d_{X_t}(x'', z'')\le M_{\bar{K}}. 
$$
Thus, the triangle $\Delta''$ is $M_{\bar{K}}/2$-thin, which concludes the proof in Case 4.    
 This also completes the proof of Theorem \ref{small-ladder}. \qed 

\begin{rem}
Since the paths $c(x,y)$ are uniformly proper (Corollary \ref{wd0}), by Part (b) of Corollary \ref{cor:bowditch}, up to a reparameterization, they are  uniform quasigeodesics in $L_K(\al)$. 
\end{rem}

\section{Hyperbolicity of general ladders}\label{sec:hyperbolicity-of-ladders}

In this section we prove that all ladders contained in $X$ are uniformly hyperbolic. The idea is to decompose the given ladder as a union of subladders, each of which is carpeted and then make use of the quasiconvex chain-amalgamation, 
Theorem \ref{thm:hyp-tree}. 

We define $\bar{K}$ and $\kappa$ as in the previous section, Definition \ref{defn:carpeted ladder}, \eqref{eq:barK.4} and \eqref{eq:kappa.4}. 
Recall that each $K$-ladder comes equipped with a family of $K$-sections $\Sigma_\bullet$. These are the sections which will be used in the next proposition. 

\begin{prop}
[Vertical subdivision of general ladders] \label{vertical subdivision}
Fix numbers $K$ and $C$  such that $M_{\bar{K}}\le C$, and suppose that ${\mathfrak X}$ satisfies the uniform $\kappa$-flaring condition. 
Consider a $(K,D,E)$-ladder\footnote{The parameters $D$ and $E$ of the ladder play  no role in the proposition.} 
 $\L=\L_K(\al)$ over a subtree $S\subset T$, where $\al= [x_u y_u]_{X_u}$. 

Then the geodesic segment  $\alpha$ can be subdivided into subsegments 
$\alpha_1,...,\alpha_n$ of lengths $l_1,...,l_n$, with the end-points $x_0=x_u, x_1,...,x_n=y_u$  
such that  the following hold:

 \begin{enumerate}

\item 
The $K$-qi sections $\Si_{x_i}\subset L= L_{K}(\alpha)$ through the points 
$x_i$ for $i=0,...,n$,  
are such that each pair $(\Sigma_{x_i}, \Sigma_{x_{i+1}})$, $i=0,...,n-1$, bounds a $(K,D,E)$-subladder $\L^i=\L_{K,D,E}(\alpha_i)\subset \L$. These subladders satisfy
$$
L^{i}\cap L^{i+1}=\Sigma_{x_i}.
$$

\item Each ladder $\L^i, i=0,...,n-1$, contains a $(K,C)$-narrow carpet $\A^i=\A_{K}(\alpha'_i)$, where 
$\al_i'\subset \al_i$ contains $x_i$ and 
$$
0\le l_i- \length(\al'_i)\le M_K\le M_{\bar{K}},$$
if $i<n-1$, while $\al'_n=\al_n$. 

\item $N^{fib}_{C/2}(\LL^i)\cap N^{fib}_{C/2}(\LL^j)= \emptyset$ provided that $|i-j|>1$.

\item In each ladder $\L^i$, the pair of sections $(\Sigma_{x_i},\Sigma_{x_{i+1}})$ is $B=B_{\ref{vertical subdivision}}(K,C)$-cobounded, for $i=0,\cdots, n-2$. Moreover, the projection of  $\Sigma_{x_{i+1}}$ to $\Sigma_{x_i}$ is uniformly close to the point $x_{w_i}$, where 
$$
[x_{w_i} y_{w_i}]_{X_{w_i}}= A_{w_i}
$$ 
is the narrow end of the carpet $\A^i$. 
\end{enumerate} 
\end{prop}
\proof  
We will orient the interval $\al$ from $x_u$ to $y_u$ and,  will use the corresponding natural order on $\al$. 

Suppose $x_j\in \al$ 
have been chosen and $x_i<y$. To choose $x_{i+1}$ consider the subset
$$
\Omega_{i+1}:= \{x\in \al, x > x_i:  d_{X_b}(\Sigma_x\cap \LL_b, \Sigma_{x_i} \cap \LL_b)> C, ~\forall b\in V(\pi(\Sigma_x)\cap \pi(\Sigma_{x_i}))\}. 
$$
In other words, $\Omega_{i+1}$ consists of points $x>x_i$ in $\al$ such that 
the sections $\Sigma_x$ and $\Sigma_{x_i}$ are fiberwise $C$-separated.

Now there are two possibilities. If $\Omega_{i+1}=\emptyset$ then define $n=i+1$,
$x_{i+1}=y_u$.   
This will conclude the induction. 

Otherwise, pick $x_{i+1}\in \Om_{i+1}$ such that 
$$
x_{i+1}- \inf (\Om_{i+1}) < M_K/2.     
$$

According to Corollary \ref{cor:trisecting a ladder}, 
$\Sigma_{x_i}, \Sigma_{x_{i+1}}$ bound a $(K,D,E)$-subladder 
$\L^i:=\L_{K,D,E}(\alpha_i)\subset \L$, where $\al_i=[x_i x_{i+1}]_{X_u}\subset \al$. 
By the construction, $L^i\cap L^{i-1}= \Sigma_{x_i}$. 

In order to construct a $(K,C)$-narrow carpet $\A^i=\A_K(\al'_i)\subset \L^i$, we take a subsegment 
$\al'_i= [x_i x'_{i+1}]_{X_u}\subset \al_i=[x_i x_{i+1}]_{X_u}$ such that 
$$
x'_{i+1}\notin \Om_{i+1}, d_{X_u}(x'_{i+1}, x_{i+1})<M_{K}.  
$$

Since $x'_{i+1}\notin \Om_{i+1}$, by the definition of $\Om_{i+1}$, the sections 
$\Sigma= \Sigma_{x_i}, \Sigma'=\Sigma_{x'_{i+1}}$ contain points $y, y'\in X_w$ (for some  $w\in V(S)$), 
such that 
$$
d_{X_w}(y, y')\le C. 
$$ 
Among such vertices $v\in T$ we choose one which is closest to $u$ (there might be several). 

The $K$-sections $\Sigma, \Sigma'$ contain subsections $\ga_{y}, \ga_{y'}$ 
(over the interval $\llbracket u, v\rrbracket$) 
connecting $y, y'$ to $x_i, x'_{i+1}$ respectively. 

Since $M_K\le C$ and $w$ was chosen to be closest to $u$, we have that for every vertex $t\in \llbracket u,v\llbracket$,
$$
d_{X_t}(\ga_y(t), \ga_{y'}(t))> C \ge M_K. 
$$
Thus, the sections $\ga_y, \ga'_y$ form the 
top and the bottom of a $C$-narrow $K$-carpet $\A_{K,C}(\al'_{i})$ in $\L^i$ with the narrow end 
$\beta_i= L_w\cap \A^i$.  

This proves parts 1 and 2 of the proposition. 

We next prove part 3. Suppose that $i+1< j$. By the construction, for each vertex 
$v\in \pi(\LL^i)\cap \pi(\LL^j)$ the length of the subsegment in $\LL_v$ between $\Sigma_{x_{i+1}}, \Sigma_{x_j}$ is  $> C$. Hence, the minimal fiberwise distance between 
$\LL^i\cap X_v, \LL^j\cap X_v$ is $> C$ as well. Part 3 follows since $\LL_v$ is a geodesic segment in $X_v$.

Finally to prove (4) we will use the description of the paths $c(x,y)$ constructed for the proof Theorem 
\ref{small-ladder} which are uniform quasigeodesics in the subladder $L^i=L_{K}(\alpha_i)$.

First of all, we observe that, the difference in the lengths of $\al'_i, \al_i$ is at most $M_K$ and $\al'_i$ bounds a 
$(K,C)$-narrow carpet $\A^i$ in $\L^i$. In other words, $\L^i$ is a carpeted ladder and    
Theorem \ref{small-ladder} applies in our case.

The carpet $\A^i= (\pi: A^i\to \llbracket u, w_i\rrbracket)$ is bounded by horizontal paths $\ga_{i}\subset \Sigma_{x_i}, \ga'_i\subset \Sigma_{x_{i+1}'}$ 
and the vertical paths $\al'_i, \beta_i$ (where $\beta_i$ has length $\le C$). 

Consider points $z_i\in \Sigma_{x_i}$ and $z_{i+1}\in \Sigma_{x_{i+1}}$.  
Since the fiberwise vertical separation between $\Sigma_{x_i}$ and $\Sigma_{x_{i+1}}$
 is $> C\ge M_{\bar{K}}$, we conclude that the path $c(z_i,z_{i+1})$ has to be of {\em type 2} in the terminology of 
the proof Proposition \ref{small-ladder}. 

In other words, $c(z_i,z_{i+1})$ is the concatenation of five (actually, four) subpaths: Two of these subpaths (containing $z_i, z_{i+1}$ respectively, connect $z_i$ to $\tilde{z}_i$, $z_{i+1}$ to $\tilde{z}_{i+1}$  and are contained in, respectively, $\ga_i, \ga_{i+1}$. The point $\tilde{z}_{i+1}$ is within vertical distance $\le M_{\bar{K}}$ from a point $\bar{z}_{i+1}\in top(\A^i)$, while the point $\tilde{z}_i$ actually equals $\bar{z}_i$. The vertical geodesic
$$
[\bar{z}_{i+1} \tilde{z}_{i+1}]_{L_{v_{i+1}}}, v_i= \pi(\tilde{z}_{i+1}),  
$$
is contained in $c(z_i, z_{i+1})$. The path $c(z_i, z_{i+1})$ contains one more subpath, namely, 
$$
c_{\A^i}(\bar{z}_i, \bar{z}_{i+1}),
$$
a path in the combing of $A^i$ constructed in the proof of Proposition \ref{prop:easy-one}. By the construction of $\A^i$, each vertex-space $A_v$ of $\A^i$  has length $\ge C\ge M_{\bar{K}}> M_K$. Thus, by the definition of 
 $c_{\A^i}(\bar{z}_i, \bar{z}_{i+1})$ in the proof of Proposition \ref{prop:easy-one}   this path is a concatenation of 
 subpaths contained in the top and the bottom of $\A^i$ and, {\em most importantly}, the narrow end 
 $\beta_i= [x_{w_i} y_{w_i}]_{L_{w_i}}$. 

In particular, each uniform quasigeodesic path $c(z_i,z_{i+1})$ contains the point $x_i^-:=x_{w_i}\in A^i_{w_i}\cap \Si_{x_i}$

According to Lemma \ref{lem:citerion-of-projection},  it follows that the nearest-point projections (taken in $L^i$) of $\Sigma_{x_{i+1}}$ to $\Sigma_{x_i}$ is contained in the  $R$-neighborhoods of $x_{w_i}$, where $R=R(K,C)$. Hence, by Corollary 
\ref{cor:cob-char}, the pair of sections $\Sigma_i, \Sigma_{i+1}$ is $B(K,C)$-cobounded in $L^i$. 
This concludes the proof of the proposition.  \qed 

\medskip
In Section \ref{sec:inductive} (when describing uniform quasigeodesics in ladders) 
we will need a bit more detailed information about the uniformly cobounded pair 
$(\Sigma_{x_i}, \Sigma_{x_{i+1}})$. Let $x_i^-\in \Sigma_{x_i}, x_i^+\in \Sigma_{x_{i+1}}$ be a pair of points realizing the minimal distance between these subsets in the ladder $L^i$. According to Part (4) of the proposition, $x_i^-$ is uniformly close (in terms of $K$ and $C$) to the point $x_{w_i}$. We will also need to identify the other point, $x_i^+$, up to a uniformly bounded error. 

Let $u_{i+1}\in \llbracket u, w_{i+1}\rrbracket$ benote the maximum 
(in the oriented interval $\llbracket u, w_{i+1}\rrbracket$) of the subset 
$$
\{t\in V(\llbracket u, w_{i+1}\rrbracket): d_{X_t}(\Si_{x_{i+1}}\cap X_t), A_t)\le M_{\bar{K}}\}. 
$$
For each vertex $v\in \llbracket u, u_{i+1}\rrbracket$ we observe that the 
top-most point $y'_v$ of the segment $A_{v}$ divides $L^i_v$ in two subsegments:
$$
A_{v}\cup A'_{v}, \quad A'_{v}=  [y'_{v} y_{v}]_{X_v},$$
where 
$$
y_{v}= top(L^i_{v}).$$

The first part of the next lemma will be used in Section \ref{sec:inductive}, while second part 
will be used in Section \ref{sec:Part I} (proof of Lemma \ref{lem:proj-to-detour}).

\begin{lemma}\label{lem:top-projection} 
1. $$
d(y_{u_{i+1}}, x_i^+)\le D_{\ref{lem:top-projection}}(K,C).$$ 

2. For each $v$ as above, the length of $A'_v$ is $\le R_{\ref{lem:top-projection}}(K).$
\end{lemma}
\proof 1. Take a point $y\in \Sigma_{x_{i+1}}\cap X_v$. 
Then, as we noted in the proof of the last part of the proposition, since each path $c(x,y)$ connecting (any) $x\in \Sigma_{x_i}$ to $y$ has type 2, it has to pass through a point of the segment $A'_{u_{i+1}}$, and the latter has length $\le M_{\bar{K}}$. Thus, $c(x,y)$ passes within distance $M_{\bar{K}}$ from $y_{u_{i+1}}$. It follows that $d(y_{u_{i+1}}, x_i^+)$ is uniformly bounded in terms of $K$ and $C$. 

2. This part is an application of uniform flaring. We have two $K$-qi sections $\ga_0, \ga_1$ over the interval $J=\llbracket u, u_{i+1}\rrbracket$, defined by restricting $top(\A^i)$ and $\Si_{i+1}$. The vertical separation  between these sections at the 
 end-points of $J$ is $\le M_{\bar{K}}$. Thus, by \flaring, the vertical separation between these sections elsewhere in the interval is 
 $\le \tau_{\ref{cor:super-weak flaring}}(K,M_{\bar{K}})$. Hence, we can take 
 $$
 R_{\ref{lem:top-projection}}(K)=\tau_{\ref{cor:super-weak flaring}}(K,M_{\bar{K}}). \qed $$

\begin{theorem}\label{thm:ladders-are-hyperbolic} 
For each $K$ and a hyperbolic tree  spaces ${\mathfrak X}$ 
satisfying the uniform $\kappa=\kappa_{\ref{defn:carpeted ladder}}(K)$-flaring assumption as in  
Proposition \ref{small-ladder}, and arbitrary $D$ and $E$, 
there exists  $\delta=\delta_{\ref{thm:ladders-are-hyperbolic}}(K)$ 
such that  every $(K,D,E)$-ladder $\L=\L_K(\al)\subset \X$  has  
$\delta$-hyperbolic total space $L$ with respect to the intrinsic metric of the ladder. 
\end{theorem} 
\proof  Set $C=M_{\bar{K}}$. By Proposition \ref{vertical subdivision} we have a subdivision $x_0=x_u, x_1,...,x_n=y_u$ of the segment 
$\al=[x_uy_u]_{X_u}$.  The $K$-qi sections $\Sigma_i:=\Sigma_{x_i}$ in $\L$ passing through $x_i$'s  
decompose the ladder 
$\L$ into subladders $\L^i:=\L_{K}([x_{i-1} x_i]_{X_u})$  containing $(K,C)$-narrow  carpets $\A^i$. Hence, by Proposition \ref{small-ladder}, each $L^i$ (the total space of $\L^i$) 
is $\delta_{\ref{small-ladder}}(K,C)$-hyperbolic. By the construction, 
$$
L^i\cap L^{i+1}=\Sigma_i, i=1,...,n-1.
$$
The subsets $\Sigma_i$ are $K$-sections, hence, are  
$\la_{\ref{lem:qi-preserves2}}(\delta_{\ref{small-ladder}}(K,C), K)$-quasiconvex in $L^i, L^{i+1}$. 
Furthermore, by Proposition \ref{vertical subdivision},  each pair of ladders 
$\L^{i-1}, \L^{i+1}$, $i\ge 1$, is $B_{\ref{vertical subdivision}}(K,C)$-coboun\-ded. 
Thus, arbitrary ladder $L=L_K(\al)$ is uniformly hyperbolic by Theorem \ref{thm:hyp-tree}. \qed

\chapter{Hyperbolicity of flow-spaces}\label{ch:flows} 

In this chapter we shall prove that the $k$-flow-spaces (for $k$ in a suitable range) 
of each vertex space $X_u\subset X$, are uniformly hyperbolic (with hyperbolicity constant depending on $k$) 
 provided that $\X$ is a hyperbolic tree of spaces satisfying a certain uniform flaring condition. 
The strategy of the proof is to show that:

(a) Every two points $x, y$ in $Fl_k(X_u)$ belong to a common a ladder $(K,D,E)$-ladder $L_{x,y}$ (essentially contained in  $Fl_k(X_u)$), where $D$ is a sufficiently large number, $K$ depends only on $k$ and $E$ depends on $k$ and $D$. This is done in Section \ref{sec:ubiquity}. Actually, in Section \ref{sec:ubiquity} we prove a stronger result, 
the existence of {\em tripods of ladders} connecting points $x_1, x_2, x_3\in Fl_k(X_u)$ such that in each vertex-space of ${\mathfrak Fl}_k(X_u)$ the three geodesics of these ladders form a geodesic tripod. Hyperbolicity of total spaces of such tripods of ladders is almost immediate, Section \ref{sec:Hyperbolicity of tripods families}.

(b) Since ladders are uniformly hyperbolic (as it was proven in the previous chapter), 
this appears to yield a preferred family of paths $c(x,y)$ connecting points of $Fl_k(X_u)$ (projections of uniform quasigeodesics in corresponding ladders). Hyperbolicity of tripods of ladders should then yield the uniform slimness condition for the family of paths $c(x,y)$. 
The trouble, however, is that $L_{x,y}$ is far from canonical, and, thus, it is far from clear why the paths $c(x,y)$ satisfy the fellow-traveling condition. 
If different ladders $L^1_{x,y}, L^2_{x,y}$ were at uniformly bounded minimal distance from each other in each fiber-space where both ladders are nonempty, one could 
use hyperbolicity of the union of these two ladders. Unfortunately, it is unclear why there should be a uniform bound on such minimal distance. To resolve the problem, we use the 
 construction of a {\em coarse projection} of the ladder $L^1_{x,y}$ to $L^2_{x,y}$ defined in Section \ref{sec:proj-ladders}. This projection is used in Proposition \ref{prop:ladder-amalgam} to construct a uniformly  hyperbolic subspace $Z$ in $X$ containing the two ladders.  
The coarse independence of the paths $c(x,y)$ on the choice of $L_{x,y}$ is then almost immediate, Corollary \ref{cor:ladder-amalgam}. 

This, in turn, will conclude the verification of hyperbolicity of flow-spaces, Theorem \ref{flow of one vertex space}.

\section{Ubiquity of ladders in $Fl_k(X_u)$}\label{sec:ubiquity} 

In this section we prove that for all  ({sufficiently large}) $D$ and $k$, and all $x,y\in Fl_k(X_u)$, 
there is a $(K,D,E)$-ladder  $\L$ containing $x,y$, where $K=K(k), E=3k+D$. 
Furthermore, $\L$ will be contained in the fiberwise $4\delta_0$-neighborhood of ${\mathfrak Fl}_K(X_u)$. 
We will actually prove a stronger result, about the existence of a {\em tripod} of ladders containing given three points in $Fl_K(X_u)$. 

In the next definition and in what follows, $i$ is taken modulo $3$.

\begin{defn}\label{defn:ladder-triangle}\index{tripod of ladders}
A $(K,D,E)$-{\em tripod of ladders} in $\X$ is a semicontinuous $(K,D,E)$-family $\Y$ over a subtree $S\subset T$, which is a union of three $(K,D,E)$-ladders $\L^i=(L^i\to S_i), i=1,2,3$, such that:

1. There exists a $K$-section $\Xi$ defined over the subtree 
$$
S_{123}:= S_1\cap S_2\cap S_3  
$$
and called the {\em center-section} of $\Y$, such that for each $v\in V(S_{123})$ and $i\in \{1,2,3\}$, 
$\Xi(v)=top(\LL^i)\cap X_v$.

2. For each $v\in V(S_{123})$, $e\in E(S_{123})$, the 
vertex- and edge-space $Y_v, Y_e$ of $\Y$ is a $\delta_0$-tripod in $X_v, X_e$ 
$$
Y_v= \bigcup_{i=1}^3 [x^i_v z_v]_{X_v}, Y_e= \bigcup_{i=1}^3 [x^i_e z_e]_{X_e}, z_v=\Xi(v), L^i_v= [x^i_v z_v]_{X_v}, L^i_e= [x^i_e z_e]_{X_e}. 
$$

3. If $v\in V(S_i) \setminus V(S_{123})$, then $L_v^{i+1}= L_v^{i-1}=\emptyset$ (two legs of the tripod are missing) 
and $L^i_v= L_v= [x_v y_v]_{X_v}$ (we, thus, omit the superscript $i$ in this situation). The same applies to the edges. 

We will refer to the union of bottoms of the ladders $\L^i$ as the {\em bottom} (denoted $bot(\Y)$) of the 
tripod of ladders $\Y$. 

A tripod of ladders $\Y$ is said to be {\em degenerate} if for some $i\in \{1,2,3\}$, $\L^{i+1}= \L^{i-1}= \Xi$ and, thus, $\Y$ is reduced to a single ladder, $\L^i$. 
\end{defn}

In the following proposition, $D_0=D_{\ref{cobdd-cor}}(\delta'_0, \la'_0)$,  
$$
C= 2(\la'_0+2\delta'_0 + D_{\ref{stab-qg}}(\delta'_0, L'_0)) + C_{\ref{cor:projection-2}}(\delta'_0,\la'_0), 
$$ 
$$
D_1:=\max(D_0, C_{\ref{cor:cob-char}}(\la'_0,\delta'_0,C)),$$
\begin{equation}\label{eq:D1}
D= D_{\ref{prop:existence-of-tripod-ladders}}:= D_1+ \max( 3\delta_0, 2\delta_0 + 2(\la'_0+2\delta'_0)). 
\end{equation}
Assume also that  $k=r^\wedge=(15 L'_0 r)^3\ge K_0$ and $r$  satisfies the inequality 
$$
r\ge r_1=\max(2\la'_0+ 5\delta'_0, \la +4\la'_0+8\delta'_0+ 5\delta_0),
$$
where
$$
\la=C_{\ref{lem:proj-to-tripod1}}(\delta'_0,\la'_0,L'_0). 
$$
In other words,
\begin{equation}\label{eq:K1}
k\ge k_{\ref{prop:existence-of-tripod-ladders}}=K_1:=  \max(K_0, (15 L'_0 r_1)^3).  
\end{equation}
Note that, in particular, $k\ge  \la+4\la'_0+8\delta'_0+ 5\delta_0$.

In the proposition we will be also using  the function $\kappa'=K'_{\ref{lem:E-ladder-structure}}(\kappa)$ 
defined in Lemma \ref{lem:E-ladder-structure}.

\begin{prop}
[Existence of tripods of ladders] \label{prop:existence-of-tripod-ladders}
Let $\X=(\pi: X\to T)$ be a tree of spaces satisfying Axiom {\bf H1}.  
Then for  $k$ and $D$ as above, there exist 
constants 
$$K= K_{\ref{prop:existence-of-tripod-ladders}}(k),  
E=E_{\ref{prop:existence-of-tripod-ladders}}(k)$$
such that the following holds. 

For points $x^i, i=1, 2, 3$ in ${\mathcal Fl}_k(X_u)$, we let $\ga^i:= \ga_{x^i}$ denote $k$-sections in 
$Fl_k(X_u)$ over $\llbracket u,\pi(x^i)\rrbracket$, connecting $x^i$ to $X_u$. 

Then: 

(i) There exists a $(K,D,E)$-tripod of ladders $\Y= \L^1\cup \L^2\cup \L^3$,  
centered at $u$ such that:   

1. $\YY\subset N_{5\delta_0}^{fib} {\mathcal Fl}_k(X_u)$, while $bot(\YY)\subset {\mathcal Fl}_k(X_u)$.  

2. $\ga^i\subset bot(\L^i)$, $i=1,2,3$, thus, $\ga^i\subset bot(\Y)$. 

3. If, for some $i$, $\ga^{i-1}=\ga^{i+1}$, 
then the tripod of ladders $\Y$ is degenerate and the section $\ga^{i-1}=\ga^{i+1}$ is contained in the 
center-section $\Xi$ of $\Y$. 

(ii) There exist $(K,D,E)$-ladders $L^{ij}$ containing $x^i, x^j$,  such that 
 $top(L^{ij})\subset bot(L^j)$, $bot(L^{ij})\subset bot(L^j)$, and $L^{ij}$ is 
 contained in $\delta_0$-fiberwise neighborhood of $\L^i\cup \L^j$. 
\end{prop}
\proof We first note that, according to Lemma \ref{max-qi-section}, given $x\in Fl_k(X_u)$, there exist a maximal 
$K$-section $\Sigma_x\subset Fl_k(X_u)$ through $x$, intersecting $X_u$ and containing $\ga_x$. Thus, we define the maximal $k$-sections $\Sigma_{x^i}, i=1,2,3$ through points $x^i$, intersecting $X_u$, 
and containing $\ga^i$ (if it is given), otherwise, chosen arbitrarily. 
(Note that these sections $\Sigma_{x^i}$ need not be disjoint and, in general, they have different domains in $T$.) In line with Part 3 of the proposition, if $\ga^{i-1}=\ga^{i+1}$, we require $\Si_{x^{i-1}}=\Si_{x^{i+1}}$. We define tripods $Y_v, v\in V(T)$ inductively, by induction on the distance from $v$ to $u$. 

As the base of induction, we define $x_u^i$ as $\ga^i(u)$ and $L_u^i$ as a geodesic segment $[x_u^i z_u]_{X_u}$, where  
$z_u$ is a $\delta_0$-center of the geodesic triangle $\Delta_u=\Delta x^1_u x^2_u x^3_u$ in $X_u$.

We then proceed inductively as in the proof of Lemma \ref{lem:E-ladder-structure} 
to which we refer the reader for the notation used below. Namely, assume that segments $L_v^i$ 
are defined for vertices of the subtree 
$B_n\subset T$,
$$
L^i_v= [ x^i_v z_v]_{X_v}, v\in V(B_n), Y_v= L^1_v\cup L^2_v\cup L_v^3, 
$$ 
where $z_v$ is a $\delta_0$-center of the geodesic triangle $\Delta_v= \Delta x^1_v x^2_v x^3_v$. 
(In order to simplify the notation, we allow for empty vertex and edge-spaces in ladders.) 
We let $d_{Y_v}$ denote the intrinsic path-metric on $Y_v$. 
Then the inclusion map $(Y_v, d_{Y_v})\to X_{vw}$ is an  
$L'_0(2\delta_0)$-qi embedding for each edge $[v,w]\in E(T), v\in B_n$, directed away from $u$. 
We will be also assuming (inductively) that the extremities $x^i_v$ of $Y_v$ 
belong to $Q_v=Fl_k(X_u)\cap X_v$ and, thus, $Y_v$ is contained in the 
$5\delta_0$-neighborhood of $Q_v$ taken in $X_v$. 

\medskip 
Next, we apply the modified projection $\bar{P}_{X_{vw},Y_v}$ to $X_w$. 
 The image $\bar{Y}_v$ is the closed convex hull of ${P}_{X_{vw},Y_v}(X_w)$ in the tripod $Y_v$ 
(with respect to the path-metric of $Y_v$). 
Specifically, if $\bar{Y}_v\ne\emptyset$, then it is the convex hull (taken in $Y_v$) of three points $\bar{x}^i_v, i=1, 2, 3$, 
such that $\bar{x}^i_v$ is the nearest-point projection (taken in $(Y_v, d_{Y_v})$) of $x^i_v$ to $\bar{Y}_v$.

\medskip 
Note that it is entirely possible for the center $z_v$  not to belong to $\bar{Y}_v$, in which case  $\bar{Y}_v$ is a 
segment contained in the relative interior of one of the legs of $Y_v$. 
It even can happen that $\bar{Y}_v$ is empty, if $Q_w=\emptyset$. Our next task is to analyze  implications of the 
containment $\bar{Y}_v\subset L^i_v$. 

\begin{lemma}\label{lem:leg-proj} 
Suppose that $\bar{Y}_v$ is contained in $L^i_v, d_T(u,v)=n$ 
and $e=[v,w]$ is an edge of $T$ oriented away from $u$. 
Then for $j=i\pm 1$ we have 
$$
P_{X_{vw},L_v^j}(X_w)\subset B^e(z_v, C),$$
where $B^e(z_v,C)$ is the ball with respect to the metric of $X_{vw}$ and $C=C_{\ref{lem:leg-proj}}$. 
\end{lemma}
\proof According to Corollary \ref{cor:projection-2}, for the $\la'_0$-quasiconvex subsets $V=L_v^j, U= Y_v$ in $X_{vw}$ we have  
$$
d_{X_{vw}}(P_{X_{vw},V}, P_{U,V}\circ P_{X_{vw},U})\le 
C_{\ref{cor:projection-2}}(\delta'_0,\la'_0).
$$ 
Here all the projections are taken with respect to the restriction of the metric on $X_{vw}$. 
Thus, we need to prove that $P_{U,V}(L^i_v)$ is uniformly close to the point $z_v$. This is obviously true for the intrinsic nearest-point projection $P'_{U,V}(L^i_v)$ (taken with respect to the metric $d_{Y_v}$), 
since $P'_{U,V}$ sends $L^i_v$ to $\{z_v\}$.  Therefore, we need to compare the 
intrinsic projection $Y_v\to L_v^j$ and the extrinsic projection, with respect to the metric of $X_{vw}$. 

Take some $x\in L^i_v$ and let 
$\bar{x}= P_{X_{vw},V}(x)$. According to Lemma \ref{lem:projection-1}, the geodesic $\al^*=[xz_v]_{X_{vw}}$ 
passes within distance $\la'_0+2\delta'_0$ from the point $\bar{x}$. Since the segment  $\al=[xz_v]_{X_{v}}$ is an $L'_0$-quasigeodesic in $X_{vw}$, it follows that 
$$
\Hd_{X_{vw}}(\al,\al^*)\le D_{\ref{stab-qg}}(\delta'_0, L'_0). 
$$
Hence, $\al^*$ contains a point $y$ satisfying
$$
d_{X_{vw}}(y, z_v)\le \la'_0+2\delta'_0 + D_{\ref{stab-qg}}(\delta'_0, L'_0). 
$$
 Since
$$
d_{X_{vw}}(x, \bar{x})\le d_{X_{vw}}(x, z_v),
$$
we get
$$
d_{X_{vw}}(\bar{x}, z_v)\le 2(\la'_0+2\delta'_0 + D_{\ref{stab-qg}}(\delta'_0, L'_0)). 
$$
Combining this estimate with Corollary  \ref{cor:projection-2}, and the hypothesis that  $P_{X_{vw},U}(X_w)\subset L^i_v$, 
we conclude that for each point $q\in X_w$, $P_{X_{vw},U}(q)=x\in L^i_v$ and 
\begin{align*}
d_{X_{vw}}(P_{X_{vw},L_v^j}(q), z_v)\le C_{\ref{lem:leg-proj}}:= \\
2(\la'_0+2\delta'_0 + D_{\ref{stab-qg}}(\delta'_0, L'_0)) + C_{\ref{cor:projection-2}}(\delta'_0,\la'_0). \qed  
\end{align*}

Combining the lemma with Corollary \ref{cor:cob-char}, yields: 

\begin{cor}
If $\bar{Y}_v$ is contained in $L^i_v$, then for $j=i\pm1$, the pair $(L_v^j, X_w)$ is $C_{\ref{cor:cob-char}}(\la'_0,\delta'_0,C)$-cobounded in $X_{vw}$, where $C= C_{\ref{lem:leg-proj}}$. 
\end{cor}

\begin{rem}\label{rem:cbd1}
The assumption $D\ge D_1$ made in the proposition ensures that $D\ge C_{\ref{cor:cob-char}}(\la'_0,\delta'_0,C)$. 
Thus, $\bar{Y}_v\subset L^i_v$ implies that the pair $(L_v^j, X_w)$ is $D_1$-cobounded in $X_{vw}$, hence, $D$-cobounded.  
\end{rem}

\medskip  
We now return to the construction of a family of tripods. 
Let $e=[v,w]$ be an edge directed away from $u$, $v\in B_n\subset T, w\notin B_n$.  
 There are several things which can now happen, primarily depending on the 
coboundedness  of $Y_v$ and $X_w$, but also on intersections of the sections $\Si_{x^i}$ with $X_w$.  

\medskip 
{\bf Case 1:} Suppose that the tripod $Y_v$ and $X_w$ are $D_1$-cobounded (in $X_{vw}$) and 
all three sections $\Si_{x^i}$ are disjoint from $X_w$. 
Then we set $Y_w=\emptyset$.

\begin{rem}\label{rem:cbd2}
Observe that if the pair $(Y_v,X_w)$ is $D_1$-cobounded, so are the pairs $(L_v^i, X_w)$, $i=1, 2, 3$.  
\end{rem}

\medskip 
{\bf Case 2:} Suppose that the tripod $Y_v$ and $X_w$ are not $D_1$-cobounded. We will also assume that the tripod $Y_v$ has ``all its legs,'' i.e. $L^i_v\ne \emptyset$, $i=1,2,3$.  According to Lemma \ref{lem:proj-to-tripod1}, 
\begin{equation}\label{eq:proj-to-tripod1}
\bar{Y}_v\subset N_{\la} ({P}_{X_{vw},Y_v}(X_w)),
\end{equation}
where, as before, 
$$
\la=C_{\ref{lem:proj-to-tripod1}}(\delta'_0,\la'_0,L'_0). 
$$
We now use the fact that $D_1\ge D_0= D_{\ref{cobdd-cor}}(\delta'_0, \la'_0)$.   
 Since $Y_v$ and $X_w$ are not $D_1$-cobounded, by Corollary \ref{cor:noncbd}, 
$$
{P}_{X_{vw},Y_v}(X_w) \subset N^e_{4\la'_0+8\delta'_0}(X_w)\cap Y_v,
$$
hence,
$$
\bar{Y}_v\subset N^e_{\la+4\la'_0+8\delta'_0} (X_w)\cap Y_v. 
$$
However, $Y_v\subset N^{fib}_{5\delta_0}(Q_v)$, by the inductive hypothesis. 
Hence, each point $x\in X_w$ within distance (in $X_{vw}$) 
$$
\la+4\la'_0+8\delta'_0+ 5\delta_0$$
from $\bar{Y}_v$, belongs to $Q'_w=N^e_r(Q_v)\cap X_w$, 
$r=k^\vee$ (see the definition of flow-spaces in Section \ref{sec:flow-spaces}). 
Since, by the assumption of the proposition, 
$$
{k^\vee > \la+4\la'_0+8\delta'_0+ 5\delta_0},
$$
we see that  
$$
\bar{Y}_v\subset N^e_{\la+4\la'_0+8\delta'_0} (Q_w)\cap Y_v. 
$$

\medskip
Since the extremities $\bar{x}^i_v, i=1, 2, 3$, of the (possibly degenerate) tripod $\bar{Y}_v$ are at the distance 
$4\la'_0+8\delta'_0+ 5\delta_0$ from $Q_w$, we take points 
$\tilde{x}^i_w\in Q_w$ which are nearest-point projections of $\bar{x}^i_v, i=1, 2, 3$, and, thus, 
\begin{equation}\label{eq:moving-feet} 
d_{X_{vw}}(\bar{x}^i_v, \tilde{x}^i_w)\le 4\la'_0+8\delta'_0+ 5\delta_0. 
\end{equation}
(If  $\bar{x}^i_v=\bar{x}^j_v$ then $\tilde{x}^i_w=\tilde{x}^j_w$.) 
Similarly, if $z_v$ belongs to $\bar{Y}_v$ then there exists 
$\tilde{z}_w\in Q_w$ satisfying  
$$
d_{X_{vw}}(z_v, \tilde{z}_w)\le \la+4\la'_0+8\delta'_0+ 5\delta_0. 
$$
In this case, we  define the tripod  
$\tilde{Y}_w\subset X_w$ centered at $\tilde{z}_w$ and with the legs $[\tilde{z}_w \tilde{x}^i_w]_{X_{w}}$.

The actual tripod $Y_w$, as we will see, is uniformly Hausdorff-close to $\tilde{Y}_w$.
For now, we observe that, according to Lemma \ref{lem:sub-close}:

\begin{equation}\label{eq:tYbY} 
\Hd_{X_{vw}}(\tilde{Y}_w, \bar{Y}_v)\le D_{\ref{lem:sub-close} }(\delta'_0,L'_0,4\la'_0+8\delta'_0+ 5\delta_0). 
\end{equation}

Depending on the intersections 
$\Si_{x^i}\cap X_w$, the points $\tilde{x}^i_w$ might be the vertices of the tripod $Y_w$. 
Specifically, there are four subcases:  

(a) If for some $i$, $\Si_{x^i}\cap X_w= \{x^i_w\}$, then we use the point $x^i_w$ 
 as one of the vertices of $\Delta_w$. Thus,
 $$
 d_{X_{vw}}(x^i_v, x^i_w)\le k
 $$
 in this subcase. 
 
(b) If for some $i$, $\bar{Y}_v$ is disjoint from (necessarily both) $L_v^{i\pm 1}$ 
and $\Si_{x^{i+1}}\cap X_w=\emptyset$, 
 $\Si_{x^{i-1}}\cap X_w=\emptyset$, we set $L^{i\pm 1}_w=\emptyset$. Thus, in this subcase 
 the tripod $Y_w$ will be missing two legs. This degenerate tripod will 
 be equal the oriented geodesic segment $L^i_w= L_w=[x_w y_w]_{X_w}$, where $x_w =\tilde{x}_w^{i\pm 1}$ and 
 $y_w$ will be either $\tilde{x}_w^{i}$ (if $\Si_{x^i}\cap X_w=\emptyset$) or, as in subcase (a),  $\Si_{x^i}\cap X_w=\{y_w\}$. 
 In this situation, by the construction, for 
 $$
 \hat{x}_v:= \bar{x}_v^{i\pm 1}
 $$
 $$
 d_{X_{vw}}(x_w, \hat{x}_v)\le 4\la'_0+8\delta'_0+ 5\delta_0,
 $$
 while
 $$
 d_{X_{vw}}(y_w, \hat{y}_v)\le \max(k,4\la'_0+8\delta'_0+ 5\delta_0)=k. 
 $$
 Here $\hat{y}_v=\bar{x}^i_v$ (if $\Si_{x^i}\cap X_w=\emptyset$) or $\hat{y}_v=x^i_v$ (otherwise).  
 
 \begin{rem}\label{rem:cbd3}
 In this subcase, due to our assumptions on $D_1$, both pairs $(L_v^{i\pm 1}, X_w)$ will be $D_1$-cobounded, see Remark 
 \ref{rem:cbd1}.  
 \end{rem}

 (c)  If for some $i$, $\bar{Y}_v$ is disjoint from (necessarily both) $L_v^{i\pm 1}$ 
 and for exactly one element $j\in \{i\pm 1\}$, the section $\Si_{x^j}$ 
 intersects  $X_w$, then we discard the point $\bar{x}_v^{i\pm 1}$  and 
 let $z_w=x_w^{i-1}= x_w^{i+1}$ be that point of intersection.  
 We let $x_w^i$ either be equal to the intersection point of   
 $\Si_{x^i}$ and $X_w$ (if the intersection is nonempty)   or equal to $\tilde{x}_w^i$. Thus,
 $$
 d_{X_{vw}}(z_w, \bar{x}_v^j)\le k, 
 $$ 
 while either
 $$
 d_{X_{vw}}(x^i_w, \hat{y}_v)\le \max(4\la'_0+8\delta'_0+ 5\delta_0,k)=k,
 $$
 where, as in the subcase (b), $\hat{y}_v= \bar{x}^i_v$, or $\hat{y}_v=x_v^i$.

 (d) In the ``generic'' case (i.e. when $z_v\in \bar{Y}_v$),  for each $i$ such that  
 $\Si_{x^i}\cap X_w=\emptyset$, we set $x_w^i=\tilde{x}_w^i$. (Of course, if 
  $\Si_{x^i}\cap X_w$ is nonempty, we use this intersection point as $x^i_w$, see subcase (a).)  As above, we obtain:
  $$
 d_{X_{vw}}(x^i_w, \{\bar{x}^i_v, x_v^i\})\le 
 \max(4\la'_0+8\delta'_0+ 5\delta_0,k)=k, i=1, 2, 3. 
 $$

\medskip 
Except for the subcase (b), we, thus, obtain three points $x_w^1, x_w^2, x_w^3$ spanning a 
(possibly degenerate) geodesic 
triangle $\Delta_w=\Delta x_w^1 x_w^2 x_w^3\subset X_w$.  We let $z_w$ be a $\delta_0$-center of this triangle. 
(The subcase (c) above does not cause trouble because the triangle $\Delta_w$ is degenerate and one of its sides equals the 
point $z_w$, which is, therefore, the center of $\Delta_w$.) 
Accordingly, we  define geodesic segments
$$
L^i_w:= [ x^{i}_w z^i_w]_{X_w}  
$$
and the tripod $Y_w= L^1_w\cup L_w^2\cup L^3_w$. 

The subcase (b) requires a separate discussion since the tripod $Y_w$ is missing 
 two out of its three legs.  In this situation, by the definition of the point $x_w$, 
$$
 d_{X_{vw}}(x_w, \bar{x}_v^{i\pm 1})\le \la+4\la'_0+8\delta'_0+ 5\delta_0\le k,  
 $$ 
 by the assumption on $k$ made in the proposition.

\medskip 
{\bf Case 3.} We still assume that $Y_v$ and $Q_w$ are not $D_1$-cobounded, but consider the case that $Y_v$ is degenerate and 
has only one leg, $L^i_v= L_v= [x_v y_v]_{X_v}$: The other two legs are empty. We treat this case exactly the same way as 
we treated the subcases (2b) and (2c) above: The tripod $\bar{Y}_v$ has empty intersection with the empty legs 
$L_v^{i\pm 1}$ of $Y_v$. The points $x_w, y_w\in X_w$ define the oriented segment $L_w= [x_w y_w]_{X_w}$ and 
the points $x_w, y_w$ are within distance $k$, respectively, from points $\hat{x}_v, \hat{y}_v$, where (as in 
subcase (2b)) $\hat{x}_v\in \{\bar{x}_v, x_v\}$, $\hat{y}_v\in \{\bar{y}_v, y_v\}$, 
$$
\bar{Y}_v= [\bar{x}_v\bar{y}_v]_{X_v}\subset L_v 
$$ 
By the definition, the points $\hat{x}_v, \hat{y}_v$ satisfy
$$
x_v\le \hat{x}_v\le \hat{y}_v\le y_v
$$
in the oriented segment $[x_v y_v]_{X_v}$. (Compare Lemma \ref{lem:E-ladder-structure}(a3).)

\medskip 
This concludes the construction of the segments $L^i_w$. 
We just note that  
since $Q_w\subset X_w$ is 
$4\delta_0$-quasiconvex and $x^i_w\in Q_w$ for all $i$, we get: 
$$
Y_w\subset N^{fib}_{5\delta_0}(Q_w).$$
These are the inductive assumptions we made earlier. We set
$$
\LL^i:= \bigcap_{v\in V(T)} L^i_v, i=1,2,3.  
$$
We define a subtree of spaces  $\Y\subset \X$ using the tripods $Y_v, Y_e$ as, respectively, vertex and edge-sets. The incidence maps $Y_e\to Y_v$ will be compositions of restrictions of incidence maps of $\X$ with nearest-point projections $X_v\to Y_v$. 

We, are done with the induction but it remains to verify that each $\LL^i$ satisfies the 
properties required by Lemma \ref{lem:E-ladder-structure}: This 
lemma is used to promote the unions of geodesics segments in vertex spaces of 
$X$ to the union of vertex-spaces of a ladder. We also have to show that 
$\Y$ is a $(\kappa_1,D_1,E_1)$-semicontinuous family in $\X$, 
as required by the definition of a tripod of ladders, for suitable constants $\kappa_1, E_1$.  

\medskip 
As we observed in the discussion of subcases, points $x_w^i$ satisfy  
\begin{equation}\label{eq:moving-feet3} 
d_{X_{vw}}(\hat{x}_v^i, x_w^i)\le k, 
\end{equation} 
where $\hat{x}_v^i\in \{x_v^i, \bar{x}_v^i\}$.

We next turn our attention to the center $z_w$ of $Y_w$. 
Except for the generic subcase (d) above, the tripod $Y_w$ is degenerate 
and  $z_w$ is one of its extremities, i.e. equals to one of the points $\tilde{x}_w^j$. Hence, apart from the generic case, 
as we observed while discussing nongeneric cases, 
\begin{equation}\label{eq:moving-center-0} 
d_{X_{vw}}(z_w, \bar{x}^j_v)\le k. 
\end{equation}
Note that the point $z_w$ is within uniformly bounded (in terms of $k$) distance from $z_v$ in subcase (2c) but 
might be quite far from $z_v$ in the subcase (b). The next lemma allows us to control the 
position of $z_w$ in the generic subcase (d).

\begin{lemma}\label{lem:moving-center} 
Suppose that we are in the generic subcase (d). Then 
$$
d_{X_{vw}}(z_v, z_w)\le C_{\ref{lem:moving-center}}(k).  
$$
\end{lemma}
\proof We define the geodesic triangles (in $X_v, X_w$) $\hat\Delta_v:= \Delta \hat{x}_v^1 \hat{x}_v^2 \hat{x}_v^3$. 
$\Delta_w:= \Delta {x}_w^1 {x}_w^2 {x}_w^3$ and corresponding geodesic triangles in $X_{vw}$:
$$
\hat\Delta^*_v:= \Delta_{X_{vw}} \hat{x}_v^1 \hat{x}_v^2 \hat{x}_v^3,  \quad 
\Delta^*_w:= \Delta_{X_{vw}} {x}_w^1 {x}_w^2 {x}_w^3.$$
Then the points $z_v, z_w$ are, respectively,  a $3\delta_0$-center of $\hat\Delta_v$ and $\delta_0$-center of $\Delta_w$. 
Since the sides of $\hat\Delta_v$ are $L'_0$-quasigeodesics in $X_{vw}$, stability of quasigeodesics implies that 
 $z_v$ is within distance $3\delta_0 + D_{_{\ref{stab-qg}}}(\delta'_0,L'_0)$ from all three sides of 
$\hat\Delta^*_v$. Since the respective end-points of the geodesic sides of the triangles 
$\hat\Delta^*_v, \Delta^*_w$ are within distance $k$ in $X$, it follows that  
$z_v$ is within distance $3\delta_0 + D_{\ref{stab-qg}}(\delta'_0,L'_0)+ \delta'_0+k$ from all the sides of  
$\Delta^*_w$, i.e. $z_v$ is a $3\delta_0 + D_{\ref{stab-qg}}(\delta'_0,L'_0)+ \delta'_0+k$-center of 
$\Delta^*_w$. 

 Similarly, the point $z_w$ is a $\delta_0+ D_{\ref{stab-qg}}(\delta'_0,L'_0)$-center of the triangle 
 $\Delta^*_w$ in $X_{vw}$. Thus,  by Lemma \ref{lem:centers},
$$
d_{X_{vw}}(z_w, z_v)\le D_{\ref{lem:centers}}(\delta_0,3\delta_0 + D_{\ref{stab-qg}}(\delta'_0,L'_0)+ \delta'_0+k). 
$$
Setting $C_{\ref{lem:moving-center}}(k):= D_{\ref{lem:centers}}(\delta_0,3\delta_0 + D_{\ref{stab-qg}}(\delta'_0,L'_0)+ \delta'_0+k)$ 
concludes the proof. \qed 

\begin{cor}\label{cor:moving-center} 
For every edge $e=[v,w]$ (pointing away from $u$), we have 
$$
\Hd_{X_{vw}}([z_v \bar{x}^i_v]_{X_v}, L^i_w)\le D_{\ref{cor:moving-center}}(k)= 
D_{\ref{lem:sub-close}}(\delta'_0,L'_0,C_{\ref{lem:moving-center}}(k)),$$
unless $L^i_w=\emptyset$. In any case, 
$$
\Hd_{X_{vw}}(\bar{Y}_v, Y_w)\le D_{\ref{cor:moving-center}}(k). 
$$
\end{cor}
\proof In the generic subcase (d) the first claim is an immediate application of Lemma \ref{lem:sub-close} and Lemma 
\ref{lem:moving-center}. In other subcases, both tripods $\bar{Y}_v, Y_w$ are degenerate, equal to geodesic segments whose respective end-points are within distance $k$ in $X_{vw}$, e.g. $Y_w= L^i_w$. 
Therefore, we similarly conclude that
 $$
 \Hd_{X_{vw}}([z_v \bar{x}^i_v]_{X_v}, L^i_w) \le D_{\ref{lem:sub-close}}(\delta'_0,L'_0, k). 
 $$
Observing that 
$$
k\le C_{\ref{lem:moving-center}}(k),$$
we obtain the first claim in other subcases as well. 
The second claim is an immediate corollary of the first one. \qed  

\medskip
We are now ready to verify that $\Y$ satisfies axioms of a $(K,D,E)$-semicontinuous family and that each $\L^i$ 
is the vertex set of a $(\kappa_1,D_1,E_1)$-ladder.  The constants $\kappa_1$ and $E_1$ will be computed in the end of the proof of the proposition.   

In line with the proof of Lemma \ref{lem:E-ladder-structure}, for the edges 
$e=[v,w]$ (oriented away from $u$) 
we define tripods $Y_e\subset X_e$  as
$$
T_{z_e}(x_e^1 x_e^2 x_e^3)
$$
where $x_e^i$ is a nearest-point projection of $x_w^i$ to $X_e$ in $X_{vw}$, while $z_e$ is a $\delta_0$-center of the geodesic triangle $\Delta  x_e^1 x_e^2 x_e^3 \subset X_e$. Thus, by the construction, each $Y_v, Y_e$ is a $\delta_0$-quasiconvex subset 
of the respective vertex and edge-space of $\X$.

According to Corollary \ref{cor:moving-center}, every point 
$x\in Y_w$ is within distance  
\begin{equation}\label{eq:k1} 
\kappa_1:=D_{\ref{cor:moving-center}}(k)
\end{equation}  
from a point $y\in \bar{Y}_v\subset Y_v$. Thus, every 
point of $\YY$ is connected to $X_u$ by a $k_1$-qi section. Similarly,   we get the inequality $\Hd_{X_{vw}}(Y_w, Y_e)\le \kappa_1$.  
Let's verify the inequality \eqref{eq:E-in} (see Definition \ref{defn:scfamily}) 
for a suitable value of the parameter $E$, i.e. get a uniform bound on the Hausdorff distance between $Y_w$ and the projection 
of $Y_v$ to $X_w$. First of all, since
$$
d_{X_{vw}}(Y_v, X_w)\le k= \max(4\la'_0+8\delta'_0+ 5\delta_0, k),  
$$
using Lemma \ref{lemma0-flow-space} we obtain 
\begin{align*}
\Hd_{X_{vw}}(P_{Y_v}(X_w), P_{X_w}(Y_v))\le R_{\ref{lemma0-flow-space}}(k,\la,\delta'_0)=\\ 
 2\la'_0+3\delta'_0 +k
\end{align*}
Thus, we need to estimate the Hausdorff distance between $P_{Y_v}(X_w)$ and $Y_w$. According to 
\eqref{eq:proj-to-tripod1}, 
$$
\Hd_{X_{vw}}(\bar{Y}_v, P_{Y_v}(X_w))\le \la. 
$$
Combining these inequalities with Corollary \ref{cor:moving-center}, we get:
\begin{equation}\label{eqE1} 
\Hd_{X_{vw}}(P_{X_w}(Y_v), Y_w)\le E_1:= D_{\ref{cor:moving-center}}(k) + \la + (2\la'_0+3\delta'_0 +k).  
\end{equation}
For Axiom 4 of a semicontinuous family we observe that, by the construction, if $Y_w=\emptyset$ then 
$Y_v$ and $X_w$ are $D_1$-cobounded. We conclude:

\begin{lemma}\label{lem:1tripod} 
$\Y$ is a $(\kappa_1,D_1,E_1)$-semicontinuous family containing sections $\ga^1, \ga^2, \ga^3$  
and $\YY\subset N^{fib}_{5\delta_0}({\mathcal Fl}_k(X_u))$. 
\end{lemma}


\medskip 
Next, consider the  families of intervals $\LL^i, i=1,2,3$ ($\LL^i$ is the union of geodesic intervals $L^i_v$, $v\in V(S)$). 
We will be verifying the conditions of Lemma \ref{lem:E-ladder-structure} for each $\LL^i$. The easiest thing to check is the first part of condition (a2) of the lemma, dealing with Property 4 of a semicontinuous family of spaces. Namely, by the definition of $L^i_w$, it is empty only when $L^i_v$ and $X_w$ are $D_2=D_1$-cobounded in $X_{vw}$, 
see Remarks \ref{rem:cbd2} and \ref{rem:cbd3}.

Next, consider the condition (a3) of the lemma. Suppose $L^i_w$ is nonempty, equals the oriented segment $[x_w^i z_w]_{X_{w}}$. 
There are several cases to consider, for instance, suppose we are in the generic subcase (2d). Then there exists a point 
$\hat{x}_v^i\in [x_v^i z_v]_{X_v}= L^i_v$ within distance $k$ from $x^i_w$, while 
according to Lemma \ref{lem:moving-center} 
$$
d_{X_{vw}}(z_v, z_w)\le C_{\ref{lem:moving-center}}(k).  
$$
Of course, in this case, in the oriented interval $L^i_v$, we have
$$
x_v^i\le\hat{x}_v^i\le  z_v\le z_v. 
$$
In all nongeneric cases, there are points  $\bar{x}^{i\pm 1}_v, \hat{x}^i_v\in L^i_v$ within distance $k$ from 
$z_w, x^i_w$ (subcase (2c)) or points $\hat{x}_v, \hat{y}_v\in L^i_v$ (or $L_v$) within distance $k$ from $x_w, y_w$ 
(subcase (2c) or case 3), and these points appear in the oriented interval $L^i_v$ (or $L_v$) in the correct order. 

This verifies condition (a3) of Lemma \ref{lem:E-ladder-structure} with the constant 
\begin{equation}\label{eq:k2}
k_2=\max(k, C_{\ref{lem:moving-center}}(k))=C_{\ref{lem:moving-center}}(k)  
\end{equation} 
playing the role of $K$ in Lemma \ref{lem:E-ladder-structure}.

\medskip 
Lastly, we analyze the projection of $L_v^i$ to $X_w$. Similarly to the projection of $Y_v$ to $X_w$, we have: 
$$
d_{X_{vw}}(L_v^i, X_w)\le k  
$$
and, thus,
\begin{equation}\label{LX}
\begin{aligned} 
\Hd_{X_{vw}}(P_{L_v^i}(X_w), P_{X_w}(L_v^i))\le R_{\ref{lemma0-flow-space}}(k,\la,\delta)=\\ 
 2\la'_0+3\delta'_0 +k. 
\end{aligned}
\end{equation} 
In other words,  the projection of $L_v^i$ to $X_w$  is uniformly Hausdorff-close to 
the projection of $X_w$ to $L_v^i$. Therefore, we analyze the latter using 
the arguments from the proof of Lemma \ref{lem:leg-proj}. We have four projections:
$$
P_{X_{vw},Y_v}, P_{X_{vw},L^i_v},  P_{Y_v,L^i_v}, P'_{Y_v,L^i_v}
$$
where the first three are nearest-point projections with respect to the metric of $X_{vw}$, while the last one 
 is the intrinsic nearest-point projection with respect to the metric of $Y_v$. We have, of course, 
 $$
 P'_{Y_v,L^i_v} \circ P_{X_{vw},Y_v}(X_w)= P_{X_{vw},Y_v}(X_w) \cap [\bar{x}^i_v z_v]_{X_v} 
  \subset P'_{Y_v,L^i_v}(\bar{Y}_v)= [\bar{x}^i_v z_v]_{X_v}. 
 $$
 According to Lemma \ref{lem:proj-to-tripod1}, 
 \begin{equation}\label{eq:lem:proj-to-tripod1}
 \Hd_{X_{vw}}(P_{X_{vw},Y_v}(X_w) \cap [\bar{x}^i_v z_v]_{X_v}, [\bar{x}^i_v z_v]_{X_v}) \le \la. 
 \end{equation}
As in the proof of Lemma \ref{lem:leg-proj},   
$$
d_{X_{vw}}(P_{X_{vw},L^i_v}, P_{Y_v,L^i_v}\circ P_{X_{vw},Y_v})\le 
C_{\ref{cor:projection-2}}(\delta'_0,\la'_0), 
$$ 
while for $x\in Y_v$,
$$
d_{X_{vw}}(P_{Y_v,L^i_v}(x), P'_{Y_v,L^i_v}(x))\le 2(\la'_0+2\delta'_0 + D_{\ref{stab-qg}}(\delta'_0, L'_0)). 
$$
Combining the two inequalities, we obtain that for each $q\in X_w$,
\begin{align*}
d_{X_{vw}} (P'_{Y_v,L^i_v} \circ P_{X_{vw},Y_v}(q), P_{X_{vw},L^i_v}(q))\le\\
C_{\ref{lem:leg-proj}}= C_{\ref{cor:projection-2}}(\delta'_0,\la'_0)+ 
2(\la'_0+2\delta'_0 + D_{\ref{stab-qg}}(\delta'_0, L'_0)) 
\end{align*}
Thus, taking into account \eqref{eq:lem:proj-to-tripod1}, 
$$
\Hd_{X_{vw}} ( P_{X_{vw},L^i_v}(X_w), [\bar{x}^i_v z_v]_{X_v}) \le  C_{\ref{lem:leg-proj}} +\la. 
$$
Combined with the inequality \eqref{LX}, we get: 
$$
\Hd_{X_{vw}} (P_{X_w}(L_v^i), [\bar{x}^i_v z_v]_{X_v}) \le  R_{\ref{lemma0-flow-space}}(k,\la,\delta'_0)+ 
C_{\ref{lem:leg-proj}} +\la. 
$$
Recall that by Corollary \ref{cor:moving-center}
$$
\Hd_{X_{vw}}([z_v \bar{x}^i_v]_{X_v}, L^i_w)\le D_{\ref{cor:moving-center}}(k).$$
Thus, 
\begin{equation}\label{eq:E2est}
\Hd_{X_{vw}} (P_{X_w}(L_v^i), L^i_w) \le E_2:=D_{\ref{cor:moving-center}}(k) + 
R_{\ref{lemma0-flow-space}}(k,\la,\delta'_0)+ 
C_{\ref{lem:leg-proj}} +\la. 
\end{equation} 

This concludes verification of the conditions  of Lemma \ref{lem:E-ladder-structure} and we obtain:

\begin{lemma}\label{lem:3-ladders}
For $k_2$ given by the equation \eqref{eq:k2}, $\kappa_2=k_2'=K'_{\ref{lem:E-ladder-structure}}(k_2)$, $E_2$ as in \eqref{eq:E2est},  and 
 $\la=C_{\ref{lem:proj-to-tripod1}}(\delta'_0,\la'_0,L'_0)$,  each $\LL^i$ defined earlier is the vertex-set of a 
$(\kappa_2,D_2,E_2)$-ladder $\L^i$ in $\X$. Each 
ladder $\L^i$ contains the section $\ga^i$. 
\end{lemma}

This concludes the proof of part (i) of the proposition.  We now prove part (ii). The goal, of course, is to verify the conditions 
of Lemma \ref{lem:E-ladder-structure} for the family of segments $L_v^{ij}, v\in V(S)$, $j=i+1$. We define the monotonic 
map
$$
f_{ij}: L_v^i\cup L_v^j \to L_v^{ij}= [x^i_v x^j_v]_{X_v},
$$
using Corollary \ref{lem:proj}. According to Lemma \ref{lem:3-ladders}, the map moves each point distance $\le 4\delta_0$. The points 
$$
x_v':= f_{ij}(\hat{x}_v^i), y_v':= f_{ij}(\hat{x}_v^j). 
$$
satisfy 
$$
x^i_v \le x_v'\le y'_v\le x^j_v
$$
in the oriented interval $L_v^{ij}$. Since 
$$
\max( d_{X_{vw}}(\hat{x}_v^i, x_w^i),  d_{X_{vw}}(\hat{x}_v^j, x_w^j)   )\le k_2,
$$
we obtain
\begin{equation}\label{eq:kappa}
\max( d_{X_{vw}} (x_v', x_w^i), d_{X_{vw}} (y_v', x_w^j))\le k_3:=k_2+ 4\delta_0. 
\end{equation}

Next: Since $L_v^{ij}\subset N^{fib}_{2\delta_0}(Y_v))$,
$$
\diam(P_{X_w}(L_v^{ij}))\le L'_1\cdot 3\delta_0 + \diam(P_{X_w}(Y_v))\le 3\delta_0+ D_1. 
$$
Since
\begin{equation}\label{eq:Lij}
\Hd_{X_{v}}(L_v^i\cup L_v^j, L_v^{ij})\le 2\delta_0,
\end{equation}
$$
d_{X_{vw}}(P_{L_v^i \cup L_v^j}, P_{L_v^{ij}})\le 2\delta_0 + 2(\la'_0+2\delta'_0). 
$$
Therefore, 
\begin{align}
\diam(P_{L_v^{ij}}(X_w))\le \diam( P_{L_v^i \cup L_v^j} (X_w))+ 2\delta_0 + 2(\la'_0+2\delta'_0) \le \\
\diam( P_{Y_v} (X_w))+ 2\delta_0 + 2(\la'_0+2\delta'_0) \le D_1 + 2\delta_0 + 2(\la'_0+2\delta'_0). 
\end{align}
Therefore, for 
$$
D_3:= D_1+ \max( 3\delta_0, 2\delta_0 + 2(\la'_0+2\delta'_0)), 
$$
for every boundary edge $e=[v,w]$ of $S$, the subsets $L_v^{ij}, X_w\subset X_{vw}$ are $D_3$-cobounded. 

Lastly, we estimate the Hausdorff distance between the projection of $L_v^{ij}$ to $X_w$ and $L_w^{ij}$ for every edge 
$e=[v,w]\in S$. We again use the inequality \eqref{eq:Lij} and the coarse Lispchitz property of the projection $P_{X_{vw},X_w}$:
$$
\Hd(P_{X_w}(L_v^i\cup L_v^j) , P_{X_w}(L_v^{ij}))\le L'_1\cdot 3\delta_0. 
$$
Therefore, the inequality \eqref{eq:E2est}, implies that 
$$
\Hd(P_{X_w}(L_v^{ij}), L_w^{ij})\le E_3:= L'_1\cdot 5\delta_0 +E_2.  
$$
It follows (in view of Lemma \ref{lem:E-ladder-structure}) 
that $\L^{ij}$ is the union of vertex-sets of a $(\kappa_3, D_3, E_3)$-ladder, $\kappa_3=k'_3$. 
Taking 
$$
K:= \max(\kappa_1, \kappa_2, \kappa_3), D:= \max(D_1, D_3), E:= \max(E_1, E_2, E_3),
$$ 
concludes the proof of the proposition. \qed

\begin{cor}\label{cor:existence-of-ladders} 
For $k, K, D, E$ as in the proposition, 
any two points $x, y\in \Fl_k(X_u)$ belong to a $(K,D,E)$-ladder $\L$ 
centered at $u$ and  
contained in the fiberwise 
$4\delta_0$-neighborhood of  $\Fl_k(X_u)$. Furthermore, if we are given two $k$-leaves $\ga_x, \ga_y$ in $Fl_k(X_u)$ connecting $x, y$ to $X_u$, the ladder $\L$  can be chosen to satisfy:
$$
\ga_x\subset bot(\L), \ga_y\subset top(\L). 
$$
\end{cor}

\section{Projection of ladders}\label{sec:proj-ladders}

In this section we discuss an important {\em projection} procedure which converts a pair of ladders $\L^i= (\pi: L^i\to S_i), i=1,2$, with the common center $u$  into a pair of subladders $\bar\L^i\subset \L^i, i=1, 2$ (with the same center $u$, but possibly a different set of qi sections $\bar{\Sigma}_\bullet$),  within uniformly bounded (fiberwise) Hausdorff distance from each other. This construction will be used in Section \ref{sec:Hyperbolicity of flow-spaces} to show the coarse independence of a combing path $c(x,y)$ in a ladder $L_{x,y}$ containing the given points $x, y$ on the choice of a ladder $L_{x,y}$. 

The intersection 
$\pi(\L^1)\cap \pi(\L^2)= S_1\cap S_2$ is a subtree $S\subset T$. 
This subtree will contain (but, in general will be different from) the tree $\bar{S}$ which is the common base of the ladders $\pi: \bar\L^i\to \bar{S}$. 
For each $v\in V(S)$ we let
$$
\bar{L}^1_v:= \bar{P}_{X_v,L^1_v}(L^2_v)\subset L^1_v, \bar{L}^2_v:= 
\bar{P}_{X_v,L^2_v}(L^1_v)\subset L^2_v,
$$
see Definition \ref{def:modified projection} for the definition of the modified
 fiberwise   projection $\bar{P}$. In the definition of ladders $\bar\L^i$ below, 
 the segment  $\bar{L}^i_v$ will equal the fiber $\bar{L}^i_v$ of $\bar\L^i$, unless $v\notin V(\bar{S})$.

By Lemma \ref{cobdd-cor} and Corollary \ref{cor:cob-char} we have the dichotomy:

i. Either the pair of geodesic segments $L^1_v, L^2_v\subset X_v$ is $7\delta_0$-separated (i.e. $d_{X_v}(L^1_v, L^2_v)> 7\delta_0$), 
in which case this pair is  $D_{\ref{cobdd-cor}}(\delta, \delta)= 9\delta_0$-cobounded.

ii. Or $d_{X_v}(L^1_v, L^2_v)\le 7\delta_0$ and
$$
\Hd_{X_v}(P_{L^1_v}(L^2_v), P_{L^2_v}(L^1_v)) \le 12\delta_0.$$ 
According to Remark \ref{rem:barP}, in this case 
$$
\Hd_{X_v}(P_{L^1_v}(L^2_v), \bar{L}^1_v)\le 4\delta_0$$
and
$$ 
\Hd_{X_v}(P_{L^2_v}(L^1_v), \bar{L}^2_v)\le 4\delta_0. 
$$
Combining the inequalities, we get that in the second case, 
\begin{equation}\label{eq:5}
\Hd_{X_v}(\bar{L}^1_v, \bar{L}^2_v)\le 20\delta_0. 
\end{equation}

Accordingly, for each vertex $v\in S$ such that the pair 
of geodesic segments $L^1_v, L^2_v\subset X_v$ is $7\delta_0$-separated, we 
remove from $V(S)$ all the vertices (and edges) $w\ne v$ such that $v$ is 
between $u$ and $w$; we let $\bar{S}$ denote the subtree of $S$ spanned by the remaining set of vertices of $S$. Note that if $v=u$,  $\bar{S}=\{u\}$.

We define 
$$
\bar{\LL}^i:= \bigcup_{v\in V(\bar{S})} \bar{L}^i_v. 
$$
Thus, apart from the boundary vertices $v$ of $\bar{S}$, the intervals $\bar{L}^1_v, \bar{L}^2_v$ are fiberwise Hausdorff $20\delta_0$-close, while for {\em some} boundary vertices $v$ of $\bar{S}$, both  $\bar{L}^1_v, \bar{L}^2_v$ have length $\le 9\delta_0$. We will prove (Claim \ref{claim:vwE}) that even for the boundary vertices of $\bar{S}$, the  
 Hausdorff distance between the intervals $\bar{L}^1_v, \bar{L}^2_v$ is also uniformly bounded. 

\begin{rem}
By the construction, both $\bar{\LL}^i$ contain the intersection $\LL^1\cap \LL^2$. 
\end{rem}

Thus, we obtain subsets $\bar{\LL}^i\subset \XX$ 
which are unions of geodesic segments in vertex-spaces 
$X_v, v\in V(\bar{S})$. Our next goal it to prove that these subsets are unions of vertex-spaces of ladders in $\X$. 

\medskip 
Recall that $L'_1$ is a coarse Lipschitz constant for the composition of the inclusion map $X_v\to X_{vw}$ with 
the nearest-point projection $P_{X_{vw},X_{w}}: X_{vw}\to X_w$, cf.  Notation \ref{not:edge-space-constants}. 
Set
$$
 \eps= \eps_{\ref{prop:proj-of-pairs}}(E)=
 E+ (2E+ 21\delta_0 L'_1) + L_{\ref{lip-proj}}(\delta_0,\delta_0)(2E+ 21\delta_0 L'_1 +1),  
$$
$$
 k:=2(2K+ \eps+ 20\delta_0) +K + \eps+ 20\delta_0,
$$
$$
\tilde{K}= \tilde{K}_{\ref{prop:proj-of-pairs}}(K,E)= K'_{\ref{lem:E-ladder-structure}}(k),  
$$
$$
\tilde{E}:=  \tilde{E}_{\ref{prop:proj-of-pairs}}(K,E)= 2k, 
$$
$$
\tilde{D}:= \tilde{D}_{\ref{prop:proj-of-pairs}}(D)= 
\max(10\delta_0\cdot L'_1, D, C_{\ref{cor:cob-char}}(\la'_0,\delta'_0, 21\delta_0 L'_1 )).  
$$

\begin{prop}
[Projections of pairs of ladders] \label{prop:proj-of-pairs}
If ${\L}^1, {\L}^2$ are $(K,D,E)$-ladders centered at $u$, then there exist 
$(\tilde{K}, \tilde{D}, \tilde{E})$-ladders $\bar\L^i=(\pi: \bar{L}^i\to \bar{S}), i=1,2$,  centered at $u$, where 
$\tilde{K}=\tilde{K}_{\ref{prop:proj-of-pairs}}(K,E)$, $\tilde{D}= \tilde{D}_{\ref{prop:proj-of-pairs}}(D)$, 
$\tilde{E}= \tilde{E}_{\ref{prop:proj-of-pairs}}(K,E)$, such that 
$\bar{\LL}^i = \bar{L}^i\cap \XX$ is the union of vertex-spaces of $\bar\L^i$.  
\end{prop}
\proof As with tripods of ladders, we will prove the proposition by verifying the conditions of Lemma \ref{lem:E-ladder-structure}. 
Below,  $e=[v,w]$ is an edge of $T$ oriented away from $u$, with $v\in V(\bar{S})$. 

\medskip
1. We first check part 1 of the condition (a2) in Lemma \ref{lem:E-ladder-structure}, i.e. that for every boundary edge $e$ of 
$\bar{S}$,   $e=[v,w]\notin E(\bar{S})$, $v\in V(\bar{S})$,  
both pairs of subsets $\bar{L}^i_v, X_{w}\subset X_{vw}$ are $\tilde{D}$-cobounded, $i=1,2$. 

(i) It can happen that $e=[v,w]\notin E(\bar{S})$ because the pair $L^1_v, L^2_v\subset X_v$ was $7\delta_0$-separated, i.e. both  
$\bar{L}^1_v, \bar{L}^2_v$ have length $\le 9\delta_0$. 
Since the projection $P_{X_{vw},X_w}: X_{vw}\to X_w$ is $L'_1$-coarse Lipschitz,  
diameters of the projections $P_{X_{vw},X_w}(\bar{L}^i_v)\subset X_w$ are at most 
$10\delta_0\cdot L'_1\le \tilde{D}$. (Recall that $\delta_0\ge 1$.) 

(ii) If  the pair $L^1_v, L^2_v\subset X_v$ was not $7\delta_0$-separated then the Hausdorff distance in $X_v$ between $\bar{L}^1_v, \bar{L}^2_v$ is $\le 20\delta_0$. After swapping the labels of 1 and  2 we may assume that 
there is an edge $e=[v,w]\notin E(S_1)$. Then, because $\L^1$ was a $(K,D,E)$-ladder, the 
pair $L^1_v, X_w$ is $D$-cobounded in $X_{vw}$. The same, of course applies to the pair 
$\bar{L}^1_v, X_w$, since $\bar{L}^1_v\subset L^1_v$. We need to get a coboundedness estimate for the pair 
$\bar{L}^2_v, X_w$. Since 
\begin{equation} 
\Hd_{X_v}(\bar{\LL}^1_v, \bar{\LL}^2_v)\le 20\delta_0 
\end{equation}
and the projection $P_{X_{vw},X_w}: X_v\to X_w$ is coarse $L_1'$-Lipschitz, 
the diameter the projection of $\bar{L}_v^2$ to $X_w$ is at most 
$$
L'_1(20\delta_0+1)\le 21\delta_0 L'_1 \le \tilde{D}. 
$$
The estimate 
$$
\diam_{X_{vw}}(P_{X_{vw},\bar{L}^2_v}(X_w)) \le  C_{\ref{cor:cob-char}}(\la'_0,\delta'_0, 21\delta_0 L'_1 )
$$
follows from Corollary \ref{cor:cob-char}.

This verifies part 1 of the condition (a2) in Lemma \ref{lem:E-ladder-structure}. 

\medskip 
2.  We assume now that $e=[v,w]$ is an edge of $\bar{S}$, in 
particular, $\bar{L}^1_v, \bar{L}^2_v$ are $20\delta_0$-Hausdorff close in $X_v$. 
Since $\L^1, \L^2$ satisfy the $(K,D,E)$-ladder axioms, for $P= P_{X_{vw},X_w}$,
$P(L^i_v)$ and $L^i_w$ are $E$-Hausdorff close in $X_{vw}$, $i=1, 2$. Our goal is to estimate the Hausdorff distance between $L^i_w$ and the projection of $L^i_v$ to $X_w$ ($i=1,2$). In this part of the proof we will get only half of the estimate, we will get the other half in Part 4 of the proof. 

Pick $x\in \bar{L}^1_v$. Then there exists $y\in \bar{L}^2_v$ such that $d_{X_v}(x,y)\le 20\delta_0$; 
we also have 
$$
d_{X_w}(P(x), L^1_w)\le E, \quad d_{X_w}(P(y), L^2_w)\le E.
$$ 

Thus, there exist $x'\in L^1_w$ and $y'\in L^2_w$ such that 
\begin{equation}\label{eq:1}
d_{X_w}(P(x), x')\le E,\quad d_{X_w}d(P(y), y')\le E,
\end{equation}
which in turn implies the inequality 
\begin{equation}\label{eq:2}
d_{X_w}(x', y')\le 2E+ 21\delta_0 L'_1.
\end{equation}

We next estimate  the distance $d_{X_w}(y', \bar{L}^2_w)$. 
Since the projection $P': X_w\to L^2_w$ is $L_{\ref{lip-proj}}(\delta_0,\delta_0)$-coarse Lipschitz, 
we have 
$$
d_{X_w}(P'(x'), y')\le  L_{\ref{lip-proj}}(\delta_0,\delta_0) (d_{X_w}(x',y') +1)  \le 
 L_{\ref{lip-proj}}(\delta_0,\delta_0)(2E+ 21\delta_0 L'_1 +1).   
$$

Since $x'\in L^1_w$,  $P'(x')\in \bar{L}^2_w$, and 
we obtain:
\begin{equation}\label{eq:3}
d_{X_w}(y',  \bar{L}^2_w)\le L_{\ref{lip-proj}}(\delta_0,\delta_0)(2E+ 21\delta_0 L'_1 +1). 
\end{equation}

Switching the roles of $x$ and $y$ we similarly obtain:
\begin{equation}\label{eq:4}
d_{X_w}(x',  \bar{L}^1_w)\le L_{\ref{lip-proj}}(\delta_0,\delta_0)(2E+ 21\delta_0 L'_1 +1).  
\end{equation}

Combining the equations \eqref{eq:1} and \eqref{eq:4} we obtain
$$
d_{X_w}(P(x),  \bar{L}^1_w)\le E+ L_{\ref{lip-proj}}(\delta_0,\delta_0)(2E+ 21\delta_0 L'_1 +1) 
$$
and, similarly, combining  \eqref{eq:1} and \eqref{eq:3} we obtain 
$$
d_{X_w}(P(y),  \bar{L}^2_w)\le E+ L_{\ref{lip-proj}}(\delta_0,\delta_0)(2E+ 21\delta_0 L'_1 +1).  
$$
At the same time, combing the equations \eqref{eq:1}, \eqref{eq:2} and \eqref{eq:3} we get:
$$
d_{X_w}(P(x),  \bar{L}^2_w)\le E+ (2E+ 21\delta_0 L'_1) + L_{\ref{lip-proj}}(\delta_0,\delta_0)(2E+ 21\delta_0 L'_1 +1) \le \tilde{E}. 
$$
Similarly,
$$
d_{X_w}(P(y),  \bar{L}^1_w)\le 
\eps= E+ (2E+ 21\delta_0 L'_1) + L_{\ref{lip-proj}}(\delta_0,\delta_0)(2E+ 21\delta_0 L'_1 +1). 
$$
Thus, we proved:
$$
P(\bar{L}^1_v)\cup P(\bar{L}^2_v) \subset N^e_{\eps}( \bar{L}^1_w) \cap N^e_{\eps}( \bar{L}^2_w),
$$
which gives us half of the estimate on $\Hd_{X_{vw}}( P(\bar{L}^i_v), \bar{L}^i_w)$, but also gives an upper bound on the minimal distance between $\bar{L}^1_w$ and $\bar{L}^2_w$, something which was not a priori clear from the construction. We will derive the other half of the estimate on $\Hd_{X_{vw}}( P(\bar{L}^i_v), \bar{L}^i_w)$ in Part 4 of the proof.  Before proceeding with Part 3 of the proof we establish:

\begin{claim}\label{claim:vwE}
For every vertex $w\in V(\bar{S})$
$$
\Hd_{X_{vw}}(\bar{L}^1_w, \bar{L}^2_w) \le C_{\ref{claim:vwE}}(E)= \eps+ 20\delta_0. 
$$
\end{claim}
\proof There are two cases to consider according to the definition of the projections $\bar\LL^i$:

(i) Suppose that the pair of geodesic segments $L^1_w, L^2_w$ is $7\delta_0$-separated, hence, 
each segment $\bar{L}^i_w, i=1,2$ has length $\le 9\delta_0$. As we observed in the end of the Part 2 of the proof, 
there exists $z\in X_w$ which lies in the intersection 
$$
N^e_{\eps}( \bar{L}^1_w) \cap N^e_{\eps}( \bar{L}^2_w). 
$$
Therefore, 
$$
\Hd_{X_{vw}}(\bar{L}^1_w, \bar{L}^2_w) \le \eps+ 9\delta_0. 
$$

(ii) Suppose that the pair of geodesic segments $L^1_w, L^2_w$ is not $9\delta_0$-cobounded. 
According to \eqref{eq:5}, 
$$
\Hd_{X_v}(\bar{\LL}^1_v, \bar{\LL}^2_v)\le 20\delta_0. \qedhere 
$$

\medskip 
We now return to the proof of the proposition. 

\medskip 
3. We next verify the condition (a3)  in Lemma \ref{lem:E-ladder-structure}.

Since $\L^i, i=1,2$ were $K$-ladders,  each $x_i\in \bar{L}^i_w$ is within distance $K$ (measured in $X_{vw}$) 
 from some point in $L^i_v$. Thus, setting $y_i:= P_{X_{vw},L^i_v}(x_i)$,  
$$
d_{X_{vw}}(x_i, y_i)\le 2K, i=1,2.$$ 
According to Claim \ref{claim:vwE}, given $x_1\in \bar{L}^1_w$ we can find $x_2\in \bar{L}^2_w$ such that
$$
d_{X_{vw}}(x_1, x_2)\le \eps+ 20\delta_0. 
$$
Thus,
$$
d_{X_{vw}}(y_1, y_2)\le 2K+ \eps+ 20\delta_0. 
$$
Let $\bar{y}_i$ denote a nearest-point projection of $y_i$ to $L^{3-i}_v, i=1,2$. Then the above inequality implies 
that
$$
d_{X_{vw}}(\bar{y}_i, y_i)\le 2(2K+ \eps+ 20\delta_0). 
$$
However, by the construction, $\bar{y}_i$ belongs to $\bar{L}_v^{3-i}$. Hence, by the above estimates: 
$$
d_{X_{vw}}(x_1, \bar{L}_v^{1}) \le d_{X_{vw}}(x_1,  \bar{y}_2) \le 
k:=2(2K+ \eps+ 20\delta_0) +K + \eps+ 20\delta_0.  
$$
Similarly, for each $x_2\in \bar{L}^2_w$ we have 
$$
d_{X_{vw}}(x_2, \bar{L}_v^{2}) \le k=2(2K+ \eps+ 20\delta_0) +K + \eps+ 20\delta_0.  
$$
This verifies the condition (a3)  in Lemma \ref{lem:E-ladder-structure}.

\medskip 
4. By Part 3, each $x\in \bar{L}_w^{i}$ is within distance $k$ from some $y\in \bar{L}_v^{i}$, $i=1,2$. Therefore, 
$$
d_{X_{vw}}(x, P_{X_{vw},X_w}(y))\le 2k. 
$$
In other words, 
$$
\bar{L}_w^{i}\subset N^e_{2k}(P_{X_{vw},X_w}(\bar{L}_v^i)). 
$$
Combining this with the estimate in the end of Part 2, we obtain:
$$
\Hd_{X_{vw}}( P(\bar{L}^i_v), \bar{L}^i_w) \le \max(2k, \eps)= 2k= \tilde{E}. 
$$
Since we defined $\tilde{E}$ to be $2k$, we are done with the proof of the proposition. \qed

\section{Hyperbolicity of tripods families} \label{sec:Hyperbolicity of tripods families}

\begin{prop}[Tripod families are hyperbolic]\label{tripod-bdle}
Suppose that $\X$ is a tree of hyperbolic spaces 
satisfying the uniform $\kappa_{\ref{defn:carpeted ladder}}(K)$-flaring condition. 
Then for each $(K,D,E)$-tripod family $\Y=(\pi: Y\to S)$ in $\X$,  the total space $Y$ is 
 $\delta_{\ref{tripod-bdle}}(K)$-hyperbolic. 
\end{prop}
\proof The total space $Y$ of $\Y$ is the union of total spaces $L^i$ of the ladders $\L^i$, $i=1,2,3$. The pairwise intersections of these ladders equal their triple intersection, namely, the $K$-qi section $\Xi$, which is, intrinsically, a tree. Thus, $Y$ has a structure of a tree of spaces with the vertex-spaces $L^i$, $i=1,2,3$, and $\Xi$.  
According to Theorem \ref{thm:ladders-are-hyperbolic} 
each ladder $L^i$ is $\delta_{\ref{thm:ladders-are-hyperbolic}}(K)$-hyperbolic (this is where we need the $\kappa_{\ref{defn:carpeted ladder}}(K)$-flaring condition), 
the space $Y$ is $\delta_{\ref{tripod-bdle}}(K)$-hyperbolic according to Corollary \ref{cor:finite-tree-hyp}.  \qed

\section{Hyperbolicity of flow-spaces} \label{sec:Hyperbolicity of flow-spaces}

In the following theorem, we fix $k\ge k_{\ref{prop:existence-of-tripod-ladders}}$, 
$K= K_{\ref{prop:existence-of-tripod-ladders}}(k)$.  

\begin{theorem}\label{flow of one vertex space}
Suppose that $k$, $K$ are as above,  $\X$ is a tree of hyperbolic spaces 
satisfying the uniform $\kappa_{\ref{defn:carpeted ladder}}(K)$-flaring condition. 
Then there is a function $\delta=\delta_{\ref{flow of one vertex space}}(k)$ such that
for each $u\in v(T)$, the flow space 
$Fl_k(X_u)$ is $\delta$-hyperbolic. 
\end{theorem}
\proof According to Corollary \ref{cor:existence-of-ladders}, whenever $k\ge k_{\ref{prop:existence-of-tripod-ladders}}$,  
$K= K_{\ref{prop:existence-of-tripod-ladders}}(k)$, $D= D_{\ref{prop:existence-of-tripod-ladders}}$, 
$E=E_{\ref{prop:existence-of-tripod-ladders}}$, for any two points $x, y\in  {\mathcal Fl}_k(X_u)$ 
there exists a $(K,D,E)$-ladder $\L=\L_{x,y}$ centered at $u$ and containing $x, y$, such that 
$\LL$ is contained in the fiberwise $5\delta_0$-neighborhood of ${\mathcal Fl}_k(X_u)$. 

Recall that the total space $L_{x,y}$ of the ladder $\L$ is  
$L_{\ref{cor:mitras-projection}}(K,D,E,\delta_0)$-qi embedded in $X$. 
Define $c(x,y)$ to be a projection to $Fl_k(X_u)$ of a geodesic in $L_{x,y}$ connecting $x$ to $y$. 
We note that the definition of $c(x,y)$ depends on the choice of $\L_{xy}$ which is far from canonical. Our first task is to prove that 
different choices lead to uniformly Hausdorff-close paths. 

\begin{prop}\label{prop:ladder-amalgam}
Let $\L^1=\L^1_{x,y}, \L^2=\L^2_{x,y}$ be  $(K,D,E)$-ladders containing   $x, y$. Then $L^1_{x,y}\cup L^2_{x,y}$ is contained in a 
$\delta_{\ref{prop:ladder-amalgam}}(K,D,E)$-hyperbolic subspace $Z$ in $X$. 
\end{prop}
\proof  We let $\bar{\L}^i\subset \L^i$ denote the $(\tilde{K},\tilde{D},\tilde{E})$-subladders obtained by the 
projection construction described in Section \ref{sec:proj-ladders}; see also Proposition \ref{prop:proj-of-pairs}. 
Note that $\tilde{K}\ge K, \tilde{E}\ge E, \tilde{D}\ge D$. Also, note that the subladders $\bar{\L}^i$ are nonempty since they both contain 
$x$ and $y$.

The subladders 
$\bar{L}^i$ have equal projection to $T$, which is a subtree $\bar{S}\subset S_1\cap S_2$. As usual, we extend these ladders over 
the rest of the tree $S_1\cup S_2$ by empty fibers. According to Claim \ref{claim:vwE}, the ladders $\bar{\L}^1, \bar{\L}^2$ 
are fiberwise $C_{\ref{claim:vwE}}(E)$-Hausdorff close. 
Therefore, for each $v\in V(\bar{S})$ the union $\bar{L}^1_v\cup \bar{L}^2_v$ is $C_{\ref{claim:vwE}}(E)+\delta_0$-quasiconvex 
in $X_v$. For each vertex $v\in V(S)$ we set 
$$
Z^0_v:= \hull_{\delta_0}  (\bar{L}^1_v\cup \bar{L}^2_v), 
$$ 
$$
Z^i_v:= L^i_v\cup Z^0_v, i=1,2, 
$$
and  $Z_v:= Z^1_v\cup Z^2_v$. Thus, $Z^j_v$ ($j=0,1,2$) and $Z_v$ are rectifiably connected 
$4\delta_0$-quasiconvex subsets of $X_v$ (see Lemma \ref{lem:qc-hull}). By \eqref{eq:qc-nbd}, 
$$
\Hd_{X_v}(Z^0_v, \bar{L}^1_v\cup \bar{L}^2_v)\le C_{\ref{claim:vwE}}(E)+2\delta_0
$$ 
and, hence,
\begin{equation}\label{eq:Z^0} 
\Hd_{X_v}(Z^0_v, \bar{L}^i_v)\le 2C_{\ref{claim:vwE}}(E)+2\delta_0, i=1,2. 
\end{equation}
Accordingly,
$$
\Hd_{X_v}(Z^i_v, L^i_v)\le 2C_{\ref{claim:vwE}}(E)+2\delta_0, i=1,2. 
$$
We repeat the same construction (and estimates) for all edges $e\in E(\bar{S})$. 
We, thus, obtain four subtrees of spaces $\ZZZ^j$ ($j=0,1,2$) and $\ZZZ$ in $\X$ whose vertex-spaces are, respectively  
$Z^j_v$ ($j=0,1,2$) and $Z_v$, $v\in V(S)$. The total space $Z$ of $\ZZZ$ contains both ladders $L^1, L^2$. 
We  equip $Z^j$'s and $Z$ with natural path-metrics $d_{Z^j}, d_Z$; the goal is to show that $Z$ is uniformly hyperbolic. 
The space $Z$ is the union of subsets $Z^1, Z^2$ whose intersection is $Z^0$. 
According to Corollary \ref{cor:mitras-projection},  the ladders $L^i, \bar{L}^i, i=1,2$ are 
$L_{\ref{cor:mitras-projection}}(\tilde{K},\tilde{D},\tilde{E},\delta_0)$-qi embedded in $X$. 
In view of \eqref{eq:Z^0}, the inclusion $Z^0\to Z$ is an  
$(L_{\ref{cor:mitras-projection}}(\tilde{K},\tilde{D},\tilde{E},\delta_0), 2C_{\ref{claim:vwE}}(E)+2\delta_0)$-qi embedding.

Thus, $Z=Z^1\cup Z^2$ satisfies the assumptions of Theorem \ref{thm:hyp-tree} and, therefore, is 
$\delta$-hyperbolic for some $\delta$ which depends only on $K, D, E$.  \qed

\begin{cor}\label{cor:ladder-amalgam}
Let $L^i_{x,y}, i=1,2$ be two $(K,D,E)$-ladders containing $x, y$ and $c^i(x,y)$ be projections to $Fl_k(X_u)$ of geodesics 
$[xy]^i\subset L_{x,y}^i$. Then
$$
\Hd_{X}(c^1(x,y), c^2(x,y))\le C_{\ref{cor:ladder-amalgam}}(K,D,E). 
$$ 
\end{cor}
\proof The paths $c^i(x,y)$ are uniformly close to the  
geodesic segments  
$[xy]^i\subset L_{x,y}^i$, $i=1,2$; hence, it suffices to bound the Hausdorff distance between these segments. 
Since both  $[xy]^i$ are contained in a 
$\delta_{\ref{cor:ladder-amalgam}}(K,D,E)$-hyperbolic space 
$(Z,d_Z)$ and are $L_{\ref{cor:mitras-projection}}({K},{D},{E},\delta_0)$-quasigeodesics in $Z$ with common end-points, 
by Lemma 
\ref{lem:sub-close}, we get: 
\begin{align}
\Hd_{X}([xy]^1, [xy]^2)\le  \Hd_{Z}([xy]^1, [xy]^2)\le  
\\
D_{\ref{lem:sub-close}}(\delta_{\ref{cor:ladder-amalgam}}(K,D,E),L_{\ref{cor:mitras-projection}}({K},{D},{E},0)). 
\end{align}
Corollary follows.  \qed

\medskip  
We now  check that the family of paths $c$ satisfies 
the conditions of the Corollary \ref{cor:bowditch}, characterizing hyperbolic spaces; this will conclude the proof of the theorem. 

\medskip 
{\bf Condition a1:} This follows from the fact that $c(x,y)$ is within uniformly bounded distance from 
a geodesic in $L_{x,y}$ and that $L_{x,y}$  is  $L_{\ref{cor:mitras-projection}}(K',D,E,\delta_0)$-qi embedded in $X$.

\medskip

{\bf Condition a2:}  Consider points $x, y, z\in {\mathcal Fl}_k(X_u)$. According to  Proposition \ref{prop:existence-of-tripod-ladders}, there exists a $(K,D,E)$-tripod family $\Y=(\pi: Y\to S)$ contained  the fiberwise $5\delta_0$-neighborhood of ${\mathcal Fl}_k(X_u)$ 
and containing the points $x, y, z$. Moreover, according to Part (ii) of  Proposition \ref{prop:existence-of-tripod-ladders}, 
there are $(K,D,E)$-ladders $\L_{x,y}, \L_{y,z}, \L_{z,x}$ contained in the fiberwise $\delta_0$-neighborhood of $\Y$ and containing the respective pairs of points $x, y$, etc. By Proposition \ref{tripod-bdle}, $Y$ is 
$\delta_{\ref{tripod-bdle}}(K, D, E)$-hyperbolic. The paths $c(x,y), c(y,z), c(z,x)$ are 
uniformly close to geodesics $[xy]_{L_{x,y}}, [yz]_{L_{y,z}}, [zx]_{L_{z,x}}$, which are $\kappa$-quasigeodesics in $Y$, where $\kappa$ depends only on $K, D$ and $E$.  Therefore, by the $\delta_{\ref{tripod-bdle}}(K, D, E)$-hyperbolicity of $Y$, 
for 
$$
\eps= 2D_{\ref{stab-qg}}(\delta_{\ref{tripod-bdle}}(K, D, E),\kappa)+ \delta_{\ref{tripod-bdle}}(K, D, E)
$$
we have 
$$
[xy]_{L_{x,y}}\subset N^{Y}_{\eps}([yz]_{L_{y,z}}\cup [zx]_{L_{z,x}})\subset N_{\eps}([yz]_{L_{y,z}}\cup [zx]_{L_{z,x}}), 
$$
where the first neighborhood is taken in $Y$ and the second is taken in $X$. Condition (a2) follows.  

\medskip
Lastly, by the construction, 
each $c(x,y)$ is uniformly close to an $L_{\ref{cor:mitras-projection}}(K',D,E,\delta_0)$-quasigeodesic in $X$. \qed 

\medskip

In the next corollary illustrates an application of Theorem \ref{flow of one vertex space} to proving hyperbolicity of various subspaces of $X$. As in Theorem \ref{flow of one vertex space}, we assume that $k\ge k_{\ref{prop:existence-of-tripod-ladders}}$, but set 
$K:= K_{\ref{prop:existence-of-tripod-ladders}}(k^\wedge)$. We will use similar arguments in Section \ref{sec:generalized hallways} to prove hyperbolicity of unions of pairs of flow-spaces.  
We refer to Definition \ref{defn:generalized flow-space} for the notion of a {\em generalized flow-space}, $Fl_k(\QQQ)$, used below.

\begin{cor}
[Hyperbolicity of generalized flow-spaces] \label{cor:hyp-of-gen-flows}
Suppose that $\X$ satisfies the uniform 
$\kappa_{\ref{defn:carpeted ladder}}(K)$-flaring condition. Then for 
every $k$-bundle  $\QQQ=(\pi: Q\to S)\subset \X$, the $k$-flow space 
$Fl_k(\QQQ)$ is $\delta_{\ref{cor:hyp-of-gen-flows}}(k)$-hyperbolic. 
\end{cor} 
\proof Pick $u\in V(S)$. Observe that  every $x\in Fl_k(\QQQ)$ is connected to $Fl_k(\QQQ)\cap X_u$ by a $k$-leaf $\ga_x$. 
This leaf is contained in $Fl_{k^\wedge}(Q_u)$ (see Proposition \ref{prop:flow-prop}(2)). Therefore, 
$$
Fl_k(\QQQ)\subset Fl_{k^\wedge}(Q_u). 
$$ 
In view of the uniform flaring condition of the corollary, Theorem  \ref{flow of one vertex space} applies to 
$Fl_{k^\wedge}(Q_u)$ and, hence, the latter is $\delta_{\ref{flow of one vertex space}}(k^\wedge)$-hyperbolic. 
According to Theorem \ref{thm:generalized retraction},  $Fl_k(\QQQ)$ is an   
$L_{\ref{thm:generalized retraction}}(k,k)$-coarse Lipschitz retract of $X$. Hence $Fl_k(\QQQ)$ is $\la(k)$-quasiconvex in 
$Fl_{k^\wedge}(Q_u)$ and, therefore (in view of hyperbolicity of the latter), 
is $\delta_{\ref{cor:hyp-of-gen-flows}}(k)$-hyperbolic. \qed

\chapter{Hyperbolicity of trees of spaces: Putting everything together}\label{sec:everything together}

In this chapter we finish the proof of Theorem \ref{thm:mainBF}, establishing hyperbolicity of 
trees of hyperbolic spaces, satisfying the uniform $K$-flaring condition for suitable values of $K$. The key  
is to show hyperbolicity of flow-spaces $Fl_K(X_J)$ for intervals $J\subset T$, Theorem \ref{thm:hyperbolicity over segments}. 
This is done in three steps:

{\bf Step 1.} Hyperbolicity of $Fl_K(X_J)$ for {\em special intervals} $J$ (Theorem \ref{thm:geod-flow-step1}). This is the hardest 
part of the chapter, we deal with it in Section \ref{sec:union-of-two}. An outline of this part of  proof is given  in the introduction to Section \ref{sec:union-of-two}. 

{\bf Step 2.} Hyperbolicity of $Fl_K(X_J)$, when $J$ is the union of three special subintervals (Proposition 
\ref{prop:geod-flow-step2}).   

{\bf Step 3.} Hyperbolicity of $Fl_K(X_J)$ for general intervals, which is done by subdividing $J$ as the union of subintervals $J_i$, 
each of which is a union of (at most) three special subintervals, and then using quasiconvex chain-amalgamation (Theorem \ref{thm:hyp-tree}). 

Once we are done with Theorem \ref{thm:hyperbolicity over segments}, applying quasiconvex amalgamation (Corollary \ref{cor:finite-tree-hyp})  one more time, in Proposition \ref{prop:hyperbolicity over tripds} we will prove that flow-spaces $Fl_K(X_S)$ are uniformly hyperbolic, whenever $S$ is a tripod in $T$.  We then conclude the proof of Theorem \ref{thm:mainBF} by appealing to 
Corollary \ref{cor:bowditch} one last time, by constructing a slim combing in $X$ via geodesics in flow-spaces of interval-spaces,
 see Section  \ref{sec:conclusion}.  

\section{Hyperbolicity of flow-spaces of special interval-spaces}
 \label{sec:union-of-two}

This section deals with Step 1 described in the introduction to the chapter. 
Recall that in Section \ref{sec:flow-incidence graph} we defined an interval $J=\llbracket u , v\rrbracket\subset T$ to be 
{\em special} (more precisely, $K$-{\em special}) if one of its end-points (say, $u$) has the property that $J\subset \pi(Fl_K(X_u))$.  
The main result of this section is Theorem \ref{thm:geod-flow-step1}, where we prove that the flow-space $Fl_K(X_J)$ of every special interval $J\subset T$, is uniformly hyperbolic with hyperbolicity constant depending only on $K$. The fact that this is true is not at all surprising since, assuming that $X$ is hyperbolic, a uniform neighborhood of two intersecting uniformly quasiconvex subsets is uniformly qi embedded in $X$ and, hence, is uniformly hyperbolic. However, at this stage we did not yet prove hyperbolicity of $X$ and, furthermore, we are interested in describing uniform quasigeodesics in $X$. The most difficult part of the proof is to show that for each special interval $J$, a certain uniform (depending on $K$) neighborhood in $X$ of the union $Fl_{K}(X_u)\cup Fl_{K}(X_v)$ is 
uniformly properly embedded in $X$ and uniformly hyperbolic (Corollary \ref{cor:two-flows}). The idea  is to:

(i) embed  such union 
 in a larger ``modified'' flow-space 
$$
Fl_{K}(X_u)\cup Fl_{K}(X_v)\cup Fl_R(\HHH)
$$
(for a certain $R=R(K)$ and a metric bundle $\HHH$ over the interval $J=\llbracket u, v\rrbracket$), 

(ii) prove uniform hyperbolicity of a uniform neighborhood of this triple union using quasiconvex amalgamation (Theorem \ref{thm:hyp-tree}),

(iii)  lastly,  use the fact that a uniform neighborhood of 
the union of intersecting quasiconvex subsets of a hyperbolic space is again uniformly 
quasiconvex and hyperbolic.  

\medskip 
We then use quasiconvex amalgamation to prove the same result for the union of 
three  flow-spaces 
$$
Fl_{K}(X_u)\cup Fl_{K}(X_w)\cup Fl_{K}(X_v),$$
where each subinterval $\llbracket u, w\rrbracket$, $\llbracket w, v\rrbracket$ is special.  
Theorem \ref{thm:geod-flow-step1} is then proven by verifying that the family of paths in 
$Fl_K(X_J)$ which are geodesics in pairwise unions $N_D(Fl_{K}(X_s)\cup Fl_{K}(X_t)), s, t\in V(J)$, 
satisfy the slim combing axioms from Corollary \ref{cor:bowditch}.

\subsection{Proper embeddings of  unions of pairs of intersecting 
flow-spaces}\label{sec:proper-generalized hallways}

Assuming that $Fl_{K}(X_u)\cap X_v\neq \emptyset$, we will show that a certain uniform neighborhood of the union of two flow-spaces $Fl_{K}(X_u)\cup Fl_{K}(X_v)$ is uniformly properly embedded in $X$.  We first deal with the following easier case when $T$ is an interval, which we will identify with an interval $[0, n]\subset \RR$, $n\in \N$, and that the vertex-set of $T$ equals to the set of integer points in the interval. Recall that $M_k$ is the parameter from the definition of uniform $k$-flaring (Definition \ref{uniform flaring}).

\begin{lemma}\label{proper embedding in Z-spaces}
Suppose $\X=(\pi: X\to T)$ is a tree of hyperbolic spaces, such that the tree $T$ is an interval 
$T=[0,n]$, $K\geq K_0$ and $\X$  satisfies the uniform $k$-flaring condition for $k=(L+1)^2 K$, where $L=L_{\ref{mjproj}}(K)$. 
Assume, moreover, that  for vertices $u, v\in V(T)$, 
and a ${4\delta_0}$-quasiconvex subset $Q=Q_u\subset X_u$, we have $Fl_{K}(Q)\cap X_v\neq \emptyset$.

 Then  the fiberwise $M_k$-neighborhood of $Fl_{K}(Q)\cup Fl_{K}(X_n)$ 
is uniformly properly embedded in $X$, with the distortion function depending only on $K$. 
\end{lemma}

\proof It suffices  to show that for each $D>0$ and $x,y \in {\mathcal Fl}_{K}(Q)\cup {\mathcal Fl}_{K}(X_v)$ with $d(x,y)<D$, the intrinsic distance between 
$x$ and $y$ in the $M_k$-fiberwise neighborhood $U$ of $Fl_{K}(Q)\cup Fl_{K}(X_v)$ is bounded by a constant depending on 
$D, K$ only. 

\begin{claim}
The statement of the lemma holds for $v=n$. 
\end{claim}
\proof Without loss of generality $x\in {\mathcal Fl}_{K}(Q)\setminus {\mathcal Fl}_{K}(X_n)$ and $y\in {\mathcal Fl}_{K}(X_n)\setminus Fl_{K}(Q)$. 
(Otherwise, the claim follows from Theorem \ref{mjproj}.) In the proof we will be repeatedly using Mitra's retractions $\rho$ defined in Theorem \ref{mjproj}. 

\medskip 
{\bf Reduction 1.} We first reduce to the case where $x,y$ are in the same vertex-space and, moreover,  $x=\rho_{Fl_K(Q)}(y)$. 

Observe that, since $T$ is an interval with an extremal vertex $n$, $y\in Fl_K(X_n)$ and $Fl_{K}(X_u)\cap X_n\neq \emptyset$, $\pi(y)\in \pi (Fl_{K}(Q))$.  
In particular, by the definition of Mitra's retraction $\rho_{Fl_K(Q)}$, 
$$
x':=\rho_{Fl_K(Q)}(y) \in X_{\pi(y)}\cap {\mathcal Fl}_{K}(Q).$$

We also apply Mitra's retraction $\rho_{Fl_K(Q)}: X\to Fl_{K}(A)$ to $xy$, a geodesic in $X$ connecting $x$ to $y$. 
The image $\rho_{Fl_K(Q)}(xy)$ is  a path of length $\le D_1:=(D+1)L_{\ref{mjproj}}(K)$ in 
$Fl_{K}(Q)$ joining $x$ to $x'$. Thus, 
$$
d_{Fl_K(Q)}(x,x')\le D_1,
$$
while $d(x',y)\leq d(x,y)+d(x',x)\leq D+D_1$. Hence, we replace $x$ by $x'$ and, henceforth, assume that $x, y$ are in the
same vertex-space and $d(x,y)\le D$ (with $D:=D+D_1$). 
Note that $\pi(x)\neq n$ since $x\in Fl_{K}(Q)\setminus Fl_{K}(X_n)$. 

\medskip 
We therefore, assume from now on that $j=\pi(x)=\pi(y)$ and $x=\rho_{Fl_K(Q)}(y)$. 

\medskip
{\bf Reduction 2.} We next reduce to the case when both $x, y$ are connected by $k$-qi leaves to the same point $z\in X_n$, where $k$ is as in the statement of lemma.

Let $\gamma_y$ be a $K$-qi in $Fl_K(X_n)$ leaf joining $y$ to some $z\in X_n$. 
In view of the assumption that $x=\rho_{Fl_K(Q)}(y)$, the path 
$$
\bar{\gamma}_{x}:= \rho_{Fl_K(Q)}(\gamma_y)\subset Fl_{K}(Q)$$
 connects $x$ to $z_1=\bar{\gamma}_{x}(n)=\rho_{Fl_K(Q)}(z)\in X_n$. 
For $y_1= \rho_{Fl_K(X_n)}(x)$, the path   
$$
\gamma_{y_1}:= \rho_{Fl_K(X_n)} \circ \rho_{Fl_K(Q)}(\gamma_y)\subset X_n
$$
also connects $y_1$ to $z_1$. Since both Mitra's retractions that we used are $L=L_{\ref{mjproj}}(K)$-coarse Lipschitz, the paths $\bar{\gamma}_{x}, \gamma_{y_1}$ are 
$k$-qi leaves for $k=(L+1)^2 K$.

Observe that since projections of $Fl_{K}(Q)$ and $Fl_{K}(X_n)$ to $T$ both contain 
the interval $\llbracket j,n\rrbracket$, both Mitra's  retractions restricted to the vertex spaces $X_i, i\in \llbracket j,n\rrbracket$, amount to 
fiberwise nearest-point projections to respective flow-spaces. In particular, $y_1\in X_j$. 

We next estimate $d_{X_j}(y, y_1)$. Since $\rho_{Fl_K(X_n)}$ is $L$-coarse Lipschitz, $\rho_{Fl_K(X_n)}(y)=y$, $\rho_{Fl_K(X_n)}(x)=y_1$ and $d(x,y)\le D$,  we obtain 
$$
d_{Fl_{K}(X_n)}(y, y_1)\le (L+1)D. 
$$
At the same time, $d(x, y_1)\le d(x, y)+ d(y,y_1)\le 2D$. 
Therefore, it suffices to prove the claim for the pair of points $x, y_1$: Both are connected to $z_1\in X_n$ by $k$-leaves $\bar{\gamma}_{x}\subset   Fl_K(Q) , \ga_{y_1}\subset Fl_K(X_n)$ respectively. We now reset $z:=z_1$, $D:= 2D$. 

\medskip 
Thus, we consider the case of points $x, y\in X_j$ such that $d(x,y)\le D$,  $x\in Fl_{K}(Q) \cap X_j, y\in Fl_{K}(X_n)\cap X_j$ and there exist $k$-qi leaves 
$\ga_x, \ga_y$ in $Fl_{K}(Q)$, $Fl_{K}(X_n)$, respectively, connecting $x, y$ to a point $z\in X_n$. 
Let $t\in \llbracket j, n\rrbracket$ be the minimal vertex such that 
$$
d_{X_t}(\ga_x(t), \ga_y(t))\le M_k.$$ 
(Such $t$ exists since $d_{X_n}(\ga_x(n), \ga_y(n))=0$.)  Recall that the vertex-spaces of $X$ are $\eta_0$-uniformly properly embedded in $X$; 
in particular, $d_{X_j}(x,y)\le \eta_0(D)$. 

If $d_{X_j}(x,y)\le M_k$, we will be done by taking 
$$
D_{\ref{proper embedding in Z-spaces}}(K):= M_k,
$$
since the intrinsic distance between $x$ and $y$ in $U:=N^{fib}_{M_r}( Fl_K(Q)\cup Fl_K(X_n))$ would be $\le d_{X_j}(x,y)\le \eta_0(D)$.

Otherwise (if  $d_{X_j}(x,y)> M_k$), by the uniform $k$-flaring condition, the length of the interval 
$\llbracket j, t\rrbracket$ is at most $\tau_{\ref{prop:weak flaring}}(k,\max(\eta_0(D), M_k))$. 
Therefore, the intrinsic distance between $x, y$ in $U$ is at most
$$
k \tau_{\ref{prop:weak flaring}}(k,\max(\eta_0(D), M_k)) + M_k.  
$$
This concludes the proof of the claim. \qed

\medskip

We return to the proof of the lemma and consider the general case, when $v$ need not be equal to the extreme vertex $n$ of $T=\llbracket 0,n\rrbracket$. 
Let $x\in X_{u_0}, y\in X_{v_0}$ be points within distance $D$, $u_0, v_0\in V(T)$. 
There are two cases to consider, depending on the order of the vertices $v, u_0$ in $\llbracket 0,n\rrbracket$. 

1. Suppose that $v$ lies in the interval $\llbracket u, u_0\rrbracket$ (or $\llbracket u_0, u\rrbracket$). 
Since $x\in Fl_K(X_u)\cap X_{u_0}$ Proposition \ref{prop:flow-prop} (specifically, \eqref{i} with $r=0$) 
implies that $x\in Fl_K(X_v)$. Thus, $x, y$ both belong to $Fl_K(X_v)$ and we conclude using   
the fact that flow-spaces are uniformly quasiisometrically embedded in $X$. 

2. Suppose that $v$ does not lie in the interval between $u, u_0$ in $T$. Without loss of generality (reversing the orientation on $T$ if necessary), with respect to the order on 
the interval $T$, 
$$
\max(u, u_0)< v.  
$$
If $v_0\le v$ then we  can shorten the tree $T$ replacing it with the subinterval $\llbracket 0, v\rrbracket$ and reduce the problem to the one solved in the claim above. 
Thus, we can assume that 
$$
0\le \max(u, u_0)< v < v_0\le n. 
$$
In particular, $d_T(u_0, v_0)= d_T(u_0, v)+ d_T(v, v_0)$. 
By projecting a geodesic $xy$ in $X$ to the tree $T$, we see that $d_T(u_0, v_0)\le D$ and, therefore,
$$
\max \left( d_T(u_0, v), d_T(v, v_0)\right ) \le D. 
$$ 
Hence,  a $K$-qi leaf $\ga_y$ in $FL_K(X_v)$  connecting $y$ to $y_1\in X_v$ has length 
$\le K D$, which implies that
$$
d(y, y_1)\le    d_{FL_K(X_v)}(y, y_1) \le (K+1)D. 
$$

Consider the subtree $S=  \llbracket 0, v\rrbracket$ in $T$. According to Lemma \ref{unif-emb-subtree}, the inclusion $X_S\to X$ is an 
$\eta=\eta_{\ref{unif-emb-subtree}}$-uniformly proper map (where the function $\eta$ depends only on the parameters of the tree of spaces $\X$). 
In particular, 
$$
d_{X_S}(x,y_1)\le D_1:=\eta((K+1)D).  
$$
Thus, it suffices to estimate the distance $d_{X_1}(x, y_1)$ in the  fiberwise 
$M_k$-neighborhood of $Fl_K(Q)\cup FL_K(X_v)\cap X_S$. Since $v$ is an extreme vertex of $S$, 
this is done in the claim above. Lemma follows. \qed

\medskip
Lastly, we consider the case of a general tree $T$:

\begin{prop}\label{prop:two-flows-easy}
Suppose $\X=(\pi: X\to T)$ is a tree of hyperbolic spaces, $K\geq K_0$ 
and $\X$  satisfies the uniform $k$-flaring condition, where, as in the lemma,  
$$
k= k_{\ref{prop:two-flows-easy}}(K)=(L+1)^2 K, \quad L=L_{\ref{mjproj}}(K).$$
Assume, moreover, that  for some   ${4\delta_0}$-quasiconvex subset $Q=Q_u\subset X_u$, we have $Fl_{K}(Q)\cap X_v\neq \emptyset$. 
 Then  the $M_k$-neighborhood of $Y=Fl_{K}(Q)\cup Fl_{K}(X_v)$ in $X$ 
is $\eta_{\ref{prop:two-flows-easy},K}$-properly embedded in $X$. 
\end{prop}
\proof  Suppose 
$$
x\in Fl_{K}(X_u) \cap X_{u_0}, y\in Fl_{K}(X_v) \cap X_{v_0}
$$
satisfy  $d(x,y)\leq D$. We need to show that there is a
constant constant $D_1>0$ depending only on $D$ and $K$ 
such that the distance between $x, y$ 
in $N_{M_k}(Y)$ is at most $D_1$.

If both $x,y$ are either in $Fl_{K}(X_u)$ or $Fl_{K}(X_v)$ then certainly this is true since $K$-flows of vertex-spaces are uniformly quasiisometrically embedded in $X$.

Suppose, therefore, that 
\begin{equation}\label{eq:xy-location}
x\in {\mathcal Fl}_{K}(X_u)\setminus {\mathcal Fl}_{K}(X_v) \hbox{~~and~~} y\in {\mathcal Fl}_{K}(X_v)\setminus {\mathcal Fl}_{K}(X_u)\end{equation}
If $u, v$ span an edge in $T$ 
then our assumptions on the location of $x$ and $y$ imply that  the set of vertices $\{u_0=\pi(x), v_0=\pi(y), u, v\}$ is contained in a common interval in $T$ and, hence, the claim follows from 
Lemma \ref{proper embedding in Z-spaces}. Thus, we assume that $d_T(u, v)\ge 2$.  More generally, we can assume that 
$u_0, v_0, u, v$ do not belong to a common interval in $T$.

The same assumption \eqref{eq:xy-location}  implies that the center of $\Delta u_0 uv\subset T$ is in 
$\llbracket u,v\llbracket $ and the center of $\Delta v_0uv$ is in $\rrbracket u,v\rrbracket$,  
and if these centers are equal, then they are in $\rrbracket u,v\llbracket$.

The proof of the proposition is broken in three cases:

\begin{figure}[tbh]
\centering
\includegraphics[width=80mm]{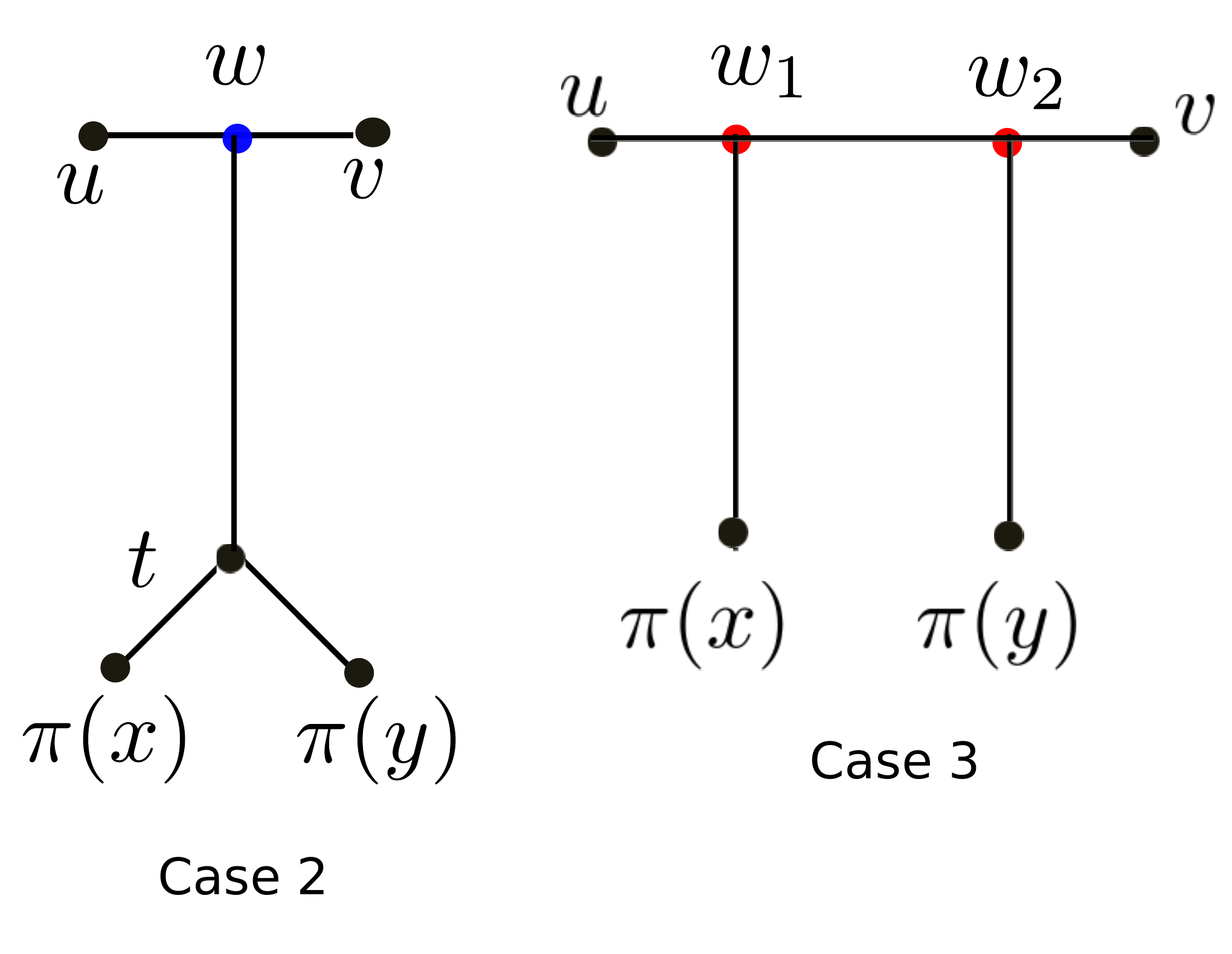}  
\caption{Cases 2 and 3}
\label{cases-2-3.fig}
\end{figure}

\medskip
{\bf Case 1:} Suppose $x, y$ are in the same vertex-space $X_t$; we let  $w$ denote the center of $\Delta uvt$.
As noted before, $w\in \rrbracket u,v\llbracket $. Let $Q= Q_w= Fl_{K}(X_u)\cap X_w$; this is a  $4\delta_0$-quasiconvex subset of $X_w$.

 For the interval $S= \llbracket v, t\rrbracket$ consider the subtree of spaces $X_S\subset X$. 
The points $x, y$ belong to the $K$-flow-spaces $Fl_K(Q), Fl_K(X_v)$ in $X_S$, and  
$$
\emptyset \ne Fl_K(X_u)\cap X_v\subset Fl_K(Q)\cap X_v. 
$$
We also have $d_{X_S}(x,y)\le \eta(D)$, where $\eta= \eta_{\ref{unif-emb-subtree}}$. 
Thus, the conclusion follows from Lemma \ref{proper embedding in Z-spaces} 
applied to the tree of spaces $X_S$ over the interval $S$.

\medskip 
{\bf Case 2:} Suppose the  triangles  $\Delta  u_0uv$ and $\Delta v_0uv$ have the common center $w$, see Figure \ref{cases-2-3.fig}. 
Let $t$ denote the center of the triangle $\Delta w u_0 v_0$. Without loss of generality, we may assume that $d_T(\pi(x),t)\leq d_T(\pi(x),\pi(y))$. 
Since $d(x,y)\le D$, we also obtain $d_T(v_0,t)\leq D$. Let $\ga_x, \ga_y$ denote, respectively, $K$-qi leaves in $Fl_{K}(X_u), Fl_{K}(X_v)$ 
connecting $x$, $y$ to $X_u, X_v$. Then, for $x_1=\ga_x(t)$ and $y_1=\ga_y(t)$ we have 
$$
d_{Fl_K(Q)}(x,x_1)\leq K D \hbox{~~and~~} d_{Fl_K(X_v)}(y, y_1)\leq KD.
$$
 In particular $d(x_1,y_1)\leq (1+2K)D$ and $x_1, y_1$ belong to the same vertex-space $X_t$. 
This reduces the proof to that of Case 1.

\medskip 
{\bf Case 3:} Suppose the triangles $\Delta u_0uv$ and $\Delta v_0uv$ have distinct centers $u_1$ and $v_1$ respectively. 
These centers necessarily belong to the interval $\llbracket u,v\rrbracket$, see Figure \ref{cases-2-3.fig}.  
As in Case $2$ we first we take two $K$-qi leaves  $\gamma_x, \ga_y$ in $Fl_{K}(X_u), Fl_{K}(X_v)$, connecting $x$, $y$ respectively to $X_u$ and $X_v$. 
Since $d(x,y)\le D$, we also have $d_T(u_0,u_1)\le D, d(v_0, v_1)\le D$.

We then replace $x$ with $x_1=\ga_x(u_1)$ and replace $y$ with $y_1=\ga_y(v_1)$, which are within distance $KD$  from $x$, respectively $y$, in $Fl_{K}(X_u), Fl_{K}(X_v)$. 
Furthermore,
$$
d(x_1, y_1)\le D(1+ 2K). 
$$
Since $x_1, y_1$ belong to the subtree of spaces $X_{uv} = \pi^{-1}(uv)$, analogously to Case 1, we conclude the proof by applying Lemma \ref{proper embedding in Z-spaces}. \qed

\begin{cor}\label{cor:two-flows-easy}
Suppose $\X=(\pi: X\to T)$ is a tree of hyperbolic spaces, $K\geq K_0$ 
and $\X$  satisfies the uniform $k$-flaring condition for 
$$
k= k_{\ref{prop:two-flows-easy}}(K).$$
Assume also that  for a ${4\delta_0}$-quasiconvex subset $Q=Q_u\subset X_u$, we have 
$Fl_{K}(Q)\cap X_v\neq \emptyset$.

 Then for each $r\ge M_k$, the   $r$-neighborhood of $Fl_{K}(Q)\cup Fl_{K}(X_v)$ in $X$ 
is $\eta_{\ref{cor:two-flows-easy},K,r}$-properly embedded in $X$.
\end{cor}

\subsection{Hyperbolicity of the union of two flow-spaces in the case of special intervals}\label{sec:generalized hallways}

The goal of this section is to prove 
that a uniform neighborhood of 
$Fl_{K}(X_u)\cup Fl_{K}(X_v)$ in $X$ is uniformly hyperbolic, with the hyperbolicity constant depending only on $K$, provided that
 $Fl_{K}(X_u)\cap X_v\neq \emptyset$ and $\X$ satisfies a suitable uniform flaring condition. 

Here is the  idea of the proof. Recall that $Fl_{K}(X_u)$ and $Fl_{K}(X_v)$ are (uniformly) hyperbolic and are (uniformly) 
qi embedded in $X$. If $X$ were hyperbolic, it would follow (since $Fl_{K}(X_u)\cap Fl_{K}(X_v)\ne \emptyset$)   
that  a uniform neighborhood of $Fl_{K}(X_u)\cup Fl_{K}(X_v)$ is uniformly hyperbolic. Hyperbolicity of $X$, of course, is 
not yet proven, so instead, we will find (see the proof of Corollary \ref{cor:two-flows-easy}) a larger subset $U=U_{r_1}$ containing 
$Fl_{K}(X_u)\cup Fl_{K}(X_v)$, which is uniformly hyperbolic and uniformly properly embedded in $X$ 
 (Proposition \ref{flow modification}).  
Thus, a suitable uniform neighborhood of  $Fl_{K}(X_u)\cup Fl_{K}(X_v)$ in $U$ is 
uniformly hyperbolic and uniformly properly embedded in $U$. From this (since $U$ is uniformly properly embedded in $X$), 
we will conclude that a suitable uniform neighborhood of  $Fl_{K}(X_u)\cup Fl_{K}(X_v)$ in $X$ is also uniformly hyperbolic 
and uniformly properly embedded in  $X$.

\begin{lemma}\label{leaves in both flows}
For all $D\ge 0, K\ge K_0$ there is a constant $K_{\ref{leaves in both flows}}=K_{\ref{leaves in both flows}}(K)\ge K$ 
such that the following holds:

Let $u, v\in T$ be vertices such $Fl_K(X_u)\cap X_v\ne \emptyset$. 
Then for each 
$$
x\in N_{D}({\mathcal Fl}_{K}(X_u))\cap N_D({\mathcal Fl}_{K}(X_v)),$$ 
there is a vertex $t\in V(T)$ and  a $K_{\ref{leaves in both flows}}$-qi section $\Sigma$ of $\pi:X\map T$ over the tripod (triangle) $S=\Delta tuv$ such that 
$$
x\in N_{3D}(\Sigma \cap X_t)
$$
and 
$$
\Sigma \subset Fl_{K}(X_u)\cup Fl_{K}(X_v). 
$$
\end{lemma}
\proof First of all, there is a vertex $t\in V(T)$ and points  $x_1\in X_t\cap Fl_{K}(X_u)$, 
$x_2\in X_t\cap Fl_{K}(X_v)$  such that 
$$
d_{X_t}(x,x_i)\leq D, i=1, 2.$$
Let $\gamma_1$ be a $K$-qi leaf in $Fl_{K}(X_u)$ connecting $x_1$ to $X_u$. 
We apply Mitra's projection $\rho=\rho_{Fl_{K}(X_v)}$ to $\gamma_1$; call the result $\ga'_1$.  Since, by the assumption of the lemma, 
$$
S\subset \pi(Fl_{K}(X_u)),
$$
$\rho$ restricted to $X_S$ amounts to the fiberwise projection to $Fl_{K}(X_v)$. In particular, 
$\pi(\ga'_1)= \llbracket u, t\rrbracket$. Hence, $\ga'_1$ is a $K_{\ref{leaves in both flows}}:= KL_{\ref{mjproj}}(K)$-qi 
section over $ \llbracket u, t\rrbracket$ whose image is contained in $Fl_{K}(X_v)$. 

Let $w\in uv$ denote the center of the triangle $S$. Since $\ga'_1$ is contained in $Fl_{K}(X_v)$, 
the point $\gamma'_1(w)$ can be joined to $X_v$ by a $K$-qi leaf $\gamma'_2$ inside $Fl_{K}(X_v)$.
Clearly, the union of these two qi leaves $\ga_1'\cup \ga_2'$ forms a 
$K_{\ref{leaves in both flows}}$-qi section $\Sigma$ over the tripod $S$. 

We have
$$
d_{X_t}(x_1, \gamma'_1(t))\le d_{X_t}(x_1,x_2)\le 2D,$$
$$
d(x,x_1)\leq D, \quad d(x, \gamma'_1(t))\leq 3D.$$
Lemma follows. \qed

 \medskip
Set 
\begin{equation}\label{eq:K2}
k=K_{\ref{leaves in both flows}}(K)   
\end{equation} 
and define $\mathcal S= \mathcal S_{k, J}$, the set of all $k$-qi sections over the interval $J=\llbracket u, v\rrbracket$. 

Assuming that $Fl_K(X_u)\cap X_v\ne \emptyset$ (which is the standing assumption of the previous and this section), 
$\mathcal S$ is nonempty since $k\ge K$ and we are assuming that
$$
Fl_K(X_u)\cap X_v\ne \emptyset. 
$$
 For each vertex $w\in V(J)$, let 
${H}_w$ denote the (fiberwise) $\delta_0$-hull of the subset $\{\gamma(w):\gamma\in \mathcal S\}$ in $X_w$.  
Define 
$$
\mathcal H:=\bigcup_{w\in V(J)} H_w.$$
Each $H_w$, of course, is a $4\delta_0$-quasiconvex subset of  $X_w$. Then Lemma \ref{lem:E-ladder-structure}(b) implies that 
${\mathcal H}$ is the union of vertex-spaces of a 
$k'= K'_{\ref{lem:E-ladder-structure}}(k)$-metric bundle ${\mathfrak H}={\mathfrak H}_{k,J}$ 
over the interval $J$, see Definition \ref{defn:flow0}.

\medskip 
As in Section \ref{sec:gfs}, we define the generalized $\kappa$-flow-spaces $Fl_{\kappa}({\mathfrak H})$ 
of the metric bundle ${\mathfrak H}$. From the definition, we recall that with each vertex $w\in V(J)$ we associate a 
subtree $T_w\subset T$ equal to the maximal subtree in $T$ containing $w$ and disjoint from all other vertices of the interval $J$. 
 
Below we will frequently use the function 
$$
\kappa\mapsto \kappa^\wedge=(15 L'_0 \kappa)^3,
$$
 defined in \eqref{eq:vee}.

\begin{lemma}\label{section5:flow intersection}
For $D\ge 0, K\ge K_0$ set $k=K_{\ref{leaves in both flows}}(K)$, $k'=K'_{\ref{lem:E-ladder-structure}}(k)$. Then the following hold:  

\begin{enumerate}

 \item Suppose that 
$$
x\in N_D({\mathcal Fl}_{K}(X_u))\cap N_D({\mathcal Fl}_{K}(X_v)).$$ 
Then there exists a vertex $t\in V(T)$ with 
$$
d(x, {\mathcal Fl}_{K}(X_u) \cap {\mathcal Fl}_{K}(X_v) \cap X_t)\le D,$$
such that the center $w$ of $\Delta uvt$ satisfies
$$
x\in N_{3D} (Fl_{k^\wedge}( H_w)\cap X_{T_w}).$$
In other words, we have
 $$
 N_D({\mathcal Fl}_{K}(X_u))\cap N_D({\mathcal Fl}_{K}(X_v))\subset 
 \bigcup_{w\in V(\llbracket u, v\rrbracket)} N_{3D} ( Fl_{k^\wedge}(H_w)\cap X_{T_w}).
 $$ 
 
\item For each $\kappa\ge k'$  
and each vertex $w\in V(\llbracket u, v\rrbracket)$ we have 
$$
Fl_{\kappa}(H_w)\cap X_{T_w}\subset Fl_{\kappa^\wedge}(X_u)\cap Fl_{\kappa^\wedge}(X_v).$$
In other words, 
$$
Fl_{\kappa}({\mathfrak H}) \subset Fl_{\kappa^\wedge}(X_u)\cap Fl_{\kappa^\wedge}(X_v).
$$

\end{enumerate}

\end{lemma} 

\proof  (1) By Lemma \ref{leaves in both flows} there is a $k=K_{\ref{leaves in both flows}}(K)$-qi section $\Sigma$ over a 
tripod $uv\cup wt$ such that
$$
d(x, \Sigma)\le 3D. 
$$
In particular, $\Sigma$ restricted to $J=\llbracket u, v\rrbracket$ is a $k$-qi
section over $J$   
hence, belongs to $\mathcal S_{k, J}$. 
Therefore, by the definition of ${\mathcal H}$, for the vertex $w\in V(J)$, the intersection 
$\Sigma\cap X_w$ belongs to $H_w$. 

Since $\Sigma(w)$ is in $H_w$, the restriction of $\Sigma$ to the interval $\llbracket t, w\rrbracket$   is a
$k$-qi leaf connecting $\Sigma\cap X_t$ to $H_w$ and, by Proposition \ref{prop:flow-prop}(2), 
$$
\Sigma\cap X_t\in Fl_{k^\wedge}({H}_w).
$$
Thus,
$$
x\in N_{3D} (Fl_{k^\wedge}(H_w) ) \cap X_t.
$$
This proves Part (1) for fiberwise neighborhoods. The proof for neighborhoods taken in $X$ is identical and we omit it.

\medskip 
(2) Given $x\in {\mathcal Fl}_{\kappa}(H_w)\cap X_{T_w}$, pick  a 
$\kappa$-qi leaf  $\gamma_x$ in the flow-space $Fl_{\kappa}(H_w)$ connecting $x$ to $H_w$. 

Since $z=\gamma_x(w)\in H_w$, and ${\mathfrak H}$ is a $k'$-metric bundle  over the interval $J= \llbracket  u, v\rrbracket$, 
there exists a $k'$-qi leaf $\ga_z$ in this bundle (over $\llbracket w, v\rrbracket$) connecting $z$ to $X_v$. 
Thus, since $\kappa\ge k'$, the point $x$ can be connected to both $X_u$ and $X_v$ by $\kappa$-qi leaves. 
According to Proposition \ref{prop:flow-prop}(2), 
$$
x\in Fl_{\kappa^\wedge}(X_u)\cap Fl_{\kappa^\wedge}(X_v). \qed 
$$

\medskip
For the rest of this section (until Corollary \ref{cor:two-flows-easy}) 
we will be working under the following  assumption (which is stronger than the one we had earlier): 

\begin{itemize}
\item 
We assume that $X_v\cap Fl_{K}(X_u)\neq \emptyset, X_u\cap Fl_{K}(X_v)\neq \emptyset$.
\end{itemize}

\begin{rem}
1. This assumption implies that for all vertices  $w\in V(\llbracket u,v\rrbracket)$,  
$$
X_w\cap Fl_{K}(X_u)\neq \emptyset,  X_w\cap Fl_{K}(X_v)\neq \emptyset. 
$$

2. The stronger assumption we are now making is not too far from the condition that $X_v\cap Fl_{K}(X_u)\neq \emptyset$ 
made earlier, since  
$$
X_u\cap Fl_{K^\wedge}(X_v)\neq \emptyset,  
$$
see Proposition \ref{prop:flow-prop}(2). We will be using this fact in the proof of Corollary \ref{cor:two-flows-easy}. 
\end{rem}

\medskip

We now fix some $K\ge K_0$, 
set, $k= K_{\ref{leaves in both flows}}(K)$ and take some $R\ge k^\wedge$.  
For an interval $J= \llbracket u , v\rrbracket\subset T$, set ${\mathfrak H}:={\mathfrak H}_{k,J}$ and  
$Y^0:= Fl_{R}({\mathfrak H})$. Since $\Fl_R(\HHH)$ is a generalized flow-space with the parameters $K_1= k'$ and $K_2= R\ge K_0$, the next lemma is a corollary of Theorem \ref{thm:generalized retraction}: 

\begin{lemma}\label{lem:Y0 coarse retract}
The inclusion map $Y^0\to X$ is a $L_{\ref{thm:generalized retraction}}(k',R)$-qi embedding. 
\end{lemma}

We also define the unions $Y^1= Fl_K(X_u)\cup Fl_{R}({\mathfrak H})$,  $Y^2= Fl_K(X_v)\cup Fl_{R}({\mathfrak H})$ and their neighborhoods 
$U^i_r:= N_r(Y^i), i=1,2$, taken in $X$.  

\begin{lemma}\label{lem:Y1-cap-Y2}
For every $r\ge 0$, 
$$
Fl_{R}({\mathfrak H}) \subset U^1_r\cap U^2_r \subset N_{3r}(Fl_{R}({\mathfrak H})), 
$$
i.e. the intersection is uniformly (in terms of $r, R$ and $K$) Hausdorff-close to $Fl_{R}({\mathfrak H})$. 
\end{lemma} 
\proof Consider $x\in U^1_r\cap U^2_r$. Thus, there exist points $x_1\in Fl_K(X_u)\cup Fl_{R}({\mathfrak H})$, 
 $x_2\in Fl_K(X_v)\cup Fl_{R}({\mathfrak H})$ at distance $\le r$ from $x$. If one of these points  
 is in $Fl_{R}({\mathfrak H})$ then $d(x, Fl_{R}({\mathfrak H}))\le r$.  
 Therefore, assume that $x_1\in Fl_K(X_u)$, $x_2\in Fl_K(X_v)$. By Lemma \ref{section5:flow intersection}(1), 
 $d(x, Fl_{R}({\mathfrak H}))\le 3r$, as required. \qed

\medskip 
Recall that $K\ge K_0$, $k= K_{\ref{leaves in both flows}}(K)$. Set $R:= k^\wedge$. 
In the next proposition, $N'$ indicates a metric neighborhood of $Y_1$ or $Y_2$ 
taken inside  $Fl_{R^\wedge}(X_u)$ or 
$Fl_{R^\wedge}(X_v)$ respectively. 
The most useful part of the  proposition is (2): Part (1) is used only to prove (2). 

\begin{prop}\label{flow modification}
Assume that $K\ge K_0$ and $\X$ satisfies the uniform $\kappa$-flaring condition for 
$$
\kappa=\max(k_{\ref{prop:two-flows-easy}}(K), \kappa_{\ref{defn:carpeted ladder}}(R^\wedge)). 
$$
Then there exist  
$\delta'_{\ref{flow modification}}=\delta'_{\ref{flow modification}}(K,C)$, 
$\delta_{\ref{flow modification}}=\delta_{\ref{flow modification}}(K)$, 
$L'_{\ref{flow modification}}=L'_{\ref{flow modification}}(K,C)$, 
$C_{\ref{flow modification}}= C_{\ref{flow modification}}(K)$, 
and a function $\eta_{\ref{flow modification}}=\eta_{\ref{flow modification},K}$ 
such that the following hold.   
\begin{enumerate}

\item For each $C\ge C_{\ref{flow modification}}$, both $U'_1=N'_{C}(Y^1)$ and 
$U'_2=N'_{C}(Y^2)$,  equipped with the induced path-metrics,  
are  $\delta'_{\ref{flow modification}}(K,C)$-hyperbolic and $L'_{\ref{flow modification}}(K,C)$-qi embedded in $X$. 

\item For  
$$
r:= \max(C_{\ref{flow modification}},M_{k_{\ref{prop:two-flows-easy}}(K)}),$$ 
the union 
$$
U_r:= U_r^1\cup U_r^2$$ 
equipped with the induced path-\-met\-ric, is $\delta_{\ref{flow modification}}(K,R)$-hyperbolic and 
$\eta_{\ref{flow modification}}$-unifor\-mly properly embedded in $X$. 
\end{enumerate}
\end{prop}
\proof (1)  We will only prove the claim for $U'_1$ 
 since the proof for $U'_2$ is obtained by relabelling. 
Recall that
$$
 Y^1=Fl_{K}(X_u)\cup Y^0. 
$$
Since $R\ge k'$, the definition of $Y^0$  and Lemma \ref{section5:flow intersection}(2) imply that 
$$
Y^0=Fl_{R}({\mathfrak H})\subset  Fl_{R^\wedge}(X_u).  
$$
Furthermore, since $R^\wedge\ge {R\ge  k}\ge K$, $Fl_{K}(X_u)\subset Fl_{R^\wedge}(X_u)$. Thus,
$$
Y^1\subset Fl_{R^\wedge}(X_u).
$$ 
Recall that we are assuming that $\X$ satisfies the uniform  {$\kappa_{\ref{defn:carpeted ladder}}(R^\wedge)$}-flaring condition. 
Therefore, Theorem \ref{flow of one vertex space} applies and the flow-space $Fl_{R^\wedge}(X_u)$ is 
$\delta=\delta_{\ref{flow of one vertex space}}(R^\wedge)$-hyperbolic. 

By Lemma \ref{lem:Y0 coarse retract}, $Y^0$ is $L_{\ref{thm:generalized retraction}}(k',R)$-qi embedded in $X$, while
$Fl_K(X_u)$ is  $L_{\ref{cor:mjproj}}(K)$-qi embedded in $X$ according  to Corollary \ref{cor:mjproj}.

Hence, for 
$$L:= 2\max(L_{\ref{thm:generalized retraction}}(k',R), L_{\ref{cor:mjproj}}(K))$$
and 
$$\la= \la_{\ref{lem:qi-preserves2}}(\delta, L),$$ 
both $Y^0$ and $Fl_K(X_u)$ are $\la$-quasiconvex in $Fl_{R^\wedge}(X_u)$. 
Moreover, these subsets have nonempty intersection (containing at least ${H}_u$). Thus, their 
union is $\la+\delta$-quasiconvex in  $Fl_{R^\wedge}(X_u)$.  

By Lemma \ref{lem:QCunion}(1), since  
$$
C\ge C_{\ref{flow modification}}\ge \la+\delta, 
$$ 
the $C$-neighborhood $U'_1$ of 
$Y^1=Fl_K(X_u)\cup Y^0$ in $Fl_{R^\wedge}(X_u)$ 
(equipped with the induced path-metric) 
is $\delta+C$-quasiconvex in 
$Fl_{R^\wedge}(X_u)$. 

Furthermore, by Part (2) of the same lemma, since $C$ was taken to be $\ge 2\la+4\delta$, 
the inclusion map $U'_1\to Fl_{R^\wedge}(X_u)$ is a 
$(1,6\delta+C)$-quasiisometric embedding, 
where $U'_1$ is equipped with the induced path-metric. 
Thus, $U'_1$ is $\delta_{\ref{lem:qi-preserves}}(\delta, 6\delta+C)$-hyperbolic, 
see Lemma \ref{lem:qi-preserves}. Moreover, the inclusion map
$U_1'\to X$ is an $L'_{\ref{flow modification}}$-qi embedding with 
$$
L'_{\ref{flow modification}}= L_{\ref{cor:mjproj}}(R^\wedge)(6\delta+C). 
$$

 This proves Part (1), where we use:
 $$
 \begin{aligned}
 \delta=\delta_{\ref{flow of one vertex space}}(R^\wedge), L=2\max(L_{\ref{thm:generalized retraction}}(k',R), L_{\ref{cor:mjproj}}(K)),
 \\ 
 \la= \la_{\ref{lem:qi-preserves2}}(\delta, L), 
 C_{\ref{flow modification}}= 2\la+4\delta,  \delta'_{\ref{flow modification}}= \delta_{\ref{lem:qi-preserves}}(\delta, 6\delta+C). 
 \end{aligned}
 $$
 
 \medskip 
(2) In the proof of Part (1) we were using hyperbolicity of the union of two quasiconvex subsets in a hyperbolic space. 
In Part (2), such ambient hyperbolic space is unavailable, so we will use hyperbolicity of 
{\em pairwise quasiconvex amalgams} of hyperbolic  spaces (Theorem \ref{thm:hyp-tree},  quasiconvex amalgamation), to prove hyperbolicity of $U_r$.  
We will be using Part (1) with $C=r$. 
 We have, of course,
$$
Y_i\subset U'_i= N'_r(Y_i) \subset U^i_r= N_r(Y_i), i=1, 2. 
$$
Since each inclusion map $U'_i\to U^i_r$ is an $L$-quasiisometry, $L'=\max(L_{\ref{cor:mjproj}}(K), r)$, 
Lemma \ref{lem:qi-preserves}  implies that $U^i_r$ is 
$\delta'_{\ref{lem:qi-preserves}}(\delta'_{\ref{flow modification}}, L')$-hyperbolic.

Since  $Y^0$ is $L_{\ref{thm:generalized retraction}}(k',R)$-qi embedded in $X$, it is 
$\la_{\ref{lem:qi-preserves2}}(\delta'_{\ref{flow modification}}, L_{\ref{thm:generalized retraction}}(k',R))$-quasiconvex in $U^i_r, i=1,2$. Since 
$Y^0$  is $3r$-Hausdorff close to 
 $U^1_r\cap U^2_r$ (Lemma \ref{lem:Y1-cap-Y2}), 
 it follows that  $U^1_r\cap U^2_r$ is  {$\la'$}-quasiconvex in 
 both $U^1_r, U^2_r$ where 
 $$\la'= 3r+\delta'_{\ref{flow modification}}  \la_{\ref{lem:qi-preserves2}}(\delta'_{\ref{flow modification}}, L_{\ref{thm:generalized retraction}}(k',R)).$$ 
 
 Thus, we are in position to 
apply Theorem \ref{thm:hyp-tree} (quasiconvex amalgams), and conclude that  the union $U_r=U^1_r\cup U^2_r$ is 
$\delta_{\ref{flow modification}}(K,D,R)$-hyperbolic. 
It remains to prove that $U_r$ is uniformly properly embedded in $X$, thereby proving half of (2).

\medskip 
Let $x_1, x_2$ be points in $U_r$. We need to estimate $d_{U_r}(x_1, x_2)$ in terms of $d(x_1,x_2)$. 
It suffices to consider the case when $x_i\in U^i_r$, $i=1,2$, since $U'_1, U'_2$ are both $L'_{\ref{flow modification}}$-qi embedded in $X$ 
and are within Hausdorff distance $r$ from $U^1_r, U^2_r$ respectively. Let $y_i\in Y^i$ be points at distance $\le r$ from $x_i$, $i=1,2$; 
$d(y_1,y_2)\le d(x_1,x_2)+2r$. Since $Y^1= Fl_{K}(X_u)\cup Y^0, Y^2= Fl_{K}(X_v)\cup Y^0$, using the fact that $Y^0$ is 
$L_{\ref{thm:generalized retraction}}(k',R)$-qi embedded in $X$, the problem reduces to the case when 
$y_1\in Fl_{K}(X_u), y_2\in Fl_{K}(X_v)$. Recall that according to Corollary \ref{cor:two-flows-easy} (applied to $Q_u=X_u$), 
since 
$$
r\ge M_{k_{\ref{prop:two-flows-easy}}(K)},
$$
the $r$-neighborhood of the union $Fl_{K}(X_u)\cup Fl_{K}(X_v)$ is 
$\eta_{\ref{cor:two-flows-easy},K}$-properly embedded in $X$. 
Therefore,  
$$
N_{r} (Fl_{K}(X_u)\cup Fl_{K}(X_v))\subset N_r(Y^1)\cup N_r(Y^2)=U_r, 
$$
and we obtain 
$$
d_{U_r}(x_1,x_2)\le d_{N_{M_k} (Fl_{K}(X_u)\cup Fl_{K}(X_v))}\le \eta_{\ref{prop:two-flows-easy},K}(d(x_1,x_2)+ 2r). 
$$
This concludes the proof of Part (2) and, hence, of the proposition.   \qed

\medskip
We can now prove the main result of this section:

\begin{cor}\label{cor:two-flows} 
Suppose that $K\ge K_0$ and  $u, v$ are vertices in $T$ such that $Fl_K(X_u)\cap X_v\ne\emptyset$. 
Set $k_{\ref{cor:two-flows}}:= K_{\ref{leaves in both flows}}(K^\wedge)$, 
$R_{\ref{cor:two-flows}}:= k_{\ref{cor:two-flows}}^\wedge$ and 
assume that $\X$ satisfies the uniform $\kappa$-flaring condition for 
$$
\kappa=\kappa_{\ref{cor:two-flows}}(K)=\max(k_{\ref{prop:two-flows-easy}}(K^\wedge), 
\kappa_{\ref{defn:carpeted ladder}}(R_{\ref{cor:two-flows}}^\wedge)). 
$$
Then there exist 
$\delta=\delta_{\ref{cor:two-flows}}(K)$, 
$D= D_{\ref{cor:two-flows}}(K)$, 
and a function $\eta=\eta_{\ref{cor:two-flows},K}$ 
such that the following hold:

The $D$-neighborhood $N_D(Fl_K(X_u)\cup Fl_K(X_v))$ (with the induced path-metric) in 
$X$ is $\delta$-hyperbolic and $\eta$-properly embedded in $X$. 
\end{cor}
\proof According to Proposition \ref{prop:flow-prop}, for $K_1= K^\wedge$,
$$
Fl_{K_1}(X_v)\cap X_u\ne \emptyset. 
$$
Of course, we still have
$$
Fl_{K_1}(X_u)\cap X_v\supset Fl_{K}(X_u)\cap X_v \ne \emptyset. 
$$
Therefore, Proposition \ref{flow modification}(2) applies and we get that for $r_1=r(K_1)$ as in the proposition, $U=U_{r_1}$ is 
$\delta(K_1)$-hyperbolic and $\eta_{K_1}$-properly embedded in $X$. 

Since both $Fl_K(X_u), Fl_K(X_v)$ are $L=L_{\ref{cor:mjproj}}(K)$-qi embedded in $X$ (hence, in $U$), they are 
$\la_1= \la_{\ref{lem:qi-preserves2}}(\delta(K_1), L)$-quasiconvex in $U$ (Lemma \ref{lem:qi-preserves2}). 
Set
$$
D:= \max(2\la_1 +4\delta(K_1), M_{k_{\ref{prop:two-flows-easy}}(K)}).
$$
By Lemma \ref{lem:QCunion}, the  
$D$-neighborhood $V$ of $Fl_K(X_u)\cup Fl_K(X_v)$ in $U$ is $(2\la_1 +5\delta(K_1))$-quasiconvex in 
$U$. By the same lemma, $V$ (equipped with its path-metric) is $(10\delta(K_1)+2\la_1)$-qi embedded in $U$. 
Hence, $V$ is $\delta_{\ref{lem:qi-preserves}}(\delta(K_1), 10\delta(K_1)+2\la_1)$-hyperbolic (see Lemma \ref{lem:qi-preserves}).

Note that  the Hausdorff distance in $X$ between $V$ 
(the $D$-neighborhood of  $Fl_K(X_u)\cup Fl_K(X_v)$ in $U$) and $N_D(Fl_K(X_u)\cup Fl_K(X_v))$ (the $D$-neighborhood in $X$) is $\le D$.  Since $U$ is $\eta_{K_1}$-properly embedded in $X$ and the inclusion map 
$V\to U$ is a $(10\delta(K_1)+2\la_1)$-qi embedding, the composition $V\to N_D(Fl_K(X_u)\cup Fl_K(X_v))$ is 
$\zeta$-proper for
$$
\zeta(t)= (10\delta(K_1)+2\la_1) \eta(t) + (10\delta(K_1)+2\la_1)^2. 
$$
 Corollary \ref{cor:hyp-th} now implies that  $N_D(Fl_K(X_u)\cup Fl_K(X_v))$ is 
$\delta_{\ref{cor:two-flows}}(K)$-hyperbolic for
$$
\delta_{\ref{cor:two-flows}}(K)= \delta_{\ref{lem:qi-preserves}}(\delta_{\ref{lem:qi-preserves}}(\delta(K_1), 10\delta(K_1)+2\la_1), 
\zeta(2D+1)).
$$
It remains to estimate the distortion of $N_D(Fl_K(X_u)\cup Fl_K(X_v))$ in $X$. Since $D\ge M_{k_{\ref{prop:two-flows-easy}}(K)}$, 
Corollary \ref{cor:two-flows-easy} applies and $Fl_K(X_u)\cup Fl_K(X_v)$ is 
$\eta_{\ref{cor:two-flows-easy},K,D}$-properly embedded in $X$. Therefore, its $D$-neighborhood is 
$\eta$-properly embedded in $X$ for 
$$
\eta_{\ref{cor:two-flows},K}(t)=\eta(t)= 2D+ \eta_{\ref{cor:two-flows-easy},K,D}(t+2D). \qed 
$$

\begin{rem}\label{rem:K*}
This corollary gives us the value of $K_*$ in Theorem \ref{thm:mainBF}: 
\begin{equation}\label{K_*}
K_*= \max(k_{\ref{prop:two-flows-easy}}(K^\wedge), 
\kappa_{\ref{defn:carpeted ladder}}(R_{\ref{cor:two-flows}}^\wedge)),
\end{equation}
where $K=K_0$. 
 \end{rem}

 \subsection{Hyperbolicity of flow-spaces of special interval-spaces $Fl_K(X_J)$} 

In this section we  conclude the proof of Theorem \ref{thm:geod-flow-step1}. We will also prove uniform hyperbolicity 
of  flow-spaces $Fl_K(X_J)$, whenever $J\subset T$ is a union of three special intervals. 
  
\medskip
For the next proposition, set $D=D_{\ref{cor:two-flows}}(K)$ and  $\kappa=\kappa_{\ref{cor:two-flows}}(K)$. 

\begin{prop}\label{union of three flows}
Assume that $K\ge K_0$ and $\X$ satisfies the uniform  $\kappa$-flaring condition. 
Let $v_0, v_1, v_2$ be vertices of $T$ such that $v_0\in \llbracket v_1,v_2\rrbracket$. 
We will assume that each subinterval 
$J_i= \llbracket v_0,v_i\rrbracket$, $i=1,2$, is special.  
 Then the $D$-neighborhood (taken in $X$) 
$$
Z:=N_D(Fl_K(X_u) \cup  Fl_K(X_v)\cup  Fl_K(X_w))$$ 
(equipped with the path-metric induced from $X$) is $\delta=\delta_{\ref{union of three flows}}(K)$-hyperbolic and 
$\eta=\eta_{\ref{union of three flows},K}$-properly embedded in $X$.
\end{prop}
\proof Define $Y_i=N_D(Fl_K(X_{v_i}))$,  $Z_i:= N_D(Y_0\cup Y_i)$, $i=0, 1, 2$. 
Thus, $Z=N_D(Y_1\cup Y_2)$ and, furthermore, $N_D(Y_0)=Z_0=Z_1\cap Z_2$ and $Z_0$ separates $Z_1$, $Z_2$ in $X$: 
Every path $c$ connecting a point of $Z_1$ to a point of $Z_2$, has to intersect $Z_0$, see Proposition \ref{prop:flow-prop}(1).

\medskip 
{\bf 1: $Z=Z_1\cup Z_2$ is hyperbolic.} The hypothesis of the proposition implies that both $Z_1$ and $Z_2$ 
satisfy the assumptions of   Corollary \ref{cor:two-flows}. Thus, 
each $Z_i$ is $\delta_{\ref{cor:two-flows}}(K)$-hyperbolic 
and  $\eta_{\ref{cor:two-flows},K}$-properly embedded in $X$. Since $Fl_K(X_{v_0})$ is 
$L_{\ref{cor:mjproj}}(K)$-qi embedded in $X$, 
it is $\la=\la_{\ref{lem:qi-preserves2}}(\delta_{\ref{cor:two-flows}}(K), L_{\ref{cor:mjproj}}(K))$-quasiconvex in 
$Z_1, Z_2$. Hence, $Z_0=N_D(Fl_K(X_{v_0}))$  is $\la+D+2\delta_{\ref{cor:two-flows}}(K)$-quasiconvex in $Z_1, Z_2$.
Theorem \ref{thm:hyp-tree} (for quasiconvex amalgams)  now implies that $Z$ is $\delta=\delta_{\ref{union of three flows}}(K)$-hyperbolic.

\medskip 
{\bf 2: $Z_1\cup Z_2$ is uniformly properly embedded in $X$.}  The proof is similar to that of Proposition \ref{flow modification}(2). 
Take points $x_1\in Z_1, x_2\in Z_2$. Then the separation property mentioned earlier, implies that each geodesic $z_1z_2$ in $X$ has to intersect $Z_0$ at some $z_0$. In particular, $\max(d(z_1, z_0), d(z_0,z_2))\le d(z_1,z_2)$. 

By Proposition \ref{prop:two-flows-easy}, for $i=1,2$, 
$$
d_{Z_i}(z_i,z_0)\le 2D+\eta_{\ref{prop:two-flows-easy},K}(d(z_i,z_0))\le  2D+\eta_{\ref{prop:two-flows-easy},K}(d(z_1,z_2)),
$$
and, therefore,
$$
d_Z(z_1,z_2)\le  \eta(d(z_1,z_2))=  \eta_{\ref{union of three flows},K}(d(z_1,z_2)):= 
4D+ 2 \eta_{\ref{prop:two-flows-easy},K}(d(z_1,z_2)). \qed 
$$

\medskip

Recall that for  subtrees $S\subset T$, we defined flow-spaces $Fl_K(X_S)$, see \eqref{eq:flow-of-subtree}. 
We again assume that $K\geq K_0$,  $\X$ satisfies the uniform $\kappa_{\ref{cor:two-flows}}(K)$-flaring condition and 
set $D=D_{\ref{cor:two-flows}}(K)$.   
 
\begin{thm}\label{thm:geod-flow-step1}
 The flow-space $Fl_K(X_{J})$ of any special interval $J=\llbracket  u,v\rrbracket\subset T$ (equipped with the intrinsic path-metric)  is 
$\delta_{\ref{thm:geod-flow-step1}}(K)$-hyperbolic. 
\end{thm}
\proof  
It suffices to prove uniform hyperbolicity of the  $D$-neighborhood of $Fl_K(X_{J})$ in $X$ 
 (with the path-metric induced from $X$), the claim then will follow from the fact that 
 $Fl_K(X_{J})$ is uniformly qi-embedded in $X$ (Proposition \ref{prop:uniform embeddings for subtree flows}).  

Note that, in view of the assumption of the theorem, 
for any two vertices $t, s\in J$, at least one of the intersections is nonempty: 
$$
Fl_K(X_t)\cap X_s\ne \emptyset, \quad \hbox{ or}\quad Fl_K(X_s)\cap X_t\ne \emptyset
$$
(depending on which distance $d(t,v)$, $d(s,v)$ is larger). Thus, any triple of vertices $v_0, v_1, v_2\in J$ satisfies the assumptions of Proposition \ref{union of three flows}, and, hence,
$$
Z=N_D(Fl_K(X_{v_0})\cup Fl_K(X_{v_1})\cup Fl_K(X_{v_2}))$$
is $\delta_{\ref{union of three flows}}(K)$-hyperbolic and 
$\eta_{\ref{union of three flows},K}$-properly embedded in $X$. Since each flow-space 
$Fl_K(X_{v_i})$ is $L_{\ref{cor:mjproj}}(K)$-qi embedded in $X$, it follows that the $D$-neighbor\-hood 
of the union of any two of these flow-spaces is $\la(K)$-quasiconvex in $Z$. 

We are now in position to apply Corollary \ref{cor:bowditch}. 
For each pair of points $x, y\in {\mathcal Fl}_{K}(X_I)$ we 
define the path $c(x,y)$ in $N_D(Fl_K(X_{J}))$ to be a geodesic  between $x, y$ in  
$$
N_{D}\left(Fl_K(X_{\pi(x)})\cup Fl_K(X_{\pi(y)})\right). 
$$
In view of the uniform proper embeddedness of this union (in $X$) and the uniform hyperbolicity of the triple unions as above, 
this family of paths in $Fl_K(X_J)$ satisfies axioms of Corollary \ref{cor:bowditch}. \qed

\medskip 
We now deal with Step 2 outlined in the introduction to this chapter. 
This step is a rather direct application of Theorem \ref{thm:geod-flow-step1}: 
We apply quasiconvex amalgamation of pairs twice to show that $Fl_K(X_J)$ is hyperbolic.

\begin{prop}\label{cor:geod-flow-step1}\label{prop:geod-flow-step2}
Suppose that $J$ is an interval in $T$, that can be subdivided as a union of three special subintervals $J = J_1\cup J_2\cup J_3$,
$$
J_i=\llbracket v_i, v_{i+1}\rrbracket, i=1, 2, 3. 
$$
Then $Fl_K(X_J)$ (equipped with the intrinsic path-metric)  is 
$\delta_{\ref{cor:geod-flow-step1}}(K)$-hyperbolic. 
\end{prop}
\proof A quick way to argue is to appeal to  Corollary \ref{cor:finite-tree-hyp} since $Fl_K(X_J)$ has a structure of a hyperbolic tree of spaces with the base-tree consisting of four vertices and three edges (forming an interval of length $3$), where the vertex-spaces are $Fl_K(X_{J_i})$'s and the edge-spaces are 
$Fl_K(X_{v_i})$'s. We will give a more explicit proof following the proof of  Corollary \ref{cor:finite-tree-hyp} since  it will also provide us with a description of uniform quasigeodesics in $Fl_K(X_J)$.

\medskip 
We will use the quasiconvex amalgamation (see Section \ref{sec:qcamalgam}) 
{\em twice}: 

a.  We have
$$
Fl_K(X_{J_1}) \cap Fl_K(X_{J_2})= Fl_K(X_{v_2}). 
$$
The intersection is  $L_{\ref{cor:mjproj}}(K)$-qi embedded in $X$, hence, in $Fl_K(X_{I})$, $I=J_1\cup J_2$. Since both 
$Fl_K(X_{J_1}), Fl_K(X_{J_2})$ are $\delta_{\ref{thm:geod-flow-step1}}(K)$-hyperbolic, 
Theorem \ref{thm:hyp-tree} (for quasiconvex amalgams)  implies $\eps_{\ref{cor:geod-flow-step1}}(K)$-hyperbolicity of their union $Fl_K(X_{I})$. 

b. We have 
$$
Fl_K(X_J)= Fl_K(X_I)\cup Fl_K(X_{J_3}),
$$
and 
$$
Fl_K(X_{I}) \cap Fl_K(X_{J_3})= Fl_K(X_{v_3}). 
$$
The intersection is  $L_{\ref{cor:mjproj}}(K)$-qi embedded in $X$, hence, in $Fl_K(X_{I})$, $I=J_1\cup J_2$. Since both 
$Fl_K(X_{I}), Fl_K(X_{J_3})$ are $\delta$-hyperbolic, 
$$
\delta= \max(\delta_{\ref{thm:geod-flow-step1}}(K),\eps_{\ref{cor:geod-flow-step1}}(K)),
$$
Theorem \ref{thm:hyp-tree} (for quasiconvex amalgams)   implies $\delta_{\ref{cor:geod-flow-step1}}(K)$-hyperbolicity of their union $Fl_K(X_{J})$. \qed

\medskip 
The following result is never used afterwards, we include the proof for the sake of completeness of the picture. 

\begin{cor}[Hyperbolicity of the union of two flow-spaces: General case]
\label{cor:union of two hallways-II}
For each $K\geq K_0$, for $D=D_{\ref{cor:two-flows}}(K)$, assuming that  
$\X$ satisfies the uniform $\kappa_{\ref{cor:two-flows}}(K)$-flaring condition, 
the following holds:

If $u,v\in V(T)$ are such that $Fl_K(X_u)\cap Fl_K(X_v)\neq \emptyset$, then 
$$N_{D}(Fl_K(X_u)\cup Fl_K(X_v))$$ is an $L_{\ref{cor:union of two hallways-II}}(K)$-qi embedded, 
$\delta_{\ref{cor:union of two hallways-II}}(K)$-hyperbolic subspace of $X$. 
\end{cor}

\proof Since $Fl_K(X_u)\cap Fl_K(X_v)\neq \emptyset$, the interval $J= \llbracket u, v\rrbracket$ 
splits as the union of two special subintervals, $J_1, J_2$, there is a vertex $w$ in the interval 
$I=\llbracket u,v\rrbracket$ such that $X_w\cap Fl_K(X_u)\neq \emptyset$ 
and $X_w\cap Fl_K(X_v)\neq \emptyset$ (see Lemma \ref{lem:special-sub}). Thus, $N_{D}(Fl_K(X_u)\cup Fl_K(X_v))$ is contained in the $D$-neighborhood 
of the $\delta_{\ref{cor:geod-flow-step1}}(K)$-hyperbolic subspace, $Fl_K(X_I)$. The result now follows from Lemma \ref{lem:QCunion} on unions of quasiconvex subsets of hyperbolic spaces, combined with Proposition \ref{prop:two-flows-easy}.   \qed

\section{Hyperbolicity of flow-spaces of general interval-spaces} 
\label{sec:hyp-inerval-flow} 

This section deals with Step 3 outlined in the introduction to this chapter. Recall that according to Proposition \ref{prop:uniform embeddings for subtree flows}, for every subtree $S\subset T$ and $K\ge K_0$, the flow-space $Fl_K(X_S)$ is 
$L_{\ref{prop:uniform embeddings for subtree flows}}(K)$-qi embedded in $X$. 

\begin{thm}\label{thm:hyperbolicity over segments}
For every $K\geq K_0$, assuming that $\X$ satisfies the uniform $\kappa_{\ref{cor:two-flows}}(K)$-flaring condition, 
for each interval $J= \llbracket u,v\rrbracket \subset T$, the flow-space 
$Fl_K(X_{J})$ is $\delta_{\ref{thm:hyperbolicity over segments}}(K)$-hyperbolic.
\end{thm}
\proof  When $u=v$, the statement is established in Theorem \ref{flow of one vertex space}. 
Therefore, we consider the case of nondegenerate intervals $J$. 

 We apply the Horizontal Subdivision Lemma (Lemma \ref{lem:subdivision}) and 
its corollary (Corollary \ref{cor:subdivision}) to subdivide the interval $J$ into subintervals $J_0,...,J_n$ such that:

(i) Each  $J_i$  is the union of three $K$-special subintervals. 

(ii) Whenever $|i-j|\ge 2$,  the flow-spaces $Fl_K(X_{J_i}), Fl_K(X_{J_j})$ are $L_{\ref{mjproj}}(K)$-Lipschitz cobounded in $X$, hence (taking restrictions of Mitra's projections), in $Fl_K(X_J)$ as well.

\medskip
According to Corollary \ref{cor:geod-flow-step1}, 
each flow-space $Fl_K(X_{J_i})$   is $\delta_{\ref{cor:geod-flow-step1}}(K)$-hyperbolic. 
 
\medskip 
 Finally, we consider the union 
$$
Fl_K(X_{J})=\bigcup_{0\leq i\leq n} Fl_K(X_{J_i}). 
$$ 
By Proposition \ref{prop:flow-prop}(1), 
$$
Fl_K(X_{J_{i-1}})\cap Fl_K(X_{J_i})=Fl_K(X_{u_i}).$$
Since whenever $|i-j|\ge 2$,  the flow-spaces $Fl_K(X_{J_i}), Fl_K(X_{J_j})$ are $L_{\ref{mjproj}}(K)$-Lipschitz cobounded in 
$Fl_K(X_{J})$ and the consecutive intersections  $Fl_K(X_{u_i})$ are $L_{\ref{mjproj}}(K)$-qi embedded in 
$Fl_K(X_{J})$, Theorem \ref{thm:hyp-tree}   applies and 
the flow-space  $Fl_K(X_{J})$ is  $\delta_{\ref{thm:hyperbolicity over segments}}(K)$-hyperbolic. 
 \qed

\section{Conclusion of the proof} \label{sec:conclusion}

In this section we finish the proof of the main result of this book, Theorem \ref{thm:mainBF}. We first prove that flow-spaces $Fl_K(X_S)$ are uniformly hyperbolic, whenever $S$ is a tripod in $T$ (Corollary \ref{cor:finite-tree-hyp}). We then conclude the proof of Theorem \ref{thm:mainBF} by appealing to Corollary \ref{cor:bowditch} one last time 
by constructing a slim combing in $X$ via geodesics in flow-spaces of interval-spaces.

\begin{prop}\label{prop:hyperbolicity over tripds}
 Assume that $\X$ satisfies the uniform $\kappa_{\ref{cor:two-flows}}(K)$-flaring condition. 
Suppose $S=T_b u_1u_2u_3\subset T$ is a tripod\footnote{see Definition \ref{defn:tripod}} 
with the center $b$ and three extremities $u_1, u_2, u_3$.
Then for every $K\geq K_0$,  the flow-space $Fl_K(X_{S})$ is $\delta_{\ref{prop:hyperbolicity over tripds}}(K)$-hyperbolic. 
\end{prop}

\proof The proof is  similar to that of Proposition 
\ref{prop:geod-flow-step2}. The tripod $S$ is the union of three segments (legs) $J_i=\llbracket u_i, b\rrbracket, i=1, 2, 3$, whose pairwise intersections equal $\{b\}$. According to Theorem \ref{thm:hyperbolicity over segments}, each flow-space $Fl_K(X_{J_i})$ is 
$\delta_{\ref{thm:hyperbolicity over segments}}(K)$-hyperbolic. The intersection 
$$
\bigcap_{i=1}^3 Fl_K(X_{J_i})= Fl_K(X_b)
$$
(see Proposition \ref{prop:flow-prop}(1)) is $\delta_{\ref{flow of one vertex space}}(K)$-hyperbolic 
(Theorem \ref{flow of one vertex space}) and $L_{\ref{mjproj}}(K)$-qi embedded in $X$. Therefore, $Fl_K(X_{S})$ has structure of a tripod ${\mathfrak Y}$ of hyperbolic spaces 
$$
Y=Fl_K(X_{S})\to S'=T_b v_1v_2v_3,$$
where the tripod $S'$ has four vertices ($b, v_1, v_2, v_3$) and three edges, all incident to the vertex $b$. 
Namely, the vertex-spaces of ${\mathfrak Y}$ are $Y_b= Fl_K(X_b)$, $Y_{v_i}= Fl_K(X_{J_i})$. The three 
edge-spaces are all isomorphic to $Y_b$ and have natural incidence maps (inclusion maps)  to the vertex-spaces.  
Now, the proposition follows from Corollary \ref{cor:finite-tree-hyp}. \qed

\medskip 
We can now finish our proof of Theorem \ref{thm:mainBF}. We let 
$K=K_0$ and $K_*$ be as Notation \ref{not:K0} and Remark \ref{rem:K*} respectively 
(the constants $\delta'_0, \la'_0, L'_0$ used in 
Notation \ref{not:K0} are defined in Notation \ref{not:edge-space-constants}). 
We assume that the tree of spaces ${\mathfrak X}$ satisfies the uniform $K_*$-flaring condition.

We will once again apply Corollary \ref{cor:bowditch}, with $X_0={\mathcal X}$, the union of vertex-spaces in $X$. 
For each pair of vertices $u,v\in V(T)$, and points $x\in X_u$, $y\in X_v$, we define the path $c(x,y)$ 
in $X$ to be the intrinsic geodesic in the  
flow-space $Y=Fl_K(X_{uv})$. In order to verify the assumptions of Corollary \ref{cor:bowditch}, we observe that condition (a1) follows from the fact that each 
$Y$ is $L=L_{\ref{prop:uniform embeddings for subtree flows}}(K)$-qi embedded in $X$. 
We also conclude that each path $c(x,y)$ is an $L$-quasigeodesic in $X$. 
Condition (a2) is immediate from Proposition \ref{prop:hyperbolicity over tripds}. \qed

\medskip

\begin{cor}
Let $\X=(\pi: X\to T)$ be a tree of hyperbolic spaces (satisfying Axiom {\bf H}). Then the following conditions are equivalent:

1. $X$ is hyperbolic. 

2. $\X$ satisfies the uniform $\kappa$-flaring condition for all $\kappa\in [K_0, K_*]$. 

3. $\X$ satisfies the Bestvina-Feighn exponential flaring condition for all $\kappa\ge 1$.

4. Carpets in $\X$ are uniformly hyperbolic. More precisely, there exists a function $\delta(K, C)$, 
such that each $(K,C)$-narrow carpet in $\X$ has $\delta(K,C)$-hyperbolic total space. 

5. Ladders in $\X$ are uniformly hyperbolic. More precisely, there exists a function $\delta(K, D, E)$, 
such that each $(K,D,E)$-ladder in $\X$ has $\delta(K,D,E)$-hyperbolic total space. 
\end{cor}
\proof The implication 2 $\RA$ 1 is the content of Theorem \ref{thm:mainBF}. The converse implication 1 $\RA$ 2 is Lemma \ref{lem:hyp->uniform flaring}. In Proposition 
\ref{hyp to lin flaring} we proved the implication 1 $\RA$ 3, while the implication 3 $\RA$ 2 is proven in Lemma \ref{lem:exp}. Thus, we obtain the equivalence of 1, 2 and 3.

In order to establish equivalence of 1 with 4 and 5, observe that the uniform flaring condition is defined in terms of separation properties of pairs of $K$-sections over intervals in $T$ and every such pair of sections is contained in a 
carpet, while each carpet is contained in a ladder. Applying the implication 1 $\RA$ 2 to trees of spaces which are carpets and ladders, we conclude that uniform hyperbolicity of carpets/ladders implies the uniform $\kappa$-flaring condition for $\X$ for all $\kappa\ge 1$. \qed

\chapter{Description of geodesics}\label{ch:description-of-geodesics} 

\section{Inductive description} \label{sec:inductive} 

Let $\pi: X\to T$ be a tree of hyperbolic spaces (satisfying Axiom {\bf H}) with hyperbolic total space $X$.  
We can now give a description of geodesics in $X$, more precisely, of uniform quasigeodesics. This description is inductive/hierarchical.  The basis of induction is the fact that each $K$-qi leaf in $X$  (i.e. the image of a $K$-qi section of $\pi: X\to T$ over a geodesic segment in $T$)  
defines a uniform {\em horizontal} quasigeodesic in $X$.  Such quasigeodesics present one type of building blocks of geodesics in $X$. The second building block consists of {\em vertical quasigeodesics in $X$}: Such quasigeodesics are certain intrinsic geodesics in vertex-spaces $X_v$ of $X$. While $X_v$'s (typically) are not quasiisometrically embedded in $X$, {\em some} geodesics in $X_v$'s nevertheless are uniform  quasigeodesics in $X$, namely, ones satisfying the {\em small carpet condition}, see Proposition  \ref{prop:uncarpeted geodesics} below. We will see that general geodesics in $X$ are uniformly Hausdorff-close to alternating 
concatenations of horizontal and vertical uniform quasigeodesics.  In Section \ref{sec:brdy-acyl-case} we, furthermore, give a simpler description of  uniform quasigeodesics in a more limited class of trees of spaces, namely, {\em acylindrical} trees of spaces. 

There are several basic classes of subtrees of spaces in $\X$, which are used in description of geodesics in $X$. All of these are 
special cases of {\em semicontinuous families} (of subsets of vertex-spaces) in $\X$, see Chapter \ref{ch:4 classes}, all are uniformly 
quasiconvex subsets of $X$. Here is the list of these classes of subspaces, listed in order of increase of the complexity of their definitions:

\begin{itemize}
\item Carpets $\A\subset \X$. 

\item Metric bundles $\HHH\subset \X$.

\item Ladders $\L\subset \X$. 

\item Flow-spaces of vertex-spaces ${\mathfrak Fl}_K(X_u)\subset \X$. 

\item Flow-spaces ${\mathfrak Fl}_K(X_S)$, where $S\subset T$ is a subtree and $X_S= \pi^{-1}(S)$. 

\item Flow-spaces of metric bundles ${\mathfrak Fl}_K(\HHH)\subset \X$. 

\end{itemize}

  \begin{figure}[tbh]
\centering
\includegraphics[width=70mm]{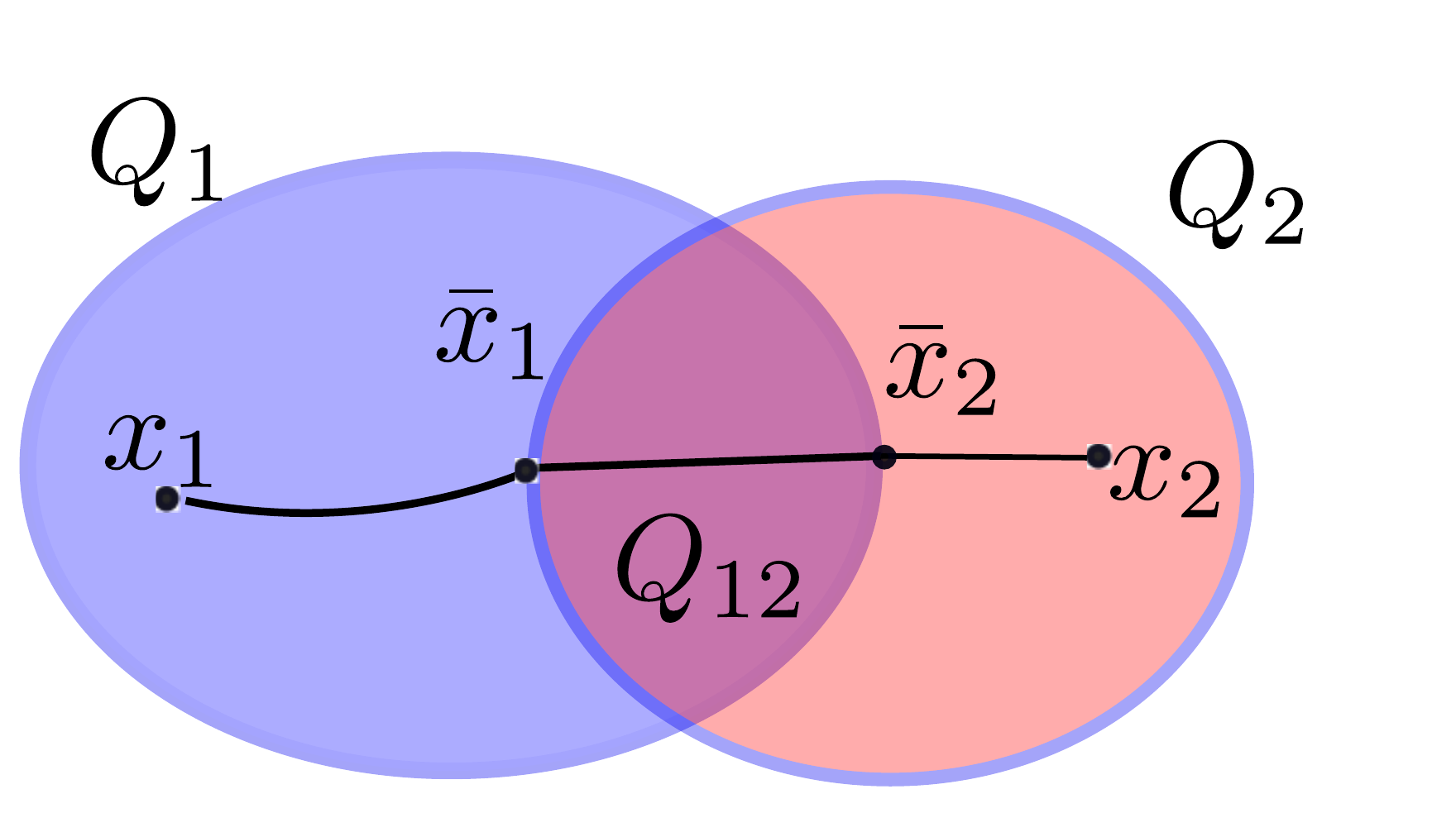}
\caption{Geodesics in amalgams}
\label{fig:amalgam}
\end{figure}

We first describe geodesics in carpets, then use those to describe geodesics in ladders, use those to describe geodesics in flow-spaces of vertex spaces. At the same time, flow-spaces of metric bundles ${Fl}_k(\HHH)$ are uniformly quasiconvex subsets of certain flow-spaces ${Fl}_K(X_u)$ (form some $K\ge k$) 
and, hence, we do not give a separate description of geodesics in the former. After Step I.4, our description of geodesics in $X$  only uses geodesics in flow-spaces ${Fl}_K(X_u)$ as building blocks. A key feature of  flow-spaces is that they are 
uniformly quasiconvex in $X$ (unlike vertex-spaces themselves). Each flow-space is itself a tree of spaces, $Fl_K(X_u)\to S_u$, where $S_u$ is a  subtree in $T$. The intersection pattern of the subtrees $S_u$ (encoded in the  flow-incidence graph $\Ga_K$) is discussed in Section \ref{sec:flow-incidence graph}; it will provide a guide for describing geodesics in $X$ in terms of geodesics in flow-spaces $Fl_K(X_u)$.

We now describe (inductively) geodesics in $X$. 
Except for the two initial steps, the rest is a repetitive use of one of the following constructions:

(a) {\bf Quasiconvex amalgams:  Section \ref{sec:qcamalgam}.}  Given two $L$-qi embedded  subsets $Q_1, Q_2$ in a $\delta$-hyperbolic space $Y$, with nonempty $L$-qi embedded intersection 
$Q_{12}=Q_1\cap Q_2$, the union 
$$
Q= Q_1\cup Q_2
$$ 
is a {\em quasiconvex amalgam} of $Q_1, Q_2$.

  \begin{figure}[tbh]
\centering
\includegraphics[width=85mm]{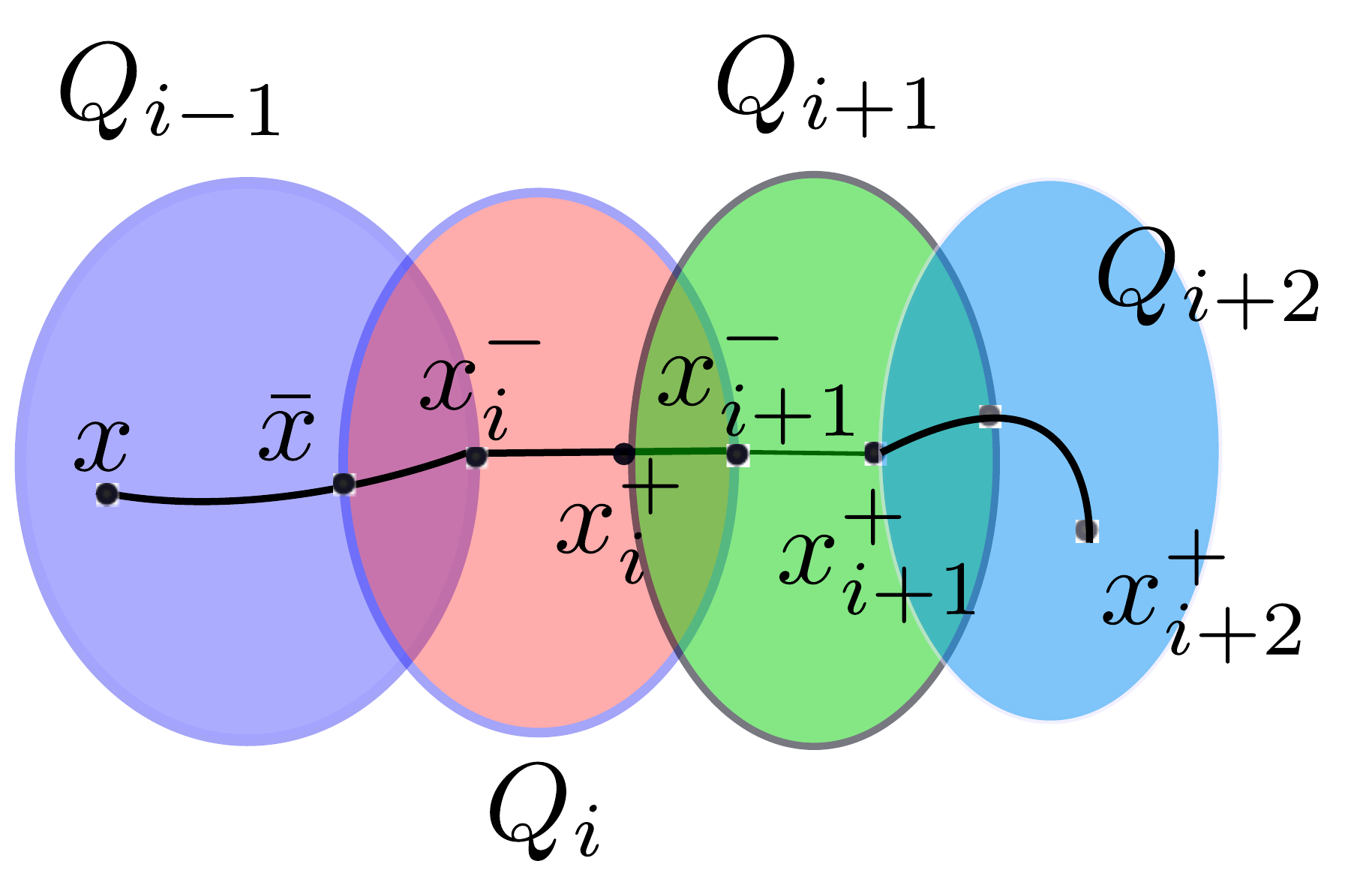}
\caption{Geodesics in chain-amalgams}
\label{fig:chain}
\end{figure}

\medskip 

(b) {\bf Pairwise cobounded quasiconvex chain-amalgamation.} We refer to Section \ref{sec:qcamalgam} for a detailed definition and description of uniform quasigeodesics in unions $Q=Q_0\cup Q_1\cup  ... \cup Q_{n}$ defining 
quasiconvex chain-amalgamation. Such amalgamations (with $n\ge 2$) will be used just in two instances in our proof. Briefly, 
paths $c$ in $Q$ are alternating concatenations of geodesics in $Q_{i,i+1}$'s and $Q_i$'s connecting points $x_i^+, x_{i+1}^+\in Q_{i,i+1}$ and $x^-_i, x^+_i\in Q_i$ where latter pairs  (up to a uniformly bounded error $C$) realize the minimal distance between $Q_{i-1,i}, Q_{i,i+1}$ in $Q_i$. See Figure \ref{fig:chain}. 

\medskip

We now begin the inductive description of uniform quasigeodesics. 
Regarding constants $C, D, E, K$ appearing below:  Ultimately, we will take $K=K_*$, 
$C=M_{\bar{K}}$, $D= D_{\ref{prop:existence-of-tripod-ladders}}$, 
$E=E_{\ref{prop:existence-of-tripod-ladders}}$. However, for instance, the description of paths in carpets works for all $K\ge 1$ and all $C\ge 0$, etc. 

\medskip 
{\bf Part 0: Geodesics in $K$-qi sections of  $\X=(\pi: X\to T)$.} 
The basis for the entire description of uniform quasigeodesics in $X$ is the fact that each $K$-qi leaf in $X$ is a $K$-quasigeodesic in $X$. 

\medskip 
{\bf Part I: Geodesics in flow-spaces of vertex-spaces $Fl_K(X_v)$.} 

These geodesics (or, rather, uniform quasigeodesics) are described in three steps.

\medskip 
{\bf Step I.1.}  {\bf Quasigeodesics in carpets: Section 
\ref{sec:Hyperbolicity of carpets}, especially, Proposition \ref{prop:easy-one}.}

  \begin{figure}[tbh]
\centering
\includegraphics[width=80mm]{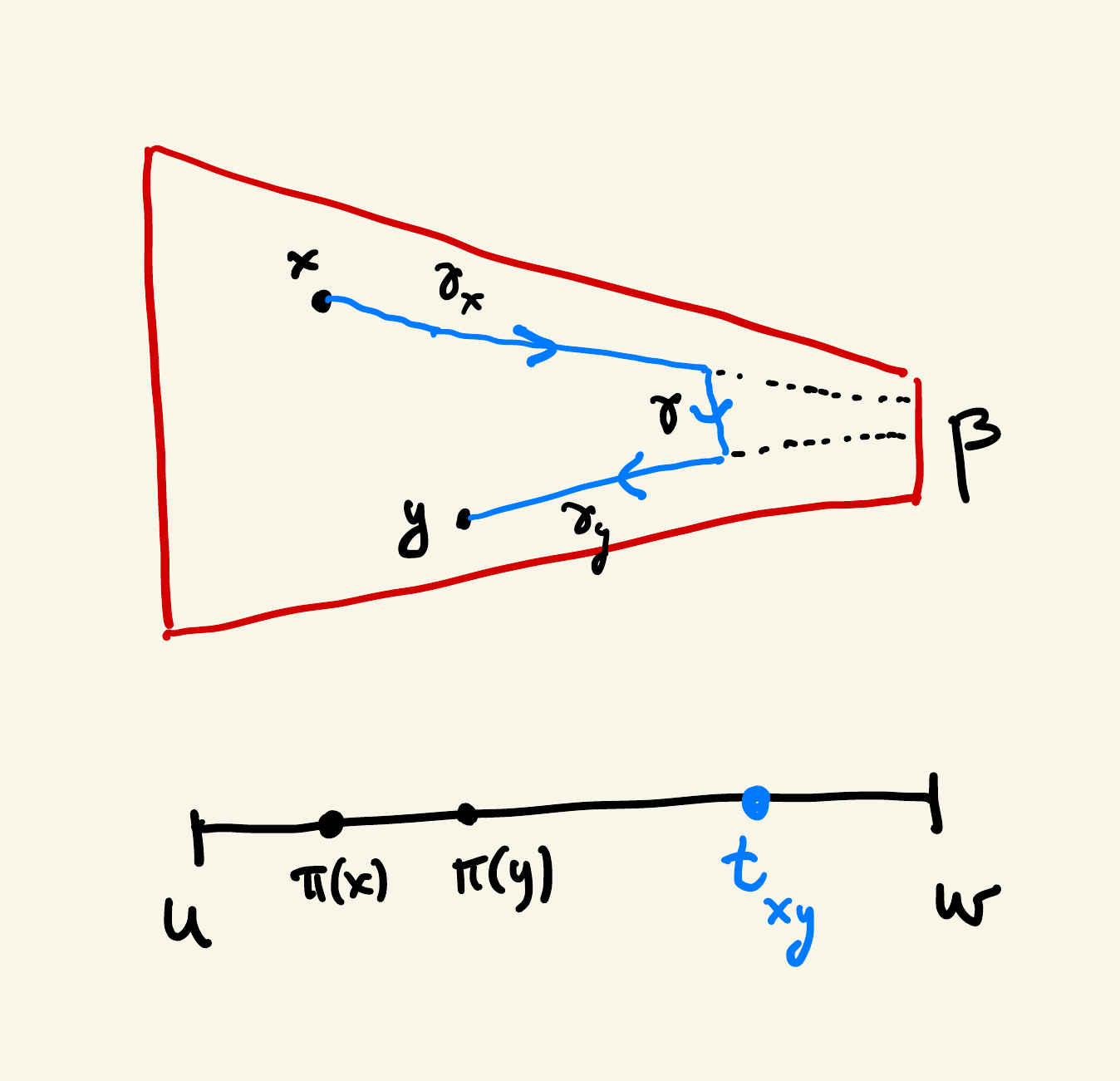}
\caption{Geodesics in carpets}
\label{fig:carpet-paths}
\end{figure}

\medskip

 Let ${\mathfrak A}=(\pi: A\to J)\subset {\mathfrak X}$ be a $(K,C)$-carpet over an interval 
$J=\llbracket u, w\rrbracket$, 
such that the end $A_w$ of  ${\mathfrak A}$ over $w$ is $C$-narrow (i.e. 
is a geodesic of length $\le C$ in $X_w$).

Then for $x, y\in \AA$  we consider $K$-sections $\ga_x, \ga_y$ over 
subintervals $\llbracket \pi(x), w\rrbracket$, $\llbracket \pi(y)w\rrbracket$ in $J$.

Let $t_{xy}\in  \llbracket w, u\rrbracket$ denote the supremum of
$$
\{ t\in \llbracket w,u\rrbracket : d_{X_t}(\ga_x(t), \ga_y(t))\le M_{K}\}. 
$$
Then an  $L_{\ref{prop:easy-one}}(K,C)$-quasigeodesic $c(x,y)$ in $A$ is defined as the concatenation of 
the section $\ga_x$ restricted to $\llbracket \pi(x), t\rrbracket $ with the vertical segment 
$\gamma=[\ga_{x}(t) \ga_y(t)]_{X_t}$, followed by the concatenation with the restriction of the section $\ga_y$ to the subinterval $\llbracket t, \pi(y)\rrbracket $. See figure \ref{fig:carpet-paths}.

\medskip 
{\bf Step I.2.} {\bf Quasigeodesics in carpeted ladders with narrow carpets: Section \ref{sec:narrow ladders}, especially 
Proposition \ref{small-ladder}.}

For a vertical geodesic segment $\al\subset X_u$, we consider a {\em carpeted ladder}, a 
$(K,D,E)$-ladder ${\mathfrak L}={\mathfrak L}_K(\al)=(\pi: A\to J)$,which contains a $(K,C)$-carpet 
${\mathfrak A}={\mathfrak A}(\al')$, where $\al'\subset \al$  is a subsegment of length 
$\ge length(\al)- M_{\bar{K}}$, where $\bar{K}$ is defined by
$$
\bar{K}:=K_{\ref{bundle-proj}}(\delta_0, K, K). 
$$  
The definition of paths $c(x,y)$ connecting points $x, y\in \LL$ is a 2-part process.

{\bf Part a: Retraction $\rho$ and paths $c_x$.} 
In Section \ref{sec:retractions} we defined a retraction $\rho: L_K(\al)\to A$: This retraction is uniformly close to 
the nearest-point projection, see Remark \ref{carpet-rem}. This retraction plays critical role in the definition of our combing of 
$L_K(\al)$. Below is a review of the definition of the retraction in terms of the structure of the tree of spaces. 

For a point $x\in \LL_v$, let $\ga_x$ be a canonical $K$-leaf in $\LL$ connecting $x$ to $\al$: Such leaves are a part of the definition of a ladder. Ultimately, which leaf one takes does not matter and the paths $c_x$ change only uniformly bounded amount if one makes a different choice. Let $t=t_x\in \llbracket u, v\rrbracket$ be the vertex farthest from $u$ such that 
$\pi(x)\in J$ and there exists a point 
$\tilde{x}\in \ga_x(t)$ for which 
$$
d_{X_t}(\tilde{x}, A_t)\le M_{\bar{K}}. 
$$
(It is possible that $t=u$ and $\tilde{x}\in A_u=\al$.)  Then define a path $c_x$ connecting $x$ to $\rho(x)=\bar{x}$ and equal to the concatenation 
$$
\ga_{x,\tilde{x}}\star [\tilde{x} \bar{x}]_{X_{t}}.$$
Here $t=t_x$, and 
$$
\ga_{x,\tilde{x}}=  \ga_x|_{\llbracket v, t\rrbracket }, 
$$
is  the subpath of $\ga_x$ connecting $x$ to $\tilde{x}$, while $\bar{x}\in A_t$ is a nearest-point projection of  
$\tilde{x}$ to $A_t$ in the vertex-space $X_t$. Since $c_x$ is uniformly Hausdorff-close to $\ga_{x,\tilde{x}}$, the 
point $\rho(x)$ essentially determines the path $c_x$.

\medskip

  \begin{figure}[tbh]
\centering
\includegraphics[width=60mm]{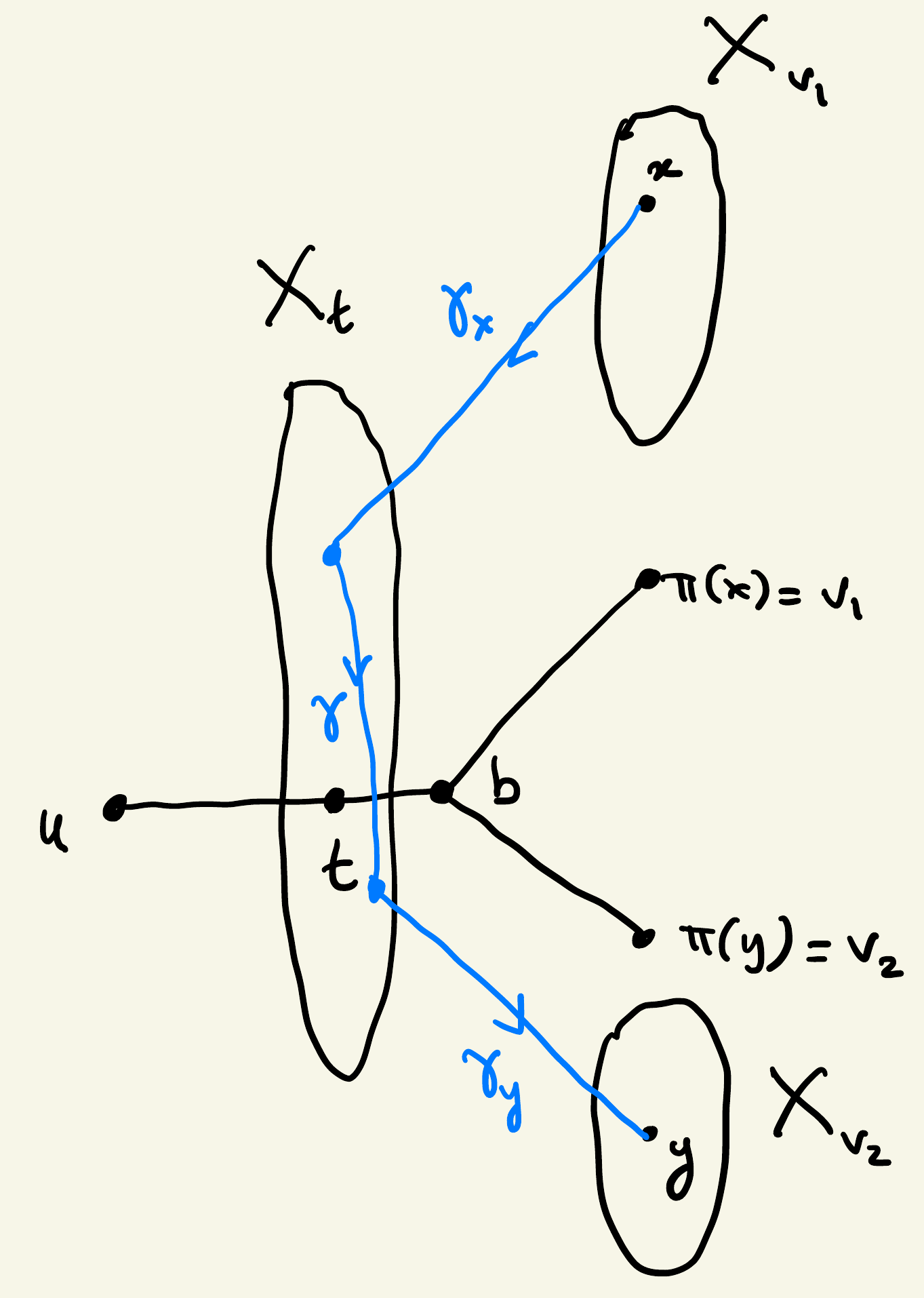}
\caption{Geodesics in ladders: Type 1}
\label{fig:type-1}
\end{figure}

{\bf Part b: Paths $c(x,y)$.} For $x,y \in \LL$ we let $b=b_{xy}$ be the center of the triangle $\triangle u\pi(x)\pi(y)$. The 
path $c(x,y)$ in $L_K(\al)$ connecting $x$ to $y$ is defined as follows.

{\bf Paths of type 1:} There exists $t\in V(\llbracket \pi(x), \pi(\bar{x})\rrbracket \cap \llbracket \pi(y), \pi(\bar{y})\rrbracket)\subset V(\llbracket u,b\rrbracket)$ such that 
$$
d_{X_t}(\gamma_x(t), \gamma_y(t))\leq M_{\bar{K}},$$
i.e. the paths $\ga_x, \ga_y$ ``come sufficiently close'' in some common vertex-space. 

Then let $t=t_{x,y}$ be the maximal vertex in $\llbracket u,b\rrbracket $ with this property. Then define $c(x,y)$ to be the concatenation of the portions of $\gamma_x$ and (the reverse of) $\gamma_y$ over $\llbracket t, \pi(x)\rrbracket $ and $\llbracket t, \pi(y)\rrbracket $ respectively with the subsegment  $\gamma\subset L_{t}$ joining their end-points. See Figure \ref{fig:type-1}.

\begin{figure}[tbh]
\centering
\includegraphics[width=70mm]{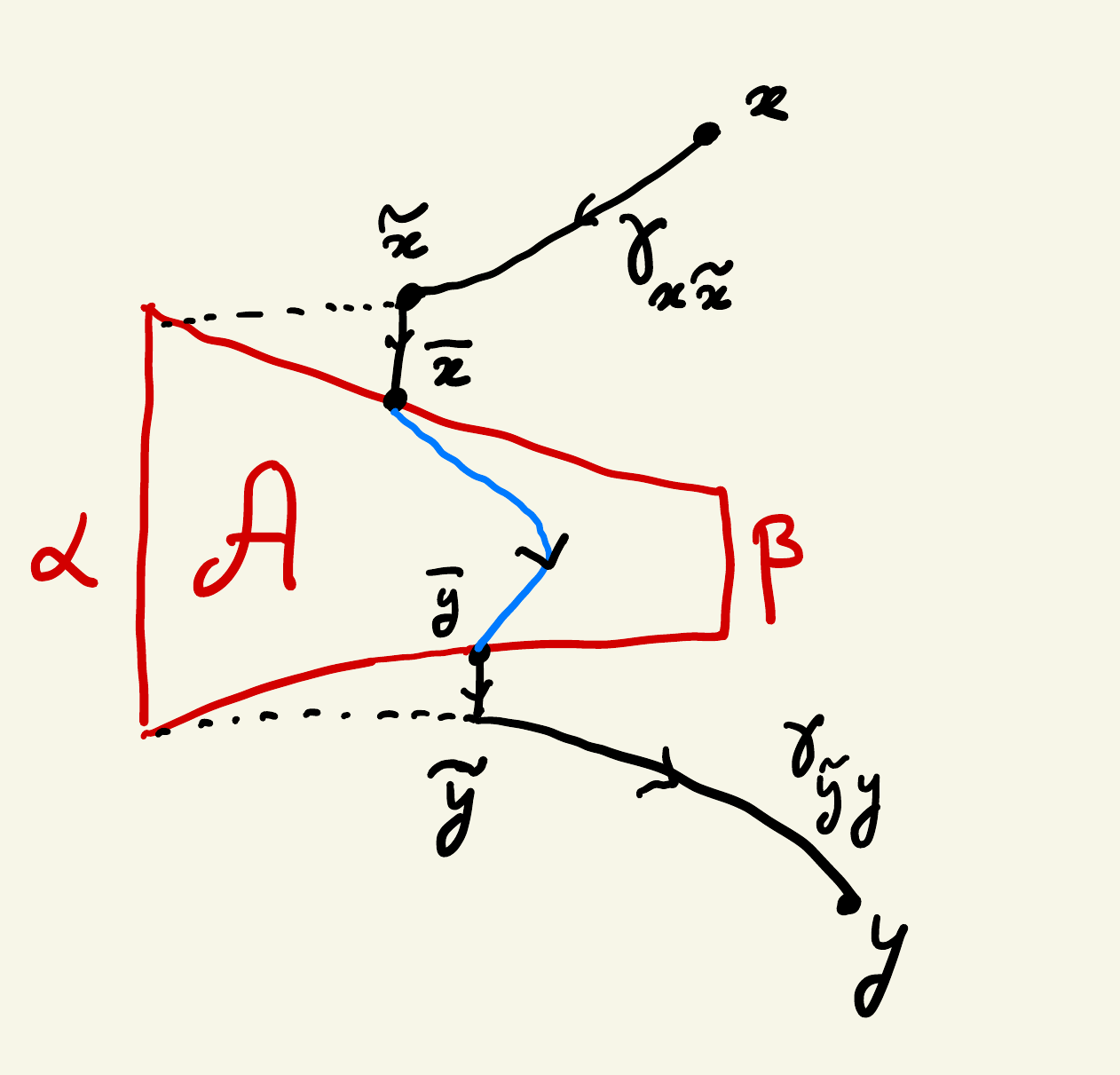}
\caption{Geodesics in ladders: Type 2}
\label{fig:type-2}
\end{figure}

\medskip 
{\bf Paths of type 2:} Suppose type 1 does not happen. Then define $c(x,y)$ to be the concatenation of $c_x$ and the reverse of 
$c_y$ with a geodesic in $A$ connecting $\rho(x)$ to $\rho(y)$. See Figure \ref{fig:type-2}. 

\begin{rem}
The geodesic in $A$ here is taken for granted, see Step I.1. In fact, instead of a geodesic one should take a path $c(\rho(x), \rho(y))$ in $A$ defined in Step I.1. Such inductive arguments will be common in what follows and we will simply say ``geodesic.'' 
\end{rem}

One more thing of importance is that the vertical segments 
$$
 [\tilde{x}\rho(x)]_{X_{u_x}}, \quad [\rho(y)\tilde{y}]_{X_{u_y}} 
$$
have length $\le M_{\bar{K}}$, hence, for paths of both types we are essentially concatenating at most three paths, two of which  are  horizontal subpaths of $c_x, c_y$ and one which lies in the carpet $A$. (The latter is also essentially a concatenation of two horizontal paths.)

 \medskip 
{\bf Step I.3.} {\bf Quasigeodesics in general $(K,D,E)$-ladders $\L(\al)$: Section \ref{sec:hyperbolicity-of-ladders}, especially Proposition \ref{vertical subdivision}.}

In this case, uniform quasigeodesics are constructed via {\em cobounded quasiconvex amalgamation}, 
amalgamating  ladders $\L_K(\al_1),..., \L_K(\al_{n+1})$ along uniformly pairwise cobounded uniformly quasiconvex subsets 
$\Si_1=\Si_{p_1}$ ,..., $\Si_{n-1}=\Si_{p_{n-1}}$, which are canonical $K$-qi sections of these ladders.  

We consider a ladder $\L=\L_K(\al)$ of a vertical geodesic $\al=[pq]_{X_u}$ and 
a pair of points $x, x'\in L_K(\al)$, connected to the end-points $p, p'$ of $\al$ by $K$-qi leaves in $L_K(\al)$.

In Proposition \ref{vertical subdivision} we prove the existence of a {\em Vertical Subdivision} of every vertical geodesic 
segment $\al\subset X_u$, $\al=[p p']_{X_u}$ into subsegments 
$$
\al_1=[p_1 p_2]_{X_u},..., \al_{n}=[p_{n} p_{n+1}]_{X_u}, p_1=p, p_{n+1}=p', 
$$
such that:

(a)  For each $i$, the segment $\al_i$ defines a {\em carpeted ladder} $\L^i=\L_K(\al_i)$ which is a 
$(K,D,E)$-subladder in $\L$ containing a $(K,C)$-narrow  carpet $\A_K(\al'_i)$ (for  $C= M_{\bar{K}}$). 

(b) The canonical sections $\Si_{i}, \Si_{{i+1}}$ (in $\L$) through the points $p_i, p_{i+1}$ are, respectively, the bottom and the top sections of $\L^i$, so that
$$
L^{i-1}\cap L^{i}=\Sigma_i
$$  

(c) The sections $\Si_{i}, \Si_{i+1}$ are $B_{\ref{vertical subdivision}}(K,C)$-cobounded unless $i=n$. 

(d) For all $i$, 
$$
0\le l_i- \length(\al'_i)\le M_K\le M_{\bar{K}}.$$

(e) The top and bottom sections $\Si_{p}, \Si_{p'}$ of $\L$ pass through the points $x, x'$.

\medskip
In particular, the ladder $L$ is a quasiconvex amalgam of its subladders $L^i$ with pairwise intersections $L^{i-1}\cap L^i=\Sigma_i$.

  \begin{figure}[tbh]
\centering
\includegraphics[width=80mm]{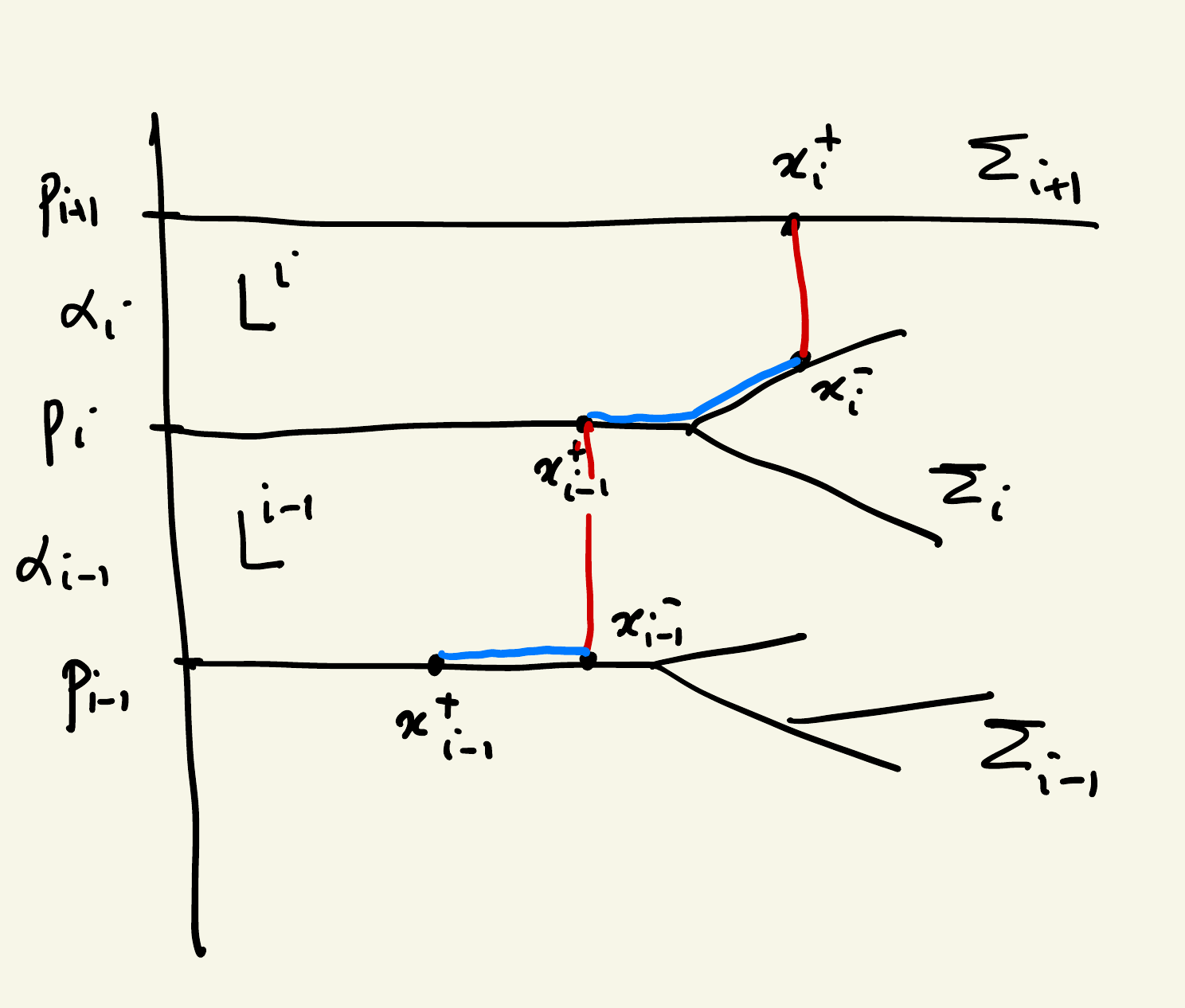}
\caption{Geodesics in general ladders: Vertical subdivision}
\label{fig:figure32}
\end{figure}

As a part of the proof of the Vertical Subdivision Proposition (Proposition \ref{vertical subdivision}), we identify (up to a uniformly bounded error) 
the nearest points $x_{i}^-\in \Si_{i}, x^+_i\in \Si_{i+1}$ between these cobounded subsets. 
Namely, by the construction, the ladder $\L^i$ contains the carpet $\A(\al'_i)$ over an interval 
$\llbracket u, w_i \rrbracket$, with the narrow end 
$$
A_{w_i}= [x_{w_i} y_{w_i}]_{X_{w_i}}.$$
Then (up to a uniformly bounded error), 
$$
x^-_i=x_{w_i}.
$$
The description of the point $x_i^+$ is more complicated, see Lemma \ref{lem:top-projection}.

The paths $c(x_{i}^-, x_i^+)$ are defined according to Step I.2 and the paths $c(x_{i}^+, x_{i+1}^-)$ are defined as in Part 0 (by using 
$K$-qi leaves in the $K$-qi section $\Si_i$).  Lastly,  according to Theorem \ref{thm:hyp-tree}, 
$$
c(x,x')= c(x, x_1^+) \star c(x_1^+, x_2^-)\star ... \star c(x^-_{n}, x^+_n) \star c(x^+_n, x'), 
$$
where $c(x, x_1^+)$, $c(x^+_n, x')$ are the uniformly quasigeodesic paths in carpeted ladders $\L^1$ and $\L^n$ respectively as  
defined in Step I.2.

\medskip 

{\bf Conclusion of Part I:}  {\bf Quasigeodesics in flow-spaces $Fl_k(X_u)$ of vertex-spaces.}  
 Any two points $x, x'\in Fl_k(X_u)$ 
belong to a common $(K,D,E)$-ladder $\L_{x,x'}= \L(\al)$, where $\al$ is a certain geodesic in $X_u$ and $K, D, E$ depend on $k$.  The ladder itself is uniformly close to a 
uniformly quasiconvex subset of $Fl_k(X_u)$. Hence, Step I.3 yields a description of uniform quasigeodesics $c(x,x')$ connecting arbitrary points $x, x'\in Fl_k(X_u)$ (the path $c(x,x')$ is the projection to $Fl_k(X_u)$ of a geodesic in $L(\al)$ connecting $x, x'$). 

\medskip
{\bf Part II:} {\bf Connecting points  in flow-spaces $Fl_k(X_J)$ for special intervals $J\subset T$.} 
An interval $J=\llbracket u , v\rrbracket\subset T$ is said to be {\em special} (more precisely, $k$-{\em special}) if some $w\in \{u, v\}$ has the property that $J\subset \pi(Fl_k(X_w))$; such $w$ is called a {\em special vertex} in $J$. 
The flow-space $Fl_k(X_J)$ then is called a {\em special flow-space} (more precisely, $k$-{\em special} flow-space). Similarly, an interval $J\subset T$ is {\em $k$-semispecial} if it is the union of two special intervals meeting at a vertex. Accordingly, for such $J$, 
the flow-space $Fl_k(X_J)$ is {\em $k$-semispecial}. Uniform quasigeodesics in $Fl_k(X_J)$ in this setting are described using pairwise quasiconvex amalgams of hyperbolic spaces.

\medskip
{\bf Step II.4. Uniform quasigeodesics in special flow-spaces:  
Section \ref{sec:generalized hallways}, especially Proposition \ref{flow modification}, Corollary \ref{cor:two-flows}.}
 
We assume that $J=\llbracket u , v\rrbracket\subset T$ is a $K$-special interval in $T$. 
The flow-spaces $Fl_K(X_u)$, $Fl_K(X_v)$ are uniformly quasiconvex in $X$, hence, a uniform $D$-neighborhood of this union 
in $X$ is hyperbolic. The constant $D$ is defined in Corollary \ref{cor:two-flows}. 
On this step, we describe uniform quasigeodesics in $X$ connecting points of $Fl_K(X_u)\cup Fl_K(X_v)$ and uniformly close to 
such a union.  

We first define a certain metric bundle $\HHH$ over $J$. We then combine the generalized flow-space $Y_0=Fl_R({\mathfrak H})$ 
of this bundle (for some $R=R(K)$) with the flow-spaces $Fl_K(X_u)$, $Fl_K(X_v)$ to get 
``modified flow-spaces'' $Y_1= Fl_K(X_u)\cup Fl_R({\mathfrak H})$,  $Y_2= Fl_K(X_v)\cup Fl_R(\HHH)$. 
The union 
$$
U=U_1\cup U_2= N_D(Y_1)\cup N_D(Y_2)
$$
is uniformly qi embedded in $X$,  while the intersection $U_1\cap U_2$ is uniformly quasiconvex and uniformly 
Hausdorff-close to $Y_0$. Moreover,  $U_1, U_2$ are uniformly Hausdorff-close to subsets in  
$Z_1=Fl_{K_1}(Y_1), Z_2=Fl_{K_1}(Y_2)$ for some (computable $K_1\ge K$). 

Hence, we are in the setting of the quasiconvex amalgamation of a pair (Section \ref{sec:qcamalgam}). To connect $x\in Fl_K(X_u), y\in Fl_K(X_v)$ 
by a uniform quasigeodesic uniformly close to $Fl_K(X_u)\cup Fl_K(X_v)$ we proceed as in Section \ref{sec:qcamalgam}:  

We first project (in $Fl_{K_1}(Y_1), Fl_{K_1}(Y_2)$ respectively) the points $x, y$ to points 
$\bar{x}, \bar{y}\in Y_0= Fl_K({\mathfrak H})$ (one can also use projections in $U$ or in $X$, the 
difference will be uniformly bounded) 
and then connect  $\bar{x}, \bar{y}$ in $Y_0$. The concatenations 
$$
[x \bar{x}]_{Z_1}\star [\bar{x} \bar{y}]_{Y_0} \star [\bar{y} y]_{Z_2}  
$$
are uniform quasigeodesics in $X$ and are uniformly close to $Fl_K(X_u)\cup Fl_K(X_v)$.

\medskip

{\bf Step II.5: Uniform quasigeodesics in semispecial  flow-spaces: 
Proposition \ref{union of three flows} and 
Theorem \ref{thm:geod-flow-step1}.}  

Suppose now a semispecial interval $S=\llbracket u, v\rrbracket= I\cup J$, where both 
$I=  \llbracket u, w\rrbracket, J= \llbracket w, v\rrbracket$  are special intervals. Then for $s\in V(I)$ and $t\in V(J)$,  
uniform quasigeodesics in $Fl_K(X_S)$ connecting points $x\in Fl_K(\llbracket s, w\rrbracket), y\in Fl_K(\llbracket w, t\rrbracket)$ 
are described as follows. The union 
$$
Fl_K(\llbracket s, w\rrbracket)\cup Fl_K(\llbracket w, t\rrbracket)
$$
is a quasiconvex amalgam over the flow-space $Fl_K(X_w)$. Therefore, according to the description of uniform quasigeodesics in quasiconvex amalgams, 
we first project (using nearest-point projections in $Z_1=N_D(Fl_K(X_s)\cup Fl_K(X_w))$, resp. in  $Z_2=N_D(Fl_K(X_t)\cup Fl_K(X_w))$) points $x$ (resp. $y$) to $\bar{x}\in Fl_K(X_w)$ (resp. $\bar{y}\in Fl_K(X_w)$), and then take the concatenation
$$
[x \bar{x}]_{Z_1} \star [\bar{x} \bar{y}]_{Fl_K(X_w)} \star [\bar{y} y]_{Z_2}.
$$
The first and the last geodesics in this concatenation are from Step II.4, while the middle one is from the conclusion of Part I.

\medskip 
{\bf Step II.6: Uniform quasigeodesics in triple unions of special flow-spaces: 
Corollary \ref{cor:geod-flow-step1}.}  

Suppose that $J$ is an interval in $T$, that can be subdivided as a union of three special subintervals, $J = J_1\cup J_2\cup J_3$,
$$
J_i=\llbracket v_i, v_{i+1}\rrbracket, i=1, 2, 3. 
$$
Set $I:= J_1\cup J_2$. Then geodesics in $Fl_K(X_J)$ are described by applying quasiconvex amalgamation of pairs {\em twice}: Once to 
$$
Fl_K(X_I)= Fl_K(X_{J_1\cup J_2}),
$$
which is the amalgam of $Fl_K(J_1), Fl_K(J_2)$ over $Fl_K(X_{v_2})$, and then once more, to 
$$
Fl_K(X_{I\cup J_3}), 
$$
which is the amalgam of $Fl_K(I), Fl_K(J_3)$ over $Fl_K(X_{v_3})$.

\medskip
{\bf Part III:} {\bf Connecting general points in $X$.} 

\medskip

{\bf Step III.7: Horizontal subdivision.} 

Any two points $x, y\in \XX$ belong to the interval flow-space $Fl_K(X_J)$, where $J= \llbracket u,v\rrbracket \subset T$, 
$x\in X_u, y\in X_v$. Since the flow-space $Fl_K(X_J)$ is (uniformly) quasiconvex in $X$, it suffices to describe a uniform 
quasigeodesic in $Fl_K(X_J)$ connecting $x$ to $y$. The key ingredient of this part is the {\em Horizontal Subdivision Lemma} (Lemma \ref{lem:subdivision}).  
This lemma gives a subdivision of the interval $J$ into subintervals $J_i=\llbracket u_i, u_{i+1}\rrbracket, i=1,...,n$, such that:

1.  $Fl_K(J_i)\cap Fl_K(J_j)=\emptyset$ whenever $|i-j|\ge 2$. 

2.  Each interval $J_i$ is subdivided in three special subintervals. 

This represents $Fl_K(X_J)$ as a (uniformly) pairwise cobounded quasiconvex chain with quasiconvex subsets 
$Q_i= Fl_K(X_{J_i})$, $i=1,....,n$, whose consecutive intersections 
$$
Q_{i-1}\cap Q_{i}= Q_{i-1,i}= Fl_K(X_{u_i})
$$
separate in the union $Fl_K(X_J)$ as required by a cobounded quasiconvex chain-amalgamation. Thus, uniform quasigeodesics in 
$Fl_K(X_J)$ are described according to Section \ref{sec:qcamalgam}, see the discussion of quasiconvex chain-amalgamation given early in this section.

This concludes our description of uniform quasigeodesics between points of $X$.

\section{Characterization of vertical quasigeodesics}  

In this section we use the description of uniform quasigeodesics in $X$ to characterize vertical geodesics in $X$ (i.e. geodesics in 
vertex-spaces $X_u$) which are quasigeodesics in $X$. We assume that $\X= (\pi: X\to T)$ is a tree of spaces satisfying the assumptions of Theorem \ref{thm:mainBF}, equivalently, a tree of hyperbolic spaces with hyperbolic total space $X$. Set $K=K_0$ 
and let $C=M_{\bar{K}}$, $D= D_{\ref{prop:existence-of-tripod-ladders}}$, 
$E=E_{\ref{prop:existence-of-tripod-ladders}}$.  

Suppose that $\al$ is  (finite or infinite) geodesic in $X_u$. We will say that $\al$ satisfies the {\em $R$-small carpet condition} \index{small carpet condition}
 if  
whenever $\al'\subset \al$ is a subsegment which bounds a $(K, C)$-narrow carpet  $\A=\A({\alpha'})\subset \X$, 
$\A= (\pi: A\to \llbracket u, w\rrbracket)$, with fiberwise distances between $top(\A), bot(\A)$ at least $M_K$, then 
we have 
$$d_T(u,w)= \length \pi(A)\le R.$$ 
In view of Corollary \ref{cor:contraction}, for such $\al$'s,  the lengths of subsegments $\al'\subset \al$ bounding $(K, C)$-narrow carpets $\A({\alpha'})$, are uniformly bounded.

\medskip
The main result of this section is that a vertical geodesic $\al$ satisfies the small carpet condition 
if and only if it is quasigeodesic in $X$. More precisely:

\begin{prop}\label{prop:uncarpeted geodesics} 
1. Each vertical geodesic $\al$ satisfying the $R$-small carpet condition is an 
$L_{\ref{prop:uncarpeted geodesics}}(R)$-quasigeodesic in $X$. 

2. If $\al$ fails  the $R$-small carpet condition for all $R$, then $\al$ is not a quasigeodesic in $X$. 
 \end{prop}
\proof 1. Let $\L=\L_K(\al)$ denote a $(K,D,E)$-ladder in $X$ based on the segment $\al$. 
Since $L_K(\al)$ is uniformly qi embedded in $X$ (see Corollary \ref{cor:ladder-retraction}), 
it suffices to show that $\al$ is a uniform quasigeodesic in  $L_K(\al)$.
We follow the description of uniform quasigeodesics geodesics in $L_K(\al)$ given in Step I.3 in the previous section. 
We subdivide the segment $\al=[pp']_{X_v}=[p_1 p_{n+1}]_{X_u}$ into subsegments 
$$
\al_i= [p_i p_{i+1}]_{X_u}, i=1,...,n. 
$$
Each $\al_i$ (except for possibly $i=n$) contains a subsegment $\al_i'$ such that
$$
M_K\le  \length(\al_i)- M_K\le  \length(\al'_i)$$
and  there exists a $(K,C)$-carpet $\A^i= \A_{(K,C)}(\al_i)\subset \L$; 
the bottom and the top of $\L^i$ are $K$-qi sections $\Si_i, \Si_{i+1}\subset \L^i$ of a subladder $\L^i\subset \L$.  Moreover, for 
each vertex $t$ in the interval $\pi(A^i)$,
$$
d_{X_t}(top(A^i)_t, bot(A^i)_t)\ge M_K. 
$$
 
For each $i$ we mark the points $x_i^-\in \Si_{i}, x_{i}^+\in \Si_{{i+1}}$ which (up to a uniformly bounded error) are the 
nearest points in $L^i$ between these subsets and consider geodesics  
$$
\beta_i=[x_i^- x_{i}^+]_{L^i}, 
$$
$$
\ga_i^-:= \ga_{x_{i-1}^+, x^-_i}= [x_{i-1}^+ x^-_i]_{\Si_i}     \subset \Si_{i},
$$
$$
\ga_i^+:= \ga_{x_{i}^+, x^-_{i+1}}= [x_{i}^+ x^-_{i+1}]_{\Si_{i+1}}\subset \Si_{i+1}. 
$$
The path $c(p,p')$ then is (up to a uniformly bounded error) equal the concatenation
$$
c= ... \star \ga_i^- \star \beta_i \star \ga_i^+ \star ...
$$
Our goal is to show that the path $c$ is uniformly Hausdorff-close to $\al$.  
Because the lengths of  projections  $\pi(A_i)$ are uniformly bounded, the geodesics $\beta_i$ are uniformly close to 
the vertical geodesics $\al_i$. For the same reason, the paths $\ga_i^\pm$ are uniformly short as well. Hence, each concatenation 
$$
\ga_i^- \star \beta_i \star \ga_i^+
$$ 
is $r=r_{\ref{prop:uncarpeted geodesics}}(R)$-Hausdorff close to the vertical geodesic segment $\al_i$. 
It follows that $\Hd(\al,c)\le r$ as well.  Now, the first statement of the proposition 
follows from Lemma \ref{lem:quasigeodesic-paths}. 

\medskip 
2. Suppose that $\al$ contains a sequence of subsegments $\al_i$ each bounding a $(K,C)$-carpet 
$\A^i= \A_{(K,C)}(\al_i)= (\pi: A^i\to J_i=\llbracket u, w_i\rrbracket)$, such that 
$$
\lim_{i\to\infty} d_T(u, w_i)=\infty,   
$$
and for each vertex $t$ in the interval $J_i$,
$$
d_{X_t}(top(A^i)_t, bot(A^i)_t)\ge M_K. 
$$
Therefore, the concatenation $c_i$ of the bottom of $\A^i$, the 
narrow end $\beta_i$ and the top of $\A^i$ is an $L_{\ref{prop:easy-one}}(K,C)$-quasigeodesic in $A^i$. Since
$$
d(\al_i, \beta_i)\ge d_T(u, w_i), 
$$
the Hausdorff distances between $\al_i$ and $c_i$ diverge to infinity. Morse lemma and hyperbolicity of $X$ then imply that 
$\al_i$'s cannot be uniform quasigeodesics in $X$. Therefore, $\al$ is not a quasigeodesic in $X$ either. \qed

\section{Visual boundary and geodesics in acylindrical trees of spaces} \label{sec:brdy-acyl-case}

In this section we specialize our discussion of  geodesics in trees of hyperbolic spaces to the  
case of {\em $(M,K,\tau)$-acylin\-drical trees} of hyperbolic  spaces satisfying Axiom {\bf H},  
 see Definition \ref{defn:acylindrical}; the constant $K$ will be taken equal 
 $K=K_*=K_{\ref{thm:mainBF}}(\delta_0,L_0)$ 
 although, many arguments will go through for smaller values of $K$. Besides uniform quasigeodesics we will also describe the ideal boundary of $X$. In the group-theoretic setting, when $X$ is the Cayley graph of the fundamental group of an acylindrical graph of hyperbolic groups with quasiconvex edge-subgroups, 
 such description of the boundary is due to Dahmani, \cite{MR2026551} (who also gave a description in the relatively hyperbolic case).

\begin{lemma}\label{acyl lemma 1}
Suppose $\pi:X\map T$ is a $(M, K,\tau)$-acylindrical tree of hyperbolic  spaces satisfying Axiom {\bf H}
for some $M$ and $\tau$. Then (1) $X$ is hyperbolic. (2) For  subtrees $S\subset T$, the subspaces $X_S$ are uniformly qi embedded in $X$.
\end{lemma} 

\proof (1) Recall (Section \ref{sec:acyl trees}), that  $(M,K,\tau)$-acylindricity implies uniform $K$-flaring, hence, by Theorem \ref{thm:mainBF}, hyperbolicity of $X$. 

(2) Let $v\in V(T)$. Then, by Proposition \ref{prop:uncarpeted geodesics}, every geodesic in $X_v$
is a uniform quasigeodesic in $X$. Thus $X_v$ is uniformly qi embedded in $X$. 
Consequently $X_v$ is uniformly quasiconvex too. Then it follows that $X_S$ is also uniformly quasiconvex in $X$ for the
following reason. Suppose a geodesic segment $\beta$ of $X$ joining a pair of points in $X_S$ is not entirely contained in $X_S$.
Then the closure $\beta_w$ of each connected component of $\beta\setminus X_S$ joins two points of $X_w$ for some $w\in S$.
Since vertex-spaces $X_w$ are uniformly quasiconvex in $X$, the geodesic $\beta_w$ is uniformly close to $X_w$. It follows that 
the entire $\beta$ is uniformly close to $X_S$. 

Finally we know that $X_S$ is uniformly properly embedded in $X$ by Proposition \ref{unif-emb-subtree}.
Hence, $X_S$ is uniformly qi embedded in $X$ by Lemma \ref{lip-proj} and Lemma \ref{qc plus prop emb}.
\qed

\medskip

{\bf Description of quasigeodesics.} Suppose that $\X$ is $(M,K,\tau)$-acylindrical,  $x\in X_u, y\in X_v$, $u,v\in V(T)$. Without loss of generality, $\tau\in \NN$. 
Since, by the above lemma, $X_{\llbracket u, v\rrbracket}$ is uniformly qi embedded in $X$, it suffices to describe 
uniform quasigeodesics connecting $x, y$ in $X_{\llbracket u, v\rrbracket}$. Hence, from now on, we will assume that the tree $T$ is an interval $\llbracket u, v\rrbracket$. 
Let
$$
t_0=u, t_1,...,t_m=v 
$$
denote the consecutive vertices in the interval $\llbracket u, v\rrbracket$. Define $A_{i-1}:=X_{e_it_{i-1}}=f_{e_i t_{i-1}}(X_{e_i})\subset X_{i-1}$ and $B_i:=X_{e_it_{i}}=f_{e_i t_{i}}(X_{e_i})\subset X_{t_i}$
for $1\leq i\leq m-1$. We inductively construct points $x^+_i, x^-_i\in X_{t_i}$ 
for $0\leq i\leq m$ as follows. 

Set $x^{-}_0= x$ and $x^{+}_0=P_{X_{t_0}, A_0}(x)$. Now, suppose that $x^+_i, x^-_i\in X_{t_i}$
are already defined for some $i<m$. Then we define 
$x^{-}_{i+1}$ to be an arbitrary point of $f_{e_{i+1} t_{i+1}}(f^{-1}_{e_{i+1}, t_i}(x^+_i))$ and
$x^{+}_{i+1}:= P_{X_{t_{i+1}}, A_{i+1}}(x^{-}_{i+1})$. We define $x^+_m:= y$.

Lastly, we define the path $\gamma(x,y)$ as the concatenation of the segments $[x^-_i x^+_i]_{X_{t_i}}$ for $0\leq i\leq m$,
and unit segments in $X_{t_i t_{i+1}}$ joining each pair $x^+_i, x^-_{i+1}$ for $0\leq i\leq m-1$.

\begin{prop}\label{acyl finite geodesics}
The path $\gamma(x,y)$ is a uniform quasigeodesic in $X$ joining $x,y$. Moreover, $x^-_m$ is uniformly close to  the nearest point
projection of $x$ to $X_v$. 
\end{prop}
\proof The proof is based on several lemmata. We first prove the claim for constants depending on $m$ and then eliminate this dependence.

\begin{lem}\label{lem:acyl finite geodesics-1}
1. The path $\gamma(x,y)$ is an $L_{\ref{lem:acyl finite geodesics-1}}(m)$-quasigeodesic. 

2. The distance between $x^{-}_m$ and the projection of $x$ to $X_{t_m}$ inside $X_{\llbracket u, v\rrbracket}$ is bounded by 
$D_{\ref{lem:acyl finite geodesics-1}}(m)$. 
\end{lem}
\proof In view of uniform quasiconvexity  in $X=X_{\llbracket u, v\rrbracket}$ of vertex-spaces and of subintervals of spaces $X_{\llbracket t_i, t_j\rrbracket}$, the statements follow from, respectively, Parts 1 and 2 of Lemma \ref{lem:two-projections-2}. \qed 

\medskip
Below we set $R= R_{\ref{proj-geod}}(\delta_0,\la_0)=\la_0+5\delta_0$ and note that $K\ge R+1$.

\begin{lem}\label{lem:acyl finite geodesics-2}
Each pair of vertex spaces $X_u, X_v$ satisfying $d_T(v,w)=m\ge \tau$,  is $C_{\ref{lem:acyl finite geodesics-2}}(m)$-cobounded in $X$. 
\end{lem}
\proof We start with a pair of points $x, z\in X_u$ and inductively project them to points $x^-_i, z^-_i\in X_{t_i}$, $1\le i\le m$, 
using the notation above.

According to Lemma \ref{proj-geod}(1), for each $i$ one of two things happens, where $\delta_0$ is the hyperbolicity constant of $X_{t_i}$ and $\la_0$ is the quasiconvexity constant of $A_i$ in $X_{t_i}$: 

(a) $d_{X_{t_i}}(x^+_i, z^+_i)\le D_{\ref{proj-geod}}(\delta_0,\la_0)= 2\la_0+7\delta_0$. 

(b) $[x^+_i z^+_i]_{X_{t_i}}\subset N_R([x^-_i z^-_i]_{X_{t_i}})$, where, as above, 
$R= R_{\ref{proj-geod}}(\delta_0,\la_0)=\la_0+5\delta_0$.  

If (a) occurs for some $i$, then the distance between $x^-_i, z^-_i\in X_{t_i}$ is uniformly bounded; hence, the distance between 
$x^-_m, z^-_m\in X_{t_m}=X_v$ is uniformly bounded as well, with bound depending on $m$. Suppose, therefore, that (b) holds for every $i$. Then for $k=R+1$, there exists a pair of $k$-sections $\gamma_0, \gamma_1$ over $\llbracket u, v\rrbracket$ such 
that 
$$
\gamma_0(t_i)= x_i^-, \quad \gamma_1(t_i)= z_i^-, i=0,...,m. 
$$
Since, by the assumption, $K\ge k$, the acylindricity condition implies that $d_{X_v}(x_m^-, z_m^-)\le M$. Combining this with Lemma \ref{lem:acyl finite geodesics-1}(2), we obtain a uniform (in terms of $m$) bound on the distance between the projections of $x, z$ to $X_v$.  \qed

\medskip 
We are now ready to prove the proposition. We let $N:= \lfloor \frac{m}{\tau}\rfloor$ and define 
the subintervals 
$$
J_i=\llbracket v_i, v_{i+1}\rrbracket, v_i= t_{i\tau}, i=0,...,N,
$$
in $\llbracket u, v\rrbracket$. These subintervals cover $\llbracket u, v\rrbracket$ except for the subinterval 
$J_{N+1}= \llbracket v_{N+1}, v\rrbracket$ whose length is $<\tau$. Set
$$
x_i:= x^-_{i\tau}, i=0,...,N. 
$$
Since the subspaces $X_{J_i}$ are uniformly quasiconvex in $X$ and each pair $X_{J_{i-1}}, X_{J_{i+1}}$ is uniformly cobounded, it follows from Theorem \ref{thm:hyp-tree} on quasiconvex chain-amalgams that the concatenation
$$
[x_0 x_1] \star ... \star [x_N x_{N+1}] \star [x_{N+1} y] 
$$
is a uniform quasigeodesic in $X$. Furthermore, by Lemma \ref{lem:acyl finite geodesics-1}(1), each path  
$$
\gamma(x_i, x_{i+1}), i=0,...,N, \gamma(x_N,y), 
$$  
is also a uniform quasigeodesic. Thus, the entire path $\gamma(x,y)$ is uniformly quasigeodesic as well. \qed

\begin{defn}
 We shall refer to the (uniform) quasigeodesics of the type described in Proposition \ref{acyl finite geodesics} as {\em
HV (horizontal-vertical) quasigeodesics} in what follows. 
\end{defn}

Up to a uniform error, these HV quasigeodesics describe all finite geodesics in $X$. 
Our next goal is to extend this description to the rays in $X$. We will do so under the extra assumption that $X$ is a proper metric space.  Of course, as before, we also suppose that $\pi:X\map T$ is a $(M,K,\tau)$-acylindrical tree of hyperbolic spaces satisfying Axiom {\bf H}.

 \medskip

Fix $v_0\in V(T)$ and $x_0\in X_{v_0}$. We will  describe (quasi)geodesic
rays in $X$ starting from $x_0$.  
First of all, for every $v\in V(T)$ we fix once and for all an HV uniform quasigeodesic $\gamma_v$ joining to $x_0$ to a point
$x_v\in X_v$ where:

 (1) $x_v$ is uniformly close to a nearest point projection of $x_0$ to $X_v$ as we obtained in 
Proposition \ref{acyl finite geodesics}; 

(2) for each $w\in \llbracket v_0,v\rrbracket$ we have $\gamma_w\subset \gamma_v$. 

One defines $\gamma_v$ by induction on $d(v_0,v)$. Note that for each vertex $w\in  \llbracket v_0,v\rrbracket$, the preimage of $w$ in $\gamma_v$ 
under the projection $\pi$ is an interval. In this situation, we will say that $\gamma_v$ {\em projects monotonically} to $\llbracket v_0,v\rrbracket$.

\medskip 
Armed with this,  we can now describe quasigeodesic rays in $X$.

{\bf Rays of type 1:} Let $v\in V(T)$ and $\xi_v\in \geo X_v$. Let $\alpha_v$ be a geodesic ray in $X_v$ joining $x_v$ to $\xi_v$.
Then $\rho_v:=\gamma_v\star \alpha_v$ is a uniform quasigeodesic in $X$ by Proposition \ref{acyl finite geodesics}.

{\bf Rays of type 2:} On the other hand, suppose that $c=v_0 \eta$ is a geodesic ray in $T$ joining $v_0$ to $\eta\in \geo T$.
Then the uniform quasigeodesic paths $\gamma_{c(n)}$ combine to form a uniform quasigeodesic ray $\gamma_{\eta}$ in $X$ which projects monotonically to $c$.

For the following proof, we recall a notion from \cite[Chapter III.H, p. 429]{bridson-haefliger}: A {\em generalized geodesic ray} in a metric space $Z$ is either a geodesic ray or a finite subinterval.

\begin{prop}\label{ray in acyl tree}
Any (quasi)geodesic ray in $X$ starting from $x_0$ is asymptotic to a quasigeodesic ray of exactly one of the above two types.
\end{prop} 

\proof Suppose that $\rho$ is a geodesic ray in $X$ emanating from $x_0$. 
For each $n\in \NN$ let $\alpha_n=\beta(x_0, \rho(n))$ be the arc-length parametrized HV 
quasigeodesic segment (discussed after Lemma \ref{acyl lemma 1}),  
joining $x_0$ to $\rho(n)$. Note that $\al_n$ projects monotonically to the interval $ \llbracket v_0,v_n\rrbracket$, 
where $v_n= \pi(\rho(n))$. 
Since $X$ is a proper metric space, the sequence of uniform quasigeodesic segments $(\alpha_n)$ subconverges to a 
(uniform) quasigeodesic ray, say $\alpha$, asymptotic to $\rho(\infty)$. Since each $\al_n$ projects to an interval in $T$, the projection of $\al$ to $T$ is a generalized geodesic ray. 

\medskip
There are two cases to consider. 

{\bf Case 1:} $\pi(\rho)$ has finite diameter. In this case the image of $\pi\circ \alpha$ is a finite geodesic segment in $T$,
say $\llbracket v_0,v\rrbracket$. It follows from the nature of paths $\al_n$ that $\alpha\cap X_v$ is a geodesic ray. Now by Proposition \ref{acyl finite geodesics}, in order to connect $x_0$ to 
$\al(\infty)$, we can first connect $x_0$ to its nearest-point projection $x_1$ in $X_v$ and then connect $x_1$ to $\al(\infty)\in \geo X_v$ 
by a geodesic ray in $X_v$. The concatenation of such two paths is a HV quasigeodesic of type 1.

\medskip
{\bf Case 2:} $\pi(\rho)$ has  infinite diameter. In this case the image of $\pi\circ \alpha$ is a geodesic ray  $c= v_0\eta$ in $T$. 
We shall show that $\alpha=\rho_{\eta}$. It suffices to prove the claim for finite subsegments in $\al$.

Since each $\alpha_n$ is an HV path and $\pi\circ \alpha_n$ 
converges to a geodesic ray in $T$ it is clear that there is $v_1\in V(T)$ adjacent to $v_0$ such that for all large $n$, 
$\gamma_{v_1}\subset \alpha_n$. Then we can run the same argument by replacing $x_0$ by $x_{v_1}$. Thus by induction on $\pi\circ \alpha(n)$
we are done. By the construction, each subsegment of $\al$ connecting $x_0$ to $X_{v_n}$, $v_n=c(n)$, is the limit of a sequence of HV quasigeodesic segments connecting $x_0$ to $X_{v_n}$. Thus, the limit is again an HV quasigeodesic segment connecting $x_0$ to $X_{v_n}$. 
\qed

\begin{lemma}
Suppose $c: [0,\infty)\map T$ is a geodesic ray with  $v_n=c(n), n\in \NN,$ and $c(\infty)=\eta$. Then every sequence 
$x_n\in X_{v_n}$ converges to $\gamma_{\eta}(\infty)\in \geo X$.  
\end{lemma}
\proof Let $\gamma_n$ denote the concatenation of the path $\gamma_{v_n}$ and the geodesic $\llbracket x_{v_n}, x_n\rrbracket_{X_{v_n}}$. These concatenations are uniform quasigeodesics in $X$. Clearly, the sequence  $\{\gamma_n\}$ converges uniformly on compacts to $\gamma_{\eta}$. \qed

\medskip 
Our next goal is to describe the ideal boundary $\geo X$ in terms of ideal boundaries of vertex-spaces $\geo X_v$ and the ideal boundary of the tree $T$. Our description is similar to the one given by Dahmani, \cite[section 2]{MR2026551}, in the setting of graphs of hyperbolic groups.

We set  
$$
\tilde\ZZ:= \coprod_{v\in V(T)} \geo X_v$$
and define a relation $\sim$ on $\tilde\ZZ$ as 
follows:

If $\xi_u\in\geo X_u$ and $\xi_v\in \geo X_v$ then $\xi_u\sim \xi_v$ iff $\xi_u$ belongs to the ideal boundary flow $Fl(\{\xi_v\})$  of  
$\xi_v$ in $\geo X$, see Section \ref{sec:Boundary flows}. It follows from the discussion of ideal boundary flows in Section \ref{sec:Boundary flows} that any two boundary flows are either equal or disjoint, hence, $\sim$ is  an equivalence relation. 
We let $p: \tilde\ZZ \to \ZZ:=\tilde\ZZ/_\sim$ be the quotient map.   
We will topologize $\ZZ$ later on. For now we note that, in general, this topology {\em is not} the quotient topology of the coproduct topology on $\tilde\ZZ$ (cf. Lemma \ref{lem:Z-topology}), but the topology will restrict to the standard topology on each $\geo X_v$. 

The space $\ZZ$ will be homeomorphic to one part of $\geo X$. 
The other part will be homeomorphic to $\geo T$ (with its natural topology of the ideal boundary of a tree). There are 
natural maps from both $\geo  T$ and $\ZZ$ to $\geo X$ defined, respectively, as follows:

1. $f_{\geo T}: \eta\mapsto \gamma_{\eta}(\infty)$, $\eta\in \geo T$. 

2. Since each $X_v, v\in V(T)$, is qi embedded in $X$, we have topological embeddings 
$$
q_v: \geo X_v\embed \geo X. 
$$
By combining these embeddings, we obtain a map 
$$
q: \coprod_{v\in V(T)} \geo X_v \to \geo X. 
$$
Next, for two vertical  geodesic rays $\al_u$ in $X_u$ and $\al_v$ in $X_v$ ($u, v\in V(T)$) 
we have equivalence of the following statements:

(a) $q_u(\al_u(\infty))=q_v(\al_v(\infty))$.

(b) $\Hd(\al_u, \al_v)<\infty$. 

(c)  The boundary flow-space $Fl(\{\al_u(\infty)\})$ contains $\al_v(\infty)$. 

While the equivalence of (a) and (b) is immediate, the equivalence of (b) and (c) follows from  
Lemma \ref{lem: flow condition}(2). 

Thus, the maps $q$ satisfies the property that $z_1\sim z_2$ iff $q(z_1)=q(z_2)$. In particular, the map $q$ descends to a  
map $f_T: \ZZ\map \geo X$. Combining the maps $f_{\geo T}$, $f_{T}$ we obtain a map 
$$
f: \geo T \sqcup \ZZ\to \geo X. 
$$

\begin{lem}
The map $f$ is a bijection. 
\end{lem}
\proof 1. That $f_{\ZZ}$ is injective is immediate from the equivalence of (a) and (c) above. 

2. If $\eta_1, \eta_2\in \geo T$ are distinct, then the rays $v_0 \eta_1, v_0\eta_2$ in $T$ diverge, hence, 
their lifts $\gamma_{\eta_1}, \gamma_{\eta_2}$ in $X$ diverge as well. Hence, 
$\gamma_{\eta_1}(\infty)\ne \gamma_{\eta_2}(\infty)$. It follows that $f_{\geo T}$ is injective. 

3. Proposition \ref{ray in acyl tree} directly implies that the images of $f_{\geo T}, f_{T}$ are disjoint and their union is the entire 
$\geo X$. \qed 

\medskip
We next topologize the disjoint union of $\ZZ$ and $\geo T$ by defining a basis of this topology as follows.  Given a vertex 
$v\in T$, we let the {\em shadow} $Sh_v$ of $v$ in $V(T)\cup \geo T$ denote the set consisting of all $\xi\in \geo T$ and $w\in V(T)$ such that the rays $v_0\xi$ and segments $v_0w$  contain $v$. 

(A) Let $\eta$ be a point in $\geo T$ and let $v$ be a vertex in the ray $v_0\eta$. Then define
$$
U_{v,\eta}:= (Sh_v\cap \geo T) \sqcup \bigcup_{w\in Sh_v \cap V(T)} p(\geo X_v)\subset \geo T \sqcup \ZZ.  
$$
These subsets will form a basis of neighborhoods of $\eta$ in $\geo T \sqcup \ZZ$. 

\medskip 
(B) Let $[\zeta]$ be a point in $\ZZ$,  $R>0$. For each $\eta\in \geo T$ and $\zeta\in \geo X_v$ in the equivalence class $[\zeta]$ define the point $z_\eta\in X_v$ as $x_v$ if the ray $v_0\eta$ does not contain $v$, and so that $[x_v z_\eta]_{X_v}$ is the 
maximal  subsegment of $\gamma_\eta$ contained in $X_v$. Similarly, define $z_\eta$ for points $\eta\in \geo X_w$: If 
$\llbracket v_0, w\rrbracket$ is disjoint from $v$, then set $z_\eta:=x_v$, otherwise it is the maximal subsegment of the HV ray 
$\rho_\eta=\gamma_w\star x_w \eta$ contained in $X_v$.

(B1) We define $U^1_{R,[\zeta]}$ as the set of points $\eta\in \geo T$ 
such that there exists a representative $\zeta\in \geo X_v$ of $[\zeta]$ for which 
$(\eta.z_v)_{x_v}> R$.

(B2) Similarly, define 
$$
U^2_{R,[\zeta]}:= p(\{\eta\in \geo X_v: (\zeta.z_v)_{x_v}> R , \zeta\in [\zeta] \} ). 
$$

Here in both (B1) and (B2) where the Gromov-product is taken in $X_v$. 
 Set 
$$
U_{R,[\zeta]}:= U^1_{R,[\zeta]}\cup U^2_{R,[\zeta]}. 
$$

We leave it to the reader to check the following properties satisfied by the collection of {\em basic subsets} $U_{v,\eta}, U_{R,[\zeta]}$ defined above:

\begin{lemma}
1. $\eta\in U_{v,\eta}, [\zeta]\in U_{R,[\zeta]}$ for every $\eta, \zeta$, and $v, R$ as above. 

2. For every two basis subsets $U', U''$ as above containing a point $\xi\in \geo T\sqcup \ZZ$, there exists another basis subset $U$ containing $\xi$ such that $U\subset U'\cap U''$. 

3. For any two distinct points $\xi', \xi''\in  \geo T\sqcup \ZZ$, there exist disjoint basic subsets $U', U''$ containing  $\xi', \xi''$ respectively. 
\end{lemma}

It follows that the collection of basic sets $U_{v,\eta}, U_{R,\zeta}$ defines a Hausdorff topology $\tau$ on $\geo T\sqcup \ZZ$. From now on, we will equip $\geo T\sqcup \ZZ$ with this topology. We will see soon (Proposition \ref{prop:f-homeo}) that the topology $\tau$ is compact and metrizable. The next lemma, where we use metrizability of $(\ZZ,\tau)$,  
describes some properties of $(\ZZ,\tau)$: 

\begin{lem}
\label{lem:Z-topology}
i. For every sequence of distinct vertices $v_n\in V(T)$ the sequence of compacts $p(\geo X_{v_n})$ subconverges to a point in 
$(\geo T\sqcup \ZZ,\tau)$.

ii. For every finite diameter subtree $S\subset T$, the restriction of $\tau$ to $p(\coprod_{v\in V(S)} \geo X_v)\subset \ZZ$ 
is compact. 
\end{lem}
\proof i.  
After extraction, there are two cases which may occur:

{\bf Case i1}. There is a vertex 
$v_0\in V(S)$ such that all the edges $e_n=[v_0, w_n]$ in the segments $\llbracket v_0, v_n\rrbracket$ 
are pairwise distinct. Pick a base-point $x_0\in X_{v_0}$. Then, again by properness of $X$, the sequence of distances from $x_0$ to $X_{e_nv_0}$ (in $X_{v_0}$) diverges to infinity. It follows that, after further extraction, the sequence of subsets $X_{e_nv_0}$ converges to a point $\xi\in \geo X_{v_0}$. By the description of the topology $\tau$, the sequence $(x_n)$ converges to $p(\xi)\in \ZZ$ as well. 

{\bf Case i2}. The sequence $(v_n)$ converges to a point $\eta\in \geo T$. Then, by the definition of the topology $\tau$, every sequence $(x_n)$ in $X_{v_n}$ converges to $\eta$. 

\medskip 
ii. Since $(\ZZ,\tau)$ is metrizable, it suffices to prove sequential compactness. Consider a sequence $x_n\in X_{v_n}$ where $v_n\in V(S)$.  
There are two cases to consider:

{\bf Case ii1}. Suppose that, after extraction, all the vertices $v_n$ are distinct. Then, by Part i (Case i1), the sequence $(x_n)$ subconverges to a point in $p(\geo X_v)$ for some $v\in V(S)$. 

{\bf Case ii2}. The same vertex $v$ appears in the sequence $(v_n)$ infinitely many times. Then, after extraction, $v_n=v$ and the sequence 
$(x_n)$ subconverges to a point in $\geo X_v$ (recall that $X$ is assumed to be proper, which implies properness of $X_v$).  \qed

\medskip 
We can now prove the main result of this section:

\begin{prop}\label{prop:f-homeo}
The map $f: \geo T\sqcup \ZZ\to \geo X$ is a homeomorphism. 
\end{prop}
\proof Since $\geo X$ is compact, $\geo T\sqcup \ZZ$ is Hausdorff and $f$ is a bijection, it suffices to prove continuity of $f^{-1}$. 

{\bf Case 1.} Suppose that $\eta_n\in \geo T$ are such that the sequence $(\gamma_{\eta_n}(\infty))$ in $\geo X$ converges to 
$\gamma_\eta(\infty)$,  $\eta\in \geo T$. This implies that there exists a constant $D$ such that for every $R\ge 0$ and all sufficiently large $n$, the point $\gamma_\eta(R)$ lies in the $D$-neighborhood of $\gamma_{\eta_n}$. The same, therefore, holds for the rays $v_0\eta$, $v_0\eta_n$. That implies convergence $\eta_n\to\eta$.  

{\bf Case 2.} Suppose that $\xi_n\in \geo X_{w_n}$ is a sequence converging to $\xi=\rho_\eta(\infty)$. Then, for the same reason as in Case 1, the sequence $(w_n)$ converges to $\eta$ in $T\cup \geo T$. It then follows by the definition of a neighborhood basis at $\eta$ 
in $\ZZ\cup \geo T$ that $\xi_n\to \eta$. 

{\bf Case 3.} Suppose that $\eta_n$ is a sequence in $\geo T$ such that the corresponding sequence $\xi_n:= \gamma_{\eta_n}(\infty)$ converges to some $\xi\in \geo X_v$. Then, after extraction, the sequence of HV (uniformly quasigeodesic) rays 
$\gamma_{\eta_n}$ converges to an HV ray $\gamma$ asymptotic to $\xi$. This limiting ray necessarily has the form of a concatenation $\rho_{\xi'}\gamma_{w}\star \alpha_w$, where $\alpha_w$ is a geodesic ray in a vertex space $X_w$. 
Since the rays 
$\gamma$, $\gamma_\eta$ are at finite Hausdorff distance from each other, it follows that $p(\xi)=p(\xi')$, i.e. $[\xi]=[\xi']$.
Furthermore, let $[x_v z_n]_{X_v}$ is the maximal subinterval of $\gamma_{\eta_n}$ contained in $X_v$. Then 
$$
\lim_{n\to\infty} (\xi'.z_n)_{x_v}=\infty. 
$$
Now, it follows from the definition 
of a neighborhood basis of $[\xi]$ in $\ZZ\sqcup \geo T$ (see the description of neighborhoods $U^1_{[\xi],R}$) that 
$$
\lim_{n\to\infty} \eta_n= [\xi] 
$$ 
in the topology of $\ZZ\sqcup \geo T$. 

{\bf Case 4.} The proof in the last case, when $\xi_n\in \geo X_{v_n}$ is a sequence converging to $\xi\in \geo X_v$ in $\geo X$ is similar to Case 3 and is left to the reader.  \qed

\chapter{Cannon--Thurston maps} \label{ch:CT}

In this chapter, as an application of the description of uniform quasigeodesics in trees of hyperbolic spaces, we establish an existence theorem for {\em Cannon--Thurston} maps (CT-maps) between ideal boundaries of trees of hyperbolic spaces induced by the inclusion maps of subtrees of spaces, $X_S=Y\to X$, Theorem \ref{thm:mainCT}. The proof of this theorem occupies most of the chapter.  Once this theorem is proven, we investigate the associated  {\em Cannon--Thurston laminations}. In Section \ref{sec:CTfibers} we  identify the CT-lamination $\La(Y,X)$, while in Section \ref{sec:CT-lamination} we relate $\La(Y,X)$ to the collection of {\em ending laminations} of $Y$ in $X$. 
We conclude  the chapter with Section \ref{sec:Group-theoretic applications} 
where we discuss group-theoretic application our results on CT-maps and CT-laminations. In particular, in Section \ref{sec:non-Anosov} we construct examples of undistorted surface subgroups of $PSL(2,\C)\times PSL(2,\C)$ which are not Anosov. (The proof of non-distortion is an application of Theorem \ref{thm:cut-paste}.)

\section{Generalities of Cannon--Thurston maps} \label{sec:CT-generalities}

If $X, Y$ are nonempty geodesic hyperbolic spaces and $f: Y\to X$ is a qi embedding, then $f$ induces a (continuous) embedding of Gromov-boundaries, $\geo Y\to \geo X$: A sequence $(y_n)$ in $Y$ is a Gromov-sequence if and only if $(f(y_n))$ is a Gromov-sequence in $X$ (see e.g. \cite[Theorem 5.38]{MR2164775} or \cite[Exercise 11.109]{Drutu-Kapovich}). The problem of existence of Cannon--Thurston maps concerns the existence of such an extension in the setting of  {\em uniformly proper maps}.

The original motivation for Cannon--Thurston maps comes from the group theory: Given a hyperbolic subgroup $H$ of a hyperbolic group $G$ (with the inclusion map $\iota=\iota_{H,G}$) or, more generally, a homomorphism with finite kernel $\phi: H\to G$ between two hyperbolic groups, 
one says that $\phi$ {\em admits a Cannon--Thurston map} if there exists a $\phi$-equivariant continuous map (a CT-map)
$$
\geo \phi: \geo H\to \geo G. 
$$
Surprisingly, CT-maps for group homomorphisms exist quite often (see \cite{MR3126566, MR3177384, mahan-pal, MR3497262, MR3652816, MR4077662}, etc.)  but not always (see \cite{MR3143716}). 
One of the earliest examples (justifying the terminology) of existence of CT-maps is due to Cannon and Thurston, \cite{MR2326947}: 
If $M$ is a compact hyperbolic 3-manifold fibered over the circle with the fiber $F$, then the natural embedding $\iota: \pi_1(F)\to \pi_1(M)$ admits a CT-map. 

One can further generalize the setting of CT-maps to the non-equivariant one: 

\begin{defn}\index{Cannon--Thurston map}
Let $X, Y$ be  geodesic Gromov-hyperbolic spaces and $f: H\to G$ is a coarse Lipschitz  map. Then $f$ is said to 
{\em admit a CT-map}  or {\em admits a CT-extension}, if there exists a  map 
$$
\geo f: \geo Y\to \geo X
$$
such that $f\cup \geo f: \bar{Y}= Y\cup \geo Y\to \bar{X}=X\cup \geo X$ is continuous on $\geo Y$.  
\end{defn}

In order to connect to the group-theoretic situation, one takes $X, Y$ to be Cayley graphs of the groups $G$ and $H$ respectively (where the finite generating set of $G$ contains that of $H$), and $f: Y\to X$ the inclusion map induced by the inclusion $\phi=\iota: H\to G$. Since $f$ is equivariant, the existence of a CT-extension of $f$ would imply the existence of a continuous equivariant 
map $\geo \phi$.

In the setting when $Y$ is a subspace of $X$ with the induced path-metric and the identity embedding $f=\iota_{Y,X}$, 
we will use the notation $\D_{Y,X}$ for the CT-map $\geo \iota_{Y,Y}$ (if it exists). The same notation $\D_{H,G}$ will be used for the CT-map $\geo \iota_{H,G}$, where $\iota_{H,G}: H\embed G$ is the identity embedding of a hyperbolic subgroup $H$ of a hyperbolic group $G$.

\medskip 
A criterion for the existence of CT-maps between hyperbolic metric spaces was established by Mahan Mitra in 
\cite[Lemma 2.1]{mitra-trees}.  Recall that $(x.y)_p$ denotes the Gromov-product in a metric space.

\begin{theorem}
[Mitra's Criterion] \label{thm:Mitra's Criterion}
Let $f: Y\to X$ be a coarse Lipschitz proper map of proper geodesic hyperbolic metric spaces.   
Then a Cannon--Thurston map $\geo f: \geo Y\to \geo X$ exists if and only if for some (each) $y_0\in Y$ 
the function
$$
t\mapsto \inf \{ (f(y_1). f(y_2))_{f(y_0)}:  y_1, y_2\in Y \hbox{~~are such that~~} (y_1. y_2)_{y_0} \ge t\} 
$$
is proper. 
\end{theorem}

Half of this theorem (most relevant for us) also holds for non-proper and non-geodesic hyperbolic spaces. The result was first proven in \cite{Krishna-Sardar}, we include a proof since it is quite simple and for the sake of completeness: 

\begin{prop}\label{prop:Mitra's Criterion}
Suppose that $f: Y\to X$ is a coarse Lipschitz 
map of hyperbolic spaces in the sense of Gromov, such that for each $p\in Y$ and each pair of sequences $y_n, y'_n\in Y$, 
$$
\lim_{n\to\infty} (y_n. y'_n)_p=\infty \RA 
\lim_{n\to\infty} (f(y_n). f(y'_n))_{f(p)}=\infty.
$$  
Then $f$ admits a Cannon--Thurston map. 
\end{prop}
\proof Note that the assumption in the proposition implies that the map $f$ is metrically proper: If a sequence $y_n\in Y$ diverges to infinity in the sense that $d(p, y_n)\to \infty$, then the sequence $f(y_n)$ also diverges to infinity. To prove properness one takes a sequence $y'_n=y_n$.

The definition of the Cannon--Thurston extension of the map $f$ is a natural one. The assumption in the proposition states that the image under $f$ of each Gromov-sequence $(y_n)$ in $Y$ is also a Gromov-sequence in $X$. Thus, one defines the extension by the formula:
$$
\geo f([y_n])= [f(y_n)],
$$
where $(y_n)$ is a Gromov-sequence. A diagonal argument shows that this map is well defined and, hence, is partially continuous: 
If a sequence $y_n\in Y$ converges to $\zeta\in \geo Y$, then the sequence $f(y_n)$  converges to $\geo f(\zeta)$. 
From this, by the density of $Y$ in $\bar{Y}$, it follows that $f\cup \geo f$ is continuous at $\geo Y$. \qed 

\begin{rem}
The converse to this proposition also holds if $Y$ is proper: If a map $f: Y\to X$ admits a CT-extension, then for every pair of sequences $(y_n), (y'_n)$ in $Y$, 
$$
\lim_{n\to\infty} (y_n. y'_n)_p=\infty \RA 
\lim_{n\to\infty} (f(y_n). f(y'_n))_{f(p)}=\infty.
$$  
The reason is that such sequences $(y_n), (y'_n)$ have to subconverge to points $\xi, \xi'\in \geo Y$ (in view of properness of $Y$). Then, necessarily, $\xi=\xi'$. 
Since $f$ admits a CT-extension to the points $\xi, \xi'$, the sequences $(f(y_n)), (f(y'_n))$ have to converge to $\geo f(\xi)$, implying that 
$$
\lim_{n\to\infty} (f(y_n). f(y'_n))_{f(p)}=\infty.
$$  
\end{rem}

In the next   lemma we will use the notion of the {\em relative boundary} $\geo (A, X)$ of a subset $A$ of a hyperbolic space $X$, see Definition \ref{not:rel-bdry}. 

\begin{lemma}\label{lem:image-of-CT}
Suppose that $Y$ is a proper metric space, and $f: Y\to X$ admits a CT-extension. Then 
$$
\geo f(\geo Y)= \geo(f(Y), X).
$$
\end{lemma}
\proof By the continuity of the CT-extension, if a sequence $y_n\in Y$ converges to $\xi\in \geo Y$, then the sequence $(f(y_n))$ converges to $\geo f(\xi)$. Conversely, suppose that a sequence $x_n=f(y_n)\in f(Y)$ converges to a point $\eta\in \geo(f(Y), X)$. By the properness of $Y$, the sequence $(y_n)$ subconverges to some $\xi\in \geo Y$ and, again by continuity of the CT-extension, $\geo f(\xi)=\eta$. \qed

\begin{lemma}[Functoriality of CT-maps] 
If $f: X\to Y, g: Y\to Z$ are coarse Lipschitz maps which admits CT-extensions, then their composition also does and 
$$
\geo (g\circ f)= \geo g \circ \geo f. 
$$
\end{lemma}
\proof Let $(x_n)$ be a Gromov-sequence representing $\xi\in \geo X$. Then its image, the sequence $(y_n)$ in $Y$ is also a Gromov-sequence in $Y$ representing the point $\eta=\geo f(\xi)$. Applying the same reasoning to the  map $g$, we conclude that $(g(y_n))$ is a Gromov sequence in $Z$ representing $\zeta=\geo g(\eta)$. Thus, $\geo g \circ \geo f$ defines the CT-extension of $g\circ f$. \qed

\medskip 
Since for $\delta$-hyperbolic spaces in the sense of Rips, the Gromov-product is comparable to the  distance to a suitable geodesic (see Lemma \ref{lem:gromov product}),   
Proposition \ref{prop:Mitra's Criterion} can be reformulated as

\begin{prop}\label{prop:Mitra's Criterion1}
Suppose that $f: Y\to X$ is a coarse Lipschitz  
map of hyperbolic spaces in the sense of Rips, such that for each pair of sequences $y_n, y'_n\in Y$, if  
$$
\lim_{n\to\infty} d_Y(p, [y_n y'_n]_Y)=\infty\ RA 
\lim_{n\to\infty} d_X(f(p),  [f(y_n) f(y'_n)]_{X})=\infty.
$$  
Then $f$ admits a Cannon--Thurston map. 
\end{prop}

The existence of a CT-map, of course, does not imply its injectivity, and the notion of a Cannon--Thurston lamination (introduced by Mitra in \cite{MR1445392}) is motivated by this lack of injectivity:

\begin{defn}\label{defn:CT-lamination} \index{Cannon--Thurston lamination}
Suppose that $f: Y\to X$ is a map of hyperbolic spaces which admits a CT-extension $\geo f$. 
The {\em Cannon--Thurston lamination} (the CT-lamination) of  $f: Y\to X$ is the (closed) subset $\La(f)$ of $\geo^{(2)}Y$ consisting of unordered pairs of distinct points $\{\xi,\eta\}$ such that $\geo f(\xi)=\geo f(\eta)$. 
In the case when $Y$ is a subset of $X$ and $f$ is the inclusion map $Y\to X$, we will use the notation $\La(Y,X)$ for the CT-lamination. A geodesic $\al\subset Y$ connecting points $\xi, \eta$ with $\{\xi,\eta\}\in \La(f)$, is called a {\em leaf} of the CT-lamination $\La(f)$. 
\end{defn}

Note that, in view of the fact that the map $\geo f$ is continuous, the lamination $\La(f)$ is a closed subset of $\geo^{(2)} Y$. 

\begin{rem}\label{rem:CT-remark}
1. The above definition of $\La(f)$ requires existence of a CT-map. However, one can extend this definition to the general case as follows (see \cite[section 2]{MR1610757}, \cite[Definition 3.1]{Mj-Rafi}). We say that a point $\{\xi, \xi'\}\in \geo^{(2)}Y$ belongs to $\La(f)$ if there exist sequences $(y_n), (y'_p)$ in $Y$ converging to $\xi, \xi'$ respectively, such that 
$$
\lim_{n\to\infty} (f(y_n). f(y'_n))_{f(p)}=\infty.
$$ 

2. As it was noted in \cite{Mj-Rafi}, in the general setting, $\La(f)$, {\em a priori} is not closed in $\geo^{(2)}Y$.  If the map $f$ admits a CT-extension, then the two definitions of $\La(f)$ agree. 

3. If $H$ is a hyperbolic subgroup of a hyperbolic group $G$ and $f$ is the inclusion map $H\to G$ then $\La(f)=\emptyset$ if and only if $H$ is quasiconvex in $G$, see 
\cite[Lemma 2.1]{MR1610757}. 
\end{rem}

\medskip 
Mitra proved in \cite{mitra-trees} that if $G$ is a hyperbolic group isomorphic to the fundamental group of a finite graph ${\mathcal G}$ of hyperbolic groups satisfying the conditions of the Bestvina-Feighn Combination Theorem, 
 then for each vertex-group $G_v$ of ${\mathcal G}$, 
the Cannon--Thurston map for the inclusion homomorphism $G_v\to G$ exists. More generally, he proves:

\begin{theorem}\label{thm:mitra-trees}
If $X\to T$ is a tree of hyperbolic spaces with hyperbolic total space $X$, 
then for every vertex space $X_v$ the inclusion map $X_v\to X$ admits a CT-map.  
\end{theorem}

Later on, in the paper by Mj and Pal \cite{mahan-pal}, this result was extended to the relatively hyperbolic setting; we will discuss the extension in the next chapter of the book.

Our goal is to generalize Theorem \ref{thm:mitra-trees} to the case of fundamental
groups of subgraphs of graphs of groups and, more generally, to 
inclusion maps $Y\to X$ of subtrees of spaces in a tree of hyperbolic spaces $\Y\subset \X$, 
satisfying the conditions of Theorem \ref{thm:mainBF}. The main result of this chapter is:

\begin{thm}\label{thm:mainCT}
Let $\X=(\pi: X\to T)$ be a tree of hyperbolic spaces with hyperbolic total space $X$. Then for every subtree $S\subset T$, the inclusion map $X_S\to X$ admits a CT-extension. 
\end{thm}

\medskip
The most difficult part of the proof   is to relate, for points $x, y\in Y$, the  geodesics $[xy]_X$ in $X$ to the geodesics $[xy]_Y$ in $Y$. This is done in Section \ref{sec:Cut-and-replace} in the form of 
a ``cut-and-replace'' theorem (Theorem \ref{thm:cut-paste}). 
Once this theorem is established, the existence of a CT-map is an almost immediate consequence of  Proposition \ref{prop:Mitra's Criterion1} (see Theorem \ref{thm:ECT}). 
The main tool in our proof of Theorem \ref{thm:cut-paste} is the description of geodesics in hyperbolic  trees of  spaces given in the previous chapter. 
As this description is inductive in nature (a 7-step process), the proof of the cut-and-replace theorem 
follows the same inductive process (but we will only need 6 steps). While our proof follows in main Mitra's proof in 
 \cite{mitra-trees}, we have to deal with some substantial complications; in fact, we will derive Mitra's theorem (Theorem \ref{thm:mitra-trees}) as an easy application of the first part of our proof, see the end of Section \ref{sec:Part I}. 
 
In the proof of Theorem \ref{thm:cut-paste} we will be using the fact that $\X$  satisfies  the uniform $K$-flaring condition for all $K\ge 1$ (see Lemma \ref{lem:hyp->uniform flaring}). 
Of course, if $\X$ satisfied the $K$-flaring condition, so does $\Y\subset \X$.

\section{Cut-and-replace theorem}\label{sec:Cut-and-replace}

\subsection{Definitions and notations}
 
 Suppose that $\X=(\pi: X\to T)$ is a tree of spaces (not necessarily hyperbolic), containing a subtree of spaces $\Y= (\pi: X_S=Y\to S)$. In the next section we shall prove that the inclusion $Y\map X$ admits a Cannon--Thurston map 
$\geo Y\map \geo X$ provided that $X$ is hyperbolic. To prove this result we need to compare the $Y$-geodesics 
$[xy]_Y$ to $X$-geodesics $[xy]_X$ joining  pairs points $x, y\in Y$.
When the points are in the same vertex space in $Y$, this is done by Mahan Mitra in \cite{mitra-trees} by constructing 
ladders.  In general, up to a uniformly bounded error, the relation between $X$-geodesics and $Y$-geodesics is given by a {\em cut-and-replace procedure} described below.

\medskip 
For each (continuous) path $c: I\subset \RR\to X$ in $X$ with $c(\partial I)\subset Y$, we define the following modification,  
a {\em cut-and-replace procedure}, transforming $c$ to a new path $\hat{c}=c_S: I\to Y$.  

\begin{defn}\index{detour subpath} \index{modified path $\hat{c}$}\label{defn:detour} 
For a closed subinterval $J= [s, s']\subset I$, we say that the restriction $\zeta=c|_J$  
 is a {\em detour subpath} in $c$, if $c(\partial J)\subset Y$, while $c(J- \partial J)$ is disjoint from $Y$. 
 Thus, the points $x=c(s), x'=c(s')$ belong to a common vertex-space $X_t\subset Y$. We then replace each detour subpath $\zeta$ in 
 $c$ by the corresponding $X_t$-geodesic $\hat\zeta=[x x']_{X_t}$, called a {\em replacement segment}.
 \footnote{While this $X_t$-geodesic is, in general, non-unique, if vertex spaces are uniformly hyperbolic (which will be the case in all our examples) the ambiguity is uniformly bounded and we will ignore it.}  
We let $\hat{c}=c_S$ denote the resulting path $I\to Y$. 
\end{defn}

\begin{figure}[tbh]
\centering
\includegraphics[width=60mm]{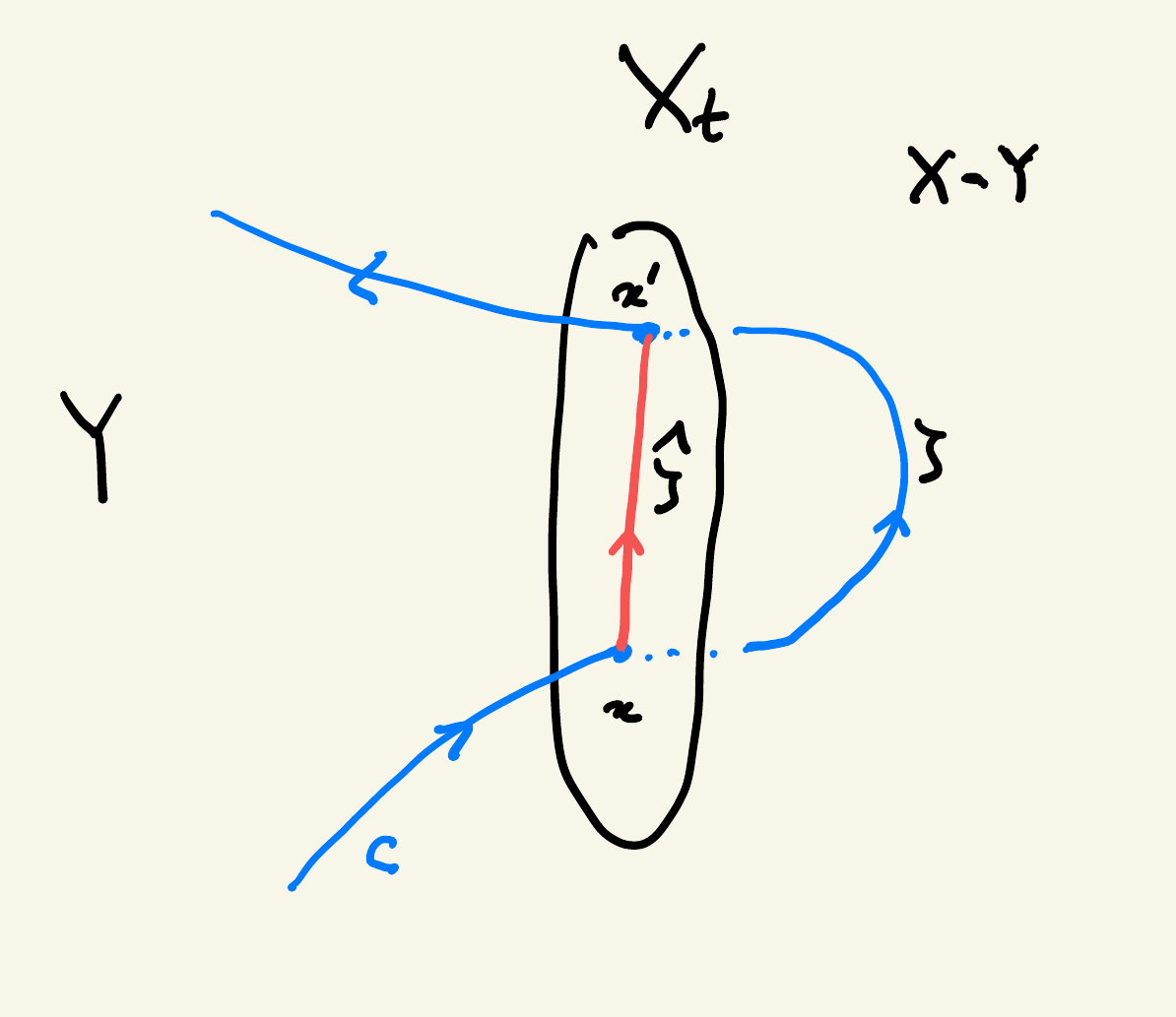}
\caption{Detours and path-modification}
\label{detour.fig}
\end{figure}

\begin{comment}
earlier?}As a general rule, given a path $c: I\to X$ and two points $x, y$ in $c$, we let $c(x,y)$ be the subpath in $c$ between $x$ and $y$. To be more precise, instead of specifying $x$ and $y$ we have to specify points $s, t \in I$ such that $c(s)=x, c(t)=y$, but in practice, our notation will be unambiguous.  \end{comment} 

\begin{defn}\index{consistent pair} 
Let $\phi$ be a continuous $\La$-quasigeodesic in $X$ with the end-points in $Y$. We will say that $\phi$ is $\La'$-{\em consistent} if 
$\hat\phi$ is a $\La'$-quasigeodesic in $Y$. We will say that a pair of points $(y, y')\in Y^2$ is $\theta$-{\em consistent}, 
where $\theta$ is a function $[1,\infty)\to [1,\infty)$,  
if every continuous $\La$-quasigeodesic $\phi$ in $X$ connecting $y$ to $y'$ is $\theta(\La)$-consistent.  
We say that a subset of $Y^2$ is (uniformly) 
consistent if it consists of  $\theta$-consistent pairs for some function $\theta$. We will say that a pair of points $(y,y')\in Y^2$ is 
$\La'$-{\em weakly consistent} if $y, y'$ are connected by {\em some} $\La'$-consistent continuous quasigeodesic in  $X$. 
Similarly, we will say that a subset of $Y^2$ is {\em uniformly weakly consistent} if it consists of $\La'$-weakly consistent 
pairs of points for some $\La'$. 
\end{defn}

While it is clear that every consistent pair of points is also $\La'$-consistent for some $\La'$, it is far from clear that if 
a subset of $Y^2$ consisting of $\La'$-consistent pairs is uniformly consistent. The issue is that while any two uniform 
quasigeodesics $\phi, \phi'$ connecting $y$ to $y'$ are uniformly close, we do not yet know if the same holds 
for $\hat\phi, \hat\phi'$. We will prove the claim, establishing equivalence of weak consistency and consistency, at the end of Part I of the proof of Theorem \ref{thm:cut-paste}. Before that, we will be only proving consistency of specific quasigeodesics, namely  slim combing paths $c(y,y')$ described  in Section \ref{sec:inductive}.

\medskip 
We will frequently use the following simple observation:

\begin{rem}\label{rem:trivial1}
Suppose that a pair $(x,y)\in Y^2$  in $X$ is $\La'$-consistent. Then every pair of points 
$x', y'\in N_r([xy]_X)\cap Y$ is  $\theta'$-consistent, where $\theta'$ depends only on $\theta'$ and $r$. In particular, perturbing the points $x, y$ by a uniformly bounded amount, we do not loose consistency of the pair. 
\end{rem}

\begin{defn} 
Given two  subsets $Z', Z''$ of a Gromov-hyperbolic space $Z$, we say that $y\in Z$ is an $R$-{\em transition point} \index{transition point}
between  $Z', Z''$ in $Z$, if every geodesic in $Z$ connecting $Z', Z''$ passes within distance $R$ from $y$. We say that a  finite sequence of points 
$z'=z_1, z_2,...,z_{n+1}=z''$ in a $\delta$-hyperbolic space $Z$ is  {\em $R$-straight} if the geodesic $\gamma=z'z''\subset  Z$ \index{straight sequence}
passes through some points $x_i\in B(z_i,R)$ in the order $x_2 ,..., x_n$. 

In relation to (sub)trees of spaces $Y\subset X$, we will talk about $X$-transition points and $Y$-transition points. We will say that $y\in Y$ is an  $R$-transition point between $Z', Z''\subset Y$ if it is $R$-transition point for $Z', Z''$ regarded as subsets in {\em both} 
$Y$ and $X$. 
\end{defn}

When dealing with sets of finite sequences, we will refer to those as {\em uniformly straight} if they are $R$-straight for some uniform value of $R$. We use a similar terminology for transition points.

\begin{example}\label{ex:trivial}
Suppose that $Z, Z'$ are $\la$-quasiconvex and $C$-cobounded subsets in a $\delta$-hyperbolic geodesic space $X$. 
Let $\beta$ be a shortest geodesic between $Z$ and $Z'$. Then every point of $\beta$ is an $R$-transition point between 
$Z$ and $Z'$ with $R=R(C,\la, \delta)$.  Conversely, every $R$-transition point between such $Z, Z'$ is $D$-uniformly 
close to a point in $\beta$, where $D=D(C, B, \delta, R)$. 
\end{example}

The importance of the concept of a transition point comes from another simple observation. Let $\X$ be a tree of spaces and 
$\Y\subset \X$ be a  subtree of spaces. 

\begin{lemma}\label{lem:transitions} 
Suppose that $Z_1, Z_2\subset X$ and $y\in Y$ is an $R$-transition point 
between $Z_1, Z_2$. Assume that 
the set of pairs $(z_i,y)$, $z_i\in Z_i$, $i=1,2$, is $\theta$-consistent. 
Then the set of pairs $(z_1, z_2)\in Z_1\times Z_2$ is  
$\theta'$-consistent where $\theta'=\theta'(L,R)$. 
\end{lemma}
\proof Since we are dealing with quasigeodesics, we can as well consider $L$-quasigeo\-de\-sics $c$ in $X$ connecting point $z_1\in Z_1$ to  $z_2\in Z_2$ and passing through $y$ (such exist due to the assumption that $y$ is an $R$-transition points in $X$). Such $c$ is a concatenation $c_1\star c_2$, where $c_1$ connects $z_1$ to $y$. Thus, since $y\in Y$, 
$$
\hat{c}= \hat{c}_1\star \hat{c}_2.
$$
Each of the subpaths $\hat{c}_i$ is a $\theta(L)$-quasigeodesic by the consistency assumption for the pairs $(z_i,y)$. We will estimate the qi constant of $\hat{c}$ in $Y$.  
Let $\gamma$ be a geodesic in $Y$ connecting two points $a, b$ in $\hat{c}$; the only interesting case to consider is when $a\in \hat{c}_1, b\in \hat{c}_2$. 
Let $\hat{c}(a,b)$ be the portion of $\hat{c}$ between $a$ and $b$. Since $y$ is an $R$-transition point in $Y$ between $z_1, z_2$, 
the path $\gamma$ has to pass within distance $r=r(R,L)$ from $y$. Subdividing $\gamma$ as a  concatenation $\gamma_1\star \gamma_2$, where $\gamma_1$ connects $z_1$ to $x\in B(y,r)$, $i=1,2$, we see that 
$$
\length (\hat{c}(a,b))\le (\theta(L) + r)\length(\gamma),
$$
as required by a uniform quasigeodesic.\footnote{Here and in what follows we repeatedly use the notation $\phi(p,q)$ for a subpath in a path $\phi$ between the points $p, q$ in $\phi$, see Section \ref{sec:Metric notions}.} \qed 

\medskip
This lemma generalizes to the case of higher number of transition/concatenation points. Since the proof is similar, we leave it to the reader:  

\begin{lemma}\label{lem:trivial3}
Suppose that $z'=z_1, z_2,...,z_{n+1}=z''$ is an $R$-straight sequence in both $X$ and $Y$, and that each pair of points 
$(z_i, z_{i+1}), i=1,...,n$, is $\theta$-consistent. Then the pair $(z', z'')$ 
is $\theta'$-consistent with the function $\theta'$ depending only on the hyperbolicity constant of $Y$ and on $\theta$.   
\end{lemma}

\medskip 
We are now ready to state the main technical result of this chapter:

\begin{theorem}[Cut-and-replace $X$-quasigeodesics to get $Y$-quasigeodesics] \label{thm:cut-paste}
Suppose that $\X= (\pi: X\to T)$ is a tree of hyperbolic spaces satisfying the uniform $K_*$-flaring condition, $S\subset T$ is a subtree and $Y=X_S$. Then $Y\times Y$ is $\theta$-consistent with $\theta$ depending only on the parameters of $\X$. 
\end{theorem}

We  break the proof of Theorem \ref{thm:cut-paste} in three  parts and each part in several steps, where we prove this theorem in  special cases and then use these special cases to prove the general case in the last part.

\section{Part I: Consistency of points in vertex flow-spaces}\label{sec:Part I}

Suppose that $k\ge k_{\ref{prop:existence-of-tripod-ladders}}$, $K= K_{\ref{prop:existence-of-tripod-ladders}}(k)$, and that 
$\X$ is a tree of hyperbolic spaces 
satisfying the uniform $\kappa_{\ref{defn:carpeted ladder}}(K)$-flaring condition. (As it was noted earlier, in Lemma 
\ref{lem:hyp->uniform flaring}, uniform $\kappa$-flaring holds for all $\kappa$'s if $X$ is hyperbolic.)  The main result of this section is:

\begin{prop}\label{prop:oneflow}
For every $u\in V(S)$, the set of pairs 
$(y, y') \in Fl_{Y,k}(X_u) \times Fl_{Y,k}(X_u)$ is $\theta= \theta_{\ref{prop:oneflow},k}$-consistent, 
with $\theta$ depending only on $k$ and the parameters of $\X$. 
\end{prop}
\proof For most of the proof we will be only proving weak uniform consistency, working  with the quasigeodesic 
paths $c=c_X(y,y')$ in the fiberwise $4\delta_0$-neighborhood of  
$\Fl_{X,k}(X_u)\subset \Y$ given by the slim combing of $\Fl_{X,k}(X_u)$ described in Section \ref{sec:inductive}.

According to Corollary \ref{cor:existence-of-ladders}, there exists a $(K,D,E)$-ladder $\L=\L_X(\al)\subset \Y$ contained in the fiberwise $4\delta_0$-neighborhood of  $\Fl_{X,k}(X_u)\subset \Y$, containing $y$ (resp. $y'$) in its bottom (resp. top), where 
$K, D, E$ depend only on $k$ and $\al$ is a geodesic in $X_u$. Recall that, by the very definition, $c_X(y,y')=c_{\L}(y,y')$ is contained in the ladder $L_X(\al)$. 
Define $\L_Y= \L\cap \Y$: It follows from the definition of a ladder that $\L_Y$ is a $(K,D,E)$-ladder  in $\Y$. It also follows from  the definition of the modification $c\mapsto \hat{c}$, 
that $\hat{c}$ is contained in $L_Y$ (up to a uniformly bounded error which we will ignore).

 \begin{comment}
\begin{figure}[tbh]
\centering
\includegraphics[width=150mm]{fig1-1.pdf}
\caption{}
\label{F1.fig}
\end{figure}
\end{comment}

There are several cases to consider, according to the construction of uniformly quasigeodesic paths $c_{\L}(y,y')$, 
depending on the properties of the ladder $\L$ and  location of the points $y, y'$.

\subsection{Part I.1: The points $y, y'$ belong to a $(K, M_{\bar{K}})$-narrow carpet 
$\A=\A_X= (\pi: A_X\to \llbracket u,w\rrbracket)\subset \L_X$}

\medskip 
 The carpet $\A$ contains the ``subcarpet'' $\A_Y= \A_X\cap \Y$, $\llbracket u, w'\rrbracket= \pi(A_Y)$. 
Note that $\A_Y$ is a $(K,C)$-carpet where $C$ the length of the ``narrow end'' $A_{w'}$ of $\A_Y$, but we cannot bound $C$ 
(from above) in terms of $k$. According to Corollary \ref{cor:X-to-bundle}, we have  the coarse  $L_{\ref{cor:X-to-bundle}}$-Lipschitz retraction $\rho_\A: X\to A_X$, where  $L_{\ref{cor:X-to-bundle}}$ depends only on $K$, $D$ and $E$ (hence, only on $k$). The restriction of this retraction to $Y$ is a retraction to $A_Y$;  in particular, 
$A_Y$ is $L_{\ref{cor:X-to-bundle}}$-qi embedded in $Y$.

\begin{lemma}\label{lem:I.1}
The pairs $(y,y')$ are $\La_{\ref{lem:I.1}}$-weakly consistent, where $\La_{\ref{lem:I.1}}$ depends only on $k$. 
\end{lemma}
\proof  We let $\ga_{y}, \ga_{y'}$ denote the $K$-leaves in $\A_X$ connecting, respectively, $y, y'$ to points of $A_w$. 
Let $v\in \llbracket u,w\rrbracket$ be the infimum of all vertices $t$ in $\llbracket u,w\rrbracket$ such that 
$$
d_{X_t}(\ga_y(t), \ga_{y'}(t))\le M_{\bar{K}}. 
$$
Then $c=c_\A(y,y')$ is the concatenation of the subpath $\ga_y$ restricted to $\llbracket u,v\rrbracket$, followed by the vertical geodesic $[\ga_y(v) \ga_{y'}(v)]_{X_{v}}$ and then followed by $\ga_{y'}$ restricted to 
 $\llbracket v,\pi(y')\rrbracket$. The path $\hat{c}$ is a similar concatenation $c_1\star c_2\star c_3$ 
 except $v$ is replaced by the vertex $v'$ which is 
 the minimum of $\{v, w'\}$ in the oriented interval $\llbracket u,w\rrbracket$ (the paths $c_1, c_3$ are contained in $\ga_y, \ga_{y'}$ respectively and $c_2$ is contained in $A_{v'}$). But this path is exactly 
 the path $c_{\A_Y}(y, y')$ as defined in Step I.1 in Section \ref{sec:inductive}, or in the proof of Proposition \ref{prop:easy-one}. 
 
\begin{figure}[tbh]
\centering
\includegraphics[width=100mm]{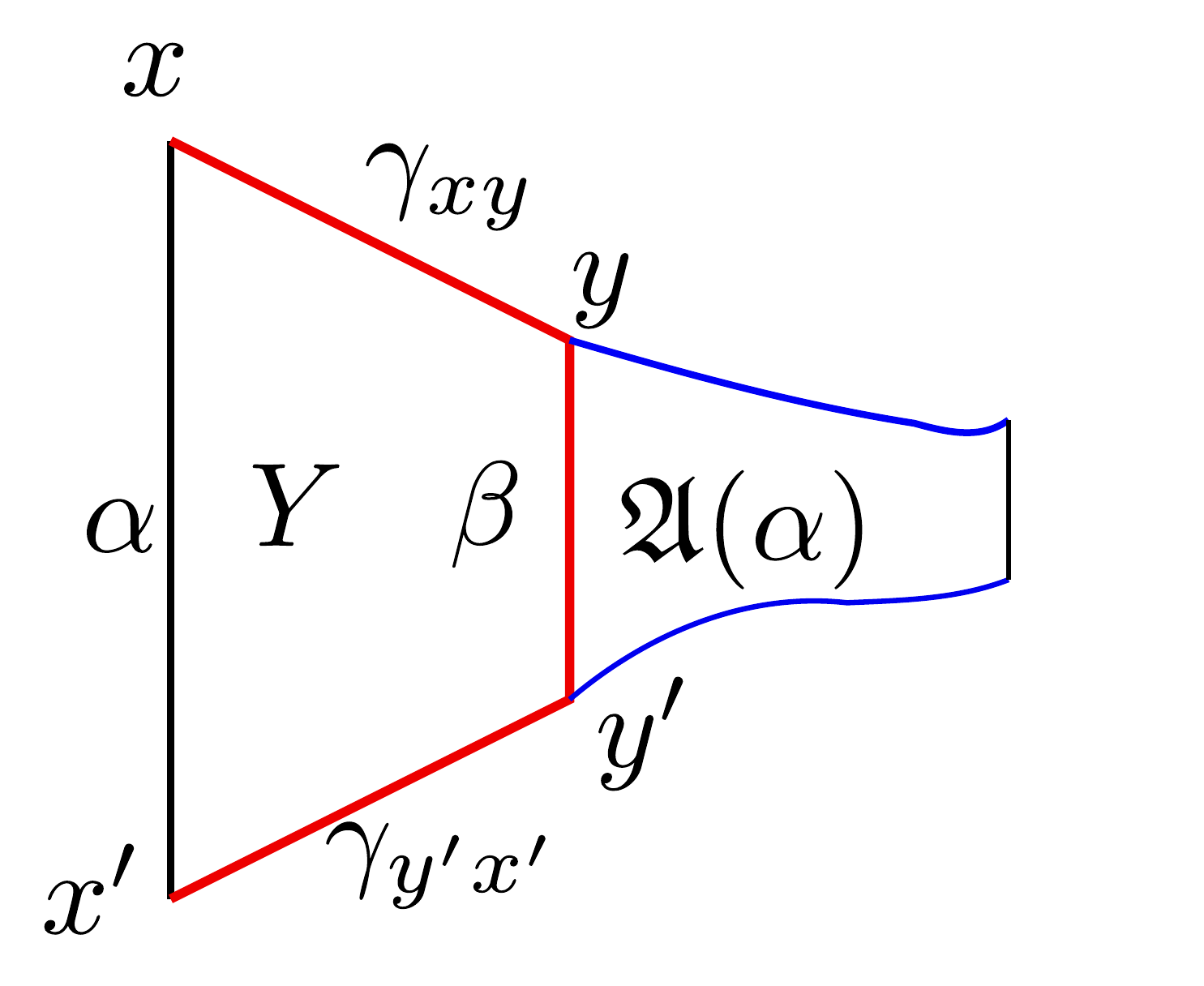} 
\caption{}
\label{F5.fig}
\end{figure}

 The quasigeodesic constant of $c_{\A_Y}(y, y')$ {\em a priori} depends on both $K$ and $C$. 
 However, according to Remark \ref{rem:dependence on C}, the dependence on  $C$ appears only in 
 the proof of Lemma \ref{lem:easy-up}, establishing uniform bounds on distortion of paths  
 $c_{\A_Y}(y, y')$ in $A_Y$. It remains, therefore, to get a uniform distortion bound depending only on $k$. Take a pair of 
 points $a, b\in \hat{c}$. There are several cases to consider depending on the location of the points $a, b$, we will treat just one since the rest are done by the same argument: We will assume that $a\in c_2, b\in c_3$. It suffices to bound the length of $\hat{c}(a,b)$ (between $a, b$) in terms of   $d_{\A}(a,b)$, equivalently,  in terms of the length of $c(a,b)$ since the latter is a uniform (with quasigeodesic constant depending only on $k$) quasigeodesic in $A_X$. The latter path is a concatenation of $c(a,a')$ and $c(a',b)$, 
 where $a'$ is in $A_{v'}$ and $c(a',b)= c_3$.  See Figure \ref{F5.fig}. We have
$$
\length(c_3)\le \length (c(a,b)),
$$
while
$$
\length(c_2)= d_{A_{v'}}(a,a') \le \eta(d_X(a,a'))\le \eta(\length(c(a,b))),    
$$
where $\eta=\eta_{\ref{unif-emb-subtree}}$. Thus,
$$
\length( \hat{c}(a,b))\le  \length (c(a,b)) + \eta(\length(c(a,b)),
$$
 thereby providing the required distortion bound depending only on $k$. \qed 
 
\subsection{Part I.2: The points $y, y'\in Y$ belong to a carpeted $X$-ladder $\L= \L_X(\al)$}

In this part of the proof we assume that $\L$ contains a $(K,M_{\bar{K}})$-narrow carpet 
$\A=\A_X(\al')$,  where $\al$ is a geodesic segment in $X_u$, $u\in S$, and $\al'\subset \al$  is a subsegment of length 
$\ge \length(\al)- M_{\bar{K}}$. The path $c=c_{\L}(y,y')$ can be of one of two types (see Section \ref{sec:inductive}, Step I.2, 
for the definition of the types of slim combing paths in a ladder):  

1. If $c$ is of type 1, then the assumption that $y, y'$ are in $Y$ and the definition of type 1 paths imply that $c=\hat{c}$, so there is nothing to prove. 

2. Suppose that $c$ is of type 2. Then $c$ is the concatenation $c_1\star c_2\star c_3$, where $c_2$ is contained in the 
carpet $\A$, $c_1, c_3$ are contained in $Y$  and, thus,  
$$
\hat{c}= c_1\star \hat{c}_2 \star c_3. 
$$
Below we will use the same notation for the subcarpet 
$\A_Y= \A\cap \Y$ as in Part I.1.

\begin{lemma}\label{lem:I.2}
The pairs $(y,y')$ are $\La_{\ref{lem:I.2}}$-weakly consistent, where $\La_{\ref{lem:I.2}}$ depends only on $k$. 
\end{lemma}
\proof 
The paths $c_1, c_3$ are uniformly quasigeodesic in $L_X$ (since $c$ is) while $\hat{c}_2$ is uniformly quasigeodesic in $A_Y$ according to Lemma \ref{lem:I.1}. Since $L_Y$ is uniformly hyperbolic, in order to prove the lemma, it suffices to verify that 
$\hat{c}$ is uniformly proper in $L_Y$ (see Lemma \ref{lem:quasigeodesic-paths}). The proof is similar to the one in 
Lemma \ref{lem:I.1}. We will prove uniform properness of $\hat{c}$ in $L_X$. As in the proof of 
Lemma \ref{lem:I.1}, we only consider the most representative case, of points $a\in c_2, b\in c_3$: We need to bound 
$\length(\hat{c}(a,b))$ in terms of $\length(c(a,b))$. The path $c(a,b)$ is the concatenation $c(a,a')\star c_3(a',b)$, where $a'$ is the concatenation point of $c_2$ and $c_3$. According to 
Lemma \ref{lem:I.1},
\begin{align*}
\length(\hat{c}(a,a'))\le  \La_{\ref{lem:I.1}} d_{\A_Y}(a,a') \le \La_{\ref{lem:I.1}} \eta_{\ref{unif-emb-subtree}}(d_{\A}(a,a')) \\
\le   \La_{\ref{lem:I.1}} \eta_{\ref{unif-emb-subtree}}(\length(c(a,b))), 
\end{align*} 
while 
$$
\length(c_3))\le \length(c(a,b)).  
$$
Lemma follows. \qed

\subsection{Part I.3: General ladders}\label{sec:StepI.3}

  Suppose that $u\in S$, $\al= [pp']_{X_u}\subset Y$, $\L=\L_X(\al)\subset \X$ is a $K$-ladder, 
$\L_Y= \L\cap \Y$. Our goal is to prove uniform consistency of  paths $c_\L$ connecting points   
$y\in bot(\L)\cap Y$, $y'\in top(\L)\cap Y$. Recall that according to Proposition Proposition \ref{vertical subdivision} (on vertical subdivision), 
we have a subdivision of $\al$ into subintervals $\al_i= [p_{i} p_{i+1}]_{X_u}$, subintervals $\al'_i\subset \al_i$, and a collection of $K$-qi sections $\Sigma_i$ in $\L$ through the points $p_i$ dividing $\L$ into subladders $\L^i= \L(\al_i)$ containing $(K, M_{\bar{K}})$-narrow carpets $\A^i= \A(\al'_i)$. We also defined points $x^\pm_i$ in the sections $\Si_i^-=\Si_i, \Si_i^+=\Si_{i+1}$ bounding $\L^i$ such that the combing paths in $L$ (connecting $y, y'$) pass through the points $x^\pm_i$. Each section 
$\Si_i$ is defined over some subtree $T_i\subset T$. 
The intersections $\Si_{i,Y}= \Si_i\cap Y$ project to subtrees $S_i\subset T_i$.

  \begin{figure}[tbh]
\centering
\includegraphics[width=100mm]{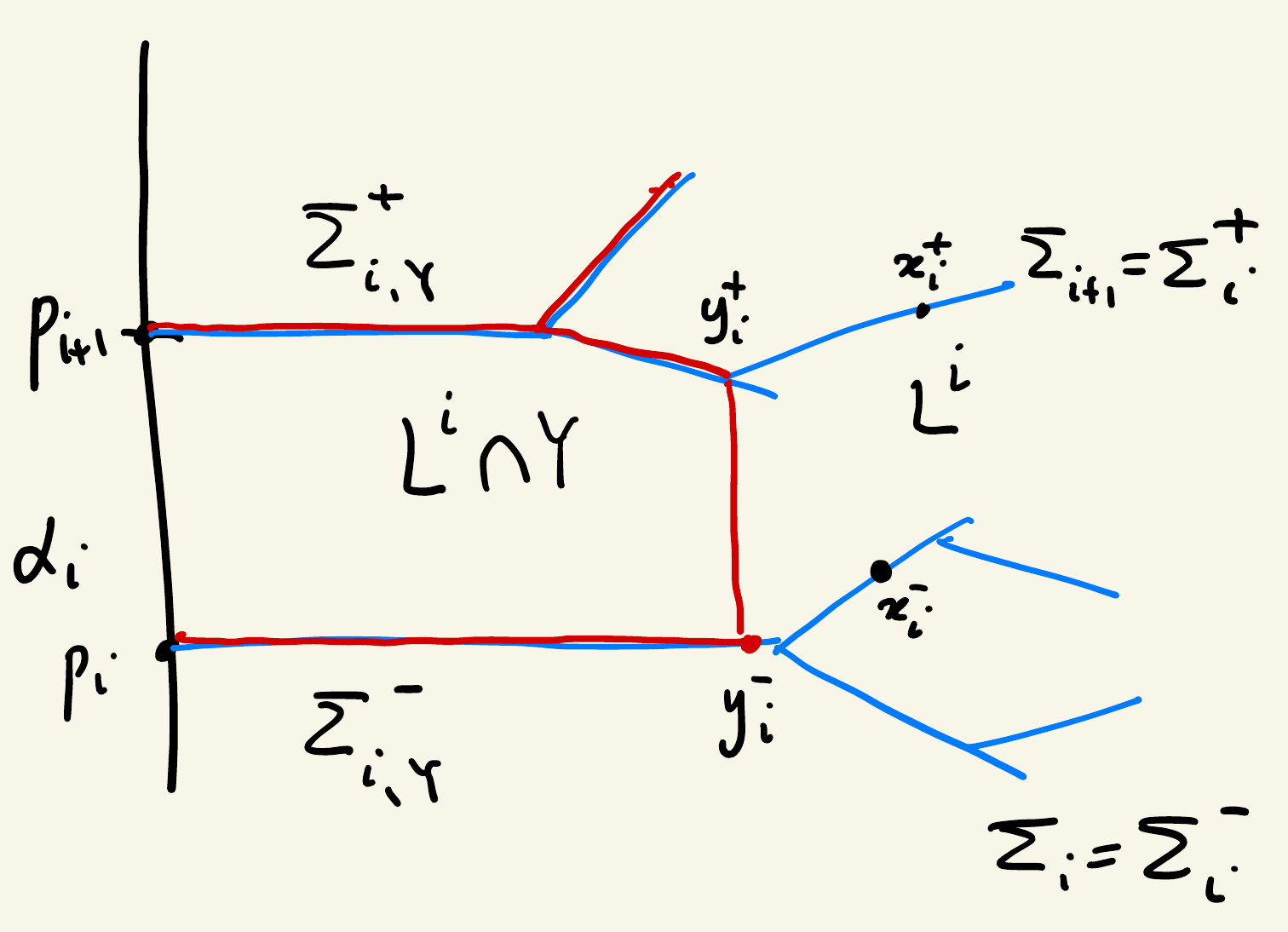}
\caption{Ladder $L^i$ and transition points.}
\label{figure30}
\end{figure}

\medskip 
We next find uniform transition points in $Y$ 
between the sections $\Si_{i,Y}$. We let $v_i^\pm= \pi(x^\pm _i)$ and  define vertices 
$w_i^\pm$ as nearest-point projections of 
 $v_i^\pm$'s to the subtrees $S_i, S_{i+1}$, where the projection is taken inside the trees 
 $T_i, T_{i+1}$ respectively.  
 Set 
 $$
 y_i^\pm := \Si^\pm_{i,Y}\cap X_{w_i^\pm}.  
 $$
 See Figure \ref{figure30}.

Since each path $c_{\L}$ connecting {\em arbitrary} points $y^\pm \in \Si_i^\pm$ is a concatenation of paths in $\Si_i^\pm$ connecting $y^\pm$ and $x_i^\pm$ and of the path $c_{\L^i}(x_i^-, x_i^+)$, we conclude that whenever $y^\pm$ are both in $Y$, the path 
$\hat{c}_{\L}(y_i^-, y_i^+)$ passes through the points $y_i^\pm$. 
Since the above paths are uniformly quasigeodesic in the ladders $L^i_Y$ (Part I.2), 
we see that the points $y_i^\pm$ are at uniformly bounded (in terms of $k$) distance $R$ from the nearest-point projections (in the ladder $L^i_Y$) of $\Si_i^-\cap Y$ to $\Si_i^+\cap Y$ and vice versa, for all $i=1,...,n-1$. For $i=n$, we have that the point $y^-_n$ within distance $R$ from the projection (in the ladder $L^n_Y$) of  $y'$ to the section $\Si_{n,Y}=\Si_n\cap Y$. 
Furthermore, by Part I.2, the paths  
$\hat{c}(y_{i-1}^+, y_i^+)$ are $\La_{\ref{lem:I.2}}$-quasigeodesics in $L_Y^i$. 
Thus, Theorem \ref{thm:hyp-tree}  implies that the alternating concatenation of the paths $\hat{c}_{\L}(y_i^-, y_i^+)$ in $L_Y^i$'s and 
of the paths $c_{\Si_{i,Y}}(y_{i-1}^+, y_i^-)$ in $\Si_{i,Y}$'s, is a $\La_{\ref{sec:StepI.3}}$-quasigeodesic in $L_Y$ connecting $y$ to $y'$.

For each vertex $v$ of $\pi(\A^i)$ we break the geodesic segment $L^i_v\subset X_v$ as a concatenation of two subsegments: $A^i_v$ (in the carpet $\A^i$) and $\beta^i_v$.

\begin{lemma}\label{lem:beta} 
For  $w= w_i^+=\pi(y_i^+)$ the length of $\beta=\beta^i_w$ is $\le C_{\ref{lem:beta}}(K)$.  
\end{lemma}
\proof If the length of $\beta$ is $\le M_{\bar{K}}$ then we are done. Otherwise, let $J$ be the largest  
subinterval in $\pi(A^i)$ containing $v$ such that for all vertices $s\in J$ the length of the subinterval $\beta_s^i$ is 
 $>M_{\bar{K}}$. Since $\beta_u^i$ has length $\le M_{\bar{K}}$ and for  $v=v_i^+=\pi(x^+_i)$ the length of $\beta^i_v$ is 
 also  $\le M_{\bar{K}}$, uniform $K$-flaring implies that the length of $J$ is $\le \tau=\tau(K, M_{\bar{K}})-2$ and, thus,  by Lemma 
 \ref{lem:growth-of-flare}
 $$
\length (\beta)\le C(k):=a^\tau  (M_{\bar{K}} + b),
 $$
where $a=L'_0$ and  $b=2L_0' K$.  \qed 

\medskip 
At this point, if we knew that for all $i$'s the vertices $\pi(x_i^+), \pi(x^-_{i+1})$ are separated in $T$ by 
the subtree $S_{i+1}= \pi(\Si^+_i)$, then the paths $c_{\Si_{i}}(x_{i-1}^+, x_i^-)$ would have to contain the subpaths $c_{\Si_{i,Y}}(y_{i-1}^+, y_i^-)$. This would imply that $\hat{c}$ is  a concatenation of the uniform quasigeodesics
$$
c_{\L_Y}(y_i^-, y_i^+) \star c_{\Si_{i,Y}}(y_{i}^+, y_{i+1}^-), i=1, 2,..., 
$$
and then we would be done with the proof of the proposition. However, this (the separation property) need not be the case. What we know, however, from the description of the points $x^\pm_i, y_i^\pm$, is that in the oriented interval $\pi(\A_i)$ either 
$$
u\le \pi(y_i^+)\le \pi(y_i^-)\le \pi(x_i^+)\le \pi(x_i^-)
$$
or 
$$
u\le \pi(y_i^+)=\pi(x_i^+) \le \pi(y_i^-)\le \pi(x_i^-). 
$$

  \begin{figure}[tbh]
\centering
\includegraphics[width=70mm]{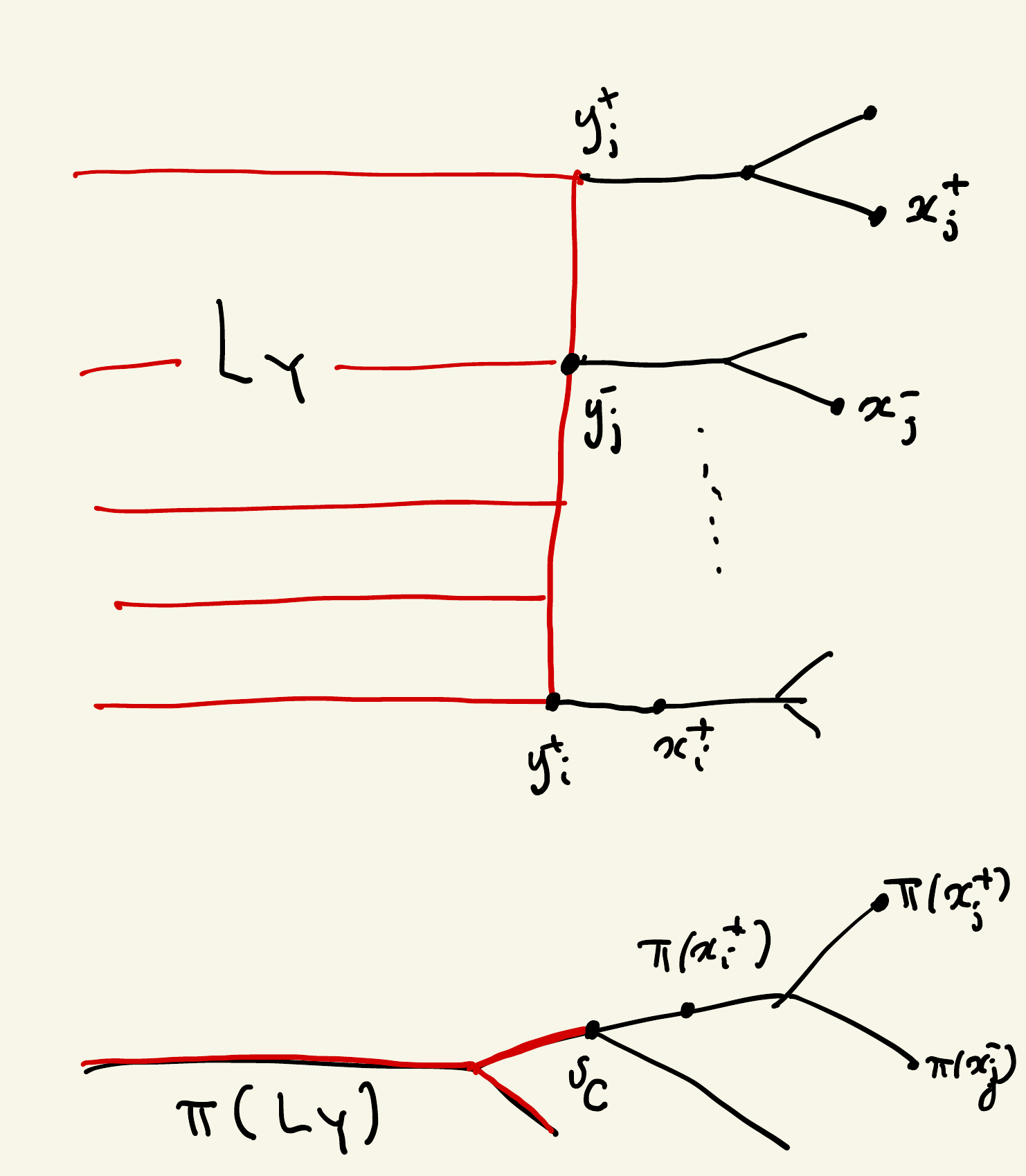}
\caption{A clique $\CC=\{x_i^+,...,x_j^-, x_j^+\}$.}
\label{figure31}
\end{figure}

Moreover, since $\L, \L_Y$ are ladders, $\pi(x_i^+), \pi(x^-_{i+1})$ are separated in $T$ by 
the subtree $S_{i+1}$ if and only if they are separated in $T$ by the subtree $\pi(L_Y)= \pi(\L)\cap S$. 
In particular, for each $i$ either $\pi(x_i^\pm)$ are in the same component of $T-\pi(L_Y)$ or $x_i^+$ lies in $Y$.  In other words, we find that $y_i^+$ is an $X$-transition point between $y, y'$ if either  $y_i^+=x_i^+$ or if $\pi(x_i^+)$ and 
$\pi(x_{i+1}^-)$ lie in distinct components of $T-\pi(\L_Y)$. Thus, we group the points $x^\pm_i$ into {\em cliques} of the form
$$
{\CC}=\{x_i^+, x_{i+1}^-, x_{i+1}^+,...,x_j^+\},
$$
or 
$$
{\CC}=\{x_i^-, x_i^+, x_{i+1}^-, x_{i+1}^+,...,x_j^+\},
$$
etc.,  consisting of maximal collections of consecutive points in $c$ of the form $x^\pm_l$, 
  whose projections to $T$ are not separated by  $\pi(L_Y)$. Here it is understood that if one of the two vertices $s, t$ lies in  a subtree $T'\subset T$, then $s$ and $t$ are separated by this subtree. We will denote a clique as 
 $$
 (x_i^\pm,...,x_j^\pm), 
 $$
 where $x_i^\pm, x_j^\pm$ are the  first and the last element of a clique $\CC$, listed in the order of their 
 appearance in the path $c=c(x,x')$.  Elements of a clique have the property that their projections to $T$ followed by the projection to $\pi(L_Y)$ equal to the same vertex $v=v_{\CC}$, and the corresponding points in $Y$ 
 $$
 y_i^\pm,...,y_j^\pm  
 $$
 all belong to the same vertical geodesic segment $L_v$, $v=v_{\CC}$. Furthermore, for each $l\in [i, j]$ 
$$
y^+_l=y^-_{l+1}, 
$$
since they both belong to $L_v\cap \Si_{l+1}$. 

\medskip 

The first and the last points $x_i^\pm, x_j^\pm$ of a clique $\CC$ determine points $y_i^\pm, y_j^\pm\in L_v\cap c$ which 
break $c$ as a concatenation of three subpaths $c_1\star c_2 \star c_3$, where $c_2= c(y_i^\pm, y_j^\pm)$. Thus, each clique defines a decomposition 
$$
\hat{c}= \hat{c}_1\star \hat{c}_2 \star \hat{c}_3, 
$$
where $\hat{c}_2= [y_i^\pm y_j^\pm]_{X_v}\subset L_v$.  The latter path is a concatenation of the subsegments 
$[y_l^- y_l^+]_{X_v}\subset L_v$. Each of these subsegments, in turn, 
 is a concatenation of two subsegments: $A^l_v$ (the narrow end of 
the carpet $\A_Y^l\subset \L_Y^i$) and a subsegment $\beta_v^l$ of length $\le C(k)$, see Lemma \ref{lem:beta}. 
Since the segment 
$A^l_v$ is a uniform quasigeodesic in $L^i_Y$ (see Lemma \ref{lem:I.1}), it follows that  
 $[y_l^- y_l^+]_{X_v}$ is a uniform quasigeodesic in $L^i_Y$.  
 
 Thus, the entire path $\hat{c}$ is broken as an alter\-nating conca\-tena\-tion of 
 uniform  $Y$-quasigeo\-desics  $[y_l^- y_l^+]_{X_v}$ connecting $\Si_{l,Y}, \Si_{l+1,Y}$ and of (possibly degenerate) horizontal $K$-qi leaves connecting  $y^+_i, y^-_{i+1}$ inside  $\Si_{l+1,Y}$. The points   $y_l^- , y_l^+$, up to a uniformly bounded error, realize the shortest distance in $L^i_Y$ between   $\Si_{l,Y}, \Si_{l+1,Y}$. 
 
 We can now finish the proof of Proposition \ref{prop:oneflow} for the slim combing paths $c$: The path $\hat{c}$ satisfies the conditions of Theorem \ref{thm:hyp-tree} and, hence, is a uniform quasigeodesic in $Y$.

\medskip 
To conclude: 

\begin{lemma}\label{lem:I.3} 
The set of pairs $(y,y')\in Fl^Y_k(X_u)\times Fl^Y_k(X_u)$ is weakly special. More precisely, there exists $\La_{\ref{lem:I.3}}(k)$ 
such that each slim combing path $c$ in $N^{fib}_{4\delta_0}Fl^X_k(X_u)$ connecting $y$ to $y'$ satisfies the property that $\hat{c}$ 
is a $\La_{\ref{lem:I.3}}(k)$-quasigeodesic in $Y$. 
\end{lemma}
\proof Points $y, y'$ belong to bottom/top of a $k$-ladder $L(\al)$ which is uniformly close to  $Fl^Y_k(X_u)$. The combing path $c=c(y,y')$ in this ladder satisfies the property 
that $\hat{c}$ is a  $\La_{\ref{lem:I.3}}(k)$-quasigeodesic in $Y$.  \qed

\medskip 
Lastly, we prove consistency for arbitrary uniform quasigeodesics $\phi$ 
connecting points of the given vertex-flow-space, i.e. uniform consistency of points in an arbitrary flow-space 
$Fl_k(X_u)$, i.e. prove Proposition \ref{prop:oneflow} in full generality. 

Since $Fl_k(X_u)$ is $\delta_{\ref{flow of one vertex space}}(k)$-hyperbolic, for each $\La\ge 1$, 
each $\La$-quasigeodesic $\phi$ in $Fl_k(X_u)$ connecting $y, y'$ is within Hausdorff distance $D(k,\La)$ 
from a combing path $c=c_L(y,y')$ contained in a ladder $\L=\L_X(\al)\subset N_{4\delta_0}^{fib} Fl_k(X_u)$, 
where $y, y'\in L_X(\al)$.

\medskip 
We will be using the notation from the proof of Lemma \ref{lem:I.3}. 
Suppose that $x\in \XX$ is a point in $c\subset L_K(\al), \al\subset X_u$, $x\notin L_Y$. Then $x$ belongs to one of the subpaths $c(y_i^\pm, y_j^\pm)$ determined by a clique $\CC$ and $\pi(x)$ is a vertex in a subtree of $T$ separated from $\pi(L_Y)$ by the vertex $v=v_{\CC}\in \pi(L_Y)$.  In particular, 
$$
d(x, X_S)= d(x, X_v)\ge d(\pi(x), v) 
$$
and taking intersection with $L_v$ of the canonical $K$-qi section $\Si_x\subset L_X(\al)$, we obtain a point $y''\in L_v\subset L_Y$ within distance $K d(x, X_S)$ from $x$. It follows that every point $z\in X_S$ within distance $R$  from $x\in c$, lies within distance 
$(K+1) R$ from a point $\hat{z}=y''$ in $\hat{c}\cap L_v$.  In particular, each intersection point of $\phi$ with $X_S$ is within distance 
$D(k,\La)+ (K+1) R$ from a point  in $\hat{c}\cap L_v$.  Furthermore, by the construction, the map $z\mapsto \hat{z}$ is monotonic: If $z_1$ appears before $z_2$ in $c$, then $\hat{z}_1$   appears before $\hat{z}_2$ in $\hat{c}$. 
It now follows that for each $L$-quasigeodesic $\phi$, the path $\hat{\phi}$ is $\hat{\La}$-quasigeodesic in $Y$. \qed 

\medskip 
 This concludes Part I of the proof of Theorem \ref{thm:cut-paste}.

\medskip 
The next result is an immediate corollary of Proposition \ref{prop:oneflow}:

\begin{cor}\label{cor:detour} 
Suppose that $p, p'$ are point in a vertex-space $X_v\subset Y$ such that a $\La$-quasigeodesic $\phi$ in $X$ connecting $p$ to 
$p'$ intersects $X_v$ only at its end-points. Then the vertical geodesic $[pp']_{X_v}$ is a 
$\La'= \La'_{\ref{cor:detour}}(\La)$-quasigeodesic in $Y$. 
\end{cor}

As another application of Part I, we also obtain a theorem which is essentially due to Mitra, \cite{mitra-trees}:

\begin{thm}\label{thm:connecting points of a vertex space}
There exists a constant $R=R(\La)$ depending only on the tree of spaces $X\to T$ such that for every vertex $v\in V(T)$ and any pair of points $y, y'\in X_u$,  and any $\La$-quasigeodesic $\phi$ in $X$ with end-points 
$y, y'$, the intersection
$$
X_u\cap  \phi 
$$
is contained in the $R$-neighborhood of the vertical geodesic $[yy']_{X_u}$. 
\end{thm} 
\proof We will apply Proposition \ref{prop:oneflow} to the subtree $S=\{u\}$. By the proposition, each intersection point $z\in X_u\cap  \phi $ is within distance $D(k,\La)+ (K+1) R$ from a point  in the geodesic segment 
$L_u\subset [yy']_{X_u}$, where $K=K(k), R=R(k,\La)$ and $k$ can be taken to be uniform, say, $k=K_0$. \qed  

\begin{figure}[tbh]
\centering
\includegraphics[width=100mm]{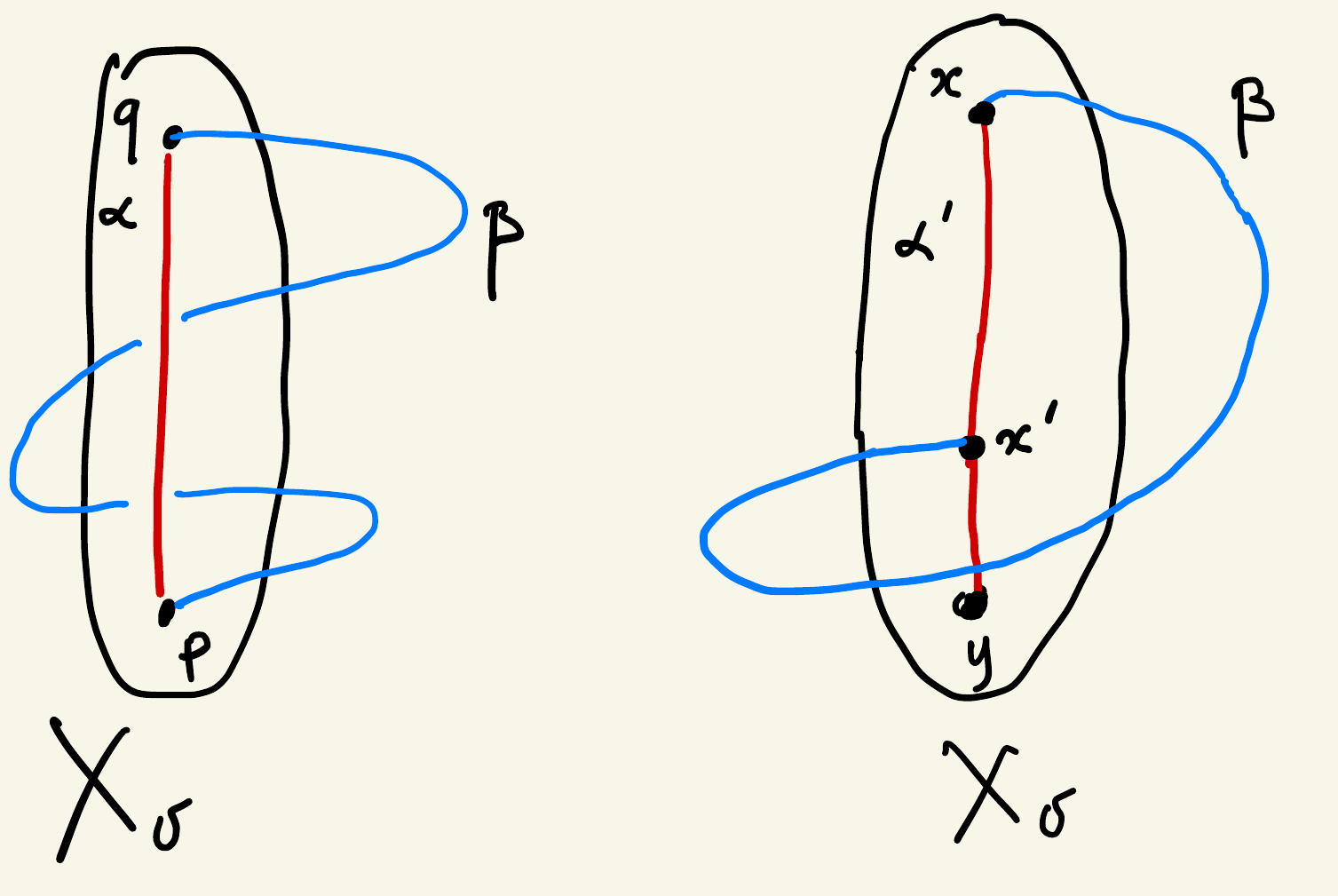}
\caption{}
\label{figure14}
\end{figure}

\begin{cor}\label{cor:geodesic_crossing_one_vertex_space}
There exists $R=R_{\ref{cor:geodesic_crossing_one_vertex_space}}(r)$ 
such that for every vertex $v\in V(T)$ and any pair of points $p, q\in X_v$ and geodesics $\al=[pq]_{X_v}$, $\beta=[pq]_{X}$, 
we have  
$$
N_r(X_v)\cap \beta  \subset N_R (\al). 
$$
In particular, if $\al'=[xx']_{X_v}\subset \al=[xy]_{X_v}$ and $y\in N_r(\beta)$, $\beta=[xx']_X$, then $d(y, x')\le R$. (See Figure \ref{figure14}.) 
\end{cor}

As another corollary, we obtain Mitra's theorem on the existence of CT-maps (Theorem \ref{thm:mitra-trees}): 

\medskip
{\em Proof of  Theorem \ref{thm:mitra-trees}.}  We will use Mitra's Criterion for the inclusion map $f: X_u\to X$ (Proposition \ref{prop:Mitra's Criterion1}). We prove that $f$ satisfies the assumption of Proposition \ref{prop:Mitra's Criterion1} by arguing the contrapositive. 

Fix a base-point $y\in Y$ and consider  a sequence of vertical geodesics $\al_n=[x_n x'_n]_{X_u}$. 
Assume that each $X$-geodesic $\beta_n=[x_n x'_n]_X$ has nonempty intersection with the $r$-ball $B(y,r)\subset X$ for 
some fixed $r$. After replacing  $\beta_n$'s with uniform quasigeodesics $\phi_n$ in $X$, we ensure  that $y\in \phi_n$ for all $n$. 
Therefore, $y\in \hat{\phi}_n$  as well. Since the Hausdorff distance between $\hat{\phi}_n$ and $\al_n$ is uniformly bounded, 
the minimal distances between $y$ and $\al_n$'s are uniformly bounded too. Thus, the assumption of 
Proposition \ref{prop:Mitra's Criterion1}  is satisfied and, hence,    the inclusion map $f: X_u\to X$ admits a  CT-extension. \qed

\medskip 
 In Part II of the proof of Theorem \ref{thm:cut-paste} we will need several technical results regarding projections to $X$-geodesics connecting points in $X_u$, $u\in V(S)$.  These results occupy the rest of this section. We consider a $(K,D,E)$-ladder $\L=\L_X(\al)$, where $p, p'\in X_u$ and 
 $\al= [p p']_{X_u}$.  We will be investigating the nearest-point projection (taken in $X$) 
of points $y\in \al= [pp']_{X_u}$ to a $\La$-quasigeodesic geodesic $\phi$ in $X$ containing a detour path $\zeta=\phi(p, p')$ 
connecting the points $p, p'$. Most of the discussion deals with the case $\phi=\zeta$. Since the path $\zeta$ is uniformly Hausdorff-close to the combing path $c=c_X(p,p')$, in order to understand the projection of $y$ to $\zeta$, it suffices to analyze the projection of $y$ to $c$ (see Corollary \ref{cor:proj-to-close-subsets}).

\begin{lemma}\label{lem:proj-to-detour} 
Suppose that $p, p'$ belong to a common vertex-space $X_u$.   
Let $\{\al_i\}$ be a vertical subdivision of $\al=[pp']_{X_u}$. Then:

1. For 
each $y\in \al_i=[p_i p_{i+1}]_{X_u}$ the projection $\bar{y}= P_{X,c}(y)\in c= c_X(p,p')$ is 
uniformly close to a point $\bar{y}'$ in the subladder $L^i=L_K(\al_i)$ determined by $\al_i$:
$$
d_X(\bar{y}, \bar{y}')\le C=C_{\ref{lem:proj-to-detour}}(K, D, E, \La). 
$$

2. The point $\bar{y}'$ can be chosen to lie in the carpet  $A^i=A(\al_i')\subset L^i$, where $\al'_i\subset \al_i$ is as in Proposition \ref{vertical subdivision}. 

3. The point $\bar{y}'\in A^i$ can be chosen so that  
$y, \bar{y}'$ are connected by a (canonical in $A^i\subset L^i$) $K$-qi leaf contained in $A^i$. 
 \end{lemma}
\proof 1. Connect $y$ to a point $z\in c\cap L^i$ by a geodesic $[yz]_X$. 
Since $L^i$ is $\la$-quasiconvex in $X$, $[yz]_X$ lies in the $\la$-neighborhood of $L^i$. On the other hand, by Lemma \ref{lem:projection-1},  $[yz]_X$, as any geodesic connecting $y$ to $c$, has to pass uniformly close to $\bar{y}$, namely, within distance  $\la'+3\delta_X$, where $\la'$ is the quasiconvexity constant (in $X$) of the path $c$. 
 It follows that $\bar{y}$ lies distance $C=\la+\la'+3\delta_X$ 
 from a point $\bar{y}'\in L^i$.

\medskip 

2. By Part (1), $\bar{y}\in c$ lies within distance $C$ from a point $\bar{y}'\in L^i$. If $\bar{y}\notin L^i$, take the 
smallest subpath in $c$ connecting $\bar{y}$ to a point $z\in L^i$. Then (by the construction of the path $c$) the point 
$z$  realizes (up to a uniform additive error) the minimal distance from $\bar{y}$ to $L^i$. Hence, by replacing $C$ with another uniform constant $C'$, we can assume that  $\bar{y}'\in c_i:=c\cap L^i$, $d_X(\bar{y}, \bar{y}')\le C'$.

  \begin{figure}[tbh]
\centering
\includegraphics[width=100mm]{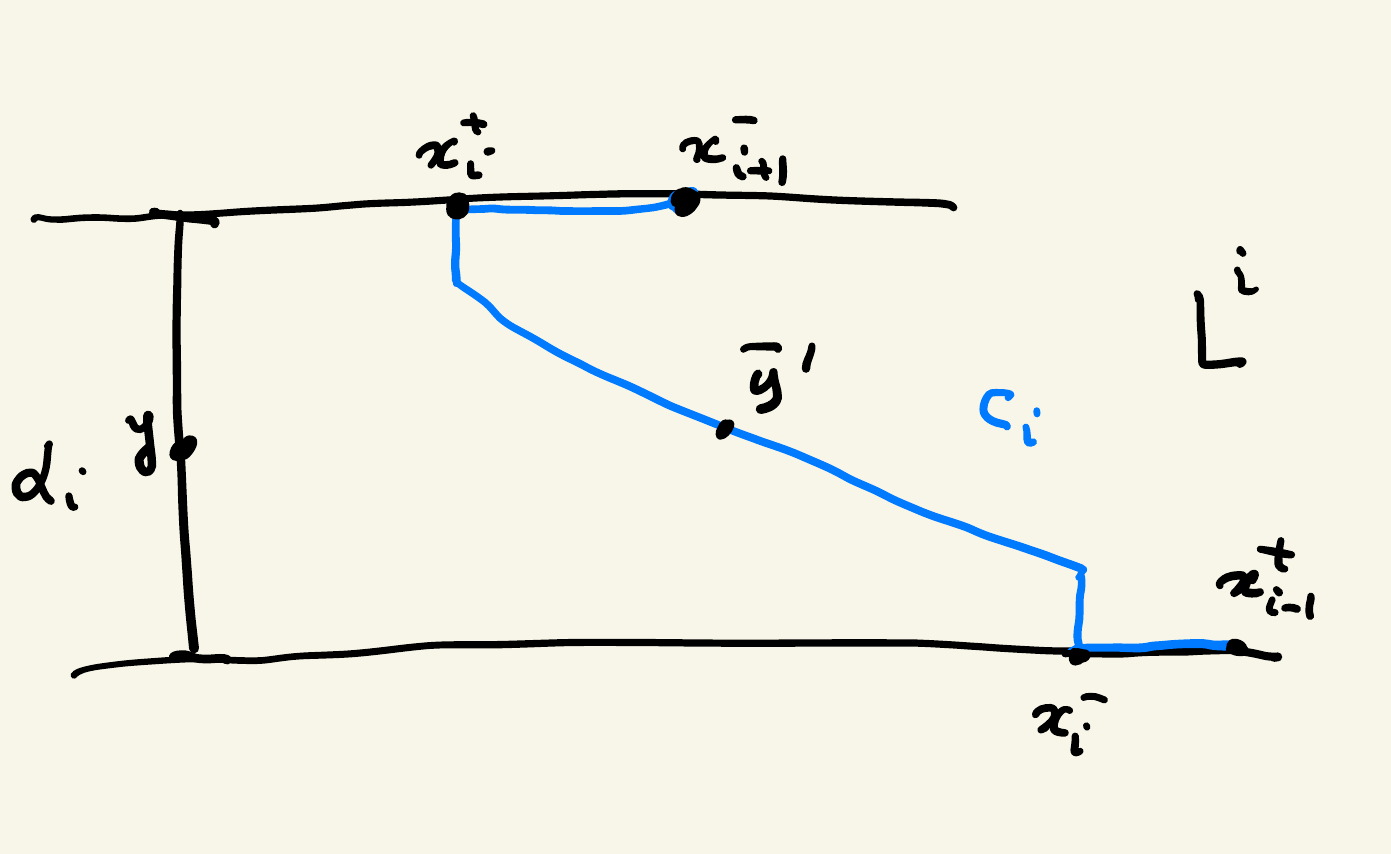}
\caption{Path $c_i$}
\label{c_i.fig}
\end{figure}

Note that $c_i$ is the concatenation 
$$
c(x_{i-1}^+, x_i^-)\star c(x_i^-, x_i^+) \star c(x_i^+, x_{i+1}^-), 
$$
see Figure \ref{c_i.fig}. The middle subpath lies in $A^i$ (except for a vertical subpath of length $\le M_{\bar{K}}$); 
therefore, there are just two cases we have to consider:

(a) $\bar{y}'\in c(x_{i-1}^+, x_i^-)\setminus A^i$. Then, by the description of combing paths in $L^i$ (see Section \ref{sec:inductive}, 
Step I.2), the path $c(\bar{y}', y)$ contains 
a subpath contained in $bot(L^i)\cap c(x_{i-1}^+, x_i^-)\subset c_i\subset c$ 
and connecting $\bar{y}'$ to a point $A^i$. Since $\bar{y}'$ is uniformly close to the projection of $y$ to $c$ 
(and since $c(\bar{y}', y)$ is uniformly quasigeodesic), the length of that subpath of $c(\bar{y}', y)$ is uniformly bounded 
and, hence, $\bar{y}'$ is uniformly close to a point in $c_i\cap A^i$, as required by the second statement of the lemma.

\begin{figure}[tbh]
\centering
\includegraphics[width=100mm]{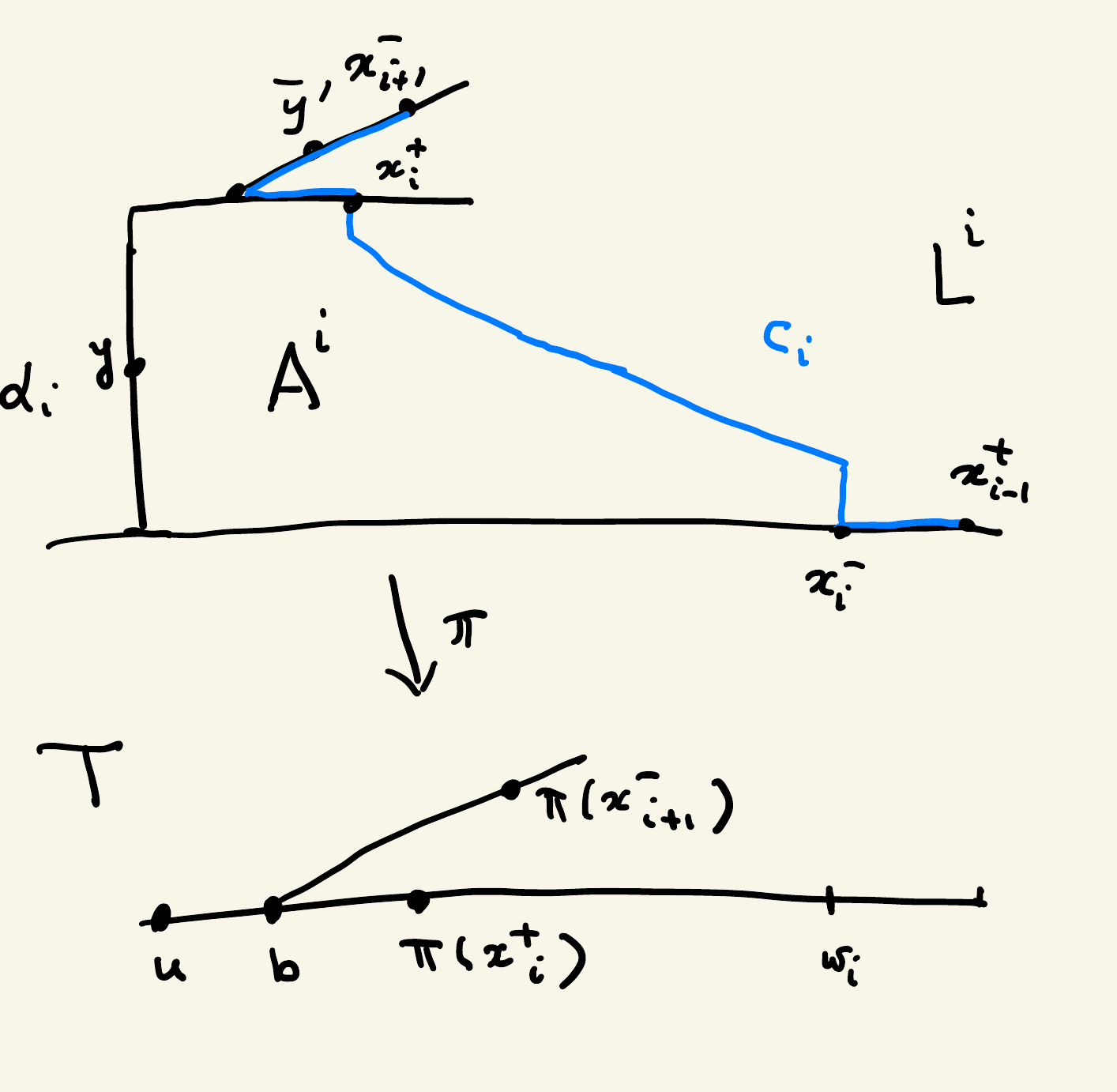}
\caption{Path $c_i$}
\label{c_i2.fig}
\end{figure}

\medskip 
(b) $\bar{y}'\in c(x_i^+, x_{i+1}^-) \setminus A^i$, see Figure \ref{c_i2.fig}. The proof is similar to Case (a): 
The path $c(\bar{y}', y)$ starts with a subpath $c_{\bar{y}'}\subset L^i$ connecting $\bar{y}'$ to $A^i$ 
(again, see Section \ref{sec:inductive}, Step I.2). Since the length of the part of 
 $c_{\bar{y}'}$ contained in $c_i$ has to be uniformly bounded (as in Case (a)), by changing the location of 
 $\bar{y}'\in c_i\cap top(L^i)$ by a uniformly bounded amount, 
 we can assume that for 
 $$
 b= \barycenter(\Delta u \pi(x_i^+) \pi(x_{i+1}^-))= \barycenter(\Delta u w_i \pi(x_{i+1}^-)),
 $$
 the vertex $v= \pi(\bar{y})$ lies in the interval 
 $$
 \llbracket b, \pi(x_i^+) \rrbracket \subset    \llbracket u, \pi(x_i^+)\rrbracket, 
 $$
 where $\pi(A^i)= \llbracket u, w_i\rrbracket$. Therefore, by Lemma \ref{lem:top-projection}(2), the 
  vertical distance from $\bar{y}'$ to the top of $A^i_v$ is bounded by $R_{\ref{lem:top-projection}}(K)$. 
 This concludes the proof of (2).

\medskip 
3. Thus, we assume that $\bar{y}'\in c_i\cap A^i$. If $y\in \al_i \setminus \al'_i$, then $\bar{y}'\in top(A^i)$ and $c(\bar{y}', x_i^+)$ 
has uniformly bounded length and, without loss of generality, $c(\bar{y}', x_i^+)$  is a vertical segment of length $\le M_{\bar{K}}$. The canonical $K$-section $\Si_y\subset L^i$ crosses the vertical interval  $c(\bar{y}', x_i^+)$ at some point $\bar{y}''$ and 
$d(\bar{y}', \bar{y}'')\le M_{\bar{K}}$. Thus, from now on, we assume that $y\in \al'_i\subset A^i$. 

The intersection $c_i\cap A^i$ consists of a path in $top(A^i)$, the narrow end of $A^i$ (which has length $\le M_{\bar{K}}$) 
and of a path contained in $bot(A^i)$ which equals 
$$
c(x_{i-1}^+, x_i^-)\cap bot(A^i). 
$$
Since every point in the narrow end of $A^i$ is connected to each point of $\al'_i$  by a $K$-qi section, we have two cases to consider: 

\medskip 
(a) 
$$
\bar{y}'\in c(x_{i-1}^+, x_i^-)\cap bot(A^i). 
$$ 
By the definition of the path
$$
c(\bar{y}', y),
$$
it is a concatenation of a horizontal subpath contained in $c(x_{i-1}^+, x_i^-)\cap bot(A^i)$, followed by a vertical path of length $\le M_{\bar{K}}$ and, then by a $K$-qi leaf
$$
\ga_{\bar{y}'', y}\subset A^i, 
$$
connecting $\bar{y}''$ to $y$. (Note that $\bar{y}''$ does not, in general, lie in $c$, it is just within distance $\le M_{\bar{K}}$ from a point in $c$.) The first horizontal subpath is entirely contained in $c_i\subset c$, hence, has to have uniformly bounded length.  Thus, 
$\bar{y}'$ is uniformly close to the point $\bar{y}''\in A^i$ connected to $y$ by a $K$-qi section.  This concludes the proof in case (a). 

(b) 
$$
\bar{y}'\in c(x_{i-1}^+, x_i^-)\cap top(A^i). 
$$ 
This case is similar to (a): The path  $c(\bar{y}', y)$ is a concatenation of a subpath contained in $c$, followed by a short vertical subpath and, then, by a $K$-qi leaf connecting a point $\bar{y}''\in A^i$ to $y$. The first two subpaths are uniformly short, hence, we are done. \qed

\begin{cor}\label{cor:proj-to-detour}
If $d(\bar{y}, p)\le R$, then
$$
d_X(y,p)\le R'=R'_{\ref{cor:proj-to-detour}}(R, K, C), 
$$  
where $C=C_{\ref{lem:proj-to-detour}}$. 
\end{cor}
\proof Taking  $\bar{y}'$ is as in Part (3) of Lemma \ref{lem:proj-to-detour}, we obtain 
$$
d_T(u, \pi(\bar{y}'))\le d_X(p, \bar{y}')\le R +C. 
$$
Since $y, \bar{y}'$ are connected by a $K$-qi leaf over the interval $\llbracket u, \pi(\bar{y}')\rrbracket$,   
$$
d_X(y, \bar{y}')\le K d_T(u, \pi(\bar{y}'))\le  K(R +C), 
$$
which, in turn, implies the inequality $d_X(y,p)\le (K+1)(R +C)$. \qed 

\medskip
Lastly, we turn to the projection of $y\in [pp']_{X_u}$ to a $\La$-quasigeodesic $\phi$ (with end-points in $Y$) containing $\zeta=\phi(p, p')$ as a detour subarc:

\begin{lemma}\label{lem:proj-detour-line}
$$
d(P_{X,\phi}(y), \zeta)\le D_{\ref{lem:proj-detour-line}}= D_{\ref{lem:proj-detour-line}}(K,D, E, \La).$$ 
\end{lemma}
\proof The result is an application of Corollaries \ref{cor:projection-2} and \ref{cor:proj-to-detour}. 
The subsets $U=\phi, V=\zeta$ are $D_{\ref{Morse}}(\delta_X,\La)$-quasiconvex in $X$. 
Assume that  $\bar{y}:= P_{X,\phi}(y)$ is not in the arc $\zeta$. 
In view of the Morse Lemma (Lemma \ref{Morse}), for each $z\in  \zeta$ the geodesic $[\bar{y} z]_X$ passes within distance $D_{\ref{Morse}}(\delta_X,\La)$ 
from $p$ or from $p'$ (depending which component of $\phi \setminus \zeta$ the point $\bar{y}$ belongs to). We will assume that it is $p$ rather than $p'$. 
Therefore, 
by Lemma \ref{lem:projection-1}(i),  the projection of $\bar{y}$ is within distance 
$2D_{\ref{Morse}}(\delta_X,\La)$ from $p$ or from $p'$. But, according to Corollary  \ref{cor:projection-2},
$$
d_X(P_{U,V}\circ P_{X,U}(y), P_{X,V}(y))\le C_{\ref{cor:projection-2}}(\delta_X,D_{\ref{Morse}}(\delta_X,\La)). 
$$
Hence, 
$$
d_X(P_V(y), p) \le  
R:=C_{\ref{cor:projection-2}}(\delta_X,D_{\ref{Morse}}(\delta_X,\La)) + 2D_{\ref{Morse}}(\delta_X,\La).  
$$
Applying Corollary \ref{cor:proj-to-detour}, we get
$$
d_X(y, p)\le R'_{\ref{cor:proj-to-detour}}(D),
$$
and, therefore, $d_X(\bar{y}, p)\le D_{\ref{lem:proj-detour-line}}:=2R'_{\ref{cor:proj-to-detour}}(R)$. \qed 

\begin{cor}\label{cor:proj-detour-line}
$$
d(P_{X,\phi}(y), P_{X,\zeta}(y))\le D_{\ref{cor:proj-detour-line}}(K,D, E, \La).$$ 
\end{cor}
\proof If $P_{X,\phi}(y)\in \zeta$ then $P_{X,\phi}(y)= P_{X,\zeta}(y)$. 
Assume, therefore, that 
$P_{X,\phi}(y)\notin \zeta$. According to lemma, $P_{X,\phi}(y)$ is within distance $D_{\ref{cor:proj-detour-line}}$ from 
a point $q\in \zeta$. Therefore,
$$
d_X(y, q)\le d(y, \zeta)+ D_{\ref{cor:proj-detour-line}}. 
$$
By Corollary \ref{cor:almost-npp},
$$
d_X(q, P_{X,\zeta}(y))\le D_{\ref{cor:proj-detour-line}}+ 2\la + 4\delta_X,
$$
where $\la=D_{\ref{Morse}}(\delta_X,\La)$ is the quasiconvexity constant of $\zeta\subset X$. \qed

\section{Part II: Consistency in semispecial flow-spaces} \label{sec.II}

\subsection{Part II.4: Projections in special flow-spaces} 

In this part of the proof (which is the most difficult part of Section \ref{sec.II}) we are assuming that vertices 
$u, v\in V(S)$ define a $K$-special interval $J= \llbracket u, v\rrbracket\subset S$.  In order to simplify the notation we set
$$
F^X_w= Fl^X_{K}(X_w), \quad F^Y_w= Fl^Y_{K}(X_w), 
$$
for vertices $w\in V(S)$. Observe that $F^X_v\cap X_u\ne \emptyset \iff F^Y_v\cap X_u\ne \emptyset$ and 
$F^X_u\cap X_v\ne \emptyset \iff F^Y_u\cap X_v\ne \emptyset$, i.e. it does not matter in what space ($X$ or $Y$) our interval $J$ is special. In the proofs below, it does not matter which of the above intersections is nonempty.  
In this section we will relate the $X$-projection $\bar{x}$ of $x\in F^Y_v$ to $F^X_u$ and the $Y$-projection $\bar{\bar{x}}$ of $x$ to $F^Y_u$. At the same time, we will establish uniform consistency of some classes of pairs $(x, x')\in F^Y_v\times F^Y_u$. These results will be the key for proving uniform consistency of pairs of points in semispecial flow-spaces, which will be done in the next section. (This will apply, of course, to points in special flow-spaces as well.) After altering $\bar{x}$ and $\bar{\bar{x}}$ by uniformly bounded distance, we can assume that these points belong to vertex-spaces of $X$ and $Y$ respectively.

Recall (see Lemma \ref{lem:projection-1}) that the concatenation
$$
\phi= [x \bar{x}]_X\star [\bar{x} x']_X
$$
is a $\La$-quasigeodesic in $X$, where $\La$ depends only on the quasiconvexity constant $\la_X$ of $F^X_u$ and the hyperbolicity 
constant of $X$ (i.e. only on the parameters of $\X$ and $K$).  If $\bar{x}\notin Y$, we let 
$\zeta=\phi(\hat{x}, \hat{x}')$ denote the detour subpath in $\phi$ containing $\bar{x}$ and connecting 
points $\hat{x}, \hat{x}'$ which belong to the same vertex-space $X_t\subset Y$.  
In the case $\bar{x}\in Y\cap X_t$, we declare $\zeta$ to be  degenerate, equal to the singleton $\{\bar{x}\}$, and, accordingly, set 
$$
\hat{x}= \hat{x}'= \bar{x}. 
$$

\begin{rem}
Even if $\bar{x}$ is far away from $Y$, the points $\hat{x}, \hat{x}'$ depend quite a bit on the quasigeodesic $\phi$ 
(i.e. on the choice of geodesics $[x \bar{x}]_X, [\bar{x} x']_X$) 
connecting $x, x'$. The notation $\hat{x}, \hat{x}'$, therefore, is ambiguous and, in truth, should contain the symbol $\phi$, which we omit for the ease of the notation.   
\end{rem}

\begin{figure}[tbh]
\centering
\includegraphics[width=80mm]{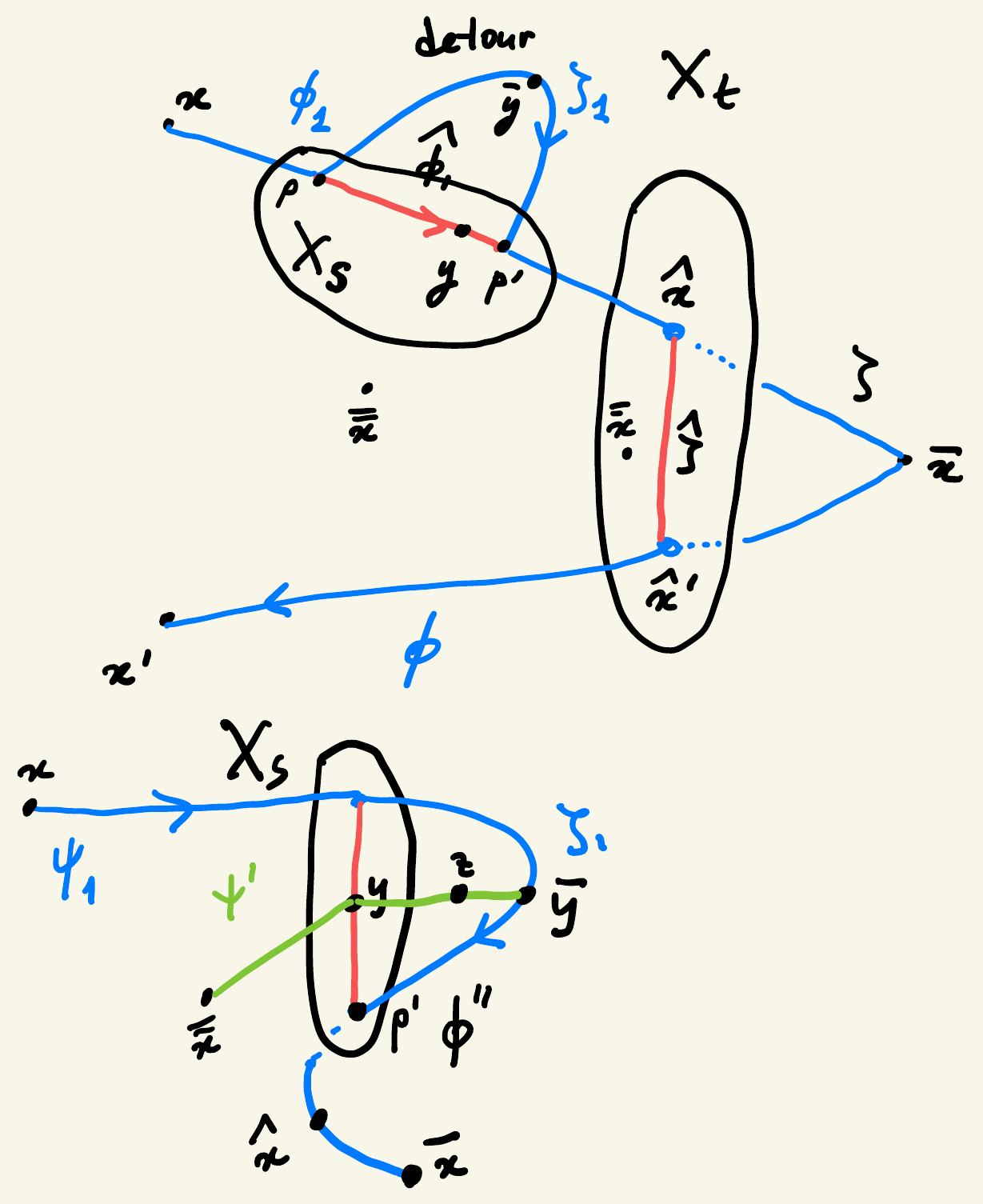}
\caption{Two projections.}
\label{two-projections.fig}
\end{figure}

The key result of Part II.4 is:

\begin{lemma}\label{lem:key-II.4}
The point $\bar{\bar{x}}$ lies within distance $R_{\ref{lem:key-II.4}}(K,\La)$ 
from a point of the segment $\hat{\zeta}=[\hat{x} \hat{x}']_{X_t}$, no matter what $x'$ and $\phi$ are.  
\end{lemma}
\proof Consider the subpath $\phi'=\phi(x,\hat{x}')$ connecting $x$ to $\hat{x}'$ and passing through 
$\hat{x}$: 
$$
\phi'= \phi_1\star \zeta. 
$$
(The subpath $\phi_1$ is geodesic since it is contained in $[x\bar{x}]_X$.) See Figure \ref{two-projections.fig}. 
Then 
$$
\hat{\phi}'= \hat{\phi}_1\star \hat{\zeta}.
$$  
We already know (by Proposition \ref{prop:oneflow} and Corollary \ref{cor:detour}) that $\hat{\phi}_1$ and 
$\hat{\zeta}$ are uniform $Y$-quasigeodesics, but we do not yet know that their concatenation is, since uniform consistency of the pairs $(x, \hat{x}')$ is not yet known, except in some special cases. We do know, however, the 
uniform consistency of the pair 
$(x, \bar{\bar{x}})$: Since $F^Y_v\cap F^Y_u\ne \emptyset$, 
Lemma \ref{lemma0-flow-space} implies that the point $\bar{\bar{x}}$ lies in the $2\la_Y+3\delta_Y$-neighborhood of $F^Y_v$, 
where $\la_Y$ is the quasiconvexity constant of flow-spaces $Fl_w^Y\subset Y$. 
Therefore, the uniform consistency of the pairs $(x, \bar{\bar{x}})$ follows from Proposition \ref{prop:oneflow}.   

 Since the geodesic $[x \hat{x}']_Y$ passes uniformly close to the projection $\bar{\bar{x}}$ of $x$ to $F^Y_u$, 
by the $\delta_Y$-slimness of the geodesic triangle
$$
\Delta_Y x \hat{x} \hat{x}'\subset Y,
$$
either $\hat\phi_1$ or $\hat\zeta$ passes within the 
distance $R=R(K)$ 
from $\bar{\bar{x}}$. If it is $\hat\zeta$, we are done. Consider, therefore,  the case 
that $\hat\phi_1$ passes within distance $R$ from $\bar{\bar{x}}$. The path $\hat\phi_1$ is a concatenation of subarcs $\phi_1\cap Y$ with vertical geodesics $\hat\zeta_i$, where each $\zeta_i$ is a detour subarc in $\phi_1$.

Suppose first that $\phi_1\cap Y$ passes within distance $R$ from $\bar{\bar{x}}\in F^Y_u$. 
Since $\phi_1$ is geodesic in $X$ and  $\bar x$ is the nearest-point projection of $x$ to 
$F^X_u$, we obtain:
$$
d(\bar{\bar{x}}, \bar{x})\le 2R. 
$$ 
In particular, $d(\bar{x}, X_t)\le 2R$ as well (since $[\bar{x} \bar{\bar{x}}]_X$ has to cross $X_t$).  
By Corollary \ref{cor:geodesic_crossing_one_vertex_space}, the intersection 
$$
B(\bar{x}, 2R)\cap X_t$$ 
is uniformly close to a subset of $\hat\zeta$ and, hence, we are done, in this case as well.

It remains to analyze the harder case when for one of the detour subarcs $\zeta_1\subset \phi_1$, the path 
 $\hat\zeta_1$ passes within distance $R$ from $\bar{\bar{x}}$. Let $X_s$ denote the vertex-space of $Y$ containing the end-points $p, p'$ of $\zeta_1$ and $y\in \hat\zeta_1\subset X_s$ a point within distance $R$ from  $\bar{\bar{x}}$. 
Let $\bar{y}$ denote the nearest-point projection of $y$ to $\zeta_1$.

The point $\bar{y}$ lies within distance $D_{\ref{lem:proj-detour-line}}$ from the projection of $y$ 
to $\phi_1$ since the latter contains $\zeta_1$, see Lemma \ref{lem:proj-detour-line}. 
The concatenation $\psi=\psi_1 \star [\bar{y} {y}]_X$ 
(where $\psi_1$ is the subpath of $\phi_1$ between $x$ and $\bar{y}$) 
is a uniform quasigeodesic in $X$ connecting $x$ to the point $y$ within distance $R$ from $\bar{\bar{x}}$. 

Since $F_v^X$ is $\la_X$-quasiconvex in $X$,  $[y \bar{\bar{x}}]_X$ has length $\le R$ and $\psi$ is a uniform quasigeodesic in $X$, 
the path $\psi'=\psi\star [y \bar{\bar{x}}]_X$ (connecting $x$ to $\bar{\bar{x}}$) lies within a uniform neighborhood of $F^X_v$. 
Since  $\bar{\bar{x}}$ is in the $2\la_Y+3\delta_Y$-neighborhood of $F^Y_v\subset F^X_v$, 
by Lemma \ref{lem:projection-1}(i) the path $\psi'$  (connecting $x\in F^X_v$ to $\bar{\bar{x}}$ which is 
uniformly close to $F^X_v$) passes, at some point $z\in \psi'$, uniformly close to $\bar{x}$.

Where could this point $z$ be? If $z$ lies in $\psi_1$, then, since $\bar{y}$ is between $z$ and $\bar{x}$ in $\phi$,  
the length of  the entire subpath $\phi''$ of $\phi$ between $\bar{x}$ and $\bar{y}$ is uniformly bounded. 
Similarly, if $z$ lies in $[y \bar{\bar{x}}]_X$ then,  since the concatenation of $\phi''$ with  $[\bar{y} y]_X$ is also a uniform quasigeodesic, the length of $\phi''$ is uniformly bounded as well. The path $\phi''$ includes a subpath (contained in $\zeta_1$) 
between $\bar{y}$ and $p'$, and we conclude that the distance between $\bar{y}$ and $p'$ is uniformly bounded by some constant 
$D$.  Therefore, by Corollary \ref{cor:proj-to-detour}, the distance between $y$ and $p'$ is 
$\le   D'=D'_{\ref{cor:proj-to-detour}}(D)$, implying a uniform upper bound on the distance from $\bar{\bar{x}}$ to $\hat{x}\in \phi''$. Therefore, in this case again, we see that $\bar{\bar{x}}$ is uniformly close to a point (namely, $\hat{x}$)  in the segment 
$$
\hat{\zeta}= [\hat{x} \hat{x}']_{X_t}. 
$$
Of course, in this case the distance between $\bar{\bar{x}}$ and $\bar{x}$ is also uniformly bounded as well. \qed

\medskip
In fact, we can pin down the location of a point $y\in \hat\zeta$ within distance 
$R_{\ref{lem:key-II.4}}(K)$ from  $\bar{\bar{x}}$ a bit further. Namely, the set of points $\hat{x}'$ in the lemma is precisely the $4\delta_0$-quasiconvex subset $Q_t\subset X_t$ equal to the intersection 
$$
X_t\cap Fl_u^X. 
$$
Since $\bar{\bar{x}}$ is uniformly close to a point in each of the geodesics $[\hat{x} \hat{x}']_{X_t}$, $\hat{x}'\in Q_t$, 
by applying Corollary \ref{cor:projection-1} we conclude that $\bar{\bar{x}}$ is uniformly close to a point in the geodesic segment 
$[\hat{x} y]_{X_t}$, where $y$ is the projection (taken in $X_t$) of $\hat{x}$ to $Q_t$. 

 At the same time, since $\bar{\bar{x}}$ belongs to $F^Y_u$,
 $y$ lies in the intersection of $X_t$  with the $R_{\ref{lem:key-II.4}}(K)$-neighborhood of $F^Y_u$ 
 (the neighborhood is taken in $X$). Thus, by Lemma \ref{lem:intersection-with-flow-nbd}, the point $y$ lies in 
$$
N^{fib}_{D}(Q_t)\subset X_t, 
$$
where $D=D_{\ref{lem:intersection-with-flow-nbd}}(R,K)$. We conclude:

\begin{cor}
The point $\bar{\bar{x}}$ is uniformly close to the projection $\tilde{x}$ (taken in $X_t$) of $\hat{x}$ to 
$Q_t=X_t\cap Fl_K^X(X_u)$. 
\end{cor}

\begin{figure}[tbh]
\centering
\includegraphics[width=80mm]{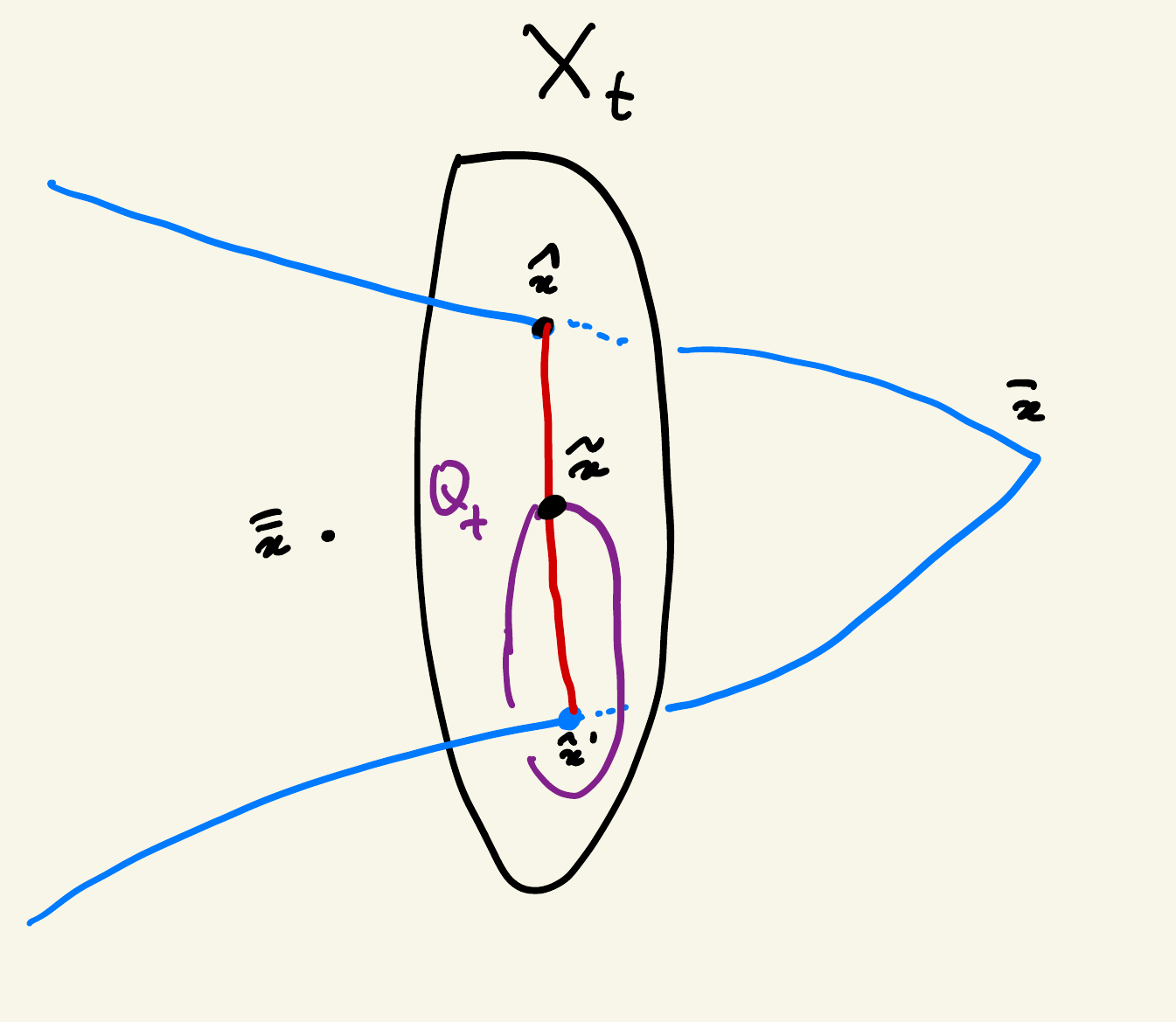}
\caption{Location of $\bar{\bar{x}}$.}
\label{two-projections1.fig}
\end{figure}

Note that, while the point $\hat{x}$ depends heavily on the choice of the path $\phi$ connecting $x, x'\in F^Y_u$,  
the point $\tilde{x}$ is canonical (up to a uniform error, depending only on $K$).  

\begin{figure}[tbh]
\centering
\includegraphics[width=80mm]{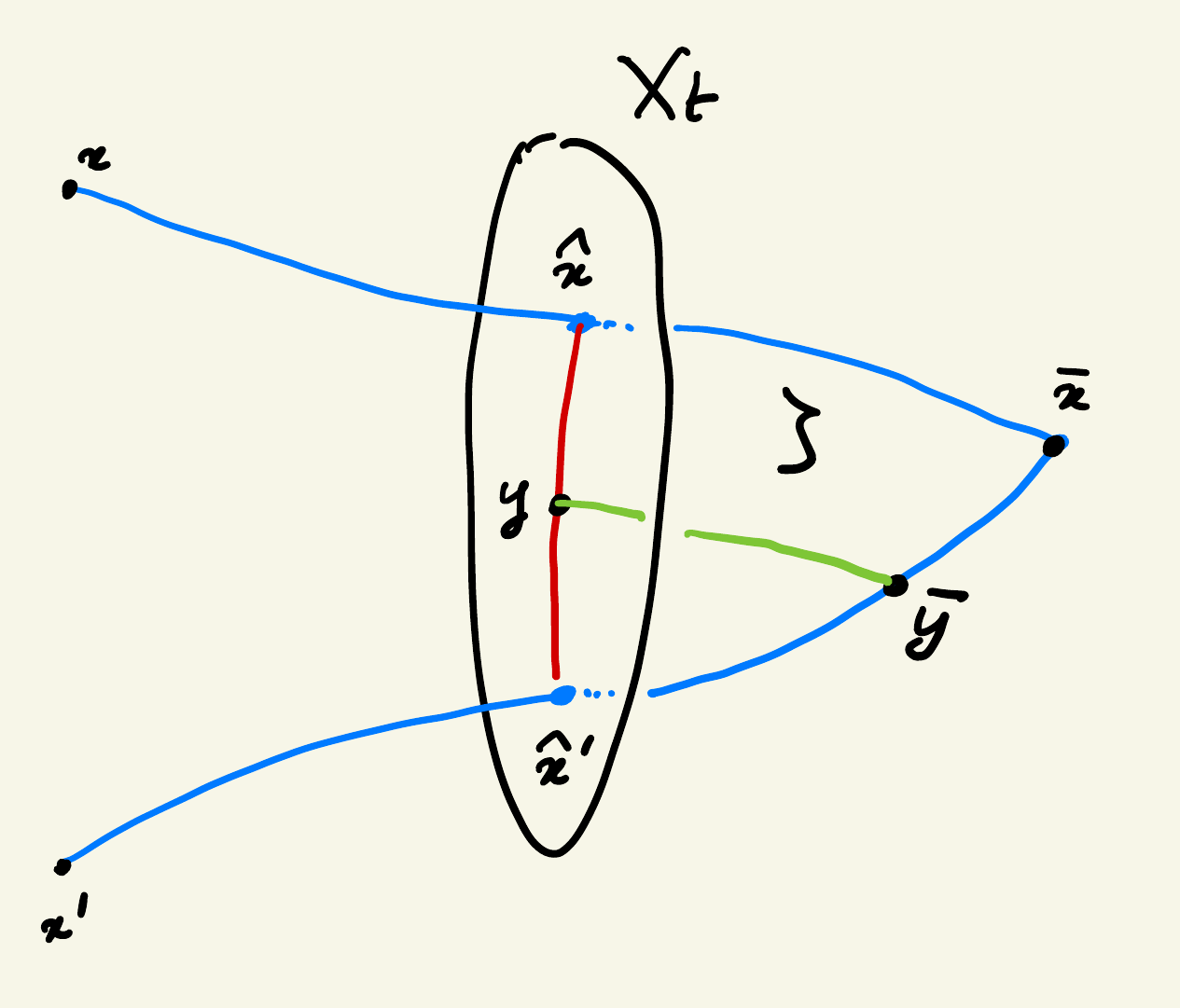}
\caption{Location of $\bar{x}$ in the detour path.}
\label{two-projections2.fig}
\end{figure}

\medskip
Another observation relating the position of the points $\bar{\bar{x}}$ and $\bar{x}$ is that, setting $y:= \tilde{x}$, if 
$\bar{y}$ denotes the $X$-projection of $y$ to the detour arc $\zeta$ (connecting $\hat{x}$ and $\hat{x}'$), then $\bar{x}$ is uniformly close to a point in the subarc $\zeta(\hat{x}, \bar{y})\subset \zeta$: 
Otherwise, as in the proof of Lemma \ref{lem:key-II.4}, the concatenation  
$[y\bar{y}]_X \star \zeta(\bar{y}, \hat{x})$ is a uniform $X$-quasigeodesic. Therefore, it has to pass uniformly close to the point $\bar{x}$. At the same time, the concatenation $[y\bar{y}]_X \star \zeta(\bar{y}, \hat{x}')$ is also a uniform $X$-quasigeodesic. Thus, $\bar{x}$ would have to be within uniformly bounded distance from both $[y\bar{y}]_X$ and $\zeta(\bar{y}, \hat{x}')$, which means that 
$\bar{x}$ is uniformly close to $\bar{y}$. We, therefore, proved (see Figure \ref{two-projections2.fig}):

\begin{lemma}\label{lem:barx-bary}
$d_X(\bar{x}, \zeta(\hat{x}, \bar{y}))\le C_{\ref{lem:barx-bary}}(K)$. 
\end{lemma}

\begin{cor}\label{cor:barx-bary}
If $d_X(\bar{x}, Y)\le r$, then $d_X(\bar{x}, \bar{\bar{x}})\le \bar{r}= \bar{r}_{\ref{cor:barx-bary}}(K,r)$. 
\end{cor}
\proof Choose $x'=y=\tilde{x}$ as above. As before, the concatenation $[x\bar{x}]_X \star [\bar{x} y]_X$ is a  $\La$-quasigeodesic in $X$ for $\La$ depending only on $K$. Take a point $p\in X_t$ within distance $r$ from $\bar{x}$. Then, according to  Lemma \ref{lem:intersection-with-flow-nbd}, 
$$
d(p, F_u\cap X_t)\le D_{\ref{lem:intersection-with-flow-nbd}}(r,K). 
$$
At the same time, by Corollary 
\ref{cor:geodesic_crossing_one_vertex_space}, the point $p$ can be chosen to lie in the $R_{\ref{cor:geodesic_crossing_one_vertex_space}}(r)$-neighborhood of  the segment $[\hat{x} y]_{X_t}$. Since 
$y=\tilde{x}$ was the projection of $\hat{x}$ to $F_u\cap X_t$ (taken in $X_t$), it follows that the distance 
from $p$ to $y$ is uniformly bounded. \qed

\medskip 
We are now ready to prove (using the notation introduced above, with the point $y:= \tilde{x}$ uniformly close to $\bar{\bar{x}}$):

\begin{lemma}\label{lem:II4}
For all $K$-special intervals $J=\llbracket u,v\rrbracket$ in $S$, and all pairs $(x,x')\in F_v^Y\times F^Y_u$:

1. The points $\hat{x}$  are uniform transition points between $x, \hat{x}'$. 

2. The pairs $(x, \hat{x}')$ are uniformly consistent. 
\end{lemma}
\proof 1. By the construction, the point $\hat{x}$ belongs to the $\La=\La(K)$-quasigeodesic 
$\phi= [x\bar{x}]_X \star [\bar{x} \hat{x}']_X$ connecting $x$ to $\hat{x}'$, implying that $\hat{x}$ is a uniform $X$-transition point between  $x, \hat{x}'$. To prove the $Y$-transition property, note that the concatenation 
$$
[x y]_Y\star [y\hat{x}']_{X_t}
$$
is a uniform $Y$-quasigeodesic (since $y$ is uniformly close to the nearest-point projection in $Y$ of $x$ to $F^Y_u$ and 
$[y\hat{x}']_{X_t}$ is a uniform quasigeodesic in $F^Y_u$). Thus, we only have to prove that $[x y]_Y$ passes uniformly close to 
$\hat{x}$. As we noted earlier, the pairs $(x,y)$ are uniformly consistent (since $x\in F^Y_v$ and 
$y$ lies in a uniform neighborhood of $F_v^Y$). The concatenation
$$
\psi=[y\bar{y}]_X\star [\bar{y} x]_X
$$ 
is a uniform quasigeodesic in $X$. The path $\hat\psi= [y \hat{x}]_{X_t} \star \hat{\psi}(\hat{x},x)$ passes through $\hat{x}$ and is 
a uniform $Y$-quasigeodesic (by the uniform consistency of the pair  $(x,y)$). The same, therefore, holds for the geodesic $[yx]_Y$. 
This implies that $\hat{x}$ is a uniform transition point between $x$ and $\hat{x}'$. 

2. The pairs $(x,\hat{x})$, $(\hat{x}, \hat{x}')$ are uniformly consistent according to Proposition \ref{prop:oneflow}, 
because the first is in $F_v^Y\times F_v^Y$ and the second is in $X_t^2$. Since $\hat{x}$ is a uniform transition point between 
$x, \hat{x}'$ according to Part 1, the pairs $(x, \hat{x}')$ are uniformly consistent (see Lemma \ref{lem:transitions}). \qed

\medskip 
This concludes  Part II.4 of the proof of Theorem \ref{thm:cut-paste}.

\subsection{Part II.5: Pairs in semispecial flow-spaces}  
Consider a $K$-semi\-special interval $J= J_1\cup J_2$ which is a union of two $K$-special intervals $J_1, J_2$ 
meeting only at a common end-point $w$. We will prove uniform consistency 
of pairs of points in $Fl^Y_K(X_J)$: This will also apply to the case of special intervals $J$ since would arbitrarily 
subdivide it into two subintervals. Observe that it suffices to prove uniform consistency of pairs of points $x_i\in F^Y_{v_i}=Fl_K^Y(X_{v_i})$, 
$v_i\in V(J_i)$, $i=1,2$. Indeed, if ${v_i}, i=1,2$ are both in, say, $J_1$, then we subdivide the interval $J_1$ further, to subintervals $J_1', J_1''$ containing $v_1, v_2$ respectively.

\begin{prop}\label{prop:II.5} 
Each pair $(x_1, x_2)$ as above is $\theta_{\ref{prop:II.5},K}$-consistent. 
\end{prop} 
\proof For $i=1,2$ consider the point 
$$
\bar{x}_i= P_{X,F^X_w}(x_i),$$
which is the nearest-point projection of $x_i$ (taken in $X$) to the flow-space $F^X_w= Fl_K^X(X_{w})$.

\begin{figure}[tbh]
\centering
\includegraphics[width=80mm]{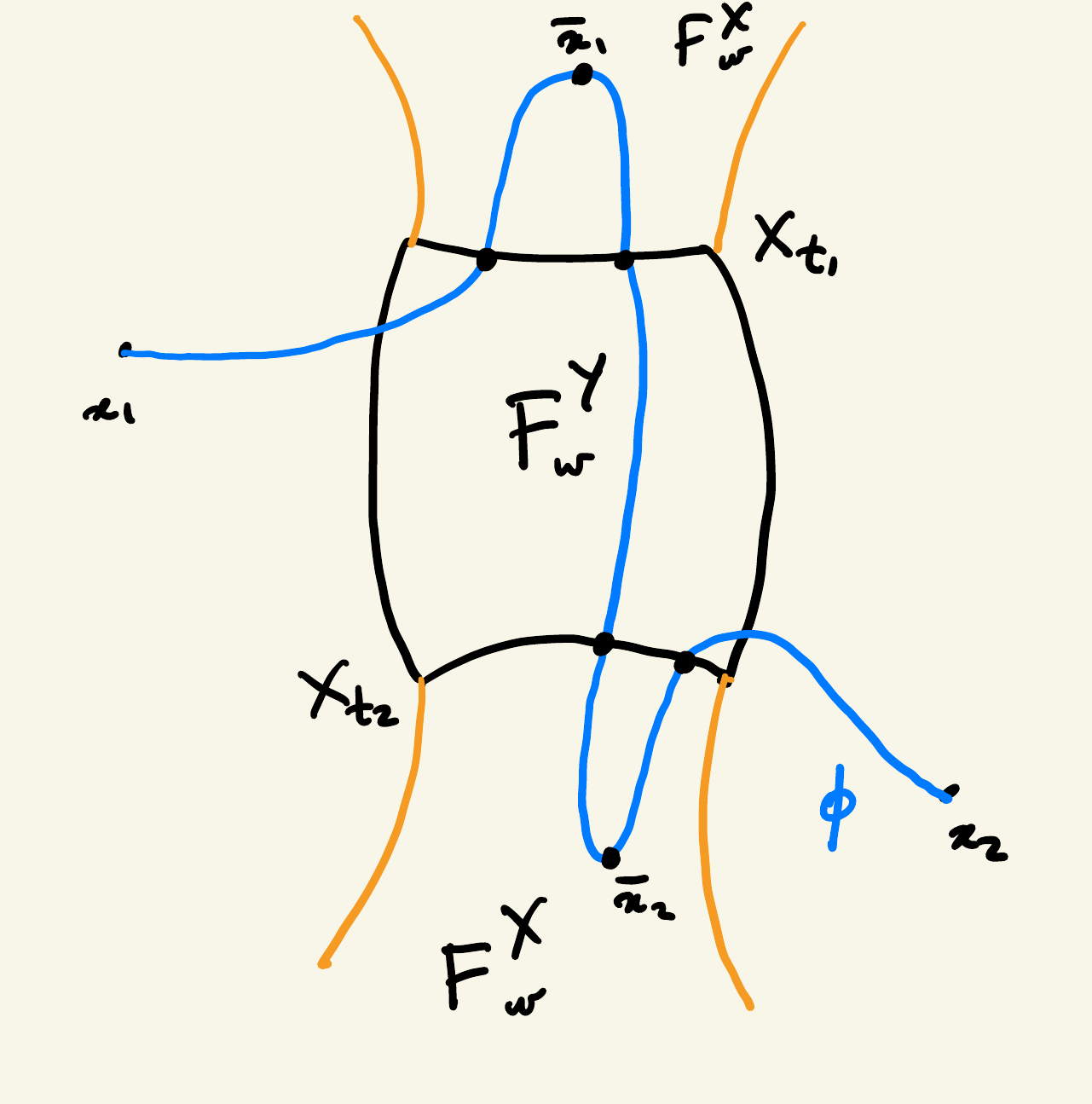}
\caption{Path $\phi$ from $x_1$ to $x_2$.}
\label{two-projections-3.fig}
\end{figure}

\begin{lemma}\label{lem:disjoint} 
The points $\bar{x}_1, \bar{x}_2$ cannot lie in the same component of $X-Y$. 
\end{lemma}
\proof Assume that both $\bar{x}_1, \bar{x}_2$ belong to $X-Y$; in particular,  
$\pi(\bar{x}_i)\ne w, i=1,2$. Then the vertex $w$ cannot separate $w_i=\pi(\bar{x}_i)$ from $v_i$, $i=1,2$: If $w$ were to separate these vertices, then the geodesic $[x_i \bar{x}_i]_X$ would cross into $X_w$ before reaching $\bar{x}_i$ and, thus, $\bar{x}_i$ would not be a nearest point to $x_i$ in $F^X_w$.  
In particular, the geodesics $\llbracket v_i , w_i\rrbracket$ lie in distinct components 
of $T-\{w\}$ and, hence, $w_1, w_2$ are separated by $w\in S$. Thus, $\bar{x}_1, \bar{x}_2$ cannot lie in the same connected component of $X-Y$.  \qed

\begin{figure}[tbh]
\centering
\includegraphics[width=80mm]{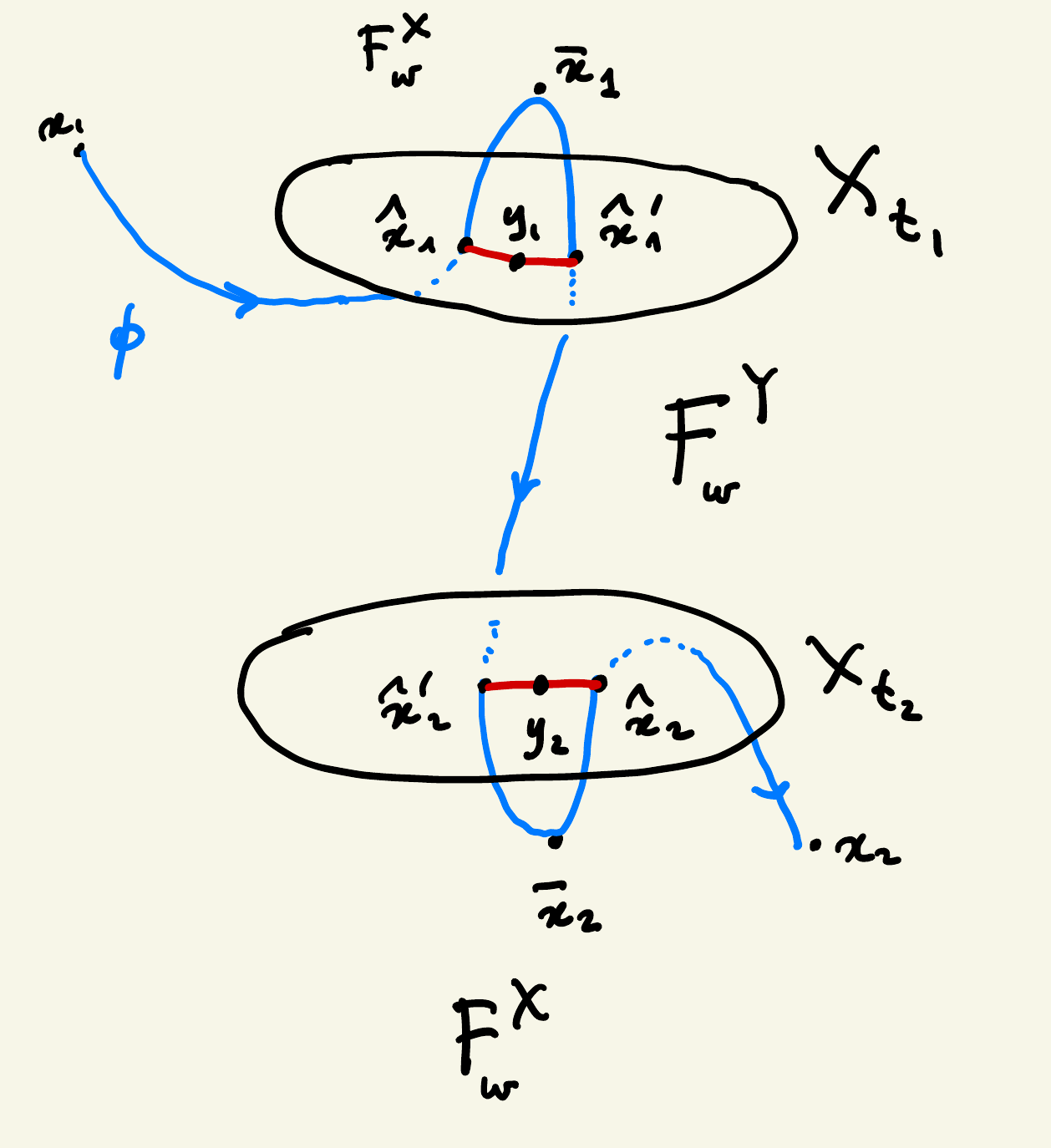}
\caption{Transition points on the path from $x_1$ to $x_2$.}
\label{two-projections-4.fig}
\end{figure}

\medskip 
As we proved in Section \ref{sec:qcamalgam}, the concatenations 
$$
\phi=[x_1 \bar{x}_1]_X \star [\bar{x}_1 \bar{x}_2]_X \star [ \bar{x}_2 x_2]_X,
$$
are $\La$-quasigeodesics,  where $\La=\La(K)$. As in the previous section, we consider  detour subpaths $\zeta_i\subset \phi$ containing $\bar{x}_i$, 
$$
\zeta_i= \phi(\hat{x}_i, \hat{x}'_i),
$$
where $\hat{x}_i, \hat{x}'_i\in X_{t_i},  i=1,2$. (See Figure \ref{two-projections-4.fig}.) By Lemma \ref{lem:disjoint}, these detour subpaths have to be disjoint, except, possibly at their end-points. The points $\hat{x}'_i$ belong to the middle portion $ [\bar{x}_1 \bar{x}_2]_X$ of 
$\phi$. We  mark points 
$$
y_i\in [\hat{x}_i \hat{x}'_i]_{X_{t_i}}, i=1,2$$ 
which are uniformly close to the projections $\bar{\bar{x}}_i=P_{Y,F^Y_w}(x_i)$, see Lemma \ref{lem:key-II.4}.

\begin{lemma}
The points $\hat{x}'_i$ are uniform transition points between $x_1$ and $x_2$, more precisely,  the sequences 
$$
x_1, \hat{x}'_1, \hat{x}'_2, x_2
$$ 
are uniformly straight in both $X$ and in $Y$. 
\end{lemma}
\proof The uniform straightness in $X$ follows from the fact that the above sequence appears in the $\La$-quasigeodesic $\phi$ in $X$ (in the correct order). To prove the uniform straightness in $Y$, note that the concatenation
$$
\psi=[x_1 y_1]_Y \star [y_1 y_2]_Y \star [y_2 x_2]_Y,
$$
is a uniform quasigeodesic in $Y$, since each $y_i$ is uniformly close to the projection of $x_i$ to the $\la_Y$-quasiconvex subset 
$F^Y_w\subset Y$  (again, see Section \ref{sec:qcamalgam}). Thus, it suffices to show that $\hat{x}'_1, \hat{x}'_2$ 
are uniform $Y$-transition points between $y_1, y_2$. As in the previous section, consider the nearest-point projections 
$\bar{y}_i\in \zeta_i$ of $y_i$. According to  Lemma \ref{lem:proj-detour-line}, there exist points 
$$
\tilde{y}_i\in [\bar{x}_i, \hat{x}'_i]_X, i=1,2, 
$$
within distance $C_{\ref{lem:barx-bary}}(K)$ from $\bar{y}_i$. The points $\bar{y}_i$  (and, hence, $\tilde{y}_i$) 
are uniformly close to the projections of $y_i$ to the quasigeodesic $\phi$ (Lemma \ref{lem:proj-detour-line}).  Hence, 
the concatenation 
$$
[y_1 \tilde{y}_1]_X\star [\tilde{y}_1 \tilde{y}_2]_X \star [\tilde{y}_2 y_2]_X 
$$
is a uniform quasigeodesic in $X$, where $[\tilde{y}_1 \tilde{y}_2]_X\subset [\bar{x}_1 \bar{x}_2]_X$. The middle geodesic in this concatenation passes through the points $\hat{x}_1', \hat{x}_2'$ (in this order), implying the uniform straightness of the sequence
$$
y_1, \hat{x}'_1, \hat{x}'_2, y_2
$$ 
in $Y$, as claimed.  \qed 

\medskip 
Now, we can finish the   Part II.5 of the proof of Theorem \ref{thm:cut-paste}. The sequence $x_1, \hat{x}'_1, \hat{x}'_2, x_2$ is uniformly straight in both $X$ and in $Y$. The pair $(\hat{x}'_1, \hat{x}'_2)\in F_w^Y\times F_w^Y$ is uniformly consistent according to Proposition \ref{prop:oneflow}, while 
both pairs $(x_1, \hat{x}'_1), (\hat{x}'_2, x_2)$ are uniformly consistent by Lemma \ref{lem:II4}(2). Therefore, by  Lemma 
\ref{lem:trivial3}, the pair $(x_1,x_2)$ is $\theta$-consistent for some function $\theta=\theta_{\ref{prop:II.5},K}$. \qed

\medskip 
This concludes Part II of the proof.

\section{Part III: Consistency in the general case} 


Suppose that $x, y \in \YY$ belong to vertex-spaces $X_u, X_v$, respectively, $u, v\in V(S)$. 
Fix $K=K_0$. Using Lemma \ref{lem:subdivision} and Theorem \ref{thm:hyperbolicity over segments}, 
 we obtain a {\em horizontal subdivision} of 
the interval $J= \llbracket u, v\rrbracket$, 
into subintervals $J_i=\llbracket u_i, u_{i+1}\rrbracket, i=1,...,n$, such that the pairs of distinct 
flow-spaces $F^X_i:= Fl^X_K(X_{u_i})$ have disjoint projections to $T$ and, hence, are 
$C=C_{\ref{thm:hyperbolicity over segments}}$-cobounded, unless, possibly 
$i=n$; $u=u_1, v=u_{n+1}$. (The same, of course also holds for the $Y$-flows.) 
As in Lemma \ref{lem:subdivision}, we also define vertices $u_i'', u_{i+1}'\in J_i$ such that

1.   $u'_i, u_i$ span an edge $e_i$  in $T$ (except, possibly, for $i=n+1$, in which case we can have $d_T(u'_{n+1}, u_{n+1})\le 1$). 

2. Each subinterval $\llbracket u'_i, u_{i}\rrbracket, \llbracket u_i, u''_{i}\rrbracket,  \llbracket u''_i, u'_{i+1}\rrbracket$ is $K$-special. 
In particular, the subinterval $J'_i= \llbracket u_{i}, u'_{i+1}\rrbracket$ is semispecial.

For each $i$ we define pairs of points $(x_i'', x'_{i+1})\in \FF^X_i\times \FF^X_{i+1}$ realizing the minimal distance between the subsets $\FF^X_i, \FF^X_{i+1}$ of $X$ unless $i=n$ and $u'_{n+1}=u_{n+1}$, 
in which case we take $x''_n$ to be the projection of $y$ to $\FF^X_n$.

\begin{figure}[tbh]
\centering
\includegraphics[width=80mm]{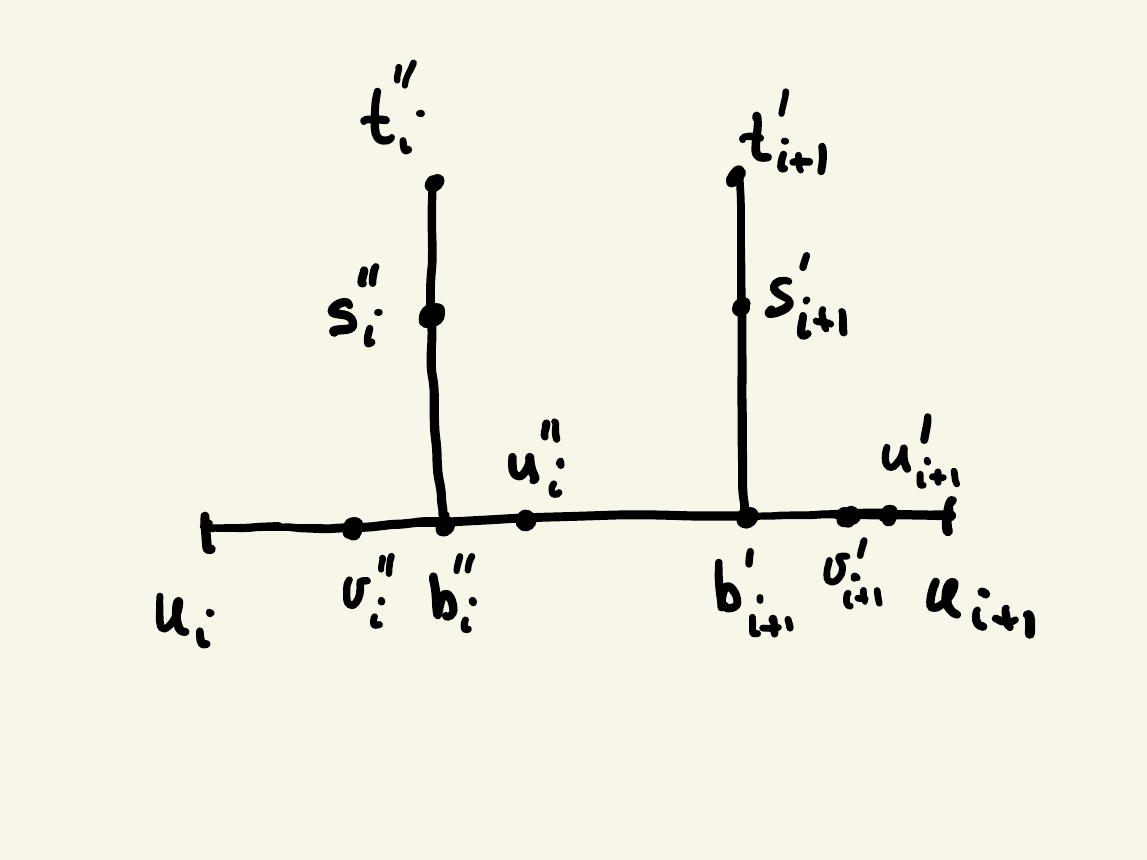}
\caption{Projections to $J$.}
\label{subdivision-projections.fig}
\end{figure}

\medskip 
For each $i$ define the vertices $t'_i:= \pi(x'_i)$, $t''_i:= \pi(x''_i)$ and let $s_i', s''_i$ denote their respective projections to the subtree $S\subset T$. 

\begin{lemma}\label{lem:semis}
Let $b'_i$ (resp. $b''_i$, $v'_i$, $v''_i$) denote the projection to $J$ of the vertices $t'_i$ (resp. $t''_i$, $\pi(y'_i)$, $\pi(y''_i)$). Then:

1. The projections $b'_i, v'_i$  lie 
in the subinterval $\rrbracket u''_{i-1}, u'_i\rrbracket$. 

2. The projections $b''_i, v''_i$   lie in the subinterval  
 $\llbracket u_{i}, u''_i\rrbracket$.  
\end{lemma} 
\proof We will prove the claim for the vertex $v'_i$, since the rest is proven by the same argument. 

First of all, 
$$
\pi(F^Y_i)\subset \pi(F^X_i)
$$
and, by the definition of the horizontal subdivision 
of the interval $J$,  $\pi(F^X_i)$ intersects $J$ in a subinterval of 
$\rrbracket u''_{i-1}, u''_i\rrbracket$. Thus, $v'_i$ belongs to $\rrbracket u''_{i-1}, u''_i\rrbracket$. 
Suppose, for the sake of a contradiction, that $v'_i$  is  
in the interval $\llbracket u_{i}, u''_i\rrbracket$. Then each geodesic connecting $y''_{i-1}$ to $y'_i$ goes through the edge-space $X_{e}$ of the edge 
$$
e= [u'_i, u_i],
$$
before reaching $x'_i$. But then, this geodesic also passes through the subset
$$
X_{e u'_i}\subset X_{u'_i}. 
$$
Since $K\ge 1$, the entire subset $X_{e u'_i}$ is contained in the flow-space $Fl_K^Y(X_{u_i})$. Hence, $y'_i$ cannot possibly be 
the $Y$-projection of $y''_i$ to $F_i^Y$. A contradiction. \qed 

\begin{cor}
The pairs $(y''_i, y'_{i+1})$ are uniformly consistent. 
\end{cor}
\proof Let $v_i'', v'_{i+1}$ denote the projections of $\pi(y''_i)$, $\pi(y'_{i+1})$ to $J$. By the lemma, both  
$v_i'', v'_{i+1}$  lie in the interval $J'_i= \llbracket u_i, u'_{i+1}\rrbracket$, which is semispecial. Moreover, 
$y''_i\in Fl_K^Y(X_{u_i})$ and $y'_i\in  Fl_K^Y(X_{u'_{i+1}})$. Thus, 
$$
(y''_i, y'_{i+1})\in Fl_K^Y(X_{J'_i})
$$
and, hence, the statement is a special case of Proposition \ref{prop:II.5}. \qed  

  \begin{figure}[tbh]
\centering
\includegraphics[width=100mm]{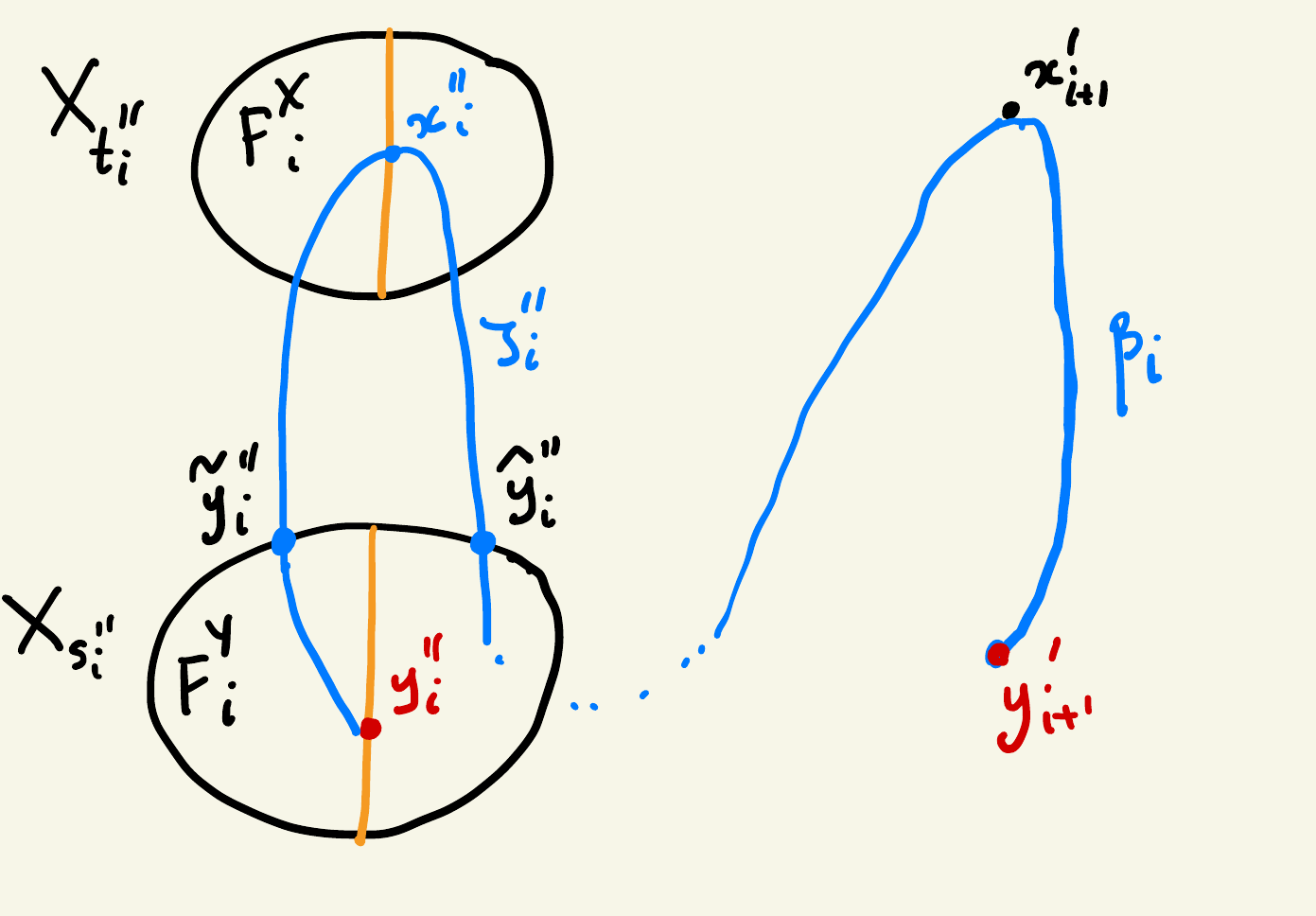}
\caption{Transition points}
\label{transitions-2.fig}
\end{figure}

\medskip 
Our next task is to relate the points $x'_i, y'_i$ and also relate the points $x''_i, y''_i$. Since $y_i'$ is in $Y$, the segment 
$[y'_i x'_i]_X$ has to cross the vertex-space $X_{s_i'}$; the same holds for the segment $[y''_i x''_i]_X$ and $X_{s_i''}$. 
We define the points  $\tilde{y}'_i\in [y'_i x'_i]_X$, $\tilde{y}''_i\in [y''_i x''_i]_X$ by the condition that the subsegments 
$$
[\tilde{y}'_i x'_i]_X\subset [y'_i x'_i]_X
$$
and 
$$
[\tilde{y}''_i x''_i]_X\subset [y''_i x''_i]_X
$$
are the smallest subsegments terminating in $x'_i, x''_i$,  
with the property that $\tilde{y}'_i\in X_{s_i'}, \tilde{y}''_i\in X_{s_i''}$. (See Figure \ref{transitions-2.fig}.) 
Thus,   
$$
Y\cap [\tilde{y}'_i x'_i]_X \setminus \{\tilde{y}'_i\} = \emptyset,
$$
and, similarly, for $[\tilde{y}''_i x''_i]_X$. 
 
Similarly, define points $\hat{y}_i'', \hat{y}_{i+1}'\in [x''_i x'_{i+1}]_{X}$ such that 
$$
[ \hat{y}'_{i} x'_{i}]_X\subset [x''_{i-1} x'_{i}]_{X}, [\hat{y}_i'' x''_i]_X\subset [x'_{i+1}x''_i]_{X}, 
$$
are the smallest subsegments terminating in $x'_i, x''_i$ with $\hat{y}'_i\in X_{s_i'}, \hat{y}''_i\in X_{s_i''}$.  Observe that the concatenations
$$
\beta_i:= [y_i'' x''_i]_X\star [x''_{i} x'_{i+1}]_{X} \star [x'_{i+1} y'_{i+1}]_X
$$
are uniform quasigeodesics in $X$ since $x''_i\in F^X_i, x'_{i+1}\in F^X_{i+1}$ realize 
the minimal distance between these two uniformly quasiconvex subsets of $X$ and $y_i''\in  F^Y_i\subset 
F^X_i, y'_{i+1}\in F^Y_{i+1}\subset F^X_{i+1}$. Both 
$$
\zeta_i'= [\hat{y}'_i x'_i]_X \cup [x'_i \tilde{y}'_i]_X, \zeta_i''= [\hat{y}''_i x''_i]_X \cup [x''_i \tilde{y}''_i]_X
$$
are detour subpaths in $\beta_i$, containing the points $x'_i, x''_i$ respectively.

\begin{lemma}
1. The pairs $(\hat{y}''_i, \hat{y}'_{i+1})\in Y^2$ are uniformly consistent. 

2. The distances $d(y'_i, \tilde{y}_i')$, $d(y''_i, \tilde{y}_i'')$ are uniformly bounded. 
\end{lemma}
\proof 
Part 1. By the previous corollary, the pair $(y''_i, y'_{i+1})$ is uniformly consistent. In particular, the path
$$
\hat{\beta}_i= \widehat{[y_i'' \tilde{y}''_i]_X} \star [\tilde{y}''_i \hat{y}''_i]_{X_{s_i''}} \star \widehat{[\hat{y}''_i \hat{y}'_{i+1}]_X} 
\star [\hat{y}'_{i+1} \tilde{y}'_{i+1}]_{X_{s'_{i+1}}} \star \widehat{[\tilde{y}'_{i+1} y'_{i+1}]_X}
$$
is uniformly quasigeodesic. Thus, the geodesic $[y_i'' y'_{i+1}]_Y$ passes within uniform distance $D$ from the points $\tilde{y}''_i, \hat{y}''_i, \hat{y}'_{i+1}, \tilde{y}'_{i+1}$ (in this order). Taking into account uniform consistency of the pairs $(y''_i, y'_{i+1})$ and Remark \ref{rem:trivial1}, we conclude uniform consistency of the pairs $(\hat{y}''_i, \hat{y}'_{i+1})$.

\medskip 
Part 2. We will estimate the second distance, $d(y''_i, \tilde{y}_i'')$, since the proof for the first one is similar. Since $F^X_i$ is $\la_X$-quasiconvex in $X$, 
and $(x'_i, y'_i)\in F_i^X\times F_i^X$, the point $\tilde{y}_i'$ lies within distance $\la_X$ from a point $q$ in  $F_i^X$. The point $q$ might not be in $F^Y_i= F^X_i\cap Y$, but it is within distance $K\la_X$ from $F^X_i\cap X_{s'_i}\subset F^Y_\la$. In Part 1 we observed that 
the geodesic $[y_i'' y'_{i+1}]_Y$ passes (at some point $p$) within uniform distance $D$ from the point $\tilde{y}''_i$. Thus, we found a point $p\in [y_i'' y'_{i+1}]_Y$ within distance $D+(K+1)\la$ from $F^Y_i$. It follows that 
$$
d_Y(p, y''_i)\le D+(K+1)\la_X
$$
and, hence,
$$
d_Y(\tilde{y}''_i, y''_i)\le 2D+(K+1)\la_X. \qed 
$$

\begin{rem}
While it is not needed for our purposes, one can prove similarly to Lemma \ref{lem:key-II.4}, that the points ${y}'_i, {y}''_i$ are uniformly close to  the projections (taken in the vertex-spaces $X_{s_i'}, X_{s_i''}$ respectively) of the points $\hat{y}_i', \hat{y}_i''$ 
to the $4\delta_0$-quasiconvex subsets 
$$
F^Y_i\cap X_{s_i'}, F^Y_i\cap X_{s_i''}
$$
respectively. This provides a description (up to a uniform error) of the points ${y}'_i, {y}''_i$ in terms of $x''_{i-1}, {x}'_i, {x}''_i$. 
\end{rem}

\begin{lemma}
The pairs $(\hat{y}'_i, \hat{y}_i'')$ are uniformly consistent. 
\end{lemma}
\proof Each interval 
$$
I_i= \llbracket s'_i, s''_i\rrbracket\subset T
$$
is contained $\llbracket t'_i, t''_i\rrbracket$ and intersects $J$ along the interval $\llbracket b'_i, b''_i\rrbracket$, which, in turn, 
contains the vertex $u_i$, see Lemma \ref{lem:semis}. Since  
$F_i^Y$ has nonempty intersection with both $X_{s'_i}, X_{s''_i}$, it follows that the interval $I_i$ is semispecial (it is the union of special subintervals $\llbracket s'_i, u_i\rrbracket$,  $\llbracket s''_i, u_i\rrbracket$). Thus, 
$$
y'_i, y''_i\in X_{I_i}\subset Fl^Y_K(X_{I_i})
$$
and, hence, the pair $(\hat{y}'_i, \hat{y}_i'')$ is uniformly consistent by Proposition \ref{prop:II.5}. \qed

\medskip
We now can finish the proof of Theorem \ref{thm:cut-paste}. For each $i$, the points $\hat{y}''_i, \hat{y}'_{i+1}$ are both $X$ and $Y$-transition points between, respectively, $x''_i, x'_{i+1}$ and $y''_i, y'_{i+1}$. The sequence
$$
x, \hat{y}''_1, \hat{y}'_2, \hat{y}''_2,..., \hat{y}'_n, \hat{y}''_n, \hat{y}'_{n+1}, y
$$
is uniformly straight in both $X$ and in $Y$ and the consecutive pairs of points in this sequence are uniformly consistent. Now, the uniform consistency of the pairs $(x,y)$ follows from Lemma \ref{lem:trivial3}. \qed

\section{The existence of CT-maps for subtrees of spaces}

We finally are ready to prove the main result of this chapter: 

\begin{theorem}\label{thm:ECT}
Let $\X=(\pi: X\to T)$ be a tree of hyperbolic spaces with hyperbolic total space $X$ and let $\Y=(\pi: Y\to S)$  be a subtree of spaces in $\X$, 
where $S\subset T$ is a subtree and $Y= \pi^{-1}(S)$. Then the inclusion map $Y\to X$ admits a CT-extension. 
\end{theorem}
\proof  We will derive  this result  from Theorem \ref{thm:cut-paste}.  Fix $p\in Y$ and suppose that $(y_n), (y'_n)$ are sequences in $Y$ such that
$$
\lim_{n\to\infty} (y_n. y'_n)^Y_{p}=\infty,
$$ 
where the superscript $Y$ refers to the Gromov-product of the intrinsic path-metric of $Y$. Equivalently, for the geodesic 
$\beta_n=[y_n y_n']_Y\subset Y$ we have
$$
\lim_{n\to\infty} d(p, \beta_n)=\infty. 
$$
For the sake of contradiction, assume that the corresponding $X$-geodesics 
$\alpha_n=[y_n y_n']_X\subset X$ all pass through a ball 
$B_X(y, R)\subset X$ for some fixed $R$ (which is equivalent to saying that the sequence $(y_n, y'_n)^X_{p}$ does not diverge to infinity), i.e. there exist 
points 
$$
q_n\in \al_n\cap B_X(p,R). 
$$ 
By Theorem \ref{thm:cut-paste}, each $\hat\alpha_n$ is a $\La$-quasigeodesic in $Y$ connecting $y_n, y_n'$, for some $\La$ independent of $n$.  
By the Morse Lemma (Lemma \ref{Morse}), 
$$
\Hd_Y(\beta_n, \hat\alpha_n)\le D=D(\delta_Y,\La).
$$ 
There are two cases which may occur:

a. $q_n\in Y$. Then (by the definition of $\hat{\al}_n$) $q_n\in \hat{\al}_n$, hence, 
$$
d_X(p, \beta_n)\le D +R
$$
for all $n$, which is a contradiction.

b. $q_n\notin Y$. Let $\zeta_n=[z_n z'_n]_X\subset \al_n$ be a detour subpath in $\al_n$ containing $q_n$; the end-points $z_n, z'_n$ of $\zeta_n$ belong to a vertex-space $X_{v_n}\subset Y$. Since $p\in Y$, $q_n\notin Y$, the vertex-space $X_{v_n}$ separates $p$ from $q_n$ and, hence, the geodesic $[p q_n]_X$ has to pass through $X_{v_n}$ at some point $p_n\in X_{v_n}$. We, obviously, have
$$
d_X(p, p_n)\le R.
$$
Thus, by Corollary \ref{cor:geodesic_crossing_one_vertex_space}, points $p_n$ are all uniformly close (within distance $r$ depending only on the parameters of $\X$) to the geodesic $[z_n z'_n]_{X_{v_n}}\subset \hat{\al}_n$. Therefore,
$$
d_X(p, \beta_n)\le D + R+ r,
$$
 which is again a contradiction. \qed

\section{Fibers of CT-maps}\label{sec:CTfibers}

Let $\X=(\pi: X\to T)$ be a  tree of hyperbolic spaces with hyperbolic total space $X$. According to Theorem \ref{thm:ECT}, 
for every subtree $S\subset T$ and $Y=X_S$, the inclusion map 
$$
f_{Y,X}: Y=X_S\to X$$
admits a Cannon--Thurston extension 
$$
\D_{Y,X}:=\geo f_{Y,X}: \geo Y\to \geo X.
$$
 
For the rest of the chapter, we will be working under the extra assumptions that $X$ is a proper metric space, i.e. that an abstract tree of spaces $\X$ admits a proper total space. Note that if some edge-spaces are non-discrete, the  total space $X$ defined in the proof of Theorem \ref{thm:existence-of-trees} need not be proper even if $V(T)=\{v, w\}$ 
 and  $X_v, X_w, X_{[v,w]}$ are proper. However, this $X$ will be proper (a locally finite metric graph) under the 
following conditions:

1. Each vertex-space $X_v$ is a locally-finite graph (with the standard graph-metric).  

2. All edge-spaces $X_e$ are discrete and $X_{ev}=f_{ev}(X_e)\subset V(X_v)$ for 
each edge $e$ and incident vertex $v$.

3. Each finite subset of each vertex-space $X_v$ meets only finitely many subsets of the form $X_{ev}$, where $e$'s are edges incident to $v$.

\begin{prop}\label{prop:CT-fibers1} 
Suppose that $\X$ has a proper total space $X$. Then there exists a number $r$, depending only on $\X$ and $X$, for which  the following holds. For each $Y=X_S\subset X$, if $\xi_\pm\in \geo Y$ are two distinct points such that $\D_{Y,X}(\xi_+)=\D_{Y,X}(\xi_-)$,  
then there exists a vertex $v\in S$ and two ideal boundary points $\xi_\pm'\in \geo X_v$ such that:

1. $\D_{X_v,Y}(\xi_\pm')= \xi_\pm$. 

2. The $Y$-geodesic $\eta: \RR\to Y$, connecting $\xi_+, \xi_-$, is contained in the $r$-neighborhood of 
a vertex space $X_v, v\in S$, and, moreover, is $r$-Hausdorff-close to a vertical geodesic $\al$ in $X_v$. 
\end{prop}
\proof 
Set $x_n:=\eta(n)$, $n\in \Z$. Define subarcs of the geodesic $\eta$ by 
$$
\eta_n:= [x_{-n} x_{n}]_Y, n\in \NN,  
$$
and geodesics $\beta_n= [x_{-n} x_{n}]_X$. 

Since $\D_{Y,X}(\xi_+)= \D_{Y,X}(\xi_-)$, we have that 
$$
\lim_{n\map \infty} d_X(x_0, \beta_n)  = \infty. 
$$
Thus, $x_0$ cannot be uniformly close to $\beta_n\cap Y$.  At the same time, since $x_0$ belongs to each $\eta_n$ 
(which, by Theorem \ref{thm:cut-paste},  is uniformly close to the $\La$-quasigeodesic $\hat\beta_n$ in $Y$), 
the point $x_0$ lies in the $C$-neighborhood of  
one of the replacement arcs $\hat\zeta_n=[p_n q_n]_{X_{v_n}}\subset \hat\eta_n$, where $\La$ and $C$ depend 
only on the parameters of $\X$. In particular, $d_X(x_0, X_{v_n})\le C$ for all $n$. Since $X$ is assumed to be proper,  
there are only finitely many vertex-spaces in $X$ which can intersect $B(x_0, C)$. 
It follows that there is an infinite subset $M\subset {\mathbb N}$ and a vertex $v\in S$ 
 such that for all $m\in M$, $v_m=v$. At the same time, since the distances $d(x_0, \beta_n)$ diverge to infinity, 
 we also have
 $$
 \lim_{n\to\infty} d(x_0, p_n)=\infty, \quad  \lim_{n\to\infty} d(x_0, q_n)=\infty. 
 $$
Hence, after passing to a further subsequence, the sequences $(p_m)_{m\in M}$, $(q_{m})_{m\in M}$, 
diverge in $X_v$ to two ideal boundary points $\xi'_\pm$.   

By continuity at infinity of the CT-map $\D_{X_v, Y}: \geo X_v\to \geo Y$, it follows that 
$$
\D_{X_v,Y}(\xi_\pm')= \xi_\pm. 
$$   
After passing to a subsequence once more, we can assume that the sequence of geodesics $\hat\zeta_n\subset X_v$ converges to a complete geodesic $\zeta\subset X_v$ asymptotic to the points $\xi'_\pm$. Since geodesics $\hat\zeta_n$ are 
$\La$-quasigeodesics in $Y$, so is their limit $\al$. Since the $Y$-geodesic $\eta$ is also asymptotic to $\xi_\pm$, 
it follows that $\eta, \al$ are within Hausdorff distance $r:=2D_{\ref{Morse}}(\delta_Y,\La)$,  
which depends only on the parameters of $\X$. It also follows that  $\eta$ is contained in $N_r(X_v)$. \qed

\medskip
In the following addendum to this proposition we describe more precisely (up to a uniform error) 
the nature of geodesic segments $\beta_{-m,n}$ connecting the points $x_{-m}, x_n$, $n, m\in \NN$, 
in the setting of the  theorem. 
Since the points $x_{-m}, x_n$ belong to the $r$-neighborhood of $X_v$, we will be considering instead of $\beta_{-m,n}$'s 
the uniform  quasigeodesics $c_{-m,n}$ (from the slim combing of $X$ described in Section \ref{sec:inductive}) 
connecting points $y_{-m}, y_n$, $m, n\in \NN$, where  $y_i:= \al(i), i\in \Z$.  
We will see that the path $c_{-m,n}$ first diverges away from $X_v$ (in the metric of $X$) at a  linear speed and then converges back to $X_v$ at a  linear speed. We refer the reader to Section \ref{sec:inductive} (step I.1) 
for the description of uniform quasigeodesics connecting 
points in narrow carpets used in the proof of the next proposition.

\begin{rem}\label{rem:carpet paths} 
In what follows, we will frequently use the following notation. Let $\A=\A(\al)$ be a $(K,C)$-narrow carpet with the narrow end $\beta\subset X_w$, bottom and top sections $\gamma_-, \gamma_+$. Then we define the  path $c=c_{\A}$ as the concatenation
$$
\gamma_- \star \beta\star \gamma_+
$$ 
connecting the bottom and and the top points of the segment $\al\subset X_u$. Such paths will be uniform (with qi constants depending on $K$ and $C$) 
quasigeodesics in $\A$ and, hence, in $X$, as 
long as for each vertex $v\in X_v$, the length of $A_v$ is $\ge M_{\bar{K}}$, see Section \ref{sec:Hyperbolicity of carpets} or Section \ref{sec:inductive}, Step I.1. 
\end{rem}

\begin{prop}\label{prop:CT-fibers2} 
1. For all  $n>0, m>0$, the segment $\al_{-m,n}= [y_{-m} y_n]_{X_v}$ bounds a certain $(K,C)$-narrow carpet 
$\A_{-m,n}=\A_K(\al_{-m,n})$ in $X$ over the interval $J_{-m,n}= \llbracket v, w_{-m,n}\rrbracket$.  Here $K$ and $C$ depend only on the parameters of $\X$. 

2. For   $n>0, m>0$, the uniform quasigeodesic $c_{-m,n}$ in $A_{m,n}\subset X$ connecting points 
$y_{-m}, y_n$ is the concatenation 
$$
\ga_{-m} \star \beta_{-m,n}\star \ga_n, 
$$
where $\beta_{-m,n}$ is a vertical geodesic in $X_{w_{-m,n}}$ of length $\le C$ and $\ga_{-m}, \ga_n$ are 
$K$-qi sections over the interval $J_{-m,n}= \llbracket v, w_{-m,n}\rrbracket$. 
\end{prop}
\proof Set $K:= K_0$ (defined in Notation \ref{not:K0}). Consider maximal $K$-qi sections $\Si_{-m}, \Si_n$ in $X$ through the points $y_{-m}, y_n$. If, for some $n_0, m_0$, 
these sections have vertical separation $\ge M_{K}$ everywhere, then they are uniformly 
cobounded in $X$, see \ref{sec:inductive}. 
In this situation, for all 
$n\ge n_0, m\ge m_0$, the path $c_{-m,n}$ (and, hence, 
the geodesic $\beta_{-m,n}=[x_{-m} x_n]_X$) 
has to come uniformly close to a pair of points $x^-\in \Si_{-m_0}, x^+\in \Si_{n_0}$, contradicting the assumption that 
$$
\lim_{n\to\infty, m\to\infty} d(y_0, \beta_{-m,n})=\infty. 
$$

Thus, for each pair $(m, n)$ there is a vertex $w_{m,n}$ and a pair of $K$-qi leaves $\ga_n\subset \Si_n, \ga_{-m}\subset \Si_{-m}$ over an interval  
$J_{m,n}= \llbracket v, w_{m,n}\rrbracket$, with vertical separation $\ge M_K$, such that 
$$
d_{X_{w_{m,n}}}(\ga_{-m}(w_{m,n}), \ga_{n}(w_{m,n}))\le M_{{K}}. 
$$
The uniform quasigeodesic $c_{-m,n}$ connecting points $y_{-m}, y_n$ then is defined as the concatenation
$$
\ga_{-m} \star \beta_{-m,n}\star \ga_n, \quad  \beta_{-m,n}= [\ga_{-m}(w_{-m,n}) \ga_{n}(w_{-m,n})]_{X_{w_{-m,n}}}.
$$
(This argument is similar to the proof of Proposition \ref{prop:uncarpeted geodesics}.)  \qed

\medskip 
We can now give a complete description of pairs of distinct points in the fibers of the CT-map $\D_{Y,X}$. 
Recall that the metric space $X$ is assumed to be proper.

\begin{thm}\label{thm:CT-fibers} 
There are constants $K, C$ depend only on the parameters of $\X$ and a function $D=D(k)$ such that the following hold:

1. Suppose that $\xi^\pm$ are distinct points in $\geo Y$ such that $\D_{Y,X}(\xi^-)=\D_{Y,X}(\xi^+)$. Then there exists a vertex-space 
$X_u\subset Y$ and a complete geodesic $\al: \RR\to X_u$, 
which is a uniform quasigeodesic in $Y$ asymptotic to $\xi^\pm$, such that the intervals $\al_{-m,n}= [\al(-m) \al(n)]_{X_u}\subset \al$, bound $(K,C)$-narrow carpets $\A(\al_{-m,n})$ in $X$ for all $m>0, n>0$. 

2. Conversely, if $X_u$ is a vertex-space of $\Y$, $\al\subset X_u$ is a complete geodesic asymptotic to distinct points $\xi^\pm\in \geo Y$, such that each subinterval 
$\al_{-m,n}$ as above bounds a $(K,C)$-narrow carpet $\A(\al_{-m,n})$ in $X$, then $\D_{Y,X}(\xi^-)=\D_{Y,X}(\xi^+)$.

3. Suppose that $X_u$ is a vertex-space of $\Y$, $\al\subset X_u$ is a complete geodesic asymptotic to distinct points $\xi^\pm\in \geo Y$. Then 
$\D_{Y,X}(\xi^-)\ne \D_{Y,X}(\xi^+)$ if and only if for some (equivalently, every) $k\ge 1$, there exist points $x, y\in \al$ and maximal $k$-qi sections $\Si_x, \Si_y$ over, possibly different, subtrees $T_x, T_y$ in $T$  through the points $x, y$ such that the vertical separation between $\Si_x, \Si_y$ over every vertex of $T_x\cap T_y$ is $\ge D$. 
\end{thm}
\proof We take $K=K_0$ and $C=M_{\bar{K}}$. 

1. The first part of the theorem is the content of Propositions \ref{prop:CT-fibers1} and \ref{prop:CT-fibers2}. 

2.  Since the lengths of the intervals $\al_{-m,n}= [y_{-m} y_n]_{X_u}$ 
diverge to $\infty$ as $m\to\infty, n\to\infty$, for all sufficiently large $m, n$, 
without loss of generality, we may assume that the vertical separation between the top and the bottom of each carpet 
$\A(\al_{-m,n})$ is $\ge M_K$. (Otherwise, since $C=M_{\bar{K}}\ge M_K$, 
we take a smaller $(K,C)$-subcarpet $\A'(\al_{-m,n})\subset \A(\al_{-m,n})$ containing no vertical intervals of length $< M_K$.) 
Now, just as in proof of Proposition \ref{prop:uncarpeted geodesics}, each 
carpet $\A(\al_{-m,n})$ defines a uniform $X$-quasigeodesic $c_{-m,n}$ connecting $y_{-m}$ to $y_n$ and 
$$
d_X(y_0, c_{-m,n})\ge \length(\pi(A(\al_{-m,n}))). 
$$ 
Thus,
$$
\lim_{m\to\infty, n\to\infty} d_X(x_0, [y_{-m} y_n]_X)=\infty
$$
and, therefore, $\D_{Y,X}(\xi^-)=\D_{Y,X}(\xi^+)$.

  \begin{figure}[tbh]
\centering
\includegraphics[width=70mm]{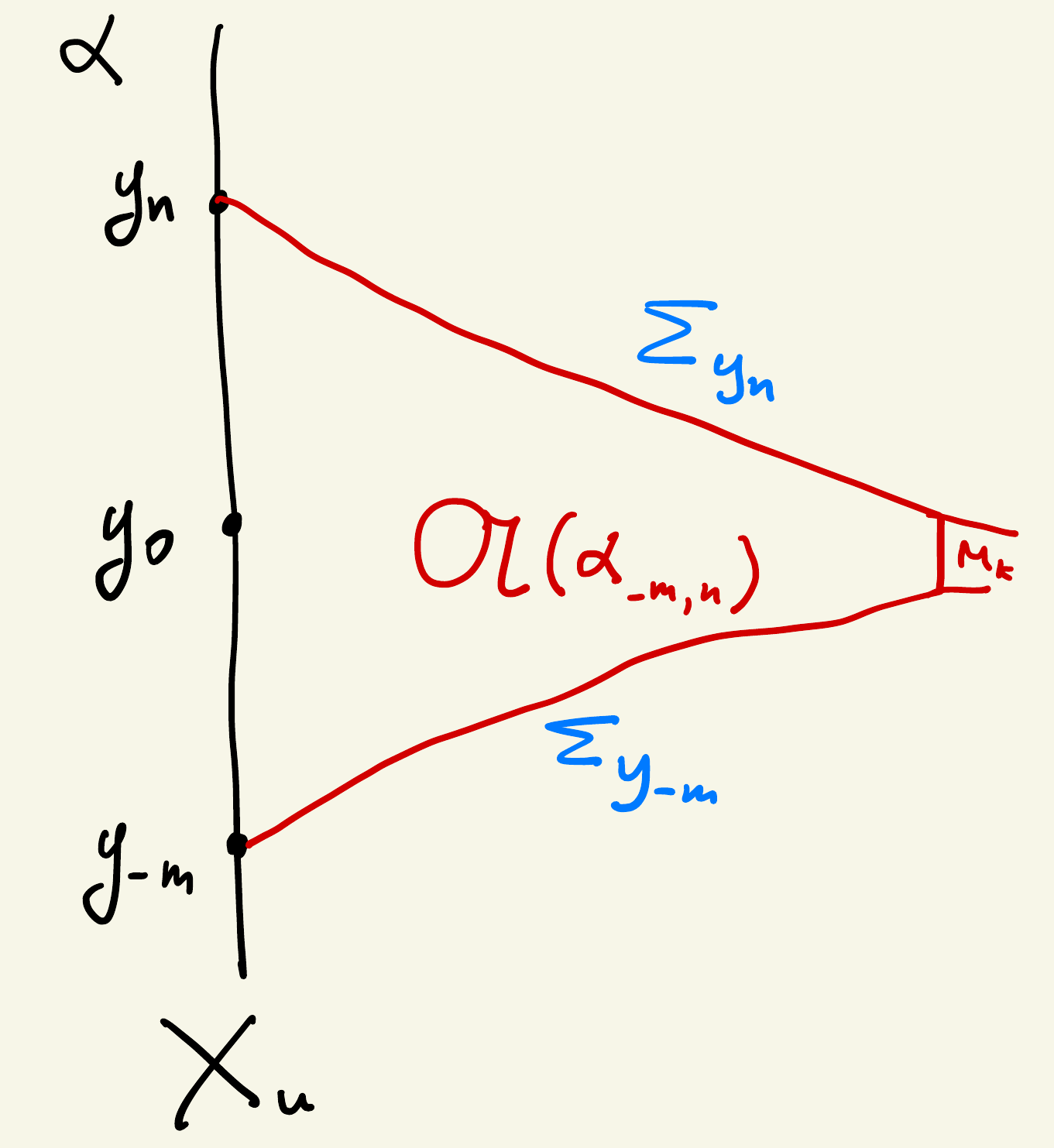}
\caption{}
\label{figure15}
\end{figure}

3. This part of the proof is similar to that of Proposition \ref{prop:CT-fibers1}. We choose $D=M_k$. 

(a) Suppose that  sections $\Si_x, \Si_y$ exist. After reparameterizing $\al$, we can assume that 
$x=\al(-m_0), y=\al(n_0), m_0>0, n_0>0$. There is a ladder $\L_{X,k}(\al)\subset \X$ containing sections $\Si_x, \Si_y$.   
Consider the combing paths $c_{-m,n}=c(\al(-m), \al(n))$ in the ladder $\L_{X,k}(\al)$ for $n\ge n_0, m\ge m_0$. These paths have to go through both sections   $\Si_x, \Si_y$ and, hence, pass uniformly close to  a pair of points $x', y'$ (independent of $m, n$) in these sections realizing the minimal fiberwise distance 
between $\Si_x, \Si_y$ (the sections are cobounded in $L_{X,k}(\al)$, see Section \ref{sec:inductive}). Thus, for $y_0=\al(0)$, the minimal distances  
$d_X(y_0, c_{-m,n})$ are uniformly bounded (from above), hence,  
$$
\lim\sup_{m, n\to\infty} (\al(-m). \al(n))_{y_0}<\infty,
$$
and, therefore, $\D_{Y,X}(\xi^-)\ne \D_{Y,X}(\xi^+)$. 

(b) Suppose that the points $x, y$ do not exist. Take arbitrary maximal $k$-qi sections  $\Si_{y_{-m}}, \Si_{y_n}$ in $X$ through the points $y_{-m}=\al(-m), y_n=\al(n)$. Then for all $m\ge 0, n\ge 0$ the minimal vertical separation between  
 $\Si_{y_{-m}}, \Si_{y_n}$ is $\le M_k$ and, hence, each interval 
 $\al_{-m,n}= [y_{-m} y_n]_{X_u}\subset \al$ bounds a $(k, M_k)$-narrow carpet  $\A(\al_{-m,n})$ in $X$. According to  Part 2,  then   $\D_{Y,X}(\xi^-)=\D_{Y,X}(\xi^+)$. \qed

\begin{rem}\label{rem:carpet-bound} 
Note that in Part 2 of the theorem we can bound from above  lengths of the intervals $\llbracket u, w_{m,n}\rrbracket= \pi(\A(\al_{-m,n}))$ in terms of lengths of the segments 
$\al_{-m,n}$: In view of the exponential flaring satisfied by $\X$ (see Lemma \ref{lem:exp} and Proposition \ref{hyp to lin flaring}) 
there exists a constant $\la$ (depending only on $K$ and the parameters of $\X$) such that 
$$
\length(\al_{-m,n}) \ge \la^{t} C, \quad t= d_T(u, w_{m,n}).
$$
Thus, 
$$
d_T(u, w_{m,n}) \le \tau_{m,n}:= \log_{\la} (\length(\al_{-m,n})) - \log_{\la}(C).$$ 
\end{rem}

As an application of Part 1 of the theorem we obtain the following: 

\begin{cor}
Suppose that $\gamma$ is a geodesic in $Y$ whose projection to $T$ is unbounded. Then $\gamma$ is not a leaf of the CT-map $\D_{Y,X}$. 
\end{cor}

\begin{comment}
We should write the corresponding group theoretic application of the above theorem in
the following section named 'group theoretic applications'. Rewrite part 1 for graphs of groups and their subgraphs. (done, Corollary \ref{cor:CT-lamination-gr-of-groups}) Also interpret parts 2 and 3 for semidirect products $H\rtimes_\phi \Z$: Two pseudo-orbits starting at $\al(\pm n)$ come uniformly close to each other. Can one use just orbits instead of pseudo-orbits? 

\end{comment}

\section{Boundary flows and CT laminations} \label{sec:Boundary flows and CT laminations}

In this section we will be using the notion of {\em ideal boundary flows} $Fl_t$ (where $t$'s are vertices and edges of the tree $T$) which are defined and discussed in Section \ref{sec:Boundary flows}.

\medskip
Before proving the next proposition we will need two definitions which will be discussed in much greater detail in Section \ref{sec:CT-lamination}. 
Recall that $\geo(Z,X)$ denotes the {\em limit set} (relative ideal boundary) of a subset $Z\subset X$. Recall also that $\La(Y,X)$ denotes the {\em Cannon--Thurston lamination} 
of a hyperbolic subspace $Y$ in a hyperbolic space $X$, see Definition \ref{defn:CT-lamination}. 

\begin{defn}\label{defn:CT-xi}
Consider a point $\xi\in \geo T$ and a geodesic ray  $v\xi$ in $T$ joining a vertex $v\in T$ to $\xi$. Define {\em $\xi$-relative ideal boundary} 
$$\geo^{\xi}(X_v,X):=\{
\eta\in \geo(X_v,X): \exists  
\mbox{~~a qi section $\gamma$ over}\,\, v\xi,~~ \gamma(\infty)=\eta\}$$
and 
$$
\La^\xi(X_v,X)=\{ \{z^-, z^+\}\in \La(X_v,X): \D_{X_v,X}(z^\pm)\in \geo^{\xi}(X_v,X)\}.
$$
\end{defn}

Note that, by the definition of $\La(X_v,X)$, $\D_{X_v,X}(z^+)=\D_{X_v,X}(z^-)$. 

\medskip
While the definition of $\La^\xi(X_v,X)$ at this point looks rather unmotivated, in the next section we will prove that it equals the 
the {\em $\xi$-ending lamination} $\La(X_v, X_{v\xi})$. Examples of points in $\geo^{\xi}(X_v,X)$ are given by points $\ga(\infty)$ for qi sections $\gamma$ that are limits (as $n\to\infty$) of bottom sections of carpets $\A(\al_{-m_0,n})$ appearing in the proof of Theorem \ref{thm:CT-fibers}. 

\begin{comment}

However, $\geo^{\xi}(X_v,X)$ might be nonempty even if $X_v$ 
is quasiconvex in $X$. {\mini For instance, if $\X$ is the tree of spaces corresponding to a nontrivial finite graph of hyperbolic groups satisfying Axiom {\bf H} and the acylindricity condition, then for each $X_v$ there will be some $\xi$ with $\geo^{\xi}(X_v,X)\ne \emptyset$, while each $X_v$ is quasiconvex in $X$. } 
{\mini We will prove in Proposition \ref{flow of CT leaf}(4) below that $\D_{X_v,X} (\Lambda^\xi(X_v,X))\subset \geo^{\xi}(X_v,X)$, i.e. $\geo^{\xi}(X_v,X)$ contains the image of the 
$\xi$-ending lamination under the CT-map. }

{\mini in general, the inclusion is strict, give an example: $X_v$ is a line whose ideal end-points $\eta_\pm$ 
have different CT-images; $T$ is a ray; half-line $y_0\eta_+$ of $X_v$ is covered by carpets with $\diam\to\infty$, while the other half-line has vertical decomposition into infinitely many intervals, in particular, there is an infinite $K$-section $\gamma$ through $y_0$. Then $\xi_+=\gamma(\infty)$ but $\xi_+$ is not in the image of the CT-lamination. Maybe even give example with hyperbolic plane. Also, connect to conical/nonconical limit points and note that $\D_{X_v,X} (\Lambda^{\xi}(X_v,X))$ consists of nonconical limit points while $X_v$ might be undistorted, while $\geo^{\xi}(X_v,X)$ nonempty for some $\xi$} 

\end{comment}

\begin{prop}\label{flow of CT leaf}
Fix a point $\xi\in \geo T$ and a pair of distinct points $z^+, z^+\in \partial_{\infty} X_v$ 
such that $\{z^+,z^+\}\in \Lambda^{\xi}(X_v,X)$. Let $L_v=\al\subset X_v$ be a biinfinite geodesic 
asymptotic to the points $z^\pm$. 

(1) For all vertices and edges $t$ in the ray $v\xi$ we have 
$$
Fl_t(\{z^\pm\})= \{z^\pm_t\}\ne \emptyset.  
$$

(2) For each vertex/edge $t$ in the ray $v\xi$  
let $L_t$ be a biinfinite geodesic in $X_t$ connecting 
the points $z^\pm_t$. Then the collection of such geodesics forms the union of vertex/edge sets of a metric bundle $\L=(\pi: L\to v\xi)\subset \X$ over  $v\xi$, which is also a $K'$-ladder for some $K'$.

(3) There are constants $K_1, C_1$ such that for each $n$, 
the segment $\al_{-n,n}=[\alpha(-n)\alpha(n)]_{X_v}$ bounds a $(K_1,C_1)$-narrow carpet 
${\mathfrak B}^n\subset \L$, where $\La$ is the ladder from (2).

(4) Every qi section over $v\xi$ contained in $L$ is asymptotic to 
$z=\D_{X_v,X}(z^-)= \D_{X_v,X}(z^+)\in \geo X$. In particular, any two such qi sections are 
at a finite Hausdorff distance from each other. 

(5) $\partial_{\infty} L$ is the singleton $\{z\}$, in particular, $\geo(L, X)= \{z\}$. 
\end{prop}
\proof We let $\Y=(\pi: Y\to v\xi)$ denote the restriction of $\X$ to the ray $v\xi$. We will need: 

\begin{lem}
$\D_{X_v,Y}(z^+)= \D_{X_v,Y}(z^-)$.
\end{lem} 
\proof We consider a sequence of $(K,C)$-narrow carpets $\A^n=\A(\al_{-n,n})= (\pi: A^n\to \llbracket v, w_n\rrbracket)$ in $\X$ given by Theorem \ref{thm:CT-fibers}(1) and bounded by 
the intervals  
$$
\al_{-n,n}= [\al(-n) \al(n)]_{X_u}\subset \al.$$
Set 
$$
\llbracket v, v_n\rrbracket= v\xi\cap \llbracket v, w_n\rrbracket. 
$$

i. Let us first verify that $\lim_{n\to\infty} v_n=\xi$. Suppose not. After passing to a subsequence, the the sequence $(v_n)$ would have to be constant. For sufficiently large each $n$, 
there is a point $x_n$ in the narrow end of $\A^n$ and a  $K$-section $\gamma_{x_n}$ in $A^n$ 
connecting $x_n$ to $\ga(v)$. The concatenation $\ga_{x_n}\star \ga$ is then a uniform quasigeodesic in $X$ (since $v_n$ is fixed). But then $\lim_{n\to\infty} x_n\ne \ga(\infty)$ in $\geo X$, contradicting the assumption that $\{z^+,z^-\}\in \Lambda^{\xi}(X_v,X)$. 

ii. The carpets $\A^n$ define uniform quasigeodesics $c_{-n,n}$ in $X$ connecting points $\al(\pm n)$ (see the proof of Theorem \ref{thm:CT-fibers}). 
Applying the cut-and-replace procedure to the paths $c_{-n,n}$ with respect to the subtree of spaces $Y\subset X$, we obtain uniform quasigeodesics $\hat{c}_{-n,n}$ in $Y$ 
which project to the intervals $\llbracket v, v_n\rrbracket$. Since $\lim_{n\to\infty} d_T(v,v_n)=\infty$, it follows that 
$$
d_Y(\al(0), \hat{c}_{-n,n})=\infty,
$$
implying that $\D_{X_v,Y}(z^+)= \D_{X_v,Y}(z^-)$. \qed

\medskip
We now begin the proof of the proposition. In view of the lemma, we can replace $\X$ with $\Y$.

(1) Since $\D_{X_v,Y}(z^+)= \D_{X_v,Y}(z^-)$, 
Theorem \ref{thm:CT-fibers}(1) implies that for each  $n>0$ the interval  
$$
\al_{-n,n}= [\al(-n) \al(n)]_{X_u}\subset \al$$
bounds a $(K,C)$-narrow carpet $\A^n=\A(\al_{-n,n})= (\pi: A^n\to \llbracket v, w_n\rrbracket)$ in 
$\Y$ for some $K, C$ depending only on the parameters of $\X$, $w_n\in v\xi$. Since for each vertex $t\in v\xi$ 
$$
\Hd_{Y}(A^n_v, A^n_t)\le K d_T(v,t),
$$
it follows that the entire geodesic $L_v$ is contained in the $K d_T(v,t)$-neighborhood of the vertex-space $X_t$. Lemma \ref{lem: flow condition} now implies that 
$$
Fl_t(\{z^\pm\})= \{z^\pm_t\}\ne \emptyset. 
$$
This proves (1).  

\medskip 
(2) By the construction, for each pair of vertices $s, t\in v\xi$ with $d_T(s,t)=1$, the geodesics $L_s, L_t$ are asymptotic to the same pair of points in $X_{st}$. Therefore, $\Hd(L_s, L_t)\le D$ for 
some $D$ depending only on $L'_0$ and the hyperbolicity constant $\delta'_0$ of $X_{st}$. The same applies to the geodesic $L_e\subset X_e$, $e=[s,t]$. This proves that the union of $L_t$'s, $t\in V(v\xi)$, $L_e$'s, $e\in E(v\xi)$ forms a metric bundle $\L$ in $\Y$. This metric bundle has structure of a $K'$-ladder according to Lemma \ref{lem:E-ladder-structure}.   

\medskip 
(3) As noted in the proof of (1), we already have the carpets $\A^n$ in $\Y$ bounded by the segments $\al_{-n,n}$. The trouble is that these carpets need not be contained in the ladder $\L$.  However, according to Theorem \ref{thm:mitras-projection}, there is a uniformly coarse Lipschitz projection $\nu: \Y\to \L$ which, for every vertex $t\in v\xi$, equals  to the restriction of the nearest-point projection $P_{X_t, L_t}$. Taking the corresponding modified projection $\bar{P}_{X_t, L_t}(A_t^n)$ (see Definition \ref{def:modified projection}) we then obtain a collection of subsegments $B_t^n\subset L_t$ satisfying axioms of a $(K_1,C_1)$-narrow carpet over $\llbracket v, w_n\rrbracket$, which we denote ${\mathfrak B}^n$.

\medskip 

(4) Suppose that $\gamma$ is a $k$-qi section of $\L$ over the ray $v\xi$, $p=\gamma(v)\in \al$. 
Then $p$ belongs to the segment $\al_{-n,n}$ for all sufficiently large $n\ge n_p$. We let $c_n$ denote the uniform quasigeodesic in $Y$ connecting the end-points of $\al_{-n,n}$ and equal to the 
concatenation of the two horizontal boundary sections of ${\mathfrak B}^n$ along with
its narrow end (the vertical geodesic segment in $X_{w_n}$). Since 
both sequences $(\al(n)), (\al(-n))$ converge to the same point $z\in \geo Y$, it follows that the sequence of quasigeodesics $(c_n)$ also converges to $z$. For each  
$n\ge n_p$, there exists a vertex $v_n\in \llbracket v, w_n\rrbracket$ such that 
$$
d_{X_{v_n}}(\gamma(v_n), x_n)\le K_1+K', 
$$
where $x_n$ lies in $c_n$. Clearly, 
$$
\lim_{n\to\infty} d_Y(\alpha(0), x_n)= \infty. 
$$
In particular, the sequence $(x_n)$ converges to $\gamma(\infty)$. Since the sequence $(x_n)$ also converges to $z$, we obtain $\gamma(\infty)=z$.   

(5) The proof of this part is similar to that of (4). Consider a geodesic ray $\beta$ in $L$, $\beta(0)=\al(0)=p$. Then for each $n$, there is a point $x_n$ in $\beta$ within distance $K'$ from the path $c_n$ in ${\mathfrak B}^n$ defined as in the proof of (4). Again, $d(p, x_n)\to \infty$ and, hence, 
$\beta(\infty)= \lim_{n\to\infty} c_n=z$.  \qed 
 
\medskip 
We now can relate boundary flows to CT-laminations: 

\begin{prop}\label{prop:CT surjective}
Suppose that $\X=(\pi: X\to T)$ is a tree of hyperbolic spaces with hyperbolic total space $X$, $\Y=(\pi: Y\to S)\subset \X$ is a subtree of spaces 
and $v\in V(S)$  satisfies the following conditions:

1. $\geo X_w\subset Fl(\geo X_v)$ for each vertex $w\in V(S)$. 

2.  $\D_{X_v,X}(\geo X_v)=\geo X$. 
 
\noindent Then the  CT-map $\D_{X_v, Y}: \partial_{\infty} X_v\map \partial_{\infty} Y$ is also surjective.
\end{prop}
\proof We claim that each $z\in \geo Y$ belongs to $\D_{X_v,Y}(X_v)$. Since $\D_{X_v,X}$ is surjective, 
there is $z_1\in \geo X_v$ such that $\D_{X_v,X}(z_1)=\D_{Y,X}\circ \D_{X_v,Y}(z_1)=\D_{Y,X}(z)$.
Thus, for $z'= \D_{X_v,Y}(z_1)\in \geo Y$, $\D_{Y,X}(z')=\D_{Y,X}(z)$. 

If $z'=z$, then we are done. If not, then by the description of the fibers of the CT-map $\D_{Y,X}$ given in Proposition \ref{prop:CT-fibers1}, 
there exists a vertex-space $X_w\subset Y$ such that a geodesic $\beta=zz'\subset Y$ asymptotic to  $z, z'$  is Hausdorff-close to a geodesic $\al\subset X_w$. In particular, 
$z\in \geo(X_w,Y)$. The first assumption of the proposition implies that $z\in \geo(X_v,Y)$. Since $\geo(X_v,Y)= \D_{X_v, Y}(\geo X_v)$ (see Lemma \ref{lem:image-of-CT}), 
the claim follows.  \qed

\section{Cannon--Thurston lamination and ending laminations}
\label{sec:CT-lamination}

In this section we shall significantly expand on Theorem \ref{thm:CT-fibers}(1); many of our results are generalizations of the ones proven by Mitra in  \cite{MR1445392}. 
Throughout this section we will assume that $\X =(\pi:X\map T)$ is a tree of hyperbolic spaces with proper hyperbolic total space $X$. 

To motivate the discussion, we recall Thurston's notion of the {\em ending laminations} in the setting of hyperbolic 3-manifolds. (We refer the reader for a detailed overview of end-invariants of hyperbolic 3-manifolds to Minsky's surveys \cite{minsky-cdm} and \cite{MR2044543}.) 
For simplicity of the discussion, we consider a noncompact complete connected hyperbolic 3-manifold $M$ with finitely-generated fundamental group   that does not split as a nontrivial free product and such that $M$ has positive injectivity radius. The manifold $M$ contains a (unique up to isotopy)   compact submanifold with smooth boundary $M_c$ (the {\em compact core} of $M$), such that the complement $M\setminus \int(M_c)$ is homeomorphic to $\partial M_c \times \RR_+$. The group $G=\pi_1(M_c)\cong \pi_1(M)$ is hyperbolic and the assumption that it does not split as a free product implies that for each surface component $S\subset \partial M_c$, the inclusion map $S\to M_c$ is $\pi_1$-injective. 

Each component $E= S\times \RR_+$ of $M\setminus M_c$ is an {\em end} of $M$; the surface $S$ is a component of $\partial M_c$; it is a compact surface which admits a hyperbolic metric (which we fix from now on).  For each end $E= S\times \RR_+$ one defines an {\em ending lamination} $\la=\la(E)$ of $E$, which is a certain nonempty compact subset of $S$, equal to a disjoint union of complete geodesics in $S$. Lifting $\la$ to the 
universal covering space of $S$, one obtains a $\pi_1(S)$-invariant closed subset of $\tilde{S}\cong \H^2$ equal to the disjoint union of geodesics. Each geodesic $\beta$ in $\tilde\la$ is uniquely determined by an unordered pair $\{\xi^+,\xi^-\}\in \geo^{(2)}\H^2$, the ideal boundary points such that $\beta= \xi^- \xi^+$. We, thus, identity $\la$ with a 
$\pi_1(S)$-invariant closed subset of $\geo^{(2)}\H^2$ consisting of such pairs. Consider a component $\tilde{E}$ of the preimage of $E$ in $\H^3$ (the universal covering space of $M$); the boundary surface $\tilde{S}$ of $\tilde{E}$ is a copy of the universal covering space of $S$. It was proven by Minsky \cite{minsky-jams} (under the above assumption on the injectivity radius of $M$) and Mj \cite{MR3652816} in full generality, that each inclusion map $\tilde{S}\to \tilde{E}$ has a CT-map; these CT-maps combine in a CT-map for the inclusion $\tilde{M}_c\to \tilde{M}=\H^3$ (where $\tilde{M}_c$ is the universal covering space of $M_c$ equipped with the pull-back Riemannian metric). Each ending lamination $\la(E)$ is $\pi_1(S)$-equivariantly homeomorphic to the CT-lamination $\La(\tilde{S},\tilde{E})$ and the union of $G$-orbits of these laminations in $\geo \tilde{M}_c$ is the CT-lamination $\La(\tilde{M}_c, \H^3)$. 

We now relate this discussion to trees of spaces. For each end $E$, the space $\tilde{E}$ (with its intrinsic Riemannian path-metric) is $\pi_1(E)$-equivariantly quasiisometric to the total space of a certain metric bundle $\X_E= (E\to \RR_+)$, with vertex and edge-spaces isometric to the hyperbolic plane. 
(This bundle structure is implicit in \cite{minsky-jams}. It is obtained via pull-back of the universal bundles over the Teichm\"uller spaces of boundary surfaces of $M_c$.)  
Putting these spaces together, we obtain a tree of spaces 
$\X=(X\to T)$ (on which $G=\pi_1(M)$ is acting) which has a distinguished vertex $v$ fixed by $G$. The total space $X$ of $\X$ is isometric to $\H^3$. 
The tree $T$ is a union of geodesic rays; the intersection of any two distinct rays in this collection is the vertex $v$. Thus, the $G$-orbits of ending laminations $\la(E)$ in 
$\geo^{(2)}G$ can be described as CT-laminations 
$$
\La(X_v,X_{v\xi})\subset \geo^{(2)} X_v,
$$ 
where, $\xi$'s are the ideal boundary points of $T$ and $X_{v\xi}$ is the total space of the pull-back of $\X$ to the ray $v\xi$ in $T$. The result stated above, relating ending  laminations of $M$ with the CT-lamination $\La(\tilde{M}_c, \H^3)$ can then be restated as:
$$
\La(X_v, X)= \bigcup_{\xi\in \geo T} \La(X_v, X_{v\xi}).
$$

In the context of general trees of hyperbolic spaces, points at infinity $\xi\in \geo T$ play the role of ends of the hyperbolic manifolds and, accordingly, ending laminations are 
defined as CT-laminations $\La(X_v, X_{v\xi})$.  The main goal of this section is to prove an analogue (actually, a sharper version) of the above equality in the setting of more general trees of spaces,  
Theorem \ref{thm:CT lamination} below. In particular, we will also prove that for each $\xi\in \geo T$, the ending lamination  $\La(X_v, X_{v\xi})$ equals the subset $\La^\xi(X_v, X)\subset \La(X_v,X)$ defined in the previous section. 
 This alternative interpretation of 
the ending lamination  $\La(X_v, X_{v\xi})$ will be used in several places, e.g. proof of Theorem \ref{thm:CT lamination}, Parts (4) and (5) and proof of 
Proposition \ref{conical-singleton}.

\medskip 
The results below are motivated by similar results
obtained in \cite[section 6.2]{Krishna-Sardar}; our  notation and  proofs are similar (see also \cite{MR1445392, MR3186670}).

\begin{lemma}\label{lemma1: CT lam}
Suppose that $\X=(\pi:X\map T)$ is a tree of hyperbolic spaces with proper hyperbolic total space $X$. 
We fix a vertex $v\in V(T)$ and $K\ge 1$. 

1. For every $\xi\in \geo T$, there is a geodesic ray $\rho$ joining $v$ to $\xi$.

2. There is $K_1$ depending on $K$ and the parameters of $\X$ such that the following holds. 
Let $(w_n)$ be a sequence of points in $V(T)\cup \partial_{\infty} T$ and let $(\gamma_n)$ be a sequence of $K$-qi sections of $\X$ of over the geodesic  $vw_n$. Suppose 
that the sequence $(\gamma_n(v))$ belongs to a bounded subset $B$ of $X$ and the sequence $(\gamma_n(w_n))$ converges to a point  
$\eta\in \partial_{\infty} X$. Then the sequence 
$(w_n)$ converges to a point $\xi\in \geo T$ and there is a $K_1$-qi section $\gamma$ over the geodesic $v\xi$ such that $\eta=\gamma(\infty)$ and $\gamma(v)\in B$.
\end{lemma}
\proof 1. Let $(w_n)$ be a Gromov-sequence of vertices in $T$ representing the point $\xi$. Then 
$$
\lim_{m,n\to\infty} d_T(v, w_mw_n)=\infty. 
$$
It follows that the union of geodesic segments $vw_n, n\in \NN$, is a locally finite subtree $S\subset T$. Therefore, the sequence of segments $vw_n$ subconverges to a geodesic ray $\rho$ in $S$ emanating from $v$. In order to prove that $\rho$ joins $v$ to $\xi$ we note that for each $m\in \N$ and all sufficiently large $n$, the points $t_m=\rho(m)$ satisfy 
$$
d_T(v, t_mw_n)= d_T(v, t_m)=m.  
$$
Thus, the sequence $(t_m)$ is a Gromov-sequence equivalent to $(w_n)$. 

2. Pick a base-point $x_0$ in a bounded subset $B\subset X$ of diameter $D$ containing all the points $\gamma_n(v)$, e.g. we can take $x_0=\gamma_1(v)$. 

Since the sequence $(\gamma_n(w_n))$ converges  $\eta\in \geo X$, we have 
$$
\lim_{n\to\infty} d_T(v, w_n)=\infty. 
$$
Moreover, since the sequence of geodesic segments $\gamma^*_i=[x_0 \gamma_n(w_n)]_X$ {\em coarsely converges} to a  geodesic ray $x_0\xi$ (see \cite[Definition 8.32]{Drutu-Kapovich}), there is a constant $C$ (depending only on $K, D$ and the hyperbolicity constant of $X$) such that for each $R$, there is a number $n_0$ such that for all $m, n\ge n_0$ the Hausdorff-distance between $\gamma_m\cap B(x_0,R)$ and $\gamma_n\cap B(x_0,R)$ is $\le C$.  Since $\gamma_i$'s are sections over geodesic segments 
$vw_i$ in $T$, it follows that 
$$
\lim_{i\to\infty} \sup \{R: vw_m\cap B(v,R)= vw_n\cap B(v,R) ~~\forall m, n\ge i\}=\infty. 
$$
In particular,  $(w_n)$ is a Gromov-sequence in $T$ converging to some $\xi\in \geo T$ and by Part 1 of the lemma, the   sequence of segments $(vw_n)$ converges to the ray $v\xi$.

Furthermore, in view of properness of $X$, the sequence of $K$-qi sections $(\ga_n)$ 
subconverges to a $K$-qi section $\gamma$ over the geodesic ray $v\xi$ in $T$ 
(this is a coarse version of the Arzela--Ascoli Theorem, cf. \cite[Proposition 8.34]{Drutu-Kapovich}). \qed

 \medskip

Fix $K\geq K_0$, pick a vertex $v\in V(T)$ and let $\Y=(\pi: Y\to S)$ be a $(K,D,E,\la)$-semicontinuous family in $\X$ relative to a vertex $v\in S\subset T$, see Definiton 
\ref{defn:scfamily}. The following proposition is motivated by the results of \cite{Krishna-Sardar} and \cite[Proposition 8.2]{MR3186670}:

\begin{prop}\label{boundary of K-flow}
We have
$$
\geo Y= U:=\geo(X_v,Y)\cup \left(\bigcup_{\xi\in \partial_{\infty} S} \geo^{\xi} (X_v, X)\right).$$
\end{prop}

\proof Recall that, according to Theorem \ref{thm:mitras-projection}, $Y$ is qi embedded in $X$. Thus, we will identify $\geo Y$ with a subset of $\geo X$. 
Since $U$ is obviously contained in $\geo Y$, we only have to prove that every point $z\in \geo Y$ lies in $U$. Fix a base-point $x\in X_v$. 
Suppose that $x_n\in Y_{v_n}$ is a sequence of points converging to $z$. 
Let $\gamma_n$ be a $K$-qi section in $Y$ over $vv_n$, joining $x_n$ to $y_n\in Y_v$.  

(i) Suppose first that $(y_n)$ is a bounded sequence. Then by Lemma  \ref{lemma1: CT lam}(2), 
the sequence $(v_n)$ converges to  some $\xi\in \geo T_v$ 
and the sequence $(\gamma_n)$ coarsely converges to a $K$-qi section $\gamma$ in $Y$ over the ray 
$v\xi$, so that $x_n\map \gamma(\infty)$. Thus $z=\lim_{n\to\infty} x_n\in  \geo^{\xi} (X_v,X)\subset U$ in this case.

\medskip 
(ii) Consider now the case when $(y_n)$ is an unbounded sequence. After extraction, we can assume that $(y_n)$ converges to some $z'\in \geo (X_v,X)= \geo(X_v,Y)$ 
(since $Y$ is quasiconvex in $X$). We claim that $z'=z$. Since each $\ga_n$ is a $K$-qi section, it suffices to show that $d(x, \ga_n)\to \infty$. Suppose that the sequence  
 $d(x, \ga_n)$ is bounded, and $p_n=\gamma_n(v_n), v_n\in V(\llbracket v, \pi(y_n)\rrbracket)$, is a sequence such that  
 $d_X(x, p_n)\le C$ for all $n$.  Since $d_T(v,v_n)=d_T(\pi(x), \pi(p_n))\le d_X(x, p_n)\le C$, it follows that 
 $$
 d_X(y_n, p_n)\le CK, d_X(x, y_n)\le C + CK,
 $$
 contradicting the assumption that the sequence $(y_n)$ is unbounded. Thus, 
 $$
 z=z'\in \geo(X_v,Y)\subset U. \qedhere$$

\begin{cor} For any ladder $\mathfrak L=(\pi: L\map \pi(L))$ in $\X$, centered at a vertex $v\in V(T)$, we have 
$$
\geo (L, X)= \bigcup_{\xi\in \geo \pi(L)}\, \geo^{\xi} (L_v,X).$$
\end{cor}
\proof Since $L_v$ is a finite geodesic segment, $\geo(L_v, L)=\emptyset$, and, thus, the corollary is an immediate consequence of Proposition \ref{boundary of K-flow}. 
\qed

\medskip 
We now return to the discussion of properties of ending laminations $\La(X_v, X_{x\xi})$ and their relation to the CT-laminations $\La(X_v,X)$ and their subsets $\La^\xi(X_v,X)$ defined in the previous section (see Definition \ref{defn:CT-xi}).

\begin{thm}[Properties of ending laminations]\label{thm:CT lamination}
Suppose that $X=(\pi:X\map T)$ is a tree of hyperbolic spaces with hyperbolic and proper total space $X$. 
There exists $K$ depending only of the parameters of $\X$ and the hyperbolicity constant of $X$ such that the following hold: 
\begin{enumerate}
\item 
Let $v$ be a vertex of $T$ and let $\alpha:\RR\map X_v$ be a complete geodesic in $X_v$ such that $\{\alpha(-\infty), \alpha(\infty)\}\in \La(X_v,X)$.  
Then for both $z\in \{\al(\pm\infty)\}$, there exists  point $\xi\in \geo T_z$ such that for each  $p\in \al$, there is a $K$-qi section $\gamma$ over $v\xi$ satisfying $\gamma(v)=p$ and $\gamma(\infty)=\D_{X_v,X}(\alpha(\infty))$. In other words, for each point $\{z^-,z^+\}\in \La(X_v,X)$, there exists $\xi\in \geo T$ such that 
$$
\D_{X_v,X}(z^\pm)\in \La^\xi(X_v,X) 
$$
and, thus,
$$
\La(X_v,X)= \bigcup_{\xi\in \partial_{\infty} T_v} \La^\xi(X_v,X).
$$

\item 
For each $\xi\in \geo T$ we have
$$
\Lambda^{\xi}(X_v, X)=\Lambda(X_v, X_{v\xi})=\Lambda(X_v, Fl_{K}(X_v)\cap X_{v\xi}).
$$

\item 
Each $\La^\xi(X_v,X)$  
is a closed subset of $\partial^{(2)}_{\infty} X_v$.

\item 
Suppose $\xi_1\neq \xi_2\in \partial_{\infty} T$, and $\alpha_1, \alpha_2$ complete geodesics in $X_v$
such that $\{z_i^-, z_i^+\}=\{\alpha_i(-\infty), \alpha_i(\infty)\} \in \Lambda^{\xi_i}(X_v, X)$, $i=1,2$. Then the subsets 
$\{z_1^-, z_1^+\}$, $\{z_2^-, z_2^+\}$ of $\geo X_v$ are disjoint. In particular,  the ending laminations $\Lambda^{\xi_1}(X_v,X)$, $\Lambda^{\xi_2}(X_v,X)$ are disjoint and the point 
$\xi$ in (1) is uniquely determined by $\{z^-,z^+\}\in \La(X_v,X)$.

\item 
Ending laminations $\La^\xi$ depend upper semicontinuously\footnote{It is shown in \cite[section 7]{MR3557464} that in general $\La^\xi$ does not depend continuously on $\xi$.} on $\xi$: 
Suppose that $\xi_n\map \xi$ in $\geo T_v$, $\{z^+_n, z^-_n\}\in \Lambda^{\xi_n}(X_v, X)$ and $\{z^+_n, z^-_n\}\map \{z^+,z^-\}\in \geo^{(2)} X_v$.
Then $\{z^+,z^-\}\in \Lambda^{\xi}(X_v, X)$. 
 
\item 
If $\xi_1\neq \xi_2\in \partial_{\infty} T_v$, then any two leaves $\al^i$ of $\Lambda^{\xi_i}(X_v, X)$, $i=1,2$, are uniformly cobounded in $X_v$. Namely,   
given $D>0$  there exists $R=R(D)$ (independent of $\al^1, \al^2$ but possibly depending on $\xi_1, \xi_2$)  
such that $\alpha^1\cap N_D(\alpha^2)$ has diameter  $\le R$. 

\end{enumerate}
\end{thm}
\proof
(1) The proof follows the argument in  the proof of Proposition \ref{flow of CT leaf}(1) (or Part (3) of that proposition). By Theorem \ref{thm:CT-fibers}(1), 
for each  $n>0$ the interval $$ \al_{-n,n}= [\al(-n) \al(n)]_{X_u}\subset \al$$
bounds a $(K,C)$-narrow carpet $\A^n=\A(\al_{-n,n})$ in $X$, which, in turn, defines a uniform quasigeodesic $c_n= c(\A^n)$ in  $\A(\al_{-n,n})$ connecting the points $\al(-n), \al(n)$. 
Take any point $p\in \al$. Then for all $n\ge n_p$, $p$ belongs to the segment $\al_{-n,n}$. Since $\A^n$ is a $K$-metric bundle, there exists a $K$-section $\ga^n$ over 
the segment $\llbracket v,w_n\rrbracket= \pi(A^n)$. The end-point $x_n=\ga^n(w_n)$ belongs to $c_n$, which implies that the sequence $(x_n)$ converges to the limit point 
$\D_{X_v,X}(\al(-\infty))= \D_{X_v,X}(\al(\infty))$.  Then the existence of the point $\xi$ follows from Lemma \ref{lemma1: CT lam}(2): Proposition \ref{flow of CT leaf}(1) implies that $\xi\in \geo T_z$, $z=\al(\pm\infty)$. The rest of the assertions of Part (1) follow  immediately from the definition of $\La^\xi(X_v,X)$. 

\medskip 
(2) The inclusion $\La^\xi(X_v,X)\subset \La(X_v,X_{v\xi})$ is clear from the definition of $\La^\xi(X_v,X)$. The opposite inclusion is a direct consequence of Part (1) of the theorem. The inclusion $\Lambda(X_v, Fl_{K}(X_v)\cap X_{v\xi})\subset \Lambda(X_v, X_{v\xi})$ for every $K\ge K_0$ is clear from the fact that $Fl_{K}(X_v)$ is quasiconvex in $X$. Suppose that $\{z^-,z^+\}\in \Lambda(X_v, X_{v\xi})$.  Proposition \ref{flow of CT leaf}(1) implies that for each vertex $t\in v\xi$, $Fl_t(\{z^\pm\})\ne \emptyset$.  
By Lemma \ref{lem:flow-flow}(1), there exists $K$ such that 
for each vertex/edge $t$ in $v\xi$ there exists a biinfinite geodesic $L_t\subset Fl_K(X_v)$ asymptotic to the points $Fl_t(\{z^\pm\})$. By Proposition \ref{flow of CT leaf}(2), these geodesics form a ladder $\L=(\pi: L\to v\xi)$ in $X$. By Part 5 of the same proposition, $\geo L$ is a singleton, which implies 
$$
\D_{X_v,Fl_K(X_v)\cap X_{v\xi}}(z^+)= \D_{X_v,Fl_K(X_v)\cap X_{v\xi}}(z^-). 
$$
In other words, 
$$
\{z^-,z^+\}\in \La(X_v, Fl_K(X_v)\cap X_{v\xi}).$$

(3) In view of Part (2), the claim follows from the fact that the CT-lamination $\La(X_v, X_{v\xi})$ is closed in $\geo^{(2)} X_v$. 

\medskip 
(4) By Proposition \ref{flow of CT leaf}(2) and (4) there are qi sections $\gamma_i$ over $v\xi_i$ asymptotic to $\D_{X_v,X}(z^\pm_i)$, $i=1,2$.
Since $\xi_1\neq \xi_2$, $\Hd(v\xi_1, v\xi_2)=\infty$, which, in turn, implies that 
$\Hd(\gamma_1, \gamma_2)=\infty$. 
 Thus $\gamma_1(\infty)\neq \gamma_2(\infty)$,  
whence $\{z_1^-, z_1^+\}\cap \{z_2^-, z_2^+\}=\emptyset$.

\medskip 
(5) Since $\D_{X_v,X}(z^+_n)= \D_{X_v,X}(z^-_n)$ and $\D_{X_v,X}$ is continuous, we have $\D_{X_v,X}(z^+)=\D_{X_v,X}(z^-)$.
Thus $\{z^+,z^-\}\in \Lambda(X_v,X)$. Since $z^+\ne z^-$, geodesics $\al_n$ in $X_v$ connecting  the points $z_n^\pm$ all intersect a certain bounded subset of $X_v$. Hence, we can parameterize these geodesics so that the sequence $(\al_n(0))$ is bounded in $X_v$. In view of properness of $X_v$, by the Arzela--Ascoli theorem, the sequence of geodesics $\al_n$ subconverges to a geodesic $\al$ in $X$. This geodesic is necessarily asymptotic to the 
points $z^\pm$, see e.g. \cite[Theorem 11.104]{Drutu-Kapovich}. 
Since $\{z_n^+, z_n^-\}$ belongs to $\La^\xi(X_v,X)$, there exist uniform qi sections $\gamma_n$ over $v\xi_n$ connecting $\al_n(0)$ to 
$\D_{X_v,X}(z^\pm)$. By the continuity of the CT-map, $z_n\map z$ implies that $\D_{X_v,X}(z_n)\map \D_{X_v,X}(z)$. Accordingly, $\gamma_n(\infty)\map \D_{X_v,X}(z)$. 
Hence, by Lemma \ref{lemma1: CT lam}(2) 
there is a qi section $\gamma$ over $v\xi$ such that $\gamma(\infty)=\D_{X_v,X}(z^\pm)$, which implies that $\{z^+, z^-\}$ belongs to $\La^\xi(X_v,X)$.

\medskip 

(6) Note that this statement is a strengthening of Part (4) since that part is equivalent to the statement that the leaves $\al^1, \al^2$ are cobounded in $X_v$. Note also that, since vertex-spaces of $\X$ are uniformly properly embedded in $X$, the following two properties are equivalent for subsets $Y^1, Y^2\subset X_v$:

(i) There exists a function $R_v(D)$ such that $\diam_{X_v}(Y^1\cap N^{X_v}_D(Y^2))\le R_v(D)$. 

(ii) There exists a function $R(D)$ such that $\diam_{X}(Y^1\cap N_D(Y^2))\le R(D)$.

\medskip 
Set  $\al^i(\pm \infty)= z_i^\pm$, $i=1, 2$. Let $x_i^\pm \in \al_i,  i=1, 2$, be points such that 
\begin{equation}\label{eq:qua}
d_{X_v}(x^\pm_1, x^\pm_2)\leq D. 
\end{equation}
Our goal is to get an upper bound (in terms of $D$) on the distances $d_{X_v}(x_i^+, x_i^-)$, $i=1,2$.

a. We first consider the special case when the rays $v\xi_1, v\xi_2$ intersect only at the vertex $v$. We have subtrees of spaces $\Y^i= (\pi: X_{v\xi_i}\to v\xi_i)$ in $X$, $i=1, 2$. 
Since for $i=1, 2$, 
$$
\{z^-_i, z^+_i\}\in \La^{\xi_i}(X_v, X)= \La(X_v, X_{v\xi_i}),
$$
there is a sequence of $(K,C)$-narrow carpets 
$$
\A^{i,n}=\A(\al^i_{-n,n}) \subset \Y^i, n\in \N,$$ 
where $\al^i_{-n,n}$ is the subinterval in $\al^i$ between $\al^i(-n), \al^i(n)$. In particular, for all sufficiently large $n$, $x_i^\pm\in \al^i_{-n,n}$, $i=1,2$. Connect $x^\pm_i$ to 
the narrow end of $\A^{i,n}$ by a $K$-section $\ga^\pm_i$ in $\A^{i,n}$, $i=1,2$. Since $d_{X_v}(x_1^\pm, x^\pm_2)\le D$, 
both concatenations
$$
\phi^-:= \ga^-_1\star [x^-_1 x^-_2]_{X_v} \star \ga^-_2, \phi^+:= \ga^+_1\star [x^+_1 x^+_2]_{X_v} \star \ga^+_2
$$
are $k$-quasigeodesics in $X$ with $k$ depending only on $K$ and $D$. The respective end-points of these quasigeodesics are at most $C$-apart from each other. It follows that $d_{X_v}(x_i^+, x_i^-)\le R=R(k,\delta_X)$, cf. Lemmata \ref{lem:sub-close} and \ref{lem:hyp->uniform flaring}.  

\medskip 
b. We now consider the general case: The rays $v\xi_1, v\xi_2$ intersect along a finite subinterval $vw$ (this subinterval is finite since $\xi_1\ne \xi_2$). Since 
$\{z^-_i, z^+_i\}\in \La^{\xi_i}(X_v, X)$, there exist vertical geodesics $\beta^1, \beta^2$ in $X_w$ within uniformly bounded (in terms of $d_T(v,w)$) 
Hausdorff distance from $\al^1, \al^2$ respectively, see Proposition \ref{flow of CT leaf}. By Part (a), the geodesics $\beta^1, \beta^2$ are uniformly cobounded in $X_w$. 
It follows that $\al^1, \al^2$ are uniformly cobounded as well.  \qed

\begin{cor}
For each vertex $v\in T$ and 
$$
z\in \geo(X_v,X)\setminus (\bigcup_{\xi\in \geo T} \geo^{\xi}(X_v,X)),$$ 
the preimage $\D_{X_v,X}^{-1}(z)$ is a singleton. 
\end{cor}
\proof By Lemma \ref{lem:image-of-CT}, $\geo(X_v,X)= \D_{X_v,X}(\geo X_v)$. Thus, all we need is to show that $|\D_{X_v,X}^{-1}(z)|\le 1$.  
As we noted in the previous section, for each $\xi\in \geo T$, 
$$
\D_{X_v,X}(\La^\xi(X_v,X))\subset \geo^\xi(X_v, X). 
$$
Hence, $z\notin \D_{X_v,X}(\La^\xi(X_v,X))$ for any $\xi\in \geo T$. However, according to Theorem \ref{thm:CT lamination}(1),  
$$
\La(X_v,X)= \bigcup_{\xi\in \partial_{\infty} T_v} \La^\xi(X_v,X), 
$$
which means that there is no $\{z^+, z^-\}\in \La(X_v,X)$ satisfying $\D_{X_v,X}(z^\pm)=z$. Thus, $\D_{X_v,X}^{-1}(z)$ contains at most one point. \qed 

\begin{cor}\label{cor:nonempty point-flow}
1. If for each $\xi\in \geo T$, $\La^\xi(X_v,X)=\emptyset$, then $\La(X_v,X)=\emptyset$, i.e. the CT-map $\D_{X_v,X}$ is 1-1. 

2. If $z^\pm \in \geo X_v$ are distinct points such that $\D_{X_v,X}(z^+)=\D_{X_v,X}(z^-)$, then there exists $\xi\in \geo T$ such that for every vertex $w\in v\xi$, $Fl_w(\{z^\pm\})\ne \emptyset$. 
\end{cor}
\proof 1. The first claim is a direct consequence of   Theorem \ref{thm:CT lamination}(1). 

2. Since $\D_{X_v,X}(z^+)=\D_{X_v,X}(z^-)$, $\{z^+, z^-\}\in \La(X_v,X)$. By Theorem \ref{thm:CT lamination}(1), there  exists $\xi\in \geo T$ such that 
$\{z^+, z^-\}\in \La^\xi(X_v,X)$. Now the claim follows from  Proposition \ref{flow of CT leaf}(1). \qed

\section{Conical limit points in trees of hyperbolic spaces}

In this section we consider trees of hyperbolic spaces $\X=(\pi: X\to T)$ with proper total space $X$ and 
discuss the relation between conicality for limit points of  subtrees of spaces $Y\subset X$ and the CT-maps $\D_{Y,X}$. Namely, identifying $\La(Y,X)$ with a subset $\Si(Y,X)$ of 
$\geo Y$ equal  
$$
\bigcup_{\{z^+, z^-\}\in \La(Y,X)} \{z^+, z^-\},
$$
we'll see that $\partial_{Y,X}(\Si(Y,X))$ is disjoint from the conical limit set of $Y$ in $\geo X$. 

\smallskip
The next definition is motivated by the notion of conical limit points of group actions on hyperbolic spaces, see Definition \ref{defn:conical-limit}, as well as 
Definition 11.93 and Section 11.13.4 in \cite{Drutu-Kapovich}. 

\begin{defn} 
\index{conical limit point}
Suppose $X$ is an arbitrary hyperbolic geodesic metric space and $Y\subset X$. Then a point $\xi\in \geo(Y, X)\subset \partial_{\infty} X$
is called a {\em conical limit point} of $Y$ if for some (any) (quasi)geodesic $\alpha\subset X$ asymptotic to $\xi$ there is
$R>0$ and a sequence of points $\{y_n\}$ in $N_R(\alpha)\cap Y$ converging to $\xi$. The set of conical limit points of $Y$ is called the {\em conical limit set} of $Y$ in $\geo X$. 
\end{defn}

Thus, if $Y$ is an orbit $Gx$ of an isometric proper action $G\acts X$, then $\xi$ is a conical limit point of $Y$ if and only if it is a conical limit point of the $G$-action on $X$.

\begin{prop}\label{conical-singleton}
Suppose $\X=(\pi: X\to T)$ is a tree of hyperbolic spaces with proper and hyperbolic total space, and $\Y=(\pi: Y\to S)\subset \X$ is a subtree of spaces.
 Let $\D_{Y,X}: \geo Y \to \geo X$ be the CT-map. If $\eta \in \geo(Y,X)$ is a conical limit point of $Y$, 
then $|\D_{Y,X}^{-1}(\eta)| = 1$.  
\end{prop} 
\proof If not, then there are distinct points $z_\pm\in \geo Y$ such that 
$\D_{Y,X}(z_-)=\D_{Y,X}(z_+)=\eta$. Consider a geodesic $\beta$ in $Y$ 
asymptotic to the points $z_\pm$. 
By Proposition \ref{prop:CT-fibers1}, there exists a 
vertex space $X_v\subset Y$ and a complete geodesic $\alpha\subset X_v$ such that 
$\Hd(\al, \beta)<\infty$. Let $z_\pm'= \al(\pm \infty)$. 
It follows that $\D_{X_v,X}(z'_\pm)=\eta$. By Theorem 
\ref{thm:CT lamination}(1), there is a point $\xi\in \geo T$ and a qi section $\gamma$ over the ray $v\xi$, such that 
$\gamma(\infty)=\eta$. We claim that $\xi\in \geo T\setminus \geo S$. If not, then 
the ray $v\xi$ is contained in the subtree $S$.  
By Theorem \ref{thm:CT lamination}(2), 
$$
\{z'_-, z'_+\}\in \La^{\xi}(X_v,X)= \La(X_v, X_{v\xi})\subset \La(X_v, Y). 
$$
But then $\D_{X_v,Y}(z'_\pm)= z_\pm$ and $z_+=z_-$, contradicting our assumption that the points $z_\pm$ are distinct. 
 It then follows that $\lim_{n\to \infty} d_X(\gamma(n),Y)=\infty$. Since $\ga(\infty)=\eta$ and $\ga$ is a quasigeodesic in $X$, 
this contradicts the hypothesis that $\eta$ is a conical limit point of $Y$ and 
proves the proposition. \qed

\begin{rem}
This proposition is a geometric counterpart of the following group-theoretic result: 
If $H$ is a hyperbolic subgroup of a hyperbolic group $G$, and the CT-map $\D_{H,G}$ exists, and a limit point $z$ of $H$ in $G$ is conical, then $|\D_{H,G}^{-1}(z)|=1$. The converse to this implication is false, see \cite{MR3488025}. 
\end{rem}

The following conjecture is motivated by \cite{Kapovich-Liu}:

\begin{conj}
Suppose that $H<G$ is a hyperbolic subgroup of a hyperbolic group, the CT-map $\D_{H,G}$ exists, but $H$ is not quasiconvex in $G$. Then there is a continuum of (nonconical) limit points of $H$ in $G$ whose preimages under  $\D_{H,G}$ are not singletons. 
\end{conj}


\section{Group-theoretic applications} \label{sec:Group-theoretic applications}

In this section we collect group-theoretic applications of our existence results for CT-maps.

\subsection{Maps to products and examples of undistorted subgroups in $PSL(2, \CCC)\times PSL(2, \CCC)$}\label{sec:non-Anosov}

Set $H:= PSL(2, \CCC)$ and $G:=H\times H$. We equip $G$ with a left-invariant  Riemannian metric and the corresponding left-invariant distance function $d_G$. A finitely generated subgroup $\Gamma< G$ is said to be {\em undistorted} if the inclusion map
$$
(\Ga, d_\Ga)\to (G, d_G)
$$
is a qi embedding, where $d_\Ga$ is a word metric on $\Ga$. Since $G$ acts properly, isometrically and transitively on $\H^3\times \H^3$, a subgroup $\Ga< G$ is undistorted if and only if for (some/every) point $x\in \H^3\times \H^3$ the orbit map $\ga\mapsto \ga x$ is a qi embedding of $\Ga$ (with its word metric) into 
$\H^3\times \H^3$. 

An element $h\in H$ is called {\em parabolic} if it has precisely one fixed point in the Riemann sphere. 
An element $g=(h_1, h_2)\in G$ is called {\em semisimple} if neither component $h_1$ nor $h_2$ is a parabolic element of $H$. We will not attempt to define here Anosov subgroups of $G$, it suffices to say that each Anosov subgroup $\Ga< G$ is Gromov-hyperbolic and for 
one of the factors $H_\pm$ of $G=H\times H=H_+\times H_-$, the projection to $\Ga$ to $H_\pm$ has finite kernel and convex-cocompact image. 
Moreover, each Anosov subgroup is undistorted in $G$. We refer to the reader to \cite{MR2981818, MR3888689, MR3736790} for the detailed definitions. O.~Guichard constructed in \cite{Guichard-thesis} (see also \cite{MR3608719}) an example of an undistorted non-Anosov free subgroup $\Ga< G$. The subgroup in his example contained non-semisimple elements. (Its projections to both factors were geometrically finite with parabolic elements, we refer the reader to \cite{MR1317633} for definitions of geometric finiteness.) Let $S$ be a closed connected oriented hyperbolic surface with the fundamental group $\pi$.

\begin{thm}
There exists an undistorted subgroup $\Ga< G$  isomorphic to $\pi$, such that every element of $\Ga$ is semisimple, 
but $\Ga$ is not Anosov.  
\end{thm}
\proof Let $c$ be a complete geodesic in the Teichm\"uller space $T(S)$ of $S$, such that the projection of $c$ to the moduli  space of $S$ is bounded. For instance, $c$ can be taken to be the unique invariant geodesic ({\em axis}) in $T(S)$ of a pseudo-Anosov homeomorphism $h$ of $S$.  
The asymptotics of $c$ in positive/negative directions are described by two transversal geodesic laminations $\la^\pm$ on $S$, called {\em ending laminations}: Such laminations contain no closed geodesics and each component of $S\setminus \la^\pm$  is simply-connected (see \cite{klarreich-el}). In the example where $c$ is the axis of a 
pseudo-Anosov homeomorphism $h$, the laminations $\la^\pm$ are stable/unstable laminations of $h$ (see e.g. \cite{Casson}).  

Take $\la^\pm$, the ending laminations of a pseudo-Anosov homeomorphism of $S$ or, more generally, any two {\em transversal}  ending geodesic laminations on $S$. There exist discrete embeddings $\rho_\pm: \pi \to H$ such that the image $\Ga_\pm$ of each $\rho_\pm$ is a singly-degenerate subgroup of $H$ without parabolic elements such that $\la^\pm$ is the {\em ending lamination} of the geometrically infinite end\footnote{As it is customary in 3-dimensional topology we will be conflating ends and their neighborhoods.} $E^\pm$ of the hyperbolic manifold $M^\pm=\H^3/\Ga_{\pm}$. Furthermore, there exists a  
discrete embedding $\rho_0: \pi\to H$ such that the group $\Ga_0=\rho_0(\pi)$ is {\em doubly-degenerate group}, whose  quotient manifold $M_0=\H^3/\Ga_0$ has two ends $E_0^\pm$ with the ending laminations $\la^\pm$. We refer the reader to \cite{Ohshika-2009} for proofs of more general existence theorems of this type (which are generalizations of Thurston's double limit theorem). The ends $E^\pm$ are bilipschitz homeomorphic to the ends $E^\pm_0$ of the manifold $M_0$ (see \cite{minsky-jams} or \cite{minsky-elc2} for more general results).  The manifolds $M_0, M^\pm$ have injectivity radii bounded from below and $M^\pm$ has structure of a metric bundle over $\R$, whose fibers are uniformly bilipschitz to the surface $S$. For instance, in the case when $\la^\pm$ are stable/unstable laminations of a pseudo-Anosov homeomorphism $h$, the manifold $M_0$ is isometric to a cyclic covering space of the mapping torus of $h$, equipped with the unique hyperbolic metric.

We then obtain a discrete and faithful representation
$$
\rho: \pi\to \Ga< G= H\times H, \rho(\ga)=(\rho_+(\ga), \rho_-(\ga)). 
$$ 
By the construction, the image of each element of $\pi$ is a semisimple element of $G$. Since the projections $\Ga_\pm$ of $\Ga$ to the factors $H_\pm$ of 
$H\times H$ are geometrically infinite, the representation $\rho$ is not Anosov. It remains to prove that $\Ga$ is undistorted in $G$, i.e. that the map $\rho$ is a qi embedding. This qi embedding condition can be reformulated as follows. Consider the closed convex hulls $C^\pm\subset \H^3$ of the limit sets of the subgroups $\Ga_\pm< H$. The group $\pi$ acts properly discontinuously, isometrically and cocompactly on the boundaries of these convex hulls. Accordingly, we obtain quasiisometries
$$
\H^2\to \partial C^\pm,
$$
where the targets are equipped with intrinsic path-metrics. The quotient manifolds $C^\pm/\Ga_\pm$ are isometric to the ends $E^\pm\subset M^\pm$. Since 
$C^\pm$ are isometrically embedded in $\H^3$, $\Ga$ is qi embedded in $G$ if and only if the map
$$
f: \H^2\to \H^3\times \H^3
$$  
given by the composition of the isometries $f_\pm: \H^2\to \partial C^\pm$ with the inclusion maps $\partial C^\pm\to C^\pm\to \H^3$, is a qi embedding. Since the ends $E^\pm$ are  bilipschitz homeomorphic to the ends $E^\pm_0$ of the manifold $M_0$, $f$ is a qi embedding if and only if the following holds:

For some (every) $\Ga_0$-invariant embedded simply-connected hypersurface $\Sigma\subset \H^3$ separating $\H^3$ into components $\Sigma^\pm$ (equipped with the induced path-metrics), the inclusion maps $\Sigma\to \Sigma^\pm$ combine to a qi embedding 
$$
\Sigma\to \Sigma^- \times \Sigma^+. 
$$ 
Since $\H^3$ has a $\Ga$-invariant structure of a metric bundle $\X= (\pi: X\to T=\R)$ with fibers uniformly qi to $\H^2$, we just need to prove that for some (every) vertex 
$v\in T$ and the ideal boundary points $\xi_\pm$ of $T=\R$, the inclusion maps $X_v\to X_{T_\pm}=X_{v\xi_\pm}$ combine to a qi embedding
$$
\Phi: X_v\to X_{T_+}\times X_{T_-},
$$ 
where $T_\pm= v\xi_\pm$, a half-line. We will prove that $\Phi$ is indeed a qi embedding (and even more)  below, Proposition \ref{prop:KS}. 

\medskip 
Suppose $\X=(\pi: X\to T)$ is a tree of hyperbolic metric spaces with hyperbolic total space $X$ and let $v\in V(T)$ be a vertex of finite degree $n\ge 2$, with edges $e_1,...,e_n$ incident to $v$. For each $i=1,...,n$, let $T_i$ denote the subtree in $T$ which is the union of subintervals of the form $\llbracket v, w\rrbracket$, containing the edge $e_i$. 
We then obtain the subtrees of spaces $\X_{T_i}=(\pi: X_{T_i} \to T_i)$ in $\X$.  For each $i$, we let 
 $f_i: X_v\map X_{T_i}$ denote the inclusion map. We  equip the product 
$Q=\prod_{1\leq i\leq n} X_{T_i}$ with the $\ell_1$-metric
$$
d_Q(p, q)= \sum_{i=1}^n d_{X_{T_i}}(p_i, q_i),  \quad p=(p_1,...,p_n), q=(q_1,...,q_n). 
$$
(One can also use the $\ell_2$-metric, the {\em product} metric: The two metrics are qi to each other.)  

In what follows, we take $K=K_*$, $D= D_{\ref{prop:existence-of-tripod-ladders}}$, 
$E=E_{\ref{prop:existence-of-tripod-ladders}}$, depending on the parameters of the tree of spaces $\X$.

\medskip
The next proposition is a generalization of a result from \cite{Krishna-Sardar}, where it was proven in the case when $\X$  is a
 metric bundle: 

\begin{prop}\label{prop:KS}
Under the above assumptions, the diagonal map 
$$
\Phi: X_v\map Q=\prod_{1\leq i\leq n} X_{T_i}, \quad x\mapsto (f_i(x)),
$$ 
is a qi embedding.
\end{prop}
\proof We note that the inclusion maps $X_v\map X_{T_i}$ are all $1$-Lipschitz. Hence, the diagonal map 
$\Phi$ is  $n$-Lipschitz. The proof of the proposition is divided in two cases. 

\medskip 
{\bf Case 1:} Suppose that $n=2$. Consider a pair of points $x, y\in X_v$ and let $\L=\L(\al)=\{\pi:L\map \pi(L)\}$ be a $(K,D,E)$-ladder centered 
at $v$, with $\al=[xy]_{X_v}= L_v$. For   $i=1,2$, we have the $(K,D,E)$-ladders $\L^i=\{L^i\map \pi(L)\cap T_i\}$ in $\X_{T_i}$,  obtained by pull-back 
of the ladder $\L$ to the subtree of spaces $\X_{T_i}$. 

Let $c=c_L(x,y)$ be a combing path in $L$ connecting $x$ to $y$. We let $\hat{c}_1, \hat{c}_2$ be the paths in $L_1, L_2$ respectively, obtained from $c$ via the 
cut-and-replace procedure with respect to the inclusions $L_1\to L, L_2\to L$  (see Definition \ref{defn:detour}). Note that, since $L^i$ is a ladder in $X_{T_i}$, it 
is qi embedded in  $X_{T_i}$. Moreover, according to Theorem \ref{thm:cut-paste}, both $\hat{c}_1, \hat{c}_2$ are (uniform) $\kappa$-quasigeodesics in 
$X_{T_i}$. (Actually, this fact is established in Part I of the proof of Theorem \ref{thm:cut-paste}.) 

We will now estimate the length of  $\al$ from above in terms of the distance between $\Phi(x), \Phi(y)$ in $Q$. We claim that 
the segment $\al$ is contained in the union  $\hat{c}_1\cup \hat{c}_2$. Indeed, by the definition of combing paths $c=c_L$ in $L$, there exists a finite monotonic sequence 
$x_0=x, x_1,..., x_m=y$ in $\al$ such that $c$ is the concatenation of paths $c(x_i, x_{i+1})$ between points $x_i, x_{i+1}$, such that (after switching the roles of $L_1, L_2$ if necessary), 
$c(x_i, x_{i+1})$ is contained in $L_1$ for odd $i$ and is contained in $L_2$ for even $i$. Now, it follows from   the definition of 
the cut-and-replace procedure that $[x_i x_{i+1}]_{X_v}\subset \al$ is contained in $\hat{c}_1$ for each odd $i$ and is contained in $\hat{c}_2$ for even $i$.  

Thus, 
\begin{align*}
d_{X_v}(x,y)= \length(\al)\le \length(\hat{c}_1) +  \length(\hat{c}_2) \le \\
(\kappa +1) \left(d_{X_{T_1}}(x,y) + d_{X_{T_2}}(x,y)\right) = 2(\kappa +1)d(\Phi(x), \Phi(y)).  
\end{align*}
It follows that $\Phi$ is a qi embedding.

\medskip 
{\bf Case 2:} Suppose that $n\geq 3$. Consider two points $x, y\in X_v$. Observe that for $p=(x,x), q=(y,y)\in X_{T_1}\times X_{T_2}$, we have 
$$
d_Q((\underbrace{x,....,x}_{n \hbox{\ times}}), (\underbrace{y,...,y}_{n \hbox{\ times}}))\ge d_{X_{T_1}\times X_{T_2}}(p,q). 
$$
Therefore, Case 1 implies that the diagonal embedding $\Phi: X_v\to Q$ is a qi embedding. \qed

\medskip 
This concludes the proof of the theorem as well. \qed 

\begin{question}
{In the example given in this theorem, is the subgroup $\Ga< G$ a coarse Lipschitz retract of $G$?} 
\end{question}

Note that Anosov subgroups of semisimple Lie groups are coarse Lipschitz retracts, see \cite{bordif}.

 \subsection{CT-maps for hyperbolic graphs of groups}

In this section, ${\mathcal G}'$ is a finite graph of hyperbolic groups satisfying Axiom {\bf H}, with the underlying connected graph $\Ga'$ and the Bass--Serre tree $T'$. 
We will also assume that the group   $G'= \pi_1({\mathcal G}')$ is hyperbolic.

\medskip
We first prove the existence of CT-maps for some classes of hyperbolic subgroups $G< G'$. 

\begin{prop}\label{prop:CT-map-for-subgroups} 
Suppose that $\Ga\subset \Ga'$ is a connected subgraph and ${\mathcal G}\subset {\mathcal G}'$ is the subgraph of groups obtained by restricting 
${\mathcal G}'$ to $\Ga$ (see Section \ref{sec:generalities}), with $G= \pi_1( {\mathcal G})$. Then the subgroup $G< G'$ admits a CT-map $\geo G\to \geo G'$.  
\end{prop}
\proof Let $T, T'$ denote the Bass--Serre trees of  ${\mathcal G}$, ${\mathcal G}'$. The embedding of graphs of groups ${\mathcal G}\embed {\mathcal G}'$ induces 
 a   $G$-equivariant embedding $T\embed T'$. Since the subgraph of groups 
${\mathcal G}\subset {\mathcal G}'$ is obtained by the restriction, for each vertex $v$ and edge $e$ of the subtree $T'$, the stabilizer of $v$ (resp. $e$) in $G$ equals its stabilizer in $G'$. Thus, the tree of spaces $\X= (\pi: X\to T)$ corresponding to the graph of groups ${\mathcal G}$ is obtained as the pull-back of the tree of spaces 
$\X'=(\pi: X'\to T')$, $X=X'_{T}$. Since the groups $G, G'$ are naturally quasiisometric to the spaces $X, X'$ (via respective orbit maps) the existence of a CT-extension for the embedding $G'\to G$ is equivalent to that of the embedding $X'\to X$. Since the existence of a CT-map for the inclusion $X'\embed X$  is the content of Theorem \ref{thm:mainCT}, the proposition follows. \qed 

\medskip
The next theorem shows that one does not need to restrict to subgraphs of ${\mathcal G}'$ to obtain subgroups with CT-maps: 

\begin{thm}\label{thm:CT-maps-for-groups}
Assume  that $G< G'=\pi_1({\mathcal G}')$ is a subgroup preserving a subtree $T\subset T'$ such that the quotient graph $T/G$ is finite and 
 that  the vertex and edge stabilizers of this action on $T$ are quasiconvex in the respective subgroups of $G'$: $G_v< G'_v$ and $G_e< G'_e$ are quasiconvex for all $v\in V(T), e\in E(T)$.  Then the subgroup $G$ is hyperbolic and the inclusion map $G\to G'$ admits a CT-map.  
\end{thm}
\proof We will use Proposition \ref{prop:retract-so-subgroup}: As in the proof of the proposition we observe that the $G$-action on $T$ defines a graph-of-groups decomposition of $G$: 
$\pi_1({\mathcal G})=G$, and the graph of groups ${\mathcal G}$ satisfies Axiom {\bf H} (in view of the quasiconvexity assumptions in the theorem). We let 
$\X= (\pi: X\to T)$ and $\X'= (\pi: X'\to T')$ denote the trees of spaces corresponding to the graphs of groups ${\mathcal G}, {\mathcal G}'$ respectively. 

Since $G'$ is hyperbolic, so is $X'$ and, hence, $\X'$ satisfies the proper flaring condition. We have a $G$-equivariant 
relatively retractive morphism of trees of spaces $\X\to \X'$ ($h: X\to X'$, over the inclusion $T\to T'$). The proper flaring condition for $\X'$ then implies the proper flaring condition for $\X$, 
hence, $X$ and, thus, $G$,  is also hyperbolic.  We let $\Y= (\pi: X'_T\to T)$ denote the restriction of the tree of spaces $\X'$ to $T$. The quasiconvexity assumption for the subgroups 
$G_v< G'_v, G_e< G'_e$, $v\in V(T), e\in E(T)$, implies that for each $v\in V(T)$, and edge $e=[v,w]\in E(T)$, 
the  $G_v$-orbit of $X'_{ev}$ is {\em locally finite} in $X'_v$, see Lemma \ref{lem:geo->lf}. Thus, Proposition \ref{prop:retract-so-subgroup} implies that the map  
$h: X\to Y=X'_T$ is a qi embedding. According to Theorem \ref{thm:mainCT}, the inclusion $Y\to X'$ admits a CT-map $\D_{Y,X'}$. Composing it with the boundary map of the qi embedding 
$X\to Y$, we obtain a CT-map $\geo h$ for the map $h: X\to X'$. Since $G$ acts geometrically on $X$ and $G'$ acts geometrically on $X'$ (see Section \ref{sec:Group actions} for the definition and Lemma \ref{lem:Milnor--Schwarz Lemma}), we conclude from the existence of $\D_{Y,X'}$ the existence of a CT-map for the subgroup $G< G'$. \qed

\begin{example}
Let $G'= F\star_{\varphi}$ be a hyperbolic group which is the descending HNN extension of a finitely generated free group $F$ via an injective endomorphism $\varphi: F\to F$. Then $G= F\star_{\varphi^n}$ is a hyperbolic subgroup of $G'$ and  
the embedding $G\to G'$ admits a CT-map. 
\end{example}

For a boundary vertex $v$ of a subtree $T\subset T'$, we let  $T'(v)\subset T'$ denote  the maximal subtree of $T'$ containing $v$ and disjoint from the rest of the vertices of $T$.  
Thus, if $g$ is an automorphism of $T$ fixing $v$ and preserving $T$, it preserves the subtree $T'(v)$ as well.

\begin{thm}
Assume that $G< G'$ are as in Theorem \ref{thm:CT-maps-for-groups} and that for each boundary vertex $v$ of $T$ in $T'$,  the stabilizer $G_v<G$ acts $k$-acylindrically on 
the subtree $T'(v)\subset T'$. Then $G$ is a quasiconvex subgroup of $G'$. 
\end{thm}
\proof According to Remark \ref{rem:CT-remark}(2), in order to prove the quasiconvexity of $G$ in $G'$ it suffices to show that the CT-map $\D_{G,G'}$ is injective, i.e. that the CT-lamination $\La(G,G')$ is empty. For the sake of contradiction, suppose that $\La(G,G')\ne \emptyset$.  As we observed in the proof of Theorem \ref{thm:CT-maps-for-groups}, 
 the action of $G$ on the space $Y=X'_T$ is quasiconvex.  Since $\La(G,G')\ne \emptyset$,  follows that there exists a pair of distinct limit points $z_\pm$ of $G$ in $\geo Y$ with equal images under $\D_{Y,X'}$. By Proposition \ref{prop:CT-fibers1}, there is a biinfinite vertical geodesic $\al\subset X'_v$ (for some $v\in V(T)$) which is a quasigeodesic in $Y$, such that  $z_\pm =\al(\pm \infty)$. 

\begin{lem}
Suppose that $\rho$ is a geodesic ray in a vertex-space $X'_v, v\in T$,  which is also a quasigeodesic ray in $Y$, such that $z=\rho(\infty)$ is a limit point of $G$ in $\geo Y$. Then 
$z$ is a (conical) limit point of the action of $G_v$ on $X'_v$.  
\end{lem}
\proof Since the $G$-action on $Y$ is quasiconvex, the limit point $z$ is a conical limit point (see Definition \ref{defn:conical-limit}). Thus, there is a sequence $g_i\in G$ and a constant $r$ 
such that for $x=\rho(0)$, $d(g_i x, \rho)\le r$, and $\lim_{i\to\infty} g_i(x)=z$. At the same time, the $G$-orbit of $X'_v$ in $Y$ is locally finite (since for $g\in G$, $gX'_v= X'_{gv}$ and 
each compact in $X$ intersects only finitely many vertex-spaces). Since $g_i(x)\in N_r(X'_v)$,  Proposition \ref{prop:proj-to-qc-action} implies that $g_i(x)\in N_R(G_v x)$ for some $R$ 
independent of $i$. Hence, for $h_i\in G_v$ such that $d(g_i(x), h_i(x))\le R$, we obtain $\lim_{i\to\infty} h_i(x)=z$ in $X'_v$. \qed

\medskip
We now return to the proof of the theorem. By the lemma, the points $z_\pm$ are limit points of the $G_v$-action on $X'_v$. 
Since $\D_{Y,X'}(z_-)= \D_{Y,X'}(z_+)$ and $z_\pm\in \geo X'_v\subset \geo Y$ are  limit points of the $G_v$-action on $X'_v$,  
Theorem \ref{thm:CT lamination}(1) implies that there is a point $\xi\in \geo T' \setminus \geo T$ 
such that $\{z_-, z_+\}\in \La^\xi(X'_v,X')$. By Proposition \ref{flow of CT leaf}(2), for each vertex $w\in V(v\xi)$, 
$Fl_w(\{z_\pm\})\ne \emptyset$. Thus, according to Lemma \ref{lem: flow condition}, for each $x\in X_{v}$ and vertex $w\in v\xi$, 
the pair of subsets $G_{v}x, X'_w$  is not cobounded in $X_{vw}$.  Lemma \ref{lem:trivial stabs} then implies that the $G_v$-stabilizer of the interval $J_w=\llbracket v, w\rrbracket$ is infinite 
for each $w\in V(v\xi)$. Since  $\xi\in \geo T' \setminus \geo T$, the intersection of  the ray $v\xi$ with the subtree $T$ 
is a finite interval $\llbracket v, v'\rrbracket$. The vertex $v'$ is a boundary vertex of $T$ in $T'$. 
Let $G_{J_w}$ denote the $G$-stabilizer of the interval $J$. Consider a vertex $w\in v\xi$ such that $v'\in V(\llbracket v, w\rrbracket)$ and, moreover,
$d_T(v', w) > k$. Then the infinite subgroup $G_J\cap G_v$ fixes the interval $\llbracket v', w\rrbracket\subset  T'(v')$ of length $>k$, contradicting 
 the hypothesis that the group $G_{v'}$ acts $k$-acylindrically on the subtree $T'(v')$.  \qed

\begin{rem}
Assume that $G< G'$ are as in Theorem \ref{thm:CT-maps-for-groups}.  

1. The subgroup $G$ is at most exponentially distorted in $G'$ since $Y$ is at most exponentially distorted in $X'$ (Corollary \ref{cor:exp-dist}) and the orbit map $o_y: G\to Gy\subset Y$ is a qi embedding. 

2. In the setting of  Theorem \ref{thm:CT-maps-for-groups}, we can drop the assumption of hyperbolicity for $G'$, but assume that for each boundary vertex $v\in T$ and the boundary edge $e=[v,w]$ the stabilizer $G_e$ is finite. Then  $G$ is a coarse Lipschitz retract of $G'$, since for each boundary edge $e$ the projection of $G_v$ to $G'_e$ is uniformly bounded, cf. 
Theorem \ref{thm:left-inverse} and Proposition \ref{prop:retract-so-subgroup}.

3. In \cite{MR3143716} Baker and Riley construct examples of finitely generated free subgroups $G$ in certain  hyperbolic groups $G'$ such that the CT-maps for 
the inclusions $G\to G'$ do not exist. Their groups $G'$ free-by-cyclic, hence, are isomorphic to fundamental groups of graphs of groups ${\mathcal G}'$ satisfying Axiom {\bf H}. However, in their examples, the intersections of $G$ with vertex/edge subgroups of ${\mathcal G}'$ are not finitely generated. 
\end{rem}

The next result is a direct group-theoretic application  of Proposition \ref{prop:CT-fibers1} regarding the nature of CT-laminations for subgraphs of groups:

\begin{cor}\label{cor:CT-lamination-gr-of-groups}
Suppose that $G< G'$ are as in Proposition \ref{prop:CT-map-for-subgroups} and $\xi^\pm\in \geo G$ are distinct points which have the same image in $\geo G'$ under the  
CT-map $\geo G\to \geo G'$. Then there exists a vertex $v\in T$ and a pair of points $\xi^\pm_v\in \geo G_v$ such that:

1. $\D_{G_v,G}(\xi_v^\pm)=\xi^\pm$.  

2. Geodesics in $G_v$ connecting $\xi_v^\pm$ are uniform quasigeodesics in $G$. 
\end{cor}

Here is another application, this time of Theorem \ref{thm:CT-fibers}. Consider a hyperbolic group $H$, an automorphism $f: H\to H$ and the semidirect product $G=H\rtimes_f \Z$. The next result describes when a geodesic in $H$ is a leaf of the CT-lamination $\La(H,G)$. For the formulation of the result we will use the notion of {\em pseudo-orbits} of the automorphism $f$, Definition \ref{defn:pseudo-orbit}. 

\begin{cor}
A geodesic $\al$ in the Cayley graph of $H$ is a leaf of $\La(H,G)$ if and only if the following holds for some numbers $K\ge 1$ and $C\ge 0$:

There exist $K$-pseudo-orbits $(y^\pm_i)$, of $h_{\pm n}=\al(\pm n)\in H$ under the automorphism $f$ which approach each   other within distance $C$. More precisely, there exist 
$i=i(n)$ such that 
$$
d_H(y^+_i, y^-_i)\le C. 
$$
\end{cor}
\proof We let $\X=(\pi: X\to \T=\R)$ denote the tree of spaces corresponding to the graph-of-groups structure on $G$ given by the 
HNN-extension of $H$ via the automorphism $f$. 

Suppose that $\A^n=(\pi: A^n\to \llbracket 0, w_n\rrbracket)$, $w_n=i=i(n)$, is a $(K,C)$-narrow carpet in $G$ bounded by 
$[h_{-n} h_{n}]_{X_0}$ and $\ga_n^\pm$ denote the $K$-qi sections corresponding to the top/bottom of the carpet $\A^n$. 
By the definition of a $(K,C)$-narrow carpet,   
$$
d_{t^{i}H}(\ga^-_n(i), \ga_n^+(i))\le C. 
$$   
As it was explained in Section \ref{sec:semi-direct-flaring}, for each $n$, the sequences $y^\pm_j=\ga^\pm_n(j)$ are 
precisely the partial pseudo-orbits of $f$ in $H$ through the points $h_{\pm n}$. Thus, the claim is a direct consequence of 
Theorem \ref{thm:CT-fibers}. \qed

\begin{rem}
One can show that $\al$ as above is a leaf of $\La(H,G)$ if and only if there are sequences $h_{\pm n}$ in $H$ converging to $\al(\pm \infty)$ (but not necessarily of the form $\al(\pm n)$), such that the $f$-orbits of $h_{\pm n}$ 
approach each   other within distance $C$, where $C$ is a uniform constant depending only on the group $H$, its generating set, and $f$. 
\end{rem}

\subsection{Miscellaneous results}

The next proposition is a partial converse to Proposition \ref{conical-singleton}:  

\begin{prop}
Suppose $\GG$ is a finite graph of hyperbolic groups satisfying Axiom {\bf H} and $G=\pi_1(\GG)$ is
hyperbolic. Suppose $X\map T$ is the tree of spaces associated to this graph of groups. If $v$ is a vertex of $T$ and $z\in \geo X_v$ is such that the subtree 
$T_z= \pi(Fl(\{z\}))$ contains no geodesic rays, then $\D_{X_v,X}(z)$ is a conical limit point of $X_v$ in $X$.
\end{prop}
\proof By the definition, $T_z\subset T$ is the subtree whose vertex set consists of those vertices $w\in V(T)$ for which $Fl_w(z)\neq \emptyset$.

\begin{lemma}
The subtree $T_z$ is finite.  
\end{lemma}
\proof Since $T_z$ contains no rays, it suffices to prove that the tree $T_z$ is locally finite. Consider a vertex $w\in T_z$ and the collection of edges $e_i, i\in I,$ 
in $T_z$ incident to $w$. Let $Fl_w(\{z\})=\{z'\}$, $z'\in \geo X_w$. Then, according to the definition of the boundary flow in Section \ref{sec:Boundary flows}, for each edge $e_i$ we have $Fl_{e_i}(\{z'\})\ne \emptyset$.  
It follows that $z'\in \geo (G_{e_i}, G_w)$ for each $i\in I$. Since $\GG$ is a finite graph of groups, there are only finitely many $G_w$-conjugacy classes of edge-stabilizers  $G_{e_i}< G_w$. At the same time, 
$$
z'\in \bigcap_{i\in I}  \geo (G_{e_i}, G_w) \ne \emptyset,  
$$ 
hence (since each subgroup $G_{e_i}$ is quasiconvex in $G_w$) each intersection $G_{e_i}\cap G_{e_j}$ is an infinite subgroup of $G_w$, 
see e.g. Lemma 2.6 in \cite{MR1389776}. The main theorem in  \cite{MR1389776} states that {\em quasiconvex subgroups of hyperbolic groups have finite width}. Without defining {\em width of subgroups} here, we only note that, as a consequence of this finiteness theorem, if $H_i, i\in I$, is a collection of pairwise distinct quasiconvex subgroups of a hyperbolic group $H$ which belong to finitely many $H$-conjugacy classes and $|H_i\cap H_j|=\infty$ for all $i, j\in I$, then $I$ is finite. Applying this result in our setting, with the ambient hyperbolic group $H$ equal $G_w$ and quasiconvex subgroups $H_i$ equal to the edge-subgroups $G_{e_i}$,  we conclude that   
 the set $I$ is finite. Thus, $T_z$ is a locally finite tree and, hence, is finite. \qed

\medskip
We can now prove the proposition. In the proof it will be convenient to assume that each edge-space $X_e$ of $\X$ is discrete, cf. introduction to Section \ref{sec:CTfibers}. 
Let $\beta$ be a ray in $X$ asymptotic to the point $\D_{X_v,X}(z)$. Suppose for a moment that the intersection $\beta\cap X_{T_{z}}$ 
is bounded. Then there exists a boundary edge $e=[v',w]$ of $T_z$, $w\notin V(T_z), v'\in V(T_z)$, 
such that an unbounded subray $\beta'$ of $\beta$ projects to the subtree $T_{w,z}\subset T$ which is the maximal subtree of $T$ containing $w$ and disjoint from $T_z$. 
(Here are are using the discreteness assumption on the edge-spaces of $\X$.) Let $K$ be as in  
Proposition \ref{flow of CT leaf} and let $\al\subset X_v$ be a geodesic ray  asymptotic to $z$. 

Since $w\notin T_z= \pi(Fl(\{z\}))$, the intersection $Fl_K(\al)\cap X_w$ is bounded, cf. Lemma 
\ref{lem: flow condition}. Recall that the flow-space $Fl_K(\al)$ is quasiconvex in $X$. Hence, we will identify $\geo Fl_K(\al)$ with a subset of $\geo X$. 
We claim that 
$$
\D_{\al, Fl_K(\al)}(z)\ne \beta(\infty) \in \geo Fl_K(\al). 
$$
Indeed, the assumption that $\pi(\beta')\subset T_{w,z}$ implies that each geodesic in  $Fl_K(\al)$ connecting points of $\al$ to that of $\beta'$ has to pass through 
$Fl_K(\al)\cap X_w$, i.e. within distance $D$ from $p=\al(0)$, where $D=\Hd(\{p\}, Fl_K(\al)\cap X_w)$. 
But this means that the sequences $(\al(n)), (\beta(n))$ cannot define the same ideal boundary point of $\geo Fl_K(\al)$. Quasiconvexity of $Fl_K(\al)\subset X$ implies that 
$\D_{X_v,X}(z)\ne \beta(\infty)$, which is a contradiction.  Thus, $\beta$ contains an unbounded sequence of points $(x_n)$ contained in $Y:=Fl_K(\al)\cap X_{T_z}$. Since (by the lemma) 
the subtree $T_z$ is finite, the subset $Y$ is Hausdorff-close to the ray $\al\subset X_v$. In other words, the point  
$\D_{X_v,X}(z)= \beta(\infty)$ is a conical limit point of $X_v$. \qed

\medskip 
We conclude the chapter with a proposition that deals with the case of {\em nonhyperbolic} graphs of hyperbolic groups and relates this lack of hyperbolicity to various notions discussed earlier, such as boundary flow-spaces and unbounded sequences of carpets:  

\begin{prop}
Suppose $\mathcal G$ is a finite graph of hyperbolic groups satisfying Axiom {\bf H}; let $\X=  (\pi: X\map T)$ denote the corresponding tree of metric spaces. 
If $G=\pi_1(\mathcal G)$ is not hyperbolic, then there is a vertex $v\in V(T)$ 
such that $\pi(Fl(\partial_{\infty} X_v))$ contains a geodesic ray. 
\end{prop}
\proof Let $G$ be the fundamental group of $\mathcal G$.
By Corollary \ref{cor:super-weak flaring}(2) there is a constant $D>0$, a sequence of intervals $I_n=\llbracket -t_n,t_n\rrbracket\subset T$
of length at least $n$  and a sequence $(\Pi^n)$ of pairs of $\kappa$-qi sections $(\gamma^n_0, \gamma^n_1)$ over
$I_n$ ($n\in \NN$), such that $\Pi^n_{max}\leq D$, but $\Pi^n_0\map \infty$. Since the $G$-action is cofinite on $T$ 
and cocompact on $X$, and the map $\pi: X\map T$ is $G$-equivariant, after extraction, 
we may assume that for each $n$ the midpoint of $I_n$ is a fixed vertex $v\in V(T)$ and the midpoint of $[\gamma^n_0(0) \gamma^n_1(0)]_{X_v}$ is within unit distance from a fixed 
point $x\in X_v$.
After passing  to a further subsequence, we may assume that the sequence of segments 
$[\gamma^n_0(0) \gamma^n_1(0)]_{X_v}$ converges to a complete 
geodesic in $X_v$. We note that for each $n$ the vertical geodesic segments $[\gamma^n_0(t)  \gamma^n_1(t)]_{X_t}, t\in V(I_n)$, form a
$K$-metric bundle over $I_n$ (with $K=\kappa'$, see Lemma \ref{lem:E-ladder-structure}) 
and hence there is a $K$-qi section  
 $\gamma^n$ over $I_n$ passing through $B(x,2)$ and contained in this metric bundle.
Hence, after passing to a further subsequence, we may assume that the sequence of sections 
$(\gamma^n)$ converges to a 
complete quasigeodesic $\gamma$ in $X$ and the sequence $I_n$ converges of complete geodesic $I$ in $T$ such that $\gamma$
is a qi section over $I$, cf. the proof of Lemma \ref{lemma1: CT lam}(2). Lastly, 
Lemma \ref{lem: flow condition}(3) implies that $I\subset \pi(Fl(\partial_{\infty} X_v))$.
\qed

\chapter{Cannon-Thurston maps for relatively hyperbolic spaces} \label{ch:CTR}

The goal of this chapter is to generalize the results on Chapter \ref{ch:CT} (primarily, the existence of CT-maps and some basic facts about CT-laminations)  
in the context of {\em relatively hyperbolic spaces}. Such a 
generalization was achieved in \cite{mahan-pal} for 
the inclusion maps of vertex-spaces.

\section{Relative hyperbolicity}

\subsection{Relative hyperbolicity in the sense of Gromov} \label{sec:GRH} 

We refer the reader to \cite{farb-relhyp, bowditch-relhyp} for the background on the  theory of relatively hyperbolic spaces. Briefly, a {\em relatively hyperbolic space} is a pair $(Y, {\mathcal H})$ consisting of a geodesic metric space $(Y, d)$ together with a collection ${\mathcal H}=\{H_i: i\in I\}$ of {\em peripheral subspaces}, \index{peripheral subspace} 
which are nonempty subsets of $Y$ 
satisfying certain conditions discussed below and, in an alternative form, in Section \ref{sec:Farb->Gromov}. While this is not always required for relatively hyperbolic spaces, we will assume that each $H_i\in {\mathcal H}$ is rectifiably connected and the inclusion maps $(H_i, d_{H_i})\to (Y, d_Y)$ are uniformly proper, where $d_{H_i}$ are the intrinsic path-metrics on $H_i$'s. Given such a pair $(Y, {\mathcal H})$, one defines two new metric spaces:

\medskip

1. The {\em extended hyperbolic space} $(Y^h, d^h)=(Y^h, d_{Y^h})= {\mathcal G}(Y, {\mathcal H})$ (or the {\em horoballification} of $(Y, {\mathcal H})$), which is a path-metric space 
obtained by attaching along each $H_i$ its {\em hyperbolic cone} $H_i^h$ defined in Section \ref{sec:hyperbolic cone};  the latter are  (intrinsically) {\em uniformly hyperbolic horoballs}, with $H_i$ the boundary horosphere in $H_i^h$. Recall that each $H_i^h$ has unique ideal boundary point, called the {\em foot-point} (or the {\em ideal center}) of $H_i^h$, $\xi(H_i^h)$.

2. The {\em electric space}  
$(Y^{\ell}, d^\ell)={\mathcal E}(Y, {\mathcal H})$ 
(the {\em electrification} of $(Y, {\mathcal H})$), obtained by coning off each $H_i$, i.e. attaching to $Y$ along each $H_i$ the cone $H_i^{\ell}=  C(a_i, H_i)$ with the apex $a_i=a(H_i^\ell)$ 
within distance $1/2$ from each point in $H_i$. We refer the reader to Section \ref{sec:metric cones} for the precise definition of the metric spaces $C(a_i, H_i)$. Here we recall only that each point $x\in H_i$ is connected to $a_i$ 
by a canonical geodesic segment of length $1/2$, called  a {\em radial line-segment}. 
The set of apexes of these cones is called the {\em cone-locus} of $Y^\ell$ and denoted $a(Y^\ell)$. \index{cone-locus}

We will use the notation $\mathring{H}_i^h$ (the {\em open peripheral horoball}) and $\mathring{H}_i^\ell$ (the {\em open peripheral cone}) for the complements
$$
H_i^h \setminus H_i, \quad H_i^\ell \setminus H_i
$$
 respectively. \index{peripheral horoball}

\begin{defn}
[Relative hyperbolicity in the sense of Gromov] \index{relatively hyperbolic in the sense of Gromov}
A pair $(Y, {\mathcal H})$ is called {\em relatively hyperbolic in the sense of Gromov (GRH)} if the metric space $(Y^h, d^h)$ is hyperbolic. 
\end{defn}

\medskip 
If $(Y, {\mathcal H})$ is GRH, then the electric space $(Y^\ell, d^\ell)$ is hyperbolic. This 
is a standard fact, usually attributed to Farb, \cite{farb-relhyp}: Although the  proofs given in his paper are only in the setting of manifolds of negative curvature, they go through in greater generality. We will give a proof in Section \ref{sec:el-space} with the side-benefit of relating quasigeodesics in $(Y^h, d^h)$ to those in $(Y^\ell, d^\ell)$. 

\subsection{Extrinsic geometry of the peripheral horoballs} 

Throughout this section we will assume that $(Y, {\mathcal H})$ is GRH with the hyperbolicity constant $\delta$. Our goal is to discuss the {\em extrinsic} geometry of the peripheral horoballs $H_i^h$. Among other things, we will prove that they are uniformly quasiconvex and uniformly pairwise cobounded.

\begin{lem}\label{lem:qc horoballs} 
1. The subsets $H_i^h$  are $\la_{\ref{lem:qc horoballs}}(\delta)$-quasiconvex and $L_{\ref{lem:qc horoballs}}(\delta)$-qi embedded in $Y$.

2.  $Y$ is uniformly properly embedded in $Y^h$ with the distortion function depending only on $\delta$. 
\end{lem}
\proof 1. We observe that, by the definition of the metric on $H_i^h$, for each point $x=(z,t)\in H_i^h$ the distance from $x$ to $H_i$ equals $\log(t)$. It follows that for any two points $x_1=(z,t_1), x_2=(z,t_2) \in H_i^h$,
$$
d_{Y^h}(x_1, x_2)= |\log(t_1/t_2)|
$$
and the vertical segment in $H_i^h$ between $x_1, x_2$ is isometrically embedded in $Y^h$. In particular, the vertical rays 
in $H_i^h$ are isometrically embedded in $Y^h$. Since any two such rays $\rho_1(t), \rho_2(t)$ converge as $t\to\infty$, it follows 
that the horoballs $H_i^h$ 
are uniformly quasiconvex in $Y^h$. The fact that the peripheral subspaces $H_i^h$ are uniformly properly embedded in $Y$ implies that the horoballs $H_i^h$ are uniformly properly embedded in $Y^h$. Combined with the uniform quasiconvexity, we obtain that these horoballs are uniformly qi embedded in $Y^h$. 

2. Lastly, uniform properness of  the inclusion maps $H_i\to H_i^h$  (see Proposition 
\ref{prop:exp-distortion}) implies that $Y$ is uniformly properly embedded in $Y^h$ as well.
\qed 

\medskip
Recall that a closed subset $C$ of a geodesic metric space $X$ is called {\em strictly convex} if for any two points $x, y\in C$, every geodesic $xy\subset X$ is contained in the interior of $C$, except maybe for its end-points. For instance, closed balls and horoballs  in the classical hyperbolic space are strictly convex, while a closed hyperbolic half-space is not.  
The next lemma establishes a form of  coarse strict convexity of peripheral horoballs in the context of GRH spaces. 

\begin{lemma}\label{lem:horodiving} 
There exist $L=L_{\ref{lem:horodiving}}(K,r,\delta)$ and $R=R_{\ref{lem:horodiving}}(K,r,\delta)$ 
satisfying the following properties. Suppose that $\beta: [0,T]\to Y^h$ is a continuous $K$-quasigeodesic in $Y^h$ connecting points $z=\beta(0), y=\beta(T)$.  

1. Suppose that the image of $\beta$ is entirely contained in $N_r(H)$ for some $H \in {\mathcal H}$. 
Then $d(y,z)\le R=R_{\ref{lem:horodiving}}(K,r,\delta)$. 

2. Suppose that the points $z=\beta(0), y=\beta(T)$ both belong to $N_r(H^h)\setminus \mathring{H}^h$. 
Then either $T\le 2L$ or there exist $a\in [0, L]$ and $b\in [T-L,T]$ such that $z'=\beta(a)\in H$, $y'=\beta(b)\in H$ and the subpath $\beta(z', y')$ in $\beta(z,y)$ is entirely contained in $\mathring{H}^h$, except for the end-points $z', y'$.  
\end{lemma}
\proof 1. Let $\bar{y}, \bar{z}\in H^h$ be the images of $y, z$, respectively, under the projection $P_{Y^h,H^h}$. Thus, 
$d(y,\bar{y})\le r, d(z, \bar{z})\le r$. Let $c=c(\bar{z}, \bar{y})$ denote the combing path in $H^h$ connecting $\bar{z}$ to $\bar{y}$,  as defined in the proof of Proposition \ref{prop:hyper-cone}. By the version of stability of quasigeodesics in hyperbolic spaces 
(Lemma \ref{lem:sub-close}), the Hausdorff distance between the images of $\beta$ and $c$ is at most   
$D_{\ref{lem:sub-close}}(\delta, \max(k,K), r)$, where $k$ is the bound on qi constant of $c$ given in Remark \ref{rem:qg-in-horoball}. 

We will consider the case 
when $\bar{x}, \bar{y}$ are both in $H$ and leave the other cases to the reader as the proofs are similar.  
 According to the description of the combing paths in $H^h$, for each constant $D\ge 0$ if 
 the path  $c$ is contained in 
 the $D$-neighborhood of $H$, then $A=d_H(\bar{z}, \bar{y})\le e^D$. 
 
 
 In our case, $c$ is contained in the $D=(D_{\ref{lem:sub-close}}(\delta, \max(k,K), r)+ r)$-neighborhood of $H$ and, hence, we conclude that $d_H(\bar{z}, \bar{y})\le e^D$. It follows that 
 $$
 d(y,z)\le R_{\ref{lem:horodiving}}(K,r,\delta):= e^D + 2r. 
 $$
 
 2. We follow the arguments of Part 1, define points $\bar{y}, \bar{z}$ and the path $c$ connecting them, setting now 
 $D:=D_{\ref{lem:sub-close}}(\delta, \max(k,K), r)$, an upper bound on the Hausdorff distance between $\beta$ and $c$. 
 If  $d_H(\bar{z}, \bar{y})> e^D$, then there are points $\bar{y}', \bar{z}'$ within distance $D$ from $\bar{y}, \bar{z}$ respectively, such that the subpath $c(\bar{z}, \bar{y})$ is disjoint from the $D$-neighborhood of $H$, apart from the end-points. 
 Thus, there are points $y'', z''$ in the image of $\beta$ which are within distance $D$ from, respectively,  
 $\bar{y}', \bar{z}'$, such that the subpath $\beta(z'', y'')$ is entirely contained in $\mathring{H}^h$, except, possibly, 
 the end-points. We then take the points $y', z'$ on the subpaths $\beta(y,y''), \beta(z,z'')$ where these paths cross into $H^h$ and such that $\beta(z', y')$ is entirely contained in $\mathring{H}^h$, except for the end-points $x', y'$.  Lemma follows. \qed

\begin{cor}\label{lem:strict-convexity}
Suppose that $y, z\in N_r(H^h)$ and the geodesic $[yz]_{Y^h}$ is disjoint from $\mathring{H}^h$. Then $d(y,z)\le 2L_{\ref{lem:horodiving}}(1,r,\delta)$. 
\end{cor}

\begin{cor}\label{cor:deep-pieces}
Suppose that $H_i, H_j$ are distinct elements of ${\mathcal H}$, $\beta: [0,T]\to Y$ is a $K$-quasigeodesic, $[s_i, t_i], [s_j, t_j]$ are subintervals of length $> 2L_{\ref{lem:horodiving}}(K,r,\delta)$ in $[0,T]$ such that $\beta(s_i)\in H_i, \beta(t_i)\in H_i, \beta(s_j)\in H_j, \beta(t_j)\in H_j$. Then $[s_i, t_i]\cap [s_j, t_j]=\emptyset$. 
\end{cor}

   \begin{figure}[tbh]
\centering
\includegraphics[width=70mm]{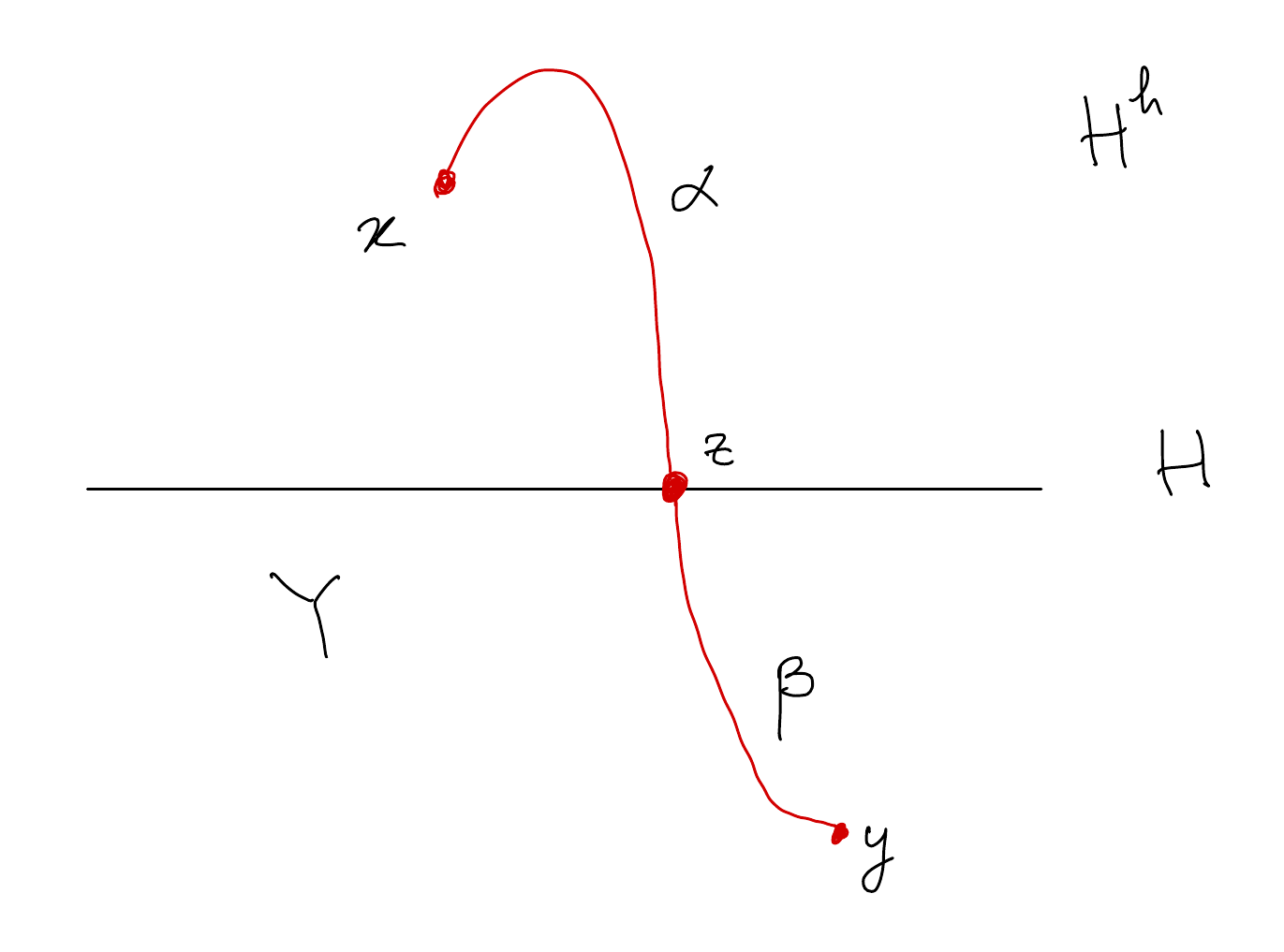}  
\caption{}
\label{c1.fig}
\end{figure}

\begin{lem}\label{lem:two geodesics}
Consider three points, $x\in H^h, z\in H$ and $y\in Y^h$, and $K$-quasigeodesics $\alpha=\al(z,x), \beta=\beta(y,z)$ in $Y^h$ connecting $z$ to $x$ and $y$ to $z$ respectively, such that $\beta$ is disjoint from $\mathring{H}^h$ and $\al$ is contained in $H^h$. (See Figure \ref{c1.fig}.) 
Then the concatenation $\beta\star \al$ is an 
$L_{\ref{lem:two geodesics}}(K,\delta)$-quasigeodesic in $Y^h$.  
\end{lem}
\proof We set $D:=D_{\ref{stab-qg}}(\delta, K)$, $R:=R_{\ref{lem:horodiving}}(1,D,\delta)+1$. Note that 
$D> 2\delta$, see Lemma \ref{stab-qg}.  

 In view of stability of quasigeodesics, it suffices to get an upper bound on the qi constant of the concatenation 
 $\beta^* \star \al^*$,  where  $\al^*=[zx]_{Y^h}$, $\beta^*=[yz]_{Y^h}$ (geodesics which are $D$-Hausdorff close to $\al, \beta$, respectively, where $D=D_{\ref{stab-qg}}(\delta, K)$). By the assumption of the corollary,  $\beta^*$ is contained in $N_D(Y^h\setminus H^h)$. There are two cases to consider: 

a. $d_{Y^h}(y,z)< R$ or $d_{Y^h}(x,z)< R$. Then the claim follows immediately: A concatenation of a uniform quasigeodesic with a uniformly bounded quasigeodesic is again a uniform quasigeodesic. 

b. $R\le \min(d(x,z), d(y,z))$.  Then, by Lemma \ref{lem:horodiving}(1), there exists a point $z'\in \beta^*$ at the distance $R'\le R$ from $z$ such that $z'\notin N_D(H^h)$. Take the point $x'\in \al^*\subset H^h$ at the same distance $R'$ from $z$. Then 
$$
d(z', x')\ge D> 2\delta. 
$$
By Lemma \ref{lem:g-concat}, the concatenation $\beta^*\star \al^*$ is then 
an  $L_{\ref{lem:g-concat}}(R',\delta)$-quasigeodesic in $Y^h$. \qed 

\medskip
In the next lemma we prove a generalization of this result, for concatenations of three uniform quasigeodesics.

   \begin{figure}[tbh]
\centering
\includegraphics[width=70mm]{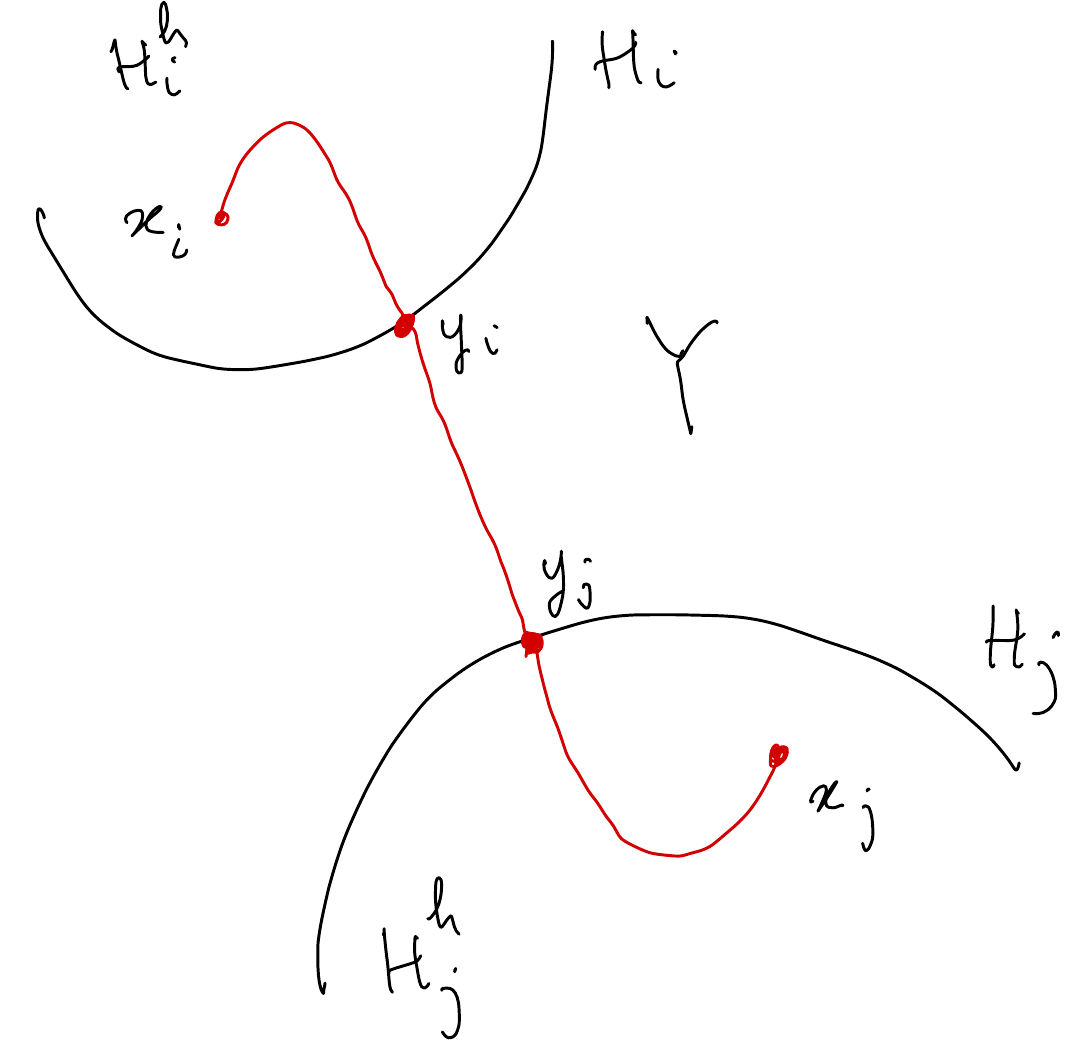}  
\caption{}
\label{c2.fig}
\end{figure}

\begin{lem}\label{lem:three geodesics}
Suppose we have three $K$-quasigeodesics $\beta(x_i,y_i)$, $\beta(y_i,y_j)$ and $\beta(y_j, x_j)$ in $Y^h$ such that the paths  
$\beta(x_k,y_k)$ are contained in the distinct peripheral horoballs $H_k^h$, and 
$\beta(y_i,y_j)\cap \mathring{H}_k^h=\emptyset$, $k= i, j$. (See Figure \ref{c2.fig}.)  Then  
 the concatenation 
 $$
 \beta(x_i,y_i)\star \beta(y_i,y_j)\star \beta(y_j, x_j)$$
 is again a (uniform) $L_{\ref{lem:three geodesics}}(K,\delta)$-quasigeodesic in $Y^h$. 
\end{lem}
\proof We first consider the concatenation $\beta(y_i,y_j)\star \beta(y_j, x_j)$ and observe that it satisfies the assumptions of 
Lemma \ref{lem:two geodesics} with the peripheral horoball $H^h$ equal $H_j^h$. Thus, $\beta=\beta(y_i,y_j)\star \beta(y_j, x_j)$ is a $K_1= L_{\ref{lem:two geodesics}}(K,\delta)$-quasigeodesic in $Y^h$.  Then we consider the concatenation $\al\star \beta$, where $\al= \beta(x_i,y_i)$. This concatenation again  satisfies the assumptions of 
Lemma \ref{lem:two geodesics} with the peripheral horoball $H^h$ equal  $H_i^h$. Lemma follows.  \qed 

\medskip
We now analyze the nearest-point projection to quasigeodesics of the type described in Lemma \ref{lem:two geodesics}:

\begin{lem}\label{lem:proj-2-gamma} 
Consider a quasigeodesic $\gamma=\beta\star \al$ as in Lemma \ref{lem:two geodesics}. Then for each point $p\in H^h$ the nearest-point projection $\bar{p}= P_{\ga}(z)$ satisfies:
$$
d(\bar{p}, \al)\le C_{\ref{lem:proj-2-gamma}}(K). 
$$
\end{lem}
\proof We continue with the notation of Lemma \ref{lem:two geodesics}. Suppose that $\bar{p}\in \beta$. Since $\gamma$ is a $K$-quasigeodesic, it is $\la=\la_{\ref{lem:qi-preserves2}}(\delta, K)$-quasiconvex; hence (see Lemma \ref{lem:projection-1}), the geodesic $p z$ in $Y^h$ passes within distance $\le D=\la+2\delta$ from $\bar{p}$. In particular, 
for
$$
r=D+\la_{\ref{lem:qc horoballs}}(\delta),
$$
$d(\bar{p}, H^h)\le r$, i.e. 
$\bar{p}\in N_r(H)$. We also have $p\in H$ and $\beta\cap \mathring{H}=\emptyset$. 
Hence, according to  Lemma \ref{lem:horodiving}(2) the distance 
between $\bar{p}$ and $z$ is $\le 2K$. \qed 

\begin{lem}\label{lem:cobounded horoballs} 
The peripheral horoballs $H_i^h, i\in I$, are uniformly pairwise cobounded in $Y^h$. 
\end{lem}
\proof Consider two distinct peripheral horoballs $H_i^h, H_j^h$ and suppose that $z_i\in H_i, z_j\in H_j$ are within distance 
$c$ from each other in $Y^h$. The points $z_i, z_j$ are connected by a path $\beta$ of length $\le K=4e^c$ in $Y$. 
Hence, $\beta$ is a $K$-quasigeodesic in $Y^h$. Let $\rho_i, \rho_j$ denote the vertical geodesic rays in $H_i^h, H_j^h$ emanating from $z_i, z_j$ respectively. Then,  by Lemma \ref{lem:three geodesics},  the 
concatenation $\al$ of the rays $\rho_i, \rho_j$ and the path $\beta$ is an $L=L_{\ref{lem:three geodesics}}(K,\delta)$-quasigeodesic in $Y^h$.

   \begin{figure}[tbh]
\centering
\includegraphics[width=70mm]{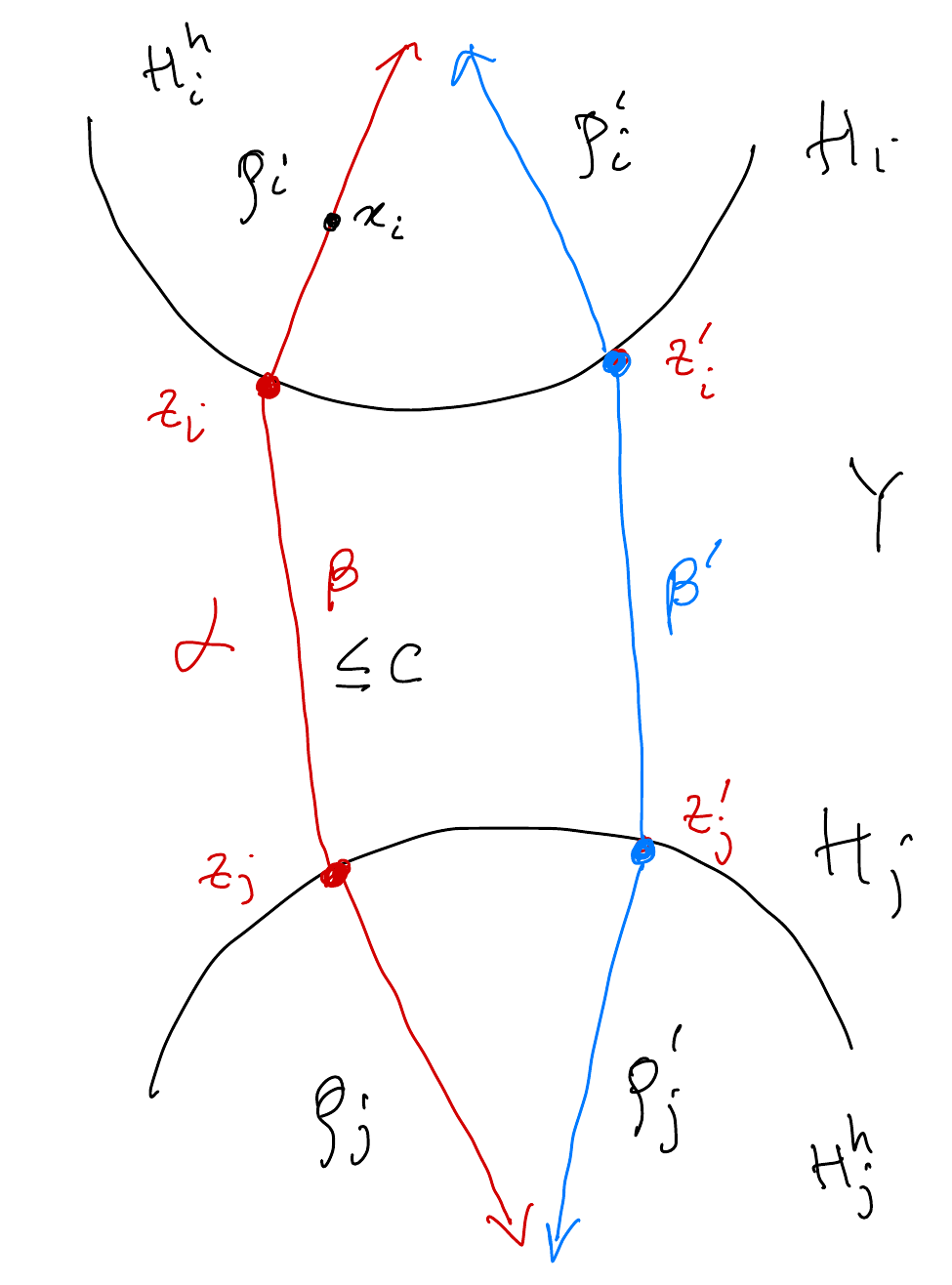}  
\caption{}
\label{c3.fig}
\end{figure}

Suppose now that we have four points $z_i, z_i'\in H_i, z_j, z_j'\in H_j$ satisfying $d_{Y^h}(z_i, z_j)\le c$, 
$d_{Y^h}(z'_i, z'_j)\le c$. We then, as above, form two $L$-quasigeodesic lines $\al, \al'$ in $Y^h$ passing through the points 
$z_i, z_j$ and $z_i', z'_j$ respectively and asymptotic to the centers $\xi(H_i^h)$, $\xi(H_j^h)$ of our peripheral horoballs. 
Thus, the $\delta$-hyperbolicity of $Y^h$, and the  fact that $\al, \al'$ are $L$-quasigeodesics in $Y^h$ asymptotic to the 
same pairs of points at infinity, imply that the Hausdorff distance between these quasigeodesics is 
$\le D=D_{\ref{lem:gen-bigons}}(L,\delta)$, see Lemma \ref{lem:gen-bigons}.  

We claim a uniform upper bound on the distance between $z_i, z_i'$. By Lemma \ref{lem:proj-2-gamma}, there exists a point $x_i\in \rho_i$ within distance 
$C=C_{\ref{lem:proj-2-gamma}}(L)$ from the nearest-point projection $\bar{z}'_i$ of $z_i$ to $\al'$. Hence, 
$$
d(z'_i, x_i)\le D'=D+C. 
$$
By the upper bound $\Hd(\al, \al')\le D$, in view of  
Lemma \ref{lem:proj-2-gamma}, there exists a point $x_i\in \rho_i$ within distance $D'=D'_{\ref{lem:proj-2-gamma}}(L,D)$ from $z_i'$.  
Since $z_i$ is the point in $H_i$ closest to $x_i$, it follows that $d(x_i,z_i)\le D'$ as well. Combining the inequalities, we get $d(z_i,z'_i)\le 2D'$. 
Similarly, $d(z_j, z'_j)\le 2D'$. Uniform quasiconvexity of the horoballs $H_i^h, H^h_j$ combined with Proposition \ref{prop:cobounded2} now implies that these horoballs are uniformly cobounded.   \qed

\begin{lem}\label{lemma:relative-limits}
Every point $\xi\in \geo Y^h$ is the limit of a sequence $(z_n)$ in $Y$, unless $\xi$ is the center of a horoball  $H^h$ with bounded horosphere $H$. 
\end{lem}
\proof   
There are two cases to consider.

1. $\xi=\xi(H^h)$ is the center of a peripheral horoball $H^h$. Fix a base-point $y_0\in H$ and the vertical ray $\rho$ 
in $H^h$ emanating from $y_0$. Since $H$ is unbounded, take a sequence $z_n\in H$ which diverges to infinity, i.e. $D_n=d_H(y_0, z_n)\to \infty$. We claim that $(z_n)$  converges to $\xi$.  Indeed, for $y_n:=(z_n, \log(D_n))\in H^h$, the combing path $c(y_0, y_n)$ in $H^h$ contains the subsegment of length $\log(D_n)$ in $\rho$. Hence, $\lim_{n\to\infty} z_n= \lim_{n\to\infty} \rho(\log(D_n))= \xi$. 

2. $\xi$ is not the center of any peripheral horoball. Then each geodesic ray $\rho=y\xi$ will cross into $Y$ along a sequence $z_n$ which diverges to infinity. Thus, $(z_n)$ converges to $\xi$. 
\qed

\subsection{Electrification and hyperbolization of quasigeodesics}\label{sec:electrification}  

In this section we describe two procedures of converting paths in $Y^h$ to paths in $Y^\ell$ and vice-versa.

\medskip
We first describe the procedure of {\em electri\-fication} of continuous hyperbolic quasigeode\-sics. 
In view of Corollary \ref{cor:deep-pieces}, for each continuous  $K$-quasigeodesic $\beta: [0,T]\to Y^h$ we obtain a maximal collection of maximal pairwise disjoint subintervals $[s_1, t_1]$,...,$[s_n,t_n]$ in $[0,T]$ 
satisfying the assumptions of  Corollary \ref{cor:deep-pieces}  
with respect to certain peripheral subspaces $H_{j_1},...,H_{j_n}\in {\mathcal H}$.  We then perform the following {\em electrification} procedure on $\beta$: 

For each subinterval $[s_i, t_i]$ we replace the restriction of $\beta$ to $[s_i,t_i]$ with the concatenation of two geodesic segments 
(of length $1/2$ each) in $H_{j_i}^\ell$ connecting $\beta(s_i), \beta(t_i)$ to the apex $a(H_{j_i}^\ell)$ of the cone $H_{j_i}^\ell$. 
In the special case when we have to deal with the subintervals $[0, t_1]$ and/or $[s_n, T]$ (i.e, $\beta(0)$ or $\beta(T)$ belongs to one of the open horoballs $\mathring{H}_j^h$), we replace $\beta|_{[0, t_1]}$ (resp. $\beta|_{[s_n, T]}$) with the unit length geodesic in $H_j^\ell$ connecting the apex $a(H_j^\ell)$ to 
$\beta(t_1)$ (resp. $\beta(s_n)$). We let $\beta^\ell$ denote the resulting path in $X^\ell$. Lastly, each subpath of $\beta$ whose domain has length $\le 2L(K)$ and which connects points of $H_i$ in $H_i^h$, is replaced by a geodesic in $H_i$ connecting the same points. Lengths of such geodesics are uniformly bounded by some constant $E=E(K)$ since horospheres $H_i$ are uniformly properly embedded in the horoballs $H_i^h$.   
The resulting map $\beta^\ell={\mathcal E}_{\mathcal P}(\beta)$ will be called the {\em electrification} of $\beta$.

We note that $\beta^\ell$ visits each cone-point $a(H^\ell_{j_i}), i=1,...,n$,  exactly once; more precisely, 
$\beta^\ell$ is {\em tight} in the following sense:

\begin{defn}
A continuous path $\gamma$ in $Y^\ell$ is {\em tight} if for each cone $H_i^\ell$, $i\in I$, the preimage 
$\gamma^{-1}(\mathring{H}^\ell_i)$ 
is a (possibly empty) interval and the restriction of $\gamma$ to this interval is 1-1.  
\end{defn}

\begin{rem}
1. The commonly used name for tight paths is {\em paths without backtracking}; we find this terminology cumbersome.  

2. The continuity assumption here is simply a matter of convenience and electrification can be defined for general 
$K$-quasigeodesics by using $s_i, t_i$'s such that $\beta(s_i), \beta(t_i)$ are within distance $K$ from the horosphere $H_i$. 

3. The electrification procedure described in \cite{MR3646028} is similar, except they replace by geodesics in $H_i^\ell$ 
the subsegments in $\beta$  with end-point in $H_i$ at distance $> 1$. This construction, while technically simpler, does not result in tight paths, which  we find undesirable. 
\end{rem}

\medskip 
More generally, we define tightness of paths in $Y^h$ as follows: 

\begin{defn}\label{defn:tight-path} 
We say that a continuous path $\beta$ in $Y^h$ is {\em tight} (relative to ${\mathcal H}$) if for each peripheral horoball \index{tight path}
$H\in {\mathcal H}$, the preimage $\beta^{-1}(\mathring{H}^h)$ is an open interval (possibly empty). 
\end{defn}

It is clear that for each tight path $\beta$, its electrification $\beta^\ell$ is a tight path in $Y^\ell$.

\medskip
We next use electrification to compare distances in $Y^h$ and in $Y^\ell$:  

\begin{lem}\label{lem:contraction} 
For any two pair of points $x, y\in Y$, ${d}^\ell(x,y)\le d_{Y^h}(x,y)$, i.e. the inclusion map $(Y, d_{Y^h})\to (Y^\ell, d^\ell)$
 is 1-Lipschitz. 
\end{lem}
\proof Recall that by the definition of the metric on the hyperbolic cones $H_i^h$, 
if $x, y$ belong to the same peripheral subspace $H_i$ and $d_{Y^h}(x,y)\le 1$, then $d_Y(x,y)\le 1$. It follows that if $\beta$ is a geodesic in $Y^h$ connecting the points $x, y\in Y$ then the length of $\beta^\ell$ is at most the length of $\beta$. \qed

\medskip 
{\bf Hyperbolization of electric geodesics.} Conversely, given any  (continuous) 
path $\beta$ in $Y^\ell$, we define its {\em hyperbolization} $\beta^h$ by replacing each subpath $\beta_H$ of $\beta$ 
connecting $x, y\in H\in {\mathcal H}$ and contained in $\mathring{H}^\ell$ except for its end-points, with the combing path 
$c(x,y)$ in $H^h$, see Section \ref{sec:hyperbolic cone}.  Recall that the paths $c(x,y)$ are uniform quasigeodesics in $H^h$ (with respect to the intrinsic metric $d^h$ on $H^h$). 
It is clear that if $\beta$ was tight, so is $\beta^h$.

In the case when $\beta$ is a quasigeodesic in $Y^\ell$, the path $\beta^h$ is called an 
{\em electro-ambient quasigeodesic}, see \cite{MR3646028}. This construction yields a collection of paths in $Y^h$ connecting points in $Y$. In Section \ref{sec:Farb->Gromov} we describe (uniform) electro-ambient quasigeodesics connecting arbitrary pairs of points in $Y^h$. 

\medskip 
The next result appears in \cite[Lemma 2.15]{MR3646028}.

\begin{lem}\label{lem:Dahmani-Mj} 
For tight\footnote{This assumption was forgotten in  \cite{MR3646028}.} uniform quasigeodesics $\beta$ in ${Y}^\ell$ connecting points of $Y$, the paths $\beta^h$ are uniform quasigeodesics in $Y^h$. More precisely, there exists a function $L=L_{\ref{lem:Dahmani-Mj} }(K)$ such that if $\beta$ is a tight $K$-quasigeodesic in $Y^\ell$ connecting points of $Y$, its hyperbolization $\beta^h$ is an $L$-quasigeodesic in $Y^h$. 
\end{lem}

An alternative and detailed proof was given by A.~Pal and A.~Kumar Singh in \cite{Pal-Kumar}; we will discuss this further in   Section \ref{sec:Farb->Gromov}.

\section{Hyperbolicity of the electric space} \label{sec:el-space}

In this section we will prove that for each relatively hyperbolic space $(Y, {\mathcal H})$, the space ${Y}^\ell$ is hyperbolic and describe uniform quasigeodesics in this space. The key result of this section is:

\begin{prop}
Electrifications $\al^\ell$ of uniform (continuous) quasigeodesics $\al$ in $Y^h$ connecting points of $Y$, 
are uniformly proper in $Y^{\ell}$. 
\end{prop}
\proof Let $\al: J=[0,T]\to Y^h$ be a continuous $K$-quasigeodesic in $Y^h$ connecting points $x, y\in Y$. Our goal is to prove that the length of the domain of $\al^{\ell}$ is bounded in terms of the distance $M:={d}^\ell(x,y)$. This will imply uniform properness 
of the maps $\al^{\ell}$ for all continuous $K$-quasigeodesics $\al$ in $Y^h$. 

Let $\beta$ be a geodesic in $Y^{\ell}$ connecting $x$ to $y$ and the length of $\beta$ is $\le M$. In particular, $\beta$ goes through at most $N$ cones $H^\ell_i$ where $N$ depends only on $M$; these cones correspond to the horoballs $B_1=H_{i_1}^h,...,B_n=H_{i_n}^h$ in $Y^h$, $n\le N$. It also follows that the length of each subsegment of $\beta$ between two distinct horoballs is at most $M$. By renumbering the horoballs, we can also label these subsegments 
$\beta_1,...,\beta_{n-1}$ so that $\beta_i$ connects $B_i$ to $B_{i+1}$. We also label $\beta_0$ the subsegment of $\beta$ between $x$ and $B_1$ and label $\beta_n$ the subsegment of $\beta$ between $B_n$ and $y$. Set
$$
B:= \beta_0\cup B_1\cup .... \cup B_n \cup \beta_n.  
$$
In particular, the subset $B\subset Y^h$  is $\la=\la(M)$-quasiconvex in $Y^h$. According to Lemma 
\ref{cobdd-cor}, there exists $D$ and 
$\epsilon$ such that for each component $C$ of $\al^{-1}(Y^h - N_D(B))$ the projection of $\al(C)$ to $B$ has diameter $\le \eps$.  
(Here $\eps, D$ depend on $\delta, K$ and $\la$, and   $D=D_{\ref{cobdd-cor}}(\delta, \la')$, where $\la'$ is the maximum of 
$\la$ and the quasiconvexity constant of the images of $K$-quasigeodesics in $Y^h$.) 
We now decompose the interval $J=[0,T]$ according to the mutual position of points $\al(t), t\in J,$ with respect  to the horoballs $B_i$ and the segments $\beta_i$. 

Let $J(\beta_j)$ denote the maximal subinterval in $J$ such that the images of the end-points of $J(\beta_j)$ under $\al$ belong to $N_D(\beta_j)$. Similarly, we define subintervals $J(B_j)$. 

1. Since $\al$ is a $K$-quasigeodesic and $\beta_j$ has length $\le M$, the length of each $J(\beta_j)$ is $\le K(2D+M+1)$. 

2. Similarly, Lemma \ref{lem:horodiving} implies that the length of each $J(B_j) \setminus \al^{-1}(B_j)$ is at most $2L_{\ref{lem:horodiving}}(K+D,0,\delta)$. Recall that by the definition of $\al^{\ell}$, whenever  $J(B_j) \setminus \al^{-1}(B_j)$ is nonempty, the interval $J(B_j)$ contributes the length 
$\le 2+ 2L_{\ref{lem:horodiving}}(K+D,0,\delta)$ to the length of the domain of $\al^{\ell}$. 

3. It now remains to estimate the length of 
$$
J':= J\setminus \bigcup_{i=0}^{n+1} (J(\beta_i) \cup J(B_i)) 
$$ 
Let $C=[s,t]$ be a component of this complement. Then 
both points $\al(s), \al(t)$ belong to $\partial N_D(B)$ and the diameter of the projection of $\al(C)$ to $B$ is $\le \eps$. Hence, the distance in $Y^h$ between $\al(s), \al(t)$ is at most $2D+\eps$. It follows that $C$ has length $\le K(2D+\eps+1)$. Since the number of components $C$ is $\le 2N+2$, it follows that the total length of $J'$ is at most
$$
K(2D+\eps+1) (2N+2). 
$$
Combining these estimates we conclude that the domain of $\al^{\ell}$ has length at most
$$
(M+1)K(2D+M+1)+ (2+ 2L_{\ref{lem:horodiving}}(K+D,0,\delta))M + K(2D+\eps+1) (2N+2). 
$$
Uniform properness of the maps $\al^{\ell}$ follows. \qed 

\medskip
As an application of this result we prove:

\begin{thm}\label{thm:hyperbolicity-of-electric-space} 
1. If $(Y, {\mathcal H})$ is relatively hyperbolic, then $Y^{\ell}$ is hyperbolic. 

2. Moreover, continuous $K$-quasigeodesics $\al$ in $Y^h$ yield $k_{\ref{thm:hyperbolicity-of-electric-space}}(K)$-quasigeodesics ${\al}^\ell$ in  $Y^{\ell}$. 
\end{thm}
\proof  We will verify that the conditions of Corollary \ref{cor:bowditch} are met. Namely, we will check that the combing of $Y^h$ by continuous $K$-quasigeodesics $\al$ results in a thin combing of $Y^{\ell}$ by paths $\al^{\ell}$ connecting points of the subset $Y^{\ell}_0\subset Y^{\ell}$ which is the union of $Y$ and the set $a(Y^\ell)$ of apexes $a_i$ of the cones $H^\ell_i$, $i\in I$.  We already know that the paths $\al^{\ell}$ are uniformly proper (Property 1 in Corollary \ref{cor:bowditch}).  
Let us verify Property 2. Consider a triple of $K$-quasigeodesics 
$\al_{x,y}, \al_{y,z}, \al_{z,x}$ connecting points $x, y, z$ in $Y^h$. Since $Y^h$ is $\delta$-hyperbolic, 
$$
\al_{x,y}\subset N_{3D_{\ref{stab-qg}}(\delta,K)+\delta}( \al_{y,z}\cup \al_{z,x}). 
$$
Let $u$ be a point of $\al_{x,y}\cap Y$ and suppose that $v$ is a point in $\al_{y,z}$ at a distance 
$d_{Y^h}(u,v)\le  3D_{\ref{stab-qg}}(\delta,K)+\delta$ from $u$. If $v$ happens to be in $Y$, then Lemma \ref{lem:contraction} implies that
$$
d^\ell(u,v)\le 3D_{\ref{stab-qg}}(\delta,K)+\delta
$$
as well. If $v=\al_{y,z}(t)$ belongs to an open horoball $\mathring{H}_i^h$ but $t$ does not lie in one of the subintervals $[s_i, t_i]$ in the domain of the path $\al_{y,z}$ which are coned-off when we define $\al^{\ell}_{y,z}$, then the distance in $Y^h$ from 
$\al_{y,z}(t)$ to $\al_{y,z}\cap H_i$ is at most $2(K+1) L_{\ref{lem:horodiving}}(K,0,\delta)$ and, thus, $d^\ell(u, \al_{y,z})\le 2(K+1) L_{\ref{lem:horodiving}}(K,0,\delta)$ as well. 
Lastly, if $t$ is one of the intervals $[s_i, t_i]$ then $p_i$ is in $\al_{y,z}$ and $d^\ell(u, p_i)\le d(u,v)+1\le 3D_{\ref{stab-qg}}(\delta,K)+\delta +1$. 

Finally, each point of $\al^{\ell}_{x,y}$ lies within unit distance (as measured in $Y^{\ell}$) 
from a point of $\al^{\ell}_{x,y}\cap Y$ and, therefore, 
for 
$$
R=3D_{\ref{stab-qg}}(\delta,K)+\delta +1 + 2(K+1) L_{\ref{lem:horodiving}}(K,0,\delta),
$$
the path $\al^{\ell}_{x,y}$ is contained in the $R$-neighborhood
 (in $Y^{\ell}$) of the union $\al^{\ell}_{y,z}\cup \al^{\ell}_{z,x}$.  \qed 

\medskip 
The following consequence of the theorem appears in \cite[Proposition 4.3]{klarreich}, see 
also \cite[Lemma 4.5 and Lemmas 4.8, 4.9]{farb-relhyp} in the setting of manifolds of 
negative curvature:

\begin{cor}\label{cor:hyperbolicity-of-electric-space} 
For $x, y\in Y$ let $\al=[xy]_{Y^\ell}$ and $\beta= [xy]_{Y^{h}}$. Then 
the Hausdorff distance in $Y$ between $\al\cap Y$ and $\beta\cap Y$ is uniformly bounded. 
\end{cor}
\proof 1. Take $z\in \al\cap Y$. Since $\al^h$ is a uniform $k$-quasigeodesic in $Y^h$ connecting the endpoints of $\beta$, 
the $k$-quasigeodesics $\al^h$ and $\beta$ are $D$-Hausdorff-close in $Y^h$, $D=D_{\ref{stab-qg}}(\delta, k)$. 
Thus, there exists $w\in \beta$ within  distance $D$ from $z$. If $w$ happens to be 
in $Y$, we are done. Suppose, therefore, that $w\in \mathring{H}$ for some $H\in {\mathcal H}$. Since $w$ is within distance $D$ 
(as measured in $Y^h$) from $Y$, and $\beta$ is a geodesic in $Y^h$, 
it follows that there exists a point $v\in \beta\cap H$ within distance $C=D+D'$ from $w$, where $D'=D_{\ref{stab-qg}}(\delta, K)$ 
and $K$ is the qi constant of the combing paths in $H^h$. Since $d_{Y^h}(v,z)\le D + C$, it follows that $d_Y(v,z)$ is also uniformly bounded by a uniform constant $E$, as $(Y,d)$ is uniformly properly embedded in $(Y^h, d^h)$, see Lemma \ref{lem:qc horoballs}(2). Hence, 
$$
\al\cap Y\subset N^Y_E(\beta\cap Y). 
$$

2. The proof of the opposite inclusion is similar, swapping the roles of $\al^h$ and $\beta$ and is left to the reader. \qed

\medskip 
We also record another application of Theorem \ref{thm:hyperbolicity-of-electric-space}:

\begin{cor}\label{cor:tight combing} 
There exists a constant $L=L_{\ref{cor:tight combing}}(\delta)$ such that every two points $x, y$ is $X$ are connected by a tight $L$-quasigeodesic in $X^\ell$. Namely, given a geodesic $c= [xy]_{X^h}$ in $X^h$, its electrification $c^\ell$ is a uniform tight quasigeodesic in $X^\ell$.  
\end{cor}

We define a coarse projection $q=q_Y: Y^h\to {Y}^\ell$ as follows: 

1. The map $q$ is the identity on $Y$. 

2. For each peripheral subspace $H_i$, $q(\mathring{H}_i^h) =\{a(H_i^\ell)\}$, 
the apex of the cone over $H_i$. 

\medskip 
\noindent The following lemma is immediate: 

\begin{lemma}
The map $q$ is $(1,1)$-coarse Lipschitz. 
\end{lemma}

The next result is a direct corollary of Theorem \ref{thm:hyperbolicity-of-electric-space}(2):

\begin{cor}
There exists a function $\la\mapsto \hat\la$ such that if $Z\subset Y^h$ is a $\la$-quasiconvex subset, then $q(Z)$ is $\hat\la$-quasiconvex in $Y^\ell$. 
\end{cor}

\begin{lemma}\label{lem:coning-projection-commutes}
The projection $q$ coarsely commutes with nearest-point projections. More precisely, there is  
$D=D_{\ref{lem:coning-projection-commutes}}(\delta,\la)$, 
where $\delta$ is the hyperbolicity constant of $Y^h$, such that for each $\la$-quasiconvex subset $Z\subset Y^h$, 
$$
d(P_{Y^\ell, q(Z)}, q\circ P_{Y^h,Z})\le D. 
$$
\end{lemma}
\proof For $y\in Y^h$, up to a uniformly bounded error, the projection  $\bar{y}=P_{Y^h,Z}(y)$ of $y$ to $Z$ is defined by the property that for each $z\in Z$ the concatenation $\beta=y\bar{y} \star \bar{y}z=\beta_1\star \beta_2$ is a uniform quasigeodesic (see Lemma \ref{lem:projection-1}). Suppose first that 
$\bar{y}\in Y$. Theorem \ref{thm:hyperbolicity-of-electric-space}(2) implies that  
the concatenation $\beta^\ell= \beta^\ell_1\star \beta^\ell_2$ is a uniform quasigeodesic in 
$\hat{Y}$. By applying Lemma \ref{lem:citerion-of-projection}, we conclude that $\bar{y}$ is uniformly close to $P_{{Y}^\ell, q(Z)}(q(y))$. 

Suppose now that $\bar{y}$ lies in some $H_i^h$; let $a_i=a(H_i^h)$ denote the apex of the cone in ${Y}^\ell$ corresponding to $H_i$, and $\beta_0\subset H^h_i$ be the maximal subsegment of $\beta$ contained in $H_i^h$ and connecting points $x_i, y_i\in H_i$.  Then $\beta$ is the concatenation $\beta_1\star \beta_0\star \beta_2$. Accordingly, 
 $\beta^\ell$ is the concatenation $\beta^\ell_1\star [x_i a_i]\star [a_i y_i]\star \beta^\ell_2$. Thus, $a_i$ is uniformly close to  
 $P_{Y^\ell, q(Z)}(q(y))$. At the same time,  $q(\bar{y})=a_i$. \qed

\medskip 
Lastly, we discuss relatively hyperbolic structures on metric spaces which are already hyperbolic, see \cite{farb-relhyp}:

\begin{thm}\label{thm:H->RH}
Suppose that $X$ is a hyperbolic space and ${\mathcal H}= \{H_i, i\in I\}$ is a collection of uniformly pairwise cobounded uniformly quasiconvex subsets in $X$. Then: 

1. The pair $(X,  {\mathcal H})$ is relatively hyperbolic. 

2. Conversely, if $X$ is hyperbolic and 
${\mathcal H}= \{H_i, i\in I\}$ is a collection of uniformly quasiconvex subsets, and $X^h$ is hyperbolic, then 
the subsets $\{H_i, i\in I\}$ are uniformly pairwise cobounded in $X$. 
\end{thm} 
 \proof 1.  We first equip $X^h$ with a structure of a tree of hyperbolic spaces $\Y=(\pi: Y\to T)$. 
 We define the tree $T$ by  taking the wedge of rays $R_i, i\in I$, where each $R_i$ is the positive half-line $[0,\infty)$ equipped with the standard simplicial structure (vertices are the nonnegative integers) and the wedge of rays is obtained by identifying $0$'es in all $R_i$'s with a single vertex $v_0\in T$. 
 For $n\in {\mathbb N}\subset R_i$ we use the notation $v_{in}$ for the corresponding vertex in $T$. 
We define vertex-spaces of  $\Y$ as $Y_{v_0}:= X$; $Y_{v_{in}}:= H_i\times \{n\}\subset H_i^h=H_i\times [1,\infty)$. The edge-space $Y_{e_{in}}$ for the edge $e_{in}=[v_{in}, v_{in+1}]$ is $H_i\times \{n+\frac{1}{2}\}$. The incidence maps 
$$
f_{e_{in},v_{in}}: Y_{e_{in}} \to Y_{v_{in}}, f_{e_{in},v_{in+1}}: Y_{e_{in}} \to Y_{v_{in+1}}
$$
 are given by 
 $$
 (y, n+\frac{1}{2})\mapsto (y, n), \quad (y, n+\frac{1}{2})\mapsto (y, n+1)
 $$
 respectively. These maps clearly are uniform quasiisometries. Thus, we obtain a tree of hyperbolic spaces $\Y=(\pi: Y\to T)$. 

  Since  the hyperbolic cones $H_i^h$ are uniformly hyperbolic (see Proposition \ref{prop:hyper-cone}), the restrictions of the tree of spaces $\Y$ to the rays $R_i$ satisfy the uniform flaring condition with flaring constants independent of $i$, see Lemma \ref{lem:hyp->uniform flaring}. Suppose that $\llbracket v, w\rrbracket\subset T$ is an interval containing vertices (different from $v_0$) 
 of different rays $R_i, R_j$. Let $\gamma_0, \gamma_1$ be $K$-qi sections of $\pi: Y\to T$ over $\llbracket v, w\rrbracket$. The assumption that  
 peripheral subspaces $H_i, H_j$ in $X$ are uniformly  cobounded implies that 
 $$
 d_{Y_{v_{i1}}}(\gamma_0( v_{i1}), \gamma_1(v_{i1}))\le C(K), \quad d_{Y_{v_{j1}}}(\gamma_0( v_{j1}), \gamma_1(v_{j1}))\le C(K). 
 $$
 In view of the uniform flaring of $\Y$ over the rays $R_i, R_j$, it then follows that 
 $$
 d_{Y_{v_{in}}}(\gamma_0( v_{i1}), \gamma_1(v_{in})) 
 $$
 is either uniformly bounded (in terms of $K$) or grows at a linear rate (as a function of $n$). The same applies to 
 $$
 d_{Y_{v_{jn}}}(\gamma_0( v_{j1}), \gamma_1(v_{jn})). 
 $$
 Hence $\Y$ satisfies the uniform flaring condition and, thus,  $Y$ is uniformly hyperbolic by Theorem \ref{thm:mainBF}. 
 
 \medskip 
 2.  This part is a consequence of Lemma \ref{lem:cobounded horoballs}.  \qed

\medskip
This theorem has a useful addendum, relating quasigeodesics in $X$ and in $X^h$. In the setting of the theorem, let $\beta$ be a $k$-quasigeodesic in $X^h$ connecting points in $X$. For each maximal subsegment $\beta_i$ contained in some $H_i^h$, we replace $\beta_i$ with a geodesic in $X$ connecting the end-points of $\beta_i$. We let $\beta_X$ denote the resulting path in $X$. 

\begin{lem}\label{lem:relativization} 
The  paths $\beta_X$ in $X$ are  $K_{\ref{lem:relativization}}(k)$-quasigeodesic. 
\end{lem} 
\proof Thinking of $X^h$ as the total space $Y$ of a tree of spaces $\Y$ as above, we note that the path $\beta_X$ obtained via the above procedure of converting $\beta$ to $\beta_X$ 
is exactly the cut-and-replace procedure in Definition \ref{defn:detour}. 
 Now, the result follows from Theorem \ref{thm:cut-paste} (actually, it follows already from Proposition 
 \ref{prop:oneflow} proven earlier by Mitra in \cite{mitra-trees}). \qed

\subsection{Equivalence of the two definitions of relative hyperbolicity}\label{sec:Farb->Gromov}

We start by reviewing Farb's definition of relative hyperbolicity. Given a metric space $Y$ and a collection of its rectifiably connected, uniformly properly embedded subspaces ${\mathcal H}= \{H_i: i\in I\}$, we get the associated electric space 
$Y^\ell={\mathcal E}(Y, {\mathcal H})$ with the metric $d^\ell$ as  described in Section \ref{sec:GRH}. 

\begin{defn}\index{weakly relatively hyperbolic space} 
The pair $(Y, {\mathcal H})$ is said to be {\em weakly relatively hyperbolic} if the metric space $(Y^\ell, d^\ell)$ is hyperbolic. 
\end{defn}

Every two points in $Y^\ell$ are  connected by a tight $L$-quasigeodesic in $Y^\ell$ for some uniform constant $L$ (for instance, one can use geodesics in $Y^\ell$). We refer to the tight quasigeodesics in $Y^\ell$ connecting points in $Y$ as {\em electric quasigeodesics} in $Y^\ell$. \index{electric quasigeodesic}

In Section \ref{sec:GRH}, given an electric quasigeodesic $\beta$ in $Y^\ell$ connecting  points $x, y\in Y$, 
we defined a path $\beta^h$ in $Y^h$ (the hyperbolization of $\beta$) connecting $x$ and $y$. We now extend this definition to connect arbitrary pairs of points $x, y$ in $Y^h$. We let $\hat{x}=q(x), \hat{y}=q(y)$ denote the projections of $x, y$ to $Y^\ell$. Let $\beta$ be an electric quasigeodesic connecting $\hat{x}$ to $\hat{y}$. In the case when $x\in \mathring{H}^h, \hat{x}=a(H^\ell)$ 
is the apex of the cone $H^\ell$, we let $x_H\in H$ denote the exit point of $\beta$ from the cone $H^\ell$; similarly, if $y\in H^h$, we let $y_H$ denote the entry point 
of $\beta$ into the cone $H^\ell$. We set $x'=x$ if $x\in Y$ and $x':= x_H$ if $x\notin Y$; similarly, we define the point $y'$. 
Then for the subpath $\beta(x',y')$ of $\beta$ between $x', y'$, we define its hyperbolization $\beta(x',y')^h$ as before, and connect $x$ to $x'$, $y$ to $y'$ by geodesics in the corresponding horoballs $H^h$ in the case when $x\ne x'$ or $y\ne y'$. 

We, thus, obtain a family of paths connecting points of $Y^h$. As it turns out (see Theorem \ref{thm:Pal-Kumar}), under suitable assumptions, the resulting paths define a slim combing of $Y^h$. In the next definition, $D=D(K)$. 

\begin{defn}\index{quasigeodesics tracking each other} 
Two electric $K$-quasigeodesics $\al_1, \al_2$ with the same end-points are said to have {\em the same $D$-intersection pattern} with respect to the collection of cones 
$H^\ell_i, i\in I$, or, simply, {\em $D$-track each other}, if the following conditions hold for all $i\in I$:

1. If one path, say, $\al_1$, contains $a(H_i^\ell)$ but then other (namely, $\al_2$)  does not, then $d_{H_i}(\al_1(s), \al_1(t))\le D$, where $\al_1^{-1}(\mathring{H}_i^\ell)$ is the open interval $(s,t)$.   

2. Suppose that, for some $i\in I$, $\al_j^{-1}(\mathring{H}_i^\ell)=(s_i,t_i)\ne \emptyset, j=1,2$. Then 
$$
\max \{d_{H_i}(\al_1(s_1), \al_2(s_2)),  d_{H_i}(\al_1(t_1), \al_2(t_2)) \} \le D. 
$$
The function $D(K)$ is the {\em tracking function} of $(Y, {\mathcal H})$. \index{tracking function}
\end{defn}

\begin{defn}\index{space relatively hyperbolic in Farb's sense}
A pair $(Y, {\mathcal H})$ is  {\em  relatively hyperbolic in Farb's sense (FRH)} if it is weakly relatively hyperbolic and for every $K\ge 1$ there exists $D=D(K)$ such that 
any two electric $K$-quasigeodesics with the same end-points $D$-track each other. 
\end{defn}

The next theorem relating the two definitions was proven by A.~Pal and A.~Kumar Singh in   \cite{Pal-Kumar} (for relatively hyperbolic groups the equivalence of two definitions was known earlier, cf. \cite{bowditch-relhyp}): 

\begin{thm}\label{thm:Pal-Kumar} 
If $(Y, {\mathcal H})$ is FRH, then $Y^h$ is uniformly\footnote{With respect to the tracking function and the hyperbolicity constant of $Y^\ell$.} 
hyperbolic and the hyperbolization of uniform electric quasigeodesics yields uniform quasigeodesics in $Y^h$. In particular,  
$(Y, {\mathcal H})$ is GRH. 
\end{thm}

\begin{rem}
A.~Pal and A.~Kumar Singh in   \cite{Pal-Kumar} assume that the subsets $H_i, i\in I$, are uniformly separated in $Y$. In our setting one achieves uniform separation by replacing $Y$ with the space $Y'$ obtained by attaching the products $H_i\times [0,1]$, $i\in I$, along the subsets $H_i$, and replacing the subsets $H_i\subset  Y$ with $H_i\times \{1\}$, $i\in I$. 
\end{rem}

To conclude the discussion, we note a relation between  uniform tight quasigeodesics in $Y^\ell$ and uniform quasigeodesics 
in the space $Y$ itself. (This is not needed for any proofs in the book, but clarifies the overall picture.)  
The following result is proven in \cite[Lemma 8.8]{Hruska} in the context of relatively hyperbolic groups, but the proof works for general relatively hyperbolic spaces:

\begin{thm}\label{thm:Hruska} 
There are functions $D=D_{\ref{thm:Hruska}}(K,\delta)$ and $L=L_{\ref{thm:Hruska}}(K,\delta)$ 
such that if $\al$ is an electric $K$-quasigeodesic in $Y^\ell$ (connecting points $x, y$ of $Y$), then for every  
$L$-quasigeodesic $\beta$ in $Y$  also connecting $x$ and $y$, we have that 
$\al\cap Y$ is contained in the $D$-neighborhood of $\beta$, with respect to the metric $d_Y$ of $Y$.   
\end{thm}

\subsection{Morphisms of relatively hyperbolic spaces}

Given a pair of relatively hyperbolic spaces $(Y, {\mathcal H})$, $(Y', {\mathcal H}')$, a {\em relative morphism} of these pairs is a 
uniformly proper map $f: Y\to Y'$ such that: \index{morphism of relatively hyperbolic spaces}

a. For each $H'\in {\mathcal H}'$, $f^{-1}(H')$ is either empty or equals some $H\in {\mathcal H}$. 

b. For each $H\in {\mathcal H}$, there exists $H'\in {\mathcal H}'$ satisfying $f(H)\subset H'$.

\begin{rem}
1. The first condition implies that if $H_i, H_j\in {\mathcal H}$ are distinct, then $f(H_i), f(H_j)$ are not contained in the same $H'\in {\mathcal H}'$. 

2. One can relax a bit the above conditions by requiring existence of a uniform constant $D$ such that:

(a') For each $H'\in {\mathcal H}'$, $f^{-1}(H')$ either has diameter $\le D$ or is $D$-Hausdorff-close to some $H\in {\mathcal H}$.

(b') For each $H\in {\mathcal H}$, there exists $H'\in {\mathcal H}'$ satisfying $f(H)\subset N_D(H')$.

This is the approach taken in \cite{MS}. However, the two definitions are easily seen to be effectively equivalent since one can replace the peripheral horoballs $H^h, H\in {\mathcal H}$ and $(H')^h, H'\in {\mathcal H}'$, by suitable smaller subsets.   
\end{rem}

A relative morphism is said to be a {\em relative qi embedding} of the pairs if, additionally:

c. $f: Y\to Y'$ is a qi embedding.

\medskip
The qi constants of $f$  
are called the {\em parameters} of the relative qi embedding $f$.

\medskip 
Given a relative morphism $f: (Y, {\mathcal H})\to (Y', {\mathcal H}')$, the {\em coned-off map} $f^\ell: Y^\ell\to (Y')^\ell$ is defined as follows:  

i. The restriction of ${f}^\ell$ to $Y$ equals $f$. \index{coned-off map}

ii. Consider a peripheral subset $H\in {\mathcal H}$ such that $f(H)\subset H'\in {\mathcal H}'$ and let $a=a(H^\ell), a'=a((H')^\ell)$ denote the respective apexes of the cones $H^\ell, (H')^\ell$ over these peripheral subsets in $Y^\ell, (Y')^\ell$. Then ${f}^\ell(a)=a'$ and for each $x\in H$, $x':= f(x)\in H'$, the map $\hat{f}$ sends the radial segment $xa$ to $x'a'$ isometrically.  This defines $f^\ell$ on $H^\ell$ for each $H\in {\mathcal H}$.

\medskip 
Similarly, we define the {\em hyperbolic extension} $f^h$ of $f$, $f^h: Y^h\to (Y')^h$; this extension construction goes back to the work of Mostow on strong rigidity of nonuniform lattices in rank one Lie groups.  

For $x\in H\in {\mathcal H}, x'=f(x)\in H'\in {\mathcal H}'$, we consider the vertical geodesic rays 
$\rho_x: [0,\infty)\to H^h, \rho_{x'}: [0,\infty)\to (H')^h$, emanating from $x, x'$ and asymptotic to the centers of the horoball $H^h, (H')^h$. Then for $t\in [0,\infty)$, we set
$$
f^h(\rho_x(t)):= \rho_{x'}(t), \quad  x\in H. 
$$ 

\begin{lemma}\label{lemma: extended relhyp tree}
1. The hyperbolic extension $f^h$ of a  proper coarse Lipschitz map  is  again a  proper coarse Lipschitz map. 

2. If $f, f^\ell$ are   qi embeddings, then so is $f^h$ and the constants of $f^h$ depend only on the parameters of  $(Y, {\mathcal H}), (Y', {\mathcal H}')$ and qi constants of $f, f^\ell$.

3. Conversely, if $f$ and $f^h$ are qi embeddings, so is $f^\ell$.  
\end{lemma}
\proof Part 1. We will only show that $f^h$ is coarse Lipschitz and leave it to the reader to check properness. Since $Y^h$ is a path-metric space, the problem is local and we have to address it only in the horoballs $H^h$, $H\in {\mathcal H}$. 
Clearly, $f^h$ is isometric along the vertical geodesic rays in $H^h$. Suppose, therefore, that $y_1, y_2$ are points within unit distance on the same horosphere $H\times \{t\}\subset H^h, t\ge 0$ (where the distance is computed in the intrinsic metric of the horosphere). 
Thus, $y_i=\rho_{x_i}(t)$, $x_i\in H$, $i=1, 2$, and, by the definition of the metric on $H^h$,
$$
d_{H}(x_1, x_2)= e^t. 
$$ 

Since the horospheres $H, H'$ are uniformly qi embedded in $Y, Y'$, we conclude that the restriction map 
$f: (H, d_H)\to (H', d_{H'})$ is $(L,A)$-coarse Lipschitz,  
where $L, A$ depend only on the parameters of  $(Y, {\mathcal H})$ and $f$. 
Hence,  
$$
d_{H'}(f(x_1), f(x_2))\le L d_H(x_1, x_2) +A\le L e^t +A. 
$$
Again, by the definition of the metric on the horoball $(H')^h$,
$$
d(\rho_{x'_1}(t), \rho_{x'_2}(t))\le e^{-t}( L e^t +A)= L + e^{-t}A\le L+A. 
$$
This verifies that $f^h$ is uniformly coarse Lipschitz. 

Part 2. 
Note that given an ambient geodesic $\gamma$ in $Y^h$, its electrification $\gamma^\ell$ is a uniform electric quasigeodesic in ${Y}^\ell$, hence, $f^\ell$ carries it to a uniform electric quasigeodesic $\hat\gamma'$ in $(Y')^\ell$, coarsely preserving its arc-length. Then, applying Lemma \ref{lem:Dahmani-Mj}, we obtain that the hyperbolization $\gamma'= (\hat\gamma')^h$ of 
$\hat\gamma'$ is a uniform quasigeodesic in $(Y')^h$. By the construction, the map $f^h$ coarsely preserves arc-lengths of paths in $X^h$, cf. Section \ref{sec:hyperbolic cone}; and, thus, $f^h$ sends geodesics to uniform quasigeodesics coarsely preserving arc-length, hence, is a qi embedding.  

Part 3. The proof of this part is similar to the proof of Part 2, except one uses Theorem \ref{thm:hyperbolicity-of-electric-space}(2). \qed 

\begin{rem}\label{rem:MS}
This lemma is a weak form of a more general recent result by J. Mackay and A. Sisto \cite[Theorem 1.2]{MS}, who proved Part 2  without  assuming that $f^\ell$ is a qi embedding (in fact, they also relax the condition (c) in the definition of a relative qi embedding).  Part 3 of the lemma then implies that $f^\ell$ is a uniform qi embedding provided that $f$ is. 
See also the recent paper  \cite{Healy-Hruska} by Healy and Hruska with a similar result for relatively hyperbolic groups. 
\end{rem}
   
\subsection{CT maps}

Suppose that $(X, {\mathcal F})$ and   $(Y, {\mathcal H})$ are relatively hyperbolic spaces and 
$f: (Y, {\mathcal H})\to (X, {\mathcal F})$  is a morphism of pairs. We then have the hyperbolic extension of $f$, which is  a coarse Lipschitz uniformly proper map $f^h: Y^h\to X^h$ (see Lemma \ref{lemma: extended relhyp tree}), as well as the coned-off map 
${f}^\ell: {Y}^\ell\to {X}^\ell$. 
   
\begin{prop}\label{prop:rel-CT}
The map $f^h$ admits a CT-extension provided that 
the following holds: 

For some (equivalently, every) $y_0\in Y$, every $K\ge 1$, there is a function $C=C(y_0, K, D)$ such that: 
For all pairs of points $y_1, y_2\in Y$, all tight $K$-quasigeodesic $\beta^\ell$ in ${Y}^\ell$ connecting $y_1, y_2$, and all tight $K$-quasigeodesics $\al^\ell$ in ${X}^\ell$ connecting $x_1=f(y_1), x_2=f(y_2)$, 
$$
d_X(f(y_0), {\al}^\ell\cap X)\le D \Rightarrow d_Y(y_0, {\beta}^\ell\cap Y)\le C. 
$$
\end{prop}  
\proof Note that, in view of Corollary \ref{cor:hyperbolicity-of-electric-space}, the implication in the proposition  
 can be rewritten as  
$$
d_X(f(y_0), [x_1 x_2]_{X^h}\cap X)\le D \Rightarrow d_Y(y_0, [y_1 y_2]_{Y^h}\cap Y)\le C:= \phi(D),  
$$
for some function $\phi$. We will be verifying that $f^h$ satisfies Mitra's Criterion, Theorem \ref{thm:Mitra's Criterion}. 
Suppose that $y_1, y_2\in Y^h$ are such that 
$$
d_{X^h}(x_0, [x_1 x_2]_{X^h})\le D,$$
where $f(y_i)=x_i$, $i=0, 1,2$. Our goal is to get a uniform bound on the distance from $y_0$ to $[y_1 y_2]_{Y^h}$ in $Y^h$. 

{\bf Case 1.} We first suppose that $y_1, y_2$ are both in $Y$. Then, $x_1, x_2$ are both in $X$ as well. Since $X$ is uniformly properly embedded in $X^h$ and $x_0$ is in $X$, it follows that there is $x'\in [x_1x_2]_{X^h}\cap X$ within distance $D'$ from $x_0$, with respect to the metric of $X$, where $D'$ depends only on $D$ (and geometry of $(X, {\mathcal H})$, of course), cf. 
Lemma \ref{lem:strictly-convex}. Thus, the assumptions of the proposition imply that 
$$
d_Y(y_0, [y_1 y_2]_{Y^h}\cap Y)\le C',
$$
where $C'$ is a function of $D$. This, of course, implies that 
$$
d_{Y^h}(y_0, [y_1 y_2]_{Y^h})\le C',
$$
as required by Mitra's Criterion.

\medskip 
{\bf Case 2.} Suppose that both $y_1, y_2$ belong to the same $H_i^h$. In view of Lemma \ref{lem:strictly-convex}, 
we  obtain that for $k=1$ or $k=2$, $d(x_0, x_k)\le D'$, for some $D'$ depending only on $D$. 
Uniform properness of $f^h$ then implies 
that $d_{Y^h}(y_0, y_k)$ is also uniformly bounded, as required. 

\medskip 
{\bf Case 3.} Suppose that $y_1\in H_{i_1}^h, y_2\in H_{i_2}^h$  and $H_{i_1}^h\ne H_{i_2}^h$. Thus, $x_i\in F_{j_i}^h$,  
$F_{j_i}\in {\mathcal F}$, are distinct, $i=1,2$. Our assumption that 
$$
d_{X^h}(x_0, [x_1 x_2]_{X^h})\le D$$
implies that there is a point $p\in x_1x_2=[x_1 x_2]_{X^h}\cap X$ within distance $\le D'=D'(D)$ from $[x_1 x_2]_{X^h}$ 
(cf. the proof in Case 2). Let $x_i'\in F_{j_i}$, $i=1,2$, be points realizing the minimal distance in $X^h$ between $F_{i_1}, F_{i_2}$. Thus, the $R$-neighborhood of the segment $x_1 x_2$ in $X^h$ contains $x_1'x_2'$, where $R$ is a uniform constant, see Lemma \ref{cobdd-lem1}. Let $x_i''\in x_1x_2$ be points within distance $R$ from $x_i', i=1, 2$. By the uniform quasiconvexity of $H_{i_k}^h$ in $X^h$, $k=1,2$, if $p$ is not in $x_1'' x_2''=[x_1'' x_2'']_{X^h}$ then its distance to $x_1''$ or $x_2''$ is uniformly bounded (cf. Corollary \ref{lem:strict-convexity}). Thus, it suffices to consider the case when $p\in x_1'' x_2''$.   

Let $w_i, i=1,2$, be intersection points of $[y_1y_2]_{Y^h}=y_1y_2$ with $H_{i_1}^h, H_{i_2}^h$ respectively, and $z_i:= f(w_i)$, $t=1,2$. The points $z_1, z_2$ belong to the peripheral subsets $F_{j_1}, F_{j_2}$  respectively. Thus, again, the $R$-neighborhood of the segment $z_1z_2$ contains $x_1'x_2'$. Since we are assuming that $p\in x_1'' x_2''$, it follows that 
$$
p\in N_{2R+\delta}(z_1 z_2)
$$
with respect to the metric of $X^h$. Hence, $d_{X^h}(x_0, [z_1z_2]_{X^h})\le 2R+\delta+D'$. Now, we are in the setting of Case 1 and it follows that 
$$
d_{Y}(y_0, [w_1w_2]_{Y^h})\le \phi(2R+\delta+D').
$$
Therefore, 
$$
d_{Y^h}(y_0, [y_1 y_2]_{Y^h})\le d_{Y}(y_0, [w_1w_2]_{Y^h})\le \phi(2R+\delta+D'),
$$
as required. 

{\bf Case 4.} Suppose that $y_1\notin Y$ and $y_2\in Y$. This case is similar to Case 3 and we leave it to the reader. \qed

\section{Trees of relatively hyperbolic spaces} \label{sec:rel trees}

We can now describe axioms of trees of relatively hyperbolic spaces following \cite{mahan-pal}:

\begin{defn}
A tree of relatively hyperbolic spaces is an (abstract) tree of spaces $\X=(\pi: X\to T)$ where all vertex and edge spaces have structures of uniformly\footnote{I.e. each $X_v^h, X_e^h$ is $\delta$-hyperbolic for a uniform constant $\delta$.}  relatively hyperbolic spaces $(X_v, {\mathcal H}_v)$, $(X_e, {\mathcal H}_e)$, 
and the incidence maps $f_{ev}: X_e\to X_v$ are uniform relative qi embeddings.\footnote{Recall that this condition also requires the  incidence maps to be morphisms of these relatively hyperbolic spaces, 
$f_{ev}: (X_e, {\mathcal H}_e)\to (X_v, {\mathcal H}_v)$.} 
\end{defn}

\begin{rem}\label{rem:qi-cone-off} 
The definition of a tree of relatively hyperbolic spaces given in \cite{mahan-reeves} also requires the coned-off maps 
$f^\ell_{ev}: X_e^\ell\to X_v^\ell$ to be uniform qi embeddings. In view of \cite{MS}, this is a consequence of the 
assumption that the incidence maps are uniform relative qi embeddings (cf. Remark \ref{rem:MS}). 
The reader not willing to rely upon the results of \cite{MS}, can simply assume that the maps $f^\ell_{ev}$ 
are uniform qi embeddings. 
\end{rem}

 Next, given such a tree of spaces $\X=(\pi: X\to T)$, we define an equivalence relation on $X$ generated by the  following: 
 
 1. For each edge $e=[v_1,v_2]$ of $T$, every point in a peripheral subspace $H_{e}\subset X_e$ is equivalent to  all the points of the  mapping cylinders of the incidence maps $f_{ev_i}: H_e \to H_v\subset X_{v_i}, i=1,2$, in $X$. 
 
 2. All points of each peripheral subspace $H_v\subset X_v, v\in V(T)$, belong to the same equivalence class.  
 
 \medskip 
 Each equivalence class of this equivalence relation is called a {\em peripheral subspace} $P$ of $X$. Observe that for each peripheral subspace $P$ of $X$ and each vertex (resp. edge) space $X_v$  (resp. $X_e$) of $\X$, the intersection $P\cap X_v$ (resp, $P\cap X_e$) is either empty or equals one of the peripheral subspaces of 
 $X_v$  (resp. $X_e$). Given these peripheral subspaces, one defines the new space $X^P={\mathcal G}(X, {\mathcal P})$,
 by attaching hyperbolic horoballs $P_j^h$ to $X$ along all peripheral subspaces $P_j, j\in J$, of $X$. Here and in what follows, 
 ${\mathcal P}=\{P_j: j\in J\}$.

 \begin{rem}
The definition of $X^P$ should not be confused with the one, $X^h$, of the total space of a tree of spaces over $T$, with the vertex/edge spaces $X^h_v, X^h_e$.   
   \end{rem}

Similarly, for each subtree $S\subset T$, we define the space $X^P_S$ obtained by first restricting the tree of spaces $\X$ to $S$
and, thus, obtaining a tree of spaces $X_S\to S$,   and then applying the above procedure, so that   
$X^P_S=(X_S)^P$. The peripheral structure of $X_S$ is denoted ${\mathcal P}_S$.

\medskip 
Additionally, we have the space $\widehat{X}={\mathcal E}(X, {\mathcal P})$ obtained by coning-off the 
peripheral subspaces $P\in {\mathcal P}$ in $X$.

Furthermore, we define the tree of coned-off-spaces ${\X}^\ell=({\pi}: {X}^\ell\to T)$ 
with vertex/edge spaces ${X}^\ell_v, {X}^\ell_e$ and coned-off incidence maps 
${f}^\ell_{ev}:  {X}^\ell_e\to {X}^\ell_v$, which are uniform qi embeddings (see Remark \ref{rem:qi-cone-off}).

\begin{rem}
In the terminology of \cite{mahan-reeves}, ${X}^\ell$ is a {\em partially electrocuted space}. 
\end{rem}

The total space ${X}^\ell$ of $\X^\ell$ contains a {\em forest} ${\mathcal F}$ 
which is a collection of pairwise disjoint trees $F_P$ corresponding to the peripheral subspaces $P\in {\mathcal P}$: 

(i) The vertices $\nu_{iv}$ of $F_P$ are the apexes $a(H^\ell_{iv})$ of the cones over the peripheral subspaces $H_{iv}=P\cap X_v$ of the vertex-spaces $X_v, v\in V(T)$. 

(ii) Similarly, the edges $\eps_{ie}$ of $F_P$ are labeled by the apexes $a(H^\ell_{ie})$ 
of the cones over the peripheral subspaces of the edge-spaces $X_e, e\in E(T)$, where $H_{ie}=P\cap X_e$.  

(iii) A vertex $\nu_{iv}$ (corresponding to $a(H^\ell_{iv})$) of the tree $F_P$  
is incident to the edge $\eps_{ie}$ (corresponding to $a(H^\ell_{ie})$) if and only if 
the incidence map $f_{ev}$ of $\X$ sends $H_{ie}$ to $H_{iv}$. 

\medskip
The fact that each graph $F_{P}$ is a tree is a consequence of  the following lemma: 

\begin{lemma}\label{lem:qi-forest}
Under the projection $X^\ell\to T$, each  $F_{P}$ maps isomorphically to a subtree in $T$. 
\end{lemma}
\proof Connectivity of $F_{P}$ is clear, we need to verify the injectivity of the map. The problem is local and reduces to analyzing 
vertex-spaces of $\X$. The only way the map can fail to be injective is if there is an edge $e=[v,w]$ of $T$ and two distinct peripheral subspaces $H_{ei}, H_{ej}$ of $X_e$ which map to the same peripheral subspace $H_{vk}$ of $X_v$. But this contradicts the condition that $(X_e, {\mathcal H}_{e})\to (X_v, {\mathcal H}_v)$ is a morphism of relatively hyperbolic spaces.   \qed 

\medskip 
In particular, the trees $F_{P}$ are uniformly qi embedded in $X^\ell$.

\begin{rem}
There are natural projections $\theta_P: P\to F_P$  sending each nonempty intersection $P\cap X_v=H_{iv}$ to the 
 corresponding vertex $\nu_{iv}$ of $F_P$ and, for $H_{ie}=P\cap X_e, e=[v,w]$, sending the (double) mapping cylinder of 
 $(f_{ev}\cup f_{ew})|_{H_{ie}}$ to the edge $\eps_{ie}$ in such a way that each  interval $\{x\}\times [v,w], x\in X_e$ maps 
 linearly onto the edge   $\eps_{ie}$. The mapping cylinder of $\theta_P$ is then homeomorphic to the closure in $X^\ell$ 
 of the component of $X^\ell \setminus X$ containing $F_P$. We, therefore, will use the notation $Cyl(F_P, P)$ for this closure. This viewpoint of identifying $X^\ell$ with the space obtained from $X$ by attaching mapping cylinders of peripheral subsets $P\in {\mathcal P}$ to trees, is adapted from \cite{mahan-reeves} and \cite{mahan-pal}.   In line with the notion of radial segments in cones, we will refer to the projections of intervals $\{x\}\times [0,1]\subset P\times [0,1]$ to  $Cyl(F_P, P)$ as {\em radial segments} in  $Cyl(F_P, P)$. 
\end{rem}
 \medskip
The spaces ${X}^\ell$  and $\widehat{X}$  are related by a quotient map $\tau: {X}^\ell\to\widehat{X}$ which collapses each tree 
$F_{P}\in {\mathcal F}$ to a single point, the apex $a(P^\ell)\in \widehat{X}$ of the cone over the peripheral subspace $P$. 
See also the diagram in Section \ref{sec:comparison}. 
 
\medskip
We next describe a {\em relative flaring condition} for trees of relatively hyperbolic spaces. This condition consists of two parts. 

\medskip 
{\bf Part 1.} The first part of the relative flaring condition requires that the tree of hyperbolic spaces ${\X}^\ell$ satisfies one of the equivalent flaring conditions, equivalently, the space $X^\ell$ is hyperbolic (see Section \ref{sec:flare} for the detailed discussion). 

\medskip  
{\bf Part 2.} The $(\la,D)$-flaring:  

Let $\Pi=(\gamma_0, \gamma_1)$ be a pair of  $1$-sections\footnote{For each edge $e_j=[v_j, v_{j+1}]$ in 
$\llbracket u,w\rrbracket$, $\gamma_i(e_j)$ maps to $\gamma_i(v_j), \gamma_i(v_{j+1})$ under the maps 
$a(X^\ell_{e_j})\to a(X^\ell_{v_j})$ and $a(X^\ell_{e_j})\to a(X^\ell_{v_{j+1}})$.} of $\X^{\ell}$ over an interval 
$\llbracket u,w\rrbracket $ of length  
$\ge D$ in $T$, such that for each vertex $v\in V(\llbracket u,w\rrbracket)$, $\gamma_i(v)\in a(X^\ell_v)$, $i=0, 1$. 
 Then 
 $$
 \Pi_{max}\ge \la \Pi_0,
 $$  
where the distances are computed in electrified vertex-spaces.  
 
 \medskip 
 Mj and Reeves in \cite{mahan-reeves} prove the following relative form of the Bestvina-Feighn combination theorem: 
 
 \begin{thm}\label{mahan reeves thm}
 For every $\la>1, D\ge 2$,   if  $\X$ is a tree of relatively hyperbolic spaces satisfying the relative $(\la,D)$-flaring condition (both parts),  
 the space $X^P$ is hyperbolic.  
 \end{thm}

 We refer the reader to the papers by Sisto \cite{Sisto},  Gautero  \cite{Gautero-2016} and Dahmani \cite{MR2026551} 
 for related results.

\begin{lem}
Assuming that $\X$ is a tree of relatively hyperbolic spaces satisfying the relative flaring condition, 
the pair $({X}^\ell, {\mathcal F})$ is a relatively hyperbolic space. 
\end{lem}
\proof We will be using Theorem \ref{thm:H->RH}. Part 1 of the relative flaring condition (together with Theorem \ref{thm:mainBF}) 
implies that ${X}^\ell$ is hyperbolic. In view of  Lemma \ref{lem:qi-forest}, the peripheral subspaces $F_{P}\in {\mathcal F}$ are 
trees uniformly qi embedded in ${X}^\ell$. Lastly, Part 2 of the relative flaring condition implies that distinct peripheral subspaces in 
$ {\mathcal F}$ are uniformly pairwise cobounded. \qed 

\medskip
Since the pair $({X}^\ell, {\mathcal F})$ is relatively hyperbolic, we may perform the secondary cone-off construction, 
coning-off the peripheral subspaces $F_{P}\in {\mathcal F}$ of $X^\ell$. The result is a hyperbolic metric space 
$\widetilde{X}:={\mathcal E}({X}^\ell, {\mathcal F})$. There is a natural map 
$$
\theta: \widehat{X}={\mathcal E}({X}, {\mathcal P})\to \widetilde{X}={\mathcal E}({X}^\ell, {\mathcal F})
$$
 which is the identity on $X$ and sends the cone over each $P\in {\mathcal P}$ 
to the union of $Cyl(F_P,P)$ and $C(a_F,F)$ (cone over $F=F_{P}$)  
 so that each radial line segment connecting the apex $a=a(P^\ell)$ 
of $P^\ell=C(a, P)$ to a point $x\in P$, maps homeomorphically to the concatenation of the radial line segment in $C(a, F)$ 
connecting $a$ to a point $y\in F$ with the line radial segment in $C(F, P)$ connecting $y$ to $x$. 
The next lemma is immediate from the fact that each cone $P^\ell$ has unit diameter and for each $F\in {\mathcal F}$, 
the union of $Cyl(F_P, P)$ with the cone $C(a_F, F)$ has diameter $2$:  

\begin{lem}\label{lem:varphi} 
The  map $\theta$ is a quasiisometry. 
\end{lem}

 \medskip 
In the follow-up (to \cite{mahan-reeves}) paper \cite{mahan-pal}, Mj and Pal prove the following result regarding  existence of CT maps:

\begin{thm}\label{mahan pal thm}
For each vertex $v\in V(T)$, the inclusion map $X_v^h\to X^P$ admits a CT extension. 
\end{thm} 

The main result of this chapter is to generalize this theorem to subtrees $S\subset T$:

\begin{thm}\label{CT for relhyp tree}
For each subtree $S\subset T$, the inclusion map $X_S^P\to X^P$ admits a CT extension. 
\end{thm}

Note that the main tools in the  proof of Theorem \ref{mahan pal thm} in \cite{mahan-pal} were: 

\medskip 
(a) A construction of {\em derived ladders} in the induced tree of coned-off space ${X}^\ell\map T$ and 

(b) a construction of qi sections in $X$ lying inside these ladders. 

\medskip 
We include below  proofs of these results for the sake of completeness. 
(We will also need these results in order to prove Theorem \ref{CT for relhyp tree}.) 
However, we provide simplified proofs modulo the 
corresponding results for trees of hyperbolic metric spaces.

\subsection{Comparison of quasigeodesics}\label{sec:comparison} 

Let $\X=(\pi: X\to T)$ be a tree of relatively hyperbolic spaces. There are several spaces and, accordingly, several types of 
quasigeodesics associated with $\X$. In this section we discuss the relation between different types of quasigeodesics.  We begin, however, with a diagram describing spaces where these quasigeodesics live in. Recall that the tree of spaces $\X= (\pi: X\to T)$ gives rise to two other trees of spaces:
$$
\X^h= (\pi: X^h\to T) \leftsquigarrow \X \rightsquigarrow \X^\ell= (\pi: X^\ell\to T). 
$$
Both $\X^h$ and $\X^\ell$ are trees of hyperbolic spaces (satisfying Axiom {\bf H}) but only $\X^\ell$ satisfies the flaring condition. 
The following diagram describes the relation between five different spaces associated with $X$; the arrow $\theta$ is the  quasiisometry described in Lemma \ref{lem:varphi} and the map $\tau$ is the collapsing map from Section \ref{sec:rel trees}:  

 \newarrow{both}<--->
\newarrow {Into}C--->
\newarrow{Dboth}<...>

$$
\begin{diagram}
X^h &\lInto                        &X& \rInto & X^\ell& \rInto & \widetilde{X}=  
{\mathcal E}(X^\ell, {\mathcal F})\\
\dInto    &  & \id\dboth & &        & \tau\rdTo  & \theta~\dDboth \\
X^P= {\mathcal G}(X, {\mathcal P})  &\lInto  &        X & \rInto &         & &                    \widehat{X}= {\mathcal E}(X, {\mathcal P})        
\end{diagram} 
$$

\medskip 
We will discuss the following classes of quasigeodesics in these spaces: 

\begin{itemize}
\item Quasigeodesics in $X^P={\mathcal G}(X, {\mathcal P})$. 

\item Tight quasigeodesics  
 in $\widehat{X}={\mathcal E}(X, {\mathcal P})$.  

\item Tight quasigeodesics in $X^\ell$. 

\item Tight quasigeodesics in $\widetilde{X}=\widehat{X^\ell}={\mathcal E}(X^\ell, {\mathcal F})$.   

\end{itemize}

Our goal is to relate these quasigeodesics. 
We already know  
 that the first two types of quasigeodesics are uniformly Hausdorff-close to each other (when intersected with $X$, where the distance is computed via the metric $d_X$), see Theorem \ref{thm:hyperbolicity-of-electric-space} and Corollary 
 \ref{cor:hyperbolicity-of-electric-space}.  
In this section we prove that the same holds for the remaining types of quasigeodesics under a suitable tightness assumption on 
the quasigeodesics in $X^\ell$:

\begin{defn}
We say that a a continuous quasigeodesic in $X^\ell$ is {\em tight} if its intersection with every 
$\mathring{C}(F_P,P), P\in {\mathcal P}$, is either empty or is a concatenation of two radial geodesics with a geodesic in $F_P$. 
\end{defn}

By analogy with electrification of paths $\beta$ in $X^P$ (where the result is a tight path in $\widehat{X}$), we define  the  {\em partial electrification} $\beta^\ell$ of $\beta$ as follows:

\begin{defn} [Partial electrification] \label{defn:partial ell} 
Suppose that $\beta$ is a path in $X^P$ which is tight with respect to ${\mathcal P}$. For each pair of points $x, y$ in $\beta$ 
which   belong to some $P\in {\mathcal P}$ and such that the subpath $\beta(x,y)$ between points $x, y$ is contained in  $\mathring{P}^h$, except for the points $x, y$, 
we replace $\beta(x,y)$ with a geodesic in $\mathring{C}(F_P,P)$ connecting $x$ and $y$.  The resulting path $\beta^\ell$ is the partial electrification of $\beta$. 
\end{defn}

It follows from the definition that the path $\beta^\ell$ is tight in $X^\ell$. 

 Given a tight quasigeodesic $\beta$ in $X^\ell$, we can electrify  it with respect to the collection of peripheral subsets 
${\mathcal F}$ 
and obtain a tight quasigeodesic $\widetilde{\beta}= {\mathcal E}_{\mathcal F}(\beta)$ in $\widetilde{X}$, see Section \ref{sec:electrification}. Clearly, this defines an injective map 
$\beta\mapsto \widetilde{\beta}$ from the 
set of tight quasigeodesics in $X^\ell$ to those in $\widetilde{X}$, under which uniform quasigeodesics correspond to 
uniform quasigeodesics.

We also have the quasiisometry 
$$
\theta: {\mathcal E}(X, {\mathcal P})=\widehat{X}\to \widetilde{X}={\mathcal E}({X}^\ell, {\mathcal F})
$$
see Lemma \ref{lem:varphi}. This map  induces a bijection $\Theta$ between the sets of tight quasigeodesics in 
$\widehat{X}$ and quasigeodesics  of the form $\widetilde{\beta}$ in $\widetilde{X}$. Uniform quasigeodesics  again correspond to uniform quasigeodesics. Combining the two maps, we obtain a bijection
$$
\Phi: \alpha \mapsto \Theta(\alpha)= \widetilde{\beta}\mapsto \beta,$$ 
between the sets of tight quasigeodesics in ${\mathcal E}(X, {\mathcal P})$ and those in $X^\ell$. 
Moreover, the paths $\al$ and $\beta$ agree in $X$. 

We record these observations as in the following lemma:

\begin{lemma}\label{lem:Phi-qg}
There exists a bijection $\Phi$ between the sets of tight quasigeodesics in ${\mathcal E}(X, {\mathcal P})$ and those in $X^\ell$. 
Moreover, a tight path in ${\mathcal E}(X, {\mathcal P})$ is a uniform quasigeodesic if and only if the path $\beta=\Phi(\al)$ is 
in $X^\ell$. 
\end{lemma}

\medskip 
As an application, we obtain a result which appears as Lemma 1.21 in \cite{mahan-pal}: 

\begin{lemma}\label{lem:geodesic-comparison} 
Let $x, y$ be in $X$ and let $\al$ be a continuous $L$-quasigeodesic in $X^P$ between these points. 
Then there exists a tight  $L_{\ref{lem:geodesic-comparison}}(L)$-quasigeodesic 
$\beta$ in $X^\ell$ connecting $x$ and $y$ such that $\al\cap X= \beta\cap X$. 
\end{lemma} 
\proof As described in Section \ref{sec:electrification}, we convert $\al$ to 
a uniform tight quasigeodesic $\widehat{\al}$ in $\widehat{X}$, the electrification of $\al$. Then applying the map $\Phi$ as above, we obtain $\beta=\Phi(\al)$, a tight quasigeodesic in $X^\ell$ whose qi constant depends only on that of $\al$. By the construction. 
$\al\cap X= \beta\cap X$. \qed

\begin{cor}\label{cor:tight combing-1}
Any two points in $X\subset X^\ell$ are connected by a uniform tight quasigeodesic in $X^\ell$. 
\end{cor}

\medskip
Suppose that $\beta_1, \beta_2$ are tight $L$-quasigeodesics in $X^\ell$ with the same end-points in $X$. Applying the 
inverse bijection $\Phi^{-1}$ to $\beta_1, \beta_2$ we obtain tight $L'$-quasigeodesics $\al_1, \al_2$ 
in $\widehat{X}$, again connecting the same points  in $X$ and such that $\al_i\cap X=\beta_i\cap X, i=1,2$. Since, by 
Theorem \ref{thm:Pal-Kumar} the paths $\al_1, \al_2$ uniformly track each other in $X$, we obtain:

\begin{cor}\label{cor:tracking-2}
Any two tight quasigeodesics in $X^\ell$ connecting points of $X$ uniformly track each other in $X$. 
\end{cor}

\subsection{Ladders in trees of relatively hyperbolic spaces} \label{sec:relative ladders}

By Lemma \ref{lemma: extended relhyp tree}, we have a tree of hyperbolic spaces
$\X^h=(\pi^h: X^h\map T)$ such that edge-spaces are uniformly quasiisometrically embedded in vertex spaces.  

Hence, given $x, y \in X_u\subset X^h_u$ we  construct 
a $(K,D,E)$-ladder $\LL^h\subset \X^h$ of the geodesic segment $L^h_u=[xy]_{X_u^h}$,  
 as we have done in Chapter \ref{ch:4 classes}. (Note that flaring conditions were not used in this construction.) 
 
 For each vertex $v$ (resp. edge  $e$) in $\pi(L^h)\subset T$ we have the oriented vertical geodesic segment $L^h_v=[x_v y_v]_{X^h_v}$ (resp. $L^h_e=[x_e y_e]_{X^h_e}$) of the ladder $\LL^h$. Coning-off the subsegments of $L_v$ (resp. $L_e$) contained in the peripheral horoballs of  $X_v^h$ (resp. $X_e^h$) we obtain a collection of uniform quasigeodesics  in the electrified vertex/edge spaces ${X}^\ell_v$ (resp.  ${X}^\ell_e$), see Theorem \ref{thm:hyperbolicity-of-electric-space}(2).  We let 
 ${L}^\ell_v$ (resp. ${L}^\ell_e$) denote the geodesics in  ${X}^\ell_v$ (resp.  ${X}^\ell_e$) connecting the end-points of the above quasigeodesics.  
 
\begin{lem}\label{lem:hat-ladder}
The collection of segments ${L}^\ell_v$ and ${L}^\ell_e$ defines a $(\hat{K}, \hat{D}, \hat{E})$-ladder ${\LL}^\ell$ in 
${\X}^\ell$. 
\end{lem} 
 \proof The proof is based on Lemma \ref{lem:E-ladder-structure}. The collection of segments ${L}^\ell_v$ and ${L}^\ell_e$ satisfies 
 the assumptions of  Lemma \ref{lem:E-ladder-structure} because $\LL^h$ is a ladder and  because coning-off of quasigeodesics 
 coarsely commutes with the nearest-point projections, see Lemma \ref{lem:coning-projection-commutes} as well as Remark \ref{rem:semico}(viii). \qed

\medskip
We let ${L}^\ell, L^h$ denote the total spaces of the ladders $\LL^\ell$ and $\LL^h$ respectively.  We will say that the ladder 
$\LL^\ell$ is {\em derived} from the ladder $\LL^h$. 

\begin{rem}
The ladder construction given in \cite{mahan-pal} is a bit more complicated. 
\end{rem}

By the construction, each ladder is a tree of relatively hyperbolic spaces with the total space $L$ and the peripheral structure 
${\mathcal P}_L$ (given by the pull-back of the peripheral structure ${\mathcal P}$ of $X$), 
and the inclusion map $(L, {\mathcal P}_L) \to (X, {\mathcal P})$ is a morphism of relatively hyperbolic spaces (we will need only that distinct peripheral subsets map to distinct ones).   

 \medskip 
Corollary \ref{cor:ladder-retraction} (the existence of coarse Lipschitz retractions 
to ladders in hyperbolic trees of spaces),  
combined with Lemma \ref{lem:hat-ladder} implies: 

\begin{cor} \label{cor:mahan-pal}
There is a coarsely Lipschitz retraction ${X}^\ell\map {L}^\ell$ with 
Lipschitz constant depending only on the parameters $K, D, E$.
\end{cor}

Applying Corollary \ref{cor:tight combing-1}, we obtain: 

\begin{cor}\label{cor:tight combing-2}
Any two points in $L\subset L^\ell$ are connected by a uniform tight quasigeodesic in $L^\ell$. 
\end{cor}

\medskip
{\bf QI sections in ladders.} We next discuss the  construction in \cite{mahan-pal} 
of {\em vertical quasigeodesic rays} contained in ladders. We assume that 
$$
\LL^h=\LL^h([x_u y_u]_{X_u^h})$$
is a $(K,D,E)$-ladder in $\X^h$ and $\LL^\ell$ is the ladder in ${\X}^\ell$ derived  from it. 

\begin{lemma}\label{lem:mahan-pal}
Given a vertex $v\in \pi(L^h)$ and a point $z_v\in X_v\cap {L}^\ell_v$, there is a 
$k$-qi section $\sigma: \llbracket u,v\rrbracket \map X$ of $\pi: X\map T$ 
such that $\sigma(\llbracket u,v\rrbracket)\subset {L}^\ell$, 
$\sigma(v)=z_v$. Here, $k$ depends only on $K$ and $E$. 
\end{lemma} 
\proof Arguing inductively, it is clear that it suffices to prove the lemma in the case when $u, v$ span an edge $e=[u,v]$ of $T$. We will use the fact that $\LL^\ell$ is derived from the $K$-ladder $\LL^h$. Thus, there exists a point $z\in L^h_u$ within distance $K$ from $z_v$ (the distance is measured in $X_{uv}^h$). Since $z_v$ is in $X_v$, the incidence maps of $\X$ are morphisms 
of relatively hyperbolic spaces,  there exists a point $z_u\in L^h_u\cap X_u$ within uniformly bounded distance from $z_v$ and, hence, within uniformly bounded distance from $z$ in $X_{uv}$. Hence, we set $\sigma(u):= z_u$. \qed

\section{Cannon-Thurston maps for trees of relatively hyperbolic spaces}

{\em Proof of Theorem \ref{CT for relhyp tree}:} In the proof we shall use the following criterion for the existence of CT maps,
which appears as Lemma 1.29 in  \cite{mahan-pal}. 
Recall that we have a subtree $\X^\ell_S= (X^\ell\to S)$ in  the tree of hyperbolic spaces $\X^\ell=(X^\ell\to T)$. 
We fix a point $x_0\in X_{v_0}$, where $v_0\in V(S)$. 

\begin{lem}\label{prop:rel-CT-1} 
A CT map for $X^P_S\map X^P$ exists provided that the following condition holds:

For each $k\geq 1$ and $M\geq 0$, there is $N\geq 0$ 
such that for all $x,y\in X_S$, if $\gamma^\ell\subset {X}^\ell, \gamma^\ell_S\subset {X}^\ell_S$ are tight  
$k$-quasigeodesics  joining $x,y$, then $d_X(x_0, \gamma^\ell\cap X)\leq M$ implies 
$d_{X_S}(x_0, \gamma^\ell_S\cap X_S)\leq N$. 
\end{lem}
\proof Our proof closely follows the one of \cite[Lemma 1.29]{mahan-pal}. 
We will verify that the conditions of Proposition \ref{prop:rel-CT} are satisfied. 
The main difference with the setting of the lemma is  that Proposition \ref{prop:rel-CT} 
was stated in terms of coned-off quasigeodesics  
$\widehat\ga, \widehat\ga_S$ in $\widehat{X}= {\mathcal E}(X, {\mathcal P})$, 
$\widehat{X}_S= {\mathcal E}(X_S, {\mathcal P}_S)$ respectively, connecting $x, y$.  Here ${\mathcal P}_S$ is the collection 
of intersections $P\cap X_S$, $P\in {\mathcal P}$.  

However, according to Lemma \ref{lem:geodesic-comparison}, the Hausdorff distances (computed in $X$ and $X_S$ respectively) 
between $\ga^\ell\cap X,  \widehat\ga\cap X$ and between $\ga^\ell_S\cap X_S, \widehat\ga_S\cap X_S$ 
are uniformly bounded, with bounds depending only on $K$. 
With this in mind, Proposition \ref{prop:rel-CT} applies and lemma follows. \qed

\medskip 
 We now proceed proving the theorem. Let $x, y\in X_S$ be arbitrary points and  
let $\gamma^\ell= [xy]_{X^\ell}$ be a uniform tight quasigeodesic in ${X}^\ell$ joining them (see Corollary \ref{cor:tight combing-1}). 
Suppose $z\in \gamma^\ell\cap X$ is a point such that $d_X(x_0,z)\leq D$ for some $D\geq 0$.

We apply the cut-and-replace operation (see Definition \ref{defn:detour}) 
to the tree of spaces $\X^\ell$ 
and the quasigeodesic $\gamma^\ell$, transforming it to a path $\widehat{\gamma^\ell}$ in $X^\ell_S$. 
By Theorem \ref{thm:cut-paste},  the path $\widehat{\gamma^\ell}$ in $X^\ell_S$ is a uniform quasigeodesic. Tightness of $\gamma^\ell$ implies that  of $\widehat{\gamma^\ell}$. 

\medskip 
{\bf Case 1.} Suppose that $z\in X_S$. By the construction of $\widehat{\gamma^\ell}$, the point 
$z$ lies on $\widehat{\gamma^\ell}$. By 
Corollary \ref{cor:tracking-2}, any two tight uniform quasigeodesics in $X_S^\ell$ uniformly  track each other in $X_S$. Hence, 
each tight uniform quasigeodesic $\gamma^\ell_S$ in $X_S^\ell$ (connecting $x$ and $y$) passes within uniformly bounded distance (in terms of the metric of $X_S$) from the point $z$.   Thus, the implication required by Proposition \ref{prop:rel-CT-1} holds and we are done in this case. 

\medskip 
{\bf Case 2.} Suppose that $z\notin X_S$. Then there is a vertex $v\in V(S)$ and a component 
$\gamma^\ell_1$ of $\gamma^\ell\setminus {X}^\ell_S$, such that $z\in \gamma_1^\ell$ and the end-points $x_1, y_1$ 
of $\gamma_1^\ell$ belong to ${X}^\ell_v$. 

\medskip 
{\bf Subcase 2.1.} We first consider the subcase when $x_1, y_1$ both belong to $X_v$. 
Let $T_1$ be the smallest subtree of $T$ such that $\gamma^\ell_1$ is contained in ${X}^\ell_{T_1}$. In other words, 
$T_1$ is the span of a component of $T\setminus S$ and the vertex $v$. 
In the tree of spaces ${\X}^\ell_{T_1}$ we  construct a ladder  
$$
{\LL}^\ell={\LL}^\ell([x_1 y_1]_{{X}^\ell_{v}})$$ 
with the total space ${L}^\ell$, as described in Section \ref{sec:relative ladders}.

Let $\bar\gamma_1^\ell$ be a uniform tight quasigeodesic joining $x_1, y_1$ in ${L}^\ell$. 
Since the ladder ${L}^\ell$ is uniformly qi embedded in ${X}^\ell_{T_1}$ (see Corollary \ref{cor:ladder-retraction}), 
$\bar\gamma_1^\ell$ is also a uniform quasigeodesic in ${X}^\ell_{T_1}$. Moreover, since the inclusion $\LL^\ell\to {\X}^\ell_{T_1}$ corresponds to a morphism of relatively hyperbolic spaces 
$$
(L, {\mathcal P}_L)\to (X_{T_1}, {\mathcal P}_{T_1}),$$ 
the path $\bar\gamma_1^\ell$ is tight in ${X}^\ell_{T_1}$.

Applying Corollary \ref{cor:tracking-2} again, we obtain that $z\in \gamma^\ell_1\cap X_{T_1}$ is $C$-close 
(in terms of the metric of $X_{T_1}$) 
 to a point $\bar{z}\in \bar\gamma_1^\ell\cap X_w$ for some vertex  $w\in V(T_1)$, where $C$ depends only on qi constants of the original quasigeodesic $\gamma^\ell_1$. Recall that $d_X(x_0,z)\leq D$. Then $d_X(x_0, \bar{z})\le D_1:= D+C$. 

According to Lemma \ref{lem:mahan-pal}, there exists a $k$-qi section $\sigma$ over 
$\llbracket v, w\rrbracket$ with image  in ${L}^\ell\cap X$, such that  $\sigma(w)=\bar{z}$ and 
$$
\sigma(v)=z_1\in [x_1 y_1]_{{X}^\ell_{v}}.$$
Since $z_1$ belongs to $X_S$,  $d_T(v,w)\le D_1$ and $d_X(x_0, \bar{z})\le D_1$, it follows that 
$$
d_X(x_0,z_1)\le kD_1+ D_1.
$$
Thus, we can apply the same reasoning as in the case when $z$ is in $X_S$, to conclude that for every  uniform tight 
quasigeodesic $\gamma_S^\ell$ in ${X}^\ell_S$ connecting $x, y$ such that the minimal distance (in $X_S$) from 
$\gamma_S^\ell\cap X_S$ to $z$ is  uniformly bounded. 

\medskip 
{\bf Subcase 2.2.} Suppose that one of the two points $x_1, y_1$ is in $X_v$ and the other is not. After relabeling, we can assume that $x_1\in X_v$ and $y_1\notin X_v$.

The (tight) geodesic $[x_1 y_1]_{{X}^\ell_{v}}$ 
contains a maximal subpath $[y'_1 y_1]_{{X}^\ell_{v_1}}$ of length $\le 1$ contained in $X^\ell_{v}\setminus X_v$. 
The concatenation $\gamma_1'$ of this path and $\gamma^\ell_1$ is still tight, $\gamma_1'$ is again a uniform quasigeodesic 
in ${\X}^\ell_{T_1}$, and both of its endpoints are in $X_{v_1}$. Now, we repeat the argument in 
Subcase 2.1 with respect to the pair of points $x_1, y'_1\in X_v$ and the uniform tight quasigeodesic  $\gamma_1'$ connecting them.

\medskip 
{\bf Subcase 2.3.} The subcase when both points $x_1, y_1$ are not in $X_v$ but $[x_1 y_1]_{{X}^\ell_{v}}\cap X_v\ne \emptyset$ 
is similar to Subcase 2.2 and we leave it to the reader.

\medskip 
{\bf Subcase 2.4.}  $[x_1 y_1]_{{X}^\ell_{v}}\cap X_v= \emptyset$, which implies that the points $x_1, y_1$ are within unit distance 
from each other in ${X}^\ell_{v}$. Tightness of the path $\ga_1^\ell$ then implies that $\ga^\ell_1$ is disjoint from $X$, which contradicts the assumption that $z\in \ga^\ell_1\cap X$. \qed 
 
 \section{Cannon-Thurston laminations for trees of relatively hyperbolic spaces} 

Since the inclusion map $X_S^P\to X^P$ admits a CT-map, one can also define its Cannon-Thurston lamination. Unlike 
Chapter \ref{ch:CT} (specifically, Sections  \ref{sec:CTfibers}, \ref{sec:Boundary flows and CT laminations} and \ref{sec:CT-lamination}) 
where these laminations were discussed in great detail in the absolute case, here we limit our discussion to a  
(weak) analogue of Theorem \ref{thm:CT-fibers}(1), relating these CT-laminations to that of the inclusion maps $X_v^h\to X^P$ 
(Theorem \ref{thm:CT-lam for relhyp tree} below).  

\medskip 
From now on, we  assume, as we did in Theorem \ref{thm:CT-fibers}, 
that $X$ is a proper metric space. Let $S\subset T$ be a subtree. 
We then have the inclusion map $X_S^P\to X^P$ that, according to 
Theorem \ref{CT for relhyp tree}, admits a CT-map $\D_{X_S^P,X^P}$. 

\begin{thm}\label{thm:CT-lam for relhyp tree}
Suppose that $\xi^\pm$ are distinct points in $\geo X^P_S$ such that 
$$\eta=\D_{X_S^P,X^P}(\xi^-)=\D_{X_S^P,X^P}(\xi^+).$$ 

1. Then there exists a vertex-space 
$X_v\subset X_S$ and a complete geodesic $\al: \RR\to X^h_v$, 
asymptotic to $\xi_\pm$ in $X_S^P$ and asymptotic (in $X^h_v$) to points $\xi_v^\pm\in \geo X_v^h$, such that 
$$
\D_{X_v^h,X^P}(\xi_v^+)= \D_{X_v^h,X^P}(\xi_v^-)= \eta.   
$$

2. Moreover, $\al$ is a uniform quasigeodesic in $X^P_S$. 
\end{thm}

\begin{rem}
In this situation we necessarily have 
$$
\D_{X_v^h,X_S^P}(\xi_v^\pm)= \xi_\pm.   
$$
\end{rem}

\proof The proof follows the arguments of Proposition \ref{prop:CT-fibers1}. Our first task is to modify the {\em cut-and-replace} procedure (Definition \ref{defn:detour}), used in the proof of 
Proposition \ref{prop:CT-fibers1}. As before, we will be identifying spaces $X_v^h$ with their images in $X^P$, $v\in V(T)$; ditto $X_S^P$.  The modification of a path $\beta$ in $X^P$ to a path $\invbreve{\beta}$ in $X_S^P$ will be done in two steps. 

\medskip 
{\bf Step 1.} The first step simply repeats what is done in Definition \ref{defn:detour} (except that $X^P$ and $X_S^P$ are not exactly total spaces of trees of spaces!): We identify {\em primary detour subpaths $\zeta_v$} of $\beta$, 
connecting points $x_v, y_v$ of vertex-spaces $X_v^h$, $v\in V(S)$, and, apart from these end-points, lying outside of $X_S^P$. 
We then replace each detour subpath with a geodesic $[x_v y_v]_{X_v^h}$ in $X^h_v$.  These geodesic segments are the {\em primary replacement subpaths} for $\beta$. As in  Definition \ref{defn:detour}, we refer to the resulting path as $\beta_S$: Its image entirely lies in $X_S^P$. The trouble is that, unlike the absolute case, even if $\beta$ is uniformly quasigeodesic, the paths $\beta_S$ {\em are priori} are not  even uniformly proper in $X_S^P$. 

\medskip
{\bf Step 2.} We then modify $\beta_S$ with respect to the collection of peripheral subsets ${\mathcal P}_S$ of $X_S$ as follows. 
We define {\em secondary} detour subpaths of $\beta_S$ as subpaths $\zeta_P$ in $\beta_S$ connecting points of peripheral subsets $P\in {\mathcal P}_S$ of $X_S$ and, besides those points, lying entirely inside the open peripheral horoballs $\mathring{P}^h\subset X^P_S$.  
The path 
$$
\invbreve{\beta}= {\mathcal G}_{X_S,{\mathcal P}_S}(\beta_S)
$$
the {\em hyperbolization} of $\beta_S$, is obtained by replacing each secondary detour subpath $\zeta_P$ with a geodesic in $P^h$ connecting the end-points of $\zeta_P$. 

\medskip 
Our next goal is to relate quasigeodesic properties of $\beta$ to that of $\invbreve{\beta}$. 
We will be using the notion of {\em tight paths}, Definition 
\ref{defn:tight-path}. The following lemma is clear: 

\begin{lem}
Suppose that $\beta$ in $X^P$ is tight with respect to ${\mathcal P}$. 

1. For each subtree $S\subset T$, the  
cut-and-replace path $\beta_S$ is tight with respect to ${\mathcal P}_S$, and so is  the path $\invbreve{\beta}$. 

2. $(\beta_S)^\ell= (\beta^\ell)_S$, where   $\beta_S$ is the result of application Step 1 to $\beta$ as above and 
$(\beta_S)^\ell$ is its partial electrification (a tight path in $X_S^\ell$, see Definition \ref{defn:partial ell}),  while 
$(\beta^\ell)_S$ is the result of application of the cut-and-replace procedure from Definition \ref{defn:detour} to the path $\beta^\ell$  (the partial electrification of $\beta$) 
in the tree of spaces $\X^\ell=(X^\ell\to T)$. 
\end{lem}

In view of Part 2 of this lemma, we will use the notation $\beta^\ell_S$ for $(\beta_S)^\ell= (\beta^\ell)_S$.

\medskip 
Recall that the map $\Phi$ defined in Section \ref{sec:comparison} establishes a bijection between tight 
paths in ${\mathcal E}(X, {\mathcal P})$ and those in $X^\ell$; same for the corresponding map $\Phi_S$ for 
the electric space ${\mathcal E}(X_S, {\mathcal P}_S)$  and $X_S^\ell$. The following lemma is again a direct consequence of the definitions:

\begin{lem}
For every path $\beta$ in $X^P$, tight with respect to ${\mathcal P}$, the path  $\invbreve{\beta}$ can be described as follows. Set $\hat\beta_S:= \Phi_S^{-1}(\beta_S^\ell)$. Then $\invbreve{\beta}$ equals the path obtained via hyperbolization of $\hat\beta_S$ 
with respect to ${\mathcal P}_S$, i.e. 
$$
\invbreve{\beta}= {\mathcal G}_{X_S, {\mathcal P}_S}(\hat\beta_S). 
$$
\end{lem}

This lemma has an important consequence:

\begin{cor}\label{cor:inbreve} 
If $\beta$ is a tight $L$-quasigeodesic in $X$, then $\invbreve{\beta}$ is an $L_{\ref{cor:inbreve}}(L)$-
quasigeodesic in $X_S^P$.  
\end{cor}
\proof First of all, according to Theorem \ref{thm:hyperbolicity-of-electric-space}(2), the path $\beta^\ell$ is a uniform (in terms of $L$) quasigeodesic in $X^\ell$. Thus, by Theorem \ref{thm:cut-paste}, so is the path $\beta_S^\ell= (\beta^\ell)_S$, in $X_S^\ell$. By Lemma \ref{lem:Phi-qg}, the path  $\hat\beta_S$ is also. 
Lastly, by Lemma \ref{lem:Dahmani-Mj}, it follows that the path 
$\invbreve{\beta}= {\mathcal G}_{X_S, {\mathcal P}_S}(\hat\beta_S)$ is also a uniform quasigeodesic in  $X_S^P$.  \qed

\medskip 

With these preliminaries out of the way, 
we can now proceed with the proof of Theorem \ref{thm:CT-lam for relhyp tree}. Fix a base-point $x_0\in X_S$. 
Let $(x^\pm_n)$ denote  sequences in $X_S$ 
 converging in $X_S^P \cup \geo X_S^P$ to the points $\xi^\pm$ respectively, see Lemma \ref{lemma:relative-limits}. 
We let $\beta_{n}$ denote a sequence of tight (with respect to ${\mathcal P}$)  uniformly quasigeodesic paths in $X^P$ connecting the points $x^-_n, x^+_n$. According to 
Corollary \ref{cor:inbreve}, the paths $\invbreve\beta_{n}$ are uniformly quasigeodesic. Thus, since, by the assumption, $\xi^+\ne \xi^-$, it follows that there is a constant $C$ such that $d(x_0, \invbreve\beta_{n})\le C$. We let $z_n$ denote a point in 
$\invbreve\beta_{n}$ within distance $C$ from $x_0$. Without loss of generality, we may assume that these points lie in the images in $X_S^P$ of spaces $X_{v_n}^h$, $v_n\in S$. 
Since $X$ is assumed to be proper, so is $X_S$. Therefore, after extraction, we may assume that 
all points $z_n$  lie in $X^h_v$ for some vertex $v\in V(S)$.

\medskip 
{\bf Case 1:}  All points $z_n$ lie in $X_v$. Since the sequence of distances $d(x_0, \beta_n)$ diverges to $\infty$, for all but finitely many $n$'s, the points $z_n$ lie in  secondary replacement segments 
$$
[z^+_n z^-_n]_{X_{v}}= \invbreve\zeta_{vn}\subset \invbreve\beta_{n}\cap X_{v}.
$$
Recall that, by the construction of $\invbreve\beta_{n}$, for each $n$ we also have a  primary replacement segment $\zeta_{vn}=[y^+_n y^-_n]_{X^h_{v}}$ of $\beta_n$, 
containing $[z^+_n z^-_n]_{X_v}$, where $y^\pm_n$ lie on $\beta_n$. 
Since $d(x_0, \beta_n)\to\infty$, 
$$
\lim_{n\to\infty} d(x_0, y^\pm_n)=\infty. 
$$
By the properness of $X_v$, since $z_n$ is in $\zeta_{vn}$ and $d(x_0, z_n)\le C$,  after further extraction, the sequence of geodesics  
$[y^+_n y^-_n]_{X^h_{v}}$ converges to  a biinfinite geodesic $\al$ in $X^h_v$ connecting  points $\xi^\pm_v\in \geo X^h_v$.  
Using again the assumption that  $d(x_0, \beta_n)\to\infty$, and the points $x_n^\pm, y_n^\pm$ all lie on the uniform 
quasigeodesic $\beta_n$, we see that 
$$
\partial_{X^h_v,X^P}(\xi_v^+)= \partial_{X^h_v,X^P}(\xi_v^-)= \partial_{X_S^P, X^P}(\xi^+)= \partial_{X_S^P, X^P}(\xi^-). 
$$
This proves the first claim of the theorem in Case 1. 

\medskip 
It remains to show that $\al$ (or, for this matter, any uniform quasigeodesic in $X_v^h$ asymptotic to $\xi_v^\pm$) is a uniform quasigeodesic in $X^P_S$. Note that, since paths $\invbreve\beta_{n}$ are uniform quasigeodesics in $X^P_S$, the same holds for the subpaths $[z^-_n z^+_n]_{X_v}\subset \invbreve\beta_{n}$.

\medskip
{\bf Subcase 1a:} $\lim_{n\to\infty} d(x_0, z^\pm_n)=\infty$. Then the sequence of subsegments $[z^-_n z^+_n]_{X_v}$ 
(uniformly quasigeodesic in $X^P_S$) subconverges to a geodesic $\al_v$ in $X_v$ also asymptotic to $\xi_v^\pm$; it follows that $\al_v$ is a uniform quasigeodesic in $X_S^P$. 

\medskip 
{\bf Subcase 1b:} There exists a constant $D$ such that $\lim_{n\to\infty} d(x_0, z^+_n)=\infty$ and $d(x_0, z^-_n)\le D$ for all $n$. 
By properness, after further extraction, we can assume that all points $z_n^-$ lie on a fixed peripheral subset $P_v\in {\mathcal P}_{X_v}$ 
of $X_v$. We then consider the sequence of segments $[y_n^- z_n^+]_{X_v}$ instead of the segments 
$[z^-_n z^+_n]_{X_v}$  used in Case 1. 

Each segment $[y_n^- z_n^+]_{X_v}$ is a concatenation of a vertical geodesic segment 
$$
[y_n^- z_n^-]_{X_v}\subset P_v^h$$
 and the geodesic segment $[z_n^- z_n^+]_{X_v}$. Since the path $\beta_{nS}$ is tight (with respect to ${\mathcal P}_S$), its subsegment $[y_n^- z_n^+]_{X_v}$ is also tight in $X_v^h$ with respect to its peripheral structure ${\mathcal P}_{v}$. In particular, the subsegments $[z_n^- z_n^+]_{X_v}$ are all disjoint from the open peripheral horoball $\mathring{P}_v^h$. In particular, 
$[z_n^- z_n^+]_{X_v}$ is disjoint from the open peripheral horoball $\mathring{P}_S^h\subset X^P_S$ containing 
$\mathring{P}_v^h$. Using the fact that 
$[z_n^- z_n^+]_{X_v}$ and $[y_n^- z_n^+]_{X_v}$ are both uniform quasigeodesic in $X^P_S$, we apply Lemma \ref{lem:three geodesics} to conclude that the segments $[y_n^- z_n^+]_{X_v}$ are uniformly quasigeodesic in  $X^P_S$ as well. Thus, their limit 
$\al_v$, a uniform quasigeodesic in $X^h_v$ asymptotic to $\xi^\pm_v$, is also a uniform quasigeodesic in   $X^P_S$. 

\medskip
{\bf Subcase 1c:} There exists a constant $D$ such that  $d(x_0, z^\pm_n)\le D$ for all $n$. The proof is similar to the subcase 1b and we give only a sketch. We consider the sequence of geodesic segments $[y_n^- y_n^+]_{X_v}$ and break each of these as a concatenation of three geodesic segments, two of which are contained in distinct horoballs $P^h_{v\pm}$ and one has uniformly bounded length. We again apply Lemma \ref{lem:three geodesics} to conclude that the segments $[y_n^- y_n^+]_{X_v}$ are uniformly quasigeodesic in $X_S^P$. 

\medskip
{\bf Case 2.} We now assume that none of the points $z_n$ belongs to $X_S$. 
Then, after further extraction, each $z_n$ lies in a peripheral horoball $P^h$ for some 
$P\in {\mathcal P}_S$. Hence, $z_n$ belongs to  a geodesic $[p_n q_n]_{P^h}\subset \invbreve\beta_{n}$, 
where $p_n\in P_{v_n}$, $q_n\in P_{w_n}\in {\mathcal P}_{X_{w_n}}$ and 
$$
P_{v_n}= X_{v_n}\cap P, \quad P_{w_n}= X_{w_n}\cap P. 
$$
By the description of such geodesics $[p_n q_n]_{P^h}$, up to a uniformly bounded error, the path  
$[p_n q_n]_{P^h}$ is the concatenation of two vertical geodesics segments $[p_np'_n]_{X^h_{v_n}}$, 
$[q_nq'_n]_{X^h_{w_n}}$  and a unit horizontal segment in $P^h$ connecting $p'_n, q'_n$. Thus, either 
$d_{X_S^P}(z_n, p_n)\le C$ or $d_{X_S^P}(z_n, q_n)\le C$. Accordingly, 
$d_{X_S^P}(x_0, p_n)\le 2C$ or $d_{X_S^P}(x_0, q_n)\le 2C$.  
Since $p_n, q_n$ are in $X_S$, Case 2 is reduced to Case 1.

This concludes the proof of the theorem. \qed

\bibliographystyle{alpha}


\bibliography{tree-bibliography.bib}

\chapter*{List of symbols} 

\newglossaryentry{latex}
{
    name=$\Delta$,
    description={geodesic triangles in metric spaces}
}

\printglossaries


\begin{tabular}{cp{0.6\textwidth}}
$\Delta$ & geodesic triangles in metric spaces\\
  $L, \eps$ 
  & coarse Lipschitz constants\\
  $D, E, \eps, R, r$& bounds on distances\\
  $R, r$ & radii of tubular neighborhoods\\
  $\ell(\al)$, $\length(\al)$ & length of a path $\al$\\
  $c, \al, \beta, \gamma, \phi, \psi$ & a path\\
  $\hat{c}= c_S$ & path obtained from $c$ via the cut-and-replace procedure, Definition \ref{defn:detour} \\
  $\cev{c}$ & reversed path of $c$\\
  $c_1\star c_2$ &  concatenation of paths $c_1, c_2$\\
  $D_0$ & Equation \eqref{eq:D0}\\
  $K_0$ & Notation \ref{not:K0}\\
  $\kappa, k, K$ & quasiisometry and flow constants \\
  Function $K'$ & Lemma \ref{lem:E-ladder-structure}\\
  Functions $K^\vee$ and $K^\wedge$ & Equation \eqref{eq:vee}\\
  $\la$ & quasiconvexity constant \\
  $\rho, \mu, \nu$ & retractions \\
  $f, g, h$ & maps\\
    $P$, $P_{X,Y}$ & nearest-point projection\\
    $\bar{P}$ & the modified projection to a tripod\\
  $\bar{x}$ & nearest-point projection of $x$\\
  $\eta$ & a distortion function\\
    $\eta_0$ & the distortion function of vertex-spaces in a tree of spaces $X$ \\
  $V(T)$ & the vertex set of a tree $T$\\
  $e=[u,v]$ & the edge of $T$ connecting vertices $u, v$\\
  $\dot{e}$ & an edge in a tree minus its vertices \\
  $E(T)$ & the edge set of a tree $T$\\
  $\ga_x$ & qi sections over an interval or a subtree, passing through $x$\\
  $\Hd$ & Hausdorff distance\\
   \end{tabular}\\

\begin{tabular}{cp{0.6\textwidth}}  
    $\llbracket u, v \rrbracket$ or $uv$ & the geodesic segment (interval) connecting vertices $u, v\in V(T)$ in a tree\\
    $ \rrbracket u, v \rrbracket$, $\rrbracket u,v\llbracket$, $\llbracket u, v\llbracket $ & open and half-open intervals in $T$\\
    $S(u,r)$ & the sphere of radius $r$ and center $u$\\ 
      $B(u,r)$ & the closed ball of radius $r$ and center $u$\\
      $N_R(A)$ & closed $R$-neighborhood of $A$\\
      $b=\barycenter(\Delta)$ & center of a tripod in a tree\\
  ${\mathcal X}:= \coprod_{v\in V(T)} X_v$ & union of vertex-spaces of a tree of spaces\\
  ${\mathcal Fl}_K(\cdot)$, $Fl_K(\cdot)$, ${\mathfrak Fl}_K(\cdot)$ & flow-space (as a tree of spaces), the total space and the intersection with ${\mathcal X}$ respectively\\
  $Fl(X_v)$ & ideal boundary flow\\
  $\mathfrak{L}, \LL_{K,D,E}$ & a ladder, as a tree of spaces and intersection with ${\mathcal X}$ respectively\\
  ${\mathfrak A}, \AA_{K,C}$ & carpets\\
  $\La(f), \La(Y,X)$ & CT-lamination\\
  $\geo (Y, X)$ & relative ideal boundary of $Y$ in $X$\\
  $\delta_0$ & hyperbolicity constant of vertex spaces (total space and intersection with ${\mathcal X}$)\\
  $\delta'_0$ & hyperbolicity constant of spaces $X_{uv}$, $[u,v]\in E(T)$\\
  $\la_0$ & the quasiconvexity constant of $X_{eu}\subset X_u$, $e=[u,v]\in E(T)$\\
  $\la'_0$ & the quasiconvexity constant for the inclusion maps $X_u\to X_{uv}$, $[u,v]\in E(T)$\\
  $L_0$ & the quasiisometry constant for the incidence maps $f_{eu}: X_e\to X_v$\\
  $L'_0$ & the quasiisometry constant for the inclusion maps $X_v \to X_{uv}$, $e=[u,v]\in E(T)$\\
  $N^{fib}_R$ & closed fiberwise $R$-neighborhood\\
  $N^e_R$ & closed neighborhood in $X_{uv}$, where $[u,v]=e\in E(T)$\\
  $\Pi=(\ga_0,\ga_1)$& a pair of $K$-qi sections over the same interval in $T$\\
  $\Pi_0$& the girth of $\Pi$\\
  $\Pi_{\max}$& the maximal separation of the ends of $\Pi$\\
  ${\mathfrak X}=(\pi: X\to T)$ & a tree of spaces\\
  $X_S:=\pi^{-1}(S)$ & the total space of a subtree of spaces over $S\subset T$\\
  $(Y^\ell, d^\ell)={\mathcal E}(Y, {\mathcal H})$ & the electrified space \\
  $(Y^h, d^h)= {\mathcal G}(Y, {\mathcal H})$ & the horoballification\\
  $a(Y^\ell)$ & cone-locus\\
  \end{tabular}\\

\printindex

\end{document}